\documentclass[a4paper,twoside]{report}

\usepackage{fancyhdr}
\usepackage{graphicx}
\usepackage{german}
\usepackage[english]{babel}
\usepackage{xcolor}
\usepackage{color}
\usepackage[utf8]{inputenc}
\usepackage{amscd,amssymb,amsmath,amsthm,bbm}
\usepackage[mathscr]{euscript}
\usepackage{caption}

\usepackage[breaklinks=true,colorlinks=true,
linkcolor=black,urlcolor=black,citecolor=black,
bookmarks=true,bookmarksopenlevel=2]{hyperref}

\usepackage{tikz}
\usetikzlibrary{decorations,arrows}
\usetikzlibrary{decorations.pathmorphing}
\usepgflibrary{decorations.pathreplacing} 
\usepackage{pgf}
\usepackage{pgfkeys}
\usepackage{pgffor}

\usepackage{pgfpages}
\usepackage{scalefnt}
\usetikzlibrary{shapes}

\usepackage{titlesec}
\titleformat{\chapter}[hang]{\huge\bfseries}{\thechapter\hspace{7mm}}{0pt}{\huge\bfseries}

\definecolor{hellgrau}{gray}{0.7}

\allowdisplaybreaks

\newtheorem{definition}{Definition}[section]
\newtheorem{theorem}[definition]{Theorem}
\newtheorem{proposition}[definition]{Proposition}
\newtheorem{lemma}[definition]{Lemma}
\newtheorem{corollary}[definition]{Corollary}
\newtheorem{remark}[definition]{Remark}

\theoremstyle{definition}
\newtheorem{conjecture}{Conjecture}
\newtheorem{question}{Question}
\newtheoremstyle{example}
      {0.4em}				
      {0.4em}				
      {\hangindent1em}		
      {1.5em}				
      {\itshape}				
      {.}					
      { }					
      {}						

\theoremstyle{example}
\newtheorem{example}[definition]{Example}

\renewenvironment{proof}{{\bfseries Proof.}}{\hfill$\Box$}

\let\tempone\itemize
\let\temptwo\enditemize
\renewenvironment{itemize}{\tempone\addtolength{\itemsep}{-0.3\baselineskip}}{\temptwo}

\makeatletter 
\let\orgdescriptionlabel\descriptionlabel 
\renewcommand*{\descriptionlabel}[1]{%
\let\orglabel\label   
\let\label\@gobble   
\phantomsection   
\edef\@currentlabel{#1}%
\let\label\orglabel   
\orgdescriptionlabel{#1}%
} 
\makeatother 

\newenvironment{Itemize}[1][1]{\list{\color{blau}$\bullet$}{\setlength\parsep{#1\baselineskip}}}{\endlist}

\makeatletter
\newcommand{\pushright}[1]{\ifmeasuring@#1\else\omit\hfill$\displaystyle#1$\fi\ignorespaces}
\newcommand{\pushleft}[1]{\ifmeasuring@#1\else\omit$\displaystyle#1$\hfill\fi\ignorespaces}
\makeatother

\usepackage{geometry}
\usepackage[
acronym,	
toc			
]
{glossaries}

\newglossary[slg]{symbolslist}{syi}{syg}{ }

\protect

\newglossaryentry{as}{
name= $\mathscr{A}$,
description={alphabet},
sort=alphabet, type=symbolslist
}

\newglossaryentry{asG}{
name= $\mathscr{A}^G$,
description={product space $\prod_{g\in G} \mathscr{A}=\{\xi:G\to\mathscr{A}\}$},
sort=asG, type=symbolslist
}

\newglossaryentry{AT}{
name= $(A_t)_{t\in\protect\ts}$,
description={field of linear operators $A_t:\hs_t\to\hs_t$ over the family of Hilbert spaces $(\hs_t)_{t\in\ts}$},
sort=AT, type=symbolslist
}

\newglossaryentry{Nub}{
name= $\mathcal{B}_{\mathscr{W}}$,
description={branching map of the dictionary $\ws$},
sort=branchingmap, type=symbolslist
}

\newglossaryentry{csX}{
name= $\protect\cs(X)$,
description={set of closed subsets of $X$},
sort=closedsubsetsX, type=symbolslist
}

\newglossaryentry{CM}{
name= $\mathbb{C}$,
description={complex numbers},
sort=CM, type=symbolslist
}

\newglossaryentry{CGred}{
name= $\protect\CG^{\ast}_{\textit{red}}(\Gamma)$,
description={reduced groupoid $C^\ast$-algebra},
sort=CGreduced, type=symbolslist
}

\newglossaryentry{CGfull}{
name= $\protect\CG^{\ast}_{\textit{full}}(\Gamma)$,
description={full groupoid $C^\ast$-algebra},
sort=CGfull, type=symbolslist
}

\newglossaryentry{comp}{
name= $p_{\mathscr{W}}$,
description={complexity function of the dictionary $\ws$},
sort=pwcomplexityfunction, type=symbolslist
}

\newglossaryentry{DG}{
name= $\mathfrak{D}_{G}(\mathscr{A})$,
description={space of dictionaries over the group $G$ and the alphabet $\mathscr{A}$},
sort=DGspacedictionaries, type=symbolslist
}

\newglossaryentry{dH}{
name= $d_H$,
description={Hausdorff metric},
sort=dHhausdorffmetric, type=symbolslist
}

\newglossaryentry{Ext}{
name= $\protect\Ext[v]$,
description={set of extension of the pattern $[v]\in\PaG$},
sort=extensions, type=symbolslist
}

\newglossaryentry{Gamma}{
name= $\Gamma$,
description={topological groupoid},
sort=Gamma, type=symbolslist
}

\newglossaryentry{Gamma0}{
name= $\Gamma^{(0)}$,
description={unit space of the topological groupoid $\Gamma$},
sort=Gamma0, type=symbolslist
}

\newglossaryentry{Gamma1}{
name= $\Gamma^{(1)}$,
description={arrow space of the topological groupoid $\Gamma$},
sort=Gamma1, type=symbolslist
}

\newglossaryentry{Gamma2}{
name= $\Gamma^{(2)}$,
description={set of composable arrows in the topological groupoid $\Gamma$},
sort=Gamma2, type=symbolslist
}

\newglossaryentry{Gammax}{
name= $\Gamma^x$,
description={$r$-fiber of $\Gamma$ at $x\in\Gamma^{(0)}$},
sort=Gammax, type=symbolslist
}

\newglossaryentry{Gun}{
name= $\Gamma_{\mathfrak{u}}$,
description={universal groupoid of $(X,G,\alpha)$},
sort=Gun, type=symbolslist
}

\newglossaryentry{Hol}{
name= $\protect\Hol$,
description={H\"older constant},
sort=hoelderconstant, type=symbolslist
}

\newglossaryentry{IGX}{
name= $\mathcal{I}_{G}(X)$,
description={space of dynamical subsystems associated with $(X,G,\alpha)$},
sort=invariantclosedsubsets, type=symbolslist
}

\newglossaryentry{ksX}{
name= $\protect\ks(X)$,
description={set of compact subsets of the topological space $X$},
sort=kompakt1subsetsX, type=symbolslist
}

\newglossaryentry{ke}{
name= $\protect\ke$,
description={set of translation equivalent compact sets of the group $G$},
sort=kompakt2equivalence, type=symbolslist
}

\newglossaryentry{L(H)}{
name= $\mathcal{L}(\protect\hs)$,
description={algebra of bounded, linear operators on a Hilbert space $\protect\hs$},
sort=L(H), type=symbolslist
}

\newglossaryentry{MathN}{
name= $\protect\mathfrak{N}$,
description={set of normal operators on a Hilbert space $\hs$},
sort=N, type=symbolslist
}

\newglossaryentry{NM}{
name= $\mathbb{N}$,
description={natural numbers},
sort=NM, type=symbolslist
}

\newglossaryentry{normfull}{
name= $||\protect\cdot ||_{\textit{full}}$,
description={full norm},
sort=normfull, type=symbolslist
}

\newglossaryentry{normreduced}{
name= $||\protect\cdot ||_{\textit{red}}$,
description={reduced norm},
sort=normreduced, type=symbolslist
}

\newglossaryentry{oNM}{
name= $\overline{\mathbb{N}}$,
description={one-point compactification of $\mathbb{N}$},
sort=NMonepointcompactification, type=symbolslist
}

\newglossaryentry{Orb}{
name= $\protect\Orb(x)$,
description={orbit of $x$},
sort=orbit, type=symbolslist
}

\newglossaryentry{PGXi}{
name= $\protect\PG^\ast(\Gamma(\Xi))$,
description={pattern equivariant algebra of the grou\-poid $\Gamma(\Xi)$},
sort=PatEquAlg, type=symbolslist
}

\newglossaryentry{osKu}{
name= $\protect\os(K{,}[u])$,
description={base for the topology on $\as^G$},
sort=osKu, type=symbolslist
}

\newglossaryentry{PaG}{
name= $\mathit{Pat}_{G}(\mathscr{A})$,
description={set of patterns},
sort=patternset, type=symbolslist
}

\newglossaryentry{PAX}{
name= $\protect\PA(X)$,
description={set of periodically approximable dynamical subsystem of $(X,G,\alpha)$},
sort=periodicallyapproximablesubshifts, type=symbolslist
}

\newglossaryentry{pix}{
name= $\pi^x$,
description={left-regular representation at $x\in\Gamma^{(0)}$},
sort=pix-left-regularrepresentation, type=symbolslist
}

\newglossaryentry{Po}{
name= $\protect\Po$,
description={set of complex-valued polynomials},
sort=polynomialscomplex, type=symbolslist
}

\newglossaryentry{Pt}{
name= $\protect\Pt$,
description={set of real-valued polynomials up to degree $2$},
sort=polynomialsdegreetwo, type=symbolslist
}

\newglossaryentry{PtM}{
name= $\protect\Pt(M)$,
description={set of real-valued polynomials up to degree $2$ with $1$-norm bounded by $M$},
sort=polynomialsdegreetwoM, type=symbolslist
}

\newglossaryentry{RM}{
name= $\mathbb{R}$,
description={real numbers},
sort=RM, type=symbolslist
}

\newglossaryentry{MathS}{
name= $\protect\mathfrak{S}$,
description={set of self-adjoint operators on a Hilbert space $\hs$},
sort=S, type=symbolslist
}

\newglossaryentry{Sigma}{
name= $\Sigma$,
description={map $\Sigma:\protect\ts\to\protect\ks(\CM)\,,\; t\mapsto \sigma(A_t)$ associated with a field of bounded operators $(A_t)_{t\in\protect\ts}$},
sort=Sigma, type=symbolslist
}

\newglossaryentry{SPX}{
name= $\protect\mathpzc{SP}_{G}(X)$,
description={set of strongly periodic dynamical subsystem of $(X,G,\alpha)$},
sort=stronglyperiodicsubsystems, type=symbolslist
}

\newglossaryentry{Stab}{
name= $\protect\Stab_\alpha(x)$,
description={stabalizer of $x\in X$ with respect to the group action $\alpha$ of the group $G$ on the compact space $X$},
sort=stabalizer, type=symbolslist
}

\newglossaryentry{supp}{
name= $\protect\supp(\protect\fz)$,
description={support of the function $\protect\fz$},
sort=supportfunction, type=symbolslist
}

\newglossaryentry{Sub}{
name= $\protect\Sub[u]$,
description={set of subpatterns of $[u]\in\PaG$},
sort=subpatterns, type=symbolslist
}

\newglossaryentry{MathU}{
name= $\protect\mathfrak{U}$,
description={set of unitary operators on a Hilbert space $\hs$},
sort=U, type=symbolslist
}

\newglossaryentry{UFO}{
name= $\protect\us(F{,}\protect\os)$,
description={base for the Vietoris-topology (with $F$ closed and $\os$ a finite family of open sets) respectively Fell-topology (with $F$ compact and $\os$ a finite family of open sets)},
sort=UFO, type=symbolslist
}

\newglossaryentry{ws}{
name= $\mathscr{W}$,
description={dictionary},
sort=wsdictionary, type=symbolslist
}

\newglossaryentry{wsxi}{
name= $\mathscr{W}(\xi)$,
description={dictionary associated with $\xi\in\as^G$},
sort=wsdictionary1Xi, type=symbolslist
}

\newglossaryentry{wsXi}{
name= $\mathscr{W}(\Xi)$,
description={dictionary associated with the subshift $\Xi\in\SG\big(\as^G\big)$},
sort=wsdictionary2Xi, type=symbolslist
}

\newglossaryentry{Xun}{
name= $X_{\mathfrak{u}}$,
description={universal space of $(X,G,\alpha)$},
sort=Xun, type=symbolslist
}

\newglossaryentry{XalphaG}{
name= $X\rtimes_\alpha G$,
description={transformation group groupoid},
sort=XalphaG, type=symbolslist
}

\newglossaryentry{Xiws}{
name= $\Xi(\protect\ws)$,
description={subshift defined by the dictionary $\ws$},
sort=Xiws, type=symbolslist
}

\newglossaryentry{ZM}{
name= $\mathbb{Z}$,
description={integers},
sort=ZM, type=symbolslist
}

\makeglossaries

\newcommand{\Cc}{{\mathcal C}}

\newcommand{\Ff}{{\mathcal F}}
\newcommand{\Ll}{{\mathcal L}}
\newcommand{\Pp}{{\mathcal P}}

\newcommand{\CM}{{\mathbb C}}

\newcommand{\NM}{{\mathbb N}}
\newcommand{\RM}{{\mathbb R}}
\newcommand{\SM}{{\mathbb S}}

\newcommand{\ZM}{{\mathbb Z}}

\newcommand{\AG}{{\mathfrak A}}

\newcommand{\BG}{{\mathfrak B}}

\newcommand{\CG}{{\mathfrak C}}

\newcommand{\PG}{{\mathfrak P}}

\newcommand{\as}{{\mathscr A}}
\newcommand{\bs}{{\mathscr B}}
\newcommand{\cs}{{\mathscr C}}
\newcommand{\es}{{\mathscr E}}
\newcommand{\gs}{{\mathscr G}}
\newcommand{\hs}{{\mathscr H}}
\newcommand{\ks}{{\mathscr K}}
\newcommand{\os}{{\mathscr O}}
\newcommand{\us}{{\mathscr U}}
\newcommand{\ts}{{\mathscr T}}
\newcommand{\vs}{{\mathscr V}}
\newcommand{\ws}{{\mathscr W}}

\newcommand{\Gz}{\mathpzc{G}}
\newcommand{\fz}{\mathpzc{f}}
\newcommand{\gz}{\mathpzc{g}}
\newcommand{\hz}{\mathpzc{h}}
\newcommand{\sz}{\mathpzc{s}}

\newcommand{\hg}{{\widehat{g}}}

\newcommand{\hphi}{{\widehat{\phi}}}

\newcommand{\ke}{{\mathscr K}_{\sim}\!(G)}			
\newcommand{\Sub}{{\it Sub}} 						
\newcommand{\Ext}{{\it Ext}} 						

\newcommand{\supp}{\mbox{\rm supp}}                	

\newcommand{\Stab}{{\it Stab}} 						
\newcommand{\Orb}{\textit{Orb}}			   		   	

\newcommand{\pt}{(\textit{p2})}			   		   	
\newcommand{\Pt}{\Pp_{\textit{2}}}		   		   	
\newcommand{\po}{(\textit{p})}			   		   	
\newcommand{\Po}{\Pp}					   		   	

\newcommand{\oNM}{\overline{{\mathbb N}}}	   		
\newcommand{\dil}{\mbox{\rm dil}}	  	   			
\newcommand{\Hol}{\mbox{\rm Hol}}	  	   			

\DeclareMathAlphabet{\mathpzc}{OT1}{pzc}{m}{it}

\newcommand{\SG}{\mathcal{I}_{G}}		   			
\newcommand{\SZ}{\mathcal{I}_{\ZM}}	   				
\newcommand{\SZd}{\mathcal{I}_{\ZM^d}}	   			
\newcommand{\SP}{\mathpzc{SP}_{G}}		   			
\newcommand{\SPZ}{\mathpzc{SP}_{\ZM}}		   		
\newcommand{\PA}{\mathpzc{PA}_{G}}		   			

\newcommand{\PaG}{\mathit{Pat}_{G}(\as)}		   		
\newcommand{\PaZ}{\mathit{Pat}_{\ZM}(\as)}		   	
\newcommand{\PaZd}{\mathit{Pat}_{\ZM^d}(\as)}	   	
\newcommand{\DG}{\mathfrak{D}_{G}(\as)}				
\newcommand{\DZ}{\mathfrak{D}_{\mathbb{Z}}(\as)}		
\newcommand{\DZd}{\mathfrak{D}_{\mathbb{Z}^d}(\as)}	
\newcommand{\Nub}{\mathcal{B}_\ws}				   	
\newcommand{\comp}{p_\ws}						   	
\newcommand{\compX}{p_{\ws(\Xi)}}				   	

\newcommand{\Xun}{\mathit{X}_{\mathfrak{u}}}			   		
\newcommand{\Gun}{\Gamma_{\!\mathfrak{u}}}				   		
\newcommand{\asGun}{\as^G_{\mathfrak{u}}}				   		

\let\X\chi
\renewcommand{\chi}{\raisebox{3pt}{$\X$}}

\renewcommand{\chaptermark}[1]%
         {\markboth{\thechapter.\ #1}{}}
\renewcommand{\sectionmark}[1]%
         {\markright{\thesection\ #1}}         
\newcommand{\JENATitle}[9]{

  \thispagestyle{empty}

  \vspace*{\stretch{1}}
  {\parindent0cm
  \rule{\linewidth}{.7ex}}
  \begin{flushright}
    \vspace*{\stretch{1}}
    \sffamily\bfseries\Huge
    #1\\
    \vspace*{\stretch{1}}
  \end{flushright}
  \rule{\linewidth}{.7ex}

  \vspace*{\stretch{3}}
  \begin{center}
    \Huge Dissertation\\[0.2cm]
    \Large zur Erlangung des akademischen Grades\\[0.3cm]
    \Large \textbf{doctor rerum naturalium}\\[0.3cm]
    \Large vorgelegt dem Rat\\
    \Large der Fakult\"at f\"ur Mathematik und Informatik\\
    \Large der Friedrich-Schiller-Universit\"at Jena\\
    \vspace*{\stretch{1}}
    \Large von\\[0.2cm]
    \Large \textbf{#2}\\[0.2cm]
    \Large geboren am #3 in #4\\
    \vspace*{\stretch{2}}
    \includegraphics[width=1.5in]{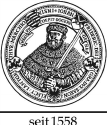}
  \end{center}

  \newpage
  \thispagestyle{empty}

  \vspace*{\stretch{1}}

  \begin{flushleft}
	\begin{tabular}{ll}
   	\Large1. Gutachter: & \Large #7 \\
    		&(Friedrich-Schiller-Universi\"at Jena,\\
    		&\, Ernst-Abbe-Platz 2, 07743 Jena, Deutschland)\\[2mm]
    	\Large 2. Gutachter: & \Large #8 \\
    		&(Georgia Institute of Technology, School of Mathematics,\\
    		&\, 686 Cherry Street, Atlanta, GA 30332-0160, USA)\\[2mm]
	\Large 3. Gutachter: & \Large #9 \\
    		&(Institut Camille Jordan, Universit\'e Claude Bernard - Lyon I,\\
    		&\, 21 avenue Claude Bernard, F-69622 Villeurbanne Cedex, Frankreich)\\[7mm]
    	\end{tabular}
    	\begin{tabular}{ll}
		\Large Tag der \"offentlichen Verteidigung: & \Large 6. Oktober 2016\Large \\	
	\end{tabular}
  \end{flushleft}

  \cleardoublepage
}

\begin{document}

\pagenumbering{roman} 
  \JENATitle
      {Spectral approximation\\ of aperiodic Schr\"odinger operators 
	  }               									
      {Dipl.-Math. Siegfried Beckus}                       			
      {11.04.1988}                             			
      {Erfurt}                             				
      {Friedrich-Schiller-Universit\"at Jena, 2016}   	
      { }                            					
      {Prof. Dr. Daniel H. Lenz}							
      {Prof. Dr. Jean V. Bellissard}						
      {Prof. Dr. Johannes Kellendonk}					

\fontsize{12}{14}
\selectfont

\markboth{Zusammenfassung in deutscher Sprache}{Zusammenfassung in deutscher Sprache}
  \addcontentsline{toc}{chapter}{\protect Zusammenfassung in deutscher Sprache}

\chapter*{Zusammenfassung in deutscher Sprache}

Die vorliegende Arbeit befasst sich mit Stetigkeitseigenschaften der Abbildung 
$$
\Sigma:\ts\to\ks(\CM)\,,
	\quad 
	t\mapsto \sigma(A_t)\,,
$$
f\"ur einen gegebenen topologischen Raum $\ts$ und eine Familie von linearen, beschr\"ankten Operatoren $A_t\in\Ll(\hs_t)\,,\; t\in\ts\,,$ definiert auf den Hilbertr\"aumen $\hs_t\,,\; t\in\ts$. Hier be\-zeichnet $\sigma(A_t)$ das Spektrum des Operators $A_t$. Der Raum $\ks(\CM)$ der kompakten Teilmengen von $\CM$ ist ausgestattet mit der durch die Hausdorffmetrik induzierten Topologie. Die Stetigkeit von kompakten Teilmengen in $\RM$ wird durch die Stetigkeit des Randes der kompakten Mengen charakterisiert. Hierbei werden interessante Effekte f\"ur selbst-adjungierte Operatoren $A_t\,,\; t\in\ts\,,$ beobachtet, bei denen sich spektrale L\"ucken schlie\ss en.

\medskip

Die Stetigkeit der Abbildung $\Sigma$ wird charakterisiert durch die Stetigkeit der Normen $\ts\ni t\mapsto\|p(A_t,A_t^\ast)\|\in[0,\infty)$ f\"ur gewisse (komplexwertige) Polynome. F\"ur selbst-adjungierte, unit\"are beziehungsweise normale Operatoren m\"ussen dabei unterschiedliche Klassen betrachtet werden. 
Falls die Operatoren alle auf einem Hilbertraum $\hs=:\hs_t\,,\; t\in\ts\,,$ definiert sind, so wird die Stetigkeit von $\Sigma$ auf eine Topologie auf dem Raum der linearen, beschr\"ankten Operatoren $\Ll(\hs)$ zur\"uckgef\"uhrt. F\"ur selbst-adjungierte und unit\"are Operatoren l\"asst sich zus\"atzlich die H\"older-Stetigkeit von $\Sigma$ durch die gleich\-m\"a\ss ige H\"older-Stetigkeit gewisser Normen der Operatoren charakterisieren, falls $(\ts,d)$ ein vollst\"andiger metrischer Raum ist. Hierbei wird die Konvergenzrate der Spektren bez\"uglich der Hausdorffmetrik verringert, welches dem m\"oglichen Schlie\ss en von spektralen L\"ucken zugrunde liegt.

\medskip

Diese neuen Resultate basieren auf der gemeinsamen Arbeit \cite{BeBe16} mit {\sc J. Bellissard} und sind motiviert durch ein abstraktes Resultat aus der Theorie der stetigen Felder von $C^\ast$-algebren. Es stellt sich heraus, dass sich $\Sigma$ genau dann stetig verh\"alt, wenn die $C^\ast$-algebren $\CG^\ast(A_t,I_t)\,,$ $t\in\ts\,,$ ein stetiges Feld von $C^\ast$-algebren definieren. Damit liefert die Arbeit \cite{LaRa99} ein Hilfsmittel, um die Stetigkeit der Abbildung $\Sigma$ zu beweisen, falls die Operatoren Elemente von Gruppoid $C^\ast$-algebren sind.

\medskip

Motiviert durch diese sehr allgemeinen Resultate werden anschlie\ss end Operatoren, assoziiert zu einem dynamischen System, bez\"uglich des Verhaltens der Spektren unter \"Anderung des zugrundeliegenden dynamischen Systems untersucht. Jedes dynamische System $(X,G,\alpha)$ definiert auf nat\"urliche Weise ein Gruppoid, welches die topologischen und dynamischen Eigenschaften kodiert. Die sogenannten verallgemeinerten Schr\"odinger\-oper\-atoren lassen sich als stetige Funktionen auf dem Gruppoid darstellen. Basierend auf der gemeinsamen Arbeit \cite{BeBeNi16} mit {\sc J. Bellissard} und {\sc G. de Nittis}, liefert die Konstruktion des sogenannten universalen dynamischen Systems $(\Xun,G,\alpha)$ und des zugeh\"origen universalen Gruppoids $\Gun$, assoziiert zu einem dynamischen System $(X,G,\alpha)$, ein stetiges Feld von Gruppoiden. Durch Verwendung des Resultats \cite{LaRa99} f\"uhrt dies zur Stetigkeit der Abbildung
$$
\Sigma:\SG(X)\to\ks(\CM)\,,
	\quad
	Y\mapsto\sigma(H_Y)\,,
$$
f\"ur alle zu $(\Xun,G,\alpha)$ assoziierten verallgemeinerten Schr\"odingeroperatoren. Hier bezeichnet $\SG(X)$ die Menge der $G$-invarianten, abgeschlossenen Teilmengen von $X$, ausgestattet mit der Hausdorfftopologie. Zus\"atzlich wird gezeigt, dass eine Folge von dynamischen Systemen $Y_n\subseteq X\,,\;n\in\NM\,,$ genau dann gegen das dynamische System $Y\subseteq X$ konvergiert, falls die Spektren $\big(\sigma(H_{Y_n})\big)_{n\in\NM}$ f\"ur alle verallgemeinerten Schr\"odingeroperatoren $H$ gegen $\sigma(H_Y)$ konvergieren. 

\medskip

Dieses Resultat zeigt den starken Zusammenhang zwischen den dynamischen Systemen und den entsprechenden Schr\"odingeroperatoren. Weiterhin liefert dies ein Hilfsmittel zur Approximation von Schr\"o\-dingeroperatoren in dem das zugrundeliegende dynamische System approximiert wird, wie im Folgendem beschrieben. 

\medskip

Im Fall von symbolischen dynamischen Systemen $(\as^G,G,\alpha)$ f\"ur eine diskrete, abz\"ahlbare Gruppe $G$ wird das Konzept eines W\"orterbuches konzeptionell eingef\"uhrt. Es stellt sich heraus, dass die Menge der W\"orterb\"ucher, ausgestattet mit der lokalen Mustertopologie, hom\"oomorph zu dem topologischen Raum $\SG\big(\as^G\big)$ ist. Damit lassen sich gewisse Subshifts von $(\as^G,G,\alpha)$ f\"ur $G=\ZM^d$ durch stark periodische Subshifts approximieren. Im Zusammenhang mit den vorangegangen Resultaten lassen sich damit alle verallgemeinerten Schr\"odingeroperatoren dieser Subshifts durch periodische Schr\"odingeroperatoren approximieren, so dass die Operatoren in der starken Operatortopologie konvergieren und gleichzeitig die Spektren der Operatoren konvergieren. 

\medskip

Diese Resultate sind besonders interessant, da mit den bisherigen Methoden haupts\"achlich Schr\"o\-dinger- und Jakobioperatoren \"uber der Gruppe $\ZM$ untersucht werden konnten. Die Charakterisierung des stetigen Verhaltens der Spektren gilt f\"ur viel allgemeinere Gruppen $G$. Zus\"atzlich lassen sich die hier verwendeten Methoden auch auf Delonemengen und deren zugeh\"origen Schr\"odingeroperatoren erweitern, c.f. \cite{BeBeNi16}. Insgesamt liefert die hier entwickelte Theorie eine M\"oglichkeit, auch h\"oherdimensionale Systeme zu untersuchen. Zus\"atzlich zeigt die vorliegende Arbeit einen starken Zusammenhang zu der Frage auf, ob ein Subshift endlichen Typs, eine dichte Teilmenge an periodischen Elementen enth\"alt. Diese Frage spielt im Gebiet der dynamischen Systeme und Tilings eine wichtige Rolle. Weitere Resultate der vorliegenden Arbeit finden sich in den Aufs\"atzen \cite{BeLeMaCh14,BeLeMaCh16}.

\cleardoublepage
\tableofcontents
\markboth{Table of Contents}{Table of Contents}
\clearpage

\markboth{Acknowledgments}{Acknowledgments}
  \addcontentsline{toc}{chapter}{\protect Acknowledgments}

\chapter*{Acknowledgments}

First of all, I would like to express my sincere thank to my adviser Daniel Lenz for his continuous and careful guidance and support during the last years. He always had an open door for my questions and he provided me with all necessary means to travel for my research. I would like to thank him for lots of fruitful and enlightening discussions and advises. Furthermore, I am grateful to him for connecting me with other researchers in the community and showing me different perspectives.

\medskip

I am especially grateful to Jean Bellissard who generously shared his mathematical experiences and knowledge with me. We had a lot of inspiring, motivating and interesting discussions. Furthermore, I am thankful for his guidance during the last years and that he offered me the support to visit him several times in Atlanta. Additionally, I thank him for providing me with the contact to David Damanik and Anton Gorodetski and giving me the opportunity to meet them.

\medskip

Also I am deeply thankful to my coworkers Giuseppe de Nittis and Felix Pogorzelski for their supervision especially in the beginning of my thesis. I thank them also for several stimulating discussions over the last years that pushed me further.

\medskip

For proofreading, I am very grateful to Juliane Beckus, Annegret Harendt, Daniel Lenz, Markus Lange, Felix Pogorzelski and Marcel Schmidt for reading parts of the thesis in advance and providing me with useful remarks and suggestions. 

\medskip

I am especially thankful to my colleagues and friends Daniel L., Daniel S., Felix, Gerhard, Marcel, Matthias, Markus, Melchior, Ren\'e, Sebastian, Therese and Xueping for creating a lively and joyful atmosphere and all the recreative breaks we had together. I would particularly single out Therese Mieth for developing with me our joint ideas. In light of this, I would like to thank the faculty of mathematics and computer science in Jena supporting our ideas for the PhD Seminar and the graduate school. Specifically, I would like to thank Dorothee Haroske, Daniel Lenz and Martin Mundhenk for their encourage and support for our ideas. In addition, I am grateful to the Graduiertenakademie Jena for their financial support and their guidance in organizational questions. Especially, I would like to single out Hanna Kauhaus who was our contact person.

\medskip

I would like to thank Tobias Hartnick for pointing out the work of de Bruijn and its connection to our work.

\medskip

In addition, I would like to express my thank to David Damanik and Anton Gorodetski inviting me to the workshop \glqq Spectral Properties of Quasicrystals via Analysis, Dynamics, and Geometric Measure Theory\grqq\, at BIRS Oaxaca in 2015. This trip was very enlightening for me and it allowed me to create new networks. Furthermore, I am grateful to Anton Gorodetski for raising stimulating questions at this workshop that led to some results presented here. Also I thank Jake Fillman whom I met in Oaxaca for rising up the question for the continuous behavior of the spectra for unitary operators and an interesting discussion about it.

\medskip

I would like to thank Antoine Julien and Franz Luef for inviting me twice to Trondheim where I had very good discussions with them. I really enjoyed the time there.

\medskip

I am thankful to Magnus Landstadt for several discussions concerning the question when a groupoid $C^\ast$-algebra is unital, especially, in the case of transformation group groupoids. In light of this I also thank Jean Renault for good ideas concerning the general case.

\medskip

In addition, I would like to thank the School of Mathematics at Georgia Institute of Technology for providing an office for me and for their support during my stays in Atlanta, especially, during my long term stay from February to April 2014. Additionally, I would like to express my gratitude to the National Science Foundation supporting several trips and stays at the Georgia Institute of Technology by the grant \glqq Spectral Properties of Aperiodic Solids\grqq, Grant No. DMS{1160962}.

\medskip

I thank the Research Training Group (1523/2), at the Friedrich-Schiller-University of Jena, Germany, for financial support for several business trips during the last years which allowed me to advertise my work and to create new projects.

\medskip

Furthermore, I am grateful to the Erwin Schr\"odinger Institute, Vienna, for support during the Summer of 2014 where parts of these results were obtained.  Additionally, I would like to thank the DAAD for sponsering my business trip to Santiago de Chile to join the International Congress of Mathematical Physics XVIII in July 2015 where also parts of my thesis were obtained.

\pagenumbering{arabic}
\setcounter{page}{1}


\chapter{Introduction}
\label{Chap1-Intro}

When do the spectra of a family of linear operators vary continuously? What is the correct topology on the space of bounded linear operators so that the spectra of normal operators vary continuously? Can the continuous behavior of the spectra be characterized quantitatively? Once identified this topology, is there a tool to check that a family of operators is continuous in this topology? These questions are addressed and answered in this thesis. In addition, our approach leads to a characterization for the continuous behavior of the spectra of Schr\"odinger operators by the continuous variation of the under\-lying dynamical system. With this at hand, periodic approximations for Schr\"odinger operators are provided.

\medskip

Let $\ts$ be a topological space and consider a family of bounded normal linear operators $A_t:\hs_t\to\hs_t\,,\; t\in\ts\,,$ over a family of Hilbert spaces $\hs_t\,,\; t\in\ts$. The map
$$
\text{\gls{Sigma}}:\ts\to\ks(\CM)
	\,,\quad t\mapsto \sigma(A_t)\,,
$$
which maps $t\in\ts$ to the spectrum $\sigma(A_t)\in\ks(\CM)$ of $A_t$, is analyzed for its (H\"older-)con\-tinuity. Here $\ks(\CM)$ denotes the set of compact subsets of the complex plane \gls{CM} equipped with the Hausdorff metric. In particular, when $\ts$ is the set of all bounded normal linear operators on a fixed Hilbert space $\hs$, we seek to identify the coarsest topology so that $\Sigma$ is continuous. 

\medskip

This thesis provides a new characterization for the continuity of the map $\Sigma$ when all the operators are self-adjoint or unitary operators by the continuity of specific norms of the operators based on joint work with {\sc J. Bellissard} \cite{BeBe16}. Following the lines of this characterization, we investigate a method to characterize the H\"older-continuity of $\Sigma$ for these operators if $(\ts,d)$ is a complete metric space, c.f. Section~\ref{Chap3-Sect-CharHolContSpect}. This project was motivated by an abstract result in the theory of continuous fields of $C^\ast$-algebras. It states that continuous normal vector fields have spectra that vary continuously. Even in the case of normal operators, a characterization of the continuity of $\Sigma$ is provided by using the machinery of continuous fields of $C^\ast$-algebras, c.f. Theorem~\ref{Chap3-Theo-PContEquivContSpectNormal}. The case of unbounded operators is treated by passing to the resolvent of the operators, c.f. Section~\ref{Chap3-Sect-CharContSpecUnbouSelfAdjOp}. 

\medskip

Next a tool is presented to verify the continuity of $\Sigma$ in concrete situations by $C^\ast$-algebraic techniques based on a result by {\sc Landsman} and {\sc Ramazan} \cite{LaRa99}, c.f. Chapter~\ref{Chap4-ToolContBehavSpectr}. It is then applied to Schr\"odinger operators associated with dynamical systems. Specifically, the continuous behavior of the spectra of generalized Schr\"odinger operators is characterized by the continuous behavior of the associated dynamical systems. This new characterization is based on joint work with {\sc J. Bellissard} and {\sc G. de Nittis} \cite{BeBeNi16} which also covers more general situations than considered here. More precisely, the continuous variation of the spectra of Schr\"odinger operators associated with Delone sets can be treated as well, c.f. the discussion in Section~\ref{Chap8-Sect-DeloneTiling}. The main ingredient for the proof of this characterization is the investigation of the concept of the universal dynamical system and the universal groupoid corresponding to a dynamical system. 

\medskip

The connection between the dynamical systems and the spectra of associated Schr\"odinger operators yields a tool to find periodic approximations for operators in terms of periodic approximations of the underlying dynamical systems. As a starting point for further research, we analyze symbolic dynamical systems associated with the group $\ZM^d$ for the existence of strongly periodic approximations, c.f. Chapter~\ref{Chap5-OneDimCase}, Chapter~\ref{Chap6-HigherDimPerAppr} and Chapter~\ref{Chap7-Examples}. 

\medskip

It is remarkable that these results about the continuity of the spectra are independent of the dimension. This is of particular interest as the known results for Schr\"odinger operators associated with quasicrystals are restricted to one-dimensional systems or systems that can be reduced to them. Furthermore, the developed theory in this thesis shows the strong connection between the existence of periodic approximations and the existence of a (dense) subset of strongly periodic elements in a subshift of finite type. The second problem is studied in the field of dynamical systems \cite{Pia08,Fio09,Hoc09,SiCo12,Coh14,CaPe15} which is also related to the so called Wang tiles \cite{Wan61,Ber66,CuKa95}.

\medskip

As described above, this thesis combines several mathematical concepts. Since the aim of this work is to provide an accessible text for a broad readership, the background and the basic concepts are presented here. Additionally, many examples highlight the different aspects of the theory. Finally, it is shown in Chapter~\ref{Chap7-Examples} how these results apply to specific examples which were studied in the literature before. These consequences extend and confirm known results by using different methods.

\medskip

\begin{figure}
\centering
\includegraphics[width=0.82\textwidth]{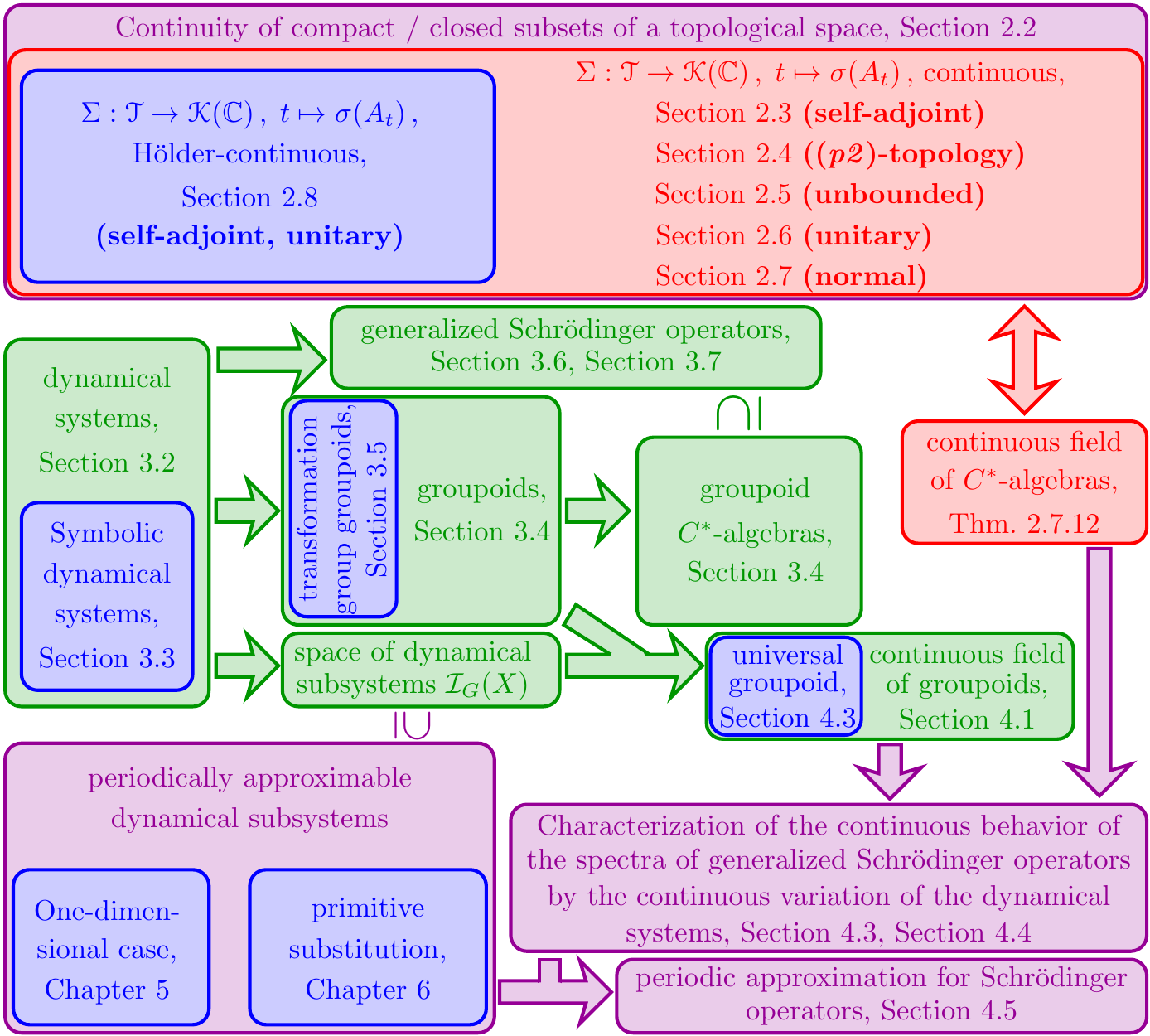}
\caption{ }
\label{Chap1-Fig-Overview}
\end{figure}

In the following, we elaborate the various mathematical theories and how they combine and apply to the study of Schr\"odinger operators in the field of Solid State Physics. An overview of this thesis is given in Figure~\ref{Chap1-Fig-Overview}

\subsection*{The Hausdorff-, Vietoris- and Fell-topology}

The spectrum of a bounded normal linear operator is a non-empty, compact subset of $\CM$. The distance of such sets is usually measured in terms of the Hausdorff metric, c.f. \cite{Hausdorff14,Hausdorff62,CastaingValadier77}. Later on we will also need the concept of a topology on the set of dynamical subsystems which are compact subsets of a compact metrizable space $X$. In order to study the continuity, only the structure of a topology is necessary. In the following we provide an overview for different topologies on compact/closed subsets of a topological space $X$ which are used throughout this work.

\medskip

In 1914, Hausdorff \cite{Hausdorff14} introduced the Hausdorff metric $d_H:\ks(X)\times\ks(X)\to[0,\infty)$ on the set $\ks(X)$ of compact subsets of a complete metric space $(X,d)$. Later, in 1922, {\sc Vietoris} \cite{Vie22} established the so called Vietoris-topology on the set $\cs(X)$ of all closed subsets of a topological space $X$. A base for this topology is given by
$$
\us(F,\os)\; 
	:= \; \big\{Y\in\cs(X)\;|\; 
		F\cap Y=\emptyset\,,\; 
		O\cap Y\neq\emptyset \text{ for all } O\in\os
	\big\}
$$
where $F\subseteq X$ ranges over all closed subsets and $\os$ ranges over all finite families of open subsets of $X$. This topology is also called finite topology \cite{Mic51} or Hausdorff-topology due to the fact that the Vietoris-topology is equal to the topology induced by the Hausdorff metric on $\ks(X)$ if $(X,d)$ is a complete metric space, c.f. \cite[Theo\-rem~II-6]{CastaingValadier77} or Theorem~\ref{Chap3-Theo-VietFellHausMetricEquiv}. In particular, this topology is independent of the metric $d$ as long as $(X,d)$ is complete.

\medskip

{\sc Fell} \cite{Fel61,Fel62} introduced a topology on $\cs(X)$ in the same spirit as Vietoris with the difference that $F\subseteq X$ only ranges over all compact subsets instead of closed subsets in the definition of the base. He was motivated by a $C^\ast$-algebraic point of view. This topology on $\cs(X)$ is called Fell-topology \cite{Bee93,BeBe16} or Chabauty-Fell-topology since {\sc Chabauty} \cite{Cha50} already investigated a special case of the Fell-topology when $X$ is a group. The advantage of the Fell-topology in comparison with the Vietoris-topology is that $\cs(X)$ is automatically compact and Hausdorff in the Fell-topology whenever $X$ is locally compact, c.f. \cite{Fel61}. Further topological properties of the Fell-topology are analyzed in \cite{Bee93}. Since closed subsets of a compact space are compact, the Fell and the Vietoris-topology coincide on compact spaces. Besides, the Fell-topology on $\cs(X)$ can be interpreted as a local version of the Vietoris-topology, c.f. Theorem~\ref{Chap3-Theo-VietFellHausMetricEquiv}. 

\medskip

The Fell-topology and the Vietoris-topology belong to the class of hit and miss topologies which are intensively studied, c.f. \cite{Kuratowski66,Kuratowski68,DiNa08}. All of the topological spaces $X$ that are considered in this thesis are metrizable and the most of them are compact (except $X=\RM$ and $X=\CM$). If $X$ is metrizable and compact, the Fell-topology, Vietoris-topology and the topology induced by the Hausdorff metric coincide for each metric $d$ on $X$ such that $(X,d)$ is complete. Thus, we call this topology Hausdorff-topology since Hausdorff was the first investigating a topology on $\ks(X)$. 

\medskip

From the philosophical point of view, continuous behavior of spectra is typically connected with the continuity of the spectral edges. Based on this idea, the notion of continuous boundaries was first introduced in \cite{BeBe16} which is also presented here. It turns out that the continuity of the boundaries is equivalent to the Fell-continuity of closed subsets in the real numbers \gls{RM}. In the case of compact sets (like for the spectrum of bounded self-adjoint operators), the continuity of the boundaries is actually equivalent to the Vietoris-continuity, c.f. Theorem~\ref{Chap3-Theo-BoundVietContCompVers}. The connection to the behavior of the boundaries (spectral edges) plays an important role in the quantitative estimates of the spectra.

\subsection*{Continuous behavior of the spectra and its relation to \texorpdfstring{$C^\ast$}{C}-algebras}

\begin{figure}
\centering
\includegraphics[width=\textwidth]{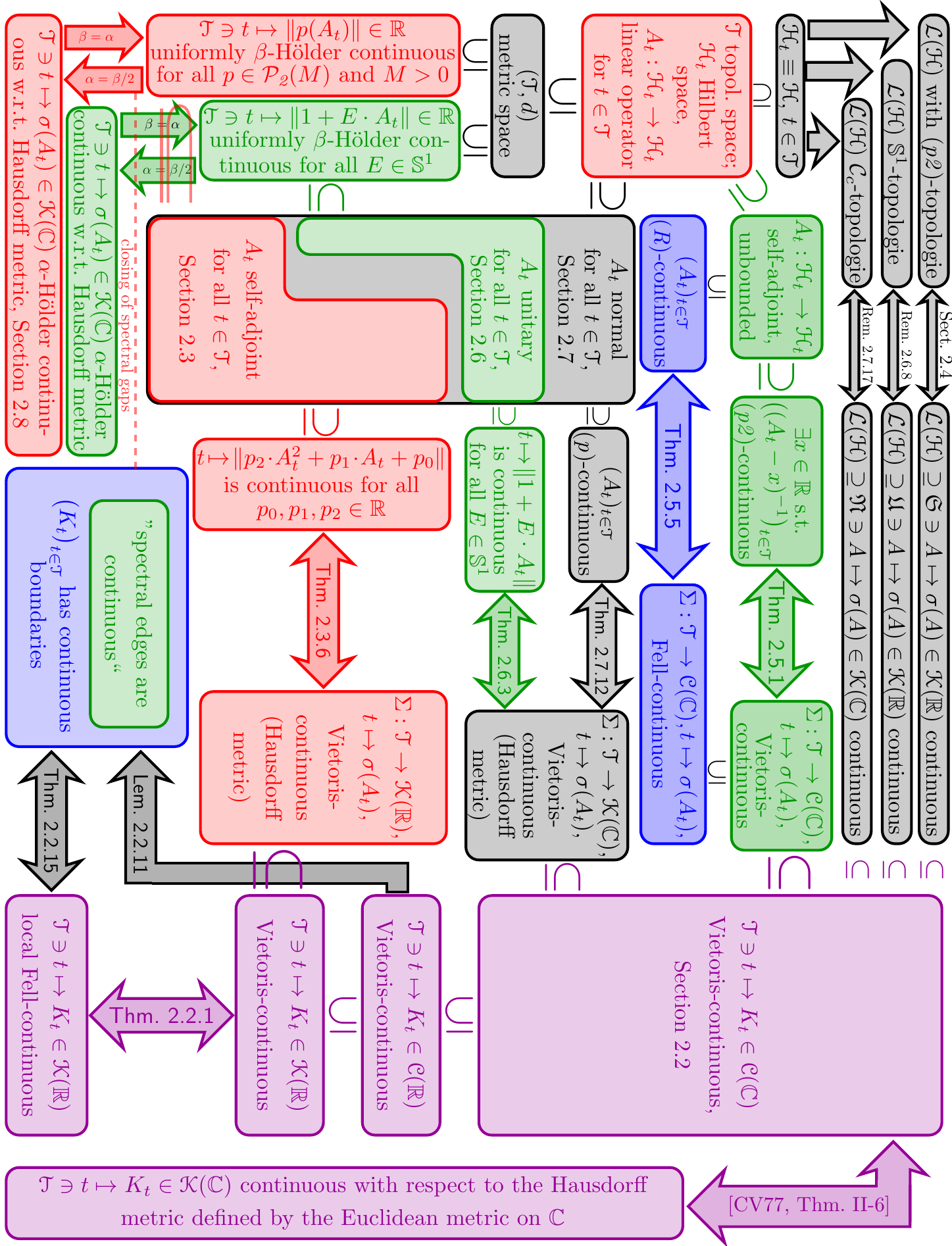}
\caption{Overview and connections between Chapter~\ref{Chap3-SpectAppr} about the continuous behavior of the spectra.}
\label{Chap1-Fig-Strategy3}
\end{figure}

The continuous dependence of spectra is interesting in different areas of mathematics and physics. On the one hand, it appears in perturbation theory. Specifically, the question arises how one can perturb an operator so that its spectrum does not vary too much. On the other hand, it is useful in approximation theory. Then one asks whether there exist approximations for an operator $A$, which is difficult to handle, by operators with known spectral theory. In this case, it is necessary to study which spectral properties are preserved by taking the limit. 
Due to the complexity of the connections of the results presented here, an overview of these results is presented in Figure~\ref{Chap1-Fig-Strategy3}. 

\medskip

Let $\hs$ be a Hilbert space and denote by \gls{L(H)} the algebra of all bounded linear operators on $\hs$. Consider the map
$$
\Sigma_\AG:\AG\subseteq \Ll(\hs)\to\ks(\CM)\,,\;
	A \; \mapsto \; \sigma(A)\,,
	\qquad\quad
	\AG\in\{\mathfrak{N},\mathfrak{S},\mathfrak{U}\}\,,
$$
where
\begin{align*}
\text{\gls{MathN}} \; 
	&:= \; \big\{ A\in\Ll(\hs)\;\big|\; A^\ast A= AA^\ast\big\}\,,
	&\pushright{(\textbf{normal})}\\
\text{\gls{MathS}} \; 
	&:= \; \big\{ A\in\Ll(\hs)\;\big|\; A= A^\ast\big\}\,,
	&\pushright{(\textbf{self-adjont})}\\
\text{\gls{MathU}} \; 
	&:= \; \big\{ A\in\Ll(\hs)\;\big|\; A^\ast A= AA^\ast=I\big\}\,.
	&\pushright{(\textbf{unitary})}
\end{align*}

Note that the study of the map $\Sigma_\AG$ does only encode the spectrum of the operators in $\AG$ as a set and not their spectral type. In order to avoid confusion with the terminology concerning continuity properties of the spectral measures, "continuous behavior of the spectra" is used for the continuity of the map $\Sigma_\AG$ instead of "continuity of the spectra".

\medskip

The algebra $\Ll(\hs)$ is usually equipped with the operator norm topology which preserves many spectral properties. In this case $\Sigma_{\mathfrak{N}}$ is continuous. Unfortunately, the convergence in the operator norm is too restrictive for many applications, see e.g. Example~\ref{Chap3-Ex-AlmostMathieu}. For instance, in the situation of Schr\"odinger operators a change of the potential yields that the operators are far apart in terms of the operator norm while the spectra can be close in the Hausdorff metric, c.f. Example~\ref{Chap4-Ex-KeineNormKonv}. Another possible choice for a topology on $\Ll(\hs)$ is the strong operator topology. Continuity of families of operators in this topology can be verified in many circumstances like for Schr\"odinger operators, c.f. Proposition~\ref{Chap2-Prop-CovFamOp-Spect}. However, the strong operator continuity has only weak implications for the spectral properties, c.f. Section~\ref{Chap3-Sect-OpTop}. The reader is also referred to \cite{Weidmann80,Wei97} for a more detailed discussion on the spectral implications for the operator norm topology and the strong operator topology. 

\medskip

One aim of this work is to provide a complete description of the (H\"older-)continuity of $\Sigma_\AG$. Specifically, the coarsest topology on $\AG\in\{\mathfrak{N},\mathfrak{S},\mathfrak{U}\}$ is identified so that $\Sigma_\AG$ is continuous. The precise description of the topology depends on the class $\AG$ of operators that are considered. The topology is called \pt-topology (Definition~\ref{Chap3-Def-P2Top}) in the self-adjoint case $\AG=\mathfrak{S}$, $\Cc_c$-topology (Remark~\ref{Chap3-Rem-PTopology}) in the normal case $\AG=\mathfrak{N}$ and $\SM^1$-topology (Remark~\ref{Chap3-Rem-S1topology}) in the unitary case $\AG=\mathfrak{U}$. They are defined by specific norms on the operators. These topologies are not vector space topologies. More precisely, the addition in $\Ll(\hs)$ is not continuous with respect to them, c.f. Example~\ref{Chap3-Ex-SumNotP2Cont} and Section~\ref{Chap3-Sect-RelP2ContOthTop}. Clearly, these topologies are coarser than the operator norm topology. On the other hand, there is no relation to the strong operator topology in general, c.f. Proposition~\ref{Chap3-Prop-StrongNotImplP2} and Proposition~\ref{Chap3-Prop-P2NotImplStrong}.

\medskip

One known tool to verify the continuous behavior of the spectra is the concept of continuous fields of Hilbert spaces and $C^\ast$-algebras. $C^\ast$-algebras have their origins in the work of {\sc Gelfand} and {\sc Naimark} \cite{GeNe43} in 1943. The idea of continuous fields was initially proposed by {\sc Kaplanski} \cite{Kap49,Kap51}, {\sc Fell} \cite{Fel60}, {\sc Tomiyama} and {\sc Takesaki} \cite{ToTa61,Tom62} and further developed by {\sc Dixmier}  and {\sc Douady} \cite{DiDo63}. The first result about the connection to the continuous behavior of the spectra goes back to {\sc Kaplanski} \cite[Lemma~3.3]{Kap51} if the operators are self-adjoint, c.f. the discussion at Theorem~\ref{Chap3-Theo-ContFieldCALgContSpectr}. {\sc Elliot} \cite{Ell82} provides a different proof which he applied to almost periodic Schr\"odinger operators. Among others, {\sc Bellissard} \cite{Bel94} used this result for an effective Hamiltonian in the rotation algebra to show the continuous variation of the spectra. This method was recently used in \cite{Man12,BeMa12,PaRi16,BeBeNi16}.

\medskip

The spectrum of a bounded operator is a compact subset of the complex plane. Thus, the Hausdorff distance of two spectra can be measured even if the operators are defined on different Hilbert spaces. Hence, it is natural to ask for the continuous behavior of the spectra even if the operators are defined on different Hilbert spaces. (This is crucial in \cite{BeBeNi16} for the generalization to Schr\"odinger operators associated with Delone sets.) Let $\ts$ be a topological space and $A_t:\hs_t\to\hs_t\,,\; t\in\ts\,,$ be a family of  bounded operators defined on a family of Hilbert spaces $(\hs_t)_{t\in\ts}$. The map
$$
\Sigma:\ts\to\ks(\CM)\,,\; 
	t\mapsto\sigma(A_t)\,,
$$
is analyzed for its continuity where $\ks(\CM)$ is equipped with the Vietoris-topology (Hausdorff metric). Or more quantitatively, the H\"older-continuity of $\Sigma$ is studied in terms of the Hausdorff metric on $\ks(\CM)$ if $(\ts,d)$ is a complete metric space.

\medskip

As shown in \cite{BeBe16}, the concept of continuous fields of $C^\ast$-algebras naturally appears in the study of continuous behavior of the spectra. More specifically, for $\AG\in\{\mathfrak{N},\mathfrak{S},\mathfrak{U}\}$, suppose that $A_t\in\AG$ holds for all $t\in\ts$. Then the following assertions turn out to be equivalent, c.f. Theorem~\ref{Chap3-Theo-P2ContEquivContSpect}, Theorem~\ref{Chap3-Theo-CharContSpectUnitary} and Theorem~\ref{Chap3-Theo-PContEquivContSpectNormal}.
\begin{description}
\item[(i)] The map $\Sigma:\ts\to\ks(\CM)\,,\; t\mapsto\sigma(A_t)\,,$ is Vietoris-continuous (Hausdorff metric on $\ks(\CM)$).
\item[(ii)] The map $\ts\ni t\mapsto\|p(A_t)\|\in\RM$ is continuous for all polynomials $p\in\Po(\AG)$.
\item[(iii)] The $C^\ast$-algebras $\CG_t$ generated by $A_t$ and the identity $I_t\in\Ll(\hs_t)$ for $t\in\ts$ define a continuous field of unital $C^\ast$-algebras.
\end{description}
Here $\Po(\AG)$ denotes one of the following sets
\begin{align*}
\Po(\mathfrak{S}) \; 
	&:= \; \big\{ p:\RM\to\RM \;\big|\; p(z)=p_0+p_1\cdot z + p_2\cdot z^2\,,\; p_0,p_1,p_2\in\RM \big\}\,,\\
\Po(\mathfrak{U}) \; 
	&:= \; \big\{ p:\CM\to\CM \;\big|\; p(z)=1+E\cdot z\,,\; E\in\SM^1 \big\}\,,\\
\Po(\mathfrak{N}) \; 
	&:= \; \big\{ p:\CM\to\CM \text{ complex-valued polynomial} \big\}\,.
\end{align*}

It is remarkable that this result mainly uses the Lemma of Urysohn, the functional calculus and the equality $\|A\|=\max_{\lambda\in\sigma(A)}|\lambda|$ for a bounded normal operator $A$ on a Hilbert space. In the case of self-adjoint or unitary operators, the Lemma of Urysohn is not needed since the geometry of the spectra as a subset of $\CM$ can be used explicitly. More precisely, the spectrum lies in a hyperplane of $\CM$ for self-adjoint operators and the spectrum is contained in the unit sphere of the complex plane for unitary operators. As discussed before, the result that normal continuous vector fields in a continuous field of $C^\ast$-algebras admit spectra that vary continuously was known before, c.f. \cite{Kap51,Ell82,Bel94}. Surprisingly, the converse implication does also hold which is proven in this thesis. Additionally, the characterizations for self-adjoint and unitary operators presented here are new in the literature and they are not based on the theory of $C^\ast$-algebras.

\medskip

Additionally, the continuity of $\Sigma$ is also equivalent to the continuity of the boundaries (spectral edges) in the self-adjoint case by the general considerations for the continuity of compact subsets of $\RM$. The continuity of the boundaries can also be investigated for unitary operators, c.f. discussion in Section~\ref{Chap3-Sect-CharContSpectUnit}.

\medskip

Following the lines of the proof for the characterization of the continuous behavior of the spectra, this delivers a characterization of the H\"older-continuous behavior of the spectra for self-adjoint or unitary operators in terms of the uniform H\"older-continuity of 
$$
\Phi_p:\ts\to[0,\infty)\,,\quad 
	t\mapsto\|p(A_t)\|\,,
$$
for specific $p\in\Po(\AG)$, c.f. Theorem~\ref{Chap3-Theo-CharHolContSpect} and Theorem~\ref{Chap3-Theo-CharHolContSpectUnitary}. Note that this approach is completely new. Astonishingly, the spectra only behave $\alpha/2$-H\"older-continuous with respect to the Hausdorff metric if the norms $\Phi_p$ behave $\alpha$-H\"older-continuous, i.e., the exponent of H\"older-continuity decreases by a factor $1/2$. Conversely, if the spectra behave $\alpha$-H\"older-continuous then the $\Phi_p$'s also behave uniformly $\alpha$-H\"older-continuous in $p$. In the self-adjoint case, the decline of the rate is due to the possibility that gaps of the spectra can close: Clearly, the $\alpha/2$-H\"older-continuity of $\Sigma$ leads also to the $\alpha/2$-H\"older-continuity of the boundaries (spectral edges), c.f. Theorem~\ref{Chap3-Lem-BoundHolCont} and Theorem~\ref{Chap3-Theo-HolConClosGap}. However, if, for certain polynomials $p\in\Po(\mathfrak{S})$, the maps $\Phi_p$ are uniformly $\alpha$-H\"older-continuous, then the gaps that stay open are also $\alpha$-H\"older-continuous, c.f. Theorem~\ref{Chap3-Lem-BoundHolCont}. Thus, the decrease of the rate of convergence is due to the possibility that gaps close. This observation is confirmed by semiclassical analy\-sis for the Almost-Mathieu operator, c.f. \cite{RaBe90,HeSj90} or Example~\ref{Chap3-Ex-AlmostMathieu}. Specifically, the norms of the Almost-Mathieu operator behave Lipschitz-continuous while the spectra only behave $1/2$-H\"older-continuous, c.f. \cite{AvSi85,AvMoSi91,Bel94,JiKr02}.

\medskip

The closing gap condition is also observed in many other models like small perturbations of the discrete Laplacian on $\ell^2(\text{\gls{ZM}})$. More specifically, for $\lambda \geq 0$, let $H_\lambda:\ell^2(\ZM)\to\ell^2(\ZM)$ be defined on $\ell^2(\ZM)$ by 
$$H_\lambda\psi(n) \; = \;
   \psi(n+1)+\psi(n-1)+\lambda V(n)\psi(n) 
$$ 
where $V:\ZM\to\RM$ takes only finitely many values. This is a specific pattern equivariant Schr\"odin\-ger operator, c.f. Theorem~\ref{Chap2-Theo-PESchrOpFinRang}. For $\lambda=0$, the spectrum of $H_\lambda$ is the interval $[-2,+2]$. The prediction provided by the {\em Gap Labeling Theorem} \cite{Bel92,Bel93} yields that for certain potentials gaps may open for $\lambda>0$, i.e., there is a gap tip at $\lambda=0$, c.f. Definition~\ref{Chap3-Def-ClosedGap}. Explicit calculations have been made on several examples, which show the opening of gaps, such as the Fibonacci sequence \cite{SiMo90,DaGo11}, Thue-Morse sequence \cite{Bel90} and the period doubling sequence \cite{BeBoGh91}. Furthermore, for Schr\"odinger operators associated with a compact metric space, the authors of \cite{AvBoDa12} proved that each gap can open by perturbing the potential.

\medskip

If a field of $C^\ast$-algebras is defined via a field of groupoids, the concept of continuous fields of groupoids leads to the continuity of the field of $C^\ast$-algebras (implying the continuity of the spectra), c.f. \cite[Theorem~5.5]{LaRa99}. This result is based on the work by {\sc Rieffel} \cite{Rie89} in the group case and a general result about $C^\ast$-algebras by Blanchard \cite{Bla96}, c.f. Theorem~\ref{Chap4-Theo-LanRamContFieldGroupoid} for details. The concept of continuous fields of groupoids was introduced by {\sc Landsman} and {\sc Ramazan} \cite[Section~5]{LaRa99}. It involves a topological groupoid $\Gamma$, a Hausdorff space $\ts$ and a continuous, open, surjective and $\Gamma$-invariant map $p:\Gamma\to\ts$. Then the preimages $p^{-1}(\{t\})\,,\; t\in\ts\,,$ define a family of disjoint, closed and open subgroupoids of $\Gamma$. 

\medskip

The standard reference for groupoid $C^\ast$-algebras is the book of {\sc Renault} \cite{Renault80} based on his thesis \cite{RenaultThesis78}. Topological groupoids can be seen as a generalization of topological groups with the difference that the composition is only partially defined, see Definition~\ref{Chap2-Def-Groupoid}. Then groups (Example~\ref{Chap2-Ex-Group}) and sets (Example~\ref{Chap2-Ex-Set}) are extremal examples of groupoids. Groupoids defined by equivalence relations (Example~\ref{Chap2-Ex-EquivalRelGroupoid}) are typical examples for groupoids. The main focus of this work is on groupoids induced by topological dynamical systems $(X,G,\alpha)$, c.f. Section~\ref{Chap2-Sect-TransformationGroupGroupoid}. The so called transformation group groupoid $X\rtimes_\alpha G$ was first introduced by {\sc Ehresmann} \cite{Ehr57}. It encodes the geometric and combinatorial properties as well as the dynamics of the underlying dynamical system. Here $X\rtimes_\alpha G$ is the generalization of the semidirect product in the group case to groupoids, c.f. Definition~\ref{Chap2-Def-TransformationGroupGroupoid}. 

\medskip

The reduced $C^\ast$-algebra associated with a topological group $G$ is defined by the left-regular representation of $\Cc_c(G)$. More precisely, elements of $\Cc_c(G)$ are viewed as convolution operators on $L^2(G,\lambda)$ where $\lambda$ is the left-invariant Haar measure on $G$. 

\medskip

A similar construction is used in the groupoid case while several obstacles have to be tackled. For instance, the analog of a Haar measure does not exist for all groupoids. The definition of the so called left-continuous Haar system on a groupoid (that is used nowadays) goes back to {\sc Renault} \cite{RenaultThesis78,Renault80} based on the works of {\sc Westman} \cite{Wes67} and {\sc Seda} \cite{Sed75,Sed76}. If a topological groupoid $\Gamma$ admits a left-continuous Haar system, $\Cc_c(\Gamma)$ is naturally equipped with a structure of a $\ast$-algebra. Like in the group case, the reduced groupoid $C^\ast$-algebra is defined via the left-regular representation of $\Cc_c(\Gamma)$. Another important class of groupoids are topologically amenable groupoids that are intensively studied in the book by {\sc Anantharaman-Delaroche} and {\sc Renault} \cite{AnantharamanRenault00}, see \cite{AnRe01} for a summary. The notion of topologically amenable groupoids has its roots in  the works by {\sc Zimmer} \cite{Zim77a,Zim77b,Zim78} and {\sc Anantharaman-Delaroche} \cite{Ana79} which is an extension of amenability for groups. For amenable groupoids, the corresponding reduced and full $C^\ast$-algebra are isomorphic, c.f. \cite{RenaultThesis78,Renault80,AnantharamanRenault00}. In contrast to the group case \cite{Hul66}, there exist groupoids where the reduced and full $C^\ast$-algebra agree while the groupoid is not topologically amenable, c.f. \cite{Wil15}. In our situation, it is crucial that these $C^\ast$-algebras coincide. This is necessary for the proof that a continuous field of groupoids defines a continuous field of $C^\ast$-algebras, c.f. Remark~\ref{Chap4-Rem-LRUseAmenable}. Compare also Figure~\ref{Chap1-Fig-Strategy2-3-4} for the connection between dynamical systems, groupoids and associated $C^\ast$-algebras.

\medskip

Having the previously discussed theorem about the characterization of the continuous behavior of the spectra in mind, the notion of continuous fields of groupoids is a good tool to prove the continuous behavior of the spectra of operators contained in groupoid $C^\ast$-algebras. A similar tool for the H\"older-continuous behavior of the spectra is still an open question which is under development, c.f. Section~\ref{Chap8-Sect-HolderCont}. Since this work shall provide an access to a broad readership, a detailed elaboration of groupoid $C^\ast$-algebras is presented based on the work \cite{Renault80} in Section~\ref{Chap2-Sect-GroupoidCalgebras}.

\subsection*{Dynamical systems and associated Schr\"odinger operators}

A triple $(X,G,\alpha)$ is called a dynamical system if $X$ is a second countable, compact space, $G$ is a discrete, countable group and $\alpha:X\times G\to X$ is a continuous action of the group $G$ on $X$, c.f. Definition~\ref{Chap2-Def-DynSyst}. The assumption that the group $G$ is discrete is made in order to guarantee that the reduced and full $C^\ast$-algebra of $(X,G,\alpha)$ have a unit which is necessary to get the continuity of the spectra, c.f. Remark~\ref{Chap3-Rem-NecessIdentityContSpect}. The set of $G$-invariant, closed subsets of $X$ is denoted by $\SG(X)$. Elements of $\SG(X)$ are called dynamical subsystems of $X$. Then the space of dynamical subsystems $\SG(X)\subseteq\cs(X)$ endowed with the Hausdorff-topology is compact, second countable and Hausdorff, c.f. Proposition~\ref{Chap2-Prop-SpaDynSyst}.

\medskip

The class of symbolic dynamical systems is of particular interest. A symbolic dynamical system is defined by a finite set $\as$ equipped with the discrete topology and a discrete, countable group $G$. Then the product space $\as^G:=\prod_{g\in G}\as$ is compact and the group $G$ acts on $\as^G$ by the shift 
$$
\alpha:\as^G\times G\to\as^G\,,\;
 (\xi,g)\mapsto \xi\big(g^{-1}\bullet\big)\,.
$$
The set $\as$ is called the alphabet. Then the elements of $\SG\big(\as^G\big)$ are called subshifts instead of dynamical subsystems. 

\medskip

In the specific case that $G=\ZM$, a so called dictionary is associated with an element of $\xi\in\as^\ZM$ being the collection of all finite subwords (patterns) of $\xi$. The concept of dictionary is useful for several applications. A partial generalization of this notion to the group $\ZM^d$ is given in \cite{LindMarcus95}. The notion of a dictionary is extended in this work to general symbolic dynamical systems while the definition of a dictionary is given irrespective of an element of $\as^G$ which is new in the literature. Specifically, a dictionary is defined by a collection of patterns satisfying a heredity and an extensibility condition, c.f. Definition~\ref{Chap2-Def-Dictionary}. The set $\DG$ of all dictionaries associated with an alphabet $\as$ and a group $G$ is naturally endowed with a topology. This new concept of local pattern topology on $\DG$ is investigated here. The local pattern topology is defined by saying that two dictionaries are closed if and only if all patterns agree on a large compact subset of $G$. This definition of the local pattern topology is in accordance with the physical intuition when two solids are close to each other. As it turns out, the space of dictionaries $\DG$ is homeomorphic to $\SG\big(\as^G\big)$. Consequently, one can switch between the notions of dictionaries and subshifts without loosing any topological properties. This approach is used to show the existence of strongly periodic approximations for subshifts (and so for Schr\"odinger operators) for the case $G=\ZM$ (Chapter~\ref{Chap5-OneDimCase}) and $G=\ZM^d$ (Chapter~\ref{Chap6-HigherDimPerAppr}).

\medskip

As mentioned before, every dynamical subsystem $(Y,G,\alpha)$ of $(X,G,\alpha)$ defines a transformation group group\-oid $\Gamma(Y):=Y\rtimes_\alpha G$ with the associated reduced $C^\ast$-algebra $\CG_{red}^\ast\big(\Gamma(Y)\big)$. Elements of $\CG^\ast_{red}\big(\Gamma(Y)\big)$ are called generalized Schr\"odinger operators. These operators are represented by a family of operators $H_Y:=(H_y)_{y\in Y}$ where $H_y\in\Ll\big(\ell^2(G)\big)\,,$ $y\in Y$. In the case of a symbolic dynamical system $(\as^G,G,\alpha)$, pattern equivariant Schr\"odinger operators are specific generalized Schr\"odinger operators that define a dense subalgebra of $\CG^\ast_{red}\big(\Gamma(\as^G)\big)$, c.f. \cite{KePu00} or Theorem~\ref{Chap2-Theo-PaEqDense}. Pattern equivariant Schr\"odinger operators are of particular interest as typically studied Schr\"odinger operators and Jacobi operators on $\ZM$ are specific pattern equivariant Schr\"odinger operators. More precisely, let $(\as^G,G,\alpha)$ be a symbolic dynamical system. Consider a finite $K\subseteq G$ and pattern equivariant functions $p_h:\as^G\to\CM\,,\; h\in K\,,$ and $p_e:\as^G\to\RM$ where $e\in G$ denotes the neutral element of the group $G$. Then the associated pattern equivariant Schr\"odinger operator $H_\xi:\ell^2(G)\to\ell^2(G)$ for a $\xi\in\as^G$ is defined by
$$
(H_\xi\psi)(g) \; 
	:= \; \left( 
			\sum\limits_{h\in K} p_h\big(\alpha_{g^{-1}}(\xi)\big) \cdot \psi(g\,h^{-1}) + \overline{p_h\big(\alpha_{(gh)^{-1}}(\xi)\big)} \cdot \psi(gh)
		\right)
		+ p_e\big(\alpha_{g^{-1}}(\xi)\big) \cdot \psi(g)
$$
where $\psi\in\ell^2(G)$ and $g\in G$, c.f. Definition~\ref{Chap2-Def-SchrOp-l2(G)}. The theory developed here applies to generalized Schr\"odinger operators while the theory is explicit for pattern equivariant Schr\"odinger operators.

\medskip

The spectrum $\sigma(H_Y)$ is equal to the union $\bigcup_{y\in Y} \sigma(H_y)$ for all generalized Schr\"odinger operators $H_Y$, c.f. \cite{Exe14} or Theorem~\ref{Chap2-Theo-SpectrGroupoidCAlg}. If the dynamical system is minimal, the spectrum of a generalized Schr\"odinger operator satisfies $\sigma(H_Y)=\sigma(H_y)=\sigma_{ess}(H_y)$ for each $y\in Y$, c.f. \cite{LeSt03-Algebras,BeLeMaCh14,BeLeMaCh16} or Theorem~\ref{Chap2-Theo-MinCharConstSpectr} and Theorem~\ref{Chap2-Theo-ConstSpectrMinimal}.

\medskip

Next, the universal dynamical system $(\Xun,G,\alpha)$ associated with a dynamical system $(X,G,\alpha)$ is defined in Section~\ref{Chap4-Sect-UnivGroupDynSyst} by following \cite{BeBeNi16}. Here $\Xun$ is the compact subspace of $\SG(X)\times X$ defined by those pairs $(Y,y)$ satisfying $y\in Y$. The universal groupoid $\Gun:=\Xun\rtimes_\alpha G$ corresponding to the universal dynamical system $(\Xun,G,\alpha)$ naturally defines a continuous field of groupoids over the topological space $\SG(X)$. Specifically, $\big(\Gun,p,\SG(X)\big)$ is a continuous field of groupoids where $p:\Gun\to\SG(X)$ is a projection. The fiber groupoids $p^{-1}\big(\{Y\}\big)$ turn out to be measurable isomorphic to $\Gamma(Y):=Y\rtimes_\alpha~G$. Hence, the $C^\ast$-algebras $\CG^\ast_{red}\big(\Gamma(Y)\big)$ and $\CG^\ast_{red}\big(p^{-1}\big(\{Y\}\big)\big)$ are isomorphic for each $Y\in\SG(X)$, c.f. Proposition~\ref{Chap4-Prop-UnivGroupContField}. Altogether, $\big(\big(\CG^\ast_{red}(\Gamma(Y))\big)_{Y\in\SG(X)},\Upsilon\big)$ with $\CG^\ast_{red}(\Gun)\subseteq\Upsilon\subseteq\prod_{Y\in\SG(X)}\CG^\ast_{red}(\Gamma(Y))$ defines a continuous field of $C^\ast$-algebras where $\Upsilon$ denotes the algebra of continuous vector fields. This leads to the continuous behavior of the spectra belonging to generalized Schr\"odinger operators as described below. 

\medskip

Before we are going into details, a short historical review about Solid State Physics is pres\-ented as well as the known results about the spectral properties of Schr\"odinger operators associated with solids. These results motivate our approach followed in this thesis.

\subsection*{Solid State Physics: a short historical review}

The theory of Solid State Physics goes back to the works \cite{Dru00a,Dru00b} of {\sc Drude} in 1900 where he was the first who described the transport properties of metals in terms of the theory of electrons. Later the microscopic structure of solids was experimentally observed by Laue in 1912 by sending $X$rays through a solid and studying the diffraction pattern. This experiment provides the first verification of crystal structures of the atoms. The reader is referred to \cite{HoBrTeWeSe92} for a historical background on Solid State Physics. 

\medskip

Solids can be classified in terms of the ordering of their atoms and mole\-cules in the space. For instance, the periodic and the Anderson model are considered. In the first case, the potential generated by the material is completely ordered and periodic with respect to the translation. In contrast, the atomic types and the locations of the atoms are given by i.i.d. random variables in the Anderson model, c.f. \cite{And58}. Solids with long range aperiodic order, i.e., ordered but not periodic, lie in between these two models. The examination of such materials has soared up since their discovery by {\sc Shechtman} \cite{SBGC84} in 1982. Such solids are called quasicystals nowadays. Since their discovery, a large machinery started creating different types of quasicrystals and studying their properties from the physical and mathematical point of view. The review articles \cite{Poo92,Ber94} provide known experimental properties of the electronic motion in quasicrystals whereas \cite{May94,Sir94} supply  theoretical approaches from the physical point of view.

\medskip

The diffraction spectrum of a quasicrystal is pure point revealing long range order of the atomic configurations. At the same time the diffraction has a symmetry which is incompatible with translation invariance of the solid which was a real surprise since the materials were created by metals such as copper, aluminum and iron, c.f. \cite{BeMaCyKlLa93,May94} or the summary in \cite[Section~2.2.1]{RiethThesis95}. The study of the corresponding Schr\"odinger operators became difficult since the Floquet-Bloch theory \cite{Flo1883,BlochThesis29} was not longer applicable. Consequently, physicists started with numerical simulations which had various difficulties corresponding to the choice of boundary conditions. In the nineties, {\sc Benza} and {\sc Sire} \cite{BeSi91} overcame this problem by the choice of suitable periodic approximations. However, a suitable mathematical ground was still missing. Stable quasicrystals were also created by changing the concentration for an alloy consisting of crystals, c.f. \cite{Rag10,TsInMa87}. Thus, periodic approximations arise naturally for aperiodic systems from the physical perspective. The work \cite{BeBeNi16} delivers a necessary and sufficient condition so that the spectra behave continuously which is also presented here. More specifically, the strong connections between the behavior of the underlying dynamical systems and the corresponding spectra of the associated Schr\"odinger operators is shown. With this at hand, periodic approximations of Schr\"odinger operators are defined by periodic approximations of the dynamical systems, see discussion below.

\medskip

From the mathematical point of view, mainly one-dimensional systems can be treated analytically by using the so called transfer matrices and the trace map formalism, see discussion below. As a consequence of semiconductor technology, physicists can create artificial structures mimicking situations of the theoreticians like one-dimensional models, see e.g. \cite{TanEtAl14} for the Fibonacci sequence. As it turns out the theoretical predictions and the experimental results agree on large accuracy.

\medskip

The first step for a mathematical theory on solids is to model solids in a suitable way. Atomic nuclei are surrounded by their electrons. This leads to a repulsion of the atoms in short distance due to the Pauli principle prohibiting electrons to occupy the same quantum state. 
A stationary description of solids does not allow arbitrary large holes if cracks are excluded. Altogether, there is a minimal and a maximal distance between atoms in a solid, c.f. \cite{BeRaSh10,BrLi11}. Consequently, the mathematical description of solids is usually modeled by Delone sets, c.f. the surveys \cite{Bel86,Bel02,Bel03} for more background. In this work, we focus on dynamical systems $(X,G,\alpha)$ where the group $G$ plays the role of a lattice which is a special case of a Delone set.

\medskip

First contributions in the mathematical study of random and quasiperiodic potentials is due to the theory of small divisor analysis by {\sc Kolmogoroff}, {\sc Arnold} and {\sc Moser} \cite{Kol53,Kol54,Arn61,Arn62,Mos62,Gal83,Chi08} in classical Hamiltonian systems, the mathematical physicists \cite{GoMoPa77} as well as multiscale analysis proposed by {\sc Fr\"ohlich} and {\sc Spencer} \cite{FrSp83}. On this ground, analytic results on the nature of the spectrum of Schr\"odinger operators with quasiperiodic potentials \cite{BLT83,Bel85} or random potentials \cite{FrMaScSp85} were proven, c.f. also the books \cite{CyconFroeseKirschSimon87,CarmonaLacroix90,PasturFigotin92}. Later, the methods of semiclassical analysis was used in combination with the study of electrons in a magnetic field, c.f. \cite{HeSj87,HeSj89,BaBeRa90,BeVi90,Bel94,CoPu15}.

\medskip

From the physical point of view, the theory of Schr\"odinger operators should be independ\-ent of the location of the origin. This leads to the study of random Schr\"odinger operators where the set of possible configurations is modeled by a probability space equipped with a dynamic. An impressive number of specific examples where analyzed for their spectral properties like pure point spectrum, purely singular continuous spectrum or Cantor spectrum, see e.g. the textbooks \cite{CyconFroeseKirschSimon87,Kirsch89,CarmonaLacroix90,PasturFigotin92,BaakeMoody00} and the surveys \cite{Sim82} for almost-periodic potentials and \cite{Bel86,Sut95,Dam00Gord,Stollmann01,Bel02,Baa02,Bel03,Dam07,DaLiQu14,BaDaGr15,DaEmGo15,KellendonkLenzSavinien15} for quasicrystals. Different classes of potentials give rise to different types of spectral properties. For instance, periodic configurations lead to band spectrum that is absolutely continuous with the help of the Floquet-Bloch theory \cite{Flo1883,BlochThesis29}. In the Anderson model, i.i.d random variables determine the values of the potential, c.f. \cite{And58}. The works \cite{FrMaScSp85,CaKlMa87} prove that Schr\"odinger operators with random potentials have pure point spectrum. The class of aperiodic potentials \-associated with quasicrystals yields spectral properties that were not expected to appear before like Cantor spectrum of Lebesgue measure zero, see discussion below.

\medskip

{\sc Connes} \cite{Connes94} provides powerful tools by introducing the concept of non-commutative geometry. This leads also to a conceptual presentation of the mathematical theory of aperiodic solids, c.f. also \cite{Bel86,Bel92,Bel03}.

\subsection*{Quasicrystals}

In the study of quasicrystals, different mathematical areas are combined like dynamical systems, spectral theory, operator algebras and non-commutative geometry. It splits up into the study of tilings and their dynamics, the diffraction and the spectral theory of the associated Schr\"odinger operators. This work focuses on the spectral theory of the corresponding Schr\"odinger operators while the dynamical part plays a crucial role. The following questions do naturally arise in the spectral theory. 
\begin{description}
\item[(Q1)\label{(Q1)}] How do the spectral properties of a Schr\"odinger operator depend on the underlying combinatorics and geometry?
\item[(Q2)\label{(Q2)}] Are the spectral properties stable under deformations of the underlying system?
\end{description}
This work focuses on the change of the spectra as a set in dependence of the combinatorics and geometry of the corresponding system. It provides an answer to the second question \nameref{(Q2)} for the behavior of the spectrum as a compact subset of $\CM$. Meanwhile, it delivers a tool to approximate systems associated with quasicrystals by periodic systems where the spectra can be computed with the help of the Floquet-Bloch theory. As previously discussed, this fits also to result in \cite{BeSi91} and that quasicrystals can be constructed by periodic approximations in physics, c.f. \cite{Rag10,TsInMa87}.

\medskip

As previously discussed, a quasicrystal is non-periodic with respect to the translation while the solid is still ordered. Up to now, there does not exist a precise mathematical definition of a quasicrystal. Typically, one supposes that the associated dynamical system is non-periodic and minimal. Minimality has strong implications for the spectrum as previously discussed. In the literature, diverse classes of quasicrystal are analyzed. For instance, system defined by a substitution or a cut and project model are considered. The case of substitutions is also analyzed here in Chapter~\ref{Chap6-HigherDimPerAppr} to define periodic approximations explicitly while the cut and project model is studied in a future project, c.f. Section~\ref{Chap8-Sect-PerApprox}.

\medskip
 
First considerations of Schr\"odinger operators associated with quasicrystals are provided in \cite{KoKaTa83,OPRSS83}. {\sc S{\"u}t{\H{o}}} \cite{Sut87,Sut89} was the first who proved that the spectrum associated with the Fibonacci subshift is a Cantor set of Lebesgue measure zero. The Fibonacci subshift is a one-dimensional system which is considered as the typical quasicrystal. Later, {\sc Bellissard}, {\sc Bovier} and {\sc Ghez} \cite{BeBoGh91,BoGh93} showed the absence of point spectrum and Cantor spectrum of Lebesgue measure zero for a large class of primitive substitutions by using the trace map formalism. This result was later generalized in \cite{LiTaWeWu02} to all primitive substitutions. Specifically, Sturmian potentials were intensively studied in the literature during the past including the Fibonacci system. The spectrum of such a Schr\"odinger operator coincides with the set where the Lyapunov exponent vanishes for irrational numbers, c.f. \cite{BIST89}. The authors of \cite{DaLiQu15} proved that the spectrum of a Schr\"odinger operator has zero Lebesgue-measure and is purely singular continuous for all pattern Sturmian potentials that arise from Toeplitz sequences. In 2002, {\sc Lenz} \cite{Len02} provides a characterization of the spectrum as a set associated with Schr\"odinger operators by properties of the Lyapunov exponent for a large class of strictly ergodic subshifts over a finite alphabet. In this setting, the author proved that the spectrum is a Cantor set of Lebesgue measure zero for aperiodic subshifts if the corresponding transfer matrices are uniform for all energies. This includes all subshifts satisfying uniform positivity of weights like Sturmian dynamical systems and primitive substitutions. More precisely, it covers previous results like \cite{BoGh93} but it also covers new subshifts like the Rudin-Shapiro substitution. Recently, this result was extended to the class of Jacobi operators in \cite{BePo13}. The result of \cite{BePo13} applies also for the Laplacian on the Schreier graphs arising from Grigorchuk's group acting on the boundary of an infinite binary tree, c.f. \cite{GrLeNa14,GrLeNa15}.

\medskip

In contrast to Cantor spectrum of Lebesgue measure zero, the absence of eigenvalues is a more delicate question. In addition to the transfer matrices techniques, the formalism of trace maps turned out be very useful in the one-dimensional case, c.f. \cite{BeBoGh91,BoGh93}. In 1998, {\sc Daminik} \cite{Dam98} verified the almost sure absence of eigenvalues for potentials with a four block structure. Upper bounds on the growth of solutions as well as the fact that the point spectrum is empty are proven by {\sc Damanik}, {\sc Lenz} and {\sc Killip} \cite{DaLe99I,DaLe99II,DaLeKi00}. 

\medskip

During the last years also the fractal dimension and the transport exponent of the spectrum were studied. {\sc Fan}, {\sc Liu}, {\sc Qu} and {\sc Wen} \cite{FaLiWe11,LiQuWe14,LiQu15} determined the fractal dimensions of the spectrum associated with Sturmian Hamiltonians in depend\-ence of their continued fraction expansion. These results are based on the general works \cite{FeWeWu97,HuRaWeWu00}. Moreover, for specific examples as the Thue-Morse sequence, the Hausdorff dimension can be bounded from below uniformly in terms of the coupling constant, c.f. \cite{LiQu15}. Recently, the authors of \cite{DaGoLiQu15} proved upper and lower bounds for the transport exponents corresponding to Sturmian Hamiltonians.

\medskip

Among others the results on the Hausdorff dimensions use that the operators can be approximated monotonically by periodic Schr\"odinger operators. Since the Lebesgue meas\-ure of the periodic operators can be estimated in high accuracy, this allows to estimate the Hausdorff dimension of the aperiodic Schr\"odinger operators of specific systems. This thesis shows the convergence of the spectra for periodic operators for a large class of one-dimensional systems including higher dimensional subshifts in $\ZM^d$ while the formalism of transfer matrices and trace maps is not used, c.f. Chapter~\ref{Chap5-OneDimCase}, Chapter~\ref{Chap6-HigherDimPerAppr}, Chapter~\ref{Chap7-Examples}.

\medskip

The reader is also referred to the surveys and lecture notes \cite{Bel86,Sut95,Dam00Gord,Stollmann01,Bel02,Baa02,Bel03,Dam07,DaLiQu14,BaDaGr15,DaEmGo15,KellendonkLenzSavinien15} for more background on the known results about spectral properties of Schr\"odinger and Jacobi operators.

\medskip

All the previously mentioned results are only valid in the one-dimensional case since they are based on the techniques of transfer matrices and trace maps. {\sc Damanik}, {\sc Gorodetski} and {\sc Solomyak} provide first results for two-dimensional systems like the square Fibonacci Hamiltonian if the Schr\"odinger operator is separable, c.f. \cite{DaGoSo13,DaGo16}. More specifi\-cally, the spectrum of the two-dimensional model turns out be the sum of the spectra of one-dimensional Schr\"odinger operators. With this at hand, the authors show that the spectrum has positive Lebesgue measure while the density of state measure is purely singular. This is verified by the general theory of sums of Cantor sets of $\RM$.

\medskip

For all of these results, the transfer matrices play a crucial role. This formalism was already used by {\sc Floquet} \cite{Flo1883} in 1883 to study ordinary differential equations with periodic coefficients. The main idea is that a solution can be generated with the knowledge of two initial parameters. The results on the fractal dimensions and the absence of eigenvalues rely additionally on the concept of trace maps. The trace map encodes dynamical properties of the system in algebraic equations. In combination, this gives a strong tool to study the spectral properties in very detail. On the other hand, this formalism is limited to the one-dimensional case besides for specific higher dimensional systems that can be decomposed into one-dimensional systems. 

\medskip

In the periodic case, the break through came with the thesis of {\sc Bloch} \cite{BlochThesis29} who used general theory of self-adjoint operators that was independent of the dimension. This created the ground for the band theory that constitutes the fundamental aspects in Solid State Physics. Later, {\sc Wannier} \cite{WannierThesis37} generalized the Fourier transform by using the translation invariance of the underlying system. This unitary transformation is nowadays called  Wannier transform which represents Schr\"odinger operators as a direct integral providing a deeper insight to the spectral properties.

\medskip

From the conceptually point of view, the known results about aperiodic Schr\"o\-dinger operators in one-dimensions show that the combination of dynamical systems with an algebraic formulation delivers a deeper insight in the spectral nature of the systems. This is also the philosophy of this thesis as described next. 

\medskip

A major breakthrough came with the works \cite{Kel95,Kel97} by {\sc Kellendonk} who proved that aperiodic solids and tilings can be treated on the same ground. Furthermore, it turned out that $C^\ast$-algebraic techniques are useful for the study of tilings. In 2003, {\sc Kellendonk} and {\sc Putnam} introduced the notion of pattern equivariant functions, c.f. \cite{Kel03}. Schr\"odinger operators with a potential arising by pattern equivariant functions are called pattern equivariant Schr\"odinger operators, c.f. Section~\ref{Chap2-Sect-PattEqSchrOp}. It turns out that the corresponding operator algebra generated by the pattern equivariant Schr\"odinger opera\-tors coincides with the groupoid $C^\ast$-algebra where the transformation group groupoid is naturally defined by the dynamical system, c.f. \cite{KePu00,LeSt03-Algebras,Whittaker05,StarlingThesis12,Sta14} or Theorem~\ref{Chap2-Theo-PaEqDense} and Theorem~\ref{Chap2-Theo-PESchrOpFinRang}. 

\medskip

Altogether, it is crucial to note for the further considerations that a pattern equivariant (generalized) Schr\"odinger operator can be represented as a continuous function on the transformation group groupoid defined by the dynamical system. This continuity property essentially leads to the continuity of the map
$$
\Sigma_H:\SG(X)\to\ks(\CM)\,,\quad
	Y\mapsto\sigma(H_Y)\,,
$$
for each suitable generalized Schr\"odinger operators $H_Y\,,\; Y\in\SG(X)$, see Theorem~\ref{Chap4-Theo-TransGroupContFieldCAlg} and Corollary~\ref{Chap4-Cor-ContSpectrContSectUnivGroup} for more details. The proof is based on the concept of universal groupoid and continuous fields of groupoids as previously mentioned. A similar strategy can be followed for Schr\"odinger operators associated with Delone sets as described in \cite{BeBeNi16}. In view of this continuity result, it is naturally to expect a H\"older-continuity of $\Sigma$ if the function that defines the Schr\"odinger operator is H\"older-continuous on the groupoid in a suitable sense. This problem remains open, see discussion in Section~\ref{Chap8-Sect-HolderCont}. A connection of the different concepts that are used for the continuity of $\Sigma_H$ is sketched in Figure~\ref{Chap1-Fig-Strategy5-6}. 

\medskip

Due to the fact that generalized Schr\"odinger operators are continuous functions on the dynamical system, the Lemma of Urysohn yields the equivalence of the convergence of dynamical subsystems and the convergence of the spectra of all associated Schr\"odinger operators. More specifically, a sequence of dynamical subsystems $Y_n\in\SG(X)\,,\; n\in\text{\gls{NM}}\,,$ converges to $Y\in\SG(X)$ in the Hausdorff-topology if and only if $\lim_{n\to\infty}\sigma\big(H_{Y_n}\big)=\sigma\big(H_Y\big)$ holds (with respect to the Hausdorff metric on $\ks(\RM)$) for each self-adjoint generalized Schr\"odinger operator, c.f. Theorem~\ref{Chap4-Theo-CharDynSystConvSpectr} and Theorem~\ref{Chap4-Theo-CharSubshiftConvSpectr}.

\medskip

Altogether, Section~\ref{Chap4-Sect-ContSpecGenSchrOp} provides an answer for Question \nameref{(Q2)} while only the spectrum as a set is analyzed and not its spectral measure. It is remarkable that these results are completely independent of the dimension and it applies to much more general operators than considered before. 

\subsection*{Periodic approximations of Schr\"odinger operators}

\begin{figure}
\begin{center}
\includegraphics[height=0.95\textheight]{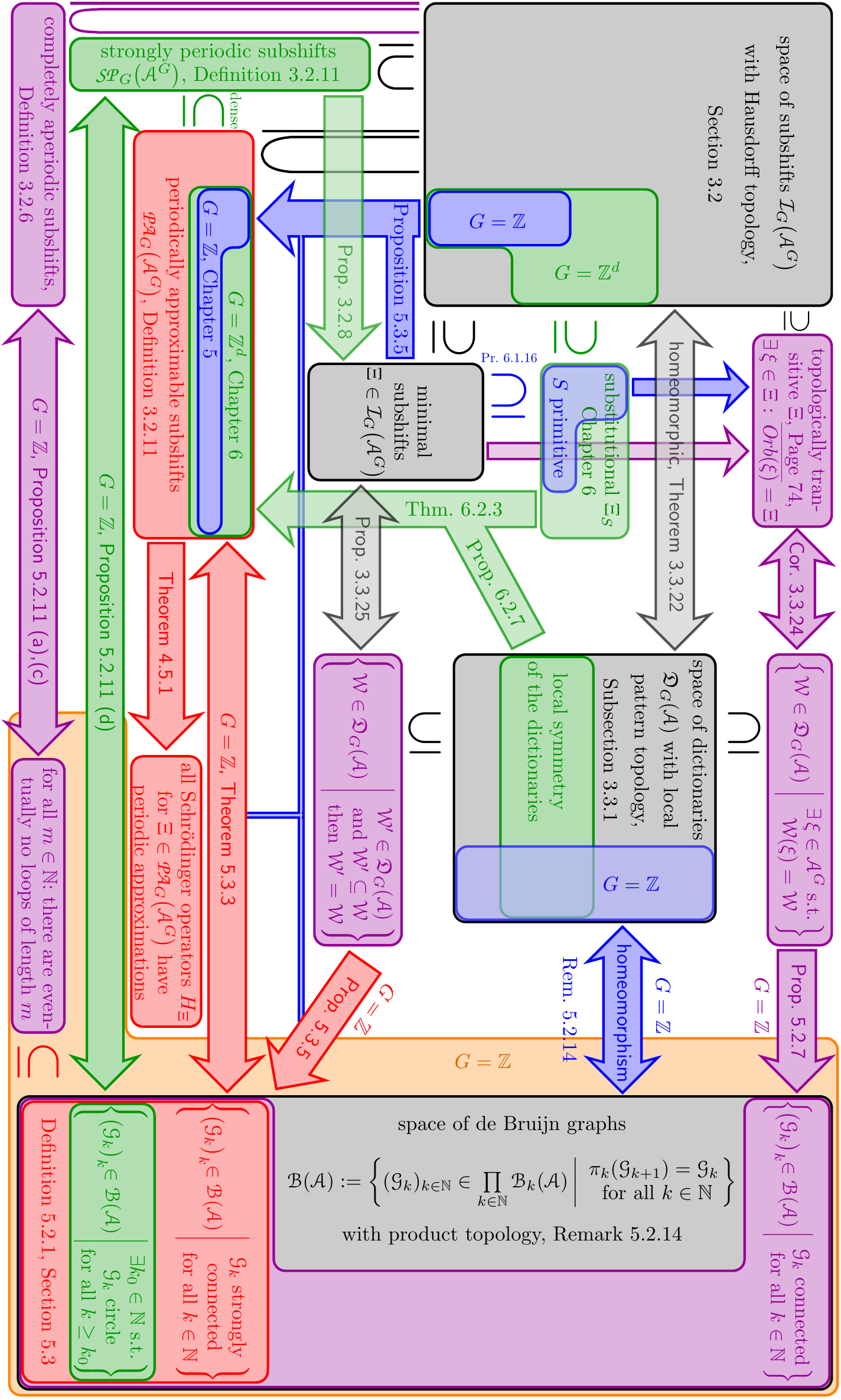}
\caption{Overview and connections between Chapter~\ref{Chap5-OneDimCase} and Chapter~\ref{Chap6-HigherDimPerAppr}.}
\label{Chap1-Fig-Strategy5-6}
\end{center}
\end{figure}

Based on the previously discussed results, an approximation for a generalized (pattern equivariant) Schr\"odinger operator can be defined by approximating the corresponding dynamical system, c.f. Theorem~\ref{Chap4-Theo-PeriodicApproximations}. More precisely, non-periodic Schr\"odinger opera\-tors can be approximated by periodic ones if the corresponding dynamical system is approximable by strongly periodic dynamical systems. This is the fundamental idea that is followed in Chapter~\ref{Chap5-OneDimCase}, Chapter~\ref{Chap6-HigherDimPerAppr} and Chapter~\ref{Chap7-Examples}. Compare Figure~\ref{Chap1-Fig-Strategy5-6} for an overview of the results. As previously mentioned, periodic approximations of aperiodic Schr\"odinger operators (corresponding to quasicrystals) are interesting from the physical and mathemati\-cal perspective.

\medskip

First the one-dimensional case is treated, i.e., $G=\ZM$. A subshift $\Xi$ of $(\as^\ZM,\ZM,\alpha)$ is approximable by periodic subshifts if and only if the so called de Bruijn graphs are strongly connected, c.f. Theorem~\ref{Chap5-Theo-ExPerAppr}. In addition, periodic approximations are defined by closed paths in these graphs. This new approach yields that all minimal subshifts are periodically approximable, c.f. Proposition~\ref{Chap5-Prop-MinAllPath}. A detailed discussion of the relation between dynamical properties and combinatorial properties of the de Bruijn graphs are discussed in Section~\ref{Chap5-Sect-DBGraphs} as well as the role of the branching vertices. The higher dimensional analog of the de Bruijn graphs are the so called Anderson-Putnam complexes, c.f. \cite{AnPu98} as well as \cite{Gah02,BeGa03,Sad03,BeBeGa06,Sad08}. There periodic approximations correspond to tori in the Anderson-Putnam complex, c.f. Section~\ref{Chap8-Sect-APComplex}. 

\medskip

If $G=\ZM^d$, new sufficient conditions are provided so that a subshift defined by a (primi\-tive) substitution is periodically approximable, c.f. Chapter~\ref{Chap6-HigherDimPerAppr}. This corresponds to the existence of local symmetries in the patterns of the subshifts, c.f. Proposition~\ref{Chap6-Prop-SuffPerApprZ2SymmPatt} and Remark~\ref{Chap6-Rem-SuffPerApprZ2SymmPatt}. The main idea is to define a periodic approximation by applying the substitution iteratively to a suitable strongly periodic element of $\as^{\ZM^d}$ and using the previously developed description of a subshift by a dictionary. Thus, the periodic approximations are given explicitly.

\medskip

According to the definition of the local pattern topology for dictionaries and Theorem~\ref{Chap2-Theo-Shift+DictSpace}, a subshift can be approximated by strongly periodic subshifts only if all the subshifts of finite type
$$
\Xi_n \;
	:= \; \big\{
		\xi\in\as^G\;|\; \ws(\xi)\cap\as^{[K_n]}\subseteq \ws(\Xi)\cap\as^{[K_n]}
	\big\}\,,
	\qquad
	n\in\NM\,,
$$ 
associated with $\Xi$ contain strongly periodic elements where $(K_n)_{n\in\NM}$ is an exhausting sequence of the group $G$ (Definition~\ref{Chap2-Def-ExhaustSeq}). The question whether a subshift of finite type contains strongly periodic elements is a current question in the field of dynamical systems, c.f. \cite{Pia08,Hoc09,SiCo12,Coh14,CaPe15}. It is known \cite{Fio09} that there are sufficient conditions for one-dimensional subshifts of finite type to contain strongly periodic elements. The proof relies on the de Bruijn graphs associated with a subshift of finite type. More specifically, closed paths in the de Bruijn graphs give rise to strongly periodic elements in the subshift of finite type, c.f. Theorem~\ref{Chap5-Theo-ExPerAppr}. In higher dimensions, the existence of strongly periodic elements in a subshift of finite type is more difficult, c.f. \cite{Hoc09}. In many situations, the set of strongly periodic elements in a subshift of finite type $\Xi$ is dense if $\Xi$ contains already one strongly periodic element, c.f. \cite{Lig03,SiCo12} or Remark~\ref{Chap2-Rem-SubsFinTypPerSpectr}.

\medskip

This question is also related to the existence of so called Wang tiles. In 1961, {\sc Wang} \cite{Wan61} was the first who asked the question whether, for a given finite set of tiles, there is a periodic tessellation of the plane with these tiles. {\sc Berger} \cite{Ber66}, a student of {\sc Wang}, provides an example of tiles that do not tile the plane in a periodic way. Such tilings are called nowadays {\em Wang tiles}. Even in the case of a lattice $\ZM^d$ for $d\geq 3$ there exist tilings, i.e., $\xi:\ZM^d\to\as$, that do not allow any periodic tiling, c.f. \cite{CuKa95}. This shows that there are subshifts that are not periodically approximable, c.f. Example~\ref{Chap5-Ex-DeBruijnNotStrongConnect} and Example~\ref{Chap5-Ex-deBruijnEventuallyNotStrongConn}.

\medskip

These results are applied to specific one-dimensional systems that are intensively studied in the literature as well as to higher dimensional subshifts in Chapter~\ref{Chap7-Examples}. Specifically, the table substitution and the Sierpinski carpet substitution over the group $\ZM^2$ are studied which both turn out to be periodically approximable.

\subsection*{Summary}

\begin{figure}
\begin{center}
\includegraphics[height=0.95\textheight]{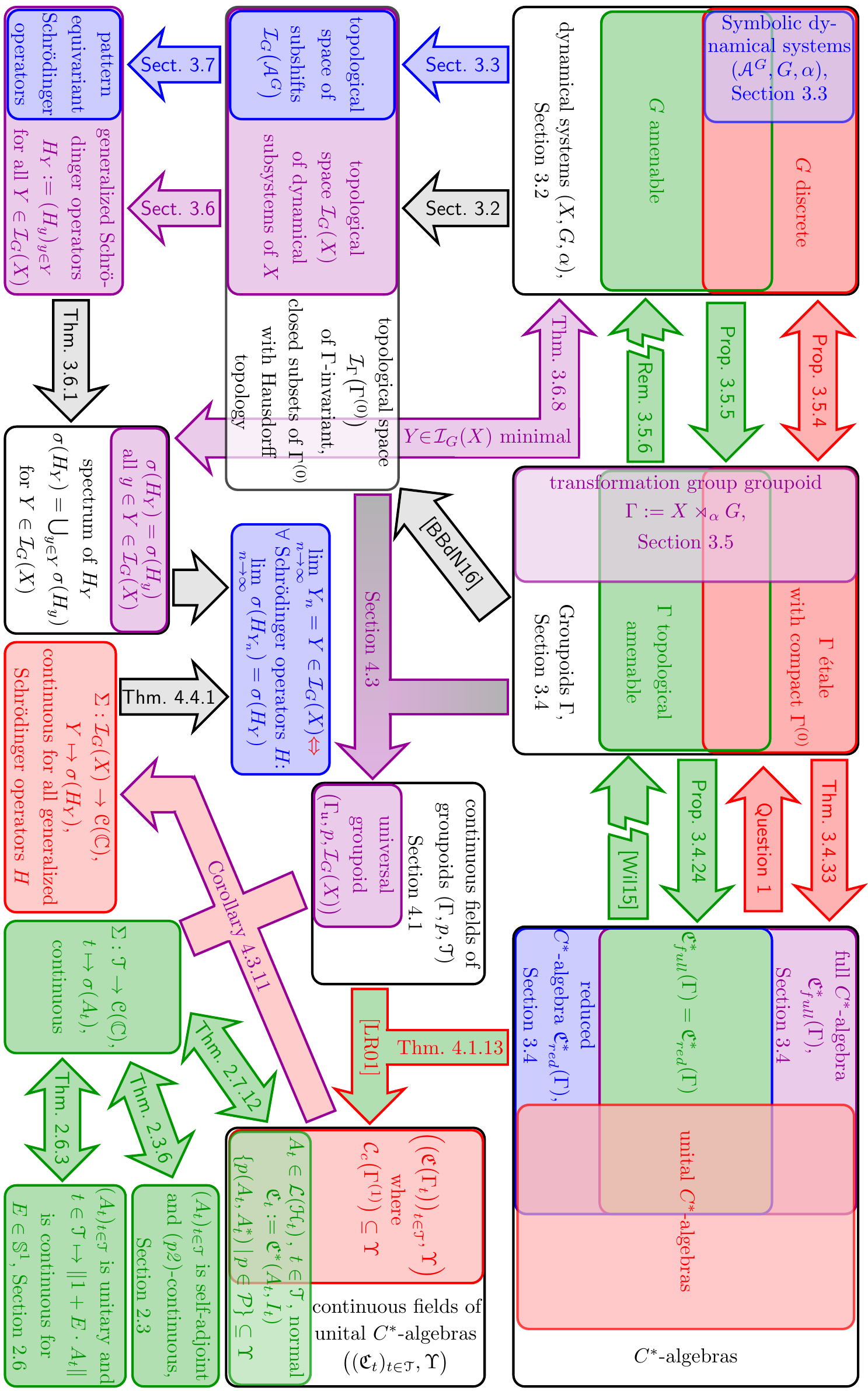}
\caption{Overview and connections between Chapter~\ref{Chap3-SpectAppr}, Chapter~\ref{Chap2-SchrOpDynSyst} and Chapter~\ref{Chap4-ToolContBehavSpectr}.}
\label{Chap1-Fig-Strategy2-3-4}
\end{center}
\end{figure}

In the following, the organization of this work is presented. For convenience of the reader, a more detailed overview for the content of Chapter~\ref{Chap3-SpectAppr}, Chapter~\ref{Chap2-SchrOpDynSyst} and Chapter~\ref{Chap4-ToolContBehavSpectr} is given in Figure~\ref{Chap1-Fig-Strategy2-3-4} while the diagram in Figure~\ref{Chap1-Fig-Strategy5-6} presents the content of Chapter~\ref{Chap5-OneDimCase} and Chapter~\ref{Chap6-HigherDimPerAppr}.

\medskip

Chapter~\ref{Chap3-SpectAppr} is devoted to the general theory of the continuous behavior of the spectra based on joint work \cite{BeBe16}. The relations of the Hausdorff metric, Vietoris-topology and Fell-topology on $\CM$ are studied in Section~\ref{Chap3-Sect-VietorisContinuity}. Furthermore, the new concept of continuous boundaries for families of closed subsets of $\RM$ is investigated and its connections to the Fell- and Vietoris-topology of these families are analyzed. In the case of self-adjoint (Section~\ref{Chap3-Sect-CharContSpectSelfAdj}) and unitary (Section~\ref{Chap3-Sect-CharContSpectUnit}) operators, the continuous variation of the spectra with respect to the Hausdorff metric is characterized by the continuity of certain norms of the operators. This leads also to a characterization of the H\"older-continuous behavior of the spectra while the rate of convergence decreases if spectral gaps close, c.f. Section~\ref{Chap3-Sect-CharHolContSpect}. In Section~\ref{Chap3-Sect-RelP2ContOthTop}, the \pt-topology on $\Ll(\hs)$ is analyzed for its relations to the operator norm topology and strong operator norm topology. Section~\ref{Chap3-Sect-CharContSpecUnbouSelfAdjOp} deals with the continuous variation of the spectra for unbounded self-adjoint operators. The case of normal operators is treated in Section~\ref{Chap3-Sect-ContFieldCAlg} by using the machinery of continuous fields of $C^\ast$-algebras.

\medskip

In Chapter~\ref{Chap2-SchrOpDynSyst}, the framework for Schr\"odinger operators induced by a dynamical system is introduced. Specifically, the basic concepts of dynamical systems (Section~\ref{Chap2-Sect-DynSystGroupoid}) and symbolic dynamical systems (Section~\ref{Chap2-Sect-SymbDynSyst}) are investigated which includes the new concept of local pattern topology on the set of dictionaries. Following the lines of \cite{Renault80}, an introduction for the theory of groupoids and their associated $C^\ast$-algebras is provided in Section~\ref{Chap2-Sect-GroupoidCalgebras}. A detailed discussion of the history and the proofs of the fundamental assertions are offered since this work is addressed also to the readers who are not familiar with the theory of groupoid $C^\ast$-algebras. Section~\ref{Chap2-Sect-TransformationGroupGroupoid} deals with those groupoids that are defined via a dynamical system which are mainly studied in this thesis. Finally, the notions of generalized Schr\"odinger operators and pattern equivariant Schr\"odinger operators are introduced in Section~\ref{Chap2-Sect-SchrOp} and Section~\ref{Chap2-Sect-PattEqSchrOp}. Spectral properties associated with dynamical properties are analyzed there.

\medskip

Chapter~\ref{Chap4-ToolContBehavSpectr} is devoted to presentation of the concept of continuous fields of groupoids as a tool to prove the continuous behavior of the spectra of operator families. Specifically, the assertion \cite[Theorem~5.5]{LaRa99} is presented in Section~\ref{Chap4-Sect-ContFieldGroupCAlg}. Section~\ref{Chap4-Sect-GroupoidIsomorphism} recalls basic relations between a measure preserving groupoid isomorphism and the $\ast$-isomorphism of the corresponding groupoid $C^\ast$-algebras that are used later. The construction of the universal dynamical system and the universal groupoid associated with a dynamical system $(X,G,\alpha)$ are investigated in Section~\ref{Chap4-Sect-UnivGroupDynSyst}. It is proven that the universal groupoid $\Gun$ defines a continuous field of groupoids which fibers over the space $\SG(X)$. This leads to the fundamental result that $\Sigma_H:\SG(X)\to\ks(\CM)$ is continuous for all suitable generalized Schr\"odinger operators $H$, c.f. Corollary~\ref{Chap4-Cor-ContSpectrContSectUnivGroup}. Based on this, the characterization of the convergence of dynamical subsystems by the convergence of the spectra of generalized Schr\"odinger operators is treated in Section~\ref{Chap4-Sect-ContSpecGenSchrOp}. Finally, these results are applied to prove the existence of periodic approximations for generalized Schr\"odinger operators in Section~\ref{Chap4-Sect-ApprPerAppr}. 

\medskip

Chapter~\ref{Chap5-OneDimCase} deals with symbolic dynamical systems over the group $\ZM$. The property that a subshift in $(\as^\ZM,\ZM,\alpha)$ can be approximated by periodic subshifts is characterized by the fact that all the associated de Bruijn graphs are strongly connected. Furthermore, it is shown that closed paths in the de Bruijn graphs give rise to periodic approximations of the subshifts. The structure of the finite de Bruijn graphs is analyzed with respect to the non-periodicity of the system and the subword complexity. Here the so called branching vertices of the de Bruijn graphs play a crucial role. As it turns out, all minimal subshifts are periodically approximable, c.f. Proposition~\ref{Chap5-Prop-MinAllPath}. 

\medskip

In Chapter~\ref{Chap6-HigherDimPerAppr} we study the class of subshifts defined by (primitive) substitutions in $\ZM^d$. While in the literature general references of higher dimensional substitutional systems are sparse, Section~\ref{Chap6-Sect-PrimBlockSubst} provides the construction and basic properties of a subshift defined by a primitive substitution. On this basis, sufficient conditions for subshifts of $(\as^{\ZM^d},\ZM^d,\alpha)$ being periodically approximable are proven in Section~\ref{Chap6-Sect-PerApprPrimBlockSubst}. It turns out that local symmetries of the patterns yield to this condition, c.f. Proposition~\ref{Chap6-Prop-SuffPerApprZ2SymmPatt}. Like in the one-dimensional case, the connections between dictionaries and subshifts is strongly used there. The advantage of substitutional subshifts is that not only the existence of periodic approximations is proven. Additionally, an approximation is defined iteratively by applying the substitution.

\medskip

This thesis is finished with Chapter~\ref{Chap7-Examples} collecting several examples that were studied in the literature. Some of the known results are covered and extended by the theory developed in the previous chapters with different methods. Most of the examples arise from a substitution which allows to define the periodic approximations iteratively based on the results of Chapter~\ref{Chap6-HigherDimPerAppr}. Finally, also two examples over the group $\ZM^2$ are analyzed. More precisely, periodic approximations are provided for the minimal, aperiodic table substitution and the aperiodic Sierpinski carpet substitution which is not minimal.

\cleardoublepage

\chapter{Characterization of the (H\"older-)continu\-ous behavior of the spectra}
\label{Chap3-SpectAppr}
\stepcounter{section}
\setcounter{section}{0}

The continuous dependence of the spectra is interesting in different areas of mathematics and physics. There are two points of view that arise naturally. On the one hand, let $H$ be an operator that is defined by some model. How can the underlying system be changed so that the spectral properties do not vary too much? On the other hand, it is common to approximate such an operator whenever the operator $H$ is difficult to handle. It is then necessary that the approximation process preserves sufficiently many spectral properties of the original object. For instance, do the spectra as a set behave well with respect to the Hausdorff metric on $\CM$ in the approximation process? Or even more quantitative, how fast do the spectra converge in the Hausdorff metric if they converge? These questions are studied in this chapter, based on joint work with {\sc J. Bellissard} \cite{BeBe16} for self-adjoint, unitary and normal operators. Depending on the class of the operator, sufficient and necessary conditions are provided which are new in the literature.

\medskip

More specifically, the main objects of this chapter are introduced next. 

\begin{definition}[Field of operators]
\label{Chap3-Def-FieldOfLinOperators}
Let $\ts$ be a topological space. A family of linear operators $A_t:\hs_t\to\hs_t,\; t\in\ts,$ on Hilbert spaces $\hs_t,\; t\in\ts,$ is called a {\em field of operators} and $(\hs_t)_{t\in\ts}$ is called a {\em field of Hilbert spaces}. If, additionally, all operators $A_t\,,\; t\in\ts\,,$ are self-adjoint/normal/unitary/bounded/unbounded the field of operators is called {\em self-adjoint/normal/unitary/bounded/unbounded}.
\end{definition}

Note that the operators \gls{AT} can be defined on different Hilbert spaces. The aim is to analyze the (H\"older-)continuity properties of the map 
$$
\Sigma:\ts\to\ks(\CM)\,,\; t\mapsto\sigma(A_t)\,,
$$ 
where $\ks(\CM)$ is equipped with the Hausdorff metric defined by the Euclidean metric on $\CM$.

\medskip

This type of problem is sometimes tackled by using continuous fields of $C^\ast$-algebras, c.f. \cite{Ell82,Bel94,PaRi16,BeBeNi16}. This work provides a new proof without this machinery whenever $(A_t)_{t\in\ts}$ is self-adjoint or unitary by following the lines of \cite{BeBe16}. Moreover, the new concept of continuous boundaries for a family of closed subsets of $\RM$ is introduced. As it turns out, the continuity of the boundaries is equivalent to the Fell-continuity for a family of closed subsets or the continuity with respect to the Hausdorff metric for a family of compact subsets, c.f. Theorem~\ref{Chap3-Theo-BoundVietContCompVers} and Remark~\ref{Chap3-Rem-BoundVietContCompVers}. The case of normal operators is discussed later by using the general tool of $C^\ast$-algebras. This concept delivers not only sufficient conditions for the continuity of $\Sigma$ as one can find in the literature. We actually prove a characterization for the continuity of $\Sigma$. If $(\ts,d)$ is a complete metric space, these conditions are extended to show a new characterization for the H\"older-continuity of $\Sigma$ for self-adjoint and unitary fields of operators. In addition, the boundaries of the spectra associated with a self-adjoint field of operators behave H\"older-continuous. Interesting phenomena are observed at gaps that close, c.f. also Example~\ref{Chap3-Ex-AlmostMathieu} below.

\medskip

As mentioned before, one tool to prove the continuous behavior of the spectra is the concept of continuous fields of $C^\ast$-algebras. It was initially proposed by {\sc Kaplanski} \cite{Kap49,Kap51}, {\sc Fell} \cite{Fel60}, {\sc Tomiyama} and {\sc Takesaki} \cite{ToTa61,Tom62} and further developed by {\sc Dixmier}  and {\sc Douady} \cite{DiDo63}, see Section~\ref{Chap3-Sect-ContFieldCAlg} for more details. The continuous behavior of the spectra in this context was first shown by {\sc Kaplanski} \cite[Lemma 3.3]{Kap51} for continuous self-adjoint vector fields. It was rediscovered by several authors over the last decades, c.f. the discussion at Theorem~\ref{Chap3-Theo-ContFieldCALgContSpectr}. For in\-stance, this technique was used for almost periodic Schr\"odinger operators \cite{Ell82}, magnetic Schr\"odinger operators \cite{Bel94,PaRi16} and Schr\"odinger operators associated with quasicrystals \cite{BeBeNi16}. Here a new approach delivers the continuity of $\Sigma$ without the abstract theory of $C^\ast$-algebras based on the work \cite{BeBe16}. However, the theory of continuous fields of $C^\ast$-algebras provides a tool to prove the continuity of $\Sigma$ based on \cite{LaRa99} which is discussed in Chapter~\ref{Chap4-ToolContBehavSpectr}. Note that such a tool is currently not available for the H\"older-continuity of $\Sigma$, c.f. Section~\ref{Chap8-Sect-HolderCont}.

\medskip

It is well-known that continuity in the operator norm guarantees the continuity of $\Sigma$ whenever the operators are defined on the same Hilbert space. Unfortunately, the convergence in operator norm is too restrictive for many applications. On the other hand, the continuity in the strong operator topology can be verified in a wider class of applications while it does not imply that $\Sigma$ is continuous, see Section~\ref{Chap3-Sect-OpTop} for more details and references therein. The following example illustrates this difficulty.

\begin{figure}[htb]
\centering
\includegraphics[scale=0.535]{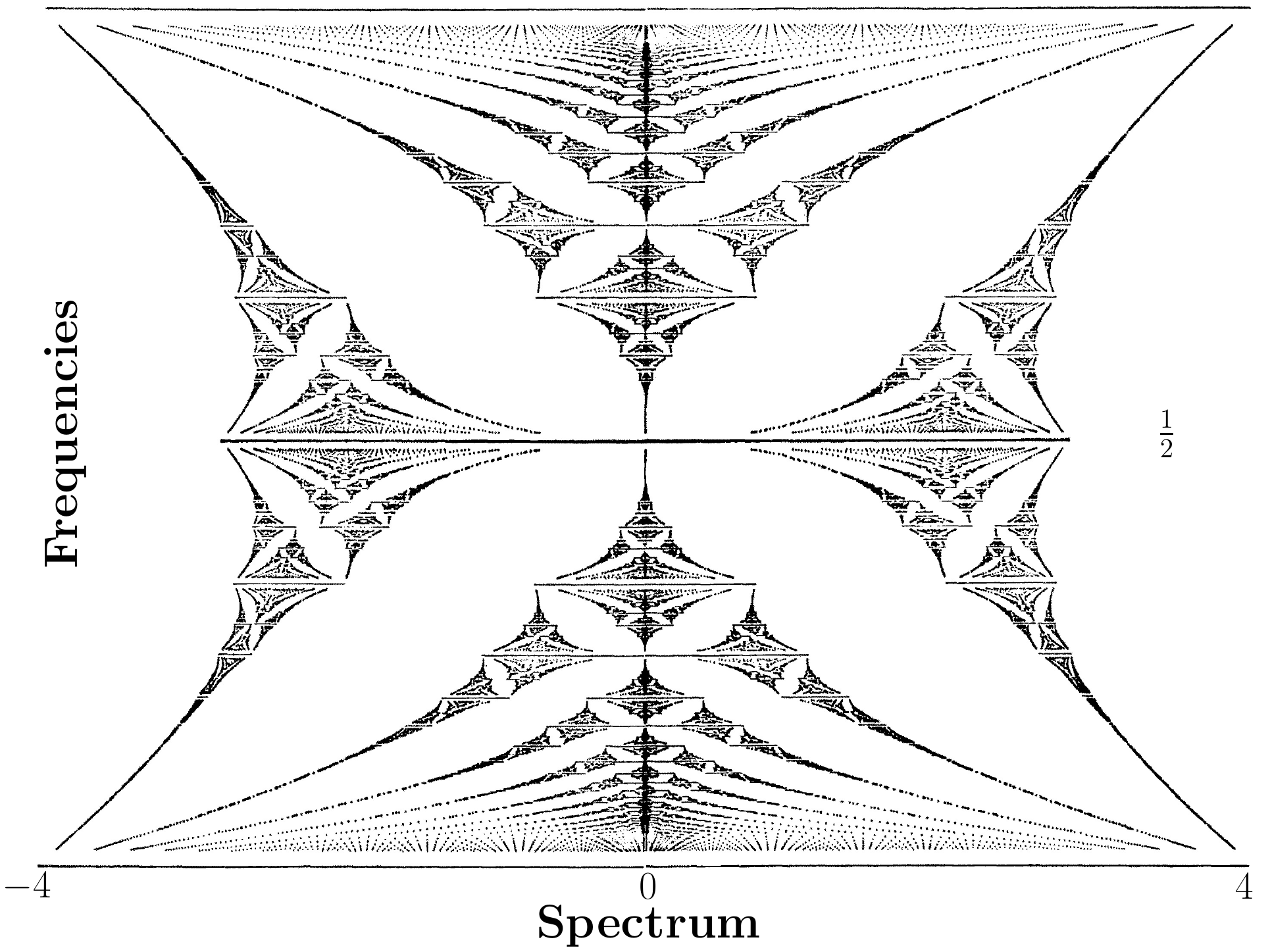}
\caption{(Hofstadter's butterfly \cite{Hof76}) The spectra of the Almost-Mathieu operator with coupling constant $\lambda=1$ for different frequencies $\alpha$.}
\label{Chap3-Fig-Hofstadter}
\end{figure}

\begin{example}
\label{Chap3-Ex-AlmostMathieu}
Let $H_\alpha$, where $\alpha \in [0,1]=:\ts$, be the Almost-Mathieu operator acting on $\hs:=\ell^2(\ZM)$ as follows
$$
H_{\alpha}\psi(n)=
 \psi(n+1)+ \psi(n-1)+ 2\lambda \cdot \cos\big(2\pi\cdot(n\alpha+ \theta)\big)\cdot \psi(n)\, ,
  \qquad
   \psi\in\ell^2(\ZM) \, ,
    	\;\;n\in\ZM \, .
$$
In this definition, $\theta$ is a fixed parameter and $\lambda>0$ is called the coupling constant. Whenever $\lambda=1$ it is also called the Harper model, c.f. \cite{Har55}. The family of operators $(H_\alpha)_{\alpha\in\ts}$ is strongly continuous in $\alpha$. On the other hand, if $\alpha\neq \beta$ are irrationally and rationally independent, it is known that $\|H_\alpha-H_\beta\|=2\lambda$. Thus, this family of operators is not norm continuous in $\alpha\in\ts$. However, it has been shown \cite{AvSi85,AvMoSi91,Bel94,JiKr02} that the boundaries of the spectra of $H_\alpha$ were (Lipschitz)-continuous as long as the gaps do not close, i.e., there do not exist gap tips, c.f. Definition~\ref{Chap3-Def-ClosedGap}. Near the points where gaps closes, the gap edges are only H\"older-continuous of exponent $1/2$, c.f. \cite{RaBe90,HeSj90}. This observation is crucial since also in the general case of a bounded self-adjoint field of operators the exponent of the H\"older-continuity decreases by a factor $1/2$ which relies on the possibility that spectral gaps can close. In addition, the spectral edge regularity of a more general class of magnetic Hamiltonians is analyzed in a recent work \cite{CoPu15}.

\medskip

\setlength\parindent{1em} 
\indent {\sc Hofstadter} \cite{Hof76} provided first numerical results on the spectra of these operators in 1976. He calculated the spectra for specific frequencies $\alpha$ and so first represented the fractal nature of the spectra, c.f. Figure~\ref{Chap3-Fig-Hofstadter}. The picture of the spectra is known by the name Hofstadter's Butterfly. Recently, two groups of physicists experimentally confirmed these numerical results independently, c.f. \cite{Pono13,DWM13}. The experi\-ments were made in graphene devices fabricated on hexagonal boron nitride substrates. The existence of gap tips is observable in Hofstadter's butterfly, for instance, for the frequency $\alpha=1/2$ at the origin. Furthermore, the figure indicates the square root behavior of the spectra in a neighborhood of $\alpha=1/2$.
\setlength\parindent{0em}
\end{example}

The main guideline that is followed here for the proof of the characterization of the (H\"older-)continuity of $\Sigma$ is described next.

\begin{description}
\item[(3.I)\label{(3.I)}] If $A$ is a bounded operator on the Hilbert space $\hs$, then the spectrum $\sigma(A)$ is contained in the compact ball $\overline{B_{\|A\|}(0)}\subseteq\CM$ of radius $\|A\|$ around the origin $0$. The unbounded case is similarly treated by passing to the resolvent.
\item[(3.II)\label{(3.II)}] Let $f:\overline{B_{\|A\|}(0)}\to\CM$ be a continuous function and $A$ be normal. Then $f(A)$ defines a normal operator by the functional calculus. Furthermore, the equality $\|f(A)\|=\sup_{\lambda\in\sigma(A)}|f(\lambda)|$ holds.
\item[(3.III)\label{(3.III)}] Due to \nameref{(3.II)}, the presence and absence of the spectrum $\sigma(A)$ can be controlled by specific norms of $A$. 
\begin{description}
\item[(3.III.1)\label{(3.III.1)}] The operator $A$ is self-adjoint. Thus, the spectrum $\sigma(A)$ is contained in $\RM$. Then the real-valued polynomials $p(z):=m^2-(z-c)^2\,,\; m,c\in\RM\,,$ of degree $2$ allow to control the presence or absence of the spectrum $\sigma(A)$ in a neighborhood of $c\in\RM$ by using \nameref{(3.II)}, c.f. Lemma~\ref{Chap3-Lem-PresenceSpectrum}.
\item[(3.III.2)\label{(3.III.2)}] The operator $A$ is unitary. Hence, the spectrum $\sigma(A)$ is contained in the sphere $\SM^1:=\{\lambda\in\CM\;|\; |\lambda|=1\}$. Then the polynomials $p(z):=1+ E\cdot z\,,\; E\in\SM^1\,,$ allow to control the presence or absence of the spectrum $\sigma(A)$ in a neighborhood of $\overline{E}\in\SM^1$ by using \nameref{(3.II)}, c.f. Lemma~\ref{Chap3-Lem-PresenceSpectrumUnitary}.
\item[(3.III.3)\label{(3.III.3)}] The operator $A$ is normal. Then Urysohn's Lemma provides the existence of a continuous function $f:\overline{B_{\|A\|}(0)}\to[0,1]$ with support in $K$ for each compact neighborhood $K$ of $c\in \overline{B_{\|A\|}(0)}$. Consequently, $f$ controls the presence or absence of the spectrum $\sigma(A)$ in the neighborhood $K$ of $c\in\RM$ by using \nameref{(3.II)}, c.f. Theorem~\ref{Chap3-Theo-ContFieldCALgContSpectr}. Since complex-valued polynomials on $\overline{B_{\|A\|}(0)}$ are dense in the space of continuous functions $f:\overline{B_{\|A\|}(0)}\to\CM$ with respect to the uniform norm $\|\cdot\|_\infty$, it suffices to control the norms of $\|p(A,A^\ast)\|$.
\end{description}
\end{description}

\medskip

The chapter is organized as follows. In the first section, a short review of typically studied operator topologies of $\Ll(\hs)$ is provided where $\Ll(\hs)$ is the set of bounded linear operators on the Hilbert space $\hs$. Of particular interests are the implications on the behavior of the spectra for the different topologies. Then Section~\ref{Chap3-Sect-VietorisContinuity} is devoted to study the continuity of closed sets with respect to the Hausdorff metric and its relation to the Vietoris and Fell-topology. Additionally, the new concept of continuous boundaries of closed subsets of $\RM$ is investigated. In Section~\ref{Chap3-Sect-CharContSpectSelfAdj}, the \pt-continuity of a self-adjoint field of operators $(A_t)_{t\in\ts}$ is introduced. It turns out that this is the correct notion so that the spectra behave continuously. The related \pt-topology on $\Ll(\hs)$ is studied for its connections with the operator norm and the strong operator topology in Section~\ref{Chap3-Sect-RelP2ContOthTop}. Then the unbounded case of self-adjoint operators is analyzed in Section~\ref{Chap3-Sect-CharContSpecUnbouSelfAdjOp}. A similar notion of continuity like in the self-adjoint case turns out to be equivalent to the continuous behavior of the spectra for a unitary field of operators, c.f. Section~\ref{Chap3-Sect-CharContSpectUnit}. Then Section~\ref{Chap3-Sect-ContFieldCAlg} extends the characterization of the continuous behavior of the spectra to normal operators. This uses the abstract technique of continuous fields of $C^\ast$-algebras and provides a tool to prove this continuity of the map $\Sigma:\ts\to\ks(\CM)$, c.f. Chapter~\ref{Chap4-ToolContBehavSpectr}. The ideas developed in Section~\ref{Chap3-Sect-CharContSpectSelfAdj} and Section~\ref{Chap3-Sect-CharContSpectUnit} apply also to functions of operators, c.f. Corollary~\ref{Chap3-Cor-CharContSpectFunctOfSelfAdj} and Corollary~\ref{Chap3-Cor-CharContSpectFunctOfUnitary}. Finally, the H\"older-continuity of $\Sigma$ is characterized for self-adjoint and unitary operators, c.f. Section~\ref{Chap3-Sect-CharHolContSpect}. Estimates of the rate of convergence of the spectra are proven as well as for the gaps of the spectra. The reader is also referred to Figure~\ref{Chap1-Fig-Strategy3} on Page~\pageref{Chap1-Fig-Strategy3} where the results and their connections are given schematically.

\medskip

It is assumed that the reader is familiar with the framework of topology. Some basics are recalled next that are used throughout. Let $X$ be a set and $\tau$ be a family of subsets of $X$. Then $\tau$ is called {\em a topology on $X$} if $X,\emptyset\in \tau$, every union of elements of $\tau$ belongs to $\tau$ and each intersection of finitely many elements of $\tau$ is an element of $\tau$. Then the tuple $(X,\tau)$ is called {\em a topological space}. The elements of $\tau$ are called {\em open} subsets of $X$. Furthermore, a subset $F\subseteq X$ is called {\em closed} if it is the complement of an open set. An {\em open neighborhood $U$} of a point $x\in X$ is an open set $U\in\tau$ such that $x\in U$. A sequence $(x_n)_{n\in\NM}$ in $X$ is called convergent to $x\in X$ if, for each open neighborhood $U$ of $x$, there exists an $n_0\in\NM$ such that $x_n\in U$ for all $n\geq n_0$. In Definition~\ref{App1-Def-BaseTopology}, the notion of a (neighborhood) base for a topological space is introduced. A topological space is called {\em first countable} if there exists, for every point, a countable neighborhood base. If there is a countable base for the whole topology $\tau$, the topological space $(X,\tau)$ is called {\em second countable}. Let $(X,\tau_X)$ and $(Y,\tau_Y)$ be topological spaces. A map $f:X\to Y$ is called {\em continuous} if the preimage $f^{-1}(U)$ is an element of $\tau_X$ for all $U\in\tau_Y$. The reader is referred to Appendix~\ref{App1-Topology} for a more detailed summary.

\medskip

In this chapter, the functional calculus for operators and elements of $C^\ast$-algebras is used. More precisely, let $\fz$ be a normal element of a unital $C^\ast$-algebra $\AG$. Then there is a unique morphism $\Phi$ of the $C^\ast$-algebra of continuous, complex-valued function on the spectrum $\sigma(\fz)$ into $\AG$ such that $\Phi(1)=\mathpzc{1}$ and $\Phi(id)=\fz$ where $id:\sigma(\fz)\to\sigma(\fz)\,,\; z\mapsto z$. Additionally, $\Phi$ is an isomorphism if the image is restricted to $\AG(\fz,\mathpzc{1})$ where $\AG(\fz,\mathpzc{1})$ denotes the sub-$C^\ast$-algebra of $\AG$ generated by the identity $\mathpzc{1}$ and $\fz$. Then, for every continuous function $\phi:\CM\to\CM$, the element $\Phi(\phi)\in\AG$ is denoted by $\phi(\fz)$. The identity $\sigma(\phi(\fz))=\phi(\sigma(\fz))$ follows by construction. In this work, polynomials of operators play an important role. Let $p(z,\overline{z})$ be a complex-valued polynomial in $z$ and $\overline{z}$ and $\phi:\CM\to\CM,\; \phi(z):=p(z,\overline{z})$. Then $\phi(\fz)$ is the usual algebraic combination $p(\fz,\fz^\ast)$. The reader is referred to \cite[Section 1.5]{Dixmier77} and \cite[Section~2.5]{Murphy90} for a more detailed discussion on the functional calculus.

\section{Operator topologies}
\label{Chap3-Sect-OpTop}

Let $\hs$ be a Hilbert space with inner product $\langle\cdot,\cdot\rangle$ and induced norm $\|\cdot\|_\hs$. For the sake of simplicity, the norm $\|\cdot\|_\hs$ is denoted by $\|\cdot\|$ whenever there is no confusion. Then $\Ll(\hs)$ denotes the $C^\ast$-algebra of bounded linear operators on $\hs$ equipped with the operator norm $\|A\|:=\sup_{\|\psi\|\leq 1}\|A\psi\|_\hs$ for $A\in\Ll(\hs)$. Different topologies are studied on $\Ll(\hs)$, see e.g. \cite{Wei97}. For instance:

\begin{description}
\item[(OT1)\label{(OT1)}] the operator norm topology induced by the operator norm $\|\cdot\|$;
\item[(OT2)\label{(OT2)}] the strong operator topology induced by the seminorms $\Ll(\hs)\ni A\mapsto\|A\psi\|_\hs\in[0,\infty)$ for all $\psi\in\hs$.
\end{description}

Let $A\in\Ll(\hs)$ and $(A_n)_{n\in\NM}\subseteq\Ll(\hs)$ be a sequence of self-adjoint operators such that $(A_n)_{n\in\NM}$ converges to $A$ in one of these two topologies. Then it is natural to ask for the spectral properties that are preserved by taking the limit. The operator norm topology safes many of the spectral properties which is well-known. For instance, the spectra $\sigma(A_n)\,,\; n\in\NM\,,$ converge to $\sigma(A)$ in the Hausdorff metric and the corresponding eigenprojections converge in norm, c.f. \cite{Wei97}. On the other hand, the convergence with respect to the operator norm is too restrictive for the most purposes, see e.g. Example~\ref{Chap3-Ex-AlmostMathieu}. This suggests to consider a weaker topology on $\Ll(\hs)$. Another possible choice for a topology is the strong operator convergence which can be veri\-fied in many circumstances. It induces a lower semi-continuity of the spectra, i.e., $\sigma(A)\subseteq\limsup_{n\to\infty}\sigma(A_n) = \bigcap_{n\in\NM}\overline{\big( \bigcup_{m=n}^\infty \sigma(A_n)\big)}$. However, it may happen in general that $\lim_{n\to\infty}\sigma(A_n)$ and $\sigma(A)$ do not agree (or even that $\big(\sigma(A_n)\big)_{n\in\NM}$ is not convergent in the Hausdorff metric). Additionally, strong operator convergence leads to the strong convergence of the spectral families. Compare \cite{Wei97} for a more detailed discussion on the spectral implications. The following example shows that strong convergence does not imply convergence of the associated spectra.

\begin{example}
\label{Chap3-Ex-StrongConvNotConvSpect}
Consider the Hilbert space $\hs:=\ell^2(\ZM)$ of square summable sequences $\psi:\ZM\to\CM$ with scalar product $\langle\psi,\varphi\rangle:=\sum_{j\in\ZM} \overline{\psi(j)}\cdot\varphi(j)$. Define the bounded operators $A_n,A\in\Ll(\hs),\; n\in\NM,$ by
$$
\begin{array}{ll}
	(A_n\psi)(j) \; 
		&:= \; \chi_{[-n,n]}(j) \cdot \psi(j) \, ,\\
	(A\psi)(j) \; 
		&:= \; \psi(j) \, ,
\end{array}
	\qquad\qquad \psi\in\ell^2(\ZM) \, 
		,\; j\in\ZM \, ,
$$
where $\chi_{[-n,n]}$ denotes the characteristic function on the interval $[-n,n]$. The operators $A_n,\; n\in\NM,$ are projections and $A$ is the identity operator on $\ell^2(\ZM)$. Thus, the equalities $\sigma(A_n)=\{0,1\},\; n\in\NM\,,$ and $\sigma(A)=\{1\}$ hold. It is not difficult to check that $(A_n)_{n\in\NM}$ converges to $A$ in the strong operator topology. On the other hand, the limit $\lim_{n\to\infty}\sigma(A_n)$ is equal to $\{0,1\}$ with respect to the Hausdorff metric on $\ks(\RM)$. 
Consequently, the spectra do not converge as $\sigma(A)\subsetneq\lim_{n\to\infty}\sigma(A_n)$.
\end{example}

Note that this semi-continuity of the spectra fails as soon as non self-adjoint operators are considered, c.f. \cite[Example~IV.3.8]{Kato95}.

\section{Continuity with respect to the Hausdorff metric}
\label{Chap3-Sect-VietorisContinuity}

In the following the relations between the Vietoris-topology, Fell-topology and the Hausdorff metric are discussed. Finally, the connection to the continuity of the boundaries associated with the closed subsets of $\RM$ is studied.

\medskip

Let $X$ be a topological space. The set of closed subsets of $X$ is denoted by $\cs(X)$. Recall the definition of the Vietoris-topology on \gls{csX} the closed subsets of $X$. For every $F\in\cs(X)$ and each finite collection $\os$ of open subsets of $X$, the sets
$$
\text{\gls{UFO}} \; := \;
	\{
		Y\in\cs(X)\;|\; 
			F\cap Y=\emptyset,\; O\cap Y\neq\emptyset \text{ for all } O\in\os 
	\} \, ,
$$
define a base for the Vietoris-topology on $\cs(X)$. A base for the Fell-topology is similarly defined while $F\subseteq X$ ranges over all compact subsets of $X$. A family $K_t\in\cs(X),\; t\in\ts,$ of closed sets, indexed by a topological space $\ts$, is called {\em Vietoris-continuous} if the map $\ts\ni t\mapsto K_t\in\cs(X)$ is continuous with respect to the Vietoris-topology on $\cs(X)$. Analogously, the {\em Fell-continuity} of a family $(K_t)_{t\in\ts}$ of closed sets is defined. The set of compact subsets of a topological spacer $X$ is denoted by \gls{ksX}.

\medskip

Let $(X,d)$ be a complete metric space. A topology on the space $\ks(X)\subseteq\cs(X)$ of compact subsets of $X$ can also be defined by the Hausdorff metric. In detail, the {\em Hausdorff metric $d_H:\ks(X)\times\ks(X)\to[0,\infty)$} is defined by
$$
\text{\gls{dH}}(F,K) \; := \;
	\max\left\{
		\sup\limits_{x\in F}\inf\limits_{y\in K} d(x,y),\;
		\sup\limits_{y\in K}\inf\limits_{x\in F} d(x,y)
	\right\} \, ,
	\qquad
	F,K\in\ks(X) \, .
$$

Then the space $(\ks(X),d_H)$ gets a complete metric space. 

\medskip

In the following, the connections between the Fell, Vietoris-topology and the Hausdorff metric are given. The equivalence of (i) and (ii) is well-known see \cite[Theorem II-6]{CastaingValadier77} which is even true for closed subsets $\cs(X)$. The equivalence to (iii) strongly depends on the fact that compact sets are considered and $(X,d)$ is assumed to be locally compact. Assertion (iii) states that the Vietoris-continuity of compact sets is a local Fell-continuity. The fact that the compact sets locally stays in a bounded region in the sense of (iii) is used in the proof of Theorem~\ref{Chap3-Theo-P2ContEquivContSpect} and Theorem~\ref{Chap3-Theo-CharContSpectUnitary}. Note that the topology induced by the Hausdorff metric might be different to the Vietoris-topology if the metric $(X,d)$ is not complete while $d$ induces the same topology on $X$.

\begin{theorem}
\label{Chap3-Theo-VietFellHausMetricEquiv}
Let $\ts$ be a topological space and $(X,d)$ be a locally compact, complete metric space. Consider a family $K_t\in\ks(X),\; t\in\ts,$ of compact subsets of $X$. Then the following assertions are equivalent.
\begin{itemize}
\item[(i)] The map $\ts\ni t\mapsto K_t\in\ks(X)$ is continuous with respect to the Hausdorff metric. More precisely, for all $t_0\in\ts$, the equation
$$
\lim\limits_{t\to t_0} d_H(K_t,K_{t_0}) \;
	= \; 0
$$
holds.
\item[(ii)] The family $(K_t)_{t\in\ts}$ is Vietoris-continuous.
\item[(iii)] For every $t_0\in\ts$, there exist a neighborhood $U_0\subseteq\ts$ of $t_0$ and a compact subset $K\subseteq X$ such that $K_t\subseteq K$ for $t\in U_0$ and $(K_t)_{t\in U_0}$ is Fell-continuous in the space $\ks(K)\subsetneq\ks(X)$.
\end{itemize}
\end{theorem}

\begin{proof}
(i)$\Leftrightarrow$(ii): This is proven in \cite[Theorem II-6]{CastaingValadier77}.

\vspace{.1cm}

(ii)$\Rightarrow$(iii): Consider a $t_0\in\ts$. Let $U_1,\ldots,U_m$ be open subsets that cover $K_{t_0}$ while the closure $\overline{U_j}$ is compact for each $1\leq j\leq m$. The existence of the desired $U_1,\ldots,U_m$ is guaranteed by the compactness of $K_{t_0}$ and since $X$ is locally compact. Altogether, $K:=\bigcup_{j=1}^m \overline{U_j}$ defines a compact neighborhood of $K_{t_0}$. Since $(K_t)_{t\in\ts}$ is Vietoris-continuous, there is an open neighborhood $U_0\subseteq\ts$ of $t_0$ such that $K_t\in \us(F,\{X\})$ where $F:=\overline{X\setminus K}$. Hence, $K_t$ is contained in $K$ for all $t\in U_0$. Obviously, $(K_t)_{t\in U_0}$ is Fell-continuous in $\ks(K)$.

\vspace{.1cm}

(iii)$\Rightarrow$(ii): Consider a $t_0\in\ts$. Let $U_0\subseteq\ts$ be open and $K\subseteq X$ be compact so that (iii) holds. Fix $F\in\cs(X)$ and $\os$ a finite family of open subsets of $X$ such that $K_{t_0}\in\us(F,\os)$. Since $K_{t_0}\subseteq K$, it follows $O\cap K\neq\emptyset$ for all $O\in\os$. Define $\tilde{\os}:=\{O\cap K\;|\; O\in\os\}$ being a finite family of non-empty, open subsets in $K$ and set $\tilde{F}:=F\cap K\in\ks(K)$. By construction, $K_{t_0}$ is an element of the open set $\us(\tilde{F},\tilde{\os})$ in the Fell-topology of $\ks(K)$. Thus, by (iii) there exists an open neighborhood $V_0\subseteq U_0$ of $t_0$ such that 
$$
K_t\in
	\us(\tilde{F},\tilde{\os}) \; 
		:= \; \left\{ K\in\ks(K) \;|\; K\cap\tilde{F}=\emptyset, \; K\cap O\neq\emptyset \text{ for all } O\in\tilde{\os}		
		\right\}
	,\qquad t\in V_0 \, .
$$
It is not difficult to check that then $K_t\in\us(F,\os)$ for all $t\in V_0$ proving that $(K_t)_{t\in\ts}$ is Vietoris-continuous.
\end{proof}

\medskip

Theorem~\ref{Chap3-Theo-VietFellHausMetricEquiv} draws a connection between the Vietoris-topology and the Fell-topol\-ogy. More precisely, the Fell-topology is coarser than the Vietoris-topology but in general these topologies do not agree. On the other hand, it is clear that in compact spaces $X$ the Vietoris-topology and the Fell-topology coincide which is used for the implication (iii) to (ii). The following example shows the difference of these topologies. Specifically, for a non-compact space $X$, a Fell-continuous family of sets $K_t\in\cs(X),\; t\in\ts,$ is not continuous with respect to the Hausdorff metric, in general. 

\begin{example}
\label{Chap3-Ex-VietFinerThanFell}
Let $X:=\RM$ be the Euclidean space and define the compact set $K:=[0,1]$ and the sequence of compact sets $K_n:=K\cup\{n+1\},\; n\in\NM$. Now let $F\subseteq\RM$ be compact and $\os$ be a finite family of open sets such that $K\in\us(F,\os)$. There is an $n_0\in\NM$ such that $n+1\not\in F$ for all $n\geq n_0$ by the compactness of $F$. Thus, $K_n$ is an element of $\us(F,\os)$ for all $n\geq n_0$. Consequently, $(K_n)_{n\in\NM}$ converges to $K$ in the Fell-topology. On the other hand, $d_H(K,K_n)=n$ tends to infinity if $n\to\infty$. Hence, $(K_n)_{n\in\NM}$ does not converge to $K$ in the Hausdorff distance and in the Vietoris-topology.
\end{example}

Let $X,Y$ be topological spaces. Then a function $\phi:X\to Y$ is called {\em closed} if the image of closed sets is closed. Such a map naturally induces a map $\hphi:\cs(X)\to\cs(Y)$ defined by $K\in\cs(X)\mapsto \phi(K)\in\cs(Y)$. It turns out that if, additionally, $\phi$ is continuous, then $\hphi$ is Vietoris-continuous. 

\begin{proposition}[\cite{Fil98,BeBe16}]
\label{Chap3-Prop-ConClosFuncHausCont}
Let $X,Y$ be locally compact spaces and consider a continuous, closed function $\phi:X\to Y$. Then, the map $\hphi:\cs(X)\to\cs(Y),\; K\mapsto \phi(K),$ is Vietoris-continuous.
\end{proposition}

\begin{proof}
Since $\phi$ is closed the function $\hphi$ is well-defined. Let $K_0\in\cs(X)$ and $\us(F,\os)$ be an open neighborhood of $\hphi(K_0)$ in $\cs(Y)$ where $F\subseteq Y$ is closed and $\os$ is a finite family of open subsets of $Y$. Define the set 
$$
W \; 
	:= \; \us\big(\phi^{-1}(F), \big\{\phi^{-1}(O)\;|\; O\in \os\big\}\big)  \, .
$$
Since $\phi$ is continuous, $\phi^{-1}(O)$ is open for every $O\in \os$ and $\phi^{-1}(F)$ is closed. Hence, $W$ is an open set in $\cs(X)$. Note that $\phi^{-1}(O)$ is non-empty for all $O\in\os$ since $\hphi(K_0)\in\us(F,\os)$, namely the image of $\phi$ intersects $O\in\os$. Pick a $K\in W$, i.e., the intersection $K\cap \phi^{-1}(F)$ is empty while $F\cap \phi^{-1}(O)$ is non-empty for all $O\in\os$. This implies that $\phi(K)\cap F=\emptyset$. In fact, otherwise, there would be an $x\in K$ with $\phi(x)\in F$ leading to $x\in \phi^{-1}(F)$ a contradiction. In addition, if $O\in\os$ and if $x\in K\cap \phi^{-1}(O)\neq \emptyset$, then $\phi(x)\in \phi(K)\cap O$ follows proving that $\phi(K)\cap O\neq \emptyset$ holds for all $O\in\os$. Altogether, $\hphi(K)$ is an element of $\us(F,\os)$. Hence, $W$ is contained in the inverse image $\hphi^{-1}(\us(F,\os))$. Finally, it is not difficult to check that $K_0\in W$. Thus, $\hphi$ is continuous. Note that the equation $W=\hphi^{-1}(\us(F,\os))$ can be actually verified. 
\end{proof}

\medskip

Next the concept of continuous boundaries is investigated for elements of $\cs(\RM)$. Let $K\in\cs(\RM)$. Then the boundary of $K$ is defined by $\partial K:=K\setminus K^\circ$ where $K^\circ$ denotes the interior of the set $K$. Whenever $(K_t)_{t\in\ts}$ is a family of compact sets, the Vietoris-continuity is characterized by the continuity of the boundaries defined in Definition~\ref{Chap3-Def-ContBoundaries} below. Of particular interest are gaps that close at a specific $t_0\in\ts$. They play an important role in Section~\ref{Chap3-Sect-CharHolContSpect}. 

\medskip

The supremum of a closed set $K$ can be equal to $\infty$ and, similarly, the infimum can be $-\infty$. Since the continuity of these elements are studied, the space $[-\infty,\infty]$ is considered. It is a compactification of $\RM$ and the topology is given by the Euclidean topology on $\RM$ and the sets $(C,\infty]$ and $[-\infty,C)$ for $C>0$ define a neighborhood base of $\infty$ and $-\infty$, respectively. Consider the supremum $\sup:\cs(\RM)\to[-\infty,\infty]\,,\; K\mapsto\sup(K):=\sup\{x\in K\}\,,$ and the infimum $\inf:\cs(\RM)\to[-\infty,\infty]\,,\; K\mapsto\inf(K):=\inf\{x\in K\}$. For $K\in\cs(\RM)$, $\sup(K)$ and $\inf(K)$ are elements of the boundary $\partial K$ if they are finite. These two elements are separately considered for the study of the continuity of the boundaries. 

\begin{proposition}
\emph{(\cite[Proposition~1]{BeBe16})}
\label{Chap3-Prop-SupInfVietCont}
The supremum $\sup:\cs(\RM)\to[-\infty,\infty]$ and the infimum $\inf:\cs(\RM)\to[-\infty,\infty]$ are Vietoris-continuous.
\end{proposition}

\begin{proof}
It is not difficult to check that the multiplication by $-1$ is Vietoris-continuous on $\cs(\RM)$. Consequently, it suffices to verify the Vietoris-continuity of the supremum as the identity $\inf(K)=-\sup(-K)$ holds for $K\in\cs(\RM)$.

\vspace{.1cm}

Let $K_0\in \cs(\RM)$ be such that $\lambda_0:=\sup(K_0)<\infty$ and $\varepsilon>0$. Since $K_0$ is closed, it follows that $\lambda_0\in K_0$. Define the closed set $F_\varepsilon:=[\lambda_0+\varepsilon,\infty)$ and the open set $O_\varepsilon=(\lambda_0-\varepsilon,\infty)$. Then $\us(F_\varepsilon, \{O_\varepsilon\})$ is an open neighborhood of $K_0$. Consider a closed set $K\in \us(F_\varepsilon, \{O_\varepsilon\})$. Since $K\cap F_\varepsilon=\emptyset$, it follows that $\sup(K)<\lambda_0 +\varepsilon$. Furthermore, the estimate $\lambda_0-\varepsilon <\sup(K)$ holds as $K\cap O_\varepsilon\neq \emptyset$. Thus, $K\in \us(F_\varepsilon, \{O_\varepsilon\})$ yields $|\sup(K)-\sup(K_0)|< \varepsilon$. Since $\varepsilon>0$ was arbitrary, the Vietoris-continuity of the supremum at $K_0\in\cs(\RM)$ follows.

\vspace{.1cm}

Let $K_0\in \cs(\RM)$ be such that $\sup(K_0)=\infty$. For each constant $C>0$, define the open set $O_C:=(C,\infty)$. Then $\us(\emptyset,\{O_C\})$ is an open neighborhood of $K_0$. By construction, every $K\in\us(\emptyset,\{O_C\})$ satisfies $O_C\cap K\neq\emptyset$ implying $\sup(K)>C$. Hence, the supremum is Vietoris-continuous at $K_0\in\cs(\RM)$.
\end{proof}

\medskip

Following the lines of \cite{BeBe16}, the notions of a gap, gap tip and continuous boundaries are discussed.

\begin{definition}[Gap]
\label{Chap3-Def-Gap}
Consider a closed set $K\in\cs(\RM)$. Then an open interval $(a,b)\subseteq\RM$ where $-\infty< a< b <\infty$ is called a {\em gap} if $a,b\in K$ and $(a,b)\cap K$ is empty. An $x\in\RM$ is called a {\em gap edge of $K$} whenever there is a gap $(a,b)$ of $K$ such that either $x=a$ or $x=b$.
\end{definition}

It is worth noticing that the open intervals $(-\infty,\inf(K))$ and $(\sup(K),\infty)$ are not gaps in terms of Definition~\ref{Chap3-Def-Gap}. It follows from definition that the boundary of a set $K$ is the set of gap edges together with $\sup(K)$ and $\inf(K)$ if $\sup(K)$ and $\inf(K)$ are finite.

\medskip

For a family of closed sets $(K_t)_{t\in\ts}$ indexed by a topological space $\ts$, it may happen that gaps close, c.f. Example~\ref{Chap3-Ex-GapClos}.

\begin{figure}[htb]
\centering
\includegraphics[scale=1.3]{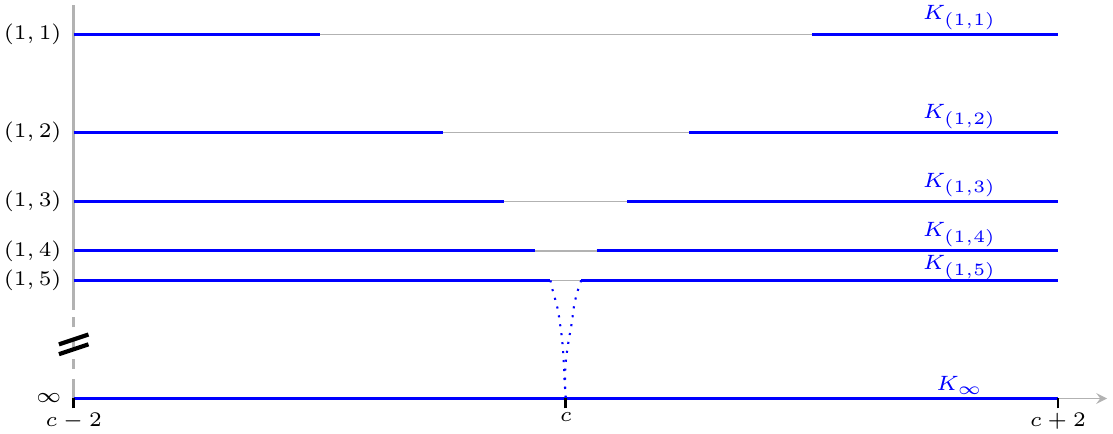}
\caption{Gaps that close at $c\in K_\infty$. The blue lines are the subsets $K_{(1,m)}\,,\; m\in\NM\,,$ and $K_\infty$ of $\RM$.}
\label{Chap3-Fig-ClosGap}
\end{figure}

\begin{example}
\label{Chap3-Ex-GapClos}
Let $\ts$ be the one-point compactification of $\NM\times\NM$, i.e., $\ts:=\NM\times\NM\cup\{\infty\}$ where a neighborhood basis of $\infty$ is given by $\{(n,m)\;|\; n,m\geq N\}$ for $N\in\NM$. Let $c\in\RM$ and, for $t\in\ts$, define the compact subset of $\RM$ by
$$
K_t\; := \; 
	\begin{cases}
		[c-2,c+2], \qquad &t=\infty \text{ or } t=(n,m) \text{ with } n\neq 1 \, ,\\
		[c-2,c-\frac{1}{m}]\cup [c+\frac{1}{m},c+2], \qquad &t=(1,m) \, .
	\end{cases}
$$
Due to construction, $d_H(K_{(n,m)},K_\infty)$ is equal to $0$ if $n\neq 1$ and otherwise it is equal to $2/m$. Hence, the family of compact sets $(K_t)_{t\in\ts}$ is continuous with respect to the Hausdorff metric $d_H$. Then the gaps $(c-\frac{1}{m},c+\frac{1}{m})$ of $K_{(1,m)}$ close at $c\in K_\infty$ if $m$ tends to infinity. This is illustrated in Figure~\ref{Chap3-Fig-ClosGap}. Note that in each neighborhood $U$ of $\infty\in\ts$ there are $t\in U$ with gaps $(a_t,b_t)$ of $K_t$ such that these gaps are close to $c$. The point $c$ is called a gap tip, c.f. Definition~\ref{Chap3-Def-ClosedGap} below. On the other hand, there does not exist a neighborhood $V$ of $\infty\in\ts$ such that every $t\in V$ has a gap close to $c$ as $K_{(n,m)}=K_\infty$ for $n\neq 1$.
\end{example}

In view of Example~\ref{Chap3-Ex-GapClos}, the notion of gap tip is defined in the following for a family of closed subsets of $\RM$. These gap tips play a crucial role for the H\"older-continuous behavior of the spectra as discussed in Section~\ref{Chap3-Sect-CharHolContSpect}. 

\begin{figure}[htb]
\centering
\includegraphics[scale=1.6]{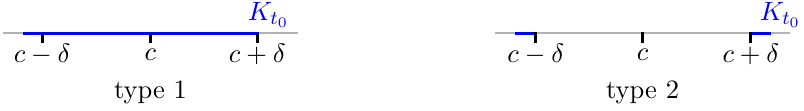}
\caption{Illustration of a gap tip $c\in K_{t_0}$ of type $1$ and $2$. The blue lines represent the set $K_{t_0}$.}
\label{Chap3-Fig-IsolGapType}
\end{figure}

\begin{definition}[Gap tip, \cite{BeBe16}]
\label{Chap3-Def-ClosedGap}
Let $\ts$ be a topological space and $K_t\in\cs(\RM),\; t\in\ts,$ be a family of closed sets. For $t_0\in\ts$, a value $c\in K_{t_0}$ is called a {\em gap tip} if there is a net $(t_\iota)_{\iota\in I}$ converging to $t_0$ and there are maps $a:I\to\RM, \; b:I\to\RM$ such that $(a_\iota,b_\iota)$ is a gap in $F_{t_\iota}$ for $\iota\in I$ and 
$$
\lim_{\iota} a_\iota 
	\; = \; c
	\; = \; \lim_{\iota} b_\iota \, .
$$
A gap tip $c\in F_{t_0}$ is called 
\begin{itemize}
\item {\em isolated of type $1$} if there is a $\delta>0$ such that $(c-\delta,c+\delta)\subseteq F_{t_0}$,
\item {\em isolated of type $2$} if there is a $\delta>0$ such that $(c-\delta,c+\delta)\cap F_{t_0}=\{c\}$,
\end{itemize}
c.f. Figure~\ref{Chap3-Fig-IsolGapType}.
\end{definition}

\begin{remark}
\label{Chap3-Rem-Net}
Since there is no requirement on the topology of $\ts$, it is necessary to use nets instead of sequences. Clearly, a net can be replaced by a sequence if $\ts$ is first countable.
\end{remark}

With this notion at hand, the continuity of the boundaries is defined which turns out to be equivalent to the Vietoris-continuity of compact sets. This concept is a bit delicate due to the possibility of closing gaps defined in \nameref{(CB3)}.

\begin{definition}[Continuous boundaries, \cite{BeBe16}]
\label{Chap3-Def-ContBoundaries}
Let $\ts$ be a topological space and $K_t\in\cs(\RM)\,,\; t\in\ts\,,$ be a family of closed subsets of $\RM$. The boundaries of $(K_t)_{t\in\ts}$ are called {\em continuous at $t_0\in\ts$} if the following assertions hold.
\begin{description}
\item[(CB1)\label{(CB1)}] The maps $\ts\ni t\mapsto\sup(K_t)\in\RM\cup\{-\infty,\infty\}$ and $\ts\ni t\mapsto\inf(K_t)\in\RM\cup\{-\infty,\infty\}$ are continuous.
\item[(CB2)\label{(CB2)}] Let $(a,b)$ be a gap of $K_{t_0}$ and $\varepsilon>0$. Then there exists an open neighborhood $U:=U(\varepsilon,a,b)$ of $t_0$ such that, for each $t\in U$, there is a gap $(a_t,b_t)$ of $K_t$ satisfying $|a-a_t|<\varepsilon$ and $|b-b_t|<\varepsilon$.
\item[(CB3)\label{(CB3)}] Let $(t_\iota)_{\iota\in I}$ be a net converging to $t_0$ and $(a_\iota,b_\iota)$ be a gap of $K_{t_\iota}$ such that $\lim a_\iota = a$ and $\lim b_\iota = b$. Then $(a,b)$ is a gap of $K_{t_0}$ that is closed whenever $a=b$.
\end{description}
The boundaries of $(K_t)_{t\in\ts}$ are called {\em continuous} if the boundaries are continuous at each $t_0\in\ts$.
\end{definition}

\begin{remark}
\label{Chap3-Rem-ContBoundariesSeq}
The nets in Definition~\ref{Chap3-Def-ContBoundaries} can be replaced by sequences if $\ts$ is first countable, c.f. Remark~\ref{Chap3-Rem-Net}.
\end{remark}

Recall that $B_\varepsilon(z)\subseteq\RM$ denotes the open interval $(z-\varepsilon,z+\varepsilon)$ for $\varepsilon>0$ and $z\in\RM$.

\begin{lemma}
\emph{(\cite[Lemma~5]{BeBe16})}
\label{Chap3-Lem-VietContImplGapCont}
Let $\ts$ be a topological space and $K_t\in\cs(\RM),\; t\in\ts,$ be a family of closed subsets of $\RM$. If $(K_t)_{t\in\ts}$ is Vietoris-continuous then the boundaries of $(K_t)_{t\in\ts}$ are continuous.
\end{lemma}

\begin{proof}
Let $t_0\in T$. Then Condition \nameref{(CB1)} follows from Proposition~\ref{Chap3-Prop-SupInfVietCont}.

\vspace{.1cm}

\nameref{(CB2)}: Let $(a,b)$ be a gap of $K_{t_0}$. Consider a $(b-a)/2>\varepsilon>0$ and define the open balls $O_a:=B_\varepsilon(a),\; O_b:=B_\varepsilon(b)$ and the closed set $F:=[a+\varepsilon,b-\varepsilon]$. Since $(b-a)/2>\varepsilon$, the inequality $a+\varepsilon<b-\varepsilon$ holds and so $F$ is well-defined. Clearly $\us(F,\{O_a,O_b\})$ is a neighborhood of $K_{t_0}$. Thus, there exists an open neighborhood $U$ of $t_0$ such that $K_t\in\us(F,\{O_a,O_b\})$ for all $t\in U$ using the Vietoris-continuity of $(K_t)_{t\in T}$. Then for $t\in U$ the maximum $a_t:=\max\big(O_a\cap K_t\big)$ and the minimum $b_t:=\min\big(O_b\cap K_t\big)$ exist as the intersections are non-empty and, additionally, $F\cap K_t=\emptyset$. Since $a_t\in O_a$ and $b_t\in O_b$, the inequalities $|a-a_t|<\varepsilon$ and $|b-b_t|<\varepsilon$ hold. Furthermore, $(a_t,b_t)$ is a gap of $K_t$ since $F\cap K_t=\emptyset$ for $t\in U$.

\vspace{.1cm}

\nameref{(CB3)}: Let $(t_\iota)_{\iota\in I}$ be a net converging to $t_0$ and $(a_\iota,b_\iota)$ be a gap of $K_{t_\iota}$ for $\iota\in I$ such that $\lim a_\iota = a$ and $\lim b_\iota = b$. Then $a\in K_{t_0}$ holds. For indeed, otherwise, there is an $r>0$ such that $\us([a-r,a+r],\{\RM\})$ is a neighborhood of $K_{t_0}$. This contradicts $K_{t_\iota}\to K_{t_0}$ and $\lim a_\iota=a$ as the second implies that there is an $\iota^{(0)}\in  I$ such that $K_{t_\iota}\not\in \us([a-r,a+r],\{\RM\})$ for $\iota\geq\iota^{(0)}$. Similarly, $b$ is an element of $K_{t_0}$. Thus, whenever $a=b$ it follows that $(a,b)$ is a gap tip of $K_{t_0}$.

\vspace{.1cm}

Only the case $a<b$ is left, namely it need to be verified that $(a,b)$ is a gap of $K_{t_0}$. Assume the contrary, i.e., $(a,b)\cap K_{t_0}\neq\emptyset$. Pick $x\in(a,b)\cap K_{t_0}$ and choose $\varepsilon>0$ such that $O_\varepsilon:=B_\varepsilon(x)\subset(a+\varepsilon,b-\varepsilon)$. Since $K_{t_0}\in\us(\emptyset,\{O_\varepsilon\})$, the Vietoris-continuity implies the existence of an $\iota^{(0)}\in I$ such that $K_{t_\iota}\in\us(\emptyset,\{O_\varepsilon\})$ for $\iota\geq\iota^{(0)}$. On the other hand, $\lim a_\iota = a$ and $\lim b_\iota = b$ lead to the existence of an $\iota^{(1)}$ such that $a_{\iota}\in B_\varepsilon(a)$ and $b_{\iota}\in B_\varepsilon(b)$ for $\iota\geq \iota^{(1)}$. Since $O_\varepsilon\subset(a+\varepsilon,b-\varepsilon)$ and $(a_{t_\iota},b_{t_\iota})\cap (a+\varepsilon,b-\varepsilon)=\emptyset$ for $\iota\geq \iota^{(1)}$, the intersection $O_\varepsilon\cap(a_{t_\iota},b_{t_\iota})$ is empty. Then by choosing an $\iota\geq\iota^{(0)},\iota^{(1)}$ this is a contradiction as $K_{t_\iota}\not\in \us(\emptyset,\{O_\varepsilon\})$ by $\iota\geq\iota^{(1)}$ while $K_{t_\iota}\in \us(\emptyset,\{O_\varepsilon\})$ by $\iota\geq\iota^{(0)}$.
\end{proof}

\medskip

Whenever $(K_t)_{t\in\ts}$ is a family of compact sets the converse of Lemma~\ref{Chap3-Lem-VietContImplGapCont} is valid.

\begin{lemma}
\emph{(\cite[Lemma~6]{BeBe16})}
\label{Chap3-Lem-GapContImplVietCont}
Let $\ts$ be a topological space and $K_t\in\ks(\RM),\; t\in\ts,$ be a family of compact subsets. Then $(K_t)_{t\in\ts}$ is Vietoris-continuous if the boundaries of $(K_t)_{t\in\ts}$ are continuous.
\end{lemma}

\begin{proof}
Consider $t_0\in T$. By compactness of $K_{t_0}$ the infimum and supremum of $K_{t_0}$ are finite. Define $m:=1+\max\{|\sup(K_{t_0})|, |\inf(K_{t_0})| \}$. It suffices to show that, for every net $(t_\iota)_{\iota\in I}$ tending to $t_0$, $K_{t_\iota}$ tends to $K_{t_0}$ in the Vietoris-topology. Let $F\subseteq\RM$ be closed and $\os=\{O_1,\ldots,O_n\}$ be a finite family of open subsets of $\RM$ such that $K_{t_0}\in \us(F,\os)$. The existence of an $\iota_{(F,\os)}\in I$ has to be shown such that $K_{t_\iota}\in\us(F,\os)$ for all $\iota\geq\iota_{(F,\os)}$. 

\vspace{.1cm}

Let $(t_\iota)_{\iota\in I}$ be a net that converges to $t_0$. According to \nameref{(CB1)}, there is an $\iota^{(0)}\in I$ such that $\inf(K_{t_\iota}),\sup(K_{t_\iota})\in[-m,m]$ for $\iota\geq\iota^{(0)}$. Thus, $K_{t_\iota}\subseteq	[-m,m]$ follows for $\iota\geq \iota^{(0)}$.

\vspace{.1cm}

Define $F_m:=F\cap[-m,m]$. Then, for $x\in F_m$, there is an $r(x)$ such that $B_{r(x)}(x)\cap K_{t_0}=\emptyset$. By compactness of $F_m$, there are $x_1,\ldots,x_l\in F_m$ with $r_k:=r(x_k)$ such that 
$$
F_m\subseteq\bigcup_{k=1}^l B_{r_k/2}(x_k) \, , 
	\qquad B_{r_k}(x_k)\cap K_{t_0}=\emptyset \, ,
		\qquad 1\leq k\leq l \, .
$$ 
Consider a $1\leq k\leq l$. Then there exists an $\iota^{(k)}\in I$ such that $\iota^{(k)}\geq\iota^{(0)}$ and for all $\iota\geq\iota^{(k)}$, the intersection $K_{t_\iota}\cap(x_k-r_k/2,x_k+r_k/2)$ is empty. For indeed, if $x_k\leq\inf(K_{t_0})$ or $x_k\geq\sup(K_{t_0})$ this follows by \nameref{(CB1)}. Whenever $\inf(K_{t_0}) < x_k < \sup(K_{t_0})$, there is a gap $(a,b)$ in $K_{t_0}$ with $a<x_k-r_k<x_k+r_k<b$. Applying \nameref{(CB2)} there exists an open neighborhood $U$ of $t_0$ such that there is a gap $(a_t,b_t)$ of $K_t$ satisfying $|a-a_t|<r_k/2$ and $|b-b_t|<r_k/2$ for all $t\in U$. Since $t_\iota\to t_0$, there is an $\iota^{(k)}\in I$ satisfying $\iota^{(k)}\geq\iota^{(0)}$ and $t_\iota\in U$ for $\iota\geq\iota^{(k)}$. Hence, the inclusion $B_{r_k/2}(x_k)\subseteq(a_{t_\iota},b_{t_\iota})$ holds for $\iota\geq\iota^{(k)}$ and so $B_{r_k/2}(x_k)\cap K_{t_\iota}=\emptyset$ follows for $\iota\geq\iota^{(k)}$. Then there is an $\iota_F\in I$ such that $\iota_F\geq \iota^{(k)}$ for all $0\leq k\leq n$. For $\iota\geq\iota_F$, the intersection $F_m\cap K_{t_\iota}$ is empty implying $F\cap K_{t_\iota}=\emptyset$ as $\iota_F\geq \iota^{(0)}$. 

\vspace{.1cm}

Now, let $O_j\in\os$. By assumption, there is an $x\in O_j\cap K_{t_0}$ and a radius $r>0$ such that $(x-r,x+r)\subseteq O_j$. Then there exists a $\tilde{\iota}^{(j)}\in I$ satisfying that $K_{t_\iota}\cap O_j\neq\emptyset$ for $\iota\geq\tilde{\iota}^{(j)}$. In fact, if $x$ is either $\inf(K_{t_0})$ or $\sup(K_{t_0})$ this immediately follows by \nameref{(CB1)}. Otherwise, assume the contrary, i.e., there is a subnet $(t_{\iota_k})$ tending to $t_0$ with $K_{t_{\iota_k}}\cap O_j=\emptyset$ for all $k$. Without loss of generality $r$ is chosen so that 
$$
\inf(K_{t_0})+r \;
	< \; x-r
		< \; x+r
			< \; \sup(K_{t_0})-r \, .
$$
By \nameref{(CB1)}, there is an $\iota'_k\geq\iota^{(0)}$ such that the intersections $K_{t_{\iota_k}}\cap[-m,x-r]$ and $K_{t_{\iota_k}}\cap[x+r,m]$ are non-empty for $\iota_k\geq\iota'_k$. Define 
$$
a_{\iota_k}:= \max\big(K_{t_{\iota_k}}\cap[-m,x-r]\big)\, ,
	\qquad 
		b_{\iota_k}:= \min\big(K_{t_{\iota_k}}\cap[x+r,m]\big) \, ,
			\qquad \iota_k\geq\iota'_k \, .
$$
Due to construction $a_{\iota_k},\, b_{\iota_k}$ are elements of $K_{t_{\iota_k}}$ and $(a_{\iota_k},b_{\iota_k})\cap K_{t_{\iota_k}}=\emptyset$. Since $(a_{\iota_k})_k$ and $(b_{\iota_k})_k$ are bounded in $\RM$, there is no loss of generality in supposing that the limits $\lim a_{\iota_k}=:a$ and $\lim b_{\iota_k}=:b$ exist otherwise pass to a subnet. Then $a\leq x-r$ and $b\geq x+r$ hold implying $(x-r,x+r)\subset(a,b)$. Since $(a,b)$ is a gap in $K_{t_0}$ by \nameref{(CB3)}, this leads to a contradiction as $x\in K_{t_0}$ and $x$ is an element of the gap $(a,b)$ at the same time. Using that $\os$ is a finite family of open sets, there is an $\iota_\os\in I$ such that $\iota_\os\geq\tilde{\iota}^{(j)}$ for all $1\leq j\leq n$. Then, for $\iota\geq\iota_\os$, the intersection $K_{t_\iota}\cap O_j$ is non-empty for all $1\leq j\leq n$.

\vspace{.1cm}

Choosing $\iota_{(F,\os)} \geq \iota_F,\iota_\os$ it follows by the previous considerations that $K_{t_\iota}\in\us(F,\os)$ for all $\iota\geq\iota_{(F,\os)}$.
\end{proof}

\begin{remark}
\label{Chap3-Rem-BoundContVietCont}
Lemma~\ref{Chap3-Lem-VietContImplGapCont} holds for all closed sets whereas Lemma~\ref{Chap3-Lem-GapContImplVietCont} is only valid for compact sets, c.f. Example~\ref{Chap3-Ex-ContBoundNotImplVietContClosSet}. 
\end{remark}

In the following, an example of a family of closed sets $(F_t)_{t\in\ts}$ is provided such that the boundaries of $(F_t)_{t\in\ts}$ are continuous while the family $(F_t)_{t\in\ts}$ is not Vietoris-continuous.

\begin{example}
\label{Chap3-Ex-ContBoundNotImplVietContClosSet}
Let $\ts:=\text{\gls{oNM}}:=\NM\cup\{\infty\}$ be the one-point compactification of $\NM$. Define iteratively the closed sets $F_1:=\RM$ and $F_n:=F_{n-1}\setminus (n,n+1/2)$ for $n\geq 2$ and set $F_\infty:=\bigcap_{n\in\NM} F_n$. The boundaries of the family $(F_t)_{t\in\ts}$ turn out to be continuous: Condition~\nameref{(CB1)} holds as $\sup(F_t)=\infty$ and $\inf(F_t)=-\infty$ for all $t\in\ts$. If $(a,b)$ is a gap of $F_\infty$, then there is an $n_0\in\NM$ such that $a=n_0$ and $b=n_0+1/2$. Thus, $(a_n,b_n):=(n_0,n_0+1/2)$ is a gap of $F_n$ for all $n\geq n_0$ and $\lim_{n\to\infty} a_n=a$ as well as $\lim_{n\to\infty} b_n=b$. Hence, \nameref{(CB2)} is satisfied. Finally, let $(a_n,b_n)$ be a gap of $F_n$ for $n\in\NM$ such that the limits $a:=\lim_{n\to\infty} a_n$ and $b:=\lim_{n\to\infty} b_n$ exist. By construction of $(F_t)_{t\in\ts}$, the sequences are eventually constant, i.e., there are $n_0,n_1\in\NM$ such that $a_n=n_0$ and $b_n=n_0+1/2$ for all $n\geq n_1$. Clearly, $(a,b)$ is a gap of $F_\infty$ implying that \nameref{(CB3)} holds. Consequently, the boundaries of $(F_t)_{t\in\ts}$ are continuous whereas $d_H(F_n,F_\infty)=1/2$ does not tend to zero if $n\to\infty$. Thus, the family of closed sets $(F_t)_{t\in\ts}$ is not Vietoris-continuous by Theorem~\ref{Chap3-Theo-VietFellHausMetricEquiv}.
\end{example}

The last two lemmata lead to the following characterization.

\begin{theorem}[\cite{BeBe16}]
\label{Chap3-Theo-BoundVietContCompVers}
Let $\ts$ be a topological space and $K_t\in\ks(\RM),\; t\in\ts,$ be a family of compact subsets. Then the following assertions are equivalent.
\begin{itemize}
\item[(i)] The family $(K_t)_{t\in\ts}$ is Vietoris-continuous.
\item[(ii)] The boundaries of $(K_t)_{t\in\ts}$ are continuous.
\end{itemize}
\end{theorem}

\begin{proof}
This is proven in Lemma~\ref{Chap3-Lem-VietContImplGapCont} and Lemma~\ref{Chap3-Lem-GapContImplVietCont}.
\end{proof}

\medskip

\begin{remark}
\label{Chap3-Rem-BoundVietContCompVers}
Following the lines of Lemma~\ref{Chap3-Lem-VietContImplGapCont} and Lemma~\ref{Chap3-Lem-GapContImplVietCont}, it is not difficult to check that the characterization of Theorem~\ref{Chap3-Theo-BoundVietContCompVers} holds also if $K_t\,,\; t\in\ts\,,$ are only closed subsets of $\RM$ while the Vietoris-topology is replaced by the Fell-topology in (i).
\end{remark}

The continuity of the boundaries of a family $(K_t)_{t\in\ts}$ of closed subsets of a suitable complete metric space $(X,d)$ can be also characterized by the continuity of certain functions see Proposition~\ref{Chap3-Prop-BoundedCompactXContIffPsi_xCont} below. This characterization is used in Section~\ref{Chap3-Sect-CharContSpectSelfAdj} and Section~\ref{Chap3-Sect-CharContSpecUnbouSelfAdjOp}.

\medskip

Let $X$ be a locally compact space. Recall the notion of the Fell-topology on the set $\cs(X)$ of closed subsets of $X$. For every compact $F\subseteq X$ and each finite family $\os$ of open subsets of $X$, the sets
$$
\us(F,\os) \; := \;
	\big\{
		Y\in\cs(X)\;\big|\; 
			F\cap Y=\emptyset,\; O\cap Y\neq\emptyset \text{ for all } O\in\os 
	\big\} 
$$
define a base for the Fell-topology. Clearly, the only difference between the Fell-topology and Vietoris-topology is that $F$ is forced to be compact in the Fell-topology and only closed in the Vietoris-topology. As discussed in Theorem~\ref{Chap3-Theo-VietFellHausMetricEquiv}, these two topologies are related if $X$ is a locally compact, complete metric space. On the other hand, the Vietoris-topology is finer than the Fell-topology, c.f. Example~\ref{Chap3-Ex-VietFinerThanFell}. The following two lemmata are analogs for the Fell-topology of Proposition~\ref{Chap3-Prop-SupInfVietCont} and Proposition~\ref{Chap3-Prop-ConClosFuncHausCont}.

\begin{lemma}
\label{Chap3-Lem-InfFellCont}
Consider the space $X:=[C,\infty)\subseteq\RM$ with induced Euclidean metric and $C\in\RM$. Then the infimum $\inf:\cs(X)\to\RM$ is Fell-continuous.
\end{lemma}

\begin{proof}
Let $K_0\in\cs(X)$ and $\varepsilon>0$. Then $\lambda_0:=\inf(K_0)$ is greater than or equal to $C$. Since $K_0$ is closed it follows that $\lambda_0\in K_0$. If $\lambda_0-C>\varepsilon$, define $F_\varepsilon:=[C,\lambda_0-\varepsilon]$ and the open set $O_\varepsilon:=(C,\lambda_0+\varepsilon)$. Then $\us(F_\varepsilon,\{O_\varepsilon\})$ defines an open neighborhood of $K_0$ in the Fell-topology. Consider a $K\in\us(F_\varepsilon,\{O_\varepsilon\})$. Since $K\cap F_\varepsilon=\emptyset$, it follows that $\inf(K)> \lambda_0-\varepsilon$. Furthermore, $K\cap O_\varepsilon\neq\emptyset$ leads to $\inf(K)< \lambda_0+\varepsilon$. Thus, $|\inf(K)-\lambda_0|<\varepsilon$ holds for all $K\in\us(F_\varepsilon,\{O_\varepsilon\})$. Whenever $\lambda_0-C\leq\varepsilon$ the same estimate holds for $K\in\us(\emptyset,\{O_\varepsilon\})$. Hence, the infimum is Fell-continuous.
\end{proof}

\medskip

Clearly, the supremum is also Fell continuous for the space $X:=(-\infty,C]$ with $C\in\RM$. 

\medskip

A function between two topological spaces is called {\em proper} if the preimage of compact sets are compact.

\begin{lemma}
\label{Chap3-Lem-FunctFellCont}
Let $X,Y$ be locally compact spaces and $\phi:X\to Y$ be a continuous, proper function. Then the map $\hat{\phi}:\cs(X)\to\cs(Y)\,,\; K\mapsto\phi(K)\, ,$ is Fell-continuous.
\end{lemma}

\begin{proof}
The proof is left to the reader since it is similar to the proof of Proposition~\ref{Chap3-Prop-ConClosFuncHausCont} with the following two additional remarks. The closed set Lemma \cite{Pal70} asserts that any continuous, proper function between locally compact spaces is also closed. Consequently, the map $\hat{\phi}$ is well-defined. Furthermore, since $\phi$ is proper the preimage of compact sets is compact. Thus, the set
$$
W\; 
	:= \; \us\big(
		\phi^{-1}(F),
		\big\{
			\phi^{-1}(O)\;|\; O\in\os
		\big\}
	\big)
$$
is an open set in the Fell-topology of $X$ for $F\subseteq Y$ compact and $\os$ a finite family of open subsets of $Y$.
\end{proof}

\medskip

With the previous considerations at hand, the following proposition delivers a characterization of the Fell-continuity for a family $(K_t)_{t\in\ts}$ of closed subsets of $\RM$ in terms of certain functions. This is used to show a characterization of the continuous behavior of the spectra in the self-adjoint case, c.f. Theorem~\ref{Chap3-Theo-P2ContEquivContSpect} and Theorem~\ref{Chap3-Theo-CharContUnbSelfAdjOpSpecGap} below.

\medskip

A complete metric space $(X,d)$ is called {\em proper} if the closed ball $\{y\in Y\;|\; d(x,y)\leq r\}$ is compact for each $x\in X\,,\; r>0$.

\begin{proposition}
\label{Chap3-Prop-BoundedCompactXContIffPsi_xCont}
Let $\ts$ be a topological space and $(X,d)$ be a proper, complete metric space. The following assertions are equivalent for a family $K_t\in\cs(X)\,,\; t\in\ts$.
\begin{itemize}
\item[(i)] The family $(K_t)_{t\in\ts}$ is Fell-continuous. 
\item[(ii)] The map 
$\Psi_x:\ts\to[0,\infty)\,,\; 
	t\mapsto\inf
		\big\{
			d(\lambda,x)\;\big|\; \lambda\in K_t
		\big\}
$
is continuous for each $x\in X$.
\end{itemize}
\end{proposition}

\begin{proof}
(i)$\Rightarrow$(ii): Define the continuous map $f_x:X\to[0,\infty)\,,\; y\mapsto d(x,y)\,,$ for each $x\in X$. Thus $f_x$ maps compact sets to compact sets. In addition, the maps $f_x\,,\; x\in X\,,$ are closed maps. For each $r>0$, the preimage 
$$
f_x^{-1}\big([0,a]\big) \;
	= \; \big\{ y\in X\;|\; d(x,y)\leq r\big\}
$$ 
is compact since $(X,d)$ is a proper, complete metric space. Thus, $f_x:X\to[0,\infty)$ is a proper map, i.e., preimages of compact sets are compact. Hence, the map $\hat{f_x}:\cs(X)\to\cs\big([0,\infty)\big)$ is Fell-continuous by Lemma~\ref{Chap3-Lem-FunctFellCont}. Then $\Psi_x$ is represented by 
$$
\ts\ni t\; \mapsto\;
	K_t \;\overset{\hat{f_x}}{\longmapsto}\;
		\big\{ d(\lambda,x)\;\big|\; \lambda\in K_t \big\}\;\overset{\inf}{\longmapsto}\;
				\inf\big\{ d(\lambda,x)\;\big|\; \lambda\in K_t \big\} = \Psi_x(t)\,.
$$
Applying Lemma~\ref{Chap3-Lem-InfFellCont} and (i), the map $\ts\ni t\mapsto\Psi_x(t)$ is continuous as a composition of continuous maps.

\vspace{.1cm}

(ii)$\Rightarrow$(i): Let $t_0\in\ts$. Consider a compact $F\in\ks(X)$ and a finite family $\os$ of open subsets of $X$ chosen so that $K_{t_0}\in \us(F,\os)$. The proof is organized as follows. First, the existence of an open neighborhood $U_F\subseteq \ts$ of $t_0$ is verified such that $K_t\cap F=\emptyset$ for all $t\in U_F$. Secondly, for each $O\in\os$, the existence of an open neighborhood $U_O$ of $t_0$ is proven such that $K_t\cap O$ is non-empty for $t\in U_O$. Then the finite intersection of the open sets $U_F$ and $U_O,\; O\in\os,$ defines an open neighborhood containing $t_0$ such that $K_t\in\us(F,\os)$ for all $t\in U_F\cap\bigcap_{O\in\os} U_O$.  

\vspace{.1cm}

Since $F$ and $K_t$ are closed and since $F\cap K_{t_0}=\emptyset$, there exists, for each given $x\in F$, an $r(x)>0$ so that $B_{r(x)}(x)\cap K_{t_0}=\emptyset$. The family of (smaller) open balls $\{B_{r(x)/2}(x)\;|\; x\in F\}$ covers $F$. By compactness of $F$, there is a finite set $\{x_1,\cdots, x_l\}\subseteq F$ with radii $r_k:=r(x_k)$ for $1\leq k\leq l$ such that
$$
F\subseteq \bigcup_{k=1}^l B_{r_k/2}(x_k) \, ,
	\qquad
	B_{r_k}(x_k)\cap K_{t_0}=\emptyset\, .
$$
Let $1\leq k\leq l$. The relation $B_{r_k}(x_k)\cap K_{t_0}=\emptyset$ is equivalent to $\Psi_{x_k}(t_0)\geq r_k$. The continuity of $\Psi_{x_k}:\ts\to[0,\infty)$ implies that there is an open neighborhood $U_k$ of $t_0$ such that $t\in U_k$ leads to $\Psi_{x_k}(t_0)\geq r_k/2$. For $U_F:=\bigcap_{k=1}^l U_k\subseteq U_0$, it follows from the previous bounds that $B_{r_k/2}(x_k)\cap K_t=\emptyset$ for $t\in U_F$ and $k\in\{1,2,\cdots,l\}$. Consequently, the intersection $F\cap K_t$ is empty for all $t\in U_F$ since $B_{r_k/2}(x_k),\, k\in\{1,2,\cdots,l\}\,,$ define a covering of $F$. 

\vspace{.1cm}

Now let $O\in\os$. The intersection $O\cap K_{t_0}\neq \emptyset$ is non-empty by assumption. For any $x\in O\cap K_{t_0}$, there is an $r:=r(x)>0$ such that $B_r(x)\subseteq O$ as $O$ is open in $X$. Since $x\in K_{t_0}$, the intersection $B_{r}(x)\cap K_{t_0}$ is non-empty. Thus, $\Psi_x(t_0)<r$ follows. Using the continuity of $\Psi_x:\ts\to[0,\infty)$, there is an open neighborhood $U_O$ of $t_0$ in $\ts$ such that the inequality $\Psi_x(t)<r$ holds for $t\in U_O$. This leads to $O\cap K_t \supseteq B_{r}(x)\cap K_t\neq \emptyset$ for $t\in U_O$ by definition of $\Psi_x$. 

\vspace{.1cm}

Since the family $\os$ is finite, the intersection $U_\os:=\bigcap_{O\in\os} U_O$ is open and contains $t_0$. Consequently, $U_\os\cap U_F\subseteq U_0$ is an open neighborhood of $t_0$. Altogether, the compact set $K_t$ satisfies $K_t\in\us(F,\os)$ for all $t\in U_\os\cap U_F$.
\end{proof}

\begin{corollary}
\label{Chap3-Cor-BoundedCompactXContIffPsi_xCont}
Let $\ts$ be a topological space, $(X,d)$ be a proper, complete metric space and $K\in\ks(X)$ be compact. Then the following assertions are equivalent for a family $K_t\subseteq K\,,\; t\in\ts\,,$ of compact subsets.
\begin{itemize}
\item[(i)] The family $(K_t)_{t\in\ts}$ is Vietoris-continuous. 
\item[(ii)] The map 
$\Psi_x:\ts\to[0,\infty)\,,\; 
	t\mapsto\inf
		\big\{
			d(\lambda,x)\;\big|\; \lambda\in K_t
		\big\}\,,
$
is continuous for each $x\in\RM$.
\end{itemize}
\end{corollary}

\begin{proof}
Clearly, a proper, complete metric space $(X,d)$ is locally compact. Thus, the equivalence follows by Theorem~\ref{Chap3-Theo-VietFellHausMetricEquiv} and Proposition~\ref{Chap3-Prop-BoundedCompactXContIffPsi_xCont}.
\end{proof}

\begin{remark}
\label{Chap3-Rem-BoundedCompactRMContIffPsi_xCont}
The space $\RM$ with the Euclidean metric $d(x,y):=|x-y|\,,\; x,y\in\RM\,,$ is a proper, complete metric space. Thus, Proposition~\ref{Chap3-Prop-BoundedCompactXContIffPsi_xCont} and Corollary~\ref{Chap3-Cor-BoundedCompactXContIffPsi_xCont} applies for $X=\RM$ which is used in Theorem~\ref{Chap3-Theo-CharContUnbSelfAdjOpSpecGap}.
\end{remark}

\section{Continuous behavior of the spectra for self-adjoint operators}
\label{Chap3-Sect-CharContSpectSelfAdj}

Following the lines of \cite{BeBe16}, a characterization for the continuous behavior of the spectra of a bounded self-adjoint field of operators $(A_t)_{t\in\ts}$ is proven, c.f. Theorem~\ref{Chap3-Theo-P2ContEquivContSpect}. More precisely, the behavior of the norms of certain polynomials up to degree two of the operators encode all the information of the spectra as a set. If only polynomials up to degree one are considered, this leads to a characterization of the continuity of the supremum and infimum of the spectra, c.f. Theorem~\ref{Chap3-Theo-CharInfSupSpectCont}. Note that the continuity of $\ts\ni t\mapsto\inf\big(\sigma(A_t)\big)\in\RM$ is usually referred to the continuity of the bottom of the spectra $(\sigma(A_t))_{t\in\ts}$ in the literature.

\medskip

The characterization of the continuous behavior of the spectra is based on the fact that the spectrum $\sigma(A)$ of a bounded self-adjoint operator $A\in\Ll(\hs)$ is compact and contained in $\RM$, see e.g. \cite[Proposition~1.3.9.]{Dixmier77}.

\medskip

Define a bounded operator $p(A)$ for a real-valued polynomial $p:\RM\to\RM$ and $A\in\Ll(\hs)$. In the following, the notation \gls{Pt} is used for the set of real-valued polynomials up to degree $2$, i.e., 
$$
\Pt \; 
	:=\; \big\{ p_2\cdot z^2 +p_1\cdot z+p_0 \;\big|\; p_2,p_1,p_0\in\RM \big\} \, .
$$
With this notion at hand, the \pt-continuity of a field of bounded operators is introduced along the lines of \cite{BeBe16}.

\begin{definition}[\pt-continuous]
\label{Chap3-Def-(P2)-Continuity}
Let $\ts$ be a topological space and $(A_t)_{t\in\ts}$ be a field of bounded operators. Then the field $(A_t)_{t\in\ts}$ is called {\em \pt-continuous} whenever all the maps $\Phi_p:\ts\to[0,\infty)\,,\; t\mapsto \|p(A_t)\|\,,$ for $p\in\Pt$ are continuous. 
\end{definition}

Note that instead of considering all polynomials of $\Pt$, it is sufficient for our purposes to define the \pt-continuity only for a specific subset of $\Pt$, c.f. Remark~\ref{Chap3-Rem-PolyP2ContEquivContSpect}.

\begin{remark}
\label{Chap3-Rem-(P2)-Continuity}
It is worth noticing that the bounded operators $A_t,\; t\in\ts,$ may be defined on different Hilbert spaces $\hs_t,\; t\in\ts$. On the other hand, the notions of operator norm continuity and strong operator continuity make only sense if all operators are defined on the same Hilbert space $\hs$. Whenever the field of Hilbert spaces $(\hs_t)_{t\in\ts}$ is constant, i.e., $\hs_t=\hs$ for all $t\in\ts$, then the \pt-continuity is expressed in terms of a topology on $\Ll(\hs)$. Details on this topology are provided in Section~\ref{Chap3-Sect-RelP2ContOthTop}.

\medskip

Note, furthermore, that there are no requirements on the topology of $\ts$.
\end{remark}

In case that $(A_t)_{t\in\ts}$ is a field of bounded self-adjoint operators, the \pt-continuity turns out to be the correct condition in terms of the continuous variation of the spectra $\sigma(A_t)\in\ks(\RM),\; t\in\ts,$ with respect to the Hausdorff metric. Lemma~\ref{Chap3-Lem-PresenceSpectrum} is the key observation. 

\medskip

\begin{figure}[htb]
\centering
\includegraphics[scale=1.2]{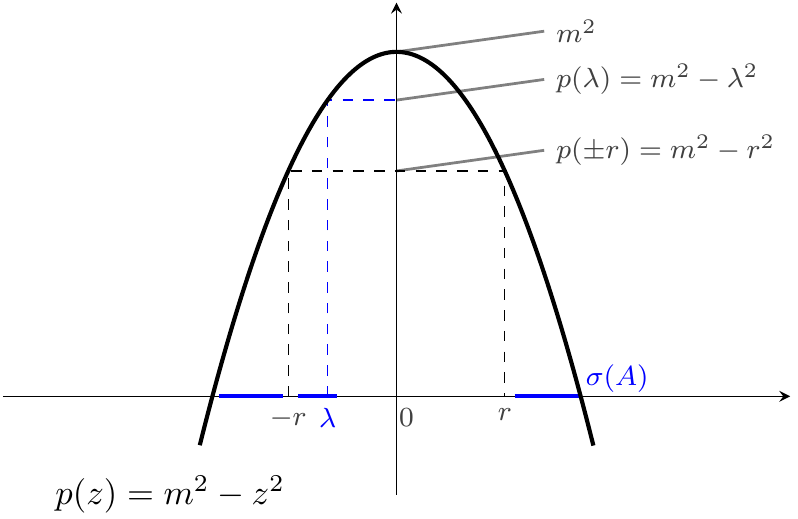}
\caption{Illustration of the proof of Lemma~\ref{Chap3-Lem-PresenceSpectrum}.}
\label{Chap3-Fig-p(z)Spectrum}
\end{figure}

In the following $B_r(x)$ denotes the {\em open ball} with radius $r>0$ around $x\in\RM$. The idea of the proof of the following assertion is sketched in Figure~\ref{Chap3-Fig-p(z)Spectrum}.

\begin{lemma}
\emph{(\cite[Lemma~1]{BeBe16})}
\label{Chap3-Lem-PresenceSpectrum}
Let $A\in\Ll(\hs)$ be a bounded self-adjoint operator on a Hilbert space $\hs$. Consider the polynomial $p(z):=m^2-z^2$ where $m\geq\|A\|$. Then the following assertions hold for every $r<m$.
\begin{itemize}
\item[(a)] The inequality $\|p(A)\|\leq m^2-r^2$ holds if and only if $B_r(0)\cap\sigma(A)=\emptyset$.
\item[(b)] The inequality $\|p(A)\|> m^2-r^2$ holds if and only if $B_r(0)\cap\sigma(A)\neq\emptyset$.
\end{itemize}
\end{lemma}

\begin{proof}
Since assertion (a) and (b) are equivalent, the proof of (a) is provided, only. The operator $A$ is self-adjoint and bounded implying that the spectrum $\sigma(A)\subseteq\RM$ is compact. As $m\geq\|A\|$ and $\|A\|\geq|\lambda|$ for all $\lambda\in\sigma(A)$, the estimates $p(\lambda)\geq m^2-\|A\|^2>0$ follows on $\sigma(A)$. Hence, the equation $\|p(A)\|=\sup_{\lambda\in\sigma(A)} p(\lambda)$ holds by using $p(\sigma(A))=\sigma(p(A))$.

\vspace{.1cm}

Then the intersection $B_r(0)\cap\sigma(A)$ is empty if and only if $|\lambda|\geq r$ for all $\lambda\in\sigma(A)$. This is equivalent to the inequality $m^2-\lambda^2\leq m^2-r^2$ for all $\lambda\in\sigma(A)$. Finally, $m^2-\lambda^2\leq m^2-r^2$ holds for all $\lambda\in\sigma(A)$ if and only if $\sup_{\lambda\in\sigma(A)}m^2-\lambda^2\leq m^2-r^2$. By the previous considerations, the identity $\|p(A)\|=\sup_{\lambda\in\sigma(A)}m^2-\lambda^2$ concludes the proof.
\end{proof}

\begin{remark}
\label{Chap3-Rem-PresenceSpectrum}
The requirement $m\geq\|A\|$ is necessary since then $m^2-A^2$ is a positive operator and so $\sup_{\lambda\in\sigma(A)}|m^2-\lambda^2|=\sup_{\lambda\in\sigma(A)}m^2-\lambda^2$ holds.
\end{remark}

\medskip

Proposition~\ref{Chap3-Prop-ConClosFuncHausCont} leads to the following result.

\begin{lemma}
\emph{(\cite[Lemma~4]{BeBe16})}
\label{Chap3-Lem-SpecContNormPhiCont}
Let $\ts$ be a topological space and $(A_t)_{t\in\ts}$ be a field of bounded normal operators. Consider a continuous and closed function $\phi:\CM\to\CM$. Then the map $\ts\ni t\mapsto\|\phi(A_t)\|\in[0,\infty)$ is continuous if the family of spectra $(\sigma(A_t))_{t\in\ts}$ is Vietoris-continuous in $\ks(\CM)$.
\end{lemma}

\begin{proof}
Let $\phi:\CM\to\CM$ be continuous and closed. Since $A_t\in\Ll(\hs_t)$ is bounded for $t\in\ts$, the spectrum $\sigma(A_t)$ is a compact subset in $\CM$ and $\phi(\sigma(A_t))=\sigma(\phi(A_t))\subsetneq\CM$ holds by the functional calculus. Then $\phi(A_t)$ is a bounded normal operator. Thus, the spectral theorem applies for $\phi(A_t)$. Consider the function $\hg:\cs(\CM)\to\cs([0,\infty))$ where $g:\CM\to[0,\infty),\; z\mapsto |z|,$ is the absolute value. Then the norm map is represented by
$$
\ts\ni t
	\longmapsto\sigma(A_t)
		\overset{\hphi}{\longmapsto}\sigma(\phi(A_t))
			\overset{\hg}{\longmapsto}|\sigma(\phi(A_t))|
				\overset{\sup}{\longmapsto}\|\phi(A_t)\|\in[0,\infty) \, .
$$
The map $g$ is continuous and closed. Hence, $\hg$ is Vietoris-continuous as well as $\hphi$ by Proposition~\ref{Chap3-Prop-ConClosFuncHausCont}. The supremum on $\cs(\RM)$ is also continuous, c.f. Proposition~\ref{Chap3-Prop-SupInfVietCont}. Finally, the Vietoris-continuity of $(\sigma(A_t))_{t\in\ts}$ implies that the norm map $t\mapsto\|\phi(A_t)\|$ is continuous as a composition of continuous maps. 
\end{proof}

\medskip

Let $(A_t)_{t\in\ts}$ be a field of bounded self-adjoint operators. Lemma~\ref{Chap3-Lem-SpecContNormPhiCont} yields to the continuity of $\ts\ni t\mapsto\|\phi(A_t)\|\in[0,\infty)$ for each continuous and closed function $\phi:\RM\to\RM$ under the assumption that $(\sigma(A_t))_{t\in\ts}$ is Vietoris-continuous in $\cs(\RM)$

\begin{theorem}
\emph{(\cite[Theorem~1]{BeBe16})}
\label{Chap3-Theo-P2ContEquivContSpect}
Let $\ts$ be a topological space and $(A_t)_{t\in\ts}$ be a bounded self-adjoint field of operators. Then the following assertions are equivalent.
\begin{itemize}
\item[(i)] The map $\Sigma:\ts\to\ks(\RM)\,,\; t\mapsto\sigma(A_t)\,,$ is continuous with respect to the Hausdorff metric on $\ks(\RM)$.
\item[(ii)] The field $(A_t)_{t\in\ts}$ is \pt-continuous.
\item[(iii)] The boundaries of $(\sigma(A_t))_{t\in\ts}$ are continuous.
\end{itemize}
\end{theorem}

\begin{figure}[htb]
\centering
\includegraphics[scale=1.6]{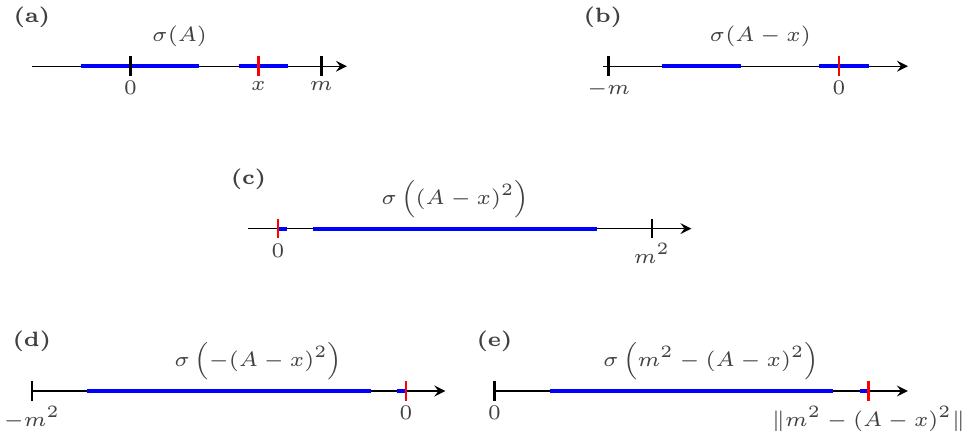}
\caption{Illustration of the utility of the polynomials of degree $2$ to control the spectra. In detail, the construction of a polynomial $p$ of degree $2$ is sketched so that $p(x)=\|p(A)\|$ holds for an $x\in\sigma(A)$. (a) The spectrum $\sigma(A)$ is contained in $[-m,m]$. (b) The spectrum is translated such that $x$ gets the origin. (c) The spectrum is folded at the new origin by taking the square. This compresses and stretches the spectrum. (d) Multiplication by $-1$ reflects the spectrum to the negative part of the real line. (e) Translate the spectrum by $m^2$ such that it is contained in the positive part of the real line. Using the fact that $\|p(A)\|=\sup_{\lambda\in\sigma(A)}|p(\lambda)|$ leads to the desired identity $p(x)=\|p(A)\|$ where $p(z):=m^2-(z-x)^2$.}
\label{Chap3-Fig-ThmContSpectSelAdjOp}
\end{figure}

\begin{proof}
Since the operators $A_t,\; t\in\ts,$ are self-adjoint and bounded, the spectra $\sigma(A_t),\; t\in\ts,$ are compact subsets of $\RM$. Then the equivalence of (i) and (iii) follows by Theorem~\ref{Chap3-Theo-BoundVietContCompVers}. According to Theorem~\ref{Chap3-Theo-VietFellHausMetricEquiv}, the continuity with respect to the Hausdorff metric is equivalent to the Vietoris-continuity.

\vspace{.1cm}

(i)$\Rightarrow$(ii): Every polynomial $p:\RM\to\RM$ of degree smaller than or equal to $2$ is continuous and closed. Then (ii) follows by Lemma~\ref{Chap3-Lem-SpecContNormPhiCont}.

\vspace{.1cm}

(ii)$\Rightarrow$(i): The main ingredient for this implication is Lemma~\ref{Chap3-Lem-PresenceSpectrum} that allows to control the absence and the presence of the spectrum by the norm of certain polynomials of the operator. Let $t_0\in\ts$. The \pt-continuity implies that $\ts\ni t\mapsto\|A_t\|\in[0,\infty)$ is continuous. Thus, there exists an open neighborhood $U_0\subseteq\ts$ of $t_0$ such that $\sup_{t\in U_0}\|A_t\|\leq \|A_{t_0}\|+1$.

\vspace{.1cm}

Consider a closed $F\in\cs(\RM)$ and a finite family $\os$ of open subsets of $\RM$ chosen so that $\sigma(A_{t_0})\in \us(F,\os)$. The proof is organized as follows. First, the existence of an open neighborhood $U_F\subseteq U_0$ of $t_0$ is verified such that $\sigma(A_t)\cap F=\emptyset$ for all $t\in U_F$. Secondly, for each $O\in\os$, the existence of an open neighborhood $U_O$ of $t_0$ is proven such that $\sigma(A_t)\cap O$ is non-empty for $t\in U_O$. Then the finite intersection of the open sets $U_F$ and $U_O,\; O\in\os,$ defines an open neighborhood containing $t_0$ satisfying $\sigma(A_t)\in\us(F,\os)$ for all $t\in U_F\cap\bigcap_{O\in\os} U_O$. In order to do so, define $m:=2\cdot \|A_{t_0}\|+1$. 

\vspace{.1cm}

For $t\in U_0$, the intersection $F\cap\sigma(A_t)$ is empty if and only if 
$$
\Big( F\cap[-\|A_{t_0}\|-1,\|A_{t_0}\|+1] \Big)\cap \sigma(A_t)=\emptyset.
$$ 
Consider the compact subset $K:=F\cap[-\|A_{t_0}\|-1,\|A_{t_0}\|+1]$ of $\RM$. Since $K$ and $\sigma(A_{t_0})$ are closed and since $K\cap\sigma(A_{t_0})=\emptyset$, there exists, for each given $x\in K$, an $0<r(x)<\|A_{t_0}\|$ so that $B_{r(x)}(x)\cap \sigma(A_{t_0})=\emptyset$. The family of (smaller) open balls $\{B_{r(x)/2}(x)\;|\; x\in K\}$ covers $K$. By compactness of $K$, there is a finite set $\{x_1,\cdots, x_l\}\subseteq K$ with radii $r_k:=r(x_k)$ for $1\leq k\leq l$ such that
$$
K\subseteq \bigcup_{k=1}^l B_{r_k/2}(x_k) \, ,
	\qquad
	B_{r_k}(x_k)\cap\sigma(A_{t_0})=\emptyset\, .
$$
Note that the estimates $\|A_{t_0}-x_k\|\leq \|A_{t_0}\|+|x_k|\leq m$ hold. Thus, using Lemma~\ref{Chap3-Lem-PresenceSpectrum}, the condition $B_{r_k}(x_k)\cap\sigma(A_{t_0})=\emptyset$ is equivalent to $\|m^2-(A_{t_0}-x_k)^2\|\leq m^2-r_k^2$. The \pt-continuity implies that there is an open neighborhood $U_k\subseteq U_0$ of $t_0$ such that $t\in U_k$ leads to $\|m^2-(A_{t}-x_k)^2\|\leq m^2-r_k^2/4$. For $U_F:=\bigcap_{k=1}^l U_k\subseteq U_0$, it follows from the previous bounds and from the Lemma~\ref{Chap3-Lem-PresenceSpectrum} that $B_{r_k/2}(x_k)\cap\sigma(A_{t})=\emptyset$ for $t\in U_F$ and $k\in\{1,2,\cdots,l\}$. Consequently, the intersection $K\cap\sigma(A_t)$ is empty as $B_{r_k/2}(x_k),\, k\in\{1,2,\cdots,l\},$ define a covering of $K$. Since, additionally, $U_F\subseteq U_0$, the intersection 
$$
\Big(\big(-\infty,-\|A_{t_0}\|-1\big]
	\cup\big[\|A_{t_0}\|+1,\infty\big)\Big)
	\cap\sigma(A_{t})\,,
	\qquad t\in U_F\,,
$$ 
is empty. Combined with $K\cap\sigma(A_t)=\emptyset$ this leads to $F\cap\sigma(A_t)=\emptyset$ for all $t\in U_F$.

\vspace{.1cm}

Now let $O\in\os$. The intersection $O\cap \sigma(A_{t_0})\neq \emptyset$ is non-empty by assumption. For each $x\in O\cap \sigma(A_{t_0})$, there is an $r:=r(x)>0$ such that $B_r(x)\subseteq O$ as $O$ is open in $\RM$. Since the inequality $|x|\leq \|A_{t_0}\|$ holds for $x\in \sigma(A_{t_0})$, it follows that $\|A_{t_0}-x\|\leq 2\|A_{t_0}\|\leq m$. Then Lemma~\ref{Chap3-Lem-PresenceSpectrum} implies that $\|m^2-(A_{t_0}-x)^2\|>m^2-r^2/4$ as $x\in\sigma(A_{t_0})$. Using the \pt-continuity, there is an open neighborhood $U_O$ of $t_0$ in $\ts$ such that the inequality $\|m^2-(A_{t}-x)^2\|>m^2-r^2$ holds for $t\in U_O$. This leads to $O\cap \sigma(A_t)\supseteq B_{r}(x)\cap \sigma(A_t)\neq \emptyset$ for $t\in U_O$. 

\vspace{.1cm}

Since the family $\os$ is finite, the intersection $U_\os:=\bigcap_{O\in\os} U_O$ is open and contains $t_0$. Consequently, $U_\os\cap U_F\subseteq U_0$ is an open neighborhood of $t_0$. Altogether, the spectrum of $A_t$ satisfies $\sigma(A_t)\in\us(F,\os)$ for all $t\in U_\os\cap U_F$.
\end{proof}

\begin{remark}
\label{Chap3-Rem-PolyP2ContEquivContSpect}
Note that not all polynomials of degree $2$ are needed to control the spectra for a fixed field $(A_t)_{t\in\ts}$ of operators. Specifically, $\big(\sigma(A_t)\big)_{t\in\ts}$ is continuous with respect to the Hausdorff metric if and only if all the maps
$$
\ts\ni t\mapsto \big\|m^2-(A_t-x)^2\big\|\,
	\qquad m\in\RM\, , \; x\in\RM\, ,
$$
are continuous. The idea of the utility of these polynomials is sketched in Figure~\ref{Chap3-Fig-ThmContSpectSelAdjOp}. Actually, it suffices to fix the parameter $m>0$ for a given field of self-adjoint operators and an open set $U\subseteq\ts$. More precisely, let $U\subseteq\ts$ be an open neighborhood of $t_0\in\ts$ such that $m:= 2\sup_{t\in U}\|A_t\|+1$ is finite. Then the parameter $x$ only has to ranges over the interval $[-m,m]$, instead of $\RM$, so that the spectra behave continuous. In this sense, the continuous behavior of the spectra only depends on one parameters $x\in[-m,m]$ for a fixed bounded self-adjoint field of operators. This is similar for fields of unitary operators, c.f. Section~\ref{Chap3-Sect-CharContSpectUnit}. In this case, the analog of the parameter $m$ is fixed by $1$ since the operator norm of unitary operators is equal to $1$. Then the continuous behavior of the spectra also only depends on one parameter, c.f. Remark~\ref{Chap3-Rem-OneParaUnitary}. Note that, for each $t_0\in\ts$, there exists an open neighborhood $U\subseteq\ts$ of $t_0$ such that $\sup_{t\in U}\|A_t\|$ is finite if $\ts\ni t\mapsto\|A_t\|\in[0,\infty)$ is continuous
\end{remark}

The following remark is important to deal with the unbounded case.

\begin{remark}
\label{Chap3-Rem-P2ContEquivContSpect}
The assertion of Theorem~\ref{Chap3-Theo-P2ContEquivContSpect} does also follow directly from Corollary~\ref{Chap3-Cor-BoundedCompactXContIffPsi_xCont}. More specifically, let $t_0\in\ts$. The \pt-continuity of $(A_t)_{t\in\ts}$ leads to the existence of an open neighborhood $U_0$ of $t_0$ satisfying $m:=\sup_{t\in U_0}\|A_t\|<\infty$. Then
$$
\big\|\big(m+|x|\big)^2-(A_t-x)^2\big\|\;
	= \; 
		\big(m+|x|\big)^2-
		\left(\inf\left.\big\{|\lambda-x| \;\right|\;
			\lambda\in\sigma(A_t)\big\}\right)^2\,,
			\qquad
			x\in\RM\,,
$$
follow. Consequently, for each $x\in\RM$, the continuity of the norm is equivalent to the continuity of the map $\Psi_x:\ts\to[0,\infty)\,,\; t\mapsto\inf\big\{|\lambda-x| \;\big|\; \lambda\in\sigma(A_t)\big\}\,,$ as $\Psi_x$ is non-negative. Then Corollary~\ref{Chap3-Cor-BoundedCompactXContIffPsi_xCont} and Remark~\ref{Chap3-Rem-BoundedCompactRMContIffPsi_xCont} imply the Vietoris continuity of $U_0\ni t\mapsto \sigma(A_t)\in\ks(\RM)$ as $\sigma(A_t)\subseteq [-m,m]$ holds for all $t\in U_0$. Since the geometric argument of Lemma~\ref{Chap3-Lem-PresenceSpectrum} provides a good intuition, the proof of Theorem~\ref{Chap3-Theo-P2ContEquivContSpect} is not presented with the help of Corollary~\ref{Chap3-Cor-BoundedCompactXContIffPsi_xCont}. The connection with the maps $\Psi_x\,,\; x\in\RM\,,$ is also used in Theorem~\ref{Chap3-Theo-(R)-Continuity} to show the Fell-continuity of the spectra associated with unbounded self-adjoint operators.
\end{remark}

The assertion of Theorem~\ref{Chap3-Theo-P2ContEquivContSpect} cannot be improved in the sense that not all polynomials up to degree are considered. More specifically, the following three examples show that the degree of the polynomials cannot be decreased and that the constant term of the polynomials is necessary, in general.

\begin{example}
\label{Chap3-Ex-Pol2Necess}
Consider the topological space $\ts:=\oNM$ and the constant field of Hilbert spaces $\hs_n:=\CM^3$ for $n\in\oNM$. Define the field of bounded self-adjoint operators on $\CM^3$ by
$$
A_n\; 
	:=\; \begin{pmatrix}
		-1 & 0 & 0\;\\
		0 & 0 & 0\\
		0 & 0 & 1
	\end{pmatrix} \, ,
	\quad
	n\in\NM \,,
	\qquad
	A_\infty \; 
		:= \; \begin{pmatrix}
			-1 & 0 & 0\;\\
			0 & \frac{1}{2} & 0\\
			0 & 0 & 1
		\end{pmatrix} \, .
$$
Clearly, the spectra are given by $\sigma(A_n)=\{-1,0,1\},\; n\in\NM,$ and $\sigma(A_\infty)=\{-1,\frac{1}{2},1\}$ and $\lim_{n\to\infty}\sigma(A_n)\neq \sigma(A_\infty)$. In order to check the continuity of the norms of $(p(A_n))_{n\in\oNM}$ the behavior at infinity in $\oNM$ has to be studied. Consider a polynomial $p(z):=p_1\cdot z+p_0$ with $p_0,p_1\in\RM$. Then the identities
\begin{gather*}
\begin{aligned} 
\|p(A_n)\|\; 
	&= \; \max\{|\pm p_1 + p_0|, |p_0|\} \; 
	&= \; |p_1| + |p_0| \, ,\qquad &n\in\NM \, ,\\
\|p(A_\infty)\|\; 
	&= \; \max\left\{
		|\pm p_1 + p_0|, 
		\left|\frac{p_1}{2} + p_0\right|
	\right\} \; 
	&= \; |p_1| + |p_0| \, ,\qquad &
\end{aligned} 
\end{gather*}
hold. Thus, the map $\oNM\ni n\mapsto\|p(A_n)\|\in[0,\infty)$ is continuous for each polynomial up to degree $1$. On the other hand, it is clear that the map $\oNM\ni n\mapsto\sigma(A_n)\in\ks(\RM)$ is not continuous at infinity, i.e., the spectra do not converge in the Hausdorff metric. Hence, polynomials of degree smaller than or equal to $1$ are not sufficient to check the continuity in this simple example. For the polynomial $q(z):=1-z^2$, the equations
$$
\|q(A_n)\| \; 
	= \; 1 \, ,
	,\qquad
	\|q(A_\infty)\|\;
		=\; \frac{3}{4} \, ,
			\qquad
			n\in\NM \, ,
$$
hold. Consequently, $\lim_{n\to\infty}\|q(A_n)\|\neq\|q(A_\infty)\|$ is satisfied implying that $(A_n)_{n\in\oNM}$ is not \pt-continuous which was expected by the previous considerations.
\end{example}

The continuity of $\oNM\ni n\mapsto\|p(A_n)\|\in[0,\infty)$ for polynomials up to degree $1$ relies on the fact that $\oNM\ni n\mapsto\sup\big(\sigma(A_n)\big)\in\RM$ and $\oNM\ni n\mapsto\inf\big(\sigma(A_n)\big)\in\RM$ are continuous in this example. Thus, this example indicates that polynomials of degree $1$ only control the supremum and the infimum of the spectra. This holds for indeed, c.f. Theorem~\ref{Chap3-Theo-CharInfSupSpectCont}.

\medskip

The following example shows that it is not sufficient to consider only polynomials of the form $p(z)=p_2\cdot z^2 + p_0$ for $p_2,p_0\in\RM$.

\begin{example}
\label{Chap3-Ex-xTermPolNecess}
Consider the topological space $\ts:=\oNM$ and the constant field of Hilbert spaces $\hs_n:=\CM^2$ for $n\in\oNM$. Define the field of bounded self-adjoint operators on $\CM^2$ by
$$
A_n\; 
	:=\; \begin{pmatrix}
		\frac{1}{2} & 0 \;\\
		0 & 1
	\end{pmatrix} \, ,
	\quad
	n\in\NM \,,
	\qquad
	A_\infty \; 
		:= \; \begin{pmatrix}
			-\frac{1}{2} & 0\;\\
			0 & -1
		\end{pmatrix} \, .
$$
The spectra are given by $\sigma(A_n)=\big\{\frac{1}{2},1\big\}\,,\; n\in\NM\,,$ and $\sigma(A_\infty)=\big\{-\frac{1}{2},-1\big\}$. Clearly, the spectra do not converge in the Hausdorff metric, i.e., $\lim_{n\to\infty}\sigma(A_n)\neq\sigma(A_\infty)$. Consider a polynomial $p(z)=p_2\cdot z^2 + p_0$ for $p_2,p_0\in\RM$. Then the identities
$$
\|p(A_n)\| \; 
	= \; \max\left\{ 
		\left|\frac{p_2}{4}+p_0 \right|,
		\left|p_2+p_0 \right|
	\right\}\;
	= \; \|p(A_\infty)\|\,,
	\qquad
	n\in\NM\,,
$$
are derived implying the continuity of $\oNM\ni n\mapsto\|p(A_n)\|\in[0,\infty)$ for all polynomials of the form $p(z)=p_2\cdot z^2 + p_0$. On the other hand, for the polynomial $q(z)=1-z$, the equations $\|q(A_n)\|=\frac{1}{2}$ and $\|q(A_n)\|=2$ are deduced. Thus, $(A_n)_{n\in\oNM}$ is not \pt-continuous.
\end{example}

In view of Example~\ref{Chap3-Ex-xTermPolNecess}, the linear part of the polynomials controls the sign of the spectra from the philosophical point of view since the square part destroy any information about the sign. 

\medskip

The following example shows that the constant part of the polynomials is also necessary.

\begin{example}
\label{Chap3-Ex-ConstTermPolNecess}
Recall the Example~\ref{Chap3-Ex-StrongConvNotConvSpect}. Consider a polynomial $p(z):=p_2 \cdot z^2 + p_1 \cdot z$ for $p_2,p_1\in\RM$ without constant term. Note that $A_n^2=A_n$ and $\|A_n\|=1$ hold for all $n\in\oNM$. Thus, the identities
$$
\|p(A_n)\|\; 
	= \; |p_2 + p_1| \cdot \|A_n\| \; = \: |p_2 + p_1|\,,
	\qquad n\in\oNM,
$$
are derived leading to the continuity of the norms $\oNM\ni n\mapsto\|p(A_n)\|\in[0,\infty)$ for all polynomials without constant term. Note that this continuity holds even for polynomials of higher degree than $2$ as long as it has no constant term. In this example, the equalities $\|q(A_n)\|=1,\; n\in\NM,$ and $\|q(A_\infty)\|=0$ hold for the polynomial $q(z):=1-z$ implying that $(A_n)_{n\in\oNM}$ is not \pt-continuous.
\end{example}

It is natural to ask whether the \pt-continuity is stable under linear combination of opera\-tors and composition of operators. Let $(A_t)_{t\in\ts}$ and $(B_t)_{t\in\ts}$ be \pt-continuous bounded self-adjoint fields of operators for a topological space $\ts$. The following two examples show that the pointwise sum $(A_t+B_t)_{t\in\ts}$ and the pointwise composition $(A_tB_t)_{t\in\ts}$ is not \pt-continuous, in general.

\begin{example}
\label{Chap3-Ex-SumNotP2Cont}
Let $\ts:=\oNM$ and $\hs:=\hs_n:=\CM^2$ for $n\in\oNM$. Define the self-adjoint operators $A_n,B_n:\hs\to\hs$ by
\begin{align*}
A_n \; := \;\begin{pmatrix}
				1 & \frac{1}{2}\\
				\frac{1}{2} & 1
		\end{pmatrix}
		\, ,\; n\in\NM\,,\qquad
	&A_\infty \; := \;\begin{pmatrix}
					\frac{1}{2} & 0\\
					0 & \frac{3}{2}
				\end{pmatrix}\,,\\
B_n \; := \;\begin{pmatrix}
				1 & -\frac{1}{2}\\
				-\frac{1}{2} & 1
		\end{pmatrix}
		\, ,\; n\in\NM\,,\qquad
	&B_\infty \; := \;\begin{pmatrix}
					\frac{1}{2} & 0\\
					0 & \frac{3}{2}
				\end{pmatrix}\,.
\end{align*}
It is easy to check that $\sigma(A_n)=\sigma(B_n)=\{\frac{1}{2},\, \frac{3}{2}\}$ holds for all $n\in\oNM$. Hence, the self-adjoint fields $(A_n)_{n\in\oNM}$ and $(B_n)_{n\in\oNM}$ are \pt-continuous. Then the equations $A_\infty+B_\infty=2\cdot A_\infty$ and $A_n+B_n=2 \cdot I\,,\; n\in\NM\,,$ follow where $I:\hs\to\hs$ is the identity. Hence, $(A_n+B_n)_{n\in\oNM}$ is not \pt-continuous since $\sigma(A_\infty+B_\infty)=\{1,\,3\}$ and $\sigma(A_n+B_n)=\{2\}\,,\;n\in\NM\,,$ hold.
\end{example}

Example~\ref{Chap3-Ex-SumNotP2Cont} shows that the \pt-topology related to the \pt-continuity is not a vector space topology, c.f. Section~\ref{Chap3-Sect-RelP2ContOthTop}.

\begin{example}
\label{Chap3-Ex-CompositNotP2Cont}
Let $\ts:=\oNM$ and $\hs:=\hs_n:=\CM^2$ for $n\in\oNM$. Define the self-adjoint operators $A_n,B_n:\hs\to\hs$ by
\begin{gather*}
A_n \; := \;\begin{pmatrix}
				1 & 0\\
				0 & 0
		\end{pmatrix}
		\, ,\; n\in\NM\,,\qquad
	A_\infty \; := \;\begin{pmatrix}
					0 & 0\\
					0 & 1
				\end{pmatrix}\,, \\				
B_n \; := \;\begin{pmatrix}
				0 & 0\\
				0 & 1
		\end{pmatrix}
		\, ,\; n\in\oNM \,.
\end{gather*}
Clearly, the equalities $\sigma(A_n)=\sigma(B_n)=\{1,\, 0\}$ hold for all $n\in\oNM$. Thus, the self-adjoint fields $(A_n)_{n\in\oNM}$ and $(B_n)_{n\in\oNM}$ are \pt-continuous. Furthermore, the equations $A_\infty B_\infty= B_\infty$ and  $A_n B_n=0\,,\; n\in\NM\,,$ hold. Hence, $(A_n B_n)_{n\in\oNM}$ is not \pt-continuous since $\sigma(A_\infty B_\infty)=\{1,\, 0\}$ and $\sigma(A_n B_n)=\{0\}\,,\; n\in\NM\,,$ hold.
\end{example}

Denote by \gls{Po} the set of complex-valued polynomials $z\mapsto p(z,\overline{z})$ in $z$ and $\overline{z}$. Consider a $p\in\Po$ and a bounded normal operator $A$ on a Hilbert space $\hs$. Then $p(A,A^\ast)$ defines a normal bounded operator on $\hs$. Whenever $A$ is self-adjoint, i.e., $A=A^\ast$, then the notation $p(A)$ is used instead of $p(A,A^\ast)$.

\begin{corollary}
\label{Chap3-Cor-AllPolNormCont}
Let $\ts$ be a topological space and $(A_t)_{t\in\ts}$ be a field of bounded self-adjoint operators. Then the following assertions are equivalent.
\begin{itemize}
\item[(i)] The field $(A_t)_{t\in\ts}$ is \pt-continuous.
\item[(ii)] The maps $\ts\ni t\mapsto\|p(A_t)\|\in[0,\infty),\; p\in\Po,$ are continuous.
\end{itemize}
\end{corollary}

\begin{proof}
Theorem~\ref{Chap3-Theo-VietFellHausMetricEquiv}, Theorem~\ref{Chap3-Theo-P2ContEquivContSpect} and Lemma~\ref{Chap3-Lem-SpecContNormPhiCont} lead to the implication (i) to (ii) whereas the converse is clear by Definition~\ref{Chap3-Def-(P2)-Continuity} and the inclusion $\Pt\subsetneq\Po$.
\end{proof}

\begin{remark}
\label{Chap3-Rem-AllPolNormCont}
Note that  the desired continuity in Corollary~\ref{Chap3-Cor-AllPolNormCont} means that $(A_t)_{t\in\ts}$ is \po-continuous, c.f. Definition~\ref{Chap3-Def-(P)-Continuity} and Theorem~\ref{Chap3-Theo-PContEquivContSpectNormal} below.
\end{remark}

In specific cases, it might be that the bounded field of operators is not self-adjoint, but it can be a function of a bounded self-adjoint field of operator. Then the continuous behavior of the spectra is reduced to the continuous behavior of the spectra of the self-adjoint field of operators.

\begin{corollary}
\label{Chap3-Cor-CharContSpectFunctOfSelfAdj}
Let $\ts$ be a topological space and $(A_t)_{t\in\ts}$ be a bounded field of operators. Suppose that there exist a field of bounded self-adjoint operators $(B_t)_{t\in\ts}$ and a continuous and injective function $\phi:X\subseteq\CM\to\CM$ such that $X$ is closed with $\sigma(B_t)\subseteq X$ for $t\in\ts$ and $A_t=\phi(B_t)$ holds for all $t\in\ts$. Then $(A_t)_{t\in\ts}$ is a bounded normal field and the following assertions are equivalent.
\begin{itemize}
\item[(i)] The map $\Sigma:\ts\to\ks(\RM)\,,\; t\mapsto\sigma(A_t)\,,$ is continuous with respect to the Hausdorff metric on $\ks(\RM)$.
\item[(ii)] The field $(B_t)_{t\in\ts}$ is \pt-continuous.
\end{itemize}
\end{corollary}

\begin{proof}
Since $(B_t)_{t\in\ts}$ is bounded and self-adjoint, $(A_t)_{t\in\ts}$ defines a bounded normal field of operators by the functional calculus.

\vspace{.1cm}

(i)$\Rightarrow$(ii): According to Theorem~\ref{Chap3-Theo-P2ContEquivContSpect} it suffices to show that $\sigma(B_t)_{t\in\ts}$ is continuous with respect to the Hausdorff metric which is equivalent to its Vietoris-continuity, c.f. Theorem~\ref{Chap3-Theo-VietFellHausMetricEquiv}. Let $t_0\in\ts$. According to Lemma~\ref{Chap3-Lem-SpecContNormPhiCont}, there exists an open neighborhood $U_0$ of $t_0$ such that $m:=\sup_{t\in U_0}\|A_t\|$ is finite. The set $K:=\overline{\bigcup_{t\in U_0}\sigma(A_t)}\subseteq X\cap \overline{B_{m}(0)}$ is by construction compact. Consequently, the map $\phi:K\to\phi(K)$ restricted to the compact set $K$ is a homeomorphism with inverse $\phi^{-1}:\phi(K)\to K$ by general topological arguments, see e.g. \cite[Satz~8.11]{Querenburg2001}. Then a general version of Tietze's Extension Theorem applies, c.f. \cite[Satz~I.8.5.4]{Schubert75} and Theorem~\ref{App1-Theo-TietzesThm}. Specifically, there exists a continuous function $\Phi:\CM\to\CM$ with compact support such that $\Phi(z)=\phi^{-1}(z),\; z\in K$. The map $\Phi$ is closed since it has compact support and it is continuous. Due to construction, the composition $\Phi\circ\phi:\CM\to\CM$ is continuous and closed. Furthermore, the identity $\Phi\big(\phi(z)\big) =z$ holds for $z\in K$ leading to $\Phi(A_t) = B_t$. With this at hand, the map $\ts\ni t\mapsto\sigma(B_t)\in\ks(\CM)$ is represented by
$$
\ts\ni t \;
	\longmapsto \; \sigma(A_t) \;
		\overset{\widehat{\Phi}}{\longmapsto} \; \Phi\big(\sigma(A_t)\big) \;
			= \; \sigma(B_t)
				\in\ks(\CM)
$$
where $\widehat{\Phi}:\cs(\CM)\to\cs(\CM)$ denotes the map $F\mapsto\Phi(F)$. The leftmost map is Vietoris-continuous by (i) and $\widehat{\Phi}$ is Vietoris-continuous by Proposition~\ref{Chap3-Prop-ConClosFuncHausCont}. Thus, the map $\ts\ni t\mapsto\sigma(B_t)\in\ks(\CM)$ is Vietoris-continuous as composition of continuous maps. Consequently, assertion (ii) follows by Theorem~\ref{Chap3-Theo-P2ContEquivContSpect}.

\vspace{.1cm}

(ii)$\Rightarrow$(i): According to Theorem~\ref{Chap3-Theo-P2ContEquivContSpect} the map $\ts\ni t\mapsto\sigma(B_t)\in\ks(\CM)$ is continuous with respect to the Hausdorff metric. By the functional calculus, the identity $\phi(\sigma(B_t))=\sigma(\phi(B_t))=\sigma(A_t)$ holds for all $t\in\ts$. Thus, Proposition~\ref{Chap3-Prop-ConClosFuncHausCont} implies the desired Vietoris-continuity of $(\sigma(A_t))_{t\in\ts}$ in $\ts$.
\end{proof}

\medskip

For certain applications, the continuity of the bottom of the spectra is important. For instance, this question is studied for Gabor frame operators associated with a Delone set, c.f. \cite{GrOrRo15,Krei16}. For this continuity, only polynomials up to degree $1$ are necessary as already indicated in Example~\ref{Chap3-Ex-Pol2Necess}.

\begin{theorem}
\label{Chap3-Theo-CharInfSupSpectCont}
Let $\ts$ be a topological space and $(A_t)_{t\in\ts}$ be a bounded self-adjoint field of operators. Then the following assertions are equivalent.
\begin{itemize}
\item[(i)]The maps $\ts\ni t\mapsto\inf\big(\sigma(A_t)\big)\in\RM$ and $\ts\ni t\mapsto\sup\big(\sigma(A_t)\big)\in\RM$ are continuous.
\item[(ii)] The maps $\ts\ni t\mapsto\|m-A_t\|\in[0,\infty)\,,\; m\in\RM\,,$ are continuous.
\end{itemize}
\end{theorem}

\begin{proof}
Let $t_0\in\ts$. Note that $\ts\ni t\mapsto\|m+A_t\|\in[0,\infty)$ is continuous by replacing $m$ by $-m$ if $\ts\ni t\mapsto\|m-A_t\|\in[0,\infty)$ is continuous.

\vspace{.1cm}

(i)$\Rightarrow$(ii): For every $m\in\RM$, the identity 
$$
\|m-A_t\| \; 
	= \; \max\big\{ 
			\big|m-\inf\big(\sigma(A_t)\big)\big|, 
			\big|m-\sup\big(\sigma(A_t)\big)\big| 
		\big\} \, ,
		\qquad t\in\ts \, ,
$$ 
holds leading to the continuity of the norms.

\vspace{.1cm}

(ii)$\Rightarrow$(i): According to (ii), the map $\ts\ni t\mapsto\|A_t\|\in[0,\infty)$ (corresponding to the choice $m=0$) is continuous. Thus, there is an open neighborhood $U_0$ of $t_0$ such that $\|A_t\|\leq m:=\|A_{t_0}\|+1$ for $t\in U_0$. Then $m\pm \lambda\geq 0$ holds for all $\lambda\in\sigma(A_t)$ and $t\in U_0$. Thus, the equations 
$$
\|m+A_t\| \;
	\; =m+\sup\big(\sigma(A_t)\big) \, , \qquad
\|m-A_t\|\; 
	= \; m-\inf\big(\sigma(A_t)\big) \, ,
$$
hold for $t\in U_0$. Then the continuity of the norms implies the continuity of the supremum and the infimum at $t_0\in\ts$.
\end{proof}

\begin{remark}
\label{Chap3-Rem-PolGra1Gra2}
According to Theorem~\ref{Chap3-Theo-CharInfSupSpectCont}, polynomials of degree $1$ control the behavior of the infimum and supremum of the spectra of a bounded self-adjoint field of operators. On the other hand, specific polynomials of degree $2$ are needed to handle also the behavior in between of the infimum and the supremum of the spectra, c.f. Remark~\ref{Chap3-Rem-PolyP2ContEquivContSpect}.
\end{remark}

\section{Relation between the (\textit{p2})-topology and other topologies on \texorpdfstring{$\Ll(\mathcal{H})$}{L(H)}}
\label{Chap3-Sect-RelP2ContOthTop}

Let $\hs$ be a Hilbert space. In Remark~\ref{Chap3-Rem-(P2)-Continuity}, it was noticed that the \pt-continuity is related to a topology on $\Ll(\hs)$. The so called \pt-topology is defined in this section and the relation to other topologies on $\Ll(\hs)$ is discussed.

\medskip

Let $X$ be a set and $\tau_1,\tau_2$ be two topologies on $X$. Then $\tau_1$ is called {\em finer} than $\tau_2$ if each $U\in\tau_2$ belongs to $\tau_1$. In this case, convergence of a sequence (or a net) $(x_n)_{n}$ in the topology $\tau_1$ implies that $(x_n)_{n}$ also converges in the topology of $\tau_2$. Both topologies $\tau_1$ and $\tau_2$ coincide if $\tau_1$ is finer than $\tau_2$ and vice versa. The reader is referred to \cite[Chapter~2]{Querenburg2001} for more details. In this section, $X$ is the space $\Ll(\hs)$ of all bounded operators over the Hilbert space $\hs$.

\medskip

Recall that $\Pt$ denotes the set of real-valued polynomials up to degree two. 

\begin{definition}[$1$-norm]
\label{Chap3-Def-1NormPol}
Consider a polynomial $p(z):=p_2\cdot z^2 + p_1\cdot z + p_0$ in $\Pt$ where $p_2,\; p_1,\; p_0\in\RM$. Then the {\em $1$-norm $\|p\|_1$ of $p$} is defined by the sum $|p_2|+|p_1|+|p_0|$. For $M>0$, the set of all $p\in\Pt$ with $\|p\|_1\leq M$ is denoted by \gls{PtM}.
\end{definition}

According to Remark~\ref{Chap3-Rem-PolyP2ContEquivContSpect}, only specific polynomials of the form $p(z):=m^2-(z-x)^2\,,\; z\in\RM\, ,$ are necessary and sufficient to control the behavior of the spectra for a given bounded self-adjoint field of operators where $m>0$ is fixed and $x\in[-m,m]$. The $1$-norm of these polynomials is estimated by
$$
\|p\|_1 \;
	= \; 1+2|x| + m^2-x^2\;
		= \; m^2+ 2 - (1-|x|)^2 \;
			\leq \; m^2+2\, .
$$
Thus, it suffices to control only the polynomials of $\Pt(M)$ for $M\geq m^2+2$. In view of that the \pt-topology is defined as follows.

\begin{definition}[\pt-topology]
\label{Chap3-Def-P2Top}
Let $\hs$ be a Hilbert space. The {\em \pt-topology on $\Ll(\hs)$} is defined by the neighborhood basis
$$
\us_{\varepsilon,M}(A_0) \,
	:= \, \left\{
		A\in\Ll(\hs) \;|\; \big| \|p(A_0)\|-\|p(A)\| \big|<\varepsilon \text{ for all } p\in\Pt(M)
	\right\}
	 \, ,\;\; \varepsilon>0 \, ,\, M>0 \, ,
$$
for $A_0\in\Ll(\hs)$. 
\end{definition}

It is not difficult to check that the sets $\us_{\varepsilon,M}(A_0),\, \varepsilon>0,\, M>0, A_0\in\Ll(\hs),$ define a basis for a topology on $\Ll(\hs)$ in terms of Definition~\ref{App1-Def-BaseTopology}. It is worth noticing that the $1$-norm of the polynomials need to be bounded in order to get a reasonable topology. Otherwise, a self-adjoint $A\in\Ll(\hs)$ is contained in an $\varepsilon$-neighborhood of $A_0$ only if $\sigma(A)=\sigma(A_0)$ for every $\varepsilon>0$. This is proved as follows: Without loss of generality, let $\lambda\in\sigma(A)\setminus\sigma(A_0)$ and $m\geq 2\cdot \max\{\|A\|,\|A_0\|\}$. Consider the polynomial $p(z):= m^2-(z-\lambda)^2$. The operators $p(A)$ and $p(A_0)$ are non-negative by definition, i.e., the spectra $\sigma(p(A))$ and $\sigma(p(A_0))$ are contained in $[0,\infty)$. Then the equalities 
$$
\|p(A)\| \;
	= \; m^2\,, 
	\qquad 
		\|p(A_0)\| \;
			= \; \sup\big\{m^2-(\mu-\lambda)^2\;\big|\; \mu\in\sigma(A_0)\big\}
$$
hold. By construction, the norm $\|p(A_0)\|$ is smaller and not equal to $m^2$. Consequently, $C:=\|p(A)\|-\|p(A_0)\|>0$ follows. Hence,
$$
\big| \|q_n(A)\| - \|q_n(A_0)\| \big| \; 
	=\; n\cdot C\,,
	\qquad
	n\in\NM\,,
$$
is derived where $q_n(z):=n\cdot p(z)$. Thus, for each $\varepsilon>0$, there exists a polynomial $q\in\Pt$ such that $\big| \|q(A)\|-\|q(A_0)\| \big|>\varepsilon$ if the $1$-norm of the polynomial $q$ is not bounded.

\medskip

In view of Theorem~\ref{Chap3-Theo-P2ContEquivContSpect}, continuity in the \pt-topology has implications for the spectrum. 

\begin{proposition}
\label{Chap3-Prop-SpecEpsClosP}
Let $A_0\in\Ll(\hs)$ be self-adjoint and $M>8(\|A_0\|+1)^2$. Then the estimate $d_H(\sigma(A_0),\sigma(A))<\varepsilon$ holds for each $0<\varepsilon<1$ and all self-adjoint $A\in\us_{\varepsilon^2/8,M}(A_0)$.
\end{proposition}

\begin{proof}
Let $1>\varepsilon>0$ and $A\in\us_{\varepsilon^2/8,M}(A_0)$ be self-adjoint. Note that the difference $\big| \|A\|-\|A_0\| \big|$ is smaller than $\varepsilon^2/8$ as $M>1$. Thus, $\sigma(A)$ and $\sigma(A_0)$ are contained in the interval $K:=[-\|A_0\|-1,\|A_0\|+1]$. Fix $m:=2\|A_0\|+2$. For $x\in K$, consider the polynomial $p_x(z):=m^2-(z-x)^2$. Then the estimates
$$
\|p_x\|_1 \; 
	= \; m^2+ 2 - (1-|x|)^2 \;
		\leq \; m^2+2 \; 
			= \; 4\cdot\|A_0\|^2+8\cdot\|A_0\|+6 \;
				< \; M\,,
				\qquad x\in K\,,
$$
follow. Thus, the inequalities
\begin{equation}
\label{Chap3-Eq-pNormEps} \tag{$\lozenge$}
\big| \|p_x(A_0)\|-\|p_x(A)\| \big| \;
	< \; \frac{\varepsilon^2}{8} \,,
	\qquad x\in K\,,
\end{equation}
hold. Now it suffices to prove that 
$$
\sigma(A)
	\subseteq U(\varepsilon,A_0):=\bigcup_{x\in\sigma(A_0)} B_\varepsilon(x)
	 \, ,\qquad
\sigma(A_0)
		\subseteq U(\varepsilon,A):=\bigcup_{x\in\sigma(A)} B_\varepsilon(x) \, .
$$
The proof follows the lines of Theorem~\ref{Chap3-Theo-P2ContEquivContSpect}. Define the closed set $F:=\RM\setminus U(\varepsilon,A_0)$. Then $\tilde{F}:=F\cap K$ is compact and so there are $x_1,\ldots,x_l\in \tilde{F}$ such that
$$
\tilde{F}\subseteq\bigcup_{k=1}^l B_{\varepsilon/2}(x_k)\, .
$$
By construction, the intersection $B_{\varepsilon}(x_k)\cap\sigma(A_0)$ is empty for $1\leq k\leq l$. Note that $x_1,\ldots,x_l\in K$. Thus, by using (\ref{Chap3-Eq-pNormEps}) and Lemma~\ref{Chap3-Lem-PresenceSpectrum}~(a), it follows that
$$
\big\|m^2-(A-x_k)^2\big\| \; 
	< \; \frac{\varepsilon^2}{8}+\big\|m^2-(A_0-x_k)^2\big\| \;
		\leq\; \frac{\varepsilon^2}{8} + m^2-\varepsilon^2
			<\; m^2-\frac{\varepsilon^2}{4} \, .
$$
Applying Lemma~\ref{Chap3-Lem-PresenceSpectrum}~(a), this implies $B_{\varepsilon/2}(x_k)\cap\sigma(A)=\emptyset$ for $1\leq k\leq l$. Hence, the intersection $\sigma(A)\cap F$ is empty leading to $\sigma(A)\subseteq U(\varepsilon,A_0)$.

\medskip

The family $B_{\varepsilon/4}(x)$ for $x\in\sigma(A_0)$ defines an open cover of $\sigma(A_0)$. By compactness there are $x_1,\ldots,x_k\in\sigma(A_0)\subseteq K$ such that $\sigma(A_0)\subseteq \bigcup_{l=1}^k B_{\varepsilon/4}(x_l)$. According to (\ref{Chap3-Eq-pNormEps}) and Lemma~\ref{Chap3-Lem-PresenceSpectrum}~(b) it follows that 
$$
\big\|m^2-(A-x_l)^2\big\|\; 
	> \; \big\|m^2-(A_0-x_l)^2\big\| -\frac{\varepsilon^2}{8} \;
		> \;  m^2-\frac{\varepsilon^2}{16}-\frac{\varepsilon^2}{8} \;
			> \;  m^2-\frac{\varepsilon^2}{4} \, .
$$
Then Lemma~\ref{Chap3-Lem-PresenceSpectrum}~(b) implies $B_{\varepsilon/2}(x_l)\cap\sigma(A)\neq\emptyset$. Hence, there is a $y_l\in\sigma(A)$ such that $B_{\varepsilon/2}(x_l)\subseteq B_{\varepsilon}(y_l)$. Thus, the inclusions $\sigma(A_0)\subseteq\bigcup_{l=1}^k B_{\varepsilon}(y_l)\subseteq U(\varepsilon,A)$ follow.
\end{proof}

\medskip

It is well-known that the operator norm topology implies the convergence of the spectra.

\begin{proposition}
\label{Chap3-Prop-NormImplP2}
The operator norm topology is finer than the \pt-topology on $\Ll(\hs)$.
\end{proposition}

\begin{proof}
According to Lemma~\ref{App1-Lem-TopFiner}, it suffices to show that, for $A_0\in\Ll(\hs)$ and $\varepsilon, M>0$, there exists an open neighborhood $U$ of $A_0$ in the operator norm topology such that $U\subseteq\us_{\varepsilon,M}(A_0)$.

\vspace{.1cm}

Let $\us_{\varepsilon,M}(A_0)$ be a neighborhood in the \pt-topology of $A_0\in\Ll(\hs)$ for $M>0$ and $\varepsilon>0$. Denote by $U\subsetneq \Ll(\hs)$ the set of all $A\in\Ll(\hs)$ satisfying
$$
\|A-A_0\| \; 
	< \; \max\left\{
			\frac{\varepsilon}{2M\cdot (\|A_0\|+1)}
			, 1
		\right\} \, .
$$ 
The set $U$ is an open neighborhood of $A_0$ in the operator norm topology. In the following, it is shown that $A\in\us_{\varepsilon,M}(A_0)$ whenever $A\in U$. Consider a $p\in\Pt(M)$ with parameters $p_2,\; p_1,\; p_0\in\RM$. Then the estimates
\begin{align*}
\big| \|p(A)\|-\|p(A_0)\| \big| \;
	&\leq \; \| p(A)-p(A_0)\| \\
		& \leq |p_2|\cdot \Big(\big\|A(A- A_0)\big\| + \big\|A_0 (A - A_0)\big\|\Big) + |p_1|\cdot \|A-A_0\|\\
			&\leq \; \left(|p_2| \cdot \big(\| A \| + \| A_0\|\big) + |p_1|\right) \cdot \| A- A_0\| \\
				&\leq \; 2 M (\|A_0\|+1) \cdot \| A- A_0\|\\
					&< \; \varepsilon
\end{align*}
hold finishing the proof.
\end{proof}

\begin{remark}
\label{Chap3-Rem-StrongFinDim}
Let $\hs$ be a finite dimensional Hilbert space. It is well-known that then the strong operator topology and the norm operator topology on $\Ll(\hs)$ coincide. Thus, if $\hs$ is finite dimensional, the strong operator topology is finer than the \pt-topology by Proposition~\ref{Chap3-Prop-NormImplP2}.
\end{remark}

This fails as soon as $\hs$ is infinite dimensional.

\begin{proposition}
\label{Chap3-Prop-StrongNotImplP2}
Let $\hs$ be infinite dimensional. Then the strong operator topology on $\Ll(\hs)$ is not finer than the \pt-topology.
\end{proposition}

\begin{proof}
If the strong operator topology were finer than the \pt-topology, the convergence of a sequence $(A_n)_{n\in\NM}$ to $A\in\Ll(\hs)$ in the strong operator topology implies the convergence in the \pt-topology. Thus, it suffices to find a sequence $(A_n)_{n\in\NM}$ converging to $A\in\Ll(\hs)$ in the strong operator topology such that the sequence does not converge in the \pt-topology.

\vspace{.1cm}

Let $(\psi_n)_{n\in\NM}$ be a countable orthonormal system in $\hs$, i.e., $\|\psi_n\|=1$ and $\langle\psi_n|\psi_m\rangle=0$ hold for $n,m\in\NM$. Denote by $U$ the closure of the linear span of this orthonormal system. Then $\hs$ is represented by the direct sum $U\oplus V$ where $V$ is the orthogonal complement of $U$. Denote by $P_V$ the projection on the space $V$. Define the operators $A_n\in\Ll(\hs)$ for $n\in\NM$ by
$$
A_n\varphi \; 
	:= \; \sum\limits_{m=1}^n 
		\langle\psi_m|\varphi\rangle \cdot \psi_m
		+ P_V\varphi \, .
$$
By standard arguments the norm $\|A_n\varphi\|$ is bounded by $\|\varphi\|$. Furthermore, $(A_n)_{n\in\NM}$ converges strongly to the identity $I$ on $\hs$. On the other hand, the operators $A_n\,,\; n\in\NM\,,$ and the identity $I$ are self-adjoint. In addition, the spectra $\sigma(A_n)=\{0,1\}\,,\; n\in\NM\,,$ do not converge to $\sigma(I)=\{1\}$ in the Hausdorff metric. Thus, $(A_n)_{n\in\NM}$ does not converge in the \pt-topology by Proposition~\ref{Chap3-Prop-SpecEpsClosP}.
\end{proof}

\medskip

The \pt-topology allows one to control the spectrum as a set, c.f. Theorem~\ref{Chap3-Theo-P2ContEquivContSpect}. On the other hand, it does not encode any information about the eigenspaces or the multiplicity of elements of the spectrum. This implies that the \pt-topology is not finer than the strong operator topology.

\begin{proposition}
\label{Chap3-Prop-P2NotImplStrong}
Let $\hs$ be a Hilbert space that is not one-dimensional. Then the \pt-topology on $\Ll(\hs)$ is not finer than the strong operator topology.
\end{proposition}

\begin{proof}
Similar to the previous proposition, it suffices to find a sequence $(A_n)_{n\in\NM}$ converging to $A\in\Ll(\hs)$ in the \pt-topology such that the sequence does not converge in the strong operator topology.

\vspace{.1cm}

Fix $\psi_1,\; \psi_2\in\hs$ such that $\langle\psi_i|\psi_j\rangle = \delta_{ij}$ for $1\leq i,j\leq 2$. Since $\hs$ is not one-dimensional, the existence of $\psi_1$ and $\psi_2$ is guaranteed. Define the operators $A,\; A_n\in\Ll(\hs)$ for $n\in\NM$ by
$$
A\varphi \; 
	:= \; \langle\psi_1|\varphi\rangle \cdot \psi_1 \, ,
	\qquad
A_n\varphi \; 
	:= \; \langle\psi_2|\varphi\rangle \cdot \psi_2\, ,
	\qquad
	\varphi\in\hs\,.
$$
Since $\langle\psi_1|\psi_2\rangle=0$, the equalities $\sigma(A_n)=\sigma(A)=\{0,1\}\,,\; n\in\NM\,,$ follow. Hence, $(A_n)_{n\in\NM}$ converges to $A$ in the \pt-topology. On the other hand, $\|(A_n-A)\psi_1\|=1$ holds for all $n\in\NM$. Thus, $\|(A_n-A)\psi_1\|$ does not tend to zero if $n\to\infty$. Consequently, $(A_n)_{n\in\NM}$ does not converge to $A$ in the strong operator topology.
\end{proof}

\medskip

By the previous considerations, it follows that the vector space $\Ll(\hs)$ equipped with the \pt-topology is not a topological vector space. More precisely, the addition is not \pt-continuous, c.f. also Example~\ref{Chap3-Ex-SumNotP2Cont}.

\begin{remark}
\label{Chap3-Rem-OneDimHS}
Let $\hs$ be a one-dimensional Hilbert space. Each bounded linear oper\-ator on $\hs$ is a multiplication operator by a constant. Thus, the spectrum of such an operator does only contain one value with multiplicity one. Then the convergence of the spectra for an operator sequence $A_n\in\Ll(\hs)\,,\; n\in\NM\,,$ to the spectrum of $A\in\Ll(\hs)$ implies that these operators converge in operator norm. Specifically, the operator norm topology, the strong operator topology and the \pt-topology coincide whenever $\hs$ is one-dimensional.
\end{remark}

\section{Continuous behavior of the spectra for unbounded self-ad\-joint operators}
\label{Chap3-Sect-CharContSpecUnbouSelfAdjOp}

Let $\ts$ be a topological space and $(A_t)_{t\in\ts}$ be a field of unbounded self-adjoint operators. Under a weak assumption on $(A_t)_{t\in\ts}$, a characterization of the continuous behavior of the spectra in the Vietoris-topology is proven, c.f. Theorem~\ref{Chap3-Theo-CharContUnbSelfAdjOpSpecGap}. Whenever this assumption is dropped, the continuity of the spectra with respect to the Fell-topology is characterized which is weaker than the Vietoris-continuity, c.f. Example~\ref{Chap3-Ex-VietFinerThanFell}.

\medskip

Let $(X,d)$ be a complete metric space. As discussed in Section~\ref{Chap3-Sect-VietorisContinuity}, the space $\ks(X)$ of compact subsets of $X$ equipped with the Hausdorff metric is a complete metric space. The induced topology on $\ks(X)$ is equal to the topology induced by the Vietoris-topology, c.f. Theorem~\ref{Chap3-Theo-VietFellHausMetricEquiv}. In the following the case, $X=\RM$ equipped with the Euclidean metric is considered and the spectra of unbounded self-adjoint operators are analyzed for their continuous behavior in the Vietoris-topology.

\medskip

Clearly, it makes no sense to consider the norm $\|A_t\|$ if $A_t$ is an unbounded operator for $t\in\ts$. The idea is to pass to the resolvent and study the continuous behavior of the norms of the resolvent. The following assertion immediately follows by Theorem~\ref{Chap3-Theo-P2ContEquivContSpect}.

\begin{theorem}
\label{Chap3-Theo-CharContUnbSelfAdjOpSpecGap}
Let $\ts$ be a topological space and $(A_t)_{t\in\ts}$ be a self-adjoint field of operators that are not necessarily bounded such that there exists an $x\in\RM$ with $x\not\in\sigma(A_t)$ for all $t\in\ts$. Then the following assertions are equivalent.
\begin{itemize}
\item[(i)] The map $\Sigma:\ts\to\cs(\RM)\,,\; t\mapsto\sigma(A_t)\,,$ is continuous in the Vietoris-topology on $\cs(\RM)$.
\item[(ii)] The family of spectra $\big(\sigma\big((A_t-x)^{-1}\big)\big)_{t\in\ts}$ is continuous in the Vietoris-topology on $\ks(\RM)$.
\item[(iii)] All the maps
$$
\ts\ni t\mapsto 
	\left\| p\left( 
		(A_t-x)^{-1}
	\right) \right\|\in[0,\infty)\, ,
		\qquad p\in\Pt\, ,
$$
are continuous, i.e., $\big((A_t-x)^{-1}\big)_{t\in\ts}$ is \pt-continuous.
\end{itemize}
\end{theorem}

\begin{proof}
Define the function $\phi:\RM\setminus\{x\}\to\RM$ by $\phi(z):=(z-x)^{-1}$. Since $x\in\RM$ is not in $\sigma(A_t)$ for all $t\in\ts$ the operators $\phi(A_t)$ are well-defined. Thus, for $t\in\ts$, the operator $(A_t-x)^{-1}$ is bounded and self-adjoint. According to the functional calculus, the equation $\sigma(\phi(A_t))=\phi(\sigma(A_t))$ holds. The function $\phi$ is continuous and invertible on $\RM\setminus\{x\}$ with continuous inverse $z\mapsto x+z^{-1}$. Consequently, Proposition~\ref{Chap3-Prop-ConClosFuncHausCont} implies that $\sigma(A_t)$ is continuous if and only if $\sigma(\phi(A_t))$ is continuous in the Vietoris-topology. Then $\sigma(\phi(A_t))$ is continuous if and only if $\ts\ni t\mapsto\|p(\phi(A_t))\|\in[0,\infty)$ is continuous for all $p\in\Pt$, c.f. Theorem~\ref{Chap3-Theo-P2ContEquivContSpect}.
\end{proof}

\medskip

If the field of self-adjoint operators does not have an $x\in\RM\setminus\sigma(A_t)$ for all $t\in\ts$, then the restriction to a subset $U\subseteq\ts$ with induced topology is helpful. More precisely, for $t_0\in\ts$ and $x\in\RM\setminus\sigma(A_{t_0})$, an open neighborhood $U$ of $t_0$ exists such that $x\in\RM\setminus\sigma(A_t)$ for $t\in U$ whenever the spectra are continuous in the Vietoris-topology on $\cs(\RM)$.

\medskip

The following statement is an easy consequence of the previous theorem.

\begin{corollary}
\label{Chap3-Cor-CharContUnbSelfAdjOpSpecGap}
Let $\ts$ be a topological space and $(A_t)_{t\in\ts}$ be a field of self-adjoint opera\-tors that are not necessarily bounded. Suppose that $(A_t)_{t\in\ts}$ is uniformly bounded from below, i.e., there exists a $C\in\RM$ such that $\inf_{t\in\ts}\inf\big(\sigma(A_t)\big)> C$. Then the spectra $\big(\sigma(A_t)\big)_{t\in\ts}$ are continuous in the Vietoris-topology if and only if all the maps $\ts\ni t\mapsto \big\|p\big((A_t-C)^{-1}\big)\big\|\in[0,\infty)\, ,\; p\in\Pt\, ,$ are continuous.
\end{corollary}

A similar version of Theorem~\ref{Chap3-Theo-CharContUnbSelfAdjOpSpecGap} holds if the resolvent is replaced by the exponential map $\RM\ni z\mapsto exp(-a\cdot z)$ for an $a>0$. 

\begin{theorem}
\label{Chap3-Theo-CharContSpectrHalbgruppe}
Let $\ts$ be a topological space and $(A_t)_{t\in\ts}$ be a self-adjoint field of operators that are not necessarily bounded. Suppose that $(A_t)_{t\in\ts}$ is uniformly bounded from below, i.e., there exists a $C\in\RM$ such that $\inf_{t\in\ts}\inf\big(\sigma(A_t)\big)> C$. Then the following are equivalent for every $a>0$.
\begin{itemize}
\item[(i)] The map $\Sigma:\ts\to\cs(\RM)\,,\; t\mapsto\sigma(A_t)\,,$ is continuous in the Vietoris-topology on $\cs(\RM)$.
\item[(ii)] The family of spectra $\big(\sigma\big(exp(-a\cdot A_t)\big)\big)_{t\in\ts}$ is continuous in the Vietoris-topology on $\ks(\RM)$.
\item[(iii)] The field of bounded self-adjoint operators $\big(exp(-a\cdot A_t)\big)_{t\in\ts}$ is \pt-continuous.
\end{itemize}
\end{theorem}

\begin{proof}
If the field $(A_t)_{t\in\ts}$ is uniformly bounded from below, $\big(exp(-a\cdot A_t)\big)_{t\in\ts}$ defines a bounded self-adjoint field of operators for each $a>0$. Thus, the equivalence of (ii) and (iii) follows by Theorem~\ref{Chap3-Theo-P2ContEquivContSpect}. The exponential function $exp:\RM\to(0,\infty)$ is continuous, open, closed and bijective, namely it has a continuous inverse which is also a closed map. Then Proposition~\ref{Chap3-Prop-ConClosFuncHausCont} leads to the equivalence of (i) and (ii).
\end{proof}

\medskip

In general, it may happen that there does not exist an $x\in\RM$ satisfying $x\not\in\sigma(A_t)$ for all $t\in\ts$ (as in Theorem~\ref{Chap3-Theo-CharContUnbSelfAdjOpSpecGap}). On the other hand, the resolvent $(A_t-z)^{-1}$ exists for all $z\in\CM\setminus\RM$ as $A_t$ is self-adjoint. This leads to the following definition which is due to \cite{BeBe16}.

\begin{definition}[$(R)$-continuous, \cite{BeBe16}]
\label{Chap3-Def-(R)-Continuity}
Let $\ts$ be a topological space and $(A_t)_{t\in\ts}$ be a field of self-adjoint operators that are not necessarily bounded. Then $(A_t)_{t\in\ts}$ is called $(R)$-continuous whenever all the maps
$$
\ts\ni t\longmapsto 
	\left\|(A_t-z)^{-1}\right\|\in[0,\infty)\, ,
		\qquad z\in\CM\setminus\RM\, ,
$$
are continuous.
\end{definition}

The $(R)$-continuity of $(A_t)_{t\in\ts}$ can be characterized by the Fell-continuity of the associated spectra by taking the considerations at the end of Section~\ref{Chap3-Sect-VietorisContinuity} into account. Recall that the Fell-continuity of the spectra is interpreted as a local Vietoris-continuity of the spectra, c.f. Section~\ref{Chap3-Sect-VietorisContinuity}.

\begin{theorem}
\emph{(\cite[Theorem~2]{BeBe16})}
\label{Chap3-Theo-(R)-Continuity}
Let $\ts$ be a topological space and $(A_t)_{t\in\ts}$ be a self-adjoint field of operators that are not necessarily bounded. Then the following are equivalent.
\begin{itemize}
\item[(i)] The map map $\Sigma:\ts\to\cs(\RM)\,,\; t\mapsto\sigma(A_t)\,,$ is Fell-continuous.
\item[(ii)] The field $(A_t)_{t\in\ts}$ is $(R)$-continuous.
\item[(ii)] The boundaries of $\big(\sigma(A_t)\big)_{t\in\ts}$ are continuous.
\end{itemize}
\end{theorem}

\begin{proof}
Let $z=x+\imath y\in\CM$ be so that $y\neq 0$. Since, for $t\in\ts$, the operator $A_t$ is self-adjoint on the Hilbert space $\hs_t$, the identities
$$
\left\| (A_t-z)^{-1} \right\|^2 \; 
	=\; \left\| \big((A_t-z)^\ast\big)^{-1} \; (A_t-z)^{-1} \right\|\;
	=\; \left\| \big((A_t-x)^2+y^2\big)^{-1} \right\|
$$
hold for every $t\in\ts$. Then by using the spectral theorem the equations 
$$
\left\|(A_t-z)^{-1} \right\|^2 \; 
	= \; \Big( 
			\left(
				\inf\big\{
					|\lambda-x| \;\big|\; \lambda\in\sigma(A_t) 
				\big\}
			\right)^2
			+ y^2 
		\Big)^{-1} \, ,
	\qquad t\in\ts \, ,
$$
are deduced. Define the maps $\Psi_x:\ts\to\RM,\; x\in\RM,$ by
$$
\Psi_x(t)\; 
	:= \; \inf\big\{
		|\lambda-x| \;\big|\; \lambda\in\sigma(A_t) 
	\big\}
	 \, ,\qquad t\in\ts \, .
$$
For each $t\in\ts$, the spectrum $\sigma(A_{t})$ is non-empty and so $\Psi_x(t)$ is finite. Thus, $(A_t)_{t\in\ts}$ is $(R)$-continuous if and only if all the maps $\Psi_x:\ts\to\RM\, , x\in\RM\, ,$ are continuous. According to Proposition~\ref{Chap3-Prop-BoundedCompactXContIffPsi_xCont} and Remark~\ref{Chap3-Rem-BoundedCompactRMContIffPsi_xCont}, the continuity of all the maps $\Psi_x\,,\; x\in\RM\,,$ is equivalent to the Fell-continuity of $(\sigma(A_t))_{t\in\ts}$. The equivalence of (i) and (iii) follows by Lemma~\ref{Chap3-Lem-VietContImplGapCont} and an adjustment of Lemma~\ref{Chap3-Lem-GapContImplVietCont}, c.f. Remark~\ref{Chap3-Rem-BoundVietContCompVers}.
\end{proof}

\begin{remark}
\label{Chap3-Rem-RContEquivFellContSpect}
According to Example~\ref{Chap3-Ex-VietFinerThanFell}, the Fell-continuity does not imply the continuity in the Vietoris-topology in general.

\vspace{.1cm}

The resolvent transforms the closed set into a compact set. Then the Vietoris-continuity of these compact sets is equivalent to the Fell-continuity of the original closed sets from the philosophical point of view.
\end{remark}

\section{Continuous behavior of the spectra for unitary operators}
\label{Chap3-Sect-CharContSpectUnit}

In this section, the continuous behavior of the spectra associated with a unitary family of operators is characterized. Like in the self-adjoint case, this extends to functions of unitary operators.

\medskip

Let $(A_t)_{t\in\ts}$ be a field of bounded self-adjoint operators over a topological space $\ts$. According to Theorem~\ref{Chap3-Theo-P2ContEquivContSpect}, the continuity of $\ts\ni t\mapsto\|p(A_t)\|\in[0,\infty)$ for specific polynomials up to degree $2$ is equivalent to the continuous behavior of the spectra. These results is based on the geometry of the spectrum as a subset of $\CM$. More precisely, the spectrum is contained in $\RM$ which is a one-dimensional hyperplane of $\CM$. The question of an analog result of Theorem~\ref{Chap3-Theo-P2ContEquivContSpect} was raised to the author by {\sc J. Fillman}.

\medskip

\begin{figure}[htb]
\centering
\includegraphics[scale=1.78]{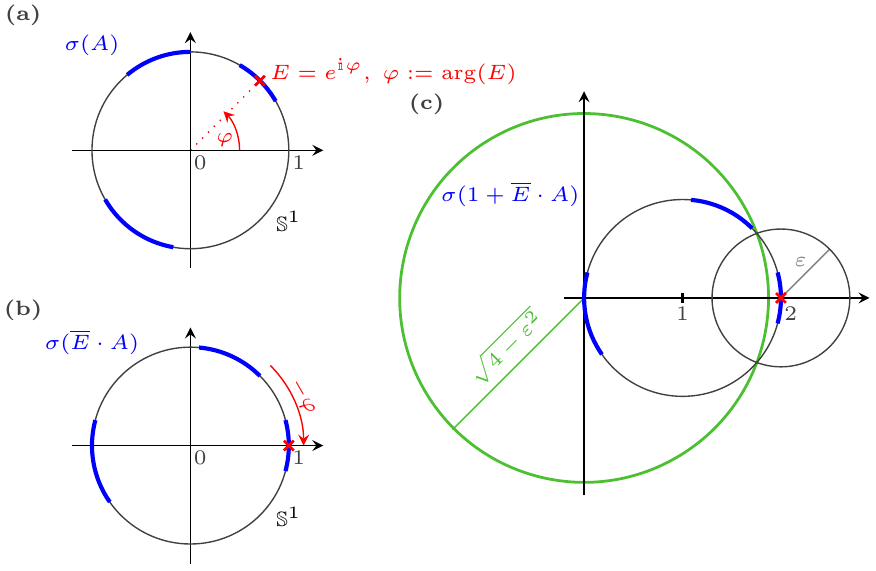}
\caption{Illustration of the proof of Lemma~\ref{Chap3-Lem-PresenceSpectrumUnitary}. Let $A\in\Ll(\hs)$ be a unitary operator with spectrum $\sigma(A)\subseteq\SM^1$, c.f. (a). Each $E\in\SM^1$ is equal to $e^{\mathbbm{i}\varphi}$ for one $\varphi\in[0,2\pi)$ (polar coordinates). Then the multiplication by its complex conjugate is a rotation of the set $\sigma(A)$ by the angle $-\varphi$, c.f. (b). Adding $1$ is a shift of the spectrum $\sigma(A)$ by $+1$, c.f. (c). Then $\|1+\overline{E}\cdot A\|=2$ holds if and only if $E\in\sigma(A)$.}
\label{Chap3-Fig-SpectrumUnitary}
\end{figure}

A similar observation appears in case of unitary operators. In detail, the spectrum of a unitary operator $A\in\Ll(\hs)$, i.e., $A$ satisfies $A^\ast A= AA^\ast =I$ where $I$ is the identity operator $\hs$, is contained in the sphere $\SM^1:=\{z\in\CM\;|\; |z|=1\}$, see e.g. \cite[Proposition~1.3.9.]{Dixmier77}. The key observation is, like in Lemma~\ref{Chap3-Lem-PresenceSpectrum}, the link between the norms of polynomials of the operator and the spectrum. The idea of the proof is illustrated in Figure~\ref{Chap3-Fig-SpectrumUnitary}.

\begin{lemma}
\emph{(\cite[Lemma~2]{BeBe16})}
\label{Chap3-Lem-PresenceSpectrumUnitary}
Let $A\in\Ll(\hs)$ be a unitary operator on a Hilbert space $\hs$. Consider an $E\in\SM^1$. Then the following assertions hold for every $r<2$.
\begin{itemize}
\item[(a)] The inequality $\|1+\overline{E}\cdot A\|\leq \sqrt{4-r^2}$ holds if and only if $B_r(E)\cap\sigma(A)=\emptyset$.
\item[(b)] The inequality $\|1+\overline{E}\cdot A\|> \sqrt{4-r^2}$ holds if and only if $B_r(E)\cap\sigma(A)\neq\emptyset$.
\end{itemize}
\end{lemma}

\begin{proof}
Since assertion (a) and (b) are equivalent the proof of (b) is provided, only. The operator $A$ is unitary implying that the spectrum $\sigma(A)\subseteq\SM^1$ is compact and that $A$ is a bounded operator. Note that the operator norm is bounded by $1$. Let $E\in\SM^1$, then, the equalities 
\begin{align*}
\|1+\overline{E}\cdot A\|\; 
	&\;\; = \;\; \max_{\lambda\in\sigma(A)} |1+ \overline{E}\cdot\lambda|\\
		&\;\; = \;\; \max_{\lambda\in\sigma(A)} \sqrt{(1+\overline{E}\cdot\lambda)\cdot (1+E\cdot \overline{\lambda})}\\
			&\underset{|\lambda|=1}{\overset{|E|=1}{=}} \; \max_{\lambda\in\sigma(A)} \sqrt{2+\overline{E}\lambda + E\cdot\overline{\lambda}}\\
				&\;\; = \;\; \max_{\lambda\in\sigma(A)} \sqrt{4-|E-\lambda|^2}
\end{align*}
hold where the last equality follows by a short computation 
$$
|E-\lambda|^2 \;
	= \; (\overline{E}-\overline{\lambda})\cdot (E-\lambda)\;
		= \; 2+\overline{E}\lambda + E\cdot\overline{\lambda}\, .
$$
Now the desired equivalence is checked as follows. If the intersection $B_r(E)\cap\sigma(A)$ is non-empty, there is a $\lambda\in\sigma(A)$ such that $|E-\lambda|<r$. Thus, the estimate $\|1+\overline{E}\cdot A\|>\sqrt{4-r^2}$ follows by the previous considerations.

\vspace{.1cm}

On the other hand, whenever $\|1+\overline{E}\cdot A\|>\sqrt{4-r^2}$ holds, then $|E-\lambda|^2<r^2$ is derived for at least one $\lambda\in\sigma(A)$. Consequently, $|E-\lambda|<r$ is deduced and so $\lambda\in B_r(E)\cap\sigma(A)$ follows.
\end{proof}

\begin{remark}
\label{Chap3-Rem-PresenceSpectrumUnitary}
As pointed out in Remark~\ref{Chap3-Rem-PolyP2ContEquivContSpect}, the value $1$ in the polynomial $p(z)=1+\overline{E}\cdot z$ plays the same role as the parameter $m$ for self-adjoint fields of operators.
\end{remark}

With this at hand, the following assertion is proved in analogy to Theorem~\ref{Chap3-Theo-P2ContEquivContSpect}.

\begin{theorem}
\emph{(\cite[Theorem~6]{BeBe16})}
\label{Chap3-Theo-CharContSpectUnitary}
Let $\ts$ be a topological space and $(A_t)_{t\in\ts}$ be a unitary field of operators. Then the following assertions are equivalent.
\begin{itemize}
\item[(i)] The map $\Sigma:\ts\to\ks(\CM)\,,\; t\mapsto\sigma(A_t)\,,$ is continuous with respect to the Hausdorff metric on $\ks(\CM)$.
\item[(ii)] For all $E\in \SM^1$, the map $\ts\ni t\mapsto\|1+ E\cdot A_t\|\in[0,\infty)$ is continuous.
\end{itemize}
\end{theorem}

\begin{proof}
(i)$\Rightarrow$(ii): This follows by Lemma~\ref{Chap3-Lem-SpecContNormPhiCont}.

\vspace{.1cm}

(ii)$\Rightarrow$(i): Let $t_0\in\ts$. According to Theorem~\ref{Chap3-Theo-VietFellHausMetricEquiv}, it suffices to prove the Vietoris-continuity of the map $\ts\ni t\mapsto\sigma(A_t)\in\ks(\CM)$. Consider a closed $F\in\cs(\CM)$ and a finite family $\os$ of open subsets of $\CM$ such that $\sigma(A_{t_0})\in\us(F,\os)$. The proof is organized as follows. First, (a) the existence of an open neighborhood $U_F$ of $t_0$ is verified such that $\sigma(A_t)\cap F=\emptyset$ for all $t\in U_F$. Secondly, (b) for each $O\in\os$, the existence of an open neighborhood $U_O$ of $t_0$ is proven such that $\sigma(A_t)\cap O\neq\emptyset$ for $t\in U_O$. Then the finite intersection of the open sets $U_F$ and $U_O,\; O\in\os,$ defines an open neighborhood containing $t_0$ such that $\sigma(A_t)\in\us(F,\os)$ for all $t\in \bigcap_{O\in\os} U_O\cap U_F$.

\vspace{.1cm}

(a) As the operators $A_t,\; t\in\ts,$ are unitary, their spectra is contained in $\SM^1$. Thus, the intersection $\big(\CM\setminus(F\cap \SM^1)\big)\cap \sigma(A_t)$ is empty for all $t\in\ts$. Consequently, it suffices to find a neighborhood $U_F$ of $t_0$ such that $\sigma(A_t)\cap\big(F\cap \SM^1\big)$ is empty for $t\in U_F$. For each $x\in F\cap\SM^1$, there exists a radius $r(x)>0$ such that the open ball $B_{r(x)}(x)$ in $\CM$ does not intersect $\sigma(A_{t_0})$. By compactness of $F\cap\SM^1$, there are $x_1,\ldots,x_l\in F\cap\SM^1$ with radii $r_k:=r(x_k)<2$ for $1\leq k\leq l$ such that 
$$
F\cap\SM^1\subseteq \bigcup_{k=1}^l B_{r_k/2}(x_k)\,
	,\qquad 
	B_{r_k}(x_k)\cap\sigma(A_{t_0})=\emptyset\, .
$$
Fix one $1\leq k\leq l$. Due to construction, $x_k\in\SM^1$ holds. Thus, using Lemma~\ref{Chap3-Lem-PresenceSpectrumUnitary}~(a), the condition $B_{r_k}(x_k)\cap\sigma(A_{t_0})=\emptyset$ is equivalent to $\|1+\overline{x_k}\cdot A_{t_0}\|\leq \sqrt{4-r_k^2}$. Then the continuity of the norms $\|1+\overline{x_k}\cdot A_{t_0}\|$ in $t\in\ts$, required in (ii), leads to the existence of an open neighborhood $U_k$ of $t_0$ such that $\|1+\overline{x_k} \cdot A_t\|\leq \sqrt{4-r_k^2/4}$ holds for $t\in U_k$. Hence, the intersection $B_{r_k/2}(x_k)\cap\sigma(A_t)$ is empty for each $t\in U_k$, c.f. Lemma~\ref{Chap3-Lem-PresenceSpectrumUnitary}~(a). Since $B_{r_k/2}(x_k),\; 1\leq k\leq l,$ is a covering of $F\cap\SM^1$, the intersection $\sigma(A_t)\cap F$ is empty for all $t\in U_F:=\bigcap_{k=1}^l U_k$. Clearly, $U_F$ defines an open neighborhood of $t_0$.

\vspace{.1cm}

(b) Let $O\in\os$. The intersection $O\cap \sigma(A_{t_0})$ is non-empty. For $x\in O\cap \sigma(A_{t_0})$, there is an $r:=r(x)>0$ such that $B_r(x)\subseteq O$ since $O$ is open in $\CM$. Lemma~\ref{Chap3-Lem-PresenceSpectrum}~(b) applies as $x\in\SM^1$, i.e., the estimate $\|1+\overline{x}\cdot A_{t_0}\|>\sqrt{4-r^2/4}$ follows. Then the map $\ts\ni t\mapsto\|1+\overline{x}\cdot A_t\|\in[0,\infty)$ is continuous by (ii). Thus, there is an open neighborhood $U_O$ of $t_0$ in $\ts$ such that the inequality $\|1+\overline{x}\cdot A_t\|>\sqrt{4-r^2}$ is derived for $t\in U_O$. Consequently, the intersection $B_{r}(x)\cap \sigma(A_t)\subseteq O\cap \sigma(A_t)$ is non-empty for $t\in U_O$ using Lemma~\ref{Chap3-Lem-PresenceSpectrum}~(b). 

\vspace{.1cm}

Thanks to construction the intersection $U:=\bigcap_{O\in\os} U_O\cap U_F$ defines an open neighborhood of $t_0$ and the spectrum of $A_t$ satisfies $\sigma(A_t)\in\us(F,\os)$ for all $t\in U$.
\end{proof}

\medskip

Due to the one-dimensional nature of the sphere $\SM^1$, the notion of gaps of closed subsets of $\SM^1$ can be defined in view of Definition~\ref{Chap3-Def-Gap}. Then the concept of continuous boundaries is extendable to a family of closed subsets $K_t\subseteq\SM^1\,,\; t\in\ts\,,$ indexed by a topological space $\ts$. This will lead to a similar assertion like in Theorem~\ref{Chap3-Theo-BoundVietContCompVers} for a family of Vietoris-continuous, closed subsets $K_t\subseteq\SM^1\,,\; t\in\ts$. More precisely, the sets $(K_t)_{t\in\ts}$ are Vietoris-continuous if and only if the boundaries of $(K_t)_{t\in\ts}$ are continuous.

\begin{remark}
\label{Chap3-Rem-OneParaUnitary}
(i) In accordance with Remark~\ref{Chap3-Rem-PolyP2ContEquivContSpect}, the continuous behavior of the spectra $\big(\sigma(A_t)\big)_{t\in\ts}$ of a field of unitary operators $(A_t)_{t\in\ts}$ only depends on one parameter $E\in\SM^1$. In this case, the analog of $m$ is fixed by $1$ since unitary operators have norm $1$.

\vspace{.1cm}

(ii) Like in the self-adjoint case (Remark~\ref{Chap3-Rem-P2ContEquivContSpect}), Theorem~\ref{Chap3-Theo-CharContSpectUnitary} also follows from Corollary~\ref{Chap3-Cor-BoundedCompactXContIffPsi_xCont} since the equations
$$
\|1 + E\cdot A_t\| \;
	= \; \sqrt{
		4-\big(
			\inf\big\{|E-\lambda|\;\big|\; \lambda\in\sigma(A_t)\big\}
		\big)^2}\,,
		\qquad
		E\in\SM^1\,,\; t\in\ts\,,
$$
are derived by the proof of Lemma~\ref{Chap3-Lem-PresenceSpectrumUnitary}.
\end{remark}

Recall that $\Po$ denotes the set of all complex-valued polynomials $z\in\CM\mapsto p(z,\overline{z})$ in $z$ and $\overline{z}$. Then $p(A,A^\ast)$ defines a bounded normal operator for all normal operators $A\in\Ll(\hs)$ on a Hilbert space $\hs$ by the functional calculus.

\begin{corollary}
\label{Chap3-Cor-AllPolNormContUnitary}
Let $\ts$ be a topological space and $(A_t)_{t\in\ts}$ be a field of unitary operators. Then the following assertions are equivalent.
\begin{itemize}
\item[(i)] For all $E\in \SM^1$, the map $\ts\ni t\mapsto\|1+E\cdot A_t\|\in[0,\infty)$ is continuous.
\item[(ii)] The maps $\ts\ni t\mapsto\|p(A_t,A_t^\ast)\|\in[0,\infty)\,,\; p\in\Po,$ are continuous.
\end{itemize}
\end{corollary}

\begin{proof}
The implication (i) to (ii) follows by Theorem~\ref{Chap3-Theo-CharContSpectUnitary} and Lemma~\ref{Chap3-Lem-SpecContNormPhiCont} whereas the converse is clear by definition.
\end{proof}

\begin{remark}
\label{Chap3-Rem-AllPolNormContUnitary}
Note that  the desired continuity in Corollary~\ref{Chap3-Cor-AllPolNormContUnitary} means that $(A_t)_{t\in\ts}$ is \po-continuous, c.f. Definition~\ref{Chap3-Def-(P)-Continuity} and Theorem~\ref{Chap3-Theo-PContEquivContSpectNormal} below.
\end{remark}

\begin{corollary}
\label{Chap3-Cor-CharContSpectFunctOfUnitary}
Let $\ts$ be a topological space and $(A_t)_{t\in\ts}$ be a field of operators. Suppose that there exist a field of unitary operators $(B_t)_{t\in\ts}$ and a continuous, injective function $\phi:X\subseteq\CM\to\CM$ such that $X$ is closed with $\SM^1\subseteq X$ and $A_t=\phi(B_t)$ holds for all $t\in\ts$. Then $(A_t)_{t\in\ts}$ is a normal bounded field and the following assertions are equivalent.
\begin{itemize}
\item[(i)] The family of spectra $(\sigma(A_t))_{t\in\ts}$ is continuous with respect to the Hausdorff metric.
\item[(ii)] For $E\in \SM^1$, the map $\ts\ni t\mapsto\|1+E\cdot B_t\|\in[0,\infty)$ is continuous.
\end{itemize}
\end{corollary}

\begin{proof}
Since $(B_t)_{t\in\ts}$ is unitary, $(A_t)_{t\in\ts}$ defines a bounded normal field of operators by the functional calculus. The proof is left to the reader as it is similar to the one of Corollary~\ref{Chap3-Cor-CharContSpectFunctOfSelfAdj}. The closed set $K$ in this proof is equal to $\SM^1$ since the inclusion $\sigma(B_t)\subseteq\SM^1$ holds for $t\in\ts$.
\end{proof}

\begin{remark}
\label{Chap3-Rem-S1topology}
As explained in Section~\ref{Chap3-Sect-RelP2ContOthTop}, the \pt-continuity is related to a topology on $\Ll(\hs)$ that is called \pt-topology. The \pt-topology is the coarsest topology on $\Ll(\hs)$ so that $\Sigma:\Ll(\hs)\to\ks(\CM)\,,\; A\mapsto\sigma(A)\,,$ is continuous if restricted to the self-adjoint operators $\mathfrak{S}\subseteq\Ll(\hs)$. Similarly, there exists a topology on $\Ll(\hs)$ so that $\Sigma$ is continuous if $\Sigma$ is restricted to the unitary operators $\mathfrak{U}\subseteq\Ll(\hs)$. Specifically, a neighborhood basis for $A_0\in\Ll(\hs)$ is given by
$$
\us_{E,\varepsilon}(A_0) \;
	:= \; \Big\{
		A\in\Ll(\hs) \;\Big|\; \big| \|1+E\cdot A_0\|-\|1+E\cdot A\| \big|<\varepsilon
	\Big\}\,,
	\qquad
	E\in\SM^1\,,\; \varepsilon>0\,.
$$
The topology on $\Ll(\hs)$ defined by this base is called {\em $\SM^1$-topology}. Note that the assertions of Proposition~\ref{Chap3-Prop-NormImplP2}, Proposition~\ref{Chap3-Prop-StrongNotImplP2} and Proposition~\ref{Chap3-Prop-P2NotImplStrong} also hold if the \pt-topology is replaced by the $\SM^1$-topology and unitary operators are considered instead of self-adjoint operators.
\end{remark}

\section{Continuous behavior of the spectra for normal operators}
\label{Chap3-Sect-ContFieldCAlg}

The continuity of a family of spectra is sometimes proven by using continuous fields of Hilbert spaces and $C^\ast$-algebras. The concept of continuous fields of algebras goes back to the work of {\sc Kaplansky} \cite{Kap49,Kap51}, {\sc Fell} \cite{Fel60,Fel61}, {\sc Tomiyama, Takesaki} \cite{ToTa61,Tom62} and {\sc Dixmier-Douady} \cite{DiDo63}. In contrast to fiber bundles, a continuous field of Hilbert spaces may have pairwise non-isomorphic fibers. The field of irrational rotation algebras over the interval $[0,1]$ is a typical example for a continuous field of $C^\ast$-algebras, c.f. \cite{Rie81,Bel94}. It turns out that the Almost-Mathieu operator is an element of this algebra which implies the continuity of the spectra, c.f. Example~\ref{Chap3-Ex-AlmostMathieu}.

\medskip

One of the most powerful consequences of the theory of continuous fields of $C^\ast$-algebras is that a normal continuous vector field (of operators) has spectra varying continuously with the parameter. In the following this result is discussed in more detail. Using, additionally, Lemma~\ref{Chap3-Lem-SpecContNormPhiCont} a characterization of the continuity of the spectra of a bounded normal field of operators is presented and proven. This result is based on Urysohn's Lemma, the functional calculus and the fact that polynomials defined on a compact set $K\subseteq\CM$ are uniformly dense in the space of continuous functions on $K$. For convenience of the reader, a short introduction to the concept of continuous fields of $C^\ast$-algebras is provided. For a more detailed discussion on fields of algebras the reader is referred to the monographs by {\sc Naimark} \cite{Naimark68} (translation \cite{Neumark90}), {\sc Dixmier} \cite{Dixmier69} (translation \cite{Dixmier77}) and {\sc Doran-Fell} \cite{DoranFell88}. 

\medskip

\label{Page-BanachAlgebera}
In the following, a short review of the notion of a $C^\ast$-algebra is presented. For more details on $C^\ast$-algebras we refer the reader to \cite{Rickart60,Dixmier69,Dixmier77,Dixmier81,Palmer94,Dixmier96,Murphy90,Palmer01}. A vector space $\AG$ over $\CM$ equipped with a multiplication $\star:\AG\times\AG\to\AG$ is called an {\em algebra} if the equations
\begin{description}
\item[(1)] $\;\;\;\fz\star (\gz\star \hz) \;\;\, = \; (\fz\star \gz)\star \hz \, , \hfill\textbf{(associativity)}$
\item[(2)] $\begin{array}{rl}
	(\fz + \gz)\star \hz \; 
		&= \; \fz\star \hz + \gz\star \hz \, ,\\
	\fz\star (\gz + \hz) \; 
		&= \; \fz\star \gz + \fz\star \hz \, ,\\
	\lambda (\fz\star \gz) \; 
		&= \; (\lambda \fz)\star \gz \; 
			= \; \fz\star (\lambda \gz) \, ,
\end{array}
\hfill \textbf{(distributivity)}
$
\end{description}
hold for all $\fz,\gz,\hz\in\AG$ and $\lambda\in\CM$. If, additionally, the equality $\fz\star\gz=\gz\star\fz$ is valid for all $\fz,\gz\in\AG$, the algebra $\AG$ is called {\em commutative}. Otherwise, $\AG$ is called {\em non-commutative}. An algebra $\AG$ is called {\em unital} whenever there exists an element $\mathpzc{1}\in\AG$ such that $\mathpzc{1} \star \fz = \fz\star \mathpzc{1} = \fz$ for all $\fz\in\AG$. Note that such an element $\mathpzc{1}\in\AG$ is unique if it exists and is called the {\em unit} or {\em identity}. An algebra $\AG$ equipped with a norm $\|\cdot\|:\AG\to [0,\infty)$ is called {\em normed algebra} if the inequality $\|\fz\star \gz\|\leq \|\fz\|\cdot\|\gz\|$ holds for all $\fz,\gz\in\AG$. Whenever the algebra $\AG$ is unital, it is additionally required that $\|\mathpzc{1}\|=1$. If a normed algebra $(\AG,\|\cdot\|)$ is complete it is called a {\em Banach algebra}. A map $^\ast:\AG\to\AG$ on an algebra $\AG$ is an {\em involution} if the equations
\begin{align*}
(\fz+\lambda \gz)^\ast \; 
	&= \; \fz^\ast + \overline{\lambda} \gz^\ast \, ,\\
(\fz\star \gz)^\ast \; 
	&= \; \gz^\ast\star \fz^\ast \, ,\\
(\fz^\ast)^\ast \; 
	&= \; \fz \, ,
\end{align*}
are valid for all $\fz,\gz\in\AG$ and $\lambda\in\CM$. Note that $\overline{\lambda}$ denotes the complex conjugation in $\CM$. Then an algebra $\AG$ with an involution is called an {\em $\ast$-algebra} or an {\em involutive algebra}. Given an $\ast$-algebra $\AG$, an element $\fz\in\AG$ is called 
\begin{itemize}
\item {\em normal} if $\fz\star\fz^\ast=\fz^\ast\star\fz$ holds;
\item {\em self-adjoint} if $\fz^\ast=\fz$ holds.
\end{itemize}

\medskip

Let $\AG$ be a unital Banach algebra. Then $\fz\in\AG$ is called {\em invertible} if there exists a $\gz\in\AG$ such that $\fz\star\gz=\gz\star\fz=\mathpzc{1}$. The inverse $\gz\in\AG$ is unique and is denoted by $\fz^{-1}$. The {\em spectrum} of an element $\fz\in\AG$ is defined by 
$$
\sigma(\fz)\; 
	:= \; 
		\left\{
			\lambda\in\CM \; |\;
				\lambda\, \mathpzc{1}-\fz \text{ is not invertible}
		\right\} \, .
$$

It is well-known, for a $C^\ast$-algebra $\AG$, that the spectrum $\sigma(\fz)$ is real if $\fz\in\AG$ is a self-adjoint element, c.f. \cite[Proposition~1.3.9.]{Dixmier77}.

\begin{definition}[$C^\ast$-algebra]
\label{Chap3-Def-C-algebra}
Let $(\AG,\|\cdot\|)$ be a Banach algebra with an isometric involution, i.e., $\|\fz^\ast\|=\|\fz\|$  holds for all $\fz\in\AG$. Then $\AG$ is called a {\em $C^\ast$-algebra} if all $f\in\AG$ satisfy
$$
\|\fz\|^2 \; 
	= \; \|\fz\star \fz^\ast\| \, .
$$
\end{definition}

Note that the $C^\ast$-norm $\|\cdot\|$ of a $C^\ast$-algebra $\AG$ is unique, i.e., if there exists another norm $|||\cdot|||$ on $\AG$ satisfying the identities $|||\fz|||^2=|||\fz\star\fz^\ast|||$ for $\fz\in\AG$ then the norms $|||\cdot|||$ and $\|\cdot\|$ coincide, c.f. \cite[Corollary~2.1.2]{Murphy90}.

\begin{example}
\label{Chap3-Ex-LinBounOpCalgebra}
For a Hilbert space $\hs$, the set of all bounded linear operators on $\hs$ is denoted by $\Ll(\hs)$. An involution on $\Ll(\hs)$ is given by the adjoint of an operator. Furthermore, a multiplication is given by the composition of operators. Together with the operator norm, $\Ll(\hs)$ is a $C^\ast$-algebra. This $C^\ast$-algebra is unital where the unit is given by the identity operator on $\hs$. Furthermore, it is a non-commutative algebra if the dimension of the Hilbert space $\hs$ is greater than or equal to $2$. Note that the definition of spectrum in a $C^\ast$-algebra coincides with the usual notion of spectrum of an operator. Every closed subalgebra of $\Ll(\hs)$ defines a $C^\ast$-algebra as well. For instance, the closure of the algebraic span of an elements of $\Ll(\hs)$ defines a sub-$C^\ast$-algebra. According to the Theorem of Gelfand-Naimark-Segal, each unital $C^\ast$-algebra is faithfully represented on $\Ll(\hs)$ for a convenient Hilbert space $\hs$, c.f. \cite{GeNe43,Seg47}.
\end{example}

\begin{example}
\label{Chap3-Ex-ContFunctCalgebra}
Let $X$ be a topological space which is not necessarily compact. Then the set $\Cc_0(X)$ of all continuous functions that vanish at infinity defines a commutative $C^\ast$-algebra with the pointwise multiplication and the complex conjugation as involution. The $C^\ast$-norm is given by the uniform norm $\|\fz\|_\infty:=\sup_{x\in X}|\fz(x)|,\; \fz\in\Cc_0(X)$. This $C^\ast$-algebra is unital if and only if $X$ is compact. In contrast to that the set of all continuous bounded functions $\Cc_b(X)$ is unital as the constant function $1$ is continuous and bounded. For $\fz$ an element of one of these $C^\ast$-algebras, the spectrum is given by the image of the function $\fz$.
\end{example}

Let $\CG_t,\; t\in\ts,$ be a family of $C^\ast$-algebras indexed by a topological space $\ts$. A {\em vector field} is an element $\fz:= (\fz_t)_{t\in\ts}$ of the product space $\CG:=\prod_{t\in\ts} \CG_t$. The algebraic operations on $\CG$ are defined pointwise, i.e.,
\begin{gather*}
\begin{aligned} 
(\fz_t)^\ast_{t\in\ts} \;
	&:= \; (\fz_t^\ast)_{t\in\ts}
	 \, ,\qquad
	&(\fz_t)_{t\in\ts} + (\gz_t)_{t\in\ts} \; 
		&:= \; (\fz_t+\gz_t)_{t\in\ts}
	 \, ,\\
\lambda\cdot (\fz_t)_{t\in\ts} \;
	&:=\; (\lambda\cdot \fz_t)_{t\in\ts}
	 \, ,\qquad
	&(\fz_t)_{t\in\ts} \star (\gz_t)_{t\in\ts} \; 	
		&:= \; (\fz_t \star \gz_t)_{t\in\ts} \, ,
\end{aligned} 
\end{gather*}
where $(\fz_t)_{t\in\ts},\; (\gz_t)_{t\in\ts}$ are vector fields and $\lambda\in\CM$. Following the lines of \cite[Chapter~10]{Dixmier77}, a continuous field of unital $C^\ast$-algebras is defined. Recall that a $C^\ast$-algebra is called unital if it admits a unit which is also called identity.

\begin{definition}[Continuous fields of unital $C^\ast$-algebras]
\label{Chap3-Def-ContFieldCAlg}
Let $\ts$ be a topological space and $(\CG_t)_{t\in\ts}$ a family of unital $C^\ast$-algebras. For $\Upsilon\subseteq\CG:=\prod_{t\in\ts}\CG_t$, the tuple $((\CG_t)_{t\in\ts},\Upsilon)$ is called a {\em continuous field of unital $C^\ast$-algebras over $\ts$} if the following assertions hold.
\begin{description}
\item[(CFC1)\label{(CFC1)}] The set $\Upsilon$ is a (complex) linear subspace of $\CG$ being closed under multiplication $\star:\CG\times\CG\to\CG$ and involution $\ast:\CG\to\CG$. Additionally, for every $t_0\in\ts$, there exist an open neighborhood $U_0\subseteq\ts$ of $t_0$ and an $\fz\in\Upsilon$ such that $\fz_t=\mathpzc{1}_t$ for $t\in U_0$ where $\mathpzc{1}_t$ is the unit of $\CG_t$.
\item[(CFC2)\label{(CFC2)}] The set $\Upsilon_t:=\{\fz_t \;|\; \fz\in \Upsilon \}$ is dense in $\CG_t$ for each $t\in\ts$.
\item[(CFC3)\label{(CFC3)}] The map $\ts\ni t \mapsto \|\fz_t \| \in [0,+\infty)$ is continuous for each $\fz\in \Upsilon$.
\item[(CFC4)\label{(CFC4)}] A vector field $\fz$ is an element of $\Upsilon$ if and only if, for each $t_0 \in\ts$ and every $\varepsilon>0$, there is a $\gz\in \Upsilon$ and a neighborhood $U_0$ of $t_0$ such that $\|\fz_t - \gz_t \| \leq \varepsilon$ for $t\in U_0$.
\end{description}
An element $\fz\in \Upsilon$ is called a {\em continuous vector field}.
\end{definition}

Combining \nameref{(CFC1)} and \nameref{(CFC4)}, the vector field $(\mathpzc{1}_t)_{t\in\ts}$ is contained in $\Upsilon$ where $\mathpzc{1}_t$ is the unit of $\CG_t$ for $t\in\ts$. Note that Condition~\nameref{(CFC3)} is obviously satisfied for this vector field as $\|\mathpzc{1}_t\|=1$ holds for $t\in\ts$.

\begin{remark}
\label{Chap3-Rem-ContFieldCAlg}
In \cite[Definition~10.3.1]{Dixmier77}, it is not required that the $C^\ast$-algebras $\CG_t\, ,\; t\in\ts\, ,$ are unital, c.f. \cite[10.3.5]{Dixmier77} as well. Thus, also the second part of \nameref{(CFC1)} does not appear in the definition of {\sc Dixmier}, namely that there exist, for every $t_0\in\ts$, an open neighborhood $U_0\subseteq\ts$ of $t_0$ and an $\fz\in\Upsilon$ such that $\fz_t=\mathpzc{1}_t$ for $t\in U_0$ where $\mathpzc{1}_t$ is the unit of $\CG_t$. This additional assumption is made so that the polynomials in Lemma~\ref{Chap3-Lem-PolContVectField} are allowed to have a constant part. In view of Remark~\ref{Chap3-Rem-ContSpectrNonUnitalAlg} and Remark~\ref{Chap3-Rem-NecessIdentityContSpect}, this is necessary for our purposes.
\end{remark}

Let $\Lambda\subseteq \prod_{t\in\ts} \CG_t$ be such that the properties \nameref{(CFC1)}-\nameref{(CFC3)} are satisfied. In accordance with the following result, such a subset $\Lambda$ is called {\em generating family}.

\begin{proposition}[\cite{Dixmier77}, Proposition~10.2.3]
\label{Chap3-Prop-GenerFamiContField}
Let $\ts$ be a topological space and $(\CG_t)_{t\in\ts}$ be a family of unital $C^\ast$-algebras. Consider a subset $\Lambda\subseteq \prod_{t\in\ts} \CG_t$ satisfying \nameref{(CFC1)}-\nameref{(CFC3)}, namely $\Lambda$ is a generating family. Then there exists a unique subset $\Upsilon\subseteq \prod_{t\in\ts} \CG_t$ such that $\Lambda\subseteq\Upsilon$ and $((\CG_t)_{t\in\ts},\Upsilon)$ is a continuous field of unital $C^\ast$-algebras.
\end{proposition}

The following observation holds even for non-unital $C^\ast$-algebras if the polynomials do not have a constant term.

\begin{lemma}
\label{Chap3-Lem-PolContVectField}
Let $((\CG_t)_{t\in\ts},\Upsilon)$ be a continuous field of unital $C^\ast$-algebras over $\ts$ and $\fz:=(\fz_t)_{t\in\ts}\in\Upsilon$ a continuous vector field. Then, for every complex valued polynomial $z\in\CM\mapsto p(z,\overline{z})\in\CM$ in $z$ and $\overline{z}$, the vector field $\big(p(\fz_t,\fz_t^\ast)\big)_{t\in\ts}$ is continuous. 
\end{lemma}

\begin{proof}
The vector field $(\mathpzc{1}_t)_{t\in\ts}$ is continuous by \nameref{(CFC4)} and the second part of \nameref{(CFC1)}. Now, the desired assertion follows by \nameref{(CFC1)} since each vector field $(p(\fz_t,\fz_t^\ast))_{t\in\ts}$, for a polynomial $p$, is a finite algebraic combination of the continuous vector fields $(\fz_t)_{t\in\ts}$, $(\fz_t^\ast)_{t\in\ts}$ and the unit vector field $(\mathpzc{1}_t)_{t\in\ts}$.
\end{proof}

\medskip

A version of the following assertion can be found in \cite{Dixmier77}. Since the statement plays a crucial role for the continuous behavior of the spectra, the proof is provided for convenience of the reader.

\begin{proposition}[\cite{Dixmier77}, Proposition~10.3.3]
\label{Chap3-Prop-ContFuncContVectFie}
Let $((\CG_t)_{t\in\ts},\Upsilon)$ be a continuous field of unital $C^\ast$-algebras over $\ts$ and $\fz\in\Upsilon$ a continuous vector field which is normal, i.e.,  $\fz_t\fz_t^\ast=\fz_t^\ast\fz_t$ for all $t\in\ts$. Then the vector field $(\phi(\fz_t))_{t\in\ts}$ belongs to $\Upsilon$ for every continuous $\phi:\CM\to\CM$.
\end{proposition}

\begin{proof}
Since $\fz_t\in\CG_t$ is a normal element for $t\in\ts$, the functional calculus applies, c.f. \cite[Section 1.5]{Dixmier77}. Thus, $\phi(\fz_t)$ defines an element of $\CG_t$. Let $t_0\in\ts$ and $U_0$ be an open neighborhood of $t_0$ such that $m:=\sup_{t\in U_0}\|\fz_t\|<\infty$. The existence of $U_0$ is guaranteed by the continuity of the norms $\|\fz_t\|$ in $t\in\ts$. The complex-valued polynomials defined on $B_m(0):=\{z\in\CM\;|\; |z|\leq m\}$ are dense in the unital $C^\ast$-algebra of continuous, complex-valued functions on $B_m(0)$ equipped with the supremum norm. Thus, for a continuous function $\phi:\CM\to\CM$, there is a sequence of polynomials $z\in\CM\mapsto p_n(z,\overline{z})\in\CM,\; n\in\NM,$ in $z$ and $\overline{z}$ such that $(p_n)_{n\in\NM}$ converges uniformly on $B_m(0)$ to $\phi$. Lemma~\ref{Chap3-Lem-PolContVectField} implies that $\big(p_n(\fz_t,\fz_t^\ast)\big)_{t\in\ts}$ is a continuous vector field. By the uniform convergence, there is, for each $\varepsilon>0$, an $n(\varepsilon)\in\NM$ such that 
$$
\big\|\phi(\fz_t)-p_n(\fz_t,\fz_t^\ast)\big\| \; 
	=\; \sup_{\lambda\in\sigma(\fz_t)} \big|\phi(\lambda)-p_n(\lambda,\overline{\lambda})\big|\;
	\leq \; \big\|\phi-p_n\big\|_\infty
	<\; \varepsilon\,,
	\quad
	t\in U_0\,,\;
	n\geq n(\varepsilon)\,.
$$ 
Since $t_0$ was chosen arbitrarily, Condition~\nameref{(CFC4)} implies that the vector field $\big(\phi(\fz_t)\big)_{t\in\ts}$ is continuous. Hence, $\big(\phi(\fz_t)\big)_{t\in\ts}\in\Upsilon$ follows.
\end{proof}

\medskip

In \cite[Proposition~10.3.3]{Dixmier77} the additional assumption $\phi(0)=0$ is made. In our case, this assumption is not needed as the unit belongs to the continuous vector fields, c.f. Remark~\ref{Chap3-Rem-ContFieldCAlg} and \cite[10.3.5]{Dixmier77}.

\medskip

With this at hand, the abstract result of the continuous variation of the spectra is proven. This result has been re-proven over the years by several authors. First {\sc Kaplanski} \cite[Lemma~3.3]{Kap51} proved it for self-adjoint operators as well as {\sc Naimark} \cite{Naimark68}, c.f. \cite[Chapter V.\S 26.4.V]{Neumark90}. {\sc Dixmier} \cite{Dixmier69} gave a proof for normal elements, c.f. \cite[Proposition 10.3.6]{Dixmier77}. In 1982, {\sc Elliot} \cite{Ell82} applied this result to almost periodic Schr\"odinger operators $H_n,\; n\in\oNM,$ to prove the convergence of the corresponding spectra. On this basis, {\sc Bellissard} \cite{Bel94} applied this technique, among others, to an effective Hamiltonian contained in the rotation algebra. Versions of this general result can be also found in \cite{Man12,BeMa12}. There the outer and inner continuity are nothing but the continuity with respect to the Vietoris-topology. The proof presented here follows the lines of \cite[Proposition~1]{Bel94}.

\begin{theorem}[Continuity of the spectra, \cite{Kap51} et al.]
\label{Chap3-Theo-ContFieldCALgContSpectr}
Let $((\CG_t)_{t\in\ts},\Upsilon)$ be a continuous field of unital $C^\ast$-algebras over $\ts$ and $\fz:=(\fz_t)_{t\in\ts}\in\Upsilon$ a continuous vector field which is normal, i.e.,  $\fz_t\fz_t^\ast=\fz_t^\ast\fz_t$ for all $t\in\ts$. Then the map $\Sigma:\ts\to\ks(\CM)\,,\; t\mapsto\sigma(\fz_t)\,,$ is continuous with respect to the Hausdorff metric where $\sigma(\fz_t)\subseteq\CM$ is the spectrum of $\fz_t$.
\end{theorem}

\begin{figure}[htb]
\centering
\includegraphics[scale=1.45]{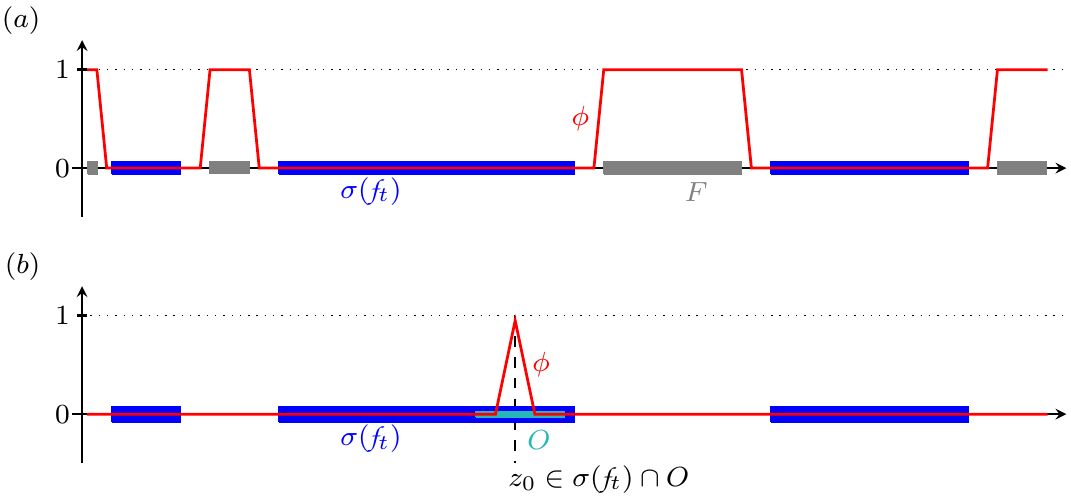}
\caption{Illustration of the proof of Theorem~\ref{Chap3-Theo-ContFieldCALgContSpectr}. For sake of presentation the figure is drawn in one-dimension and not in $\CM$. (a) The function $\phi:\CM\to[0,1]$ satisfies that it vanishes on $\sigma(\fz_t)$ and is equal to $1$ on $F$. (b) The function $\phi:\CM\to[0,1]$ satisfies that $\phi(z_0)=1\, , \; \phi(z)=0$ if $|z-z_0|\geq\varepsilon$ and $\phi(z)<1$ if $0<|z-z_0|<\varepsilon$.}
\label{Chap3-Fig-ThmContSpectNormalOp}
\end{figure}

\begin{proof}
Let $t_0\in\ts$. According to Theorem~\ref{Chap3-Theo-VietFellHausMetricEquiv}, it suffices to prove the Vietoris-continuity. Let $F\in\cs(\CM)$ and $\os$ a finite family of open subset of $\CM$ be such that $\sigma(\fz_{t_0})$ belongs to the open set $\us(F,\os)$. In the following, the functional calculus as well as the Lemma of Urysohn are used to show the existence of an open neighborhood $U$ of $t_0$ such that $\sigma(\fz_t)\in\us(F,\os)$ for $t\in U$.

\vspace{.1cm}

The intersection $\sigma(\fz_{t_0})\cap F$ is empty by assumption. Let $\phi:\CM\to [0,1]$ be continuous such that $\phi(z)=1$ for $z\in F$ and $\phi(z)=0$ for $z\in\sigma(\fz_{t_0})$, c.f. Figure~\ref{Chap3-Fig-ThmContSpectNormalOp}~(a). Such a $\phi$ exists by Urysohn's Lemma that applies for $\CM$. Since $\fz$ is normal, the element $(\phi(\fz_t))_{t\in\ts}$ is a well-defined vector field. By Proposition~\ref{Chap3-Prop-ContFuncContVectFie}, this vector field is continuous, namely $(\phi(\fz))_{t\in\ts}\in\Upsilon$. By construction, $\phi$ vanishes on $\sigma(\fz_{t_0})$ and so $\phi(\fz_{t_0})=0$ holds. By the continuity of the norms $\|\phi(\fz_t)\|$ in $t\in\ts$, there is an open neighborhood $U(F)$ of $t_0$ such that if $t\in U(F)$ then $\|\phi(\fz_t)\|<1/2$. As $\fz_t$ is a normal element of $\CG_t$, it follows that $\|\phi(\fz_t)\|=\sup\big\{ |\phi(z)|\,\big|\, z\in \sigma(\fz_t)\big\}$. Hence, the intersection $\sigma(\fz_t)\cap F$ is empty for $t\in U(F)$ since $\phi(z)=1$ holds for all $z\in F$.

\vspace{.1cm}

For an open $O\in\os$, the intersection $\sigma(\fz_{t_0})\cap O$ is non-empty. Let then $z_0\in \sigma(\fz_{t_0})\cap O$ and $\varepsilon>0$ be small enough so that the open disk $B_\varepsilon(z_0)$ is contained in $O$. Due to construction, the intersection $\sigma(\fz_{t_0})\cap B_\varepsilon(z_0)$ is non-empty. Applying Urysohn's Lemma, there is a continuous function $\phi:\CM\to [0,1]$ such that $\phi(z_0)=1$, $\phi(z)=0$ if $|z-z_0|\geq \varepsilon$ and $\phi(z)<1$ if $0<|z-z_0|< \varepsilon$, c.f. Figure~\ref{Chap3-Fig-ThmContSpectNormalOp}~(b). Proposition~\ref{Chap3-Prop-ContFuncContVectFie} implies that $(\phi(\fz_t))_{t\in\ts}$ is a continuous vector field. Since $\fz_t$ is normal, the identity $\|\phi(\fz_t)\|=\sup\big\{ |\phi(z)|\,\big|\, z\in \sigma(\fz_t)\big\}$ holds. Hence, $\|\phi(\fz_{t_0})\|$ is equal to one. By continuity of the norm $\|\phi(\fz_t)\|$ in $t\in\ts$, there exists an open neighborhood $U(O)$ of $t_0$ such that $\|\phi(\fz_t)\|\geq 1/2$ for $t\in U(O)$. Consequently, there is a point $z_t\in\sigma(\fz_t)$ such that $|\phi(z_t)|\geq 1/2$. Therefore, $|z_t-z_0|\leq \varepsilon$ holds by definition of $\phi$. Hence, $z_t\in  \sigma(\fz_t)\cap B_\varepsilon(z_0)\subseteq\sigma(\fz_t)\cap O$ follows for $t\in U(O)$. Thus, the intersection $\sigma(\fz_t)\cap O$ is non-empty.

\vspace{.1cm}

Altogether, the spectrum $\sigma(\fz_t)$ is contained in $\us(F,\os)$ for all $t\in U:=\bigcap_{O\in\os} U(O)\cap U(F)$  where $U$ is an open neighborhood of $t_0$ since it is a finite intersection of open neighborhoods of $t_0$.
\end{proof}

\begin{remark}
\label{Chap3-Rem-ContSpectrNonUnitalAlg}
The analog of Proposition~\ref{Chap3-Prop-ContFuncContVectFie} in \cite[Proposition~10.3.3]{Dixmier77} is stated with the additional assumption that $\phi(0)=0$. This is necessary due to the fact that the $C^\ast$-algebras are not required to be unital. Thus, the polynomials are not allowed to have a constant part. Even in this case Theorem~\ref{Chap3-Theo-ContFieldCALgContSpectr} is valid. There $0\in\CM$ plays a special role as the spectrum is defined by adjoining a unit, i.e., zero is contained in the spectrum of each element of the $C^\ast$-algebra, c.f. Remark~\ref{Chap3-Rem-NecessIdentityContSpect}. Altogether, the map $\ts\ni t\mapsto\sigma(\fz_t)\cup\{0\}\in\ks(\CM)$ is Vietoris-continuous for every normal continuous vector field $(\fz_t)_{t\in\ts}\in\Upsilon$ if $((\CG_t)_{t\in\ts},\Upsilon)$ is a continuous field of non-unital $C^\ast$-algebras.
\end{remark}

This abstract result immediately leads to the characterization of the continuity of the spectra corresponding to normal operators. The \pt-continuity was defined by controlling the norms of specific real-valued polynomials up to degree two of an operator. The set of all complex-valued polynomials $z\mapsto p(z,\overline{z})$ in $z$ and $\overline{z}$ is denoted by $\Po$. Recall the notion of a field of operators defined in Definition~\ref{Chap3-Def-FieldOfLinOperators}.

\begin{definition}[\po-continuous]
\label{Chap3-Def-(P)-Continuity}
Let $\ts$ be a topological space and $(A_t)_{t\in\ts}$ be a bounded normal field of operators on the field of Hilbert spaces $(\hs_t)_{t\in\ts}$. Then $(A_t)_{t\in\ts}$ is called \po-continuous if all the maps
$$
\Phi_p:\ts\to[0,\infty)\,,\; t\mapsto \|p(A_t,A_t^\ast)\|\,,
	\qquad p\in\Po\,,
$$
are continuous.
\end{definition}

According to Example~\ref{Chap3-Ex-LinBounOpCalgebra}, the $\ast$-algebra of bounded operators $\Ll(\hs)$ on $\hs$ with operator norm, adjoint for the involution and composition for the multiplication defines a $C^\ast$-algebra. With this notion at hand, the continuous behavior of the spectra is characterized.

\begin{theorem}
\label{Chap3-Theo-PContEquivContSpectNormal}
Let $\ts$ be a topological space and $(A_t)_{t\in\ts}$ a bounded normal field of opera\-tors on the field of Hilbert spaces $(\hs_t)_{t\in\ts}$. Then the following assertions are equivalent.
\begin{itemize}
\item[(i)] The map $\Sigma:\ts\to\ks(\CM)\,,\;t\mapsto\sigma(A_t)\,,$ is continuous with respect to the Hausdorff metric on $\ks(\CM)$.
\item[(ii)] The field $(A_t)_{t\in\ts}$ is \po-continuous.
\item[(iii)] Consider the field of $C^\ast$-algebras $(\CG_t)_{t\in\ts}$ defined by $\CG_t:=\CG^\ast(A_t,A_t^\ast,\mathpzc{1}_t)\subseteq\Ll(\hs)$ where $\mathpzc{1}_t:\hs_t\to\hs_t$ is the identity operator for $t\in\ts$. Then the set
$$
\Lambda\; 
	:=\; \big\{
		\big(p(A_t,A_t^\ast)\big)_{t\in\ts} \;|\; p\in\Po
	\big\}
	\subseteq \prod_{t\in\ts}\CG_t \, .
$$
is a generating family, i.e., there exists a unique $\Upsilon\subseteq\prod_{t\in\ts}\CG_t$ such that $\Lambda\subseteq\Upsilon$ and $((\CG_t)_{t\in\ts},\Upsilon)$ defines a continuous field of unital $C^\ast$-algebras.
\end{itemize}
\end{theorem}

\begin{proof}
(i)$\Rightarrow$(ii): This follows from Lemma~\ref{Chap3-Lem-SpecContNormPhiCont}. 

\vspace{.1cm}

(ii)$\Rightarrow$(iii): Define the $C^\ast$-algebra $\CG_t:=\CG^\ast(A_t,A_t^\ast,\mathpzc{1}_t)$ where $\mathpzc{1}_t:\hs_t\to\hs_t$ is the identity operator. Specifically, $\CG_t$ is defined by the closure in $\Ll(\hs_t)$ of the algebra generated by $A_t,\; A_t^\ast$ and $\mathpzc{1}_t$. Since the product of polynomials in $\Po$ is again in $\Po$, the space $\Lambda$ is closed under multiplication. Additionally, it is closed under involution as the adjoint of $p(A_t,A_t^\ast)$ is still a polynomial in $A_t, A_t^\ast$ and $\mathpzc{1}_t$ for $t\in\ts$. Every finite algebraic combinations of $A_t,\; A_t^\ast$ and $\mathpzc{1}_t$ is equal to $p(A_t,A_t^\ast)$ for a specific polynomial $p\in\Po$. Hence, $\Lambda_t:=\{\fz_t\;|\; \fz\in\Lambda\}$ is a dense linear subspace of $\CG_t$ and so the set $\Lambda$ fulfills \nameref{(CFC1)} and \nameref{(CFC2)}. Note that the vector field $(\mathpzc{1}_t)_{t\in\ts}$ of units is contained in $\Lambda$ as $p$ may be chosen constant to $1$. Condition \nameref{(CFC3)} is satisfied by (ii). Consequently, $\Lambda$ is a generating family. Then Proposition~\ref{Chap3-Prop-GenerFamiContField} leads to the existence of a unique $\Upsilon\subseteq\prod_{t\in\ts} \CG_t$ such that $\Lambda\subseteq\Upsilon$ and $\big((\CG_t)_{t\in\ts},\Upsilon\big)$ defines a continuous field of unital $C^\ast$-algebras. 

\vspace{.1cm}

(iii)$\Rightarrow$(i): Clearly, $(A_t)_{t\in\ts}\in\Lambda$ holds as we can choose $p(z,\overline{z}):=z,\; z\in\CM$. Hence, $(A_t)_{t\in\ts}$ is a normal continuous vector field. Then Theorem~\ref{Chap3-Theo-ContFieldCALgContSpectr} implies the desired continuity of the spectra.
\end{proof}

\begin{remark}
\label{Chap3-Rem-NecessIdentityContSpect}
In Theorem~\ref{Chap3-Theo-PContEquivContSpectNormal}, it is necessary to require that the $C^\ast$-algebras $\CG_t:=\CG(A_t,A_t^\ast,\mathpzc{1}_t)\, , \; t\in\ts,$ contain the identity $\mathpzc{1}_t:\hs_t\to\hs_t$ and that the vector field $(\mathpzc{1}_t)_{t\in\ts}$ is continuous. More precisely, let $\AG_t$ be the $C^\ast$-algebra defined by $A_t$ and $A_t^\ast$. Then $\AG_t$ is a sub-$C^\ast$-algebra of $\CG_t$. The spectrum $\sigma_{\AG_t}(A_t)$ in $\AG_t$ can differ with the spectrum $\sigma_{\CG_t}(A_t)$ in $\CG_t$. Actually, the equation $\sigma_{\AG_t}(A_t)\cup\{0\}=\sigma_{\CG_t}(A_t)\cup\{0\}$ only holds, c.f. \cite[Proposition 1.3.10]{Dixmier77}. Thus, whenever the identity operators $\mathpzc{1}_t,\; t\in\ts\,,$ are not used for the definition of $\CG_t$, the spectra only behave continuous if the point $0$ is added to all the spectra $\sigma(A_t)\,,\; t\in\ts$, c.f. Remark~\ref{Chap3-Rem-ContSpectrNonUnitalAlg}.
\end{remark}

\begin{remark}
\label{Chap3-Rem-WhyCAlgebrasCOntBehav}
The proof of Proposition~\ref{Chap3-Prop-ContFuncContVectFie} and Theorem~\ref{Chap3-Theo-ContFieldCALgContSpectr} relies on Urysohn's Lemma, the functional calculus and the denseness of complex-valued polynomials in the space of complex-valued continuous functions $\phi:K\to\CM$ for $K\subseteq\CM$ bounded. This leads to the continuous behavior of the spectra by the continuity of the norms of polynomials of the operators. More precisely, the proof of Theorem~\ref{Chap3-Theo-PContEquivContSpectNormal} can be given without mentioning continuous fields of $C^\ast$-algebras. However, the theory of continuous fields of $C^\ast$-algebras provides a tool to verify the \po-continuity of a field of bounded normal operators whenever it arises by a groupoid $C^\ast$-algebra. This is discussed in Chapter~\ref{Chap4-ToolContBehavSpectr}.
\end{remark}

The spectrum of an operator is independent of the $C^\ast$-algebra that contains the operator if $0$ is excluded, see e.g. \cite[Proposition~1.3.10]{Dixmier77}. On the other hand, the spectra can contain $0$ in one $C^\ast$-algebra while it is excluded in another $C^\ast$-algebra which is shown in the following example.

\begin{example}
\label{Chap3-Ex-ProjContSpect}
Let $\hs$ be a Hilbert space and $P$ a non-trivial projection on $\hs$ and $\mathpzc{1}:\hs\to\hs$ the identity operator. Then the $C^\ast$-algebra $\AG(P)$ is nothing but the space $\{\lambda P\;|\; \lambda\in\CM\}$ as $P^2=P$ and $P^\ast=P$. Hence, $\AG(P)$ is isomorphic to $\CM$. This $C^\ast$-algebra is unital with unit $P$ since $P^2=P$. Thus, the spectrum $\sigma_{\AG(P)}(P)$ of $P$ in the $C^\ast$-algebra $\AG(P)$ is equal to $\{1\}$. On the other hand, the spectrum of $P$ in the $C^\ast$-algebra $\Ll(\hs)$ is equal to $\{0,1\}$.
\end{example}

This leads to the following observation in view of Example~\ref{Chap3-Ex-StrongConvNotConvSpect} and Remark~\ref{Chap3-Rem-NecessIdentityContSpect}. Recall that the notation $p(A)$ is used instead of $p(A,A^\ast)$ for self-adjoint $A\in\Ll(\hs)$.

\begin{example}
\label{Chap3-Ex-ContFieldProj}
Let the topological space $\ts:=\oNM$ be the one-point compactification of $\NM$ and $A_n\in\Ll(\ell^2(\ZM)),\; n\in\NM,$ projections and $A_\infty=\mathpzc{1}$ is the identity operator on $\ell^2(\ZM)$. Define for $n\in\NM$ the $C^\ast$-algebras $\CG_n:=\CG^\ast(A_n,I)$ and $\AG_n:=\CG^\ast(A_n)$. Let $p(x):=q(x)+p_0$ be a polynomial such that $q$ is a polynomial without constant term, i.e., $q(0)=0$, and $p_0\in\CM\setminus\{0\}$.  By the functional calculus of $C^\ast$-algebras, the element $p(A_n)$ in the $C^\ast$-algebra $\AG_n$ is nothing but $q(A_n)+p_0\cdot A_n$ for $n\in\NM$ as $A_n$ is the unit in $\AG_n$. On the other hand, in the $C^\ast$-algebra $\CG_n$ the element $p(A_n)$ is equal to $q(A_n)+p_0\cdot \mathpzc{1}$. By the previous considerations, there exists a unique $\Upsilon\subseteq\prod_{n\in\oNM}\AG_n$ for $\Lambda:=\{(p(A_n))_{n\in\oNM}\;|\; p\in\Po\}$ where $p(A_n)$ is defined in $\AG_n$ such that $\Lambda\subseteq\Upsilon$ and $((\AG_n)_{n\in\oNM},\Upsilon)$ defines a continuous field of unital $C^\ast$-algebras. Since $(A_n)_{n\in\oNM}$ is a continuous vector field of the continuous field of unital $C^\ast$-algebras $((\AG_n)_{n\in\oNM},\Upsilon)$, the spectrum $\sigma_{\AG_n}(A_n)$ is continuous in $n\in\oNM$. This is true as $\sigma_{\AG_n}(A_n)=\{1\}$ for all $n\in\oNM$. On the other hand, there does not exist a continuous structure on $(\CG_n)_{n\in\oNM}$ such that the vector field $(A_n)_{n\in\oNM}\in\prod_{n\in\oNM}\CG_n$ is continuous. This is shown in Example~\ref{Chap3-Ex-ConstTermPolNecess} but also follows from Theorem~\ref{Chap3-Theo-ContFieldCALgContSpectr} by a contradiction as $\sigma_{\CG_n}(A_n)=\{0,1\}$ for all $n\in\NM$ and $\sigma_{\CG_n}(A_\infty)=\{1\}$.
\end{example}

\begin{remark}
\label{Chap3-Rem-PTopology}
As explained in Section~\ref{Chap3-Sect-RelP2ContOthTop}, the \pt-continuity is related to a topology on $\Ll(\hs)$ that is called \pt-topology. The \pt-topology is the coarsest topology on $\Ll(\hs)$ so that $\Sigma:\Ll(\hs)\to\ks(\CM)\,,\; A\mapsto\sigma(A)\,,$ is continuous if restricted to the self-adjoint operators $\mathfrak{S}\subseteq\Ll(\hs)$. Similarly, there exists a topology on $\Ll(\hs)$ so that $\Sigma$ is continuous if it is restricted to the normal operators $\mathfrak{N}\subseteq\Ll(\hs)$. Specifically, a neighborhood basis for $A_0\in\Ll(\hs)$ is given by
$$
\us_{\phi,\varepsilon}(A_0) \;
	:= \; \Big\{
		A\in\Ll(\hs) \;\Big|\; \big| \|\phi(A_0)\|-\|\phi(A)\| \big|<\varepsilon
	\Big\}\,,
	\qquad
	\phi\in\Cc_c(\CM)\,,\; \varepsilon>0\,.
$$
The topology on $\Ll(\hs)$ defined by this base is called {\em $\Cc_c$-topology}. Note that the assertions of Proposition~\ref{Chap3-Prop-NormImplP2}, Proposition~\ref{Chap3-Prop-StrongNotImplP2} and Proposition~\ref{Chap3-Prop-P2NotImplStrong} also hold if the \pt-topology is replaced by the $\Cc_c$-topology and normal operators are considered instead of self-adjoint operators.
\end{remark}

For some application, the semi-continuity of the spectra is useful.

\begin{remark}
\label{Chap3-Rem-PContEquivContSpectNormal}
The spectra $(\sigma(A_t))_{t\in\ts}$ is upper semi-continuous (resp. lower semi-continuous) if the maps $\ts\ni t\mapsto\|p(A_t,A_t^\ast)\|\in[0,\infty)\,,\; p\in\Po\,,$ are upper semi-continuous (resp. lower semi-continuous), c.f. \cite{BeMa12} as well. The map $\ts\ni t\mapsto\sigma(A_t)\in\ks(\CM)$ is called
\begin{itemize}
\item[(i)] {\em upper semi-continuous} in $t_0\in\ts$ if, for each $F\in\cs(\RM)$ with $\sigma(A_{t_0})\cap F=\emptyset$, there exists an open neighborhood $U\subseteq\ts$ of $t_0$ such that $\sigma(A_t)\cap F=\emptyset$ for every $t\in U$;
\item[(ii)] {\em lower semi-continuous} in $t_0\in\ts$ if, for each open subset $O\subseteq\RM$ with $\sigma(A_{t_0})\cap O\neq \emptyset$, there exists an open neighborhood $U\subseteq\ts$ of $t_0$ such that $\sigma(A_t)\cap O\neq\emptyset$ for every $t\in U$.
\end{itemize}
This immediately follows by the proof of Theorem~\ref{Chap3-Theo-ContFieldCALgContSpectr}. A similar statement holds for the assertions in Theorem~\ref{Chap3-Theo-P2ContEquivContSpect} and Theorem~\ref{Chap3-Theo-CharContSpectUnitary}.
\end{remark}

\section{H\"older-continuous behavior of the spectra and the spectral gaps}
\label{Chap3-Sect-CharHolContSpect}

In view of Theorem~\ref{Chap3-Theo-P2ContEquivContSpect} and Theorem~\ref{Chap3-Theo-CharContSpectUnitary}, the continuity of specific norms of an operator are equivalent to the continuity of the associated spectra with respect to the Hausdorff metric on $\ks(\CM)$. In this section, a quantitative version of these results are proven for self-adjoint and unitary operators and functions of them. For self-adjoint operators, Theorem~\ref{Chap3-Theo-BoundVietContCompVers} shows that the continuous behavior of the spectra is equivalent to the continuity of the boundary. It turns out that, whenever gaps close, the rate of convergence can decrease in general. The results presented here are based on joint work \cite{BeBe16}.

\medskip

In the following several concepts of metric and connected spaces are used. For convenience of the reader, a short reminder is provided here. Let $X$ be a set. A map $d:X\times X\to [0,\infty)$ is called a {\em metric} if for every $x,y,z\in X$
\begin{itemize}
\item[(i)] $d(x,y) =0$ if and only if $x=y$, 
\item[(ii)] $d(x,y)=d(y,x)$, 
\item[(iii)] the {\em triangle inequality} $d(x,y)\leq d(x,z)+d(z,y)$ holds.
\end{itemize} 
A metric $d$ is called an {\em ultrametric} whenever the triangle inequality can be replaced by the stronger version $d(x,y) \leq \max\{d(x,z),d(z,y)\}$. A sequence $(x_n)_{n\in\NM}$ of points in $X$ is called Cauchy whenever, for each $\varepsilon >0$, there is an $n(\varepsilon)\in\NM$ such that if $n,m\geq n(\varepsilon)$ then $d(x_n,x_m)<\varepsilon$. The tuple $(X,d)$ is complete if every {\em Cauchy sequence} converges in $X$. In this case, $(X,d)$ is called a {\em complete metric space}. The induced topology by a metric $d$ on $X$ is the topology defined by the open balls $B_\varepsilon(x):=\{y\in Y\;|\; d(x,y)<\varepsilon\}$ for $x\in X$ and $\varepsilon>0$.

\begin{remark}
\label{Chap3-Rem-MetrSpaFirstCountSeq}
Obviously, complete metric spaces are first countable. Then the topology is described by sequences instead of nets, c.f. Remark~\ref{Chap3-Rem-Net}.
\end{remark}

For $\alpha>0$, let $d^\alpha:X\times X\to[0,\infty)$ be the function $(x,y)\mapsto d(x,y)^\alpha$. The inequality $(a+b)^\alpha\leq a^\alpha+ b^\alpha$ holds for $\alpha\leq 1$ and $a,b\geq 0$. Thus, $d^\alpha$ is a metric on $X$ for $0<\alpha\leq 1$. It is not difficult to check that the topology induced by $d^\alpha$ is finer than the topology of $d$ and vice versa. Hence, the topology induced by $d^\alpha$ is the same as the topology of $d$ for $0<\alpha\leq 1$. If $d$ is an ultrametric, and if $\Phi:[0,\infty)\to [0,\infty)$ is strictly monotone increasing such that $\Phi(0)=0$, then the map $d_\Phi:X\times X\to[0,\infty),\; d_\Phi:= \Phi\circ d,$ is an ultrametric on $X$ defining the same topology as $d$. In particular, for an ultrametric $d$ on $X$, the map $d^\alpha$ for $\alpha>1$ also defines a metric on $X$.

\medskip

Given $\varepsilon >0$, an $\varepsilon$-path $\gamma$ joining $x$ to $y$ is an ordered sequence $(x_0=x, x_1, \cdots , x_{n-1}, x_n=y)$ such that $d(x_{k-1},x_k)<\varepsilon$ for $1\leq k\leq n$. The length of $\gamma$ is $\ell(\gamma) = \sum_{k=1}^n d(x_{k-1},x_k)$. The topological space $X$ is {\em path connected} if for every pair $x,y\in X$ and $\varepsilon >0$, there is an $\varepsilon$-path joining them. A metric $d$ is called a {\em length metric} whenever $d(x,y)$ coincides with the infimum of the lengths of all $\varepsilon$-paths joining $x$ to $y$ for each $\varepsilon >0$, c.f. \cite{Gromov99}. In this case, $(X,d)$ is called a {\em path metric space}.

\medskip

\begin{definition}[H\"older-continuous]
\label{Chap3-Def-HolderCont}
Let $(X,d_X)$ and $(Y,d_Y)$ be two complete metric spaces. For $\alpha >0$, a function $\phi:X\to Y$ is called {\em $\alpha$-H\"older-continuous} if there is a constant $C>0$ such that $d_Y(\phi(x),\phi(y))\leq C d_X(x,y)^\alpha$ for every pair of points $x,y\in X$.  The {\em H\"older constant} is defined as 
$$
\text{\gls{Hol}}^\alpha (\phi)\;
	:= \; \sup_{x\neq y}
		\frac{d_Y(\phi(x),\phi(y))}{d_X(x,y)^\alpha} \, .
$$
Let $\Ff$ be a family of functions $\phi:X\to Y$. This family $\Ff$ is called {\em uniformly $\alpha$-H\"older-continuous} if $\sup_{\phi\in\Ff}\Hol^\alpha (\phi)$ is finite.
\end{definition}

By definition of H\"older-continuity, it immediately follows that every $\alpha$-H\"older-continuous function is continuous. Whenever $\alpha$ is equal to $1$, such a function is called {\em Lipschitz continuous}. In this case, the constant $\Hol^1 (\phi)$ is called {\em Lipschitz constant of $\phi$}. It is well-known that a function $\phi:\RM\to\RM$ has finite H\"older constant $\Hol^\alpha(\phi)$ for $\alpha>1$ if and only if $\phi$ is a constant function. However, this does not hold for general $X$. In particular, if $X$ is a Cantor set and $d_X$ is an ultrametric on $X$, the characteristic function of each clopen set is $\alpha$-H\"older-continuous for every $\alpha >0$. Note that this characteristic function is locally constant but it is not constant at all. The dilation of a function is introduced by following the lines of \cite{Gromov99}.
\begin{definition}[Dilation]
\label{Chap3-Def-Dilation}
Let $(X,d_X)$ and $(Y,d_Y)$ be two complete metric spaces. For $\alpha >0$, the {\em local dilation} of a function $\phi:X\to Y$ is defined by $\Hol^\alpha_\varepsilon(\phi):= \sup_{x\in X} \Hol^\alpha_\varepsilon(\phi)(x)$ where
$$
\Hol^\alpha_\varepsilon(\phi)(x)\; 
	:=  \; \sup \left\{\left.
			\frac{d_Y(\phi(x),\phi(y))}{d_X(x,y)^\alpha} \;\right|\;
			y\in X \text{ such that } 0<d_X(x,y)<\varepsilon
		\right\} \, .
$$
\end{definition}

By definition, the inequality $\Hol^\alpha_\varepsilon(\phi)\leq\Hol^\alpha(\phi)$ holds and, additionally, $\Hol^\alpha_\varepsilon(\phi)$ is non decreasing in $\varepsilon>0$. Thus, the limit 
$$
\dil^\alpha(\phi)(x)\; 
	:= \; \lim_{\varepsilon\downarrow 0} \Hol^\alpha_\varepsilon(\phi)(x)
$$
exists which is possibly infinite. The limit $\dil^\alpha(\phi)(x)$ is called the {\em $\alpha$-dilation of $\phi$ at $x$} whenever it exists. Then $\dil^\alpha(\phi) := \sup_{x\in X} \dil^\alpha(\phi)(x)$ is the {\em $\alpha$-dilation of $\phi$}. By construction, the estimate 
$$
\dil^\alpha(\phi)\; 
	\leq \; \Hol^\alpha(\phi)
$$
holds. If $(X,d_X)$ is a path metric space, the identity $\Hol^\alpha(\phi)=\dil(\phi)$ is derived, c.f. \cite[1.8bis Property]{Gromov99}. 

\medskip

In this section, the topology of $\ts$ is induced by a metric $d$ for which it is complete. The real line $\RM$ or the complex plane $\CM$, or any of their subsets, is endowed with the usual Euclidean metric. 

\medskip

Recall that $\Pt(M)$ for $M>0$ denotes the set of all real-valued polynomials $p:\RM\to\RM$ of the form $p(z)=p_0+p_1\cdot z+p_2\cdot z^2$ with $p_0,p_1,p_2\in \RM$ and $\|p\|_1:=|p_0|+|p_1|+|p_2|\leq M$.

\begin{definition}[\pt-$\alpha$-H\"older-continuous fields, \cite{BeBe16}]
\label{Chap3-Def-P2-Alp-Hold}
Let $(\ts,d)$ be a complete metric space and $(A_t)_{t\in\ts}$ be bounded self-adjoint a field of operators. Then $(A_t)_{t\in\ts}$ is called {\em \pt-$\alpha$-H\"older-continuous} if, for all $M>0$, the family of maps 
$$
\Phi_p:\ts\to[0,\infty)\,,\; 
	t\mapsto \|p(A_t)\|\,,
	\qquad
	p\in \Pt(M)\,,
$$ 
is uniformly $\alpha$-H\"older.
\end{definition}

Clearly, a \pt-H\"older-continuous field of operators $(A_t)_{t\in\ts}$ is  also a \pt-continuous field of operators in terms of Definition~\ref{Chap3-Def-(P2)-Continuity}.

\begin{remark}
Let $(A_t)_{t\in\ts}$ be a bounded self-adjoint field of operators. In view of Remark~\ref{Chap3-Rem-PolyP2ContEquivContSpect}, it suffices for the investigation of the continuous behavior of the spectra at $t_0\in\ts$ to consider only all the polynomials $p(z):=m^2-(x-z)^2$ where $m>0$ is fixed and $x\in[-m,m]$. Here the constant $m:=\sup_{t\in U}\|A_t\|$ depends on the field $(A_t)_{t\in\ts}$ and a suitable open neighborhood $U\subseteq\ts$ of $t_0$ so that $m$ is finite. By restricting to the set $U$, Definition~\ref{Chap3-Def-P2-Alp-Hold} can be relaxed by only requiring that the family of maps
$$
U\ni t\mapsto \|m^2-(A_t-x)^2\|\in[0,\infty)\,,\qquad x\in[-m,m]\,,
$$
is uniformly $\alpha$-H\"older-continuous. In this section, the choice of $U=\ts$ is possible since $m:=\sup_{t\in\ts}\|A_t\|$ is always assumed to be finite.
\end{remark}

Let $(X,d_X)$ be a complete metric space. Recall that $\ks(X)$ is equipped with the Hausdorff metric $d_H$ defined by
$$
d_H(K,F) \; := \;
	\max\left\{
		\sup\limits_{x\in K}\inf\limits_{y\in F} d(x,y),\;
		\sup\limits_{y\in F}\inf\limits_{x\in K} d(x,y)
	\right\} \, ,
	\quad K,\, F\in\ks(X)\, .
$$
Then $(\ks(X),d_H)$ is a complete metric space, c.f. \cite{CastaingValadier77}. In view of Theorem~\ref{Chap3-Theo-VietFellHausMetricEquiv} and Proposition~\ref{Chap3-Prop-ConClosFuncHausCont}, the following assertion holds.

\begin{proposition}
\emph{(\cite[Proposition~3]{BeBe16})}
\label{Chap3-Prop-HolConClosFuncHausHolCont}
Let $(X,d_X)$ and $(Y,d_Y)$ be two complete metric spaces and $\phi:X\to Y$ be an $\alpha$-H\"older-continuous function. Then the map 
$$
\hphi:\ks(X)\to\ks(Y)\,,\; 
	F\mapsto \phi(F)\,,
$$
is $\alpha$-H\"older-continuous with the same H\"older constant $\Hol^\alpha(\phi)=\Hol^\alpha(\hphi)$ where $\ks(X)$ and $\ks(Y)$ are equipped with the corresponding Hausdorff metric.
\end{proposition}

\begin{proof}
As $\phi$ is continuous, it maps compact sets to compact sets. Thus, the map $\hphi$ is well-defined. Let $K,F\in\ks(X)$. Then the estimate
$$
\inf_{y\in F} d_Y(\phi(x),\phi(y))\;
	\leq \; \Hol^\alpha(\phi)\cdot \inf_{y\in F} d_X(x,y)^\alpha
$$
holds for each $x\in K$. By taking the supremum over $x\in K$ and interchanging the role of $K$ and $F$ the desired result follows.
\end{proof}


\begin{proposition}
\label{Chap3-Prop-MaxMinLipCont}
The maximum $\max:\ks(\RM)\to \RM$ and the minimum $\min:\ks(\RM)\to \RM$ are Lipschitz continuous function with Lipschitz constant $1$.
\end{proposition}

\begin{proof}
Let $F_1,F_2\in\ks(\RM)$. Then $\max(F_j)$ exists for $j=1,2$ as $F_j$ is bounded. Without loss of generality, assume that $\max(F_1)\leq\max(F_2)$. Thus, the estimate
$$
\big|\max(F_2) - \max(F_1)\big| \;
 	= \; \max_{\mu\in F_2}\; \min_{\lambda\in F_1} \big( \mu-\lambda \big) \;
 		\leq \; d_H(F_1,F_2)
$$
follows by using the definition of the Hausdorff metric. The case of the minimum is similarly treated.
\end{proof}

\medskip

According to Definition~\ref{Chap3-Def-P2-Alp-Hold}, the supremum
$$
C_M \;
	:=\; \sup \left.\big\{
		\Hol^\alpha(\Phi_p)\;\right|\; \|p\|_1\leq M
	\big\}
$$
is finite for every $M>0$ and a given \pt-$\alpha$-H\"older-continuous field of operators $(A_t)_{t\in\ts}$. With this notion at hand, a quantitative control on the behavior of the spectra is characterized.

\newpage
\begin{theorem}
\emph{(\cite[Theorem~3]{BeBe16})}
\label{Chap3-Theo-CharHolContSpect}
Let $(\ts,d)$ be a complete metric space and $(A_t)_{t\in\ts}$ be a bounded self-adjoint field of operators such that $m:=\sup_{t\in\ts}\|A_t\|<\infty$.
\begin{itemize}
\item[(a)] If $(A_t)_{t\in\ts}$ is a \pt-$\alpha$-H\"older-continuous field, then the map 
$$
\Sigma:\ts\to\ks(\RM)\,,\; 
	t \; \mapsto \; \sigma(A_t)\,,
$$
is $\alpha/2$-H\"older-continuous with respect to the Hausdorff metric on $\ks(\RM)$. In this case, the H\"older constant of $\Sigma$ is bounded by $\sqrt{C_{4m^2+2}}$.
\item[(b)] If the map $\Sigma:\ts\to\ks(\RM)\,,\; t\mapsto\sigma(A_t)\,,$ is $\alpha$-H\"older-continuous with respect to the Hausdorff metric on $\ks(\RM)$ with H\"older constant $C$, then $(A_t)_{t\in\ts}$ is a \pt-$\alpha$-H\"older-continuous field. For $M>0$, the H\"older constant $C_M$ of the family of maps 
$$
\Phi_p:\ts\to[0,\infty)\,,\;
	t \; \mapsto \; \|p(A_t)\|
		\,,\qquad
		p\in\Pt(M),
$$
is bounded by $(2m+1)\cdot M\cdot C$
\end{itemize}
\end{theorem}

\begin{proof}
(a): Let $s,t\in\ts$. According to the definition of the Hausdorff metric, it suffices to show 
$$
\inf_{\mu\in\sigma(A_t)} |\lambda-\mu|\; 
	\leq \;\sqrt{C_{4m^2+2}} \cdot d(s,t)^{\frac{\alpha}{2}} \, ,
		\qquad \lambda\in\sigma(A_s) \,,
$$
since the case of interchanging $s$ and $t$ can be similarly treated. As the left hand side is zero whenever $\lambda\in\sigma(A_t)$, there is no loss of generality in supposing that $\lambda\in \sigma(A_s)\setminus\sigma(A_t)$. The estimate 
$$
\|(A_t-\lambda)^2\| \; 
	\overset{\text{self-}}{\underset{\text{adjoint}}{=}}\; \|A_t-\lambda\|^2 \;
	\leq \; \big(\|A_t\| + |\lambda|\big)^2 \;
	\leq \; 4m^2
$$
is deduced. Hence, $4m^2-(A_t-\lambda)^2$ is a positive operator. Moreover, for the polynomial $p(z)=4m^2-(z-\lambda)^2$, the estimate
$$
\|p\|_1 \; 
	= \; 1+2|\lambda|+4m^2-\lambda^2 \;
		=\; 4m^2+2-(1-|\lambda|)^2 \;
			\leq \; 4m^2+2
$$
is derived. The fact that $\lambda\in\sigma(A_s)$ leads to $\|p(A_s)\|=\|4m^2-(A_s-\lambda)^2\|=4m^2$. On the other hand, since $\lambda\not\in\sigma(A_t)$, the identity 
$$
\|p(A_t)\| \; 
	= \; \sup_{\mu\in\sigma(A_t)} |p(\mu)| \;
		= \; 4m^2-\left(\inf_{\mu\in\sigma(A_t)} |\mu-\lambda|\right)^2
$$
follows by using the spectral theorem and that the operator $p(A_t)$ is positive. Consequently, the estimate
$$
\left(\inf_{\mu\in\sigma(A_t)} |\lambda-\mu|\right)^2\;
	= \; \big| \|p(A_t)\| - \|p(A_s)\| \big| \;
		\leq \; C_{4m^2+2}\cdot d(s,t)^\alpha
$$
is deduced. Hence, the map $\Sigma$ is $\alpha/2$-H\"older-continuous with H\"older constant bounded by $\sqrt{C_{4m^2+2}}$.

\vspace{.1cm}

(b): Since $m:=\sup_{t\in\ts}\|A_t\|$ is finite, the spectrum $\sigma(A_t)$ is contained in the compact set $K:=[-m,m]\subseteq\RM$ for each $t\in\ts$. Let $M>0$ and consider a polynomial $p\in\Pt(M)$, i.e., $p(z)= p_0 + p_1\cdot z + p_2\cdot z^2$ with $p_0,p_1,p_2\in\RM$ satisfying $\|p\|_1:=|p_0|+|p_1|+|p_2|\leq M$. Then, for $x,y\in K$, the estimate
$$
\big| p(x)-p(y) \big| \; 
	\leq\; \big( 
			|p_1|+ |p_2|\cdot | x+y |
		\big)\cdot |x-y| \;
			\leq \; (2m+1)M\cdot|x-y|
$$
holds and so $p$ is Lipschitz continuous on $K$ with Lipschitz constant $(2m+1)M$. The absolute value $g:\RM\to\RM,\; z\mapsto |z|,$ is a Lipschitz continuous function with Lipschitz constant $1$. Proposition~\ref{Chap3-Prop-HolConClosFuncHausHolCont} applies to $\widehat{p}$ and $\hg$ as $p$ and $g$ are Lipschitz-continuous maps. Thus, $\widehat{p}$ and $\hg$ are Lipschitz continuous as well as $\max:\ks(\RM)\to \RM$, c.f. Proposition~\ref{Chap3-Prop-MaxMinLipCont}. 
Recall that the norm map $\ts\ni t\mapsto\|p(A_t)\|\in[0,\infty)$ is represented by
$$
\ts\ni t\mapsto \sigma(A_t)
    \stackrel{\widehat{p}}{\mapsto} \sigma(p(A_t))
     \stackrel{\hg}{\mapsto} |\sigma(p(A_t))|
      \stackrel{\max}{\mapsto} \|p(A_t)\|\, .
$$
Thus, the estimate
$$\big| \|p(A_t)\|- \|p(A_s)\|\big|
	\leq (2m+1)M\cdot d_H(\sigma(A_s),\sigma(A_t))
		\leq (2m+1)M\cdot C\cdot d(s,t)^\alpha
$$
holds for all $s,t\in\ts$ by the previous considerations. Then the supremum $C_M:=\sup\big\{\Hol^\alpha(\Phi_p)\;\big|\; p\in\Pt(M)\big\}$ is bounded by $(2m+1)M\cdot C$. Consequently, the fami\-ly of maps $\ts\ni t\mapsto\|p(A_t)\|\in[0,\infty)\,,\; p\in\Pt(M)\,,$ is uniformly $\alpha$-H\"older-continuous. Hence, $(A_t)_{t\in\ts}$ is a \pt-$\alpha$-H\"older-continuous field.
\end{proof}

\begin{remark}
\label{Chap3-Rem-CharHolContSpect}
In assertion (i), the family of spectra $\sigma(A_t)\,,$ $t\in\ts\,,$ only varies $\alpha/2$-H\"older-continuous and not $\alpha$-H\"older-continuous. According to Example~\ref{Chap3-Ex-AlmostMathieu} it is known that this result cannot be improved, in general. Specifically, whenever an isolated gap tip closes, it can, a priori, close with rate $\alpha/2$, c.f. Theorem~\ref{Chap3-Theo-HolConClosGap}. In Example~\ref{Chap3-Ex-AlmostMathieu}, it was discussed that the norms of the Almost -Mathieu operator behave Lipschitz continuous while isolated gap tips close only with rate $1/2$, i.e., these gaps are $1/2$-H\"older-continuous. In this case, also the spectra behave only $1/2$-H\"older-continuous. In view of that the result of Theorem~\ref{Chap3-Theo-CharHolContSpect} is optimal.
\end{remark}

According to the observation discussed in Remark~\ref{Chap3-Rem-CharHolContSpect}, the behavior of the boundary is crucial. Especially, those gaps that are close are of particular interest. The behavior of the boundary is studied next, in more detail, for general bounded self-adjoint operators.

\begin{theorem}
\emph{(\cite[Theorem~4]{BeBe16})}
\label{Chap3-Lem-BoundHolCont}
Let $(A_t)_{t\in\ts}$ be a \pt-$\alpha$-H\"older-continuous field of bounded self-adjoint operators such that $m:=\sup_{t\in\ts}\|A_t\|<\infty$.
\begin{itemize}
\item[(a)] The supremum $\ts\ni t\mapsto\sup\big(\sigma(A_t)\big)\in\RM$ and the infimum $\ts\ni t\mapsto\inf\big(\sigma(A_t)\big)\in\RM$ are $\alpha$-H\"older-continuous with H\"older constant less than $C_{m+1}$.
\item[(b)]For $t_0\in\ts$ and a gap $(a_{t_0},b_{t_0})$ of $\sigma(A_{t_0})$, there exist an open neighborhood $U_0$ of $t_0$ and gaps $(a_t,b_t)$ of $\sigma(A_t)$ for $t\in U_0$ such that the gaps $\big((a_t,b_t)\big)_{t\in U_0}$ are $\alpha$-H\"older-continuous with $\alpha$-dilation less than $3C_{4m^2+2}/|b_{t_0}-a_{t_0}|$, i.e.,
$$
\max\big\{|a_s-a_t|,|b_s-b_t|\big\}\;
	\leq \; \frac{3C_{4m^2+2}}{|b_{t_0}-a_{t_0}|}\;\;d(s,t)^\alpha \, ,
	\qquad
	s \, , \, t\in U_0 \, .
$$
\end{itemize}
\end{theorem}

\begin{proof}
(a): For $t\in\ts$ and $\lambda\in\sigma(A_t)$, the inequality $m\pm\lambda\geq 0$ holds. Then the equations 
\begin{gather*}
\begin{aligned}
\|m+A_t\| \; 
	&= \; \sup_{\lambda\in\sigma(A_t)} |m+\lambda| \;
		&= \; m + \sup_{\lambda\in\sigma(A_t)} \lambda \, ,\\
\|m-A_t\| \; 
	&= \; \sup_{\lambda\in\sigma(A_t)} |m-\lambda| \;
		&= \; m - \inf_{\lambda\in\sigma(A_t)} \lambda \, ,
\end{aligned}
\end{gather*}
follows by the spectral theorem. These polynomials $\RM\ni z\mapsto m\pm z\in\RM$ have $1$-norm equal to $m+1$. Applying the \pt-$\alpha$-H\"older-continuity, the desired result is obtained.

\vspace{.1cm}

\begin{figure}[htb]
\centering
\includegraphics[scale=2]{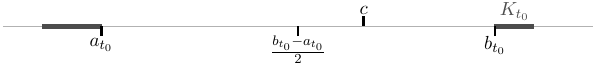}
\caption{The choice of $c\in\RM$ in Theorem~\ref{Chap3-Lem-BoundHolCont}~(b).}
\label{Chap3-Fig-ThmOpenGap}
\end{figure}

(b): Let $t_0\in\ts$ and $(a_{t_0},b_{t_0})$ be a gap of $\sigma(A_{t_0})$, i.e., $a_{t_0},\; b_{t_0}\in\sigma(A_{t_0})$ and the intersection $(a_{t_0},b_{t_0})\cap\sigma(A_{t_0})$ is empty. Subdivide the interval $(a_{t_0},b_{t_0})$ into six intervals of equal length $r:=(b_{t_0}-a_{t_0})/6$ and set $c:=a_{t_0}+4r=b_{t_0}-2r$, c.f. Figure~\ref{Chap3-Fig-ThmOpenGap}. The field $(A_t)_{t\in\ts}$ is \pt-continuous by assumption and so the boundaries of $(\sigma(A_t))_{t\in\ts}$ are continuous, c.f. Theorem~\ref{Chap3-Theo-P2ContEquivContSpect}. Thus, according to \nameref{(CB2)}, there exists an open neighborhood $U_0$ of $t_0$ such that the spectrum $\sigma(A_t)$ has a gap $(a_t,b_t)$ satisfying $|a_t-a_{t_0}|<r$ and $|b_t-b_{t_0}|<r$ for all $t\in U_0$. Consequently, for $t\in U_0$, the following inequalities hold:
\begin{gather*}
\begin{aligned}
&(1)\quad b_t-c \, 
	&&= \; b_t-b_{t_0}+b_{t_0}-c \;
	&&\Rightarrow\; r \, 
	&&< \; b_t-c \, 
	&&< \,3r\, ,\\
&(2)\quad c-a_t \, 
	&&= \; c-a_{t_0}+a_{t_0}-a_t \;
	&&\Rightarrow\; c-a_t\, 
	&&>\; 4r-r\, 
	&&= \, 3r \, .
\end{aligned}
\end{gather*}
Thus, $c$ is closer to $b_t$ than to $a_t$. Hence, if $t\in U_0$, the infimum of the spectrum of $(A_t-c)^2$ is exactly $(b_t-c)^2$. Then $c\in(a_{t_0},b_{t_0})$ leads to $|c|\leq \|A_t\|\leq m$ so that $(b_t-c)^2\leq \|(A_t-c)^2\|\leq 4m^2$. Consequently, $4m^2-(b_t-c)^2=\|4m^2-(A_t-c)^2\|$ is derived whenever $t\in U_0$. This implies
$$
\big|\|p(A_t)\|-\|p(A_s)\|\big|\; 	
	= \; \big|(b_t-c)^2-(b_s-c)^2\big| \; 
		=\; |b_t-b_s||b_t+b_s-2c| \, ,
		\qquad s,t\in U_0 \, .
$$ 
The $1$-norm of the polynomial $p(z)=4m^2-(z-c)^2$ is estimated by
$$
\|p\|_1 \; 
	= \; 1+2|c|+4m^2-c^2 \;
		= \; 4m^2+2-(1-|c|)^2 \;
			\leq \; 4m^2+2\, .
$$ 
Then (1) provides the inequality $|b_t+b_s-2c| >2r$. By the previous considerations, this leads to
$$
|b_s-b_t| \; 
	\leq \frac{3C_{(4m^2+2)}}{|b_{t_0}-a_{t_0}|}\;\;d(s,t)^\alpha
	 \, ,\qquad s,t\in U_0 \, .
$$
Changing $c$ into $a_{t_0}+2r=b_{t_0}-4r$ yields the same estimate of $|a_s-a_t|$ for $s,t\in U_0$. 
\end{proof}

\medskip

Note that the assertion Theorem~\ref{Chap3-Lem-BoundHolCont}~(b) provides an upper bound on the $\alpha$-dilation of $\big((a_t,b_t)\big)_{t\in U_0}$.

\medskip

According to Remark~\ref{Chap3-Rem-MetrSpaFirstCountSeq}, the topology of a complete metric space is described by sequences. Let $(\ts,d)$ be a complete metric space and $(A_t)_{t\in\ts}$ be a field of bounded self-adjoint operators. For $t_0\in\ts$, the notion of a gap tip $c$ of $\sigma(A_{t_0})$ is formulated as follows, c.f. Definition~\ref{Chap3-Def-ClosedGap}. There is a sequence $(t_n)_{n\in\NM}$ tending to $t_0$ and there exist maps $a:\NM\to\RM$, $b:\NM\to\RM$ such that $(a_n,b_n)$ is a gap in $\sigma(A_{t_n})$ for $n\in\NM$ and $\lim_{n\to\infty} a_n = c = \lim_{n\to\infty} b_n$. A gap tip $c$ is called isolated of type 1 (type 2) if there is a $\delta>0$ such that $(c-\delta,c+\delta)\subseteq\sigma(A_{t_0})$ $\big((c-\delta,c+\delta)\cap \sigma(A_{t_0})=\{c\}\big)$, c.f. Definition~\ref{Chap3-Def-ClosedGap}.

\begin{theorem}
\emph{(\cite[Theorem~5]{BeBe16})}
\label{Chap3-Theo-HolConClosGap}
Let $(A_t)_{t\in\ts}$ be a \pt-$\alpha$-H\"older-continuous field of bounded self-adjoint operators such that $m:=\sup_{t\in\ts}\|A_t\|<\infty$. Consider a closed isolated (type 1 or type 2) gap tip $c$ of $\sigma(A_{t_0})$ with associated sequence $(t_n)_{n\in\NM}$ tending to $t_0$ and gaps $(a_n,b_n)$ of $\sigma(A_{t_n})$ for $n\in\NM$. Then there is an $n_0\in\NM$ such that 
$$
b_n-a_n \;
	\leq \; 2\cdot\sqrt{C_{4m^2+2}}\;\cdot \;d(t_n,t_0)^{\alpha/2}
	 \, ,\qquad n\geq n_0 \, .
$$
\end{theorem}

\begin{figure}[htb]
\centering
\includegraphics[scale=1.56]{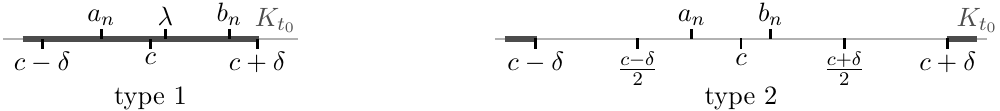}
\caption{Illustration of the proof of Theorem~\ref{Chap3-Theo-HolConClosGap}.}
\label{Chap3-Fig-ThmClosGap}
\end{figure}

\begin{proof}
The proof is organized as follows. (i) First the case of isolated gap tip of type 1 is proven. (ii) Then the assertion for an isolated gap tip of type 2 is verified.

\vspace{.1cm}

(i): Let $t_0\in\ts$ and $c\in \sigma(A_{t_0})$ be a closed isolated gap tip of type 1 with associated $\delta >0$ such that $(c-\delta,c+\delta)\subseteq \sigma(A_{t_0})$.

\vspace{.1cm}

Since $(a_n)_{n\in\NM}$ and $(b_n)_{n\in\NM}$ converge to $c$ by definition, there is an $n_0\in\NM$ such that the inequality $\max\{|a_n-c|,|b_n-c|\}<\delta$ holds for all $n\geq n_0$, c.f. Figure~\ref{Chap3-Fig-ThmClosGap}. Thus, the inequalities $c-\delta< a_n<b_n<c+\delta$ are derived if $n\geq n_0$. For $n\geq n_0$, define $\lambda_n:=(a_n+b_n)/2$. It follows that $\lambda_n \in \sigma(A_{t_0})$ and the distance of $\lambda_n$ to $\sigma(A_{t_n})$ is exactly $(b_n-a_n)/2$. Consequently, the identities
$$\|4m^2-(A_{t_0}-\lambda_n)^2\| \; 
	= \; 4m^2 \, ,
\qquad
\|4m^2-(A_{t_n}-\lambda_n)^2\| \;
	= \; 4m^2-\frac{(b_n-a_n)^2}{4} \, ,
$$
follow. The $1$-norm $\|p_n\|_1$ of the polynomial $p_n(z)= 4m^2-(z-\lambda_n)^2$ satisfies $\|p_n\|_1\leq 4m^2+2$ for $n\geq n_0$, c.f. the proof of Theorem~\ref{Chap3-Lem-BoundHolCont}~(b). Due to the \pt-$\alpha$-H\"older-continuity, the estimate
$$
(b_n-a_n)^2\;
	= \; 4\cdot \Big(\|4m^2-(A_{t_0}-\lambda_n)^2\| - \|4m^2-(A_{t_n}-\lambda_n)^2\| \Big) \;
	\leq \; 4C_{4m^2+2} \; d({t_n},t_0)^\alpha
$$
is deduced for $n\geq n_0$ leading to the desired result.

\vspace{.1cm}

(ii): Let $t_0\in\ts$ and $c\in \sigma(A_{t_0})$ be an isolated gap tip of type 2 with associated $\delta >0$ satisfying $(c-\delta,c+\delta)\cap \sigma(A_{t_0})=\{c\}$. Due to the convergence of $(a_n)_{n\in\NM}$ and $(b_n)_{n\in\NM}$ to $c$, there exists an $n_0\in\NM$ such that $a_n,\, b_n\in (c-\delta/2,c+\delta/2)$ for $n\geq n_0$, c.f. Figure~\ref{Chap3-Fig-ThmClosGap}. Consider the polynomials $p_n(z):=4m^2-(a_n-z)^2$ and $q_n(z):=4m^2-(b_n-z)^2$ for $n\in\NM$. Then the equations $\|q_n(A_{t_n})\| = \|p_n(A_{t_n})\| = 4m^2$ and
$$
\|p_n(A_{t_0})\| \;
	= \; 4m^2-(a_n-c)^2 \, ,
\qquad
\|q_n(A_{t_0})\| \;
	= \; 4m^2-(b_n-c)^2 \, ,
$$
hold. Similarly to (i), the $1$-norms of the polynomials $p_n$ and $q_n$ are uniformly bounded by $C_{4m^2+2}$, c.f. the proof of Theorem~\ref{Chap3-Lem-BoundHolCont}~(b). This leads to the estimates
$$
\begin{array}{c}
	(a_n-c)^2\;
		= \; \big|\|p_n(A_{t_0}\| - \|p_n(A_{t_n}\| \big| \;
		\leq \; C_{4m^2+2} \cdot d({t_n},t_0)^\alpha\, ,\\[0.1cm]
	(b_n-c)^2\;
		= \; \big|\|q_n(A_{t_0}\| - \|q_n(A_{t_n}\| \big| \;
		\leq \; C_{4m^2+2} \cdot d({t_n},t_0)^\alpha\, , 
\end{array}
	\qquad n\geq n_0\, ,
$$
by the required \pt-$\alpha$-H\"older-continuity. Consequently, the desired inequality
$$
|a_n-b_n| \; 
	\leq \; 2 \sqrt{C_{4m^2+2}}\cdot d({t_n},t_0)^\alpha\, ,
	\qquad n\geq n_0,
$$
is deduced.
\end{proof}

\begin{remark}
\label{Chap3-Rem-ExtensHolConClosGap}
The requirements in Theorem~\ref{Chap3-Theo-HolConClosGap} can be relaxed under few additional assumptions. If $a_n<b_n\leq c$ holds for $n\in\NM$, then it is only necessary that there exists a $\delta>0$ such that $(c-\delta,c]\subseteq\sigma(A_{t_0})$ or $(c-\delta,c)\cap\sigma(A_{t_0})=\emptyset$. The requirements change similarly if $c\leq a_n< b_n$ holds.
\end{remark}

\medskip

The closing gap condition is observed in many models. First, in the Almost-Mathieu model $(H_\alpha)_{\alpha\in[0,1]}$, c.f. Example~\ref{Chap3-Ex-AlmostMathieu}. Since it is \pt-Lipschitz continuous proven in \cite{Bel94}, the gap widths behave at least $1/2$-H\"older-continuous. A semi-classical calculation \cite{RaBe90,HeSj90} confirms this prediction and shows that the widths do not converge faster.

\medskip

Another situation where gaps close is provided by a small perturbation of the Laplacian on $\ZM$. Namely for $\lambda \geq 0$, let $H_\lambda:\ell^2(\ZM)\to\ell^2(\ZM)$ be defined on $\ell^2(\ZM)$ by 
$$H_\lambda\psi(n) \; = \;
   \psi(n+1)+\psi(n-1)+\lambda V(n)\psi(n) \, ,
$$ 
where $V(n)$ takes only finitely many values. This is a specific pattern equivariant Schr\"odin\-ger operator, c.f. Theorem~\ref{Chap2-Theo-PESchrOpFinRang} below. For $\lambda=0$, the spectrum of $H_\lambda$ is the interval $[-2,+2]$. The prediction provided by the {\em Gap Labeling Theorem} \cite{Bel92,Bel93}, indicates that for certain potentials $V$, gap may open as $\lambda$ increases from $\lambda =0$, i.e., there are gap tips at $\lambda=0$. Explicit calculations have been made on several examples, which show the opening of gaps, such as the Fibonacci sequence \cite{SiMo90,DaGo11}, Thue-Morse sequence \cite{Bel90} and the period doubling sequence \cite{BeBoGh91}. Furthermore, for Schr\"odinger operators associated with a compact metric space, the authors of \cite{AvBoDa12} proved that each gap can open by perturbing the potential.

\medskip

In Theorem~\ref{Chap3-Theo-HolConClosGap}, it was required that the gap tip is isolated of type 1 or type 2. In general, this condition cannot be relaxed as be shown by the following example, c.f. \cite[Example~1]{BeBe16}.

\begin{figure}[htb]
\centering
\includegraphics[scale=1.4]{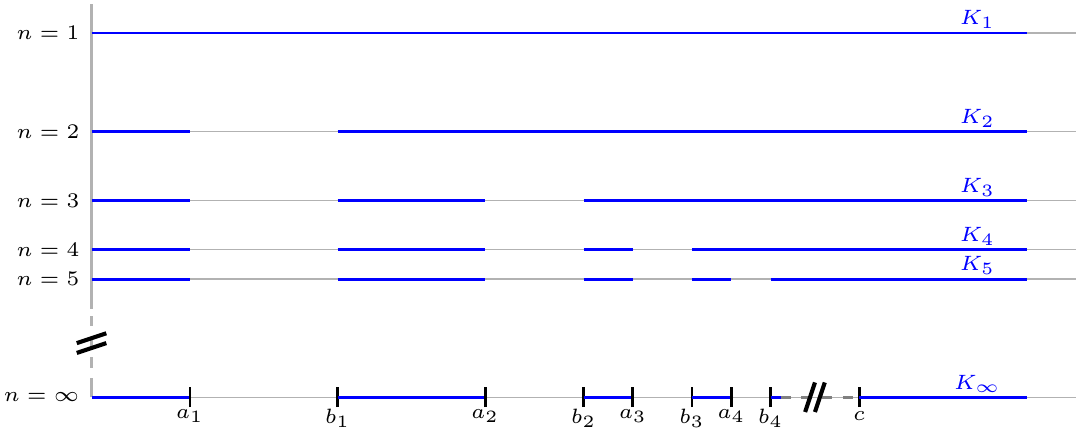}
\caption{The spectra $(K_n)_{n\in\oNM}$.}
\label{Chap3-Fig-SlowGapClos}
\end{figure}

\begin{example}
\label{Chap3-Ex-SlowClosGap}
Consider $\ts:=\oNM$ the one-point compactification of $\NM$. This topological space is locally compact, Hausdorff, second-countable. Thus, it is metrizable. Define, for $\kappa>1$, an ultrametric on $\ts$ by
$$
d(n,m)\;
	:= \; 
		\begin{cases}
			e^{-\kappa^{\min\{n,m\}}},\qquad &m\neq n \, ,\\
			0,\qquad &n=m \, .
		\end{cases}
$$
Let $(a_n,b_n)_{n\in \NM}$ be a sequence of real numbers such that $0<a_n<b_n<a_{n+1}<c$ and $\sup_{n\in\NM} a_n=c$. Thus, the limits $\lim_{n\to\infty}a_n$ and $\lim_{n\to\infty}b_n$ exist and they are equal to $c$. Now, define inductively a sequence of compact sets $(K_n)_{n\in\NM}$ associated with $(a_n,b_n)\, ,\; n\in\NM$. More precisely, set $K_0:=[0,m]$ with $m>c$ and, for $n\in\NM$, define $K_{n+1}:=K_n\setminus (a_{n+1},b_{n+1})$. Each compact subset $K_n,\; n\in\NM,$ of $\RM$ can be seen as the spectrum of a bounded self-adjoint operator $A_n$. Due to construction, the sequence $(K_n)_{n\in\NM}$ converges to $K_\infty=\bigcap_{n\in\NM} K_n$ in the Hausdorff metric and $c$ is a gap tip of $K_\infty$.

\vspace{.1cm}

\setlength\parindent{1em} 
\indent By construction, $(A_n)_{n\in\NM}$ is a field of bounded self-adjoint operators with spectra $\sigma(A_n)=K_n\,,\; n\in\NM$. If $|b_{n+1}-a_{n+1}|<|b_n-a_n|$ holds for $n\in\NM$, then $d_H(K_n,K_\infty)=|b_{n+1}-a_{n+1}|/2$. 

\vspace{.1cm}

\setlength\parindent{1em} 
\indent Now choose $(a_n)_{n\in\NM}$ and $(b_n)_{n\in\NM}$ such that there is a $C>0$ and
$$
|b_{n+1}-a_{n+1}|
		= \; C\cdot d(n,\infty)^{\alpha/2}\,,
		\qquad
		n\in\NM\,,
$$
hold for a fixed $\alpha>0$. Then $(K_n)_{n\in\NM}$ is $\alpha/2$-H\"older-continuous with constant $C/2$. Note that then $(A_n)_{n\in\NM}$ is $\alpha/2$-H\"older-continuous as well, c.f. Theorem~\ref{Chap3-Theo-CharHolContSpect}~(b). On the other hand, the identity $|b_n-a_n|= C\cdot d(n,\infty)^{\alpha/(2\kappa)}$ holds implying that the gaps width are H\"older-continuous in $n$ but with an exponent $\alpha/(2\kappa)<\alpha/2$. More specifically, there does not exist a constant $C'>0$ such that $|b_n-a_n|\leq C'\cdot d(n,\infty)^{\alpha/2}$.
\setlength\parindent{0em} 
\end{example}

In view of Theorem~\ref{Chap3-Theo-CharContSpectUnitary}, one may expect an analog of Theorem~\ref{Chap3-Theo-CharHolContSpect} for unitary operators. This is indeed the case which is discussed next.

\begin{theorem}
\emph{(\cite[Theorem~6]{BeBe16})}
\label{Chap3-Theo-CharHolContSpectUnitary}
Let $(\ts,d)$ be a complete metric space and $(A_t)_{t\in\ts}$ be a unitary field of operators.
\begin{itemize}
\item[(a)] Suppose $(A_t)_{t\in\ts}$ satisfies, for $\alpha>0$, that the family of maps $\ts\ni t\mapsto\|1+E\cdot A_t\|\in[0,\infty)\,,\; E\in\SM^1,$ is uniformly $\alpha$-H\"older-continuous with H\"older constant $C>0$. Then 
$$
\Sigma:\ts\to\ks(\CM)\,,\; 
	t \; \mapsto \; \sigma(A_t)\,,
$$
is $\alpha/2$-H\"older-continuous with respect to the Hausdorff metric on $\ks(\CM)$. In this case, the H\"older constant is bounded by $4\cdot C$.
\item[(b)] If the map $\Sigma:\ts\to\ks(\CM)\,,\; t\mapsto\sigma(A_t)\,,$ is $\alpha$-H\"older-continuous with respect to the Hausdorff metric on $\ks(\CM)$ with H\"older constant $C$, then the family of maps 
$$
\ts\ni t \; 
	\mapsto \; \|1+E\cdot A_t\|\in[0,\infty)\,,
	\qquad E\in\SM^1\,,
$$
is uniformly $\alpha$-H\"older-continuous with H\"older constant $C$.
\end{itemize}
\end{theorem}

\begin{proof}
(a): Let $s,t\in\ts$. According to the definition of the Hausdorff metric, it suffices to show 
$$
\inf_{\mu\in\sigma(A_t)} |\lambda-\mu|\; 
	\leq \;4\cdot C \cdot d(s,t)^{\frac{\alpha}{2}} \, ,
		\qquad \lambda\in\sigma(A_s) \,,
$$
since the case of interchanging $s$ and $t$ can be similarly treated. As the left hand side is zero whenever $\lambda\in\sigma(A_t)$ there is no loss of generality in supposing that $\lambda\in \sigma(A_s)\setminus\sigma(A_t)$. Then the identities 
$$
\|1+\overline{\lambda}\cdot A_s\| \; = \; 2,
	\qquad
		\|1+\overline{\lambda}\cdot A_t\| \; = \; \sup_{\mu\in\sigma(A_t)} \sqrt{4-|\mu-\lambda|^2}
$$
are derived, c.f. the computations made in Lemma~\ref{Chap3-Lem-PresenceSpectrumUnitary}. Since $|\mu-\lambda|\leq 2$, a short computation leads to
\begin{align*}
\inf_{\mu\in\sigma(A_t)}|\mu-\lambda|^2 \; 
	&= \; \Big| \|1+\overline{\lambda}\cdot A_s\|^2 - \|1+\overline{\lambda}\cdot A_t\|^2 \Big|\\
		&= \; \Big( 
				\underbrace{\|1+\overline{\lambda}\cdot A_s\|}_{\leq 2} +
				\underbrace{\|1+\overline{\lambda}\cdot A_t\|}_{\leq 2} 
			\Big) \cdot \Big|
				\|1+\overline{\lambda}\cdot A_s\| - \|1+\overline{\lambda}\cdot A_t\|
			\Big|\\
				&\leq \; 4\cdot C\cdot d(s,t)^\alpha
\end{align*}
where $\|A_s\|=\|A_t\|=1=|\lambda|$ was used in the last estimate.

\vspace{.1cm}

(b): For $t\in\ts$, the spectrum $\sigma(A_t)$ is contained in the compact set $\SM^1\subseteq\CM$. For $E\in\SM^1$ with polynomial $p:\CM\to\CM\, ,\; p(z):=1+E\cdot z\, ,$ the equalities
$$
\big| p(x)-p(y) \big| \; 
	=\; |E|\cdot |x-y| \;
			= \; |x-y|
$$
hold, namely $p$ is Lipschitz continuous on $\CM$ with Lipschitz constant $1$. Similarly, the absolute value $g:\CM\to\RM,\; z\mapsto |z|,$ is a Lipschitz continuous function with Lipschitz constant $1$. As $p$ and $g$ are closed as well, Proposition~\ref{Chap3-Prop-HolConClosFuncHausHolCont} applies. Furthermore, Proposition~\ref{Chap3-Prop-MaxMinLipCont} is applicable to the map $|\sigma(p(A_t))|\overset{\max}{\longmapsto}\|p(A_t)\|$ since the inclusion $|\sigma(p(A_t))|\subseteq [0,2]$ holds for all $t\in\ts$. Thus, by the representation 
$$
\ts\ni t\longmapsto \sigma(A_t)
    \stackrel{\hat{p}}{\longmapsto} \sigma(p(A_t))
     \stackrel{|\cdot|}{\longmapsto} |\sigma(p(A_t))|
      \stackrel{\max}{\longmapsto} \|p(A_t)\|\in[0,\infty)\, ,
$$
the norm map $\ts\ni t\mapsto\|p(A_t)\|\in[0,\infty)$ is $\alpha$-H\"older-continuous as a composition of one $\alpha$-H\"older-continuous map and three Lipschitz continuous maps. Thus, the estimate
$$\big| \|p(A_t)\|- \|p(A_s)\|\big|
	\leq d_H(\sigma(A_s),\sigma(A_t))
		\leq C\cdot d(s,t)^\alpha
$$
follows for all $s,t\in\ts$. Consequently, the family of maps $\ts\ni t\mapsto\|1+E\cdot A_t\|\in[0,\infty)$ for $E\in\SM^1$ is uniformly $\alpha$-H\"older-continuous with H\"older constant $C$.
\end{proof}

\cleardoublepage


\chapter{Generalized Schr\"odinger operators associated with dynamical systems}
\label{Chap2-SchrOpDynSyst}
\stepcounter{section}
\setcounter{section}{0}

According to Theorem~\ref{Chap3-Theo-PContEquivContSpectNormal}, the spectra of a bounded normal field of operator vary continu\-ously if and only if the field of $C^\ast$-algebras defined by these operators is continuous. Additionally, the structure of a continuous field of groupoids yields a structure of a continuous field of groupoid $C^\ast$-algebras, c.f. \cite{LaRa99} or Chapter~\ref{Chap4-ToolContBehavSpectr} for more details. Thus, it is a reasonable strategy to use continuous fields of groupoids to prove the continuous behavior of the spectra of operators. Having this in mind, this chapter provides the framework of groupoid $C^\ast$-algebras with the main focus on Schr\"odinger operators associated with dynamical systems. Additionally, a homeomorphism is defined mapping the space $\SG\big(\as^G\big)\subseteq\cs\big(\as^G\big)$ of subshifts of $\as^G$ equipped with the Hausdorff-topology onto the topological space of dictionaries $\DG$, c.f. Theorem~\ref{Chap2-Theo-Shift+DictSpace}. This new approach is intensively used to verify sufficient conditions of subshifts being periodically approximable in Chapter~\ref{Chap5-OneDimCase}, Chapter~\ref{Chap6-HigherDimPerAppr} and Chapter~\ref{Chap7-Examples}.

\medskip

This chapter is organized as follows. As an introductory to the theory, we provide a short discussion on the concepts used in the one-dimensional case in Section~\ref{Chap2-Sect-ExampleAsZMSchrOp}. Then the basic notions and properties of dynamical systems are described in Section~\ref{Chap2-Sect-DynSystGroupoid}. There the main focus is on the study of the set $\SG(X)$ of $G$-invariant, closed subsets of $X$ for a dynamical system $(X,G,\alpha)$ which is called space of dynamical subsystems. Of particular interest are symbolic dynamical systems $(\as^G,G,\alpha)$ for a discrete, countable group that are studied in Section~\ref{Chap2-Sect-SymbDynSyst}. There the notion of a dictionary is extended to general symbolic dynamical systems while a dictionary is not defined for a specific element of $\as^G$ as it is usually done in the case of $(\as^\ZM,\ZM,\alpha)$. Instead, a dictionary is defined as a set of patterns satisfying a heredity and an extensibility condition, c.f. Definition~\ref{Chap2-Def-Dictionary} for more details. Additionally, the local pattern topology on the set of dictionaries is established. As it turns out, the space $\SG\big(\as^G\big)\subseteq\cs\big(\as^G\big)$ of subshifts of $\as^G$ equipped with the Hausdorff-topology and the space of dictionaries $\DG$ are homeomorphic, c.f. Theorem~\ref{Chap2-Theo-Shift+DictSpace}. Next, Section~\ref{Chap2-Sect-GroupoidCalgebras} provides a summary of groupoids and the associated reduced and full $C^\ast$-algebra. Since this work shall be accessible for a broad readership, the proofs for the basic properties of topological groupoids and the construction of a groupoid $C^\ast$-algebra are provided here by following the lines of \cite{Renault80}. The class of transformation group groupoids associated with a dynamical system first introduced in \cite{Ehr57} is separately studied in Section~\ref{Chap2-Sect-TransformationGroupGroupoid} since the main focus in this work is on this class of groupoids. Then the notion of generalized Schr\"odinger operators is introduced and fundamental spectral properties are proven in connection with the dynamics, c.f. Section~\ref{Chap2-Sect-SchrOp}. Specifically, the constancy of the spectra of these operators is proven for all minimal dynamical systems. Additionally, minimality implies that the spectrum of a generalized Schr\"odinger operator is equal to its essential spectrum, c.f. Theorem~\ref{Chap2-Theo-ConstSpectrMinimal}. Finally, pattern equivariant Schr\"odinger operators associated with symbolic dynamical systems $(\as^G,G,\alpha)$ are introduced in Section~\ref{Chap2-Sect-PattEqSchrOp} by following the lines of \cite{KePu00,Kel03}. More precisely, for $K\subseteq G$ a finite set and $p_h:\as^G\to\CM,\; h\in K,$ and $p_e:\as^G\to\RM$ pattern equivariant functions, the pattern equivariant Schr\"odinger operators $H_\xi:\ell^2(G)\to\ell^2(G)\,,\;\xi\in\as^G\,,$ is defined by
$$
(H_\xi\psi)(g) \; 
	:= \; \left( 
			\sum\limits_{h\in K} p_h\big(\alpha_{g^{-1}}(\xi)\big) \cdot \psi(g\,h^{-1}) + \overline{p_h\big(\alpha_{(gh)^{-1}}(\xi)\big)} \cdot \psi(gh)
		\right)
		+ p_e\big(\alpha_{g^{-1}}(\xi)\big) \cdot \psi(g)
$$
where $\psi\in\ell^2(G)$ and $g\in G$. Note that if $G=\ZM$ and $K=\{-1,1\}$, the pattern equivariant Schr\"odinger operator $H_\xi:\ell^2(\ZM)\to\ell^2(\ZM)$ is exactly the Jacobi operator that is typically studied in the literature, c.f. \cite{Bel86,Sut95,Dam00Gord,Stollmann01,Bel02,Baa02,Bel03,Dam07,DaLiQu14,BaDaGr15,DaEmGo15,KellendonkLenzSavinien15}. We prove that these operators are contained in the groupoid $C^\ast$-algebra $\CG^\ast_{red}\big(\as^G\rtimes_\alpha G\big)$.

\section{Symbolic dynamical system over \texorpdfstring{$\ZM$}{ZM} and its Schr\"odinger op\-era\-tors}
\label{Chap2-Sect-ExampleAsZMSchrOp}

In this section, a short introduction to symbolic dynamical system $(\as^\ZM,\ZM,\alpha)$ is presented for a finite alphabet $\as$. Furthermore, a short review of the associated Schr\"odinger and Jacobi operator is provided. The reader is also referred to Section~\ref{Chap2-Sect-SymbDynSyst} and Chapter~\ref{Chap5-OneDimCase} for a more detailed discussion of the symbolic dynamical systems $(\as^\ZM,\ZM,\alpha)$. General references for symbolic dynamical systems $(\as^\ZM,\ZM,\alpha)$ are the textbooks \cite{Queffelec87,Fogg02,Queffelec10,BaakeGrimm13}.

\medskip

Let \gls{as} be a finite set equipped with the discrete topology. Such a set is called {\em alphabet} and its elements are called {\em letters}. Then the product $\as^\ZM:=\prod_{n\in\ZM}\as$ defines a compact, second countable, Hausdorff space, c.f. \cite{Tyc30,Cec37}. An element $\xi\in\as^\ZM$ is a map $\xi:\ZM\to\as$. It is represented by a {\em two-sided infinite word}
$$
\xi \; 
	= \; \ldots a_{-6}\, a_{-5}\, a_{-4}\, a_{-3}\, a_{-2}\, a_{-1}\, |\, a_{0}\, a_{1}\, a_{2}\, a_{3}\, a_{4}\, a_{5}\, a_{6}\ldots\;
	\in\as^\ZM	
		\,,\qquad a_i\in\as\,,\; i\in\ZM\,,
$$
where $|$ indicates the origin. A finite word $u$ with letters in $\as$ is an element of the Cartesian product $\as^n$ for an $n\in\NM$. It is represented by the concatenation $u=a_1\ldots a_n$ of letters $a_i\in\as$. For $v\in\as^n$, the {\em word length} $|v|$ of $v$ is defined by $n$. For two words $u\in\as^n$ and $v\in\as^k$, the concatenation is defined by $uv\in\as^{n+k}$. For the sake of convenience, $\epsilon$ denotes the {\em empty word} of length $|\epsilon|$ equal to zero. This is the abstract neutral element of the concatenation map, i.e., $v\epsilon = \epsilon v= v$ holds for $v\in\as^n$. Then the union $\as^\ast:=\bigcup_{k\in\NM}\as^k\cup\{\epsilon\}$ is the {\em free monoid} on the set $\as$. The set $\as^+:=\as^\ast\setminus\{\epsilon\}$ is called the {\em free semigroup} on the set $\as$.

\medskip 

Let $u\in\as^n$ and $v\in\as^k$ be arbitrary where $n\leq k$. Then $u$ is called a subword of $v$ if there exist words $v_1$ and $v_2$ such that $v=v_1uv_2$. Here, it is allowed that $v_1$ or $v_2$ is the empty word. Similarly, a subword of $\xi\in\as^\ZM$ is defined. For $n,m\in\ZM$ with $n<m$, $\xi_{[n,m]}$ denotes the finite word $\xi(n)\xi(n+1)\ldots\xi(m)$ of length $m-n+1$. With this notion at hand, a basis for the product topology of $\as^\ZM$, c.f. \cite[Definition~3.7]{Querenburg2001}, is given by the sets 
$$
\os(u,v) \; 
	:= \; \big\{ 
			\xi\in\as^\ZM\;|\; \xi|_{[-n,\ldots,m-1]}=uv 
		\big\} \, ,
	\qquad u\in\as^n\,,\; v\in\as^m \, .
$$
The group $\ZM$ naturally acts on the compact space $\as^\ZM$. Specifically, $\alpha:\ZM\times\as^\ZM\to\as^\ZM$ is defined by
$$
\big(\alpha_n(\xi)\big)(m)\; 	
	= \; \xi(m-n)\,,
		\qquad m,n\in\ZM\,,\; \xi\in\as^\ZM\,,
$$
is a continuous action. Then $(\as^\ZM,\ZM,\alpha)$ defines a topological dynamical system in terms of Definition~\ref{Chap2-Def-DynSyst} below. Let $u:=a_1\ldots a_n\in\as^n$ and $v:=b_1\ldots b_m\in\as^m$ be finite words for $n,m\in\NM$. Denote by $(u|v)^\infty$ the periodic extension of $u|v$ to a two-sided infinite word $\xi\in\as^\ZM$ such that $\xi|_{[-n+k(n+m),(m-1)+k(n+m)]}=uv$ for all $k\in\ZM$, i.e.,
$$
\xi\; 
	:= \; \ldots a_n\, b_1\ldots b_m\, a_1\ldots a_n\, | \, b_1\ldots b_m\, a_1\ldots a_n\, b_1\ldots\;
	\in\as^\ZM \, .
$$
In this notation, it is allowed that one of the words $u,v$ is the empty word $\epsilon$. Then the notation $v^\infty$ is used instead of $(\epsilon|v)^\infty = (v|\epsilon)^\infty$. Such elements are periodic with respect to the $\ZM$-action, i.e., there exists a $p\in\NM$ such that $\alpha_p(\xi)=\xi$. A subset $\Xi\subseteq\as^\ZM$ is called a {\em subshift} if $\Xi$ is closed and $\alpha$-invariant, i.e., $\alpha_n(\Xi)\subseteq\Xi$ for all $n\in\NM$. The set of all subshifts of $\as^\ZM$ is denoted by $\SZ\big(\as^\ZM\big)$. Since $\as^\ZM$ is compact, a subshift $\Xi$ is itself compact with respect to the induced topology and $(\Xi,\ZM,\alpha)$ defines a topological dynamical system. For $\xi\in\as^\ZM$, the {\em orbit $\Orb(\xi)$ of $\xi$} is defined by the set $\{\alpha_n(\xi)\;|\; n\in\ZM\}$. Then the closure $\overline{\Orb(\xi)}$ defines a subshift. Whenever $\xi$ is periodic, the orbit $\Orb(\xi)$ is finite and so it is closed. A subshift $\Xi$ is called {\em minimal} whenever $\Orb(\xi)\subseteq\Xi$ is dense for all $\xi\in\Xi$. Note that there exist subshifts $\Xi:=\overline{\Orb(\xi)}$ for $\xi\in\as^\ZM$ that are not minimal, c.f. the following Example~\ref{Chap2-Ex-OneDefect}.

\begin{example}
\label{Chap2-Ex-OneDefect}
Let $\as:=\{a,b\}$ and define $\xi\in\as^\ZM$ by 
$$
\xi(n)\; 
	:= \; 
		\begin{cases}
			b\,, \quad &n=0\,,\\
			a\,, \quad &n\neq 0\,,
		\end{cases}
		\qquad\quad n\in\ZM\,.
$$
The associated subshift $\Xi:=\overline{\Orb(\xi)}$ is given by the disjoint union $\Orb(\xi)\cup\{ a^\infty \}$. Then $\eta:=a^\infty$ satisfies $\alpha_n(\eta)=\eta$ for all $n\in\ZM$. Thus, the subshift $\Xi$ is not minimal since the letter $b$ never occur in $\eta$. This subshift is called the {\em one-defect}, c.f. Subsection~\ref{Chap7-Sect-OneDefect}.
\end{example}

Let $\as$ be a finite alphabet. Then a map $S:\as\to\as^+$ is called a {\em substitution}. A substitution $S:\as\to\as^+$ is called {\em primitive} if there is an $l_0\in\NM$ such that the letter $b$ occurs in the word $S^{l_0}(a)$ for every pair $a,b\in\as$. A substitution extends to a map $S:\as^+\to\as^+$ and to a map $S:\as^\ZM\to\as^\ZM$ by acting separately on each letter, i.e., 
\begin{gather*}
\begin{aligned}
S(v) \;
	&:= \; S(a_1)\ldots S(a_n)\,,
	\qquad 
	&v:=a_1\ldots a_n\in\as^+\,,\\
S(\xi) \; 
	&:= \;\ldots S\big(\xi(-2)\big)S\big(\xi(-1)\big) \big| S\big(\xi(0)\big)S\big(\xi(1)\big)S\big(\xi(2)\big)\ldots\,,
	\qquad
	&\xi\in\as^\ZM\,.
\end{aligned}
\end{gather*}

For $k\in\NM$, a $\xi\in\as^\ZM$ is called {\em $k$-periodic with respect to the substitution $S$} if $S^k(\xi)=\xi$. Whenever $\xi\in\as^\ZM$ is $1$-periodic with respect to a substitution $S$, then $\xi$ is called {\em fixed point of the substitution $S$}. It is well-known that $\lim_{n\to\infty}|S^n(a)|=\infty$ for each letter $a\in\as$ if $S$ is a primitive substitution and $|\as|\geq 2$. This implies that there is at least one $k$-periodic point for a primitive substitution $S$, c.f. \cite[Proposition~V.1]{Queffelec87} and Proposition~\ref{Chap6-Prop-ExKPeriodPoint}. Furthermore, it turns out that primitive substitutions give rise to minimal subshifts, c.f. \cite{Queffelec87,Fogg02,Queffelec10} or Proposition~\ref{Chap6-Prop-PrimSubstSubshZdMin}.

\medskip

Let $(\Xi,\ZM,\alpha)$ be a minimal subshift. Consider continuous maps $p:\Xi\to\RM\setminus\{0\}$ and $q:\Xi\to\RM$ that take only finitely many values. Such functions are called {\em pattern equivariant} since they depend only on the local structure of $\xi$, c.f. \cite{KePu00,Kel03} or Section~\ref{Chap2-Sect-PattEqSchrOp}. More precisely, there is an $N\in\NM$ such that $p(\xi)=p(\eta)$ whenever $\xi_{[-N,N]}=\eta_{[-N,N]}$, c.f. Definition~\ref{Chap2-Def-PattEqFunc} and Proposition~\ref{Chap2-Prop-CharPattEqFunc}. Then the family of {\em Jacobi operators} $J_\xi:\ell^2(\ZM)\to\ell^2(\ZM)\,,\; \xi\in\Xi\,,$ associated with a subshift $\Xi$ and the maps $p$ and $q$ is defined by
$$
\big(J_\xi\psi\big)(n)\;
	= \; p\big(\alpha_{-n}(\xi)\big) \cdot \psi(n-1) + p\big(\alpha_{-(n+1)}(\xi) \big)\cdot \psi(n+1) + q\big(\alpha_{-n}(\xi)\big)\cdot\psi(n)
$$
for $\xi\in\Xi\,,\; \psi\in\ell^2(\ZM)$ and $n\in\ZM$. Whenever the map $p$ is identically to one, such operators are called {\em Schr\"odinger operators}. These operators are bounded and self-adjoint, c.f. Lemma~\ref{Chap2-Lem-PatEqSchrOp}. Furthermore, the family $J_\Xi:=(J_\xi)_{\xi\in\Xi}$ is a generalized Schr\"odinger operator of finite range in terms of Definition~\ref{Chap2-Def-SchrodingerOperator}, c.f. Theorem~\ref{Chap2-Theo-PESchrOpFinRang}. Then the spectrum $\sigma(J_\Xi)$ of the family $J_\Xi$ is equal to $\bigcup_{\xi\in\Xi}\sigma(J_\xi)$, c.f. Proposition~\ref{Chap2-Prop-CovFamOp-Spect}. It is well-known that, for this family of self-adjoint operators $J_\Xi$, the spectrum is constant whenever $\Xi$ is minimal, i.e., $\sigma(J_\Xi)=\sigma(J_\xi)$ holds for all $\xi\in\Xi$. This result is proven in \cite{CyconFroeseKirschSimon87,BIST89,Jit95,Len99}. In \cite{BeLeMaCh14,BeLeMaCh16}, this result was extended to every family of operators over a dynamical system $(X,G,\alpha)$ provided that $X$ is minimal and the operators satisfy a weak continuity condition. Generalized Schr\"odinger operators satisfy this continuity condition and so the equations $\sigma(J_{\Xi})=\sigma_{ess}(J_\xi)=\sigma(J_\xi)$ hold for all $\xi\in\Xi$ if $\Xi$ is minimal, c.f. Theorem~\ref{Chap2-Theo-ConstSpectrMinimal}. Here $\sigma_{ess}(J_\xi)$ denotes the essential part of the spectrum $\sigma(J_\xi)$ for $\xi\in\Xi$.

\medskip

Self-adjoint random operators arise in the quantum mechanical treatment of disordered solids. Their study has been a key focus of mathematical physics in the last four decades. Indeed, for an impressive number of (classes of) examples, explicit spectral features (such as pure point spectrum or purely singular continuous spectrum or Cantor spectrum) are proven, see e.g. the surveys and monographs \cite{CarmonaLacroix90,PasturFigotin92,Dam00Gord,Stollmann01,Len02,BePo13,DaEmGo15}.

\medskip

One fundamental task is to study the relation between the geometric and combinatorial properties of the subshift $\Xi$ and the spectrum $\sigma(J_\Xi)$ of a given family of Schr\"odinger or Jacobi operators $J_\Xi$. This thesis provides the answer how the spectrum varies if the subshift $\Xi$ is deformed. More precisely, the continuity of the map 
$$
\SZ\big(\as^\ZM\big)\ni\Xi\; 
	\mapsto\;\sigma(J_\Xi)\in\ks(\RM)
$$
is verified. This is based on joint work with {\sc J. Bellissard} and {\sc G. de Nittis} \cite{BeBeNi16}. In detail, it is shown that a family of subshifts varies continuously with respect to the Hausdorff-topology on $\cs(\as^\ZM)$ if and only if the family of spectra associated with all Schr\"odinger or Jacobi operator varies continuously with respect to the Hausdorff metric on $\ks(\CM)$, c.f. Corollary~\ref{Chap4-Cor-CharSubshiftConvSpectrZM}. Actually, this result is proven in a much more generality. More precisely, the setting is a topological dynamical system $(X,G,\alpha)$ where $G$ is required to be an (amenable,) discrete, countable group. 

\medskip

The proof of this characterization is based on continuous fields of $C^\ast$-algebras motivated by Theorem~\ref{Chap3-Theo-PContEquivContSpectNormal}. In detail, a topological dynamical system $(X,G,\alpha)$ naturally induces a topological groupoid $\Gamma:= X\rtimes_\alpha G$. Associated with $\Gamma$, the reduced $C^\ast$-algebra is defined that contains the generalized Schr\"odinger  operators. This connection is used in Chapter~\ref{Chap4-ToolContBehavSpectr} to prove the continuous behavior of the map $\SZ\big(\as^\ZM\big)\ni\Xi\mapsto\sigma(J_\Xi)\in\ks(\RM)$ by investigating the concept of universal dynamical system and universal groupoid, c.f. Corollary~\ref{Chap4-Cor-ContSpectrContSectUnivGroup}. In view of that, the framework is investigated in this chapter.

\section{Topological dynamical systems}
\label{Chap2-Sect-DynSystGroupoid}

This section deals with the basic properties of dynamical systems. The topological space of dynamical subsystems is introduced. Furthermore, basic concepts of periodicity and minimality are established. Additionally, the set of periodically approximable dynamical subsystems is defined which plays an important role for the further considerations.

\medskip

Recall the notion of a topological group. More precisely, a {\em topological group $G$} is defined by a topological space $(G,\tau)$ equipped with a continuous map $\cdot:G\times G\to G$ such that the map $\cdot$ is associative, there exists an $e\in G$ satisfying $e\cdot g=g\cdot e = g\,,\; g\in G\,,$ and there is, for each $g\in G$, an element $g^{-1}\in G$ such that $g\cdot g^{-1}=g^{-1}\cdot g = e$. The element $e\in G$ is called {\em neutral element} and is unique. Furthermore, $g^{-1}\in G$ is called the {\em inverse of $g\in G$} and is also unique. Note that $G\times G$ is equipped with the product topology.

\begin{definition}[Dynamical system]
\label{Chap2-Def-DynSyst}
Let $G$ be a second-countable, locally compact, Hausdorff group with neutral element $e\in G$. A triple $(X,G,\alpha)$ is called a {\em topological dynamical system} if $X$ is a non-empty, second-countable, compact, Hausdorff space and $\alpha$ is a continuous action of $G$ on $X$. In detail, 
$$
\alpha:G\times X\to X,
	\qquad (g,x)\mapsto\alpha_g(x) \, ,
$$ 
is a continuous map such that:
\begin{itemize}
\item[(i)] $\alpha_e(x)=x$ for all $x\in X$,
\item[(ii)] $\alpha_{gh}(x)=\alpha_{g}(\alpha_{h}(x))$ for each $g,h\in G$ and $x\in X$.
\end{itemize}
\end{definition}

For fixed $g\in G$, the map $\alpha_g:X\to X$ defines a homeomorphism. Note that if $G$ is abelian the set $\{\alpha_g\;|\; g\in G\}$ is a commutative subset of the group of homeomorphisms on $X$.

\begin{definition}[Dynamical subsystem]
\label{Chap2-Def-SubdynSyst}
Let $(X,G,\alpha)$ be a topological dynamical system. A subset $Y\subseteq X$ is called {\em $G$-invariant} if $\{\alpha_g(y)\;|\; g\in G,\; y\in Y\}\subseteq Y$. A triple $(Y,G,\alpha_Y)$ is called {\em dynamical subsystem of $X$} if $Y\subseteq X$ is a $G$-invariant, closed, non-empty subset where $\alpha_Y:G\times Y\to Y$ is the restriction of $\alpha$ to $Y$.
\end{definition}

\begin{lemma}
\label{Chap2-Lem-DynSubsyst}
Let $(X,G,\alpha)$ be a topological dynamical system. Then, a dynamical subsystem $(Y,G,\alpha_Y)$ is a dynamical system and $\{\alpha_g(y)\;|\; g\in G,\; y\in Y\}= Y$.
\end{lemma}

\begin{proof}
Let $(Y,G,\alpha_Y)$ be a dynamical subsystem of $X$. Since $Y$ is closed, it is compact with respect to the subspace topology. Furthermore, $Y$ is second-countable and Hausdorff as a subspace of $X$. By construction, the map $\alpha_Y:G\times Y\to Y$ fulfills the requirements of Definition~\ref{Chap2-Def-DynSyst} and so $(Y,G,\alpha_Y)$ defines a dynamical system.

\vspace{.1cm}

Consider a $y_0\in Y$ and $g\in G$. Then $z:=\alpha_{g^{-1}}(y_0)$ is an element of $Y$ as $Y$  is $G$-invariant. By Definition~\ref{Chap2-Def-DynSyst}~(ii), the image $\alpha_g(z)$ is equal to $y_0$. Thus, the equation $\{\alpha_g(y)\;|\; g\in G,\; y\in Y\}= Y$ follows.
\end{proof}

\begin{definition}[Space of dynamical systems]
\label{Chap2-Def-SpaDynSyst}
Let $(X,G,\alpha)$ be a topological dynamical system. The set $\text{\gls{IGX}}\subseteq \cs(X)$ of all closed, $G$-invariant, non-empty subsets of $X$ endowed with the Hausdorff-topology induced by $\cs(X)$ is called the {\em space of dynamical systems associated with $(X,G,\alpha)$}.
\end{definition}

Note that every closed subset $Y\in\cs(X)$ is compact since $X$ is compact. Then each element $Y\in\SG(X)$ induces a dynamical subsystem $(Y,G,\alpha_Y)$. Whenever there is no confusion a dynamical subsystem is identified with $Y\in\SG(X)$. Specifically, a $Y\in\SG(X)$ is called a dynamical system. For sake of simplification, the notation $\alpha$ is used also for the restriction $\alpha_Y:G\times Y\to Y$ to the dynamical subsystem $Y\in\SG(X)$. The following assertion follows by abstract arguments of \cite{Fel62,Bee93}.

\begin{proposition}
\label{Chap2-Prop-SpaDynSyst}
Let $(X,G,\alpha)$ be a topological dynamical system. The space $\SG(X)$ is second-countable, compact and Hausdorff with respect to the induced Hausdorff-topology.
\end{proposition}

\begin{proof}
The space $X$ is second-countable, compact and Hausdorff. Hence, $X$ is a normal space, c.f. Proposition~\ref{App1-Prop-LocCompHausSecCountImplNormal}. Thus, the space $\cs(X)$ endowed with the Hausdorff-topology is second-countable, compact and Hausdorff, c.f. \cite[Theorem~3.4]{Bee93} and \cite[Theo\-rem~1]{Fel62}. Consequently, it suffices to prove that $\SG(X)$ is a closed subset of $\cs(X)$. Let $Y_n\in\SG(X),\; n\in\NM,$ be a sequence converging to $Y\in\cs(X)$ in the Hausdorff-topology. The proof is organized as follows: 
\begin{itemize}
\item[(i)] Let $(y_n)_{n\in\NM}$ be a convergent sequence to $y\in X$ such that $y_n\in Y_n,\; n\in\NM$. Then $y$ is contained in $Y$. 
\item[(ii)] For each $y\in Y$, there is a sequence $y_n\in Y_n,\; n\in\NM$, such that $\lim_{n\to\infty} y_n=y$. 
\item[(iii)] For $g\in G$ and $y\in Y$, $\alpha_g(y)$ is contained $Y$.
\end{itemize}

Assertion (iii) is proved by using (i) and (ii). It immediately leads to $Y\in\SG(X)$.

\vspace{.1cm}

(i): Let $(y_n)_{n\in\NM}$ be such a sequence converging to $y\in X$. Since $X$ is normal, there is for every neighborhood $U$ of $y$ an open set $V$ such that $y\in V\subseteq \overline{V}\subseteq U$. By convergence of the sequence $(y_n)_{n\in\NM}$, there is an $n(V)\in\NM$ such that $y_n\in V$ for all $n\geq n(V)$. Hence, the intersection $Y_n\cap V$ is non-empty for $n\geq n(V)$. Since the sequence $(Y_n)_{n\in\NM}$ converges to $Y$ in the Hausdorff-topology, the set $Y$ has to intersect $U$ which can be checked as follows. Assume that $Y\cap U=\emptyset$ implying $Y\cap\overline{V}=\emptyset$. Thus, $\us(\overline{V},X)\subsetneq\cs(X)$ defines an open neighborhood of $Y$. Consequently, there is an $n_0\in\NM$ such that $Y_n\cap \overline{V}=\emptyset$ for all $n\geq n_0$ by the convergence of $(Y_n)_{n\in\NM}$ to $Y\in\cs(X)$. This contradicts the previous considerations that $y_n\in V$ for $n\geq n(V)$. Altogether, for each neighborhood $U$ of $y$, the intersection $U\cap Y$ is non-empty. Hence, $y$ is an element of $Y$ as $Y$ is closed.

\vspace{.1cm}

(ii): For every neighborhood $U$ of $y$, there is an $n(U)\in\NM$ such that $Y_n\in\us(\emptyset,U)$ for all $n\geq n(U)$. More precisely, for all $n\geq n(U)$, the intersection $Y_n \cap U$ is non-empty. Since $X$ is second-countable, a countable neighborhood base can be chosen for $y$. Thus, a sequence $(y_n)_{n\in\NM}$ with desired properties is constructed.

\vspace{.1cm}

(iii): Let $g\in G$ and $y\in Y$. By (ii) there is a sequence $(y_n)_{n\in\NM}$ such that $\lim_{n\to\infty} y_n=y$ and $y_n\in Y_n$ for all $n\in\NM$. Since $Y_n\in\SG(X)$, the element $\alpha_g(y_n)$ is an element of $Y_n$. As $\alpha_g:X\to X$ is continuous, the limit $\lim_{n\to\infty} \alpha_g(y_n)$ exists and is equal to $\alpha_g(y)$. Thus, (i) implies that $\alpha_g(y)\in Y$ finishing the proof.
\end{proof}

\medskip

Let $(X,G,\alpha)$ be a dynamical system. For $x\in X$, the {\em associated dynamical subsystem} is defined by 
$$
X[x]:=\overline{\Orb(x)} \, ,
$$
where the set $\text{\gls{Orb}}:=\{\alpha_g(x)\;|\; g\in G\}$ is called the {\em orbit of $x$}. In general, not each dynamical system arises by the closure of an orbit $\Orb(x)$. For instance, the union $Y$ of two disjoint dynamical subsystems of $X$ (if such exist) form a dynamical subsystem while there does not exist an $x\in X$ with $Y=\overline{\Orb(x)}$. A dynamical systems $(Y,G,\alpha)$ is called\label{Chap2-Def-TopolTransit} {\em topologically transitive} if $Y=\overline{\Orb(x)}$ for an $x\in Y$.

\medskip

Let $G$ be the discrete group $\ZM$. In this case, it is natural to call an $x\in X$ periodic whenever there exists a $p\in\ZM$ such that $\alpha_p(x)=x$, c.f. Section~\ref{Chap2-Sect-ExampleAsZMSchrOp}. This is equivalent to saying that the related orbit $\Orb(x)$ is finite. However, for a more general group like $\ZM^d$, an element $x\in X$ can be periodic in a specific direction while the orbit is not finite at all. So it is necessary to distinguish these cases, in general. An element $x\in X$ is called {\em weakly periodic} whenever there exists a $g\in G$ such that $g\neq e$ and $\alpha_g(x)=x$. If, additionally, the orbit $\Orb(x)$ is finite, $x\in X$ is called {\em strongly periodic}. Note that strongly periodic elements are also called periodic \cite{Fio09,SiCo12}. However, the distinction between weakly and strongly periodic elements is standard in the literature \cite{ChereneLouis78,PiantadosiThesis06,Pia08}. By the previous considerations, it is clear that the notions of weakly and strongly periodicity coincide if $G=\ZM$.

\medskip

The notion of periodic elements is translated into the notion of subgroups of $G$. More precisely, the {\em stabilizer} \gls{Stab} of $x\in X$ is the set $\{g\in G \;|\; \alpha_g(x)=x\}$. According to Definition~\ref{Chap2-Def-DynSyst}~(i), the stabilizer $\Stab_\alpha(x)$ contains the neutral element $e\in G$ for each $x\in X$. Thus, the stabilizer contains at least one element. Furthermore, Definition~\ref{Chap2-Def-DynSyst} (ii) implies that $\Stab_\alpha(x)$ is a subgroup of $G$ for every $x\in X$. Then $x\in X$ is weakly periodic if and only if $\sharp\Stab_\alpha(x)>1$ where $\sharp A$ denotes the cardinality of a set $A$. According to the Orbit-Stabilizer Theorem, see e.g. \cite[Theorem~4.1]{PiantadosiThesis06}, an element $x\in X$ is strongly periodic if and only if the stabilizer $\Stab_\alpha(x)$ has finite index as a subgroup of a discrete countable group $G$.

\medskip

An element $x\in X$ is called {\em non-periodic} if it is not weakly periodic, i.e., the stabilizer $\Stab_\alpha(x)$ contains only the neutral element $e$. Note that in this case, the orbit $\Orb(x)$ is infinite whenever $G$ is infinite. The notion of periodicity transposes to dynamical systems which is made more precise in the following definition.

\begin{definition}
\label{Chap2-Def-AperSubs}
A topological dynamical system $(X,G,\alpha)$ is called
\begin{itemize}
\item {\em minimal} if, for all $x \in X$, the orbit $\Orb(x)$ is dense in $X$;
\item {\em weakly periodic} if each $x\in X$ is weakly periodic;
\item {\em strongly periodic} if $X$ is minimal and each $x\in X$ is strongly periodic;
\item {\em aperiodic} if there is a non-periodic element $x\in X$;
\item {\em completely aperiodic} if all $x\in X$ are non-periodic.
\end{itemize}
\end{definition}

\begin{remark}
\label{Chap2-Def-PeriodStabAbelian}
Since this work is basically concerned with the study of Schr\"odinger operators defined by dynamical systems, the class of strongly periodic dynamical subsystems is of particular interest thanks to the Floquet-Bloch theory in the abelian case. Specifically, the spectral properties of Schr\"odinger operators with periodic potentials are studied with the so called Floquet-Bloch theory \cite{Flo1883,BlochThesis29}, c.f. discussion in Section~\ref{Chap4-Sect-ApprPerAppr}. It is necessary for this theory that the corresponding stabilizer group is abelian since the Fourier transformation of groups is used. In view of that the Floquet-Bloch theory is only applicable for Schr\"odinger operators over a dynamical system $(Y,G,\alpha)$ where $Y$ is strongly periodic and the stabilizer groups $\Stab_\alpha(y)\,,\; y\in Y\,,$ are abelian. Whenever the group $G$ itself is abelian, it suffices that $Y$ is strongly periodic since the stabilizer groups are automatically abelian as subgroups of $G$. Altogether, the notion of strongly periodic subshifts has to be handled with care if the group $G$ is not abelian which is included in our situation. Note that the systems studied in Chapter~\ref{Chap5-OneDimCase}, Chapter~\ref{Chap6-HigherDimPerAppr} and Chapter~\ref{Chap7-Examples} are defined for dynamical systems where the group $G$ is abelian.
\end{remark}

\begin{proposition}
\label{Chap2-Prop-AperSubs}
Let $(X,G,\alpha)$ be a topological dynamical system where $G$ is not finite. Then the following assertions hold.
\begin{itemize}
\item[(a)] If $x\in X$ is strongly periodic then $\Orb(x)\in\SG(X)$ is a minimal dynamical subsystem, namely $\Orb(x)$ is a strongly periodic dynamical system.
\item[(b)] The dynamical system $Y\in\SG(X)$ is strongly periodic if and only if $Y$ is minimal and finite if and only if $Y$ is finite and $Y=\Orb(y)\,,\; y\in Y$.
\item[(c)] If the dynamical system $Y\in\SG(X)$ is minimal and aperiodic then it does not contain a strongly periodic element. If, additionally, $G$ is abelian, then $Y$ is completely aperiodic.
\end{itemize}
\end{proposition}

\begin{proof}
(a): Since $X$ is a Hausdorff space, a finite subset of $X$ is always closed. Thus, for a strongly periodic element $x\in X$, the orbit $\Orb(x)$ defines a dynamical subsystem that is minimal by construction.

\vspace{.1cm}

(b): This follows by (a) and the definition of strongly periodic elements.

\vspace{.1cm}

(c): Assume that $Y$ contains a strongly periodic element $x$, i.e., the orbit $\Orb(x)$ is finite. Due to minimality, it follows that $\Orb(x)$ is a dense subset of $Y$. Since $\Orb(x)$ is closed, the equation $\Orb(x)=Y$ follows. This contradicts the fact that $Y$ contains a non-periodic element.

\vspace{.1cm}

The second part is also proven by contradiction. More specifically, let $G$ be abelian and assume there is a weakly periodic element $x\in Y$. Since $Y$ is aperiodic, there is a non-periodic $y\in Y$. Then there exists a sequence $g_n\in G\,,\; n\in\NM\,,$ such that $\lim_{n\to\infty}\alpha_{g_n}(x)=y$ by minimality of $Y$. Let $g\neq e$ be such that $\alpha_g(x)=x$. The continuity of $\alpha$ implies
$$
\alpha_g(y) \; 
	=\; \lim_{n\to\infty} \alpha_{g}\big(\alpha_{g_n}(x)\big) \;
	\underset{\text{Def.~\ref{Chap2-Def-DynSyst}~(ii)}}{\overset{G\text{ abelian}}{=}}\; \lim_{n\to\infty} \alpha_{g_n}\big(\alpha_g(x)\big) \;
	=\; \lim_{n\to\infty} \alpha_{g_n}(x) \;
	=\; y\, .
$$
This is a contradiction as $y$ is non-periodic and $g\neq e$.
\end{proof}

\medskip

For convenience of the reader, some examples are presented in the following. Consider the discrete, abelian group $G=\ZM$ and a finite alphabet $\as$. Define the product space $\as^\ZM$, c.f. Section~\ref{Chap2-Sect-ExampleAsZMSchrOp} and Section~\ref{Chap2-Ssect-SpaceSubshifts}. Then $\alpha:\ZM\times\as^\ZM\to\as^\ZM$ defined by $\big(\alpha_m(\xi)\big)(n):=\xi(n-m)$ for $n,m\in\ZM$ is a continuous action of $\ZM$ on $\as^\ZM$.

\begin{example}[Non-periodic]
\label{Chap2-Ex-NonPeriodic-Fibonacci}
One of the most prominent example is the Fibonacci sequence. This sequence was already known in the year 100 by {\sc Nikomachos von Gerasa}, c.f. \cite[Footnote~86]{Lan98}. In the year 1202, the mathematician {\sc Leonardo de Pisa} (called Fibonacci) theoretically studied the evolution of rabbits by this sequence in his famous \textit{Liber abbaci}, c.f. \cite[Chapter~XII]{Boncompagni57}. The combinatorial Fibonacci sequence is constructed by a substitution. Specifically, the alphabet $\as$ has two letters, say, $a$ and $b$. Then the map $S:\as\to\as^+\,,\; a\to ab\,,\,b\to a\,,$ defines a primitive substitution, c.f. Definition~\ref{Chap6-Def-BlockPrimSubst}. Along the lines of Definition~\ref{Chap6-Def-AssDictSubshSubstZd} and Proposition~\ref{Chap6-Prop-PrimSubstSubshZdMin}, a minimal subshift $\Xi_S$ is constructed. This subshift is called the {\em Fibonacci subshift}. It is well-known that the subword complexity of the Fibonacci subshift satisfies $\comp(n)=n+1$ for all $n\in\NM$ where $\comp(n)$ is the number of subwords of the Fibonacci subshift of length $n\in\NM$, see e.g. \cite[Corollary~5.4.10]{Fogg02}. The reader is referred to \cite{LindMarcus95} and Section~\ref{Chap5-Ssect-SubwordComplexity} for more details on the subword complexity. Clearly, the subword complexity of a strongly periodic subshift of $\as^\ZM$ is uniformly bounded in $n\in\NM$, i.e., there is a constant $C>0$ such that $\comp(n)\leq C$ for all $n\in\NM$, c.f. \cite{MoHe38,MoHe40,CoHe73} and \cite{LuPl87,Lothaire02} for two-sided infinite sequences. Thus, the Fibonacci subshift is aperiodic since, for $G=\ZM$, the notions of weakly and strongly periodic coincide. Then Proposition~\ref{Chap2-Prop-AperSubs}~(c) implies that $\Xi_S$ is completely aperiodic as $\ZM$ is an abelian group and $\Xi_S$ is minimal by Proposition~\ref{Chap6-Prop-PrimSubstSubshZdMin}. Thus each $\xi\in\Xi_S$ is a non-periodic element.
\end{example}

\begin{example}[Weakly- and strongly periodic]
\label{Chap2-Ex-WeaklyStronglyPeriodic}
Let $\as:=\{a,b\}$. Define the compact space $X:=\as^\ZM\times\as^\ZM$ endowed with the product topology. Then the map $\beta:\ZM^2\times X\to X\,,\; \big( (n,m),(\xi,\eta)\big)\mapsto\big(\alpha_n(\xi),\alpha_m(\eta)\big)$ defines a continuous action of $\ZM^2$ on the space $X$. Thus, $(X,\ZM^2,\beta)$ defines a topological dynamical system. Let 
$$
\xi\;
	:= \; \ldots a\, b\, a\, b\, a\, b\, |\, a\, b\, a\, b\, a\, b\, a \ldots\in\as^\ZM
$$ 
be the strongly periodic element defined by $(b|a)^\infty$. Furthermore, let $\eta\in\as^\ZM$ be an element of the Fibonacci sequence. Then $(\xi,\eta)$ is an element of $X$ and by construction it is weakly periodic as $\alpha_{(2,0)}\big((\xi,\eta)\big)=(\xi,\eta)$. However, $(\xi,\eta)$ it is not strongly periodic since $\eta$ is non-periodic in $\as^\ZM$. On the other hand, the element $(\xi,\xi)$ is strongly periodic since the orbit $\Orb(\xi,\xi)$ is equal to $\{(\xi,\xi),\, (\alpha_1(\xi),\xi),\, (\xi,\alpha_1(\xi)),\, (\alpha_1(\xi),\alpha_1(\xi))\}$. Thus, the orbit is finite.
\end{example}

The section finishes with the notion of periodically approximable dynamical subsystems. It is discussed in Chapter~\ref{Chap4-ToolContBehavSpectr} that generalized Schr\"odinger operators of periodically approxi\-mable dynamical systems are approximated by periodic ones such that the spectra converge, c.f. Theorem~\ref{Chap4-Theo-PeriodicApproximations}. This is useful as non-periodic generalized Schr\"odinger operators are difficult to analyze whereas there are tools for periodic operators, c.f. Remark~\ref{Chap2-Def-PeriodStabAbelian}.

\begin{definition}[Periodically approximable]
\label{Chap2-Def-PeriodApprox}
Let $(X,G,\alpha)$ be a topological dynamical system. 
\begin{itemize}
\item The set of strongly periodic dynamical subsystem of $X$ is denoted by $\text{\gls{SPX}}\subseteq\SG(X)$.
\item A dynamical subsystem $Y\in\SG(X)$ is called {\em periodically approximable} if $Y\in\overline{\SP(X)}$ where the closure is taken with respect to the Hausdorff-topology in $\cs(X)$. Then the set of {\em periodically approximable} dynamical subsystems $\overline{\SP(X)}$ is denoted by \gls{PAX}.
\end{itemize}
\end{definition}

If $(X,G,\alpha)$ is a minimal dynamical system, it follows that $\SG(X)$ only contains the element $X$. Thus, the set $\PA(X)$ is empty if, for instance, $(X,G,\alpha)$ is minimal and aperiodic. On the other hand, the set $\PA(X)$ may contain a lot of elements, c.f. Proposition~\ref{Chap5-Prop-MinAllPath} and the following statement. 

\begin{theorem}[\cite{BeBeNi16}]
\label{Chap2-Theo-FinOrbDenseInvMetr}
Let $(X,G,\alpha)$ be a dynamical system satisfying that 
\begin{itemize}
\item[(a)] $\mathcal{F}:=\{x\in X\;|\; \Orb(x) \text{ finite}\}\subseteq X$ is dense;
\item[(b)] there exists a metric $d:X\times X\to[0,\infty)$ such that $(X,d)$ is a complete metric space and $d\circ\alpha=d$, i.e., $d\big(\alpha_g(x),\alpha_g(y)\big)=d(x,y)$ holds for all $x,y\in X$ and $g\in G$.
\end{itemize}
Then the set $\{Y\in\SG(X)\;|\; Y\text{ finite}\}\subseteq\SG(X)$ is dense with respect to the Hausdorff metric.
\end{theorem}

\begin{proof}
It suffices to prove for $Z\in\SG(X)$ and $\varepsilon>0$ that there exists a finite $Y\in\SG(X)$ such that 
$$
d_H(Y,Z)\;	
	:= \; \max\big\{
		\sup_{y\in Y}\inf_{z\in Z} d(y,z)\,,\;
		\sup_{z\in Z}\inf_{y\in Y} d(y,z)
	\big\}\;
	< \; \varepsilon\,.
$$

\vspace{.1cm}

Let $Z\in\SG(X)$ and $\varepsilon>0$. Then there are $x_1,\ldots,x_n\in X$ such that $d_H\big(\{x_1,\ldots, x_n\},Z\big)<\varepsilon/2$ by compactness of $X$. Due to assumption (a), there exist $y_1,\ldots, y_n\in\mathcal{F}$ satisfying $d(x_j,y_j)<\varepsilon/2$ for $1\leq j\leq n$. Consequently, the estimate
$$
d_H\big(\{y_1,\ldots,y_n\},Z\big) \;
	\leq \; d_H\big(\{y_1,\ldots,y_n\},\{x_1,\ldots,x_n\}\big) + d_H\big(\{x_1,\ldots,x_n\},Z\big) \;
	< \; \varepsilon
$$
is derived. Since $Z\in\SG(X)$ and $d\circ\alpha=d$, the estimates
$$
d_H\big(\{\alpha_g(y_1),\ldots,\alpha_g(y_n)\},Z\big) \;
	= \; d_H\big(\{\alpha_g(y_1),\ldots,\alpha_g(y_n)\},\alpha_g(Z)\big) \;
	< \; \varepsilon\,,
	\qquad
	g\in G\,,
$$
is deduced. Hence, $d_H(Y,Z)<\varepsilon$ follows for $Y:=\{\alpha_g(y_j)\;|\; g\in G\,,\; 1\leq j\leq n\}$. Clearly, $Y$ is a finite set contained in $\SG(X)$ since $y_j\in\mathcal{F}$ for $1\leq j\leq n$.
\end{proof}

\begin{remark}
\label{Chap2-Rem-FinOrbDenseInvMetr}
The assumption (a) is not too restrictive, c.f. \cite{Lig03,SiCo12}. On the other hand, the existence of an invariant metric is a strong condition. For $(\as^\ZM,\ZM,\alpha)$, this theorem shows that there does not exist an invariant metric on $\as^\ZM$ since there exists a subshift that is not periodically approximable, c.f. Example~\ref{Chap5-Ex-deBruijnEventuallyNotStrongConn}. If $(X,G,\alpha)$ is a dynamical system such that $G$ is compact, then there exists an invariant metric on $X$. More precisely, the space $X$ is metrizable and so there is a metric $d$ on $X$. Then $\tilde{d}:X\times X\to[0,\infty)$ defined by
$$
\tilde{d}(x,y) \; 
	:= \; \int_G d\big(\alpha_g(x),\alpha_g(y)\big)\; d\lambda(g)\,,
	\quad
	x,y\in X\,,
$$
is an invariant metric on $X$ where $\lambda$ is the left-invariant Haar measure of the group $G$, c.f. \cite{HaWa01}.
\end{remark}

\section{Dictionaries and symbolic dynamical systems}
\label{Chap2-Sect-SymbDynSyst}

In this section the class of symbolic dynamical systems is studied. First, the notion of dictionaries is extended to $(\as^G,G,\alpha)$ while it is important to notice that the notion of dictionary is defined without a given element of $\as^G$. Then the local pattern topology on the set of dictionaries is investigated. Basic properties of this topology are verified. Finally, it is shown that the space of subshifts is homeomorphic to the space of dictionaries.

\medskip

For the convenience of the reader, a short discussion of dictionaries associated with a finite set $\as$ and the group $\ZM$ is considered before discussing the general case. The reader is referred to \cite{LindMarcus95,BeBeNi16} for a more detailed discussion for two-sided infinite sequences. Recall the notions introduced in Section~\ref{Chap2-Sect-ExampleAsZMSchrOp}. Let $\as$ be a finite set equipped with the discrete topology. Elements of the product space $\as^\ZM=\prod_{n\in\ZM}\as$ are called {\em two-sided infinite sequences}. Together with the action $\alpha:\ZM\times\as^\ZM\to\as^\ZM\,,\; \alpha_m(\xi):=\xi(\cdot-m)\,,$ and the group $\ZM$, it forms a dynamical system. A finite word $u\in\as^n$ is called a {\em subword} of $\xi\in\as^\ZM$ if there exists an $n_0\in\ZM$ such that $\xi|_{[n_0,n_0+n-1]}=u$. Given a two-sided infinite word $\xi\in\as^\ZM$ the set $\ws(\xi)$ of subwords associated with $\xi$ is called a dictionary. It is not difficult to check that $\ws:=\ws(\xi)$ is non-empty and satisfies
\begin{description}
\item[(DZ1)\label{(DZ1)}] for all $u\in\ws$ each subword of $u$ belongs to $\ws$;
\item[(DZ2)\label{(DZ2)}] for all $u\in\ws$ there are $a,b\in\as$ such that $aub\in\ws$.
\end{description}
Now, let $\ws$ be a non-empty set of finite words satisfying \nameref{(DZ1)} and \nameref{(DZ2)} then $\ws$ is called a dictionary. Note that this definition is independent of a $\xi\in\as^\ZM$. It turns out that the set of dictionaries has a one-to-one correspondence to the set of subshifts, i.e., closed, $\ZM$-invariant subsets of $\as^\ZM$, c.f. \cite[Proposition 1.3.4.]{LindMarcus95} and Theorem~\ref{Chap2-Theo-Shift+DictSpace} for the general case. In addition, both spaces are naturally equipped with a topology so that the spaces are homeomorphic which is proven here, c.f. Theorem~\ref{Chap2-Theo-Shift+DictSpace}. 

\medskip

This idea is extended to the more general case of a countable, discrete group $G$. Since $G$ does not longer have extensions to the right and left, a more detailed study is necessary.

\begin{remark}
\label{Chap2-Rem-Dictionary}
The notion of a dictionary in $\as^\ZM$, also called language, is not new in the literature, see e.g. \cite{LindMarcus95}. In some reference a (finite) subset $\ws$ of words is also called a dictionary without assuming any additional requirements on $\ws$, c.f. \cite{More05}. The notion of a dictionary in $\as^G$, according to Definition~\ref{Chap2-Def-Dictionary}, is new in the literature to the best of the authors knowledge. Here $G$ denotes a countable, discrete group. It is the generalization of the case where $G$ is equal to $\ZM$. In this section, the local pattern topology (Definition~\ref{Chap2-Def-TopSpaceDict}) on the set of dictionaries in $\as^G$ is introduced and analyzed. This notion was not present in the literature before, even in the case $G=\ZM$.
\end{remark}

\subsection{The space of dictionaries}
\label{Chap2-Ssect-SpaceDictionaries}

Let $G$ be a discrete, countable group. Denote by $\ks(G)$ the set of non-empty, finite subsets of $G$. By discreteness of $G$, a subset $K\subseteq G$ is finite if and only if it is compact. Two finite subsets $F,K\in\ks(G)$ are called {\em $G$-equivalent $(F\sim_G K)$} if there is a $g\in G$ such that $gF=K$. It immediately follows that $\sim_G$ defines an equivalence relation, i.e., it is reflexive, symmetric and transitive. Then \gls{ke} denotes the quotient $\ks(G)/\!\!\sim_G$.

\medskip

A finite set $\as$ is called an {\em alphabet} equipped with the discrete topology. For every $K\in\ks(G)$, define the Cartesian product $\as^{K}:=\{u:K\to\as\}$. Since $K\in\ks(G)$ is a finite set, $\as^K$ is finite as well.

\begin{definition}[Pattern]
\label{Chap2-Def-Pattern}
Let $F,K\in\ks(G)$ be $G$-equivalent. Then, $v\in\as^F$ and $u\in\as^K$ are called {\em $G$-equivalent} ($v\sim_G u$) whenever there exists a $g\in G$ such that $u(g h)=v(h),\; h\in F$. The set of {\em patterns \gls{PaG}} is defined by the quotient
$$
\PaG \; 
	:= \; \{v\in\as^K\;|\; K\in\ks(G)\}/\!\!\sim_G \, .
$$
The elements of $\PaG$ are denoted by $[v]$ where the bracket indicates that it is an equivalence class. For $K\in\ks(G)$, the set $\as^{[K]}:=\as^{K}/\!\!\sim_G=\big\{[v]\;|\; v\in\as^K\big\}$ is called the {\em set of patterns with support $[K]$}. 
\end{definition}

It is worth noticing that patterns play the role of words for the symbolic dynamical system $\big(\as^\ZM,\ZM,\alpha)$. The following assertion justifies the notation $\as^{[K]}$. More precisely, if $F\in[K]\in\ke$ then the equation $\as^{[K]}=\as^{[F]}$ holds.

\begin{lemma}
\label{Chap2-Lem-EquivAK}
Let $F,K\in\ks(G)$ be $G$-equivalent, i.e., $gF=K$ for a $g\in G$. Then the sets $\as^{[F]}$ and $\as^{[K]}$ are equal. In detail, the map $\Phi_g:\as^F\to\as^K,\; u\in\as^F\mapsto u(g^{-1}\,\cdot)\in\as^K$, defines a bijection.
\end{lemma}

\begin{proof}
The map $\Phi_g$ is well-defined as $gF=K$. If $u, v\in\as^F$ are different, then there is a $h\in F$ such that $u(h)\neq v(h)$. Then $gh\in K$ and $\Phi_g(u)(gh)\neq\Phi_g(v)(gh)$ hold proving that $\Phi_g$ is injective. For $v\in\as^K$, define $\tilde{v}\in\as^F$ by $\tilde{v}(h):= v(gh),\; h\in F$. Clearly, $\Phi_g(\tilde{v})=v$ holds and so $\Phi_g$ is surjective.
\end{proof}

\medskip

Let $K\subseteq F$ be two non-empty, compact sets and $u\in\as^F$. Then the {\em restriction of the map $u:F\to\as$ to $K$} is denoted by $u|_{K}:K\to\as$.

\begin{definition}[Subpattern and extension]
\label{Chap2-Def-SubpExt}
Let $[u]\in\as^{[F]}$ be a pattern with representative $u\in\as^F$. The set $\text{\gls{Sub}}:=\{u|_K\in\as^K\;|\; K\subseteq F\}/\!\!\sim_G$ is called {\em the set of subpatterns of $[u]$}. A $[u]\in\PaG$ is called {\em extension of $[v]\in\PaG$ if $[v]\in\Sub[u]$}. The set of all extensions of $[v]\in\PaG$ is denoted by \gls{Ext}.
\end{definition}

Note that $\Sub[u]$ and $\Ext[u]$ are independent of the choice of the representative $u\in\as^F$ according to Lemma~\ref{Chap2-Lem-EquivAK}. The set $\Sub[u]$ is always finite as $F\in\ks(G)$ and $\as$ are finite. Furthermore, it follows from the definition that $\Sub[v]\subseteq\Sub[u]$ whenever $[v]\in\Sub[u]$. With this notion at hand, the formal definition of a dictionary is provided.

\begin{definition}[Dictionary]
\label{Chap2-Def-Dictionary}
A non-empty set $\text{\gls{ws}}\subseteq \PaG$ is called a {\em dictionary} if
\begin{description}
\item[(D1)\label{(D1)}] the set $\Sub[u]$ is contained in $\ws$ for all $[u]\in\ws$, \hfill\textbf{(heredity)}
\item[(D2)\label{(D2)}] for every $u\in\as^K,\; K\in\ks(G)$ with $[u]\in\ws$ and each $F\in\ks(G)$\\
satisfying $K\subseteq F$ there is a $v\in\as^F$ such that $v|_K=u$ and $[v]\in\ws$.\hfill \textbf{(extensibility)}
\end{description}
For a dictionary $\ws$, the set $\ws\cap\as^{[K]}$ is called the {\em local patterns of $\ws$ with support $[K]$}. The set of all dictionaries over the alphabet $\as$ and the group $G$ is denoted by \gls{DG}.
\end{definition}

The notation $\ws$ is chosen for a dictionary in accordance to the German translation \glqq W\"orterbuch\grqq.

\medskip

In case that $G=\ZM$ there are two directions one to the left $-1$ and one to the right $+1$. Then Condition~\nameref{(DZ2)}, which is the analog of \nameref{(D2)}, only imposes an extension of a word $u\in\ws$ to the left and to the right. Due to the loss of specific directions for a general countable group $G$, Condition~\nameref{(D2)} imposes the existence of an extension in any direction. Note that in the one-dimensional case it is assumed that the empty word is an element of a dictionary. The empty word arises naturally as the neutral element of the concatenation map. Such a map does not exists in the higher dimensional case due to the loss of two explicit directions. The empty word can also be defined in this case, but it is not useful anymore. 

\begin{lemma}
\label{Chap2-Lem-DictAllCompNon-Empty}
Let $\ws\in\DG$ be a dictionary. Then, for every $[K]\in\ke$, the intersection $\ws\cap\as^{[K]}$ is non-empty.
\end{lemma}

\begin{proof}
A dictionary $\ws\in\DG$ is non-empty by definition and so there exist an $F\in\ks(G)$ and a $u\in\as^F$ such that $[u]\in\ws$. Let $[K]\in\ke$ be with representative $K\in\ks(G)$. Then $F\cup K\in\ks(G)$ is derived and due to \nameref{(D2)}, there is a $v\in\as^{F\cup K}$ such that $[v]\in\ws$. By \nameref{(D1)}, the restriction $[v|_K]$ is a pattern in $\ws$. Thus, $[v|_K]\in\ws\cap\as^{[K]}$ follows.
\end{proof}

\medskip

Condition~\nameref{(D2)} in Definition~\ref{Chap2-Def-Dictionary} can be replaced if $G$ admits an exhausting sequence of compact sets. The notion of exhausting sequence and its implications are discussed in the following.

\begin{definition}[Exhausting sequence]
\label{Chap2-Def-ExhaustSeq}
A sequence $(K_n)_{n\in\NM}$ {\em exhausts} the group $G$ if
\begin{description}
\item[(E1)\label{(E1)}] for all $n\in\NM$, the set $K_n\subseteq G$ is non-empty and compact,
\item[(E2)\label{(E2)}] the inclusion $K_n\subseteq K_{n+1}$ holds for all $n\in\NM$,
\item[(E3)\label{(E3)}] the union $\bigcup_{n\in\NM}K_n$ is equal to $G$.
\end{description}
\end{definition}

In our setting of a discrete, countable group, an exhausting sequence exists always.

\begin{lemma}
\label{Chap2-Lem-ExExplSeq}
Let $G$ be a discrete group with countably many elements. Then there exists a sequence $(K_n)_{n\in\NM}$ that exhausts $G$.
\end{lemma}

\begin{proof}
Let $\jmath:\NM\to G$ be surjective. Define $K_1:=\{\jmath(1)\}$ and $K_{n+1}:=K_n\cup\{\jmath(n+1)\}$ for $n\in\NM$. Then $(K_n)_{n\in\NM}$ is a sequence of non-empty, compact sets that exhausts $G$.
\end{proof}

\begin{proposition}
\label{Chap2-Prop-CharDict}
Let $(K_n)_{n\in\NM}$ be an exhausting sequence of $G$ and $\ws\subseteq\PaG$ be a non-empty subset satisfying \nameref{(D1)} of Definition~\ref{Chap2-Def-Dictionary}.  Then the following assertions are equivalent.
\begin{itemize}
\item[(i)] The set $\ws$ satisfies \nameref{(D2)}, i.e., $\ws$ is a dictionary.
\item[(ii)] Let $u\in\as^K,\; K\in\ks(G)$ and $F\in\ks(G)$ be such that $[u]\in\ws$ and $K\subseteq F$. For a $g\in G$, there is a $v\in\as^{gF}$ with $v(h)=u(g^{-1}h)$ for all $h\in gK$ and $[v]\in\ws$.
\item[(iii)] For every $u\in\as^K,\; K\in\ks(G)$ with $[u]\in\ws$, there exists an $n_0\in\NM$ such that for each $n\geq n_0$ there is a $v\in\as^{K_n}$ with $K\subseteq K_n$, $v|_K=u$ and $[v]\in\ws$.
\item[(iv)] For every $u\in\as^K,\; K\in\ks(G)$ with $[u]\in\ws$, there exists an $n_0\in\NM$ such that  $K\subsetneq K_{n_0}$ and there is a $v\in\as^{K_{n_0}}$ with $v|_K=u$ and $[v]\in\ws$.
\end{itemize}
\end{proposition}

\begin{proof}
Let $u\in\as^K,\; K\in\ks(G)$ be such that $[u]\in\ws$.

\vspace{.1cm}
$(i)\Leftrightarrow(ii):$ The equivalence follows as each $v\in\as^F,\; F\in\ks(G)$ can be shifted by a fixed $g\in G$ without changing the associated equivalence class $[v]$. More precisely, define $\tilde{v}:=v(g^{-1}\,\cdot)\in\as^{gF}$ for a $g\in G$ then $[\tilde{v}]=[v]$.

\vspace{.1cm}

$(i)\Rightarrow(iii):$ According to \nameref{(E3)} there is an $n_0\in\NM$ such that $K\subseteq K_n$ for all $n\geq n_0$. Then (i) implies (iii) by setting $F:=K_n$ for $n\geq n_0$.

\vspace{.1cm}

$(iii)\Rightarrow(iv):$ This is clear.

\vspace{.1cm}

$(iv)\Rightarrow(i):$ Let $u\in\as^K$ for $K\in\ks(G)$ and $F\in\ks(G)$ be such that $K\subseteq F$. By \nameref{(E3)}, there exists an $n_0\in\NM$ such that $F\subseteq K_{n_0}$. Condition (iv) guarantees the existence of an $n_1\in\NM$ such that $K\subsetneq K_{n_1}$ and there is a $v_1\in\as^{K_{n_1}}$ with $v_1|_K=u$ and $[v_1]\in\ws$. If $n_1<n_0$, then (iv) is applied again to $v_1$. In detail, an $n_2\in\NM$ exists such that $K_{n_1}\subsetneq K_{n_2}$ and there is a $v_2\in\as^{K_{n_2}}$ with $v_2|_{K_{n_1}}=v_1$ and $[v_2]\in\ws$. As $K_{n_1}\subsetneq K_{n_2}$, \nameref{(E2)} implies $n_2>n_1$. By construction, the equation $v_2|_K=u$ holds. Repeating at most finitely many times an $N\geq n_0$ is constructed such that there exists a $v_N\in\as^{K_N}$ with $v_N|_K=u$, $[v_N]\in\ws$ and $K\subseteq F\subseteq K_N$. Set $\tilde{v}:=v_N|_F$. Since $[v_N]\in\ws$ and $\ws$ satisfies \nameref{(D1)}, the pattern $[\tilde{v}]$ is an element of $\ws$. By construction, $\tilde{v}$ is an element of $\as^F$ and the equations $\tilde{v}|_K=v_N|_K=u$ hold. As $F\in\ks(G)$ was arbitrary, Condition~\nameref{(D2)} follows.
\end{proof}

\begin{remark}
\label{Chap2-Rem-RelAZExplSeq}
The sets $K_n:=\{-n,-n+1,\ldots,n\},\; n\in\NM,$ define an exhausting sequence for the group $G:=\ZM$. Then the assertion Proposition~\ref{Chap2-Prop-CharDict} (iv) turns out to be the general analog of Condition~\nameref{(DZ2)} for dictionaries on $\as$ over the group $\ZM$.
\end{remark}

In the following, a topology is defined on $\DG$. For $[K]\in\ke$ and $U\subseteq\as^{[K]}$, define the set
\begin{equation}
\label{Chap2-Eq-TopDict}
\vs([K],U) \;
	:= \; \{\ws\in\DG\;|\; \ws\cap\as^{[K]}=U\}\subseteq\DG \, .
\end{equation}
The collection of these sets for all $[K]\in\ke$ and $U\in\as^{[K]}$ is denoted by $\bs(\DG)$. Before proving that $\bs$ defines a basis for a topology, the following key lemma is proven. It asserts that whenever $\vs([F],U_F), \vs([K],U_K)\in\bs(\DG)$ have non-empty intersection and $K\subseteq F$ then the inclusion $\vs([F],U_F)\subseteq \vs([K],U_K)$ follows.

\begin{lemma}
\label{Chap2-Lem-InclBasNonEmptInters}
Let $K,F\in\ks(G)$ be such that $K\subseteq F$. Then the inclusion $\vs([F],U_F)\subseteq \vs([K],U_K)$ holds for $U_K\subseteq\as^{[K]}$ and $U_F\subseteq\as^{[F]}$ if the intersection $\vs([F],U_F)\cap \vs([K],U_K)$ is non-empty. In particular, if $U_1,U_2\subseteq\as^{[F]}$ satisfy $\vs([F],U_1)\cap\vs([F],U_2)\neq\emptyset$ then $U_1=U_2$ follows.
\end{lemma}

\begin{proof}
Let $\ws'$ be an element of the intersection $I:=\vs([F],U_F)\cap \vs([K],U_K)$. Thus, the dictionary $\ws'$ satisfies $\ws'\cap\as^{[K]}=U_K$ and $\ws'\cap\as^{[F]}=U_F$. It suffices to prove for each $\ws\in\vs([F],U_F)$ that $\ws\cap\as^{[K]}=U_K$ is valid. The proof makes extensive use of the fact that $\ws\cap\as^{[F]}=\ws'\cap\as^{[F]}=U_F$.

\vspace{.1cm}

$\subseteq:$ Consider a $[u]\in\ws\cap\as^{[K]}$ with representative $u\in\as^K$. According to \nameref{(D2)} there exists a $v\in\as^F$ such that $v|_K=u$ and $[v]\in\ws\cap\as^{[F]}=U_F$. Thus, $[v]$ is an element of $\ws'$ and by \nameref{(D1)} the set $\Sub[v]$ is included in $\ws'$. By construction, $[u]$ is contained in $\Sub[v]$ and so $[u]\in\ws'\cap\as^{[K]}=U_K$.

\vspace{.1cm}

$\supseteq:$ Consider a $[u]\in U_K=\ws'\cap\as^{[K]}$ with representative $u\in\as^K$. Since $K\subseteq F$ holds and \nameref{(D2)} applies to $\ws'$, there is a $v\in\as^{F}$ with $v|_K=u$ and $[v]\in\ws'$. By the equation $\ws'\cap\as^{[F]} = \ws\cap\as^{[F]}$, the pattern $[v]$ is contained in $\ws$. According to the construction the pattern $[u]\in\as^{[K]}$ is a subpattern of $[v]$. Hence, the pattern $[u]$ is contained in $\ws\cap\as^{[K]}$ by \nameref{(D1)}. Thus, $U_K\subseteq\ws\cap\as^{[K]}$ follows.

\vspace{.1cm} 

Now let $F\in\ks(G)$ and $U_1,U_2\subseteq\as^{[F]}$ be so that $\vs([F],U_1)\cap\vs([F],U_2)\neq\emptyset$. Applying the first part of the lemma for $K=F$, this leads to the equation $\vs([F],U_1)=\vs([F],U_2)$. Consequently, $U_1=U_2$ is concluded since the equations $U_1=\ws\cap\as^{[F]}=U_2$ follow for $\ws\in\vs([F],U_1)=\vs([F],U_2)$.
\end{proof}

\medskip

Note that it might happen that a set $\vs([K],U)$ is empty for a $K\in\ks(G)$ and a $\emptyset\neq U\subseteq\as^{[K]}$. For instance, consider the group $G=\ZM$ with alphabet $\as:=\{a,b\}$. For $K=\{1,2\}$, the set $U:=\{ab,aa\}$ is contained in $\as^{[K]}$. For a dictionary $\ws$ containing the pattern $ab$, there exists a $u\in\as$ such that $abu\in\ws$ by \nameref{(D2)}. Thus, $bu$ is an element of $\ws$ by \nameref{(D1)}. Hence, for each dictionary $\ws$ with $ab\in\ws$, there is a $bu\in\ws\cap\as^{[K]}$ for a $u\in\as$. Consequently, $\vs([K],U)$ is empty for $K=\{1,2\}$ and $U:=\{ab,aa\}$.

\medskip

With Lemma~\ref{Chap2-Lem-InclBasNonEmptInters} at hand, it follows that the collection $\bs(\DG)$ (defined in Equation~\eqref{Chap2-Eq-TopDict}) forms a base for a topology on $\DG$. This is proven in the following proposition.

\begin{proposition}
\label{Chap2-Prop-DictBasTop}
The collection $\bs(\DG)$ of sets
$$
\vs([K],U) \; 
	:= \; \big\{\ws\in\DG\;\big|\; \ws\cap\as^{[K]}=U\big\}, \qquad [K]\in\ke,\; U\subseteq\as^{[K]} \, ,
$$
defines a base for a topology on $\DG$.
\end{proposition}

\begin{proof}
Recall the notion of a base for a topology, see e.g. Definition~\ref{App1-Def-BaseTopology}. Clearly the union $\bigcup_{U\subseteq\as^{[K]}}\vs([K],U)$ covers $\DG$ for every fixed $[K]$ in $\ke$. Thus, it suffices to show that, for each non-trivial intersection of elements of $\bs(\DG)$, there exists a non-empty element of $\bs(\DG)$ such that it is contained in the intersection. In detail, let $F,K\in\ks(G)$, $U_F\subseteq\as^{[F]}$ and $U_K\subseteq\as^{[K]}$ be non-empty such that $I:=\vs([F],U_F)\cap\vs([K],U_K)\neq\emptyset$. Consider a dictionary $\ws'\in I$. Define $\tilde{K}:=F\cup K$ and $U_{\tilde{K}}:=\ws'\cap\as^{[\tilde{K}]}$. According to Lemma~\ref{Chap2-Lem-DictAllCompNon-Empty}, the intersection $\ws'\cap\as^{[K']}$ is non-empty for all $K'\in\ks(G)$. Thus, $U_{\tilde{K}}$ is non-empty. Then the set $\vs([\tilde{K}],U_{\tilde{K}})\in\bs(\DG)$ contains at least $\ws'$. This also implies that $\vs([\tilde{K}],U_{\tilde{K}})\cap\vs([F],U_F) \neq\emptyset$ and $\vs([\tilde{K}],U_{\tilde{K}})\cap\vs([K],U_K) \neq\emptyset$. Since $F,K\subseteq\tilde{K}$, Lemma~\ref{Chap2-Lem-InclBasNonEmptInters} yields $\vs([\tilde{K}],U_{\tilde{K}})\subseteq I$.
\end{proof}

\begin{definition}[Local pattern topology]
\label{Chap2-Def-TopSpaceDict}
The topology induced by the base $\bs(\DG)$ defined in Equation~\eqref{Chap2-Eq-TopDict} is called the {\em local pattern topology on $\DG$}.
\end{definition}

A subset of a topological space is called connected if it cannot be split into two open non-empty sets $U,V$ such that $U\cap V=\emptyset$. A space is called {\em totally disconnected} if all connected components are singletons. Clearly, all connected components are singletons if $X$ is a Hausdorff space with a base of clopen sets for the topology where a set is called {\em clopen} if it is closed and open. Thus, a Hausdorff space with a clopen basis is totally disconnected. With this notion at hand, the topological properties of $\DG$ equipped with the local pattern topology are studied.

\begin{proposition}
\label{Chap2-Prop-SpacDict}
The topological space $\DG$ equipped with the local pattern topology is second-countable, compact, Hausdorff and totally disconnected. Furthermore, for every sequence $(K_n)_{n\in\NM}$ that exhausts $G$, the collection $\bs'\subseteq\bs(\DG)$ of sets
$$
\vs([K_n],U) \; 
	:= \; \big\{\ws\in\DG\;\big|\; \ws\cap\as^{[K_n]}=U\big\} \, , 
	\qquad U\subseteq\as^{[K_n]} \, , \, n\in\NM \, ,
$$
is also a base for the local pattern topology on $\DG$. In particular, the topology induced by $\bs'$ on $\DG$ is independent of the choice of the exhausting sequence $(K_n)_{n\in\NM}$.
\end{proposition}

\begin{proof} Let $(K_n)_{n\in\NM}$ be an exhausting sequence of $G$. The proof is organized as follows: 
\begin{itemize}
\item[(i)] The collection $\bs'$ defines a basis for the local pattern topology on $\DG$. 
\item[(ii)] The space $\DG$ is second-countable. 
\item[(iii)] The space $\DG$ is Hausdorff. 
\item[(iv)] The space $\DG$ is totally disconnected.
\item[(v)] The space $\DG$ is compact.
\end{itemize}

(i): That $\bs'$ defines a topology is proven analogously to Proposition~\ref{Chap2-Prop-DictBasTop}: It is clear that $\bs'$ covers $\DG$. For $n\leq m$, let $U_n\subseteq\as^{[K_n]}$ and $U_m\subseteq\as^{[K_m]}$ be such that $I:=\vs([K_n],U_n)\cap\vs([K_m],K_m)\neq\emptyset$. According to \nameref{(E2)} the inclusion $K_n\subseteq K_m$ holds. Thus, Lemma~\ref{Chap2-Lem-InclBasNonEmptInters} implies $\vs([K_m],K_m)\subseteq I$ where $\vs([K_m],K_m)\neq\emptyset$ holds by assumption.

\vspace{.1cm}

Since $\bs'\subseteq\bs(\DG)$ it suffices to prove that the topology induced by $\bs'$ is finer than the local pattern topology induced by $\bs(\DG)$. Let $\ws\in\DG$, $K\in\ks(G)$ and $U\subseteq\as^{[K]}$ be such that $\ws\in\vs([K],U)\in\bs(\DG)$. According to Lemma~\ref{App1-Lem-TopFiner}, it suffices to show the existence of a $\vs\in\bs'$ satisfying $\vs\subseteq\vs([K],U)$. By \nameref{(E3)}, there exists an $n_0\in\NM$ such that $K\subseteq K_{n_0}$. Define the subset $U_{n_0}:=\ws\cap\as^{[K_{n_0}]}\subseteq\as^{[K_{n_0}]}$ which is non-empty according to Lemma~\ref{Chap2-Lem-DictAllCompNon-Empty}. By construction, $\ws$ is an element of $\vs([K_{n_0}],U_{n_0})\in\bs'$ and Lemma~\ref{Chap2-Lem-InclBasNonEmptInters} guarantees that $\vs([K_{n_0}],U_{n_0})\subseteq\vs([K],U)$.

\vspace{.1cm}

(ii): By (i) a base for the local pattern topology is given by
$$
\bs' \;
	:= \; \bigcup\nolimits_{n\in\NM}\bigcup\nolimits_{U\subseteq\as^{[K_n]}} \{\vs([K_n],U)\} \, .
$$
Since $\as$ is finite, there are, for each fixed $n\in\NM$, at most finitely many different subsets $U\subseteq\as^{[K_n]}$. As the countable union of finite sets is countable, the set $\bs'$ is countable. Hence, the topology on $\DG$ is second-countable.

\vspace{.1cm}

(iii): Let $\ws_1,\ws_2\in\DG$ be different, i.e., there exists a $[K]\in\ke$ such that $\ws_1\cap\as^{[K]}\neq\ws_2\cap\as^{[K]}$. Thus, the sets $\vs([K],\ws_1\cap\as^{[K]})$ and $\vs([K],\ws_2\cap\as^{[K]})$ are disjoint and open. By construction, the relations $\ws_1\in\vs([K_{n_0}],U)$ and $\ws_2\in\vs([K_{n_0}],V)$ hold. Hence, $\DG$ is a Hausdorff space.

\vspace{.1cm}

(iv): It suffices to prove that the set $\vs([K],U):=\{\ws\in\DG\;|\; \ws\cap\as^{[K]}=U\}$ is clopen for $[K]\in\ke$ and $U\subseteq\as^{[K]}$. Let $\ws$ be an element of the closure $\overline{\vs([K],U)}$, i.e., each open neighborhood of $\ws$ intersects $\vs([K],U)$. As the open set $\vs([K],\ws\cap\as^{[K]})$ defines a neighborhood of $\ws$, it has to intersect $\vs([K],U)$. By Lemma~\ref{Chap2-Lem-InclBasNonEmptInters}, it follows that $\ws\cap\as^{[K]}=U$ and so $\ws\in\vs([K],U)$. Hence, $\overline{\vs([K],U)}\subseteq\vs([K],U)$ is derived implying that $\vs([K],U)$ is clopen.

\vspace{.1cm}

(v): Let $\ws_k\in\DG,\; k\in\NM,$ be a sequence. A convergent subsequence is constructed as follows: Since $\as$ is finite, the sets $\as^{K_n}$ and so $\as^{[K_n]}$ are finite for each $n\in\NM$. By induction over $n\in\NM$ the following assertion holds. There exists a function $\sigma:\NM\times\NM\to\NM$, which is defined inductively, such that for each $n\in\NM$ there exist $\emptyset\neq U_i\subseteq\as^{[K_i]},\; 1\leq i\leq n,$ and a subsequence $(\ws_{k_{\sigma(l,n)}})_{l\in\NM}$ of $(\ws_k)_{k\in\NM}$ satisfying the equations 
$$
U_i \;
	= \; \ws_{k_{\sigma(l,n)}}\cap\as^{[K_i]} \, ,
	\qquad 1\leq i\leq n \, ,
$$
hold for all $l\in\NM$. The inductive step $n\Rightarrow n+1$ is proven whereas the case $n=1$ is similarly treated.

\vspace{.1cm}

Let $(\ws_{k_{\sigma(l,n)}})_{l\in\NM}$ be a subsequence of $(\ws_k)_{k\in\NM}$ such that the above properties are fulfilled for $n\in\NM$. Since $\as^{[K_{n+1}]}$ is finite, there exists an $[u]\in\as^{[K_{n+1}]}$ such that the set 
$$
\big\{\sigma(l,n)\in\NM
	\;\big|\; [u]\in\ws_{k_{\sigma(l,n)}}\cap\as^{[K_{n+1}]}
		,\; l\in\NM
\big\}
$$ 
is infinite. Hence, a subsequence $(\ws_{k_{\sigma(l,n)_j}})_{j\in\NM}$ of $(\ws_{k_{\sigma(l,n)}})_{l\in\NM}$ is extracted such that $[u]\in\ws_{k_{\sigma(l,n)_j}}$ for all $j\in\NM$. If the equation $\ws_{k_{\sigma(l,n)_j}}\cap\as^{[K_{n+1}]}=\{[u]\}$ holds for all $j\in\NM$, set $U_{n+1}:=\{[u]\}$. Otherwise the argument is repeated (applied to this subsequence). Specifically, let $\ws_{k_{\sigma(l,n)_{j_m}}},\; m\in\NM,$ be a subsequence of $(\ws_{k_{\sigma(l,n)_j}})_{j\in\NM}$ such that there is a $[v]\neq [u]$ contained in $\ws_{k_{\sigma(l,n)_{j_m}}}$ for all $m\in\NM$. Thus, the set $\{[u],[v]\}$ is contained in the set $\ws_{k_{\sigma(l,n)_{j_m}}}$ for all $m\in\NM$. If $\{[u],[v]\}=\ws_{k_{\sigma(l,n)_{j_m}}}\cap\as^{[K_{n+1}]},\; m\in\NM,$ define $U_{n+1}:=\{[u],[v]\}$ otherwise repeat the argument. This argument has to be iterated at most finitely many times since $\as^{[K_{n+1}]}$ is finite. Thus, a set $U_{n+1}\subseteq\as^{[K_{n+1}]}$ is constructed together with a subsequence $(\ws_{k_{\sigma(l,n)_j}})_{j\in\NM}$ of $(\ws_{k_{\sigma(l,n)}})_{l\in\NM}$ such that $\ws_{k_{\sigma(l,n)_j}}\cap\as^{[K_{n+1}]}=U_{n+1}$ for all $j\in\NM$. Define $\sigma(j,n+1):=\sigma(l,n)_j$ for $j\in\NM$. Since $(\ws_{k_{\sigma(j,n+1)}})_{j\in\NM}$ is a subsequence of $(\ws_{k_{\sigma(l,n)}})_{l\in\NM}$, the equations
$$
U_i \; 
	= \; \ws_{k_{\sigma(j,n+1)}}\cap\as^{[K_i]} \, ,
	\qquad 1\leq i\leq n+1 \, .
$$
hold for all $j\in\NM$. This concludes the induction.

\vspace{.1cm}

Define the set 
$$
\ws \; := \;
	\left\{ [u]\in\PaG
		\;|\; [u]\in\Sub[u']  
			\text{ where } [u']\in \bigcup\nolimits_{n\in\NM} U_n 
	\right\}
		\subseteq\PaG \, .
$$
By definition, the set $\ws$ satisfies \nameref{(D1)}. Thus, the inclusion $\Sub[u]\cap\as^{[K_n]}\subseteq U_n$ holds for $[u]\in U_m$ and $n\leq m$.

\vspace{.1cm}

Now, it is proven that $\ws$ fulfills also \nameref{(D2)}. Let $u\in\as^K,\; K\in\ks(G)$ and $F\in\ks(G)$ with $K\subseteq F$ and $[u]\in\ws$, i.e., there exist an $n_0\in\NM$ and a $[u']\in U_{n_0}$ such that $[u]\in\Sub[u']$. Let $u'\in\as^{K_{n_0}}$ be the representative of $[u']$. According to Proposition~\ref{Chap2-Prop-CharDict}~(ii), there is no loss in generality in assuming that $K\subseteq K_{n_0}$ and $u'|_K=u$. By \nameref{(E3)}, there is an $n_1\geq n_0$ such that $F\subseteq K_{n_1}$. Due to construction of $\ws$, the pattern $[u']$ is contained in $\ws_{k_{\sigma(1,n_1)}}$. Applying \nameref{(D2)} to $\ws_{k_{\sigma(1,n_1)}}$, there is a $\tilde{v}\in\as^{K_{n_1}}$ such that $\tilde{v}|_{K_{n_0}}=u'$ and $[\tilde{v}]\in\ws_{k_{\sigma(1,n_1)}}\cap\as^{[K_{n_1}]}=U_{n_1}$. Define $v:=\tilde{v}|_F\in\as^F$. As the inclusions $K\subseteq K_{n_0}\subseteq K_{n_1}$ hold, the restriction $v|_K$ is equal to $\tilde{v}|_K=u'|_K=u$. Additionally, $[v]\in\ws$ follows since $[\tilde{v}]\in U_{n_1}\subsetneq\ws$ and $[v]\in\Sub[\tilde{v}]$ hold. This proves that $\ws$ satisfies \nameref{(D2)}.

\vspace{.1cm} 

By the previous considerations, $\ws$ is a dictionary and the family $\{\vs([K_n],U_n)\;|\; n\in\NM\}$ defines a neighborhood basis of $\ws\in\DG$. The proof is finished with a Cantor type argument that allows to extract a subsequence $(\ws_{k_n})_{n\in\NM}$ converging to $\ws$. More precisely, the subsequence is defined by $\ws_{k_n}:=\ws_{k_{\sigma(1,n)}}$ for $n\in\NM$.
\end{proof}

\medskip

The local pattern topology implies that dictionaries are close if and only if the local patterns supported on a large $[K]\in\ke$ coincide for the dictionaries. This fact is intensively used in Chapter~\ref{Chap5-OneDimCase} and Chapter~\ref{Chap6-HigherDimPerAppr} to prove the existence of strongly periodic approximation for subshifts.

\begin{corollary}
\label{Chap2-Cor-DictExhSeq}
Let $\as$ be a finite alphabet, $G$ be a countable, discrete group and $(K_n)_{n\in\NM}$ be an exhausting sequence of $G$. Then, for a dictionary $\ws\in\DG$, the sets 
$$
\vs_n(\ws) \;
	:= \; \vs\big([K_n],\ws\cap\as^{[K_n]}\big)\subseteq\DG\,,
	\qquad n\in\NM\,,
$$
define a neighborhood basis of $\ws$ in the local pattern topology. Furthermore, a sequence of dictionaries $(\ws_k)_{k\in\NM}$ converges to $\ws\in\DG$ in the local pattern topology if and only if, for every $n\in\NM$, there exists a $k_0\in\NM$ such that
$$
\ws_k\cap\as^{[K_n]}\;
	=\; \ws\cap\as^{[K_n]}\,,
	\qquad
	k\geq k_0\,.
$$
\end{corollary}

\begin{proof}
Let $\ws\in\DG$ be a dictionary. Due to Proposition~\ref{Chap2-Prop-SpacDict}, the sets
$$
\vs([K_n],U) \; 
	:= \; \big\{\ws\in\DG\;\big|\; \ws\cap\as^{[K_n]}=U\big\} \, , 
	\qquad U\subseteq\as^{[K_n]} \, , \, n\in\NM \, ,
$$
define a base for the local pattern topology. Thus, all the sets $\vs([K_n],U)$ satisfying $\ws\in\vs([K_n],U)$ define a neighborhood base of $\ws$, c.f. Definition~\ref{App1-Def-BaseTopology}. Clearly, $\ws\cap\as^{[K_n]}=U$ follows by definition if $\ws\in\vs([K_n],U)$. Hence, the sets $\vs_n(\ws)\,,\; n\in\NM\,,$ define a neighborhood basis of the dictionary $\ws$ in the local pattern topology.

\vspace{.1cm}

Recall the formal definition of convergent sequences in a topological space with a base, c.f. Definition~\ref{App1-Def-ConvergenceSequence}. Then a sequence of dictionaries $(\ws_k)_{k\in\NM}$ converges to $\ws\in\DG$ in the local pattern topology by definition if and only if, for each $n\in\NM$, there exists a $k_0\in\NM$ such that $\ws_k\in\vs\big([K_n],\ws\cap\as^{[K_n]}\big)$ for all $k\geq k_0$. Furthermore, a dictionary $\ws'\in\DG$ is an element of $\vs\big([K_n],\ws\cap\as^{[K_n]}\big)$ if and only if $\ws'\cap\as^{[K_n]} = \ws\cap\as^{[K_n]}$. Consequently, the convergence of the sequence $(\ws_k)_{k\in\NM}$ to $\ws\in\DG$ in the local pattern topology is characterized by the existence of a $k_0:=k_0(n)\in\NM$ for each $n\in\NM$ satisfying $\ws_k\cap\as^{[K_n]} = \ws\cap\as^{[K_n]}\,,\; k\geq k_0$.
\end{proof}

\begin{remark}
\label{Chap2-Rem-DictExhSeq}
It is worth noticing that Corollary~\ref{Chap2-Cor-DictExhSeq} asserts that a dictionary is uniquely defined by the union of patterns $\bigcup_{n\in\NM} \ws\cap\as^{[K_n]}\subseteq\ws$ since $\DG$ is Hausdorff and $\vs_n(\ws)\,,\; n\in\NM\,,$ defines a neighborhood base of $\ws$. In the one-dimensional case $G=\ZM$, an exhausting sequence $(K_n)_{n\in\NM}$ of $\ZM$ exists such that $K_n=\{1,\ldots,n\}$, c.f. Section~\ref{Chap5-Sect-SymbDynSystZM}. Thus, a dictionary $\ws$ is uniquely determined in this case by the set $\bigcup_{n\in\NM}\ws\cap\as^n$ since $\as^{[K_n]}=\as^n$. Note that the set $\bigcup_{n\in\NM} \ws\cap\as^n\subseteq\ws$ is usually called a dictionary in the literature.
\end{remark}

\subsection{The space of subshifts}
\label{Chap2-Ssect-SpaceSubshifts}

Let $\as$ be a finite set equipped with the discrete topology and $G$ be a discrete, countable group.
Consider the space 
$$
\text{\gls{asG}} \; 
	:= \; \prod\limits_{g\in G} \as \;
	= \; \{\xi:G\to \as\}
$$
endowed with the product topology. A basis for the product topology of $\as^G$, c.f. \cite[Definition~3.7]{Querenburg2001}, is given by the sets
$$
\text{\gls{osKu}} \; 
	:= \; \big\{ \xi\in\as^G\;\big|\; \xi|_K\in[u] \big\} \, ,
	\qquad K\in\ks(G)\,,
	\; 
	[u]\in\as^{[K]} \,,
$$
where $\xi|_K$ denotes the restriction of the map $\xi:G\to\as$ to the set $K\subseteq G$. The following lemma is similar to Lemma~\ref{Chap2-Lem-InclBasNonEmptInters} for the local pattern topology.

\begin{lemma}
\label{Chap2-Lem-SubshiftInclBasNonEmptInters}
Let $K,F\in\ks(G)$ be such that $K\subseteq F$. Then the inclusion $\os(F,[v])\subseteq \os(K,[u])$ holds for $[u]\in\as^{[K]}$ and $[v]\in\as^{[F]}$ if the intersection $\os(K,[u])\cap\os(F,[v])$ is non-empty.
\end{lemma}

\begin{proof}
Let $\xi\in\as^G$ be such that $\xi\in\os(K,[u])\cap\os(F,[v])$. Consider an $\eta\in\os(F,[v])$. Then $\eta|_F=\xi|_F$ follows by $\eta,\xi\in\os(F,[v])$ implying $\eta_K=\xi|_K$ since $K\subseteq F$. The restriction $\xi|_K$ is an element of the equivalence class $[u]$ by assumption. Hence, $\eta|_K\in [u]$ is concluded leading to $\eta\in\os(K,[u])$.
\end{proof}

\medskip

Define the map $\alpha:G\times\as^G\to\as^G$ by $\alpha_g(\xi)(h):=\xi(g^{-1}h),\; h\in G$. Recall that a topo\-logical space is called {\em totally disconnected} if the connected components are singletons. Thus, a Hausdorff space is totally disconnected if it has a clopen basis for the topology. 

\begin{proposition}[\cite{Tyc30,Cec37}]
\label{Chap2-Prop-FullShift}
Let $\as$ be a finite alphabet and $G$ a discrete, countable group. Then the space  $\as^G$ is second-countable, compact, Hausdorff and totally disconnected. Furthermore, the triple $(\as^G,G,\alpha)$ defines a topological dynamical system.
\end{proposition}

\begin{proof}
Since $G$ is countable and equipped with the discrete topology, the group is second-countable, Hausdorff and locally compact. The base for the product topology on $\as^G$ consists of clopen sets, i.e., closed and open sets. By the Hausdorff property, it follows that $\as^G$ is totally disconnected. Due to the finiteness of $\as$, the space $\as^G$ is second-countable, compact and Hausdorff \cite{Tyc30,Cec37}. By inspection the map $\alpha$ satisfies the conditions (i) and (ii) of Definition~\ref{Chap2-Def-DynSyst}. The sets of the form $\{g\}\times\os(K,[u])$ define a basis for the topology in $G\times\as^G$. Since the equation 
$$
\alpha^{-1}(\os(K,u)) \;
	= \; \bigcup\nolimits_{g\in G}\{g\}\times\os(g^{-1}K,[u])
$$
holds, the continuity of $\alpha$ follows.
\end{proof}

\medskip

The dynamical system $(\as^G,G,\alpha)$ is called the {\em full shift} or {\em Bernoulli shift}. In the literature, elements of $\SG\big(\as^G\big)$ are called {\em subshifts} instead of dynamical subsystems. Elements of $\SG\big(\as^G\big)$ are denoted by $\Xi$. This notation is chosen since $\Xi$ is the most disconnected Greek letter and all subshifts are totally disconnected, c.f. Proposition~\ref{Chap2-Prop-FullShift}. 

\begin{lemma}
\label{Chap2-Lem-xiDict}
For every $\xi\in\as^G$, the set
$$
\text{\gls{wsxi}} \; 
	:= \; \{\xi|_K\;|\; K\in\ks(G)\}/\!\!\sim_G\;\;\subseteq\PaG \, .
$$
defines a dictionary. Furthermore, the equation $\ws(\xi)=\ws(\alpha_g(\xi))$ holds for each $g\in G$.
\end{lemma}

\begin{proof}
By definition the set $\Sub[\xi|_K]=\{\xi|_F\;|\; F\subseteq K\}/\!\!\sim_G$ is contained in $\ws(\xi)$. Thus, $\ws(\xi)$ satisfies \nameref{(D1)}. Let $u\in\as^K,\; K\in\ks(G)$ be such that $[u]\in\ws(\xi)$, i.e., it exists a $g\in G$ so that $u\in[\xi|_{gK}]$. Consider an $F\in\ks(G)$ with $K\subseteq F$. Then $v(h)=u(g^{-1}h)\,,\; h\in gK\,,$ and $[v]\in\ws(\xi)$ hold for $v:=\xi|_{gF}$. Proposition~\ref{Chap2-Prop-CharDict} (ii) implies that $\ws(\xi)$ is a dictionary. 

\vspace{.1cm}

Note that each pattern in $\ws(\xi)$ is an equivalence class with respect to the shift by the group $G$. Thus, for a $g\in G$, the equality $\ws(\xi)=\ws(\alpha_g(\xi))$ follows since $\alpha_g(\xi)|_K=\xi|_{g^{-1}K}$ holds for every $K\in\ks(G)$.
\end{proof}

\medskip

The dictionary $\ws(\xi)$ is called {\em dictionary associated with $\xi\in\as^G$}.

\begin{lemma}
\label{Chap2-Lem-XiDictClosed}
Let $\Xi\in\SG\big(\as^G\big)$ be a subshift and $\xi\in\as^G$. If there exists an exhausting sequence $(K_n)_{n\in\NM}$ of $G$ such that, for each $n\in\NM$, there is an $\eta_n\in\Xi$ with $\eta_n|_{K_n}=\xi|_{K_n}$ then $\xi\in\Xi$.
\end{lemma}

\begin{proof}
Let $(K_n)_{n\in\NM}$ be an exhausting sequence of $G$ and $\eta_n\in\Xi,\; n\in\NM,$ be such that $\eta_n|_{K_n}=\xi|_{K_n}$. Then Lemma~\ref{Chap2-Lem-SubshiftInclBasNonEmptInters} implies that the family $\{\os(K_n,[\xi|_{K_n}])\;|\; n\in\NM\}$ is a neighborhood basis of $\xi$ since $(K_n)_{n\in\NM}$ is an exhausting sequence of $G$. Thus, the sequence $(\eta_n)_{n\in\NM}$ tends to $\xi$ in $\as^G$, c.f. Definition~\ref{App1-Def-ConvergenceSequence}. As $\eta_n\in\Xi,\; n\in\NM,$ and $\Xi$ is closed, $\xi\in\Xi$ is derived.
\end{proof}

\begin{lemma}
\label{Chap2-Lem-ExInfWord}
Let $\ws\in\DG$ be a dictionary. For every $[v]\in\ws$, there exists a $\xi\in\as^G$ such that $[v]\in\ws(\xi)$ and $\ws(\xi)\subseteq\ws$.
\end{lemma}

\begin{proof}
Let $[v]$ be an element of $\ws$ with representative $v\in\as^K$. Consider an exhausting sequence $(K_n)_{n\in\NM}$ of $G$. Define iteratively a sequence $v_n\in\ws\cap\as^{K_n},\; n\in\NM,$ such that there is an $n_0\in\NM$ with $[v]\in\Sub[v_n]$ for all $n\geq n_0$ and $v_n|_{K_{n-1}}=v_{n-1},\; n\in\NM\,,$ hold. The details are as follows: Proposition~\ref{Chap2-Prop-CharDict} implies the existence of an $n_0\in\NM$ and a $v_{n_0}\in\as^{K_{n_0}}$ such that $K\subseteq K_{n_0}$, $v_{n_0}|_K=v$ and $[v_{n_0}]\in\ws$. For $1\leq n\leq n_0$, define $v_n:=v_{n_0}|_{K_n}$. Then, for $n\in\NM$, the element $v_n\in\ws\cap\as^{K_n}$ is constructed by applying \nameref{(D2)} to $v_{n-1}\in\as^{K_{n-1}}$. In detail, as $K_{n-1}\subseteq K_n$ there exists a $v_n\in\as^{K_n}$ such that $v_n|_{K_{n-1}}=v_{n-1}$ and $[v_n]\in\ws$.

\medskip

The sets $\os(K_n,[v_n]),\; n\in\NM,$ are compact as closed subsets of the compact space $\as^G$, c.f. Proposition~\ref{Chap2-Prop-FullShift}. By construction of $(v_n)_{n\in\NM}$, any finite intersection $\bigcap_{j=1}^l\os(K_{n_j},[v_{n_j}])$ is equal to $\os(K_N,[v_N])\neq\emptyset$ where $N:=\max\{n_1,\ldots,n_l\}$. Thus, the intersection
$$
I \; 
	:= \; \bigcap\nolimits_{n\in\NM}\os(K_n,[v_n])
$$
is non-empty since any finite intersection of the compact sets is non-empty, c.f. \cite[Satz~8.2]{Querenburg2001}. If $\xi\in I$, then $\xi|_{K_n}=v_n$ for all $n\in\NM$. Thus, $I$ contains exactly one element. For indeed, if there were $\xi\neq\eta\in I$ then by \nameref{(E3)} there exists an $n_0\in\NM$ such that $\xi|_{K_{n_0}}\neq\eta|_{K_{n_0}}$ contradicting that $\xi,\eta\in I$. Denote this unique element in $I$ by $\xi$. 

\medskip

It is left to verify that $\ws(\xi)\subseteq\ws$: Let $[u]\in\ws(\xi)$ be with representative $\xi|_K=u$ for $K\in\ks(G)$. According to \nameref{(E3)}, there exists an $n_0\in\NM$ such that $K\subseteq K_{n_0}$. Thus, the $v_{n_0}:=\xi|_{K_{n_0}}$ fulfills $[v_{n_0}]\in\ws$ and $[u]\in\Sub[v_{n_0}]$. Then \nameref{(D1)} implies $[u]\in\ws$.
\end{proof}

\medskip

The following result is new in the literature. It was known before that, for a dictionary $\ws(\eta)$ associated with an $\eta\in\as^\ZM$, the set $\{\xi\in\as^\ZM\;|\; \ws(\xi)\subseteq\ws(\eta)\}$ defines a subshift. This idea is used to prove that the map $\Phi:\SG\big(\as^G\big)\to\DG$, defined in Theorem~\ref{Chap2-Theo-Shift+DictSpace}, is surjective in the more general case of a countable, discrete group $G$. Note that the bijectivity of $\Phi$ for $G=\ZM$ was known before, c.f. e.g. \cite[Proposition 1.3.4.]{LindMarcus95}.

\begin{theorem}
\label{Chap2-Theo-Shift+DictSpace}
Let $\as$ be an alphabet and $G$ a discrete, countable group. Then the map $\Phi:\SG\big(\as^G\big)\to\DG$ defined by
$$
\Phi(\Xi) \; 
	:= \; \text{\gls{wsXi}} \;
	:= \; \bigcup\nolimits_{\xi\in\Xi}\; \ws(\xi) \, ,
$$
is a homeomorphism. 
\end{theorem}

\begin{proof}
Since the union of dictionaries is a dictionary, Lemma~\ref{Chap2-Lem-xiDict} implies that $\ws(\Xi)$ is a dictionary for every $\Xi\in\SG\big(\as^G\big)$. By Proposition~\ref{Chap2-Prop-SpaDynSyst} and Proposition~\ref{Chap2-Prop-DictBasTop}, the space $\SG\big(\as^G\big)$ is compact and the space $\DG$ is Hausdorff. Thus, according to \cite[Satz~8.11]{Querenburg2001} it suffices to prove the continuity and the bijectivity of $\Phi$ to conclude that $\Phi$ defines a homeomorphism. In the following it is proven that $\Phi$ is (i) injective, (ii) surjective and (iii) continuous.

\vspace{.1cm}

(i): Let $\Xi_1,\Xi_2\in\SG\big(\as^G\big)$ be different. Without loss of generality, suppose that $\Xi_1\setminus\Xi_2\neq\emptyset$ holds. For $\xi\in\Xi_1\setminus\Xi_2$, there exists a $K\in\ks(G)$ such that $[\xi|_K]\not\in \Phi(\Xi_2)$ and $[\xi|_K]\in \Phi(\Xi_1)$. Otherwise, Lemma~\ref{Chap2-Lem-XiDictClosed} would imply that $\xi\in\Xi_2$. Thus, $\Phi(\Xi_1)\neq\Phi(\Xi_2)$ follows. Consequently, the map $\Phi$ is injective.

\vspace{.1cm}

(ii): Let $\ws$ be a dictionary. Define the set $\Xi(\ws)\subseteq\as^G$ by 
$$
\text{\gls{Xiws}} \; 
	:= \; \big\{\xi\in\as^G\;\big|\; \ws(\xi)\subseteq\ws\big\} \, .
$$
According to Lemma~\ref{Chap2-Lem-ExInfWord}, $\Xi(\ws)$ is a non-empty set. The set $\Xi(\ws)$ is $G$-invariant since, for each $g\in G$, the equation $\ws(\xi)=\ws(\alpha_g(\xi))$ holds. The family $\{\os(K,[\xi|_K])\;|\; K\in\ks(G)\}$ is a neighborhood base of $\xi$, c.f. Definition~\ref{App1-Def-BaseTopology}. Let $(\xi_n)_{n\in\NM}\subseteq\Xi(\ws)$ be a convergent sequence to $\xi\in\as^G$. Thus, for every $K\in\ks(G)$, there is an $n_0\in\NM$ such that $\xi_n\in\os(K,[\xi|_K])$ for all $n\geq n_0$. Hence, the equation $\xi_{n_0}|_K=\xi_K$ holds implying that $[\xi|_K]$ is an element of $\ws(\xi_{n_0})\subseteq\ws$. Consequently, the inclusion $\ws(\xi)\subseteq\ws$ follows and so $\xi\in\Xi(\ws)$ is derived. Thus, $\Xi(\ws)$ is closed. Altogether, $\Phi(\Xi(\ws))=\ws$ is concluded by Lemma~\ref{Chap2-Lem-ExInfWord} where $\Xi(\ws)$ is a subshift by the previous considerations.

\vspace{.1cm}

(iii): Consider a non-empty element $\vs([K],U)$ of the base of the local pattern topology on $\DG$ for $[K]\in\ke$ and $\emptyset\neq U\subseteq\as^{[K]}$. In order to prove the continuity of $\Phi$ it suffices to show that $\Phi^{-1}(\vs([K],U))$ is open.

\vspace{.1cm}

Since $\as$ is finite, the sets $\as^{[K]}$ and so $U=\{[u_1],\ldots,[u_l]\}$ are finite. Let $K\in\ks(G)$ be a representative of $[K]$. Define 
$$
F \; := \;
	\underset{j=1}{\overset{l}{\bigcap}} \Big(\as^G\setminus\os\big(K,[u_j]\big)\Big) \; 
		= \; \big\{\xi\in\as^G\;\big|\; \text{all }\; 1\leq j\leq l:\, \xi|_K\not\in[u_j] \big\} \, .
$$ 
Define $\os_j:=\os(K,[u_j])=\{\xi\in\as^G\;|\; \xi|_K\in[u_j] \}$ for $1\leq j\leq l$. Since the sets $\os(K,[u_j]),\; 1\leq j\leq l,$ are clopen, the set $F\subseteq\as^G$ is closed and so compact. Consider, for the finite family $\os:=\{\os_1,\ldots,\os_l\}$ of open sets the open set
$$
\us(F,\os) \;
	:= \; \big\{ \Xi\in\SG\big(\as^G\big)\;\big|\; \text{all } 1\leq j\leq l:\; \Xi\cap\os(K,[u_j])\neq\emptyset \text{ and } F\cap\Xi=\emptyset \big\}
$$
in the Hausdorff-topology of $\SG\big(\as^G\big)$. Then, the equality $\us(F,\os)= \Phi^{-1}(\vs([K],U))$ holds:

\vspace{.1cm}

$\subseteq$: It suffices to show $\Phi(\Xi)\cap\as^{[K]}=U$ for each $\Xi\in\us(F,\os)$. Let $1\leq j\leq l$. Since $\Xi\in\us(F,\os)$, the intersection $\Xi\cap\os(K,[u_j])$ is non-empty, i.e., there exists a $\xi\in\Xi$ with $\xi|_K\in[u_j]$. Thus, $[u_j]\in\ws(\xi)\cap\as^{[K]}\subseteq\Phi(\Xi)\cap\as^{[K]}$ holds leading to the inclusion $U\subseteq\ws(\Xi)\cap\as^{[K]}$. The converse inclusion is deduced by the fact that $F\cap\Xi=\emptyset$: If there were a $[v]\in\Phi(\Xi)\cap\as^{[K]}$ with $[v]\not\in U$, then there exist a representative $\tilde{K}$ of $[K]$ and a $\xi\in\Xi$ such that $\xi|_{\tilde{K}}\in[v]$. Since $\Xi$ is $G$-invariant, there is a $g\in G$ such that $K=g\tilde{K}$ and $\alpha_g(\xi)|_K=\xi|_{g\tilde{K}}\in[v]$. This contradicts the fact that $F\cap\Xi=\emptyset$. Hence, $\Phi(\Xi)\cap\as^{[K]}\subseteq U$ holds.

\vspace{.1cm}

$\supseteq$: For every $\Xi\in \Phi^{-1}(\vs([K],U))$, it follows $\Phi(\Xi)\cap\as^{[K]}=U$. Consequently, the intersections $\Xi\cap\os(K,[u_j]),\; 1\leq j\leq l,$ are non-empty by using the $G$-invariance of $\Xi$. Furthermore, if $F\cap\Xi$ were non-empty, there exists a $\xi\in\Xi$ such that $\xi|_K\not\in [u_j]$ for all $1\leq j\leq l$. Hence, $[\xi|_K]\in\ws(\Xi)\cap\as^{[K]}$ and $[\xi|_K]\not\in U$ follow being a contradiction to $\Xi\in \Phi^{-1}(\vs([K],U))$. Thus, the intersection $F\cap\Xi$ is empty proving that $\Xi\in\us(F,\os)$.
\end{proof}

\medskip

According to Theorem~\ref{Chap2-Theo-Shift+DictSpace}, a sequence of subshifts converges if and only if the local patterns of the subshifts converge. Furthermore, it immediately leads to the following consequence. 

\begin{corollary}
\label{Chap2-Cor-SGTotaDisco}
The space $\SG\big(\as^G\big)$ is totally disconnected.
\end{corollary}

\begin{proof}
This follows by Proposition~\ref{Chap2-Prop-SpacDict} and Theorem~\ref{Chap2-Theo-Shift+DictSpace}. 
\end{proof}

\medskip

Subshifts arising by the orbit closure of a $\xi\in\as^G$ are of particular interest since their associated dictionary is exactly the dictionary $\ws(\xi)$.

\begin{corollary}
\label{Chap2-Cor-DictOrbitSubshift}
Let $\as$ be a finite alphabet and $G$ a discrete, countable group and $\Xi\in\SG\big(\as^G\big)$ be a subshift. Then the following assertions are equivalent.
\begin{itemize}
\item[(i)] The subshift $\Xi$ is topologically transitive, i.e., there is a $\xi\in\as^G$ with $\Xi=\overline{\Orb(\xi)}$.
\item[(ii)] There is a $\xi\in\as^G$ such that $\ws(\xi)=\ws(\Xi)$.
\end{itemize}
\end{corollary}

\begin{proof}
(i)$\Rightarrow$(ii): Let $\xi\in\as^G$ be such that $\Xi=\overline{\Orb(\xi)}$. According to Lemma~\ref{Chap2-Lem-xiDict}, $\ws(\xi)$ defines a dictionary. The dictionary $\ws(\Xi)$ is defined by the union $\bigcup_{\eta\in\Xi}\ws(\eta)$. Thus, $\xi\in\Xi$ yields $\ws(\xi)\subseteq\ws(\Xi)$. For the converse inclusion, consider the set 
$$
\Xi'\;
	:= \; \big\{
		\eta\in\as^G\;\big |\; \ws(\eta)\subseteq\ws(\xi) 
	\big\}
	\,.
$$
In (ii) of the proof of Theorem~\ref{Chap2-Theo-Shift+DictSpace}, it is proven that $\Xi'$ is a subshift such that $\ws(\Xi')\subseteq\ws(\xi)$. By construction $\xi$ is an element of $\Xi'$ implying that $\overline{\Orb(\xi)}\subseteq\Xi'$ as $\Xi'$ is closed and $G$-invariant. Hence, the inclusion $\Xi\subseteq\Xi'$ is derived leading to $\ws(\Xi)\subseteq\ws(\Xi')\subseteq\ws(\xi)$.

\vspace{.1cm}

(ii)$\Rightarrow$(i): Let $\xi\in\as^G$ be such that $\ws(\xi)=\ws(\Xi)$. Due to Theorem~\ref{Chap2-Theo-Shift+DictSpace}, the equation $\Xi=\{\eta\in\as^G\;|\; \ws(\eta)\subseteq\ws(\Xi)\}$ holds. Thus, $\Orb(\xi)\subseteq\Xi$ follows. Let $\eta\in\Xi$ and $K\in\ks(G)$ be arbitrary. Then $[\eta|_K]$ is an element of $\ws(\Xi)=\ws(\xi)$ implying that $\os(K,[\eta|_K])\cap\Orb(\xi)\neq\emptyset$. Consequently, $\Orb(\xi)$ intersects each open neighborhood of $\eta\in\Xi$ since the sets 
$$
\os(K,[\eta|_K])\;
	:= \; \big\{\tilde{\eta}\;\big|\; \tilde{\eta}|_K\in[\eta|_K]\big\}\,,
	\qquad K\in\ks(G)\,,
$$ 
define a neighborhood base for $\eta$. Hence, $\Orb(\xi)$ is dense in $\Xi$ leading to $\Xi=\overline{\Orb(\xi)}$.
\end{proof}

\medskip

The following result is well-known in the case of $G=\ZM$, c.f. \cite[Proposition~5.1.10]{Fogg02}. For a general countable, discrete group $G$, this assertion is new in the literature.

\begin{proposition}
\label{Chap2-Prop-MinDictConst}
Let $\Xi\in\SG\big(\as^G\big)$ be a subshift. Then the following are equivalent.
\begin{itemize}
\item[(i)] The subshift $\Xi$ is minimal.
\item[(ii)] The family of dictionaries associated with $\Xi$ is constant, i.e., the equation $\ws(\xi)=\ws(\eta)$ holds for all $\xi,\eta\in\Xi$.
\end{itemize}
\end{proposition}

\begin{proof}
(i)$\Rightarrow$(ii): Let $\xi,\eta\in\Xi$. It suffices to prove that $\ws(\xi)\subseteq\ws(\eta)$. Consider a $[u]\in\ws(\xi)$ and a $K\in\ks(G)$ such that $\xi|_K\in[u]$. Define the open neighborhood $\us(K,[u])$ of $\xi$. Due to minimality, the orbit $\Orb(\eta)$ intersects the open set $\us(K,[u])$. Thus, there exists a $g\in G$ such that $\alpha_g(\eta)\in\us(K,[u])$ and so $\alpha_g(\eta)|_K\in[u]$ is derived. This yields $[u]\in\ws(\eta)$ by using $\ws(\alpha_g(\eta))=\ws(\eta)$, c.f. Lemma~\ref{Chap2-Lem-xiDict}.

\vspace{.1cm}

(ii)$\Rightarrow$(i): Let $\xi\in\Xi$. In order to prove that $\Orb(\xi)$ is dense in $\Xi$ consider an $\eta\in\Xi$ and an open neighborhood $\us(K,[\eta|_K])$ for a $K\in\ks(G)$. By (ii) the pattern $[\eta|_K]$ is an element of $\ws(\xi)$. Hence, there exists a $\tilde{K}\in\ks(G)$ such that $\tilde{K}\sim_G K$ and $\xi|_{\tilde{K}}\in[\eta|_K]$ by definition of $\ws(\xi)$. Consequently, there is a $g\in G$ with $K=g\tilde{K}$ and $\alpha_g(\xi)|_K=\xi|_{\tilde{K}}\in [\eta|_K]$. Thus, the shift $\alpha_g(\xi)$ is contained in $\us(K,[\eta|_K])$ for each $K\in\ks(G)$ implying that $\Orb(\xi)$ intersects every open subset of $\Xi$.
\end{proof}

\section{Groupoids and associated \texorpdfstring{$C^\ast$-algebras}{C-algebras}}
\label{Chap2-Sect-GroupoidCalgebras}

Following the lines of \cite{Renault80}, basic properties of groupoids are introduced and proven. Then the associated reduced and full $C^\ast$-algebras of a groupoid are constructed. The property of a groupoid being amenable and the relation to the groupoid $C^\ast$-algebra is discussed on this basis. 

\medskip

A groupoid is a generalization of a group where not each element can be composed with every other element. More precisely, a groupoid is a small category where each morphism is an isomorphism. Following \cite{RenaultThesis78,Renault80} this notion is also formulated as follows.

\begin{definition}[Topological groupoid]
\label{Chap2-Def-Groupoid}
A {\em groupoid} \gls{Gamma} is a set \gls{Gamma1} equipped with a composition $\Gamma^{(2)}\subseteq\Gamma^{(1)}\times \Gamma^{(1)}\to\Gamma^{(1)},\; (\gamma,\rho)\mapsto \gamma\rho$ and an inverse $\Gamma^{(1)}\to\Gamma^{(1)},\;\gamma\mapsto\gamma^{-1}$ such that the assertions \nameref{(G1)},\nameref{(G2)} and \nameref{(G3)} hold.
\begin{description}
\item[(G1)\label{(G1)}] If $(\gamma,\zeta),(\zeta,\varrho)\in\Gamma^{(2)}$, then $(\gamma\zeta,\varrho),(\gamma,\zeta\varrho) \in \Gamma^{(2)}$ and the equation $(\gamma\zeta)\varrho = \gamma(\zeta\varrho)$ holds.
\item[(G2)\label{(G2)}] For all $\gamma\in\Gamma^{(1)}$, the pairs $(\gamma,\gamma^{-1})$ and $(\gamma^{-1},\gamma)$ are elements of $\Gamma^{(2)}$.
\item[(G3)\label{(G3)}] If $(\gamma,\varrho)\in\Gamma^{(2)}$, then the equalities $\gamma\varrho\varrho^{-1}=\gamma$ and $\gamma^{-1}\gamma\varrho=\varrho$ hold.
\end{description}
The set $\Gamma^{(1)}$ is called the {\em arrow set} whereas \gls{Gamma2} is called the {\em set of composable arrows}. Then the {\em source} and {\em range} $s,r:\Gamma^{(1)}\to\Gamma^{(1)}$ are defined by $s(\gamma):=\gamma^{-1}\gamma$ and $r(\gamma):=\gamma\gamma^{-1}$ for all $\gamma\in\Gamma^{(1)}$. The image $\text{\gls{Gamma0}}:=r(\Gamma^{(1)})$ is called the {\em unit space} of the groupoid $\Gamma$.

\vspace{.1cm}

A groupoid $\Gamma$ is called {\em topological} if $\Gamma^{(1)}$ is a second-countable, locally compact, Hausdorff space such that the inverse and the composition are continuous. Note that $\Gamma^{(2)}$ is endowed with the subspace topology of $\Gamma^{(1)}\times\Gamma^{(1)}$.
\end{definition}

The elements $\gamma\in\Gamma^{(1)}$ are called {\em arrows} due to the interpretation that $\gamma$ is an arrow joining the source $s(\gamma)$ and the range $r(\gamma)$, c.f. Figure~\ref{Chap2-Fig-TransfGroupGroupoid} on Page~\pageref{Chap2-Fig-TransfGroupGroupoid}. The space $\Gamma^{(0)}$ is called the set of units according to the observation that $\gamma s(\gamma)=\gamma=r(\gamma)\gamma$ by \nameref{(G3)}. Note that the inversion in a groupoid is globally defined while the composition is only partially defined.

\begin{remark}
\label{Chap2-Rem-SecCountGroupoid}
The condition of the groupoid to be second-countable can be dropped in specific cases. It is used for the disintegration of representations of groupoids, c.f. \cite{AnRe01}.
\end{remark}

In the following proposition some basic statements about topological groupoids are summarized.

\newpage
\begin{proposition}[\cite{RenaultThesis78,Renault80,GoehleThesis09}]
\label{Chap2-Prop-BasGroupoid}
Let $\Gamma$ be a topological groupoid. 
\begin{itemize}
\item[(a)] The set of composable pairs $\Gamma^{(2)}$ is equal to 
$$
\{ (\gamma,\rho)\in\Gamma^{(1)}\times\Gamma^{(1)}
	\;|\; s(\gamma)=r(\rho)
\} \, .
$$
\item[(b)] The equations $r(\gamma^{-1})=s(\gamma)$ and $s(\gamma^{-1})=r(\gamma)$ hold for all $\gamma\in\Gamma^{(1)}$. In particular, the unit space $\Gamma^{(0)}:=r(\Gamma^{(1)})$ is equal to $s(\Gamma^{(1)})$.
\item[(c)] For each $(\gamma,\rho)\in\Gamma^{(2)}$, the pair $(\rho^{-1},\gamma^{-1})$ is an element of $\Gamma^{(2)}$ and $(\gamma\rho)^{-1}=\rho^{-1}\gamma^{-1}$ holds.
\item[(d)] For every $(\gamma,\rho)\in\Gamma^{(2)}$, the equations $r(\gamma\rho)=r(\gamma)$ and $s(\gamma\rho)=s(\rho)$ hold.
\item[(e)] The inverse map $^{-1}:\Gamma^{(1)}\to\Gamma^{(1)}$ is a homeomorphism. 
\item[(f)] The range $r:\Gamma^{(1)}\to\Gamma^{(1)}$ and source $s:\Gamma^{(1)}\to\Gamma^{(1)}$ are continuous.
\item[(g)] The unit space $\Gamma^{(0)}$ is a closed subspace of $\Gamma^{(1)}$. 
\item[(h)] The set of composable pairs $\Gamma^{(2)}$ is closed in $\Gamma^{(1)}\times\Gamma^{(1)}$. 
\item[(i)] The spaces $\Gamma^{(0)}$ and $\Gamma^{(2)}$ are second-countable, locally compact and Hausdorff.
\end{itemize}
\end{proposition}

\begin{proof}
(a): Let $\gamma,\rho\in\Gamma^{(1)}$ be such that $s(\gamma)=r(\rho)$. Condition~\nameref{(G2)} insures that the pairs $(\gamma,\gamma^{-1}),(\gamma^{-1},\gamma),(\rho,\rho^{-1})$ and $(\rho^{-1},\rho)$ are elements of $\Gamma^{(2)}$. Then $(\gamma,\gamma^{-1}\gamma)\in\Gamma^{(2)}$ follows by \nameref{(G3)}. By assumption the equality $\gamma^{-1}\gamma=s(\gamma)=r(\rho)=\rho\rho^{-1}$ is derived. Combining this with the previous assertion, the pair $(\gamma,\rho\rho^{-1})$ is an element of $\Gamma^{(2)}$. As $(\rho\rho^{-1},\rho)\in\Gamma^{(2)}$, Condition~\nameref{(G1)} implies $(\gamma,\rho\rho^{-1}\rho)\in\Gamma^{(2)}$. Then \nameref{(G3)} concludes the argument as $(\gamma,\rho)=(\gamma,\rho\rho^{-1}\rho)\in\Gamma^{(2)}$.

\vspace{.1cm}

Suppose $(\gamma,\rho)$ is an element of $\Gamma^{(2)}$. Condition \nameref{(G2)} yields $(\gamma,\gamma^{-1}),(\rho,\rho^{-1})\in\Gamma^{(2)}$. Using \nameref{(G1)} and \nameref{(G3)}, the equations 
$$
s(\gamma)\; 
	=\; \gamma^{-1}\gamma
		=\; \left(\gamma^{-1}\gamma\right)\left(\rho\rho^{-1}\right)
			=\; \rho\rho^{-1} \;
				= \; r(\rho)
$$
follow.

\vspace{.1cm}

(b): Let $\gamma\in\Gamma^{(1)}$. Condition \nameref{(G2)} leads to $\big((\gamma^{-1})^{-1},\gamma^{-1}\big),(\gamma^{-1},\gamma)\in\Gamma^{(2)}$. Then \nameref{(G1)} and \nameref{(G3)} imply
$$
(\gamma^{-1})^{-1}\; 
	=\; (\gamma^{-1})^{-1} \left(\gamma^{-1}\gamma\right)\; 
		=\; \big((\gamma^{-1})^{-1} \gamma^{-1}\big)\gamma \;
			=\; \gamma \, .
$$ 
Thus, the equation $r(\gamma^{-1})=\gamma^{-1}(\gamma^{-1})^{-1}=s(\gamma)$ follows. Similarly, $s(\gamma^{-1})=r(\gamma)$ is proven. As each arrow $\gamma\in\Gamma^{(1)}$ admits an inverse $\gamma^{-1}$, the equation $\Gamma^{(0)}=s(\Gamma^{(1)})$ follows by the previous considerations.

\vspace{.1cm}

(c): Let $(\gamma,\rho)$ be an element of $\Gamma^{(2)}$. According to (a) and (b), the equation $r(\gamma^{-1})=s(\rho^{-1})$ holds. Thus, $(\rho^{-1},\gamma^{-1})$ is an element of $\Gamma^{(2)}$ by (a). Then the pairs $(\gamma^{-1},\gamma)$ and $(\gamma\rho,(\gamma\rho)^{-1})$ are elements of $\Gamma^{(2)}$ by \nameref{(G2)}. Hence, \nameref{(G1)} implies $(\rho^{-1}\gamma^{-1},\gamma\rho(\gamma\rho)^{-1})\in\Gamma^{(2)}$. Applying \nameref{(G1)} and \nameref{(G3)} once more, this leads to the desired equation
$$
\rho^{-1}\gamma^{-1}\; 
	=\; \rho^{-1}\gamma^{-1}\left(\gamma\rho(\gamma\rho)^{-1}\right)\;
		=\; \left(\rho^{-1}\gamma^{-1}\gamma\rho\right)(\gamma\rho)^{-1}\;
			=\; (\gamma\rho)^{-1} \, .
$$

(d): Let $(\gamma,\rho)$ be an element of $\Gamma^{(2)}$. According to (c) and \nameref{(G3)}, the equations
\begin{align*}
r(\gamma\rho) \; &= \; \gamma\rho (\gamma\rho)^{-1} \; = \; \gamma\rho\rho^{-1}\gamma^{-1} \; = \; r(\gamma) \, ,\\
s(\gamma\rho) \; &= \; (\gamma\rho)^{-1} \gamma\rho \; = \; \rho^{-1}\gamma^{-1}\gamma\rho \; = \; s(\rho)
\end{align*}
are deduced.

\vspace{.1cm}

(e): For every $\gamma\in\Gamma^{(1)}$, it is shown in (b) that $(\gamma^{-1})^{-1}=\gamma$. Thus, $^{-1}$ is an involution and so it is a homeomorphism.

\vspace{.1cm}

(f): This statement follows by the facts that the composition of continuous maps is continuous and that the composition and the inverse map on the groupoid $\Gamma$ are continuous by definition.

\vspace{.1cm}

(g): The space $\Gamma^{(1)}$ is second-countable and so sequences describe the topology. Consider a sequence $(u_n)_{n\in\NM}$ in $\Gamma^{(0)}$ tending to $u\in\Gamma^{(1)}$. According to (c) the equation $r(u_n)=u_n$ follows immediately. Then $\lim_{n\to\infty} u_n=r(u)$ is derived by continuity of the range. As $\Gamma^{(1)}$ is Hausdorff, the equation $r(u)=u$ is deduced implying $u\in\Gamma^{(0)}$.

\vspace{.1cm}

(h): The space $\Gamma^{(1)}\times\Gamma^{(1)}$ is second-countable and so sequences describe the topology. Let $(\gamma_n,\rho_n)_{n\in\NM}$ be a sequence in $\Gamma^{(2)}\subseteq\Gamma^{(1)}\times\Gamma^{(1)}$ tending to $(\gamma,\rho)\in\Gamma^{(1)}\times\Gamma^{(1)}$. The equality $s(\gamma_n)=r(\rho_n)$ holds by (a) for all $n\in\NM$. Thus, $s(\gamma)=r(\rho)$ is deduced by the continuity of $r$ and $s$ as well as the Hausdorff property. Then (a) implies $(\gamma,\rho)\in\Gamma^{(2)}$.

\vspace{.1cm}

(i): This is clear since the product and the subspace topology preserve the desired properties.
\end{proof}

\medskip

Groups and sets are extremal examples of (topological) groupoids as discussed in the following. 

\begin{example}[Group]
\label{Chap2-Ex-Group}
Let $G$ be a (topological) group. Then $G$ is a (topological) groupoid where the set of composable pairs is $G\times G$. In this case, the unit space only contains the neutral element $e$.
\end{example}

\begin{example}[Set]
\label{Chap2-Ex-Set}
Let $X$ be a (topological) set. It gets a (topological) groupoid by defining the set of composable elements by the diagonal $\{(x,y)\in X\times X\;|\; x=y\}$. Then the composition and the inverse are defined by $(x,x)\mapsto x$ and $x^{-1}:=x$. Thus, the equations $\Gamma^{(1)}=\Gamma^{(0)}=X$ hold.
\end{example}

Typical examples for a (topological) groupoid are equivalence relations and transformation group groupoids.

\begin{example}[Equivalence relation]
\label{Chap2-Ex-EquivalRelGroupoid}
Let $X$ be a set. Consider an equivalence relation $\sim$ on $X$, i.e., it is reflexive ($x\sim x$), symmetric ($x\sim y\Leftrightarrow y\sim x$) and transitive ($x\sim y, y\sim z \Rightarrow x\sim z$) where $x,y,z\in X$. Consider the set $\Gamma^{(1)}:=\{(x,y)\in X\times X\;|\; x\sim y\}$ and the set of composable elements $\Gamma^{(2)}\subseteq \Gamma^{(1)}\times \Gamma^{(1)}$ defined by $\big\{ \big((x,y),(y',z)\big) \;\big|\; y=y' \big\}$. Then the composition is given by $\big((x,y),(y',z)\big)\in\Gamma^{(2)}\mapsto (x,z)\in\Gamma^{(1)}$ and the inverse $(x,y)\in\Gamma^{(1)}\mapsto(x,y)^{-1}:=(y,x)$. In this case the unit space is defined by the diagonal, i.e., $\Gamma^{(0)}=\{(x,y)\in\Gamma^{(1)}\;|\; x=y\}$ which is identified with $X$. It is well-known that the class of groupoids induced by an equivalence relation is characterized by the so called principal groupoids, c.f. \cite[Example~I.1.2.~c]{Renault80}. Note that $\Gamma$ is principal if the range $r:\Gamma^{(1)}\to\Gamma^{(0)}$ and the source map $s:\Gamma^{(1)}\to\Gamma^{(0)}$ are injective.
\end{example}

\begin{example}[Transformation group groupoid]
\label{Chap2-Ex-TransformationGroupGroupoid}
In analogy to the semi-direct of groups a continuous action of a group $G$ on a compact space $X$, i.e., a topological dynamical system $(X,G,\alpha)$, defines a new groupoid. Since the so called {\em transformation group groupoid} $\Gamma:=X\rtimes_\alpha G$ is mainly studied in this thesis, it is discussed separately in Section~\ref{Chap2-Sect-TransformationGroupGroupoid}.
\end{example}

It turns out that in case of a transformation group groupoid $X\rtimes_\alpha G$ the unit space is identified with $X$, c.f. Proposition~\ref{Chap2-Prop-TransGrGroupProp} below. Since this class of groupoids is mainly considered in this work, the notation $x$ and $y$ is used for elements of the unit space of a groupoid $\Gamma$ even if $\Gamma$ is not a transformation group groupoid.

\begin{definition}[$r$-fiber, $s$-fiber, $\Gamma$-orbit]
\label{Chap2-Def-r/s-Fibre-Orbit}
Let $\Gamma$ be a groupoid. Then the {\em $r$-fiber} at $x\in \Gamma^{(0)}$ is defined by 
$$
\text{\gls{Gammax}}\; 
	:=\; 
		\left\{ 
			\gamma\in\Gamma^{(1)}\;|\; 
				r(\gamma)=x  
		\right\} \;
	=\; r^{-1}(\{x\}) \, .
$$
Similarly, the $s$-fiber $\Gamma_x$ is defined by $s^{-1}(\{x\})$. Furthermore, the set $r\left(\Gamma_x \right)=s\left(\Gamma^x\right)$ is called {\em $\Gamma$-orbit} $\Orb_\Gamma(x)$ of $x\in\Gamma^{(0)}$.
\end{definition}

every locally compact, Hausdorff group admits a Haar measure, c.f. \cite{Weil40,Loomis53,Querenburg2001}. Such a group $G$ with a Haar measure admits several $C^\ast$-algebras defined via the algebra $\Cc_c(G)$ equipped with the convolution of functions coming from the group action. Two of particular interest are the maximal $C^\ast$-algebra where the completion of $\Cc_c(G)$ is taken over all unitary representations of $G$ and the reduced $C^\ast$-algebra where the completion of $\Cc_c(G)$ is taken over all left-regular representations. Since groupoids are a generalization of groups, a similar strategy is pursued. The existence of an analog of a Haar measure for groupoids is more difficult. Specifically, a topological groupoid does not necessarily admit a so called Haar system. In 1967 {\sc Westman} \cite{Wes67} introduced a notion of a left-invariant continuous system of measures. Later {\sc Seda} \cite[page 2]{Sed76} used a slightly different notion. More precisely, he requires the existence of a measure on the unit space and the continuity condition was dropped. In \cite[Theorem~2]{Sed75} and \cite[Section~4]{Sed76}, {\sc Seda} provides sufficient conditions for the desired continuity property.

\medskip

Let $\Cc_c\big(\Gamma^{(1)}\big)$ be the set of continuous functions $\fz:\Gamma^{(1)}\to\CM$ that have compact support. This means that the support $\text{\gls{supp}}:=\overline{\{\gamma\in\Gamma^{(1)}\;|\; \fz(\gamma)\neq 0\}}$ is a compact subset of $\Gamma^{(1)}$. Following \cite[Definition~I.2.2]{Renault80} a left-continuous Haar system is defined.

\begin{definition}[Left-continuous Haar system]
\label{Chap2-Def-HaarSystem}

A {\em left-continuous Haar system} on a topological groupoid $\Gamma$ is a family of Borel measures $\mu^x:\BG\big(\Gamma^{(1)}\big)\to [0,\infty],\; x\in \Gamma^{(0)},$ on the Borel-$\sigma$-algebra $\BG\big(\Gamma^{(1)}\big)$ such that the assertions {\em \nameref{(H1)},\nameref{(H2)}} and {\em \nameref{(H3)}} hold.
\begin{description}
\item[(H1)\label{(H1)}] For each $x\in\Gamma^{(0)}$, the support 
$$
\supp(\mu^x)\; :=
	\;\big\{
		\gamma\in\Gamma^{(1)}\;\big|\; \
			U\subseteq \Gamma^{(1)} \text{ open with } \gamma\in U\; 
			\Rightarrow\; \mu^x(U)>0 
	\big\}
$$
of the measure $\mu^x$ is equal to the $r$-fiber $\Gamma^x$.
\item[(H2)\label{(H2)}] For each $\fz\in \Cc_c\big(\Gamma^{(1)}\big)$, the function  $\Gamma^{(0)}\to\CM$ defined by
$$
x\; \longmapsto\; \mu^x(\fz)\;:=\;\int_{\Gamma^x} \fz(\gamma)\; d\mu^x(\gamma)
$$ 
is continuous.
\item[(H3)\label{(H3)}] The measures $\mu:=(\mu^x)_{x\in\Gamma^{(0)}}$ are left-invariant, i.e., the equation
$$
\int_{\Gamma^{s(\gamma)}} 
	\fz(\gamma\cdot\varrho)
\; d\mu^{s(\gamma)}(\varrho) \; 
		= \; \int_{\Gamma^{r(\gamma)}} 
			\fz(\varrho)
		\; d\mu^{r(\gamma)}(\varrho)
$$
holds for each $\gamma\in\Gamma^{(1)}$ and $\fz\in \Cc_c\big(\Gamma^{(1)}\big)$.
\end{description}
\end{definition}

\begin{example}
\label{Chap2-Ex-GroupHaarSyst}
Let $G$ be a locally compact, second countable, Hausdorff group. It is well-known that $G$ admits a Haar measure $\lambda$ on the Borel-$\sigma$-algebra of $G$, c.f. \cite{Weil40,Loomis53,Querenburg2001}. According to Example~\ref{Chap2-Ex-Group}, the group $G$ is a groupoid. By definition of a Haar measure the support $\supp(\lambda)$ is equal to $G$ and so \nameref{(H1)} is satisfied. Furthermore, Condition~\nameref{(H2)} holds as the unit space of $G$ is equal to $\{e\}$ where $e\in G$ is the neutral element. Finally, \nameref{(H3)} follows by the left invariance of the Haar measure $\lambda$. Thus, $\lambda$ is a left-continuous Haar system on $G$.
\end{example}

\begin{example}
\label{Chap2-Ex-SetHaarSyst}
Let $X$ be a locally compact, second countable, Hausdorff space. It defines a topological groupoid $\Gamma$ along the lines of Example~\ref{Chap2-Ex-Set}. Then the equations $\Gamma^{(1)}=\Gamma^{(0)}=X$ hold. Define the family of Borel measures $\mu^x:=\delta_x,\; x\in X,$ where $\delta_x$ is the Dirac measure at $x$. Due to definition of the composition and the inverse, the $r$-fiber $\Gamma^x$ coincides with the set $\{x\}$, namely $(\mu^x)_{x\in X}$ fulfills \nameref{(H1)}. Then $\mu^x(\fz)=\fz(x)$ holds for $x\in X$ and $\fz\in\Cc_c\big(\Gamma^{(1)}\big)$. Since $\fz$ is a continuous function, \nameref{(H2)} follows. Additionally, the identities
$$
\int_{\Gamma^{s(x)}} \fz(x\cdot y)\; d\mu^{s(x)}(y) \;
	= \; \fz(x) \;
	= \; \int_{\Gamma^{r(x)}} \fz(y)\; d\mu^{r(x)}(y)
$$
hold for every $x\in \Gamma^{(1)}$ implying that $(\mu^x)_{x\in X}$ defines a left-invariant Haar system.
\end{example}

\begin{example}
\label{Chap2-Ex-EquivalRelGroupoidHaarSyst}
Let $X$ be a set equipped with the discrete topology and an equivalence relation $\sim$ on $X$. Consider the topological groupoid $\Gamma$ defined by the equivalence relation with unit space $\Gamma^{(0)}=\{(x,x)\;|\; x\in X\}$, c.f. Example~\ref{Chap2-Ex-EquivalRelGroupoid}. Then the range and the source map are given by $r(x,y)=(x,x)$ and $s(x,y)=(y,y)$ for $(x,y)\in\Gamma^{(1)}$. The $r$-fiber $\Gamma^{(x,x)}$ is equal to $\{(x,z)\;|\; z\sim x\}$ for all $(x,x)\in\Gamma^{(0)}$. Consider the family $\mu^{(x,x)}:=\sum_{y\in X\,, y\sim x}\delta_x\times\delta_y\,,\; (x,x)\in\Gamma^{(0)}\,,$ of Borel measures on $\Gamma^{(1)}$. Clearly, the support of the Borel measure $\mu^{(x,x)}$ is equal to $\Gamma^{(x,x)}$ for $x\in X$. Hence, $\big(\mu^{(x,x)}\big)_{(x,x)\in\Gamma^{(0)}}$ satisfies \nameref{(H1)}. Since $X$ is discrete, every compactly supported, continuous function $\fz\in\Cc_c\big(\Gamma^{(1)}\big)$ has finite support. Thus, $\mu^{(x,x)}(\fz)=\sum_{y\in X\,, y\sim x} \fz(x,y)$ is a finite sum of continuous functions for each $(x,x)\in\Gamma^{(0)}$. Consequently, the map $(x,x)\in\Gamma^{(0)}\mapsto\mu^{(x,x)}(\fz)$ is continuous for $\fz\in\Cc_c\big(\Gamma^{(1)}\big)$ and so the Borel measures $\mu^{(x,x)}\,,\; (x,x)\in\Gamma^{(0)}$ fulfill \nameref{(H2)}. Finally, \nameref{(H3)} follows by a short computation
$$
\sum_{z\sim y} f\big((x,y)(y,z)\big) \;
	= \; \sum_{z\sim y} f(x,y)\;
	\overset{\sim \text{ transitiv}}{=} \; \sum_{z\sim x} \fz(x,z)\,,
	\qquad \gamma:=(x,y)\in\Gamma^{(1)} \,.
$$ 
Consequently, $\big(\mu^{(x,x)}\big)_{(x,x)\in\Gamma^{(0)}}$ is a left-continuous Haar system for the topological group\-oid $\Gamma$ defined by the equivalence relation $\sim$ on $X$.
\end{example}

\begin{example}
\label{Chap2-Ex-TransformationGroupGroupoidHaarSyst}
Let $G$ be a topological group and $X$ a topological space such that both spaces are locally compact, second countable and Hausdorff. If $\alpha$ defines a continuous action of $G$ on $X$, i.e., $(X,G,\alpha)$ is a topological dynamical system, then the transformation group groupoid $X\rtimes_\alpha G$ defines a topological groupoid, c.f. Proposition~\ref{Chap2-Prop-TransGrGroupProp} below. A left-continuous Haar system of $X\rtimes_\alpha G$ is constructed by the Haar system of $G$ and $X$ defined in Example~\ref{Chap2-Ex-GroupHaarSyst} and Example~\ref{Chap2-Ex-SetHaarSyst}. In detail, the family of Borel measures $\delta_x\times\lambda,\; x\in X,$ where $\lambda$ is the Haar measure of $G$, defines a left-continuous Haar system on the groupoid $X\rtimes_\alpha G$, c.f. Proposition~\ref{Chap2-Prop-TransGroupAmen} below.
\end{example}

\begin{remark}
\label{Chap2-Rem-ContHaarMeasSeda}
It is proven in \cite{Sed86} that for a family of Borel measures $\mu^x,\; x\in\Gamma^{(0)},$ satisfying \nameref{(H1)} and \nameref{(H3)} the continuity condition \nameref{(H2)} is equivalent to the fact that $\fz\star\gz$, defined in Proposition~\ref{Chap2-Prop-InvolutiveGroupoidAlgebra}, is a continuous function for all $\fz,\gz\in\Cc_c\big(\Gamma^{(1)}\big)$. This justifies to add the continuity condition in the definition of a Haar system since $\Cc_c\big(\Gamma^{(1)}\big)$ shall define a dense sub-$\ast$-algebra of the corresponding $C^\ast$-algebras.
\end{remark}

Recall the notion of a $C^\ast$-algebra in Definition~\ref{Chap3-Def-C-algebra}. The construction of a groupoid $C^\ast$-algebra is described next which is well-known in the literature. The interested reader is also referred to \cite{Renault80,Bun06} and references therein. For the following definition, recall that the space of bounded linear operators $\Ll(\hs)$ on a Hilbert space $\hs$ is equipped with the composition for the multiplication and the adjoint $^\ast$ of an operator for the involution, Example~\ref{Chap3-Ex-LinBounOpCalgebra}.

\begin{definition}[Representation]
\label{Chap2-Def-Representation}
Let $\AG$ be an $\ast$-algebra. Then a pair $(\pi,\hs)$ is called a {\em $\ast$-representation of $\AG$} if $\hs$ is a Hilbert space and $\pi:\AG\to\Ll(\hs)$ is a $\ast$-homomorphism, i.e.,
\begin{itemize}
\item[(i)] $\pi$ is linear;
\item[(ii)] $\pi$ is multiplicative, i.e., $\pi(\fz\star\gz)=\pi(\fz)\pi(\gz)$ holds for all $\fz,\gz\in\AG$;
\item[(iii)] $\pi$ preserves the involution, i.e., $\pi(\fz^\ast)=\pi(\fz)^\ast$ holds for all $\fz\in\AG$ where $\pi(\fz)^\ast$ is the adjoint operator of $\pi(\fz)\in\Ll(\hs)$.
\end{itemize}
A family of representations $(\pi^x,\hs_x)_{x\in X}$ is called {\em faithful} whenever the family is injective, i.e., $\fz=0$ if and only if $\pi^x(\fz)=0$ for all $x\in X$.
\end{definition}

\begin{remark}
\label{Chap2-Rem-FaithfulRepresent}
It is well-known that an injective $\ast$-homomorphism between unital $C^\ast$-algebras is automatically isometric, c.f. \cite[Theorem~3.1.5]{Murphy90}. Whenever $(\pi^x,\hs_x)_{x\in X}$ is a faithful family of representations, a faithful representation is constructed by taking the direct sum $\hs:=\oplus_{x\in X} \;\hs_x$ and $\pi:=\oplus_{x\in X} \;\pi^x$. Due to the fact that a faithful representation defines an isometric, $\ast$-homomorphism, an element $\fz\in\AG$ is invertible if and only if $\pi^x(\fz)$ is invertible for all $x\in X$ and the norm of the resolvent is uniformly bounded which is well-known, see e.g. \cite[Proposition~2.5]{NiPr15}. Clearly, this leads to a correspondence of the spectra. More precisely, a normal element $\fz\in\AG$ satisfies 
$$
\sigma(\fz) \; 
	= \; \overline{\bigcup_{x\in X}\sigma(\pi^x(\fz))}
$$
if and only if the family of representations $(\pi^x,\hs_x)_{x\in X}$ is faithful, see e.g. \cite[Proposition~2.8]{NiPr15}. It was known for specific examples that the equation holds without the closure on the right hand side. {\sc Roch} \cite{Roch2003} introduced in 2003 the concept of invertibility sufficient for a family of representation $(\pi^x,\hs_x)_{x\in X}$ being exactly those representation satisfying
$$
\sigma(\fz) \;
	= \; \bigcup_{x\in X}\sigma(\pi^x(\fz))\,,
$$
c.f. \cite[Theorem~3.6]{NiPr15}. This notion turns out to be equivalent to strictly norming representations $(\pi^x,\hs_x)_{x\in X}$, c.f. \cite[Theorem~5.7]{Roch2003}. A family of representations $(\pi^x,\hs_x)_{x\in X}$ is called {\em strictly norming} if $\|\fz\|=\sup_{x\in X}\|\pi^x(\fz)\|$ is a maximum for every $\fz\in\AG$, i.e., there exists an $x\in X$ such that $\|\fz\|=\|\pi^x(\fz)\|$. The reader is referred to \cite{NiPr15} for a more detailed discussion.
\end{remark}

\medskip

The following result is well-known and can be found in \cite[Proposition~II.1.1.]{Renault80}. Note that {\sc Renault} used a slightly different definition of the multiplication. However, the following definition is usually used in the literature.

\begin{proposition}[\cite{Renault80}]
\label{Chap2-Prop-InvolutiveGroupoidAlgebra}
Let $\Gamma$ be a topological groupoid with left-continuous Haar system $\mu:=(\mu^x)_{x\in \Gamma^{(0)}}$. Then, the vector space $\Cc_c\big(\Gamma^{(1)}\big)$ equipped with the maps
\begin{align*}
	\begin{array}{rl}
		(\fz\star \gz)(\gamma)\; 
			& := \;\displaystyle\int_{\Gamma^{r(\gamma)}} \fz(\rho)\; \gz(\rho^{-1}\cdot\gamma)\; d\mu^{r(\gamma)}(\rho) \, ,\\[0.5cm]
		\fz^*(\gamma)\; 
			& := \; \overline{\fz(\gamma^{-1})} \, ,
	\end{array}
	\qquad \fz,\gz\in\Cc_c\big(\Gamma^{(1)}\big),\; \gamma\in\Gamma^{(1)} \, ,
\end{align*}
defines an $\ast$-algebra.
\end{proposition}

\begin{proof}
The space $\Cc_c\big(\Gamma^{(1)}\big)$ gets a vector space equipped with the pointwise addition and multiplication by scalars. The proof is organized as follows.
\begin{itemize}
\item[(i)] The map $\star:\Cc_c\big(\Gamma^{(1)}\big)\times\Cc_c\big(\Gamma^{(1)}\big)\to\Cc_c\big(\Gamma^{(1)}\big)$ is well-defined and $\Cc_c\big(\Gamma^{(1)}\big)$ is an algebra with this multiplication.
\item[(ii)] The map $\ast:\Cc_c\big(\Gamma^{(1)}\big)\to \Cc_c\big(\Gamma^{(1)}\big)$ defines an involution on $\Cc_c\big(\Gamma^{(1)}\big)$.
\end{itemize}

\vspace{.1cm}

(i): Let $\fz,\gz\in\Cc_c\big(\Gamma^{(1)}\big)$. As they have compact support, the integral $(\fz\star \gz)(\gamma)$ is well-defined for every $\gamma\in\Gamma^{(1)}$. By definition the value $(\fz\star \gz)(\gamma)$ is non-zero only if there exists a $\rho\in\Gamma^{(1)}$ such that $\fz(\rho)\gz(\rho^{-1}\gamma)\neq 0$. This holds only if $\gz(\rho^{-1}\gamma)\neq 0$. Clearly, $\rho^{-1}\gamma$ is not contained in $\supp(\gz)$ for $\rho\in\Gamma^{r(\gamma)}$  if $\gamma\not\in\supp(\fz)\cdot\supp(\gz)$. Thus, the support 
$$
\supp(\fz\star \gz) \;
	:= \;
	\overline{
		\left\{ 
			\gamma\in\Gamma^{(1)} \;|\; 
				(\fz\star \gz)(\gamma)\neq 0
		\right\}
	}
$$
is contained in the compact set $\supp(\fz)\cdot\supp(\gz)$. As $\supp(\fz)$ is compact by assumption and $\supp(\fz\star \gz)$ is a closed subset, the support of $\fz\star \gz$ is compact as well. In order to prove the continuity of $\fz\star \gz$ the same trick as in \cite{Con78,Renault80} is used. According to Proposition~\ref{Chap2-Prop-BasGroupoid}~(h) the set $\Gamma^{(2)}$ is closed in $\Gamma^{(1)}\times\Gamma^{(1)}$. As $\Gamma^{(1)}\times\Gamma^{(1)}$ is second-countable, locally compact and Hausdorff, it is a normal space. Consequently, Tietze's Extension Theorem (\cite{Tie15}, \cite[Anhang~III.28.]{Ury25}) is applicable to the continuous map $(\rho,\gamma)\in\Gamma^{(2)}\mapsto \gz(\rho\gamma)$ with compact support. Specifically, there exists a continuous bounded map $\Gz:\Gamma^{(1)}\times\Gamma^{(1)}\to\CM$ such that $\Gz(\rho,\gamma)=\gz(\rho\gamma)$ whenever $(\rho,\gamma)\in\Gamma^{(2)}$. Define the map $\hz:\Gamma^{(1)}\times\Gamma^{(1)}\to\CM$ by $\hz(\rho,\gamma):=\Gz(\rho^{-1},\gamma)\fz(\rho)$. Since $^{-1}:\Gamma^{(1)}\to\Gamma^{(1)}$ is continuous, the map $\hz$ is continuous in both components. Furthermore, for each fixed $\gamma\in\Gamma^{(1)}$, the map $\hz(\cdot,\gamma):\Gamma^{(1)}\to\CM$ has compact support. It is standard to show that $\gamma\mapsto\int \hz(\rho,\gamma)\; d\mu^x(\rho)$ is continuous for each fixed $x\in\Gamma^{(0)}$. Thus, the map 
$$
(\gamma,x)\in\Gamma^{(1)}\times\Gamma^{(0)}\mapsto \int_{\Gamma} \hz(\rho,\gamma)\; d\mu^x(\rho)\in\CM
$$ 
is continuous by \nameref{(H2)}. Then the restriction of this map to the set $\{(\gamma,x)\in\Gamma^{(1)}\times\Gamma^{(0)}\;|\; r(\gamma)=x\}$ is still continuous. By construction, this restriction is equal to $\fz\star \gz$. Consequently, $\fz\star \gz$ is an element of $\Cc_c\big(\Gamma^{(1)}\big)$ and so the multiplication is well-defined.

\vspace{.1cm}

It suffices to check the associativity of $\star$ as the other conditions of an algebra immediately follow by the linearity of the integral. Let $\fz,\gz$ and $\hz$ be elements of $\Cc_c\big(\Gamma^{(1)}\big)$. Then the Theorem of Fubini-Tonelli 
for continuous, compactly supported functions amounts to
\begin{align*}
\big( (\fz\star \gz) \star \hz \big)(\gamma) \; 
	&\underset{r(\gamma)=r(\rho)}{\overset{\rho\,\in\,\Gamma^{r(\gamma)}}{=}} \; \int_{\Gamma^{r(\gamma)}} \, 
			\left(
				\int_{\Gamma^{r(\gamma)}} \,
					\fz(\zeta) \cdot \gz(\zeta^{-1}\rho) \cdot \hz(\rho^{-1}\gamma) 
				\; d\mu^{r(\gamma)}(\zeta)
			\right)
		\; d\mu^{r(\gamma)}(\rho)\\
	&\underset{r(\gamma)=r(\zeta)}{\overset{\text{Fubini}}{=}} \; \int_{\Gamma^{r(\gamma)}} \, 
			\left(
				\int_{\Gamma^{r(\zeta)}} \,
					\fz(\zeta) \cdot \gz(\zeta^{-1}\rho) \cdot \hz(\rho^{-1}\gamma) 
				\; d\mu^{r(\zeta)}(\rho) 
			\right)
		\; d\mu^{r(\gamma)}(\zeta)\\
	&\overset{\nameref{(H3)}}{\,\quad=\,\quad} \; \int_{\Gamma^{r(\gamma)}} \, 
			\left(
				\int_{\Gamma^{s(\zeta)}} \,
					\fz(\zeta) \cdot \gz\big(\zeta^{-1}(\zeta\rho)\big) \cdot \hz\big((\zeta\rho)^{-1}\gamma\big) 
				d\mu^{s(\zeta)}(\rho) 
			\right)
		d\mu^{r(\gamma)}(\zeta)\\
	&\overset{(\bigstar)}{\,\quad=\,\quad} \; \int_{\Gamma^{r(\gamma)}} \, \fz(\zeta) \cdot
			\left(
				\int_{\Gamma^{r(\zeta^{-1}\gamma)}} \, 
					\gz(\rho) \cdot \hz\big(\rho^{-1}\zeta^{-1}\gamma\big) 
				\, d\mu^{r(\zeta^{-1}\gamma)}(\rho) 
			\right) 
		d\mu^{r(\gamma)}(\zeta)\\	
	&\,\quad=\,\quad \; \int_{\Gamma^{r(\gamma)}} \, \fz(\zeta) \cdot 
			(\gz\star \hz)\big(\zeta^{-1}\gamma\big) 
		\; d\mu^{r(\gamma)}(\zeta)\\
	&\,\quad=\,\quad \; \big( \fz \star (\gz \star \hz) \big)(\gamma)
\end{align*}
for each $\gamma\in\Gamma^{(1)}$ where the equation $(\bigstar)$ follows from Proposition~\ref{Chap2-Prop-BasGroupoid}~(b) and (d).

\vspace{.1cm}

(ii): For $\fz\in\Cc_c\big(\Gamma^{(1)}\big)$, the support $\supp(\fz^\ast)$ of $\fz^\ast$ is equal to $\supp(\fz)^{-1}$. According to Proposition~\ref{Chap2-Prop-BasGroupoid}~(e) the inverse map is a homeomorphism and so $\supp(\fz^\ast)$ is a compact set. The map $f^\ast:\Gamma^{(1)}\to\CM$ is represented by the composition of $\fz$, the complex conjugation and the inverse on the groupoid. Hence, it is a composition of continuous maps if $\fz\in\Cc_c\big(\Gamma^{(1)}\big)$. Thus, $^\ast$ maps $\Cc_c\big(\Gamma^{(1)}\big)$ to $\Cc_c\big(\Gamma^{(1)}\big)$.

\vspace{.1cm}

For $\fz,\gz\in\Cc_c\big(\Gamma^{(1)}\big)$ the equations
\begin{align*}
(\fz+\lambda \gz)^\ast(\gamma)\; 
	= \; &&\overline{\fz(\gamma^{-1})} + \overline{\lambda\cdot \gz(\gamma^{-1})} \;
	= \; && \fz^\ast + \overline{\lambda}\cdot \gz^\ast(\gamma) \, ,\\
\big(\fz^\ast\big)^\ast (\gamma)\;
	= \; &&\overline{\overline{\fz\big( (\gamma^{-1})^{-1} \big)}} \;
	= \; && \fz(\gamma) \, ,
\end{align*}
hold for all $\gamma\in\Gamma^{(1)}$. By using once more Proposition~\ref{Chap2-Prop-BasGroupoid}~(c) and \nameref{(G3)}, the equalities
\begin{align*}
(\fz\star \gz)^\ast(\gamma) \; 
	&\,\;=\,\; \; \int_{\Gamma^{r(\gamma^{-1})}}
				\overline{\fz(\rho)\cdot \gz\big( \rho^{-1}\gamma^{-1} \big)}
			\; d\mu^{r(\gamma^{-1})}(\rho)\\
	&\overset{\nameref{(H3)}}{=} \; \int_{\Gamma^{s(\gamma^{-1})}}
				\overline{\fz\big( \gamma^{-1}\rho \big)\cdot \gz\big( (\gamma^{-1}\rho)^{-1}\gamma^{-1} \big)}
			\; d\mu^{s(\gamma^{-1})}(\rho)\\
	&\,\;=\,\; \; \int_{\Gamma^{r(\gamma)}}
				\gz^\ast(\rho) \cdot \fz^\ast\big( \rho^{-1}\gamma \big)
			\; d\mu^{r(\gamma)}(\rho)\\
	&\,\;=\,\; \; (\gz^\ast\star \fz^\ast)(\gamma)
\end{align*}
follow for $\fz,\gz\in\Cc_c\big(\Gamma^{(1)}\big)$ and $\gamma\in\Gamma^{(1)}$.
\end{proof}

\medskip

Let $\Gamma$ be a topological groupoid with left-continuous Haar system $(\mu^x)_{x\in\Gamma^{(0)}}$. Then $\mu_x:=\mu^x\circ^{-1},\; x\in\Gamma^{(0)},$ defines a right-continuous Haar system, where $^{-1}:\Gamma^{(1)}\to\Gamma^{(1)}$ is the inverse in the groupoid $\Gamma$. This is used to define the so called $I$-norm.

\medskip

{\sc Hahn} \cite[Page~38]{Hah78} introduced for a measured groupoid the so called $I$-norm on $\Cc_c\big(\Gamma^{(1)}\big)$. {\sc Renault} \cite[Page~50]{Renault80} changed the notation in accordance with his setting. In detail, consider a locally compact, Hausdorff groupoid $\Gamma$ with left-continuous Haar system $(\mu^x)_{x\in\Gamma^{(0)}}$. The norm $\|\cdot\|_{I,r}$ is defined by
$$
\|\fz\|_{I,r}\; 
	:=\; \sup_{x\in\Gamma^{(0)}} \mu^x(|\fz|)
	=\; \sup_{x\in\Gamma^{(0)}} \int_{\Gamma^x} |\fz(\gamma)|\; d\mu^x(\gamma)
		,\qquad \fz\in\Cc_c\big(\Gamma^{(1)}\big) \, .
$$
Note that $\mu^x(|\fz|)$ exists by definition for each $x\in\Gamma^{(0)}$. Furthermore, the map $x\mapsto\mu^x(|\fz|)$ has support $r(\supp(\fz))$. This set is compact since the range $r$ is continuous and $\supp(\fz)$ is compact. Thus, the norm $\|\fz\|_{I,r}$ is finite by \nameref{(H2)} for all $\fz\in\Cc_c\big(\Gamma^{(1)}\big)$. Similarly, the norm $\|\cdot\|_{I,s}$ is given by 
$$
\|\fz\|_{I,s}\; 
	:=\; \sup_{x\in\Gamma^{(0)}} \mu_x(|\fz|)
	=\; \sup_{x\in\Gamma^{(0)}} \int_{\Gamma^x} |\fz(\gamma^{-1})|\; d\mu^x(\gamma)
		,\qquad \fz\in\Cc_c\big(\Gamma^{(1)}\big) \, .
$$
By analog arguments as above, the norm $\|\fz\|_{I,s}$ is finite for all $\fz\in\Cc_c\big(\Gamma^{(1)}\big)$.

\begin{definition}[I-norm]
\label{Chap2-Def-Inorm}
The {\em $I$-norm} of $\fz\in\Cc_c\big(\Gamma^{(1)}\big)$ is defined by the maximum $\|\fz\|_I:=\max\{ \|\fz\|_{I,r},\; \|\fz\|_{I,s} \}$.
\end{definition}

The notation becomes clearer if one considers \'etale groupoid (Definition~\ref{Chap2-Def-Etale}) where the Haar measure is essentially the counting measure. Then the equations 
$$
\|\fz\|_{I,r}\; 
	=\; \sup_{x\in\Gamma^{(0)}}
		\sum_{x=r(\gamma)} |f(\gamma)| \cdot \mu^x(\{x\}) \, , \qquad 
\|\fz\|_{I,s}\; 
	= \; \sup_{x\in\Gamma^{(0)}}
		\sum_{x=s(\gamma)} |f(\gamma)| \cdot \mu^x(\{x\}) \, ,
$$
hold. Note that $\|\cdot\|_{I}$ defines a $\ast$-norm on $\Cc_c\big(\Gamma^{(1)}\big)$, i.e., $\|\fz\star\gz\|_I\leq\|\fz\|_I\cdot\|\gz\|_I$ and $\|\fz^\ast\|_I=\|\fz\|_I$ hold for all $\fz,\gz\in\Cc_c\big(\Gamma^{(1)}\big)$, c.f. \cite[Proposition~II.1.4]{Renault80}.

\begin{proposition}[\cite{Renault80}]
\label{Chap2-Prop-LeftRegulRepres}
Let $\Gamma$ be a topological groupoid with left-continuous Haar system $\mu:=(\mu^x)_{x\in \Gamma^{(0)}}$. Define the maps $\pi^x:\Cc_c\big(\Gamma^{(1)}\big)\to\Ll\big( L^2(\Gamma^x,\mu^x)\big),\; x\in\Gamma^{(0)},$ by
$$
\left(\pi^x(\fz\,)\psi\right)(\gamma) \; 
	:= \; \int_{\Gamma^x} \fz(\gamma^{-1}\cdot\rho)\; \psi(\rho)\;d\mu^{x}(\rho)
$$
where $\fz\in\Cc_c\big(\Gamma^{(1)}\big),\; \gamma\in\Gamma^x$ and $\psi\in L^2(\Gamma^x,\mu^x)$. Then $(\pi^x)_{x\in\Gamma^{(0)}}$ defines a family of faithful representations. In particular, $\pi^x(\fz)$ defines a bounded operator on $L^2(\Gamma^x,\mu^x)$ for each $\fz\in\Cc_c\big(\Gamma^{(1)}\big)$ and $x\in\Gamma^{(0)}$. The operator norm of $\pi^x(\fz)$ is uniformly bounded in $x\in\Gamma^{(0)}$ by $\|\fz\|_I$.
\end{proposition}

\begin{proof}
For the proof let $\fz,\gz\in\Cc_c\big(\Gamma^{(1)}\big)$, $\psi,\varphi\in L^2(\Gamma^x,\mu^x)$ where $x\in\Gamma^{(0)}$ is arbitrarily chosen. 

\vspace{.1cm}

For $\gamma\in\Gamma^x$, the function $\gz:\Gamma^x\to\CM,\; \rho\mapsto f(\gamma^{-1}\rho),$ is an element of $\Cc_c\big(\Gamma^x\big)\subseteq L^2(\Gamma^x,\mu^x)$. By H\"olders inequality the integral $\int_{\Gamma^x} \gz(\rho)\; \psi(\rho)\;d\mu^{x}(\rho)$ is well-defined and equal to $\left(\pi^x(\fz)\psi\right)(\gamma)$. Furthermore, the estimates
\begin{align*}
\langle \pi^x(\fz)\psi\;|\; \pi^x(\fz)\psi\rangle \; 
	&= \; \int_{\Gamma^x} 
		\left|
			\left( \pi^x(\fz)\psi \right) (\gamma)
		\right|^2
		\; d\mu^x(\gamma)\\
	&\leq \; \int_{\Gamma^x} 
		\left(
			\int_{\Gamma^x}
				\left| \fz \left(\gamma^{-1}\rho \right) \right|^{\frac{1}{2}} 
				\left(
					\left| \fz \left(\gamma^{-1}\rho \right) \right|^{\frac{1}{2}}
					|\psi(\rho)|
				\right)
			\; d\mu^x(\rho)
		\right)^2
		\; d\mu^x(\gamma)\\
	&\leq \; \int_{\Gamma^x} 
		\underbrace{
			\left(
				\int_{\Gamma^x}
					\left| \fz \left(\gamma^{-1}\rho \right) \right|
				d\mu^x(\rho)
			\right)
		}_{=:C_1}
		\left(
			\int_{\Gamma^x}
				\left| \fz \left(\gamma^{-1}\zeta \right) \right|
				|\psi(\zeta)|^2
			d\mu^x(\zeta)
		\right)
		d\mu^x(\gamma)\\
	&\leq \; C_1\cdot 
		\int_{\Gamma^x} 
			|\psi(\zeta)|^2 \cdot
			\underbrace{
				\left(
					\int_{\Gamma^x} 
						\left| \fz \left(\gamma^{-1}\zeta \right) \right|
					d\mu^x(\gamma)
				\right)
			}_{=:C_2}
		d\mu^x(\zeta)\\		
	&= \; C_1 \cdot C_2 \cdot \|\psi\|^2
\end{align*}
hold by applying H\"olders inequality and the Theorem of Fubini-Tonelli. Using the invariance \nameref{(H3)}, the integral $C_1$ is bounded by
$$
C_1 \;
	= \; \int_{\Gamma^{s(\gamma)}} \left| \fz(\rho) \right| d\mu^{s(\gamma)}(\rho) \; 
	\leq \; \sup\limits_{y\in\Gamma^{(0)}} \int_{\Gamma^y} \left| \fz(\rho) \right| d\mu^y(\rho) \;
	= \; \sup\limits_{y\in\Gamma^{(0)}} \mu^y(|\fz|)
	=: \; \|\fz\|_{I,r} \, .
$$
Similarly, the estimate $C_2\leq\|\fz\|_{I,s}$ holds. As discussed in the definitions of $\|\cdot\|_{I,r}$ and $\|\cdot\|_{I,s}$, these norms are finite for all $\fz\in\Cc_c\big(\Gamma^{(1)}\big)$. 
Altogether, the estimate 
$$
\|\pi^x(\fz)\psi\| \;
	\leq \; \|\fz\|_{I,r}^{\frac{1}{2}} \cdot \|\fz\|_{I,s}^{\frac{1}{2}} \cdot \|\psi\| \;
	\leq \; \|\fz\|_I \cdot \|\psi\|
$$
holds proving that $\pi^x(\fz)$ is uniformly bounded in operator norm by $\|\fz\|_I< \infty$. As the integral is linear, the operator $\pi^x(\fz)$ is linear. Thus, $\pi^x(\fz)$ is well-defined.

\vspace{.1cm}

The map $\pi^x:\Cc_c\big(\Gamma^{(1)}\big)\to\Ll\left( L^2(\Gamma^x,\mu^x) \right)$ is linear by the linearity of the integral. The multiplicativity is proven in analogy to the associativity in Proposition~\ref{Chap2-Prop-InvolutiveGroupoidAlgebra} by using \nameref{(H3)}. Specifically, for $\gamma\in\Gamma^x$ the equations
\begin{align*}
\left( \pi^x(\fz\star\gz) \psi \right) (\gamma) \;
	&= \; \int_{\Gamma^x}
			\psi(\rho)\cdot
			\left(
				\int_{\Gamma^{r(\gamma^{-1}\rho)}}
					\fz(\zeta) \cdot \gz\left( \zeta^{-1}\gamma^{-1}\rho \right)
				d\mu^{r(\gamma^{-1}\rho)}(\zeta)
			\right)
		d\mu^x(\rho)\\
	&= \; \int_{\Gamma^{r(\gamma)}}
			\fz(\gamma^{-1}\zeta)\cdot
			\left(
				\int_{\Gamma^x}
					\gz\left( \zeta^{-1}\rho \right) \psi(\rho)
				d\mu^x(\rho)
			\right)
		d\mu^{r(\gamma)}(\zeta)\\
	&= \; \Big(\pi^x(\fz)\, \big( \pi^x(\gz)\psi \big) \Big)(\gamma)
\end{align*}
are derived by using the Theorem of Fubini-Tonelli and Proposition~\ref{Chap2-Prop-BasGroupoid}~(b),(c),(d). Finally, a short computation leads to
\begin{align*}
\langle \pi^x(\fz^\ast)\psi \;|\; \varphi \rangle \;
	&= \; \int_{\Gamma^x}
			\left(			
				\int_{\Gamma^x}
					\overline{\fz^\ast\left( \gamma^{-1}\rho \right)} \cdot \overline{\psi(\rho)}
				\; d\mu^x(\rho)
			\right)
			\cdot \varphi(\gamma)
		\; d\mu^x(\gamma)\\
	&= \; \int_{\Gamma^x}
			\overline{\psi(\rho)} \cdot
			\left(			
				\int_{\Gamma^x}
					\fz\left( \rho^{-1}\gamma \right) \cdot \varphi(\gamma)
				\; d\mu^x(\gamma)
			\right)
		\; d\mu^x(\rho)\\
	&= \; \langle \psi \;|\; \pi^x(\fz)\varphi \rangle \, .
\end{align*}
Thus, $\left(\pi^x,L^2(\Gamma^x,\mu^x)\right)$ defines a representation of $\Cc_c\big(\Gamma^{(1)}\big)$. The faithfulness is proven in \cite[Proposition~II.1.11]{Renault80} by using an approximate identity. The existence of an approximate identity is provided in \cite[Proposition~II.1.9]{Renault80}.
\end{proof}

\medskip

Proposition~\ref{Chap2-Prop-LeftRegulRepres} justifies the following definition.

\begin{definition}[Left-regular representation]
\label{Chap2-Def-LeftRegulRepres}
Let $\Gamma$ be a topological groupoid with left-continuous Haar system $(\mu^x)_{x\in \Gamma^{(0)}}$. The family $\text{\gls{pix}}:\Cc_c\big(\Gamma^{(1)}\big)\to\Ll\left( L^2(\Gamma^x,\mu^x)\right),\; x\in\Gamma^{(0)},$ of representations defined by
$$
\left(\pi^x(\fz\,)\psi\right)(\gamma) \; 
	:= \; \int_{\Gamma^x} \fz(\gamma^{-1}\cdot\rho)\; \psi(\rho)\;d\mu^{x}(\rho)\,,
	\qquad
	\fz\in\Cc_c\big(\Gamma^{(1)}\big),\, \gamma\in\Gamma^x,\, \psi\in L^2(\Gamma^x,\mu^x),
$$
is called {\em left-regular representation}.
\end{definition}

\begin{remark}
\label{Chap2-Rem-LeftRegulRepresNonDegenerated}
Note that the left-regular representation defines a family of non-degen\-erated representations, i.e., the set
$$
\big\{ \psi\in L^2(\Gamma^x,\mu^x) \;\big|\;
	\pi^x(\fz)\psi =0 \text{ for all } \fz\in\Cc_c\big(\Gamma^{(1)}\big)
\big\}
$$
is equal to $\{0\}$ for each $x\in\Gamma^{(0)}$, c.f. \cite[Proposition~II.1.7]{Renault80}. 
\end{remark}

Using this representation one can define the \emph{reduced-norm} on $\Cc_c\big(\Gamma^{(1)}\big)$ given by
$$
\|\fz\|_{red}\; := \;
	\sup
		\big\{	
			\|\pi^x(\fz)\|\;\big|\; x\in\Gamma^{(0)}
		\big\},
	\qquad \fz\in\Cc_c\big(\Gamma^{(1)}\big) \, ,
$$
where $\|\pi^x(f)\|$ denotes the operator norm in $L^2(\Gamma^x,\mu^x)$. Since the operator norm on a Hilbert space is a $C^*$-norm, c.f. Example~\ref{Chap3-Ex-LinBounOpCalgebra}, it is not difficult to check that \gls{normreduced} defines a $C^\ast$-semi-norm. According to Proposition~\ref{Chap2-Prop-LeftRegulRepres}, the family $(\pi^x)_{x\in\Gamma^{(0)}}$ defines a faithful representation, c.f. \cite[Proposition~II.1.11]{Renault80}. Thus, $\|\cdot\|_{red}$ becomes a $C^\ast$-norm. Then the {\em reduced (groupoid) $C^\ast$-algebra} \gls{CGred} is the completion of $\Cc_c\big(\Gamma^{(1)}\big)$ with respect to the reduced-norm $\|\cdot\|_{red}$, c.f. \cite[Definition~II.2.8]{Renault80}.

\medskip 

Along the lines of \cite[Definition~II.1.5]{Renault80}, a representation $\pi$ of $\Cc_c\big(\Gamma^{(1)}\big)$ is called {\em bounded} whenever $\|\pi(\fz)\|\leq \|\fz\|_{I}$ for all $\fz\in\Cc_c\big(\Gamma^{(1)}\big)$ where $\|\cdot\|_I$ is the $I$-norm introduced in Definition~\ref{Chap2-Def-Inorm}. Then the {\em full-norm} is defined by
$$
\|f\|_{full}\; 
	:= \;\sup
		\big\{ 
			\|\pi(\fz)\| \;\big|\; 
			\pi \text{ bounded representation of } \Cc_c\big(\Gamma^{(1)}\big)
		\big\}
		,\qquad \fz\in\Cc_c\big(\Gamma^{(1)}\big) \, .
$$
Again \cite[Proposition~II.1.11]{Renault80} guarantees that \gls{normfull} defines a $C^\ast$-norm. The {\em full $C^\ast$-algebra} \gls{CGfull} is defined by the closure of $\Cc_c\big(\Gamma^{(1)}\big)$ with respect to the full-norm $\|\cdot\|_{full}$. By the previous considerations the following well-known assertion is deduced.

\begin{proposition}[\cite{Renault80}]
\label{Chap2-Prop-FullsubsetRedCalgebras}
Let $\Gamma$ be a topological groupoid with left-continuous Haar system $\mu:=(\mu^x)_{x\in \Gamma^{(0)}}$. Then, for all $\fz\in\Cc_c\big(\Gamma^{(1)}\big)$, the estimate
$$
\|\fz\|_{red} \; 
	\leq \; \|\fz\|_{full}
$$
holds. Specifically, there exists a $\ast$-homomorphism $\Phi:\CG^\ast_{full}(\Gamma)\to\CG^\ast_{red}(\Gamma)$ that is surjective and norm-decreasing.
\end{proposition}

\begin{proof}
Proposition~\ref{Chap2-Prop-LeftRegulRepres} implies that the left-regular representation is a bounded representation. As the supremum in the full-norm runs over all bounded representations of $\Cc_c\big(\Gamma^{(1)}\big)$, the desired estimate follows. 

\vspace{.1cm}

Consider the identity map in $\Cc_c\big(\Gamma^{(1)}\big)$ that extends to a norm-decreasing, $\ast$-homomorphism $\Phi:\CG^\ast_{full}(\Gamma)\to\CG^\ast_{red}(\Gamma)$ since $\|\fz\|_{red} \leq \|\fz\|_{full}$ holds for all $\fz\in\Cc_c\big(\Gamma^{(1)}\big)$. According to \cite[Theorem~3.1.6]{Murphy90}, the image $\AG:=\Phi\big(\CG^\ast_{full}(\Gamma)\big)$ is a $C^\ast$-subalgebra of $\CG^\ast_{red}(\Gamma)$, i.e., it is a closed $C^\ast$-algebra with respect to the reduced-norm. On the other hand, $\Cc_c\big(\Gamma^{(1)}\big)$ is contained in $\AG$ which is dense in $\CG^\ast_{red}(\Gamma)$. Hence, $\AG$ is equal to $\CG^\ast_{red}(\Gamma)$ implying that $\Phi$ is surjective.
\end{proof}

\medskip

In the case of group $C^\ast$-algebras the full $C^\ast$-algebra and the reduced $C^\ast$-algebra coincide if and only if $G$ is amenable, c.f. \cite{Hul66}. It is natural to ask for an extension for topological groupoids. In order to do so the notion of topologically amenable groupoid is introduced. However there exist groupoids that are not topologically amenable but where the full and reduced $C^\ast$-algebra coincide, c.f. \cite{Wil15}. The reader is referred to \cite{Ana16} and references therein for recent results about the question when the reduced and the full $C^\ast$-algebras are isomorphic. That the reduced and full $C^\ast$-algebra coincide plays a crucial role in Chapter~\ref{Chap4-ToolContBehavSpectr}. Note that by uniqueness of $C^\ast$-norms the norms $\|\cdot\|_{red}$ and $\|\cdot\|_{full}$ coincide whenever the reduced and full $C^\ast$-algebras coincide.

\medskip

In 1977, {\sc Zimmer} extends the notion of amenability to measured groupoids, c.f. \cite{Zim77a,Zim77b,Zim78}. {\sc Anantharaman-Delaroche} \cite{Ana79} introduced amenability for a continuous action of a locally compact group on a von Neumann algebra. Thus, {\sc Zimmer} and {\sc Anantharaman-Delaroche} provide the foundation for the extension of amenability of groups to groupoids. On this basis, the notion of topologically amenable groupoids was first introduced by {\sc Renault} \cite{Renault80}. It turns out that a topologically amenable groupoid is measurewise amenable as well, c.f. \cite[Section 3.3]{AnantharamanRenault00}. Whenever the groupoid is second-countable, locally compact with countable orbits 
both notations coincide, c.f. \cite[Theorem 3.3.7]{AnantharamanRenault00}. For a transformation group groupoid, c.f. Definition~\ref{Chap2-Def-TransformationGroupGroupoid}, this holds whenever the group $G$ is discrete, c.f. \cite{Ana87}. Clearly, this is not a necessary condition. For instance, any group can act trivially by the identity on a space $X$. In this case, each orbit contains exactly one element.

\medskip

Let $(\Gamma,\mu)$ be a topological groupoid with left-continuous Haar system $\mu$. Then \cite[Proposition~2.2.6]{AnantharamanRenault00} or \cite[Proposition~3.3]{AnRe01} provide a characterization of the property that the groupoid $(\Gamma,\mu)$ is topologically amenable. This characterization is used here as the definition.

\begin{definition}[Topologically amenable]
\label{Chap2-Def-AmenableGroupoid}
A topological groupoid $\Gamma$ with left-continu\-ous Haar system $\mu:=(\mu^x)_{x\in \Gamma^{(0)}}$ is called {\em topologically amenable} if there exists a sequence of non-negative, continuous functions $\{\fz_n\}_{n\in\NM}\subseteq\Cc_c\big(\Gamma^{(1)}\big)$ such that the following assertions hold.
\begin{description}
\item[(A1)\label{(A1)}] For all $n\in\NM$ and $x\in \Gamma^{(0)}$, the function $f_n$ is normalized, i.e., $\int_{\Gamma^x} \fz_n\; d\mu^x=1$.
\item[(A2)\label{(A2)}] The sequence of functions $m_n:\Gamma^{(1)}\to\CM$ defined by
$$
m_n(\gamma) \; 
	:= \; \int_{\Gamma^{r(\gamma)}} 
		\big| \fz_n(\gamma^{-1}\cdot\rho)- \fz_n(\rho)\big|
	\; d\mu^{r(\gamma)}(\rho)
$$ 
converges uniformly on compact subsets of $\Gamma^{(1)}$ to zero.
\end{description}
Such a sequence $\{\fz_n\}_{n\in\NM}$ of functions is called an {\em approximate invariant mean}.
\end{definition}

Let $G$ be a second-countable, locally compact, Hausdorff group with Haar measure $\lambda$. A sequence $(F_n)_{n\in\NM}$ of non-empty, compact subsets of $G$ is called {\em weak-F\o lner sequence} \label{Page-weak-Folner} if $\lambda(F_n)>0,\; n\in\NM,$ and 
$$
\lim\limits_{n\to\infty} \frac{\lambda(K F_n\bigtriangleup F_n)}{\lambda(F_n)}=0
$$
for all non-empty, compact subsets $K\subseteq G$. Here, $F\bigtriangleup K:= (F\setminus K)\cup (K\setminus F)$ denotes the symmetric difference. Then a group $G$ is called {\em amenable} whenever it has a weak-F\o lner sequence. In this case, the function $f_n:=\frac{1}{\lambda(F_n)}\chi_{F_n}$ defines an approximate invariant mean whenever the sets $F_n,\; n\in\NM,$ are open, additionally. Note that the groups of our interest are discrete and so the sets $F_n,\; n\in\NM,$ are open.

\medskip

In the case of a topologically amenable groupoid, the associated reduced $C^\ast$-algebra and full $C^\ast$-algebra coincide. The following assertion is proven in \cite[Proposition~6.1.8]{AnantharamanRenault00}, see also \cite[Proposition~4.1]{AnRe01}. Note that topological amenability of a groupoid implies measurewise amenability, c.f. \cite{AnantharamanRenault00,AnRe01}. 

\begin{proposition}[\cite{AnantharamanRenault00,AnRe01}]
\label{Chap2-Prop-Red=FullCalgebra}
Let $\Gamma$ be a topological groupoid with left-con\-tinuous Haar system $(\mu^x)_{x\in\Gamma^{(0)}}$. If the groupoid $\Gamma$ is topologically amenable, then the equation $\CG^{\ast}_{red}(\Gamma)=\CG^{\ast}_{full}(\Gamma)$ holds. In particular, the identity $\|\fz\|_{red}=\|\fz\|_{full}$ holds for all $\fz\in\Cc_c\big(\Gamma^{(1)}\big)$. 
\end{proposition}

In general, the assumption that the groupoid $\Gamma$ is second countable is necessary in Proposition~\ref{Chap2-Prop-Red=FullCalgebra}. In the case of a transformation group groupoid the result holds also if the groupoid is not second countable, c.f. \cite[Remark 4.2]{AnantharamanRenault00}.

\begin{definition}
\label{Chap2-Def-GroupoidCalgebra}
Let $\Gamma$ be a topologically amenable groupoid with left-continuous Haar system $\mu:=(\mu^x)_{x\in \Gamma^{(0)}}$. Then the {\em associated $C^\ast$-algebra} $\CG^\ast(\Gamma)$ is defined by $\CG^\ast_{red}(\Gamma)$ and the $C^\ast$-norm is given by $\|\fz\|:=\|\fz\|_{red}=\|\fz\|_{full}$ for $\fz\in\Cc_c\big(\Gamma^{(1)}\big)$.
\end{definition}

The section is finished by providing a sufficient condition so that the reduced and the full $C^\ast$-algebra are unital. In order to do so the following well-known definition is needed.

\begin{definition}[\'Etale]
\label{Chap2-Def-Etale}
Let $\Gamma$ be a topological groupoid. Then $\Gamma$ is called {\em \'etale} if the range map is a local homeomorphism, i.e., for each $\gamma\in\Gamma^{(1)}$ there is an open neighborhood $U\subseteq\Gamma^{(1)}$ of $\gamma$ such that $r:U\to r(U)$ is a homeomorphism.
\end{definition}

\begin{remark}
\label{Chap2-Rem-Etale}
Recall that the inverse map in a topological groupoid is a homeomorphism, c.f. Proposition~\ref{Chap2-Prop-BasGroupoid}~(e). Thus, the range map is a local homeomorphism if and only if the source map is a local homeomorphism following from Proposition~\ref{Chap2-Prop-BasGroupoid}~(b).
\end{remark}

\medskip

For convenience of the reader, an example of an \'etale groupoid and a non-\'etale groupoid is considered.

\begin{example}
\label{Chap2-Ex-EtalGroup}
Let $G$ be a discrete group. Then choose for every $g\in G$ the open neighborhood $\{g\}$. Clearly, the range is a homeomorphism. Hence, $G$ is an \'etale groupoid.
\end{example}

\begin{example}
\label{Chap2-Ex-NonEtalGroup}
Let $G$ be a topological group that is not discrete, namely any open neighborhood $U$ of $e$ contains a $g\neq e$. Thus, for every open neighborhood $U$ of $e\in G$ and $g\in U$, the range $r(g)$ is equal to $e$. Consequently, the range is not injective as $r(e)=e$ and so $G$ is not \'etale.
\end{example}

Let $X$ and $Y$ be topological spaces. Then $f:X\to Y$ is called {\em locally injective} if, for each $x\in X$, there exists an open neighborhood $U\subseteq X$ of $x$ such that $f:U\to Y$ is injective. Note that a local homeomorphism is always locally injective whereas the converse does not hold.

\begin{lemma}[\cite{Renault09,Tho10}]
\label{Chap2-Lem-CharUnitSpaceOpen}
Let $\Gamma$ be a topological groupoid. Then the following assertions are equivalent.
\begin{itemize}
\item[(i)] The unit space $\Gamma^{(0)}$ is open.
\item[(ii)] The range map $r$ is locally injective.
\end{itemize}
\end{lemma}

\begin{proof}
(i)$\Rightarrow$(ii): Let $\gamma\in\Gamma^{(1)}$. As $\Gamma^{(0)}$ is open in $\Gamma^{(1)}$, there exists an open neighborhood $V\subseteq\Gamma^{(0)}$ of $r(\gamma)$. The map $(\rho,\zeta)\in\Gamma^{(2)}\mapsto\rho^{-1}\zeta$ is continuous and for $\rho=\zeta=\gamma^{-1}$ its value is equal to $r(\gamma)$. Thus, there exists a $V_1\times V_2\subseteq\Gamma^{(1)}\times\Gamma^{(1)}$ open such that $\gamma\in V_1$, $\gamma\in V_2$ and the image under this map is contained in $V$ whenever the composition is defined. Recall that in the product topology the canonical projections to the components are open maps. Hence, $V_1$ and $V_2$ are open subsets of $\Gamma^{(1)}$. Consequently, $U:=V_1\cap V_2$ defines an open neighborhood of $\gamma$. Let $\rho,\zeta\in U$ be such that $r(\rho)=r(\zeta)$. Then $(\rho^{-1},\zeta)\in\Gamma^{(2)}$ is defined by Proposition~\ref{Chap2-Prop-BasGroupoid}~(a) and $\rho^{-1}\zeta$ is an element of $V\subseteq\Gamma^{(0)}$. Thus, the equalities $\zeta=\rho\rho^{-1}\zeta = \rho$ hold as $\rho^{-1}\zeta$ is a unit implying that $r:U\to\Gamma^{(0)}$ is injective.

\vspace{.1cm}

(ii)$\Rightarrow$(i): Assume the unit space $\Gamma^{(0)}\subseteq\Gamma^{(1)}$ is not open. Hence, there exists an $x\in\Gamma^{(0)}$ such that, for all open neighborhoods $V\subseteq\Gamma^{(1)}$ of $x$, the intersection $\big(\Gamma^{(1)}\setminus\Gamma^{(0)}\big)\cap V$ is non-empty. Thus, there is a sequence $\gamma_n\in\Gamma^{(1)}\setminus\Gamma^{(0)},\; n\in\NM,$ converging to $x$. More precisely, the limit $\lim_{n\to\infty}r(\gamma_n)$ is equal to $x$ by the continuity of the range, c.f. Proposition~\ref{Chap2-Prop-BasGroupoid}~(f). Let $U\subseteq\Gamma^{(1)}$ be the associated open neighborhood of $x$ such that $r:U\to\Gamma^{(0)}$ is injective. By the previous considerations, there exists an $n_0\in\NM$ such that $r(\gamma_{n_0}),\gamma_{n_0}\in U$. This contradicts the injectivity of $r:U\to\Gamma^{(0)}$ since, by construction, the relations $\gamma_{n_0}\neq r(\gamma_{n_0})$ and $r(\gamma_{n_0})=r(r(\gamma_{n_0}))$ hold.
\end{proof}

\begin{remark}
\label{Chap2-Rem-CharUnitSpaceOpen}
By Proposition~\ref{Chap2-Prop-BasGroupoid}~(e) it is worth noticing that $r$ is locally injective if and only if $s$ is so. For each \'etale groupoid  $\Gamma$, the range map is locally injective. Thus, the requirement that $\Gamma^{(0)}$ is an open subset of $\Gamma^{(1)}$ is weaker than the \'etale property.
\end{remark}

Whenever $\Gamma$ admits, additionally, a left-continuous Haar system it follows that $r$ is locally injective if and only if $r$ is a local homeomorphism. This is derived by \cite[Propo\-sition~I.2.8]{Renault80} and the following observation.

\begin{lemma}[\cite{Renault80,Tho10}]
\label{Chap2-Lem-FibreDiscrete}
Let $\Gamma$ be a topological groupoid such that $\Gamma^{(0)}$ is open. Then, for each $x\in\Gamma^{(0)}$, the $r$-fiber $\Gamma^{x}$ and the $s$-fiber $\Gamma_x$ are discrete with respect to the topology induced by $\Gamma^{(1)}$.
\end{lemma}

\begin{proof}
Let $x\in\Gamma^{(0)}$. According to Lemma~\ref{Chap2-Lem-CharUnitSpaceOpen} there exists an open set $U\subseteq\Gamma^{(1)}$ such that $r:U\to\Gamma^{(0)}$ is injective. Thus, the intersection $r^{-1}(\{x\})\cap U$ only contains the element $x$. In detail, $x$ is isolated in the fiber $\Gamma^x:=r^{-1}(\{x\})$. This argument translates to every element $\gamma\in\Gamma^x$ as follows: The map $\Gamma^x\to\Gamma^{s(\gamma)}$ defined by $\rho\mapsto\gamma^{-1}\rho$ is well-defined and a homeomorphism. Furthermore, the image of $\gamma$ is equal to $s(\gamma)$. Since $s(\gamma)$ is isolated in the topology on $\Gamma^{s(\gamma)}$, the arrow $\gamma$ is also isolated in $\Gamma^x$.
\end{proof}

\medskip

Motivated by Lemma~\ref{Chap2-Lem-FibreDiscrete} a groupoid $\Gamma$ is called {\em $r$-discrete} if the unit space $\Gamma^{(0)}$ is open, c.f. \cite[Definition I.2.6]{Renault80}.

\medskip

Note that many references assume that the groupoid is \'etale with compact unit space implying that the associated $C^\ast$-algebra is unital, c.f. the following assertion. However, by the best knowledge of the author there is no reference proving that these conditions are necessary. The necessity of these conditions is hinted without a proof in \cite{Exe14}. In the case of a transformation group groupoid, {\sc Landstad} \cite{Lan16} proved the equivalence. A similar statement for the general case of a topological groupoid is still an open question, c.f. discussion in Section~\ref{Chap8-Sect-Groupoids}.

\begin{theorem}
\label{Chap2-Theo-UnitGroupoidCalgebra}
Let $\Gamma$ be a topological groupoid with left-continuous Haar system $\mu:=(\mu^x)_{x\in\Gamma^{(0)}}$. If the groupoid $\Gamma$ is \'etale and has a compact unit space $\Gamma^{(0)}$, then the full $C^\ast$-algebra $\CG^\ast_{full}(\Gamma)$ and the reduced $C^\ast$-algebra $\CG^\ast_{red}(\Gamma)$ are unital. Furthermore, the unit $\mathpzc{1}$ is given by the characteristic function $\chi_{\Gamma^{(0)}}\in\Cc_c\big(\Gamma^{(1)}\big)$ of the unit space if, additionally, the Haar system $\mu$ is normalized, i.e., $\mu^x(\{x\})=1$ for all $x\in\Gamma^{(0)}$. 
\end{theorem}

\begin{proof}
Let $\Gamma$ be an \'etale groupoid with compact unit space. 

\vspace{.1cm}

The unit space $\Gamma^{(0)}$ is a compact subset of $\Gamma^{(1)}$ by assumption. Then the characteristic function $\chi_{\Gamma^{(0)}}:\Gamma^{(1)}\to\CM$ of the unit space defined by
$$
\chi_{\Gamma^{(0)}}(\gamma) \;
	:= \;
		\begin{cases}
			1,\quad &\gamma\in\Gamma^{(0)} \, ,\\
			0,\quad &\gamma\not\in\Gamma^{(0)} \, ,
		\end{cases}
$$
has compact support. Lemma~\ref{Chap2-Lem-CharUnitSpaceOpen} and Remark~\ref{Chap2-Rem-CharUnitSpaceOpen} imply that $\Gamma^{(0)}$ is also open as the groupoid $\Gamma$ is \'etale. Since the characteristic function on a clopen set is continuous, the map $\chi_{\Gamma^{(0)}}$ is an element of $\Cc_c\big(\Gamma^{(1)}\big)$. According to Lemma~\ref{Chap2-Lem-FibreDiscrete}, the $r$-fiber $\Gamma^x$ is discrete for all $x\in\Gamma^{(0)}$. Since $\mu^x$ has support $\Gamma^x$ for $x\in\Gamma^{(0)}$, it follows that the support of $\mu^x$ is discrete. 

\vspace{.1cm}

The map $\Gamma^{(0)}\ni x\mapsto\mu^x\big(\chi_{\Gamma^{(0)}}\big)\in[0,\infty)$ is continuous by \nameref{(H2)}. Since $\Gamma^{(1)}$ is second-countable, locally compact and Hausdorff, it is a normal space. Consequently, Tietze's Extension Theorem (\cite{Tie15}, \cite[Anhang~III.28.]{Ury25}) is applicable to the continuous map $\Gamma^{(0)}\ni x\mapsto\mu^x\big(\chi_{\Gamma^{(0)}}\big)\in[0,\infty)$ with compact support $\Gamma^{(0)}$. Hence, there exists a compactly supported, continuous function $\fz:\Gamma^{(1)}\to\RM$ such that $\fz|_{\Gamma^{(0)}}(x)=\mu^x\big(\chi_{\Gamma^{(0)}}\big)$, c.f. Theorem~\ref{App1-Theo-TietzesThm}. Define $\mathpzc{1}:\Gamma^{(1)}\to\RM$ by 
$$
\mathpzc{1}(\gamma) \; 
	= \; 	
		\begin{cases}
			\frac{1}{f(\gamma)}, \quad &\gamma\in\Gamma^{(0)} \, ,\\
			0,\quad &\gamma\not\in\Gamma^{(0)} \, .
		\end{cases}
$$
Note that the image $\fz(\gamma)$ is not zero for $\gamma\in\Gamma^{(0)}$. Thus, $\mathpzc{1}$ is well-defined and compactly supported on $\Gamma^{(0)}$. The function is continuous as $\Gamma^{(0)}$ is clopen and $f$ is continuous. Altogether, $\mathpzc{1}$ is an element of $\Cc_c\big(\Gamma^{(1)}\big)$. Then, for $\gz\in\Cc_c\big(\Gamma^{(1)}\big)$, the equations
$$
(\mathpzc{1}\star \gz)(\gamma)\; 
	= \; \int_{\Gamma^{r(\gamma)}} 
			\mathpzc{1}(\rho)\cdot \gz \left(\rho^{-1} \gamma \right) \; 
		d\mu^{r(\gamma)}(\rho) \;
	= \; \frac{\mu^{r(\gamma)}(\{r(\gamma)\})}{\fz(r(\gamma))} \cdot \gz\left( r(\gamma)^{-1} \gamma \right)\;
	= \; \gz(\gamma)
$$
are derived for all $\gamma\in\Gamma^{(1)}$. Similarly, the identity $\gz\star \mathpzc{1} = \gz$ follows. Thus, $\mathpzc{1}\in\Cc_c\big(\Gamma^{(1)}\big)$ is a unit of the full and reduced $C^\ast$-algebra since $\Cc_c\big(\Gamma^{(1)}\big)$ is a dense subset of $\CG^\ast_{full}(\Gamma)$ and $\CG^\ast_{red}(\Gamma)$. Clearly, the equations $\fz(x)=1\,,\; x\in\Gamma^{(0)}\,,$ hold if the Haar system $\mu$ is normalized. Hence, $\mathpzc{1}=\chi_{\Gamma^{(0)}}\in\Cc_c\big(\Gamma^{(1)}\big)$ is derived. 
\end{proof}

\medskip

Denote by $\Cc_0(\Gamma^{(1)})$ the set of continuous functions $\fz:\Gamma^{(1)}\to\CM$ that vanish at infinity, c.f. Example~\ref{Chap3-Ex-ContFunctCalgebra}. This set is equipped with the uniform norm $\|\fz\|_\infty:=\sup_{\gamma\in\Gamma^{(1)}}|\fz(\gamma)|$. The \'etale property of a groupoid implies that elements of the reduced $C^\ast$-algebra are represented as elements of $\Cc_0(\Gamma^{(1)})$. This result is due to {\sc Renault} \cite[Proposition~II.4.1, Proposition~II.4.2]{Renault80}.

\begin{proposition}[\cite{Renault80}]
\label{Chap2-Prop-C0functionRedCAlgebra}
Let $\Gamma$ be a topological, \'etale groupoid with left-continuous Haar system $\mu:=(\mu^x)_{x\in\Gamma^{(0)}}$. Then $\CG^\ast_{red}(\Gamma)$ is continuously embedded in the space $\Cc_0(\Gamma^{(1)})$, i.e., there exists an injective map $\jmath:\CG^\ast_{red}(\Gamma)\to\Cc_0(\Gamma^{(1)})$ such that $\|\fz\|_\infty\leq \|\fz\|_{red}$ holds for every $\fz\in\CG^\ast_{red}(\Gamma)$.
\end{proposition}

\section{Groupoids induced by a dynamical system}
\label{Chap2-Sect-TransformationGroupGroupoid}

This section is devoted to study transformation group groupoids defined by a topological dynamical system. It is proven that such a groupoid has always a left-continuous Haar system. Whenever the group is amenable it is shown that the induced transformation group groupoid is amenable as well. 

\medskip

Let $G$ and $H$ be groups and $\alpha$ an action of the group $G$ on $H$. Then a new group is defined by the semi-direct product $H\rtimes_\alpha G$. Let $(X,G,\alpha)$ be a (topological) dynamical system. Then $G$ and $X$ are also special groupoids, c.f. Example~\ref{Chap2-Ex-Group} and Example~\ref{Chap2-Ex-Set}. Meanwhile, $G$ acts on the space $X$. Like in the case of a group action on a group, a new (topological) groupoid $X\rtimes_\alpha G$ is defined. Note that this construction can be generalized to the action of a general groupoid on another groupoid where the group action on a set is a special case.

\begin{definition}[Transformation group groupoid, \cite{Ehr57}]
\label{Chap2-Def-TransformationGroupGroupoid}
Let $(X,G,\alpha)$ be a topological dynamical system. The {\em transformation group groupoid $\Gamma:=\Gamma(X):=\text{\gls{XalphaG}}$} denotes the groupoid defined by 
\begin{description}
\item[(TG1)\label{(TG1)}] the arrow space $\Gamma^{(1)}:=X\times G$ that is endowed with the product topology;
\item[(TG2)\label{(TG2)}] the set of composable arrows
$$
\Gamma^{(2)}\; :=
	\; \left.\left\{
		\big((x|g),(y|h)\big)\in\Gamma^{(1)}\times\Gamma^{(1)} \;\right|\; y=\alpha_{g^{-1}}(x) 
	\right\}
$$
and the composition map $\circ:\Gamma^{(2)}\to\Gamma^{(1)},\; \big((x|g),(\alpha_{g^{-1}}(x)|h)\big)\mapsto (x|gh)$;  
\item[(TG3)\label{(TG3)}] the inverse map $^{-1}:\Gamma^{(1)}\to\Gamma^{(1)},\; (x|g)^{-1}:=(\alpha_{g^{-1}}(x)|g^{-1})$.
\end{description}
\end{definition}

\begin{figure}[htb]
\centering
\includegraphics[scale=1.1]{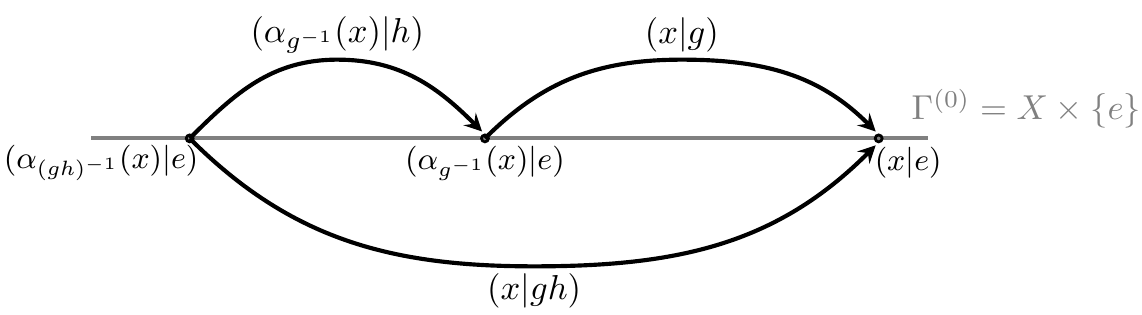}
\caption{The composition in the transformation group groupoid $\Gamma:=X\rtimes_\alpha G$.}
\label{Chap2-Fig-TransfGroupGroupoid}
\end{figure}

According to Definition~\ref{Chap2-Def-Groupoid}, a groupoid is interpreted as a graph with vertices $\Gamma^{(0)}$ and arrow set $\Gamma^{(1)}$. A $\gamma\in\Gamma^{(1)}$ is then an arrow joining the source $s(\gamma)$ with the range $r(\gamma)$. Then the composition of elements of $\Gamma^{(1)}$ is the composition of arrows, c.f. Figure~\ref{Chap2-Fig-TransfGroupGroupoid}.

\begin{remark}
\label{Chap2-Rem-OrbitGroupoidDynSystem}
In the case of a transformation group groupoid $\Gamma:=X\rtimes_\alpha G$ the $\Gamma$-orbit $\Orb_\Gamma((x|e))$, defined in Definition~\ref{Chap2-Def-r/s-Fibre-Orbit}, agrees with the orbit $\Orb(x):=\{\alpha_g(x)\;|\; g\in G\}$ of $x\in X$ in the dynamical system $(X,G,\alpha)$.
\end{remark}

\medskip

This object was first introduced by {\sc Ehresmann} \cite{Ehr57} in the year 1957. The transformation group groupoid \cite{StarlingThesis12,Sta14} is also called {\em transformation groupoid} \cite{GoehleThesis09}, {\em action groupoid} and {\em semi-direct product} \cite{Renault80,AnantharamanRenault00,AnRe01}, elsewhere. The following result is well-known, c.f. \cite{Renault80,Renault09}. 

\begin{proposition}[\cite{Renault80}]
\label{Chap2-Prop-TransGrGroupProp}
The transformation group groupoid $\Gamma:=X\rtimes_\alpha G$ defines a topological groupoid. The range and source map  $r,s:\Gamma^{(1)}\to\Gamma^{(0)}$ are given by $r(x|g)=(x|e)$ and $s(x|g)=\left(\alpha_{g^{-1}}(x)|e\right)$ for $(x|g)\in\Gamma^{(1)}$. Thus, the unit space $\Gamma^{(0)}$ is identified with the compact space $X$, i.e., these spaces are homeomorphic.
\end{proposition}

\begin{proof}
Since the topological spaces $G$ and $X$ are second-countable, locally compact and Hausdorff, the space $\Gamma^{(1)}:=X\times G$ is so as well. 

\vspace{.1cm}

\nameref{(G1)}: Let $(x|g_x),(y|g_y),(z,g_z)\in\Gamma^{(1)}$ be such that $\big((x|g_x),(y|g_y)\big),\big((y|g_y),(z|g_z)\big)\in\Gamma^{(2)}$. By definition of $\Gamma^{(2)}$, it follows that $\alpha_{g_x^{-1}}(x)=y$ and $\alpha_{g_y^{-1}}(y)=z$. Thus, the equalities 
$$
z\; 
	=\; \alpha_{g_y^{-1}}(\alpha_{g_x^{-1}}(x))\;
		=\;\alpha_{g_y^{-1}g_x^{-1}}(x)\;
			=\; \alpha_{(g_xg_y)^{-1}}(x)
$$
are derived by Definition~\ref{Chap2-Def-DynSyst}~(ii). Hence, the pair $\big((x|g_x)\circ(y|g_y),(z,g_z)\big)$ is an element of $\Gamma^{(2)}$ since $(x|g_x)\circ(y|g_y)=(x|g_xg_y)$. Consequently, Condition~\nameref{(G1)} is satisfied.

\vspace{.1cm}

\nameref{(G2)}: Assertion \nameref{(G2)} immediately follows by the definition of $\Gamma^{(2)}$ and the inverse map.

\vspace{.1cm}

\nameref{(G3)}: Note that $(x|g)^{-1}\circ(x|g)=(\alpha_{g^{-1}}(x)|e)$ and $(x|g)\circ(x|g)^{-1}=(x|e)$ hold for each $(x|g)\in\Gamma^{(1)}$. Let $\big((x|g),(y|h)\big)\in\Gamma^{(2)}$. Since $\alpha_{g^{-1}}(x)=y$, the equations
$$
(x|g)^{-1}\circ(x|g)\circ(y|h)\; 
	=\; (\alpha_{g^{-1}}(x)|e)\circ(y|h)\;
		=\; (\alpha_{g^{-1}}(x)|eh)\;
			=\; (y|h)
$$
follow. Similarly, $(x|g)\circ(y|h)\circ(y|h)^{-1}=(x|g)$ is deduced. Thus, $\Gamma$ fulfills Condition~\nameref{(G3)}.

\vspace{.1cm}

Finally, the previous considerations immediately yield the desired equalities $s(x|g)=(\alpha_{g^{-1}}(x)|e)$ and $r(x|g)=(x|e)$. Hence, the unit space $\Gamma^{(0)}=X\times\{e\}$ is identified with the topological space $X$. 
\end{proof}

\medskip

The following proposition provides a characterization of the \'etale property for a transformation group groupoid.

\begin{proposition}[\cite{Renault80}]
\label{Chap2-Prop-Etale}
Let $\Gamma:=X\rtimes_{\alpha} G$ be a transformation group groupoid. 
\begin{itemize}
\item[(a)] For each $x\in\Gamma^{(0)}$, the $r$-fiber $\Gamma^x$ is isomorphic to the topological group $G$. Specifically, $\Gamma^x$ is equal to $\{x\}\times G$.
\item[(b)] The groupoid $\Gamma$ is \'etale if and only if the group $G$ is discrete.
\end{itemize}
\end{proposition}

\begin{proof}
(a): Let $x\in\Gamma^{(0)}$. According to Proposition~\ref{Chap2-Prop-TransGrGroupProp}, the $r$-fiber $\Gamma^x$ is equal to $\{x\}\times G$. Thus, it is isomorphic to the group $G$.

\vspace{.1cm}

(b): Let $\Gamma$ be \'etale. According to Lemma~\ref{Chap2-Lem-CharUnitSpaceOpen} and Lemma~\ref{Chap2-Lem-FibreDiscrete} the $r$-fiber $\Gamma^x$ is discrete. Thus, statement (a) implies that $G$ is discrete. Conversely, let $G$ be discrete and $\gamma:=(x|g)\in\Gamma^{(1)}$ be arbitrary. Then $X\times\{g\}\subseteq\Gamma^{(1)}$ is an open neighborhood of $\gamma$. Clearly, the range $r$ is a homeomorphism if restricted to $X\times\{g\}$. Hence, $\Gamma$ is \'etale. 
\end{proof}

\medskip

The following result can already be found in \cite{Renault80,AnantharamanRenault00}. For the convenience of the reader, the proof is presented here.

\begin{proposition}[\cite{Renault80,AnantharamanRenault00}]
\label{Chap2-Prop-TransGroupAmen}
Let $(X,G,\alpha)$ be a topological dynamical system and $\lambda$ denotes the left-invariant Haar measure of the group $G$. The associated transformation group groupoid $\Gamma:=X\rtimes_\alpha G$ has a left-continuous Haar system $\mu:=(\mu^x)_{x\in X}$ defined by $\mu^x:=\delta_x\times\lambda$. If, additionally, $G$ is a discrete, amenable group then $\Gamma$ is topologically amenable.
\end{proposition}

\begin{proof}
Since $G$ is a locally compact group, there exists a left-invariant Haar measure $\lambda$, c.f. \cite{Weil40,Loomis53}. For $x\in X$, denote by $\delta_x$ the Dirac measure. Consider the measure $\mu^x:=\delta_x\times\lambda,\; (x,e)\in\Gamma^{(0)},$ defined on the Borel-$\sigma$-algebra $\Gamma^{(1)}$. In the following it is shown that (i) $\mu:=(\mu^x)_{x\in\Gamma^{(0)}}$ is a left-continuous Haar system on $\Gamma:=X\rtimes_\alpha G$ and that (ii) $(\Gamma,\mu)$ is topologically amenable if $G$ is amenable.

\vspace{.1cm}

(i): As $\lambda$ is the Haar measure on $G$ the support of $\lambda$ is $G$. Thus, $\supp(\mu^x)$ is equal to $\{x\}\times G=\Gamma^x$, c.f. Proposition~\ref{Chap2-Prop-Etale}. Thus, Condition~\nameref{(H1)} follows.

\vspace{.1cm}

Let $x\in X$ and $\fz\in\Cc_c\big(\Gamma^{(1)}\big)$. Since $\fz$ has compact support,
$$
\mu^x(\fz)\;
	=\; \int_G \fz(x|g)\; d\lambda(g)
$$
is well-defined and finite. Then the continuity of $\fz$ implies \nameref{(H2)}.

\vspace{.1cm}

The left-invariance of $\lambda$ on $G$ yields \nameref{(H3)}. More precisely, for $\fz\in\Cc_c\big(\Gamma^{(1)}\big)$ and $(x|g)\in\Gamma^{(1)}$, the equations
$$
\int_G \fz\big((x|g)\circ(\alpha_{g^{-1}}(x)|h)\big) \; d\lambda(h)\;
		=\; \int_G \fz(x|gh) \; d\lambda(h)\;
			=\; \int_G \fz(x|h)\; d\lambda(h)
$$
are derived.

\vspace{.1cm}

(ii): Recall the notion of a weak-F\o lner sequence $(F_n)_{n\in\NM}$. Since $G$ is amenable, a weak-F\o lner sequence exists. The group $G$ is discrete and the sets $F_n,\; n\in\NM,$ are compact. Hence, the sets $F_n\,,\; n\in\NM\,,$ are finite. Thus, the sets $F_n,\; n\in\NM,$ are open and compact. For $n\in\NM$, define $\fz_n:\Gamma^{(1)}\to\CM$ by
$$
\fz_n(x|g) \; 
	:= \; \frac{1}{\lambda(F_n)}\chi_{F_n}(g).
$$
By construction, it has compact support $X\times F_n$. Since the set $F_n$ is clopen for $n\in\NM$, the function $\fz_n$ is continuous, i.e., $\fz_n\in\Cc_c\big(\Gamma^{(1)}\big)$. Clearly, \nameref{(A1)} is satisfied. Condition \nameref{(A2)} follows by the F\o lner condition. In detail, for $\gamma:=(x|g)\in\Gamma^{(1)}$, the equations
\begin{align*}
m_n(\gamma) \;
	= \; &\int_{\Gamma^{r(\gamma)}} 
			\big| \fz_n(\gamma^{-1}\cdot\gamma')- \fz_n(\gamma')\big|
		\; d\mu^{r(\gamma)}(\gamma') \;\\ 
	= \; &\int_G 
			\big| \fz_n(\alpha_{g^{-1}}(x)|g^{-1} h)- \fz_n(x|h)\big| 
		\; d\lambda(h)\\
	= \; & \frac{1}{\lambda(F_n)} 
		\int_G 
			\big| \chi_{F_n}(g^{-1} h)- \chi_{F_n}(h)\big|
		\; d\lambda(h)\\
	= \; & \frac{1}{\lambda(F_n)} \lambda(gF_n\bigtriangleup F_n)
\end{align*}
hold. Thus, $\big(m_n(\gamma)\big)_{n\in\NM}$ converges uniformly to zero on compact subsets of $\Gamma^{(1)}$ since $(F_n)_{n\in\NM}$ is a weak-F\o lner sequence.
\end{proof}

\begin{remark}
\label{Chap2-Rem-ExAmenNotNec}
Note that the requirement that the group $G$ is amenable is not necessary to get a topologically amenable transformation group groupoid. More precisely, there exist non-amenable groups $G$ such that the transformation group groupoid $\Gamma:=X\rtimes_\alpha G$ is still topologically amenable. For instance, the free group $\mathbb{F}_2$ with $2$ generators acting on its boundary $\partial\mathbb{F}_2$, c.f. \cite[Example~3.8]{AnRe01}.
\end{remark}

As shown in Proposition~\ref{Chap2-Prop-TransGrGroupProp} the unit space of a transformation groupoid $\Gamma:=X\rtimes_\alpha G$ is exactly the compact space $X$. Thus, the reduced and full $C^\ast$-algebra of a transformation groupoid has a unit if the groupoid is \'etale ($G$ discrete) by Theorem~\ref{Chap2-Theo-UnitGroupoidCalgebra}. The converse implication also holds which is proven by {\sc Landstad} \cite{Lan16}.

\medskip

For the transformation group groupoid $X\rtimes_\alpha G$ with $G$ discrete, the left-regular representation acts on the space $\ell^2(G)$ as follows.

\begin{proposition}
\label{Chap2-Prop-LeftRegReprTransGrouGrou}
Let $(X,G,\alpha)$ be a topological dynamical system where $G$ is a discrete, countable group with the counting measure $\lambda$ as Haar measure. Consider the topological groupoid $\Gamma:=X\rtimes_\alpha G$ with left-continuous Haar system $\mu:=(\mu^x)_{x\in\Gamma^{(0)}},\; \mu^x:=\delta_x\times\lambda$. Then the space $L^2(\Gamma^x,\mu^x)$ is unitary equivalent to $\ell^2(G)$ and so the left-regular representation $\pi^x(\fz)$ for $\fz\in\CG^\ast_{red}(\Gamma)$ acts on $\ell^2(G)$. More precisely, for $\fz\in\CG^\ast_{red}(\Gamma)$, the equation
$$
\left(\pi^x(\fz)\psi\right)(g)\; 
	= \; \sum\limits_{h\in G} \, 
		\fz \left(
			\alpha_{g^{-1}}(x)|g^{-1} h 
		\right) \cdot \psi(h)
$$
holds for $g\in G$ and $\psi\in\ell^2(G)$. Furthermore, the identities 
$$
	\begin{array}{rl}
		(\fz\star \gz)(x|g)\; 
			& := \; \sum\limits_{h\in G} \fz(x|h)\; \cdot \gz\big(\alpha_{h^{-1}}(x)|h^{-1}g\big)\, ,\\[0.5cm]
		\fz^*(x|g)\; 
			& := \; \overline{\fz(\alpha_{g^{-1}}(x)|g^{-1})} \, ,
	\end{array}
	\qquad \fz,\gz\in\Cc_c\big(\Gamma^{(1)}\big),\; (x|g)\in\Gamma^{(1)} \, ,
$$
hold.
\end{proposition}

\begin{proof}
According to Proposition~\ref{Chap2-Prop-Etale}, the $r$-fiber $\Gamma^x$ is equal to $\{x\}\times G$ for $x\in X$. Consequently, $L^2(\Gamma^x,\mu^x)$ and $\ell^2(G)$ are unitary equivalent as $\mu^x:=\delta_x\times\lambda$ is the counting measure on $G$. The desired identities immediately follow by the definition of the Haar system $(\mu^x)_{x\in \Gamma^{(0)}}$, the left-regular representation, the multiplication and the involution, c.f. Definition~\ref{Chap2-Def-LeftRegulRepres} and Proposition~\ref{Chap2-Prop-InvolutiveGroupoidAlgebra}.
\end{proof}

\begin{remark}
\label{Chap2-Rem-SummaryPropGroup}
Altogether, the transformation group groupoid $\Gamma$ associated with a topological dynamical system $(X,G,\alpha)$ has compact unit space identified with the space $X$, c.f. Proposition~\ref{Chap2-Prop-TransGrGroupProp}. Furthermore, the $r$-fiber $\Gamma^x$ is isomorphic to the group $G$ for each $x\in X$ implying that the left-regular representation of $\Cc_c\big(\Gamma^{(1)}\big)$ represents in the algebra of bounded linear operators on the Hilbert space $\ell^2(G)$, c.f. Proposition~\ref{Chap2-Prop-Etale} and Proposition~\ref{Chap2-Prop-LeftRegReprTransGrouGrou}. For our further purposes, it is necessary that the $C^\ast$-algebras have a unit, c.f. Remark~\ref{Chap3-Rem-NecessIdentityContSpect} below. Having Theorem~\ref{Chap2-Theo-UnitGroupoidCalgebra} and Proposition~\ref{Chap2-Prop-Etale} in mind, only dynamical systems with a discrete group $G$ are considered. If a topological space is discrete and second countable, the space is automatically countable. Thus, the discrete group $G$ is countable since the group $G$ is second-countable by definition of a dynamical system. For the sake of normalization, the counting measure $\lambda=\sum_{g\in G} \delta_g$ is chosen on the group $G$ as the left-invariant Haar measure. Thus, the unit of the reduced $C^\ast$-algebra $\CG^\ast_{red}\big(X\rtimes_\alpha G\big)$ is given by the characteristic function $\chi_{\Gamma^{(0)}}\in\Cc_c\big(\Gamma^{(1)}\big)$ of the unit space $\Gamma^{(0)}=X\times\{e\}$, c.f. Theorem~\ref{Chap2-Theo-UnitGroupoidCalgebra}.
\end{remark}

\section{Generalized Schr\"odinger operators and their spectrum}
\label{Chap2-Sect-SchrOp}

This section studies the spectrum of normal elements of a groupoid $C^\ast$-algebra. Specifically, topologically amenable transformation group group\-oid are considered for a discrete, countable group $G$ equipped with the counting measure. These elements provide an equivariant, strongly continuous family of operators on $\ell^2(G)$ by the left-regular representation, c.f. Proposition~\ref{Chap2-Prop-CovFamOp-Spect}. Theorem~\ref{Chap2-Theo-MinCharConstSpectr} provides a characterization of minimality of the dynamical system $(X,G,\alpha)$ in terms of the constancy of the spectra of all normal (self-adjoint) elements.

\medskip

Let $\AG$ be a unital $C^\ast$-algebra with a faithful family of representations $(\pi^x,\hs_x)\,,\; x\in X$. As discussed in Remark~\ref{Chap2-Rem-FaithfulRepresent}, this implies that the spectrum $\sigma(\fz)$ of a normal element $\fz\in\AG$ coincides with $\overline{\bigcup_{x\in X}\sigma(\pi^x(\fz))}$. It was already proven by {\sc Renault} \cite[Proposition~II.1.11]{Renault80} that the left-regular representation of a second countable, Hausdorff groupoid is faithful, c.f. Proposition~\ref{Chap2-Prop-LeftRegulRepres}. If the groupoid is, additionally, \'etale and topologically amenable with compact unit space, {\sc Exel} proved moreover that the left-regular representation is also invertibility sufficient, c.f. \cite[Theorem 2.10]{Exe14}. Together with \cite[Theorem~3.6]{NiPr15} this leads to the following.

\begin{theorem}[\cite{Exe14,NiPr15}]
\label{Chap2-Theo-SpectrGroupoidCAlg}
Let $\Gamma$ be a topologically amenable, \'etale groupoid such that the unit space is compact. Then every normal element $\fz\in\CG^\ast(\Gamma)$ satisfies 
$$
\sigma(\fz) \;
	= \; \bigcup_{x\in\Gamma^{(0)}}\sigma\big(\pi^x(\fz)\big)\,.
$$
\end{theorem}

\begin{proof}
Let $\fz\in\CG^\ast(\Gamma)$ be normal. Thanks to \cite[Theorem~2.10]{Exe14}, $\fz$ is invertible in $\CG^\ast(\Gamma)$ if and only if $\pi^x(\fz)$ is invertible for all $x\in X$. Then \cite[Theorem~3.6]{NiPr15} implies the desired identity $\sigma(\fz)=\bigcup_{x\in \Gamma^{(0)}}\sigma\big(\pi^x(\fz)\big)$.
\end{proof}

\begin{remark}
\label{Chap2-Rem-SpectrGroupoidCAlg}
Note that an invertibility sufficient family of representation is only called sufficient in \cite{Exe14}. It is worth noticing that a topological groupoid in terms of Definition~\ref{Chap2-Def-Groupoid} is assumed to be Hausdorff which is necessary for the result of {\sc Exel} \cite{Exe14}.
\end{remark}

In the following the notion of generalized Schr\"odinger operators associated with a dynami\-cal system is introduced. The name is motivated by the fact that families of Schr\"odinger and Jacobi operators on $\ell^2(G)$ are representations of elements of the associated groupoid $C^\ast$-algebra defined by symbolic dynamical systems $X=\as^G$, c.f. Section~\ref{Chap2-Sect-PattEqSchrOp} for a detailed elaboration.

\begin{definition}[Generalized Schr\"odinger operator]
\label{Chap2-Def-SchrodingerOperator}
Let $(X,G,\alpha)$ be a topological dynamical system so that the topological groupoid $\Gamma:=X\rtimes_\alpha G$ is \'etale and topologically amenable. The family of operators $H_X:=(H_x)_{x\in X}$ defined by $H_x:=\pi^x(\fz)\in\Ll\big(\ell^2(G)\big)$ for a self-adjoint element $\fz\in\CG^\ast(\Gamma)$ is called a {\em generalized Schr\"odinger operator} where $\pi^x$ denotes the left-regular representation. 

\vspace{.1cm}

If $\fz\in\Cc_c\big(\Gamma^{(1)}\big)$, then the generalized Schr\"odinger operator $H_X$ is called of {\em finite range}. Furthermore, a generalized Schr\"odinger operator $H_X$ is called {\em periodic} if $X$ is a strongly periodic dynamical system. If $X$ is completely aperiodic, $H_X$ is called {\em non-periodic}.
\end{definition}

Note that generalized Schr\"odinger operators are assumed to be self-adjoint.

\begin{remark}
\label{Chap2-Rem-SchrodingerOperator}
The notion of finite range for an element $\fz\in\Cc_c\big(\Gamma^{(1)}\big)$ is chosen since the support of $\fz$ in the group $G$ is finite (compact subsets are finite as $G$ is discrete).
\end{remark}

\begin{proposition}
\label{Chap2-Prop-CovFamOp-Spect}
Let $(X,G,\alpha)$ be a topological dynamical system so that the topological groupoid $\Gamma:=X\rtimes_\alpha G$ is \'etale and topologically amenable. The following assertions hold for every normal element $\fz\in\CG^\ast(\Gamma)$.
\begin{itemize}
\item[(a)] The spectrum $\sigma(\fz)$ is equal to the union $\bigcup_{x\in X}\sigma(\pi^x(\fz))$.
\item[(b)] The family of operators $H_X:=(H_x)_{x\in X}$ where $H_x:=\pi^x(\fz)$ is equivariant (also called covariant), i.e., the equation
$$
H_{\alpha_h(x)} \; 
	= \; U_h H_x U_{h^{-1}}
$$
holds for all $h\in G$ where $U_h:\ell^2(G)\to\ell^2(G),\; U_h\psi(g):=\psi(h^{-1}g)$ is unitary.
\item[(c)] The map $X\ni x\mapsto \pi^x(\fz)$ is strongly continuous on $\Ll(\ell^2(G))$, i.e., the limit 
$$
\lim_{y\to x}\big\|\big(\pi^y(\fz)-\pi^x(\fz)\big)\psi\big\|
$$
is equal to zero for all $\psi\in\ell^2(G)$ and $x\in X$.
\end{itemize}
\end{proposition}

\begin{proof}
(a): According to Proposition~\ref{Chap2-Prop-TransGroupAmen} and Proposition~\ref{Chap2-Prop-Etale}, the groupoid $\Gamma$ is topologically amenable, \'etale with compact unit space. Thus, the desired identity $\sigma(\fz)=\bigcup_{x\in X}\sigma(\pi^x(\fz))$ follows by Theorem~\ref{Chap2-Theo-SpectrGroupoidCAlg}.

\vspace{.1cm}

(b): Let $x\in X$ and $h\in G$. Then, for every $\psi\in\ell^2(G)$ and $g\in G$, a short computation leads to
\begin{align*}
\left(H_{\alpha_h(x)}\psi\right)(g) \; 
	&= \; \sum\limits_{\tilde{h}\in G} \,
		\fz\left(\left. 
			\alpha_{g^{-1}}\left( \alpha_h(x) \right) \right| g^{-1}\tilde{h}
		\right) \cdot \psi(\tilde{h})\\
	&= \; \sum\limits_{\tilde{h}\in G} \,
		\fz\left(\left. 
			\alpha_{(h^{-1}g)^{-1}}(x) \right| (h^{-1}g)^{-1} \, h^{-1}\tilde{h}
		\right) \cdot \left(U_{h^{-1}}\psi\right)(h^{-1}\tilde{h})\\
	&= \; \left( \pi^x(\fz) \left(U_{h^{-1}}\psi\right) \right)(h^{-1}g)\\
	&= \; \left(U_h H_x U_{h^{-1}}\psi\right)(g)\,.
\end{align*}

(c): The \'etale property of the groupoid $\Gamma$ implies that $G$ is discrete, c.f. Proposition~\ref{Chap2-Prop-Etale}. Since $\Cc_c\big(\Gamma^{(1)}\big)\subseteq \CG^\ast(\Gamma)$ is a dense subset, it suffices to show the strong continuity for all normal elements of $\Cc_c\big(\Gamma^{(1)}\big)$. Let $\fz\in\Cc_c\big(\Gamma^{(1)}\big)$ be normal. First, (i) it is shown that $\big\|\big(\pi^x(\fz)-\pi^y(\fz)\big)\psi\big\|$ tends to zero if $x$ for all $\psi\in\Cc_c(G)$. Secondly, (ii) a $3\varepsilon$-argument leads to the desired strong continuity of the map $X\ni x\mapsto\pi^x(\fz)$.

\vspace{.1cm}

(i): Consider a $\psi\in\Cc_c(G)$. Then, for $x,y\in X$, the equation
$$
\left\|\left(\pi^x(\fz)-\pi^y(\fz)\right) \psi\right\|^2 \;
	= \; \sum\limits_{g\in G}  
		\left|
			\sum\limits_{\tilde{h}\in G} 
				\left(
					\fz\left(\left. \alpha_{g^{-1}}(x) \right| g^{-1}\tilde{h} \right) - \fz\left(\left. \alpha_{g^{-1}}(y) \right| g^{-1}\tilde{h} \right)
				\right) \cdot \psi(\tilde{h})
		\right|^2
$$
holds. Since $\fz\in\Cc_c\big(\Gamma^{(1)}\big)$ and $\psi\in\Cc_c(G)$ are compactly supported, the sums are finite. Thus, the continuity of $\fz$ implies that the field of bounded normal operators $(\pi^x(\fz))$ is strongly continuous. 

\vspace{.1cm}

(ii): Let $\varphi\in\ell^2(G)$ and $\varepsilon>0$. Since $\Cc_c(G)\subseteq\ell^2(G)$ is dense, there is a $\psi\in\Cc_c(G)$ such that $\|\psi-\varphi\|<\frac{\varepsilon}{3C}$ where $C:=\|\fz\|=\sup_{x\in X}\|\pi^x(\fz)\|<\infty$. Then choose, by (i), an open neighborhood $U\subseteq X$ of $x$ such that $\big\|\big(\pi^x(\fz)-\pi^y(\fz)\big)\psi\big\|<\frac{\varepsilon}{3}$ holds for all $y\in U$. Consequently, the estimate
\begin{align*}
\big\|\big(\pi^x(\fz)-\pi^y(\fz)\big)\varphi\big\| 
	\leq  \big\|\pi^x(\fz)\big\| \!\cdot\! \big\|\varphi-\psi\big\| + \big\|\big(\pi^x(\fz)-\pi^y(\fz)\big)\psi\big\| + \big\|\pi^y(\fz)\big\| \!\cdot\! \big\|\varphi-\psi\big\| 
	< \varepsilon
\end{align*}
is derived for all $y\in U$.
\end{proof}

\begin{remark}
\label{Chap2-Rem-CovFamOp-Spect}
Note that Proposition~\ref{Chap2-Prop-CovFamOp-Spect}~(b) implies the equalities $\sigma\big(H_{\alpha_h(x)}\big)=\sigma\big(H_x\big)$, $x\in X\,,$ where $H_x=\pi^x(\fz)$ and $\fz\in\CG^\ast(\Gamma)$ is normal. Furthermore, Proposition~\ref{Chap2-Prop-CovFamOp-Spect}~(c) yields a semi-continuity of $X\ni x\mapsto\sigma\big(H_x\big)\in\ks(\RM)$, c.f. Section~\ref{Chap3-Sect-OpTop}.
\end{remark}

In the following, it is shown that dynamical properties are encoded in spectral properties of a family of Schr\"odinger operators. More precisely, the dynamical system is minimal if and only if, for each self-adjoint (normal) generalized Schr\"odinger operator, the associated spectra are constant. This provides the first connection between a dynamical system and the spectrum of the associated generalized Schr\"odinger operators. As we will see in Section~\ref{Chap4-Sect-ContSpecGenSchrOp}, also the continuous variation of dynamical systems is characterized by the continuous behavior of the corresponding spectra for all generalized Schr\"odinger operators, c.f. Theorem~\ref{Chap4-Theo-CharDynSystConvSpectr} and Theorem~\ref{Chap4-Theo-CharSubshiftConvSpectr}. It is essential for both results that a generalized Schr\"odinger operators can be viewed as continuous functions on the dynamical system. This allows to connect the spectrum by using the Lemma of Urysohn.

\medskip

In order to do so the functional calculus of elements of operators and $C^\ast$-algebras is used. More precisely, let $\fz$ be a normal element of a unital $C^\ast$-algebra $\AG$. Then there is a unique morphism $\Phi$ of the $C^\ast$-algebra of continuous, complex-valued function on the spectrum $\sigma(\fz)$ into $\AG$ such that $\Phi(1)=\mathpzc{1}$ and $\Phi(id)=\fz$ where $id:\sigma(\fz)\to\sigma(\fz)\,,\; z\mapsto z$. Additionally, $\Phi$ is an isomorphism if the image is restricted to $\AG(\fz,\mathpzc{1})$ where $\AG(\fz,\mathpzc{1})$ denotes the sub-$C^\ast$-algebra of $\AG$ generated by the identity $\mathpzc{1}$ and $\fz$. Then, for every continuous function $\phi:\CM\to\CM$, the element $\Phi(\phi)\in\AG$ is denoted by $\phi(\fz)$. The identity $\sigma(\phi(\fz))=\phi(\sigma(\fz))$ follows by construction. In this work, polynomials of operators play an important role. Let $p(z,\overline{z})$ be a complex-valued polynomial in $z$ and $\overline{z}$ and $\phi:\CM\to\CM,\; \phi(z):=p(z,\overline{z})$. Then $\phi(\fz)$ is the usual algebraic combination $p(\fz,\fz^\ast)$. For a more detailed discussion on the functional calculus, the reader is referred to \cite[Section 1.5]{Dixmier77} and \cite[Section~2.5]{Murphy90}.

\medskip

With this at hand, the spectrum of a normal element can be characterized as a set by using Urysohn's Lemma \cite{Ury25}.

\begin{lemma}
\label{Chap2-Lem-SpectrUrysohn}
Let $\hs$ be a Hilbert space and $A\in\Ll(\hs)$ be a normal operator. Then the following assertions are equivalent for $\lambda\in\CM$.
\begin{itemize}
\item[(i)] The value $\lambda$ is an element of the spectrum $\sigma(A)$.
\item[(ii)] The normal element $\phi(A)\in\Ll(\hs)$ is non-zero for all $\phi\in\Cc_c(\CM)$ with $\phi(\lambda)=1$.
\end{itemize}
\end{lemma}

\begin{proof}
Let $\lambda\in\CM$ and $A\in\Ll(\hs)$ be normal. Denote by $\Phi$ the $\ast$-isomorphism from the continuous function on $\sigma(A)$ to the $C^\ast$-algebra generated by the identity operator $\mathpzc{1}$ on $\hs$ and $A$.

\vspace{.1cm}

(i)$\Rightarrow$(ii): Suppose $\lambda\in\sigma(A)$ and consider a $\phi\in\Cc_c(\CM)$ with $\phi(\lambda)=1$. Then $\phi$ is a non-zero element in the $C^\ast$-algebra of continuous functions on the spectrum of $A$. By the functional calculus, $\phi(A)$ is defined by the image $\Phi(\phi)$. Consequently, $\phi(A)$ is non-zero as $\phi$ does not vanish on the spectrum $\sigma(A)$ and $\Phi$ is a $\ast$-isomorphism. 

\vspace{.1cm}

(ii)$\Rightarrow$(i): Assume that $\lambda\not\in\sigma(A)$ holds. Then there is a compact neighborhood $V$ of $\lambda\in\CM$ such that $V\cap\sigma(A)=\emptyset$. Urysohn's Lemma implies the existence of a $\phi\in\Cc_c(\CM)$ satisfying $\phi|_{\CM\setminus V}\equiv 0$ and $\phi(\lambda)=1$, c.f. Proposition~\ref{App1-Prop-LemmaUrysohn}. Hence, $\phi$ vanishes on the spectrum of $A$. Thus, $\phi(A)=\Phi(\phi)=0$ is derived leading to a contradiction with (ii). 
\end{proof}

\medskip

It is well-known that the spectrum is constant for a minimal dynamical system $(X,G,\alpha)$ and a self-adjoint equivariant family of operators $H_X=(H_x)_{x\in X}$ such that $X\ni x\mapsto H_x\in\Ll\big(\ell^2(G)\big)$ is strongly continuous, c.f. \cite{CyconFroeseKirschSimon87,BIST89,Jit95,Len99,LeSt03-Algebras}. The proof provided here follows the lines of \cite[Theorem~4.3]{LeSt03-Algebras} where a similar assertion is proven for Delone dynamical systems. The proof is essentially based on the functional calculus for normal elements of a $C^\ast$-algebra and Urysohn's Lemma \cite{Ury25}.

\begin{theorem}[\cite{LeSt03-Algebras}]
\label{Chap2-Theo-MinCharConstSpectr}
Let $(X,G,\alpha)$ be a topological dynamical system where $G$ is discrete and $\Gamma:=X\rtimes_\alpha G$ is the associated transformation group groupoid. Then the following assertions are equivalent.
\begin{itemize}
\item[(i)] The dynamical system $(X,G,\alpha)$ is minimal.
\item[(ii)] For each self-adjoint (normal) element $\fz\in\CG_{red}^\ast(\Gamma)$, the spectrum of $\sigma\big(\pi^x(\fz)\big)$ is independent of $x\in X$, i.e., $\sigma(\fz)=\sigma\big(\pi^x(\fz)\big)$ for all $x\in X$.
\item[(iii)] The representation $\pi^x:\CG^\ast_{red}(X\rtimes_\alpha G)\to\Ll\big(\ell^2(G)\big)$ is faithful for every $x\in\Gamma^{(0)}$.
\end{itemize}
\end{theorem}

\begin{proof}
According to Proposition~\ref{Chap2-Prop-TransGrGroupProp}, the unit space $\Gamma^{(0)}$ is equal to $X\times\{e\}$. For simplification of the notation every element $(x,e)\in\Gamma^{(0)}$ is identified with $x\in X$.

\vspace{.1cm}

(i)$\Rightarrow$(ii): Let $\fz\in\CG_{red}^\ast(\Gamma)$ be normal and $\phi\in\Cc_c(\CM)$. Consider an $x\in X$. Due to Proposition~\ref{Chap2-Prop-LeftRegulRepres}, $\pi^x:\CG_{red}^\ast(\Gamma)\to\Ll\big(\ell^2(G)\big)$ is a $\ast$-homomorphism implying that $p\big(\pi^x(\fz),\pi^x(\fz)^\ast\big) = \pi^x\big( p(\fz,\fz^\ast)\big)$ holds for each complex-valued polynomial $p:\CM\to\CM\,,\; z\in\CM\mapsto p(z,\overline{z})$. Due to \cite[Theorem~2.1.7]{Murphy90}, $\pi^x$ is also continuous for every $x\in X$. Hence, the equality $\phi\big(\pi^x(\fz)\big) = \pi^x\big(\phi(\fz)\big)$ is derived by using that a continuous function on $\sigma(\fz)\subseteq\CM$ can be uniformly approximated by polynomials. With this at hand, the constancy of the spectrum is concluded as follows.

\vspace{.1cm}

Consider the set $X_0:=\big\{ x\in X\;|\; \pi^x\big(\phi(\fz)\big)=0 \big\}\subseteq X$. Then the set $X_0$ is $G$-invariant and closed by Proposition~\ref{Chap2-Prop-CovFamOp-Spect}~(b) and (c). Thus, $X_0$ is either empty or equal to $X$ since $(X,G,\alpha)$ is minimal by (i). 

\vspace{.1cm}

Let $x,y\in X$. Consider a $\lambda\not\in\sigma\big(\pi^x(\fz)\big)$. According to Lemma~\ref{Chap2-Lem-SpectrUrysohn}, there is a $\phi\in\Cc_c(\CM)$ satisfying $\phi(\lambda)=1$ and $\phi\big(\pi^x(\fz)\big)= 0$. Consequently, the set $X_0$ is non-empty and so $X_0$ is equal to $X$ by the previous considerations. Thus, $\phi\big(\pi^y(\fz)\big)= 0$ is deduced implying $\lambda\not\in\sigma\big(\pi^y(\fz)\big)$ by Lemma~\ref{Chap2-Lem-SpectrUrysohn}. Hence, the inclusion $\sigma\big(\pi^y(\fz)\big)\subseteq\sigma\big(\pi^x(\fz)\big)$ follows. Altogether, the equation $\sigma\big(\pi^x(\fz)\big)=\sigma\big(\pi^y(\fz)\big)$ is derived by interchanging the role of $x$ and $y$.

\vspace{.1cm}

(ii)$\Rightarrow$(iii): Let $x\in X$. Consider an $\fz\in\CG_{red}^\ast(\Gamma)$. If $\pi^x(\fz)=0$, the spectrum $\sigma\big(\pi^x(\fz)\big)$ is equal to $\{0\}$. Thus, (ii) implies that $\sigma\big(\pi^y(\fz)\big)=\{0\}$ holds for all $y\in X$. Hence, $\pi^y(\fz)=0\,,\; y\in X\,,$ follows. Consequently, $\fz=0$ is deduced since $\big(\pi^x,L^2(\Gamma^x,\mu^x)\big)_{x\in \Gamma^{(0)}}$ defines a faithful family of representations for the reduced $C^\ast$-algebra $\CG^\ast_{red}(\Gamma)$, c.f. Proposition~\ref{Chap2-Prop-LeftRegulRepres}. Since $\fz=0$ implies $\pi^y(\fz)=0$ for each $y\in X$, the previous considerations yield the equivalence of $\pi^x(f)=0$ and $f=0$. Hence, the representation $\big(\pi^x,L^2(\Gamma^x,\mu^x)\big)$ is faithful.

\vspace{.1cm}

(iii)$\Rightarrow$(i): Assume that the dynamical system $(X,G,\alpha)$ is not minimal. Then there are $x,y\in X$ such that $y\not\in\overline{\Orb(x)}\subseteq X$. The topological space $X$ is by definition second-countable, locally compact and Hausdorff. Thus, $X$ is a normal space by Proposition~\ref{App1-Prop-LocCompHausSecCountImplNormal} and so Urysohn's Lemma applies. Consequently, there exists a continuous, compactly supported function $\varphi:X\to[0,1]$ such that $\varphi(y)=1$ and $\varphi|_ {\overline{\Orb(x)}}\equiv 0$, c.f. Proposition~\ref{App1-Prop-LemmaUrysohn}. Define $\fz:=\varphi\times \delta_e$ where $\delta_e:G\to\{0,1\}$ is the Kronecker delta function at the neutral element $e\in G$. Then $\fz:\Gamma^{(1)}\to\RM$ is a continuous function with compact support. Hence, $\fz$ is a self-adjoint element of $\CG^\ast_{red}(\Gamma)$ since $\fz$ is real-valued and the support of $\fz$ contained in $X\times\{e\}$. Thanks to Proposition~\ref{Chap2-Prop-LeftRegReprTransGrouGrou}, the equations 
$$
\big(\pi^x(\fz)\psi\big)(g)\;
	= \; \sum_{h\in G} 
		\underbrace{\varphi\big(\alpha_{g^{-1}}(x)\big) }_{=0}
		\cdot \delta_e\big(g^{-1}h\big)\cdot \psi(h) \;
	= \; 0\,,
	\qquad \psi\in\ell^2(G)\,,\; g\in G\,,
$$
are derived leading to $\pi^x(\fz)=0$. On the other hand,
$$
\big\|\pi^y(\fz)\big\|^2 \;
	\geq \; \big\|\pi^y(\fz)\delta_e\big\|^2 \;
	\geq \; \big| \pi^y(\fz)\delta_e(e)\big| \;
	= \; \left|
		\sum_{h\in G} \varphi(y)\cdot \delta_e(h)\cdot \delta_e(h)
	\right| \;
	= \; 1
$$
is concluded. Thus, $\pi^y(\fz)\neq 0$ is deduced leading to $\fz\neq 0$ since $\big(\pi^x,L^2(\Gamma^x,\mu^x)\big)_{x\in \Gamma^{(0)}}$ defines a faithful family of representations for the reduced $C^\ast$-algebra $\CG^\ast_{red}(\Gamma)$, c.f. Proposition~\ref{Chap2-Prop-LeftRegulRepres}. This contradicts  assertion (iii) as $\pi^x$ is not faithful by the previous considerations.
\end{proof}

\begin{remark}
\label{Chap2-Rem-PTheory-ConstSpectr}
The constancy of the spectra was extended to general families of operators over a dynamical system $(X,G,\alpha)$ in a joint work with {\sc D. Lenz}, {\sc M. Lindner} and {\sc C. Seifert}, c.f. \cite{BeLeMaCh14,BeLeMaCh16}. The techniques that are used there differ from the methods of the papers previously mentioned. The main difference is that for non self-adjoint families of operators the strong operator topology does not longer imply the semi-continuity of the spectra, c.f. \cite[Example~IV.3.8]{Kato95}. Instead, the works \cite{BeLeMaCh14,BeLeMaCh16} take advantage of the so called $P$-continuity, some basic facts about dynamical systems and a result by {\sc Seidel} and {\sc Silbermann} \cite[Theorem~1.28]{SeSi12}. Specifically, a family of operators $(H_x)_{x\in X}$ is considered for a topological dynamical system $(X,G,\alpha)$ and $H_x$ is acting on the Banach space 
$$
\ell^p(G,Y)\;	
	= \; \left\{
			\psi:G\to Y \;\left|\; \sum_{g\in G} \|\psi(g)\|_Y^p<\infty
		\right.\right\}
$$
where $1\leq p\leq \infty$ and $Y$ is a Banach space. Following \cite{BeLeMaCh16}, $(H_x)_{x\in X}$ is called a family of operators over $(X,G,\alpha)$ if the family is equivariant and it is $P$-strongly continu\-ous. In detail, the strong continuity of the operators $(H_x)_{x\in X}$ in the self-adjoint case is replaced by the so called $P$-strong continuity. The $P$-theory has its roots in the works of {\sc Simonenko}, {\sc Pr\"ossdorf}, {\sc Roch} and {\sc Silbermann} \cite{Si68,RoSi89,ProessdorfSilbermann91}. The reader is referred to \cite{Lindner06,ChandlerLindner11,SeSi12,Sei14,BeLeMaCh14,BeLeMaCh16} for more details on the concept of $P$-strong continuity. If, additionally, $Y$ is a finite dimensional Banach space and $1<p<\infty$ holds, then this $P$-continuity turns out to be equivalent to the $\ast$-strong continuity of $(H_x)_{x\in X}$, c.f. \cite[Remark~2.5]{BeLeMaCh16}. This fact is used in the following theorem. Recall that a family of operators $(H_x)_{x\in X}$ defined on a Hilbert space $\hs$ is called $\ast$-strongly continuous if the maps $X\ni x\mapsto H_x\in\Ll(\hs)$ and $X\ni x\mapsto H_x^\ast\in\Ll(\hs)$ are both strongly continuous.
\end{remark}

{\sc Lenz}, {\sc Peyerimhoff} and {\sc Veseli\'c} \cite{LePeIv07} verified the almost sure absence of the discrete spectrum for operators affiliated to a von Neumann algebra. In \cite{BeLeMaCh14,BeLeMaCh16}, it is additionally proven that the absence of the discrete spectrum holds for all operators over a minimal dynamical system.

\begin{theorem}[\cite{BeLeMaCh14,BeLeMaCh16}]
\label{Chap2-Theo-ConstSpectrMinimal}
Let $(X,G,\alpha)$ be a minimal, topological dynamical system such that $\Gamma:=X\rtimes_\alpha G$ is topologically amenable with Haar system $\mu$ and $G$ is a discrete, countable group. For a normal element $\fz\in\CG^\ast(\Gamma)$, consider the family of operators $H_X:=(H_x)_{x\in X}$ where $H_x:=\pi^x(\fz):\ell^2(G)\to\ell^2(G)$. Then the equations 
$$
\sigma(H_X) \; 
	:= \; \sigma(\fz) \;
	= \; \sigma(H_x) \; 
	= \; \sigma_{ess}(H_x)
$$
hold for all $x\in X$.
\end{theorem}

\begin{proof}
The equation $\sigma(\fz)=\sigma(H_x)$ is an immediate consequence of Proposition~\ref{Chap2-Prop-CovFamOp-Spect}~(a) and Theorem~\ref{Chap2-Theo-MinCharConstSpectr}. Thus, only $\sigma(H_x)= \sigma_{ess}(H_x)$ is left to prove.

\vspace{.1cm}

The operators $H_x$ are defined on the Hilbert space $\ell^2(G)=\ell^2(G,\CM)$, namely $p=2$ and $Y=\CM$ in the notation of \cite{BeLeMaCh16} described above. Thus, $Y$ is finite dimensional and $1<p<\infty$ is valid. Consequently, $(H_x)_{x\in X}$ is $P$-strongly continuous if and only if the family is $\ast$-strongly continuous, c.f. \cite[Remark~2.5]{BeLeMaCh16}. According to Proposition~\ref{Chap2-Prop-CovFamOp-Spect}~(c), the family $(H_x)_{x\in X}$ is strongly continuous. Since the family of operators $(H_x)_{x\in X}$ is normal, the identities 
\begin{align*}
\big\| \big(H_x-H_y\big)\psi \big\|^2 \; 
	&\overset{\qquad\;}{=} \; \big\langle\big(H_x-H_y\big)\psi \;|\; \big(H_x-H_y\big)\psi \big\rangle\\
	&\overset{\text{normal}}{=} \; \big\langle\psi \;|\; \big(H_x-H_y\big)\big(H_x-H_y\big)^\ast\psi \big\rangle\\
	&\overset{\qquad\;}{=} \: \big\langle\big(H_x-H_y\big)^\ast\psi \;|\; \big(H_x-H_y\big)^\ast\psi \big\rangle\\
	&\overset{\qquad\;}{=} \;\big\| \big(H_x-H_y\big)^\ast\psi \big\|^2
\end{align*}
follow for each $\psi\in\ell^2(G)$. Thus, $(H_x)_{x\in X}$ is $\ast$-strongly continuous and so it is $P$-strongly continuous by the previous considerations. Together with Proposition~\ref{Chap2-Prop-CovFamOp-Spect}~(b), this implies that $H_X=(H_x)_{x\in X}$ is a family of operators over the dynamical system $(X,G,\alpha)$ in terms of \cite{BeLeMaCh16}. Then \cite[Corollary~4.5]{BeLeMaCh16} leads to the desired identity
$$
\sigma(H_X) \;
	= \; \sigma_{ess}(H_x) \;
	= \; \sigma(H_x)
	\,, \qquad x\in X\,,
$$
since $(X,G,\alpha)$ is minimal by assumption.
\end{proof}

\section{Pattern equivariant Schr\"odinger operators}
\label{Chap2-Sect-PattEqSchrOp}

The chapter is finished with a discussion about a certain class of generalized Schr\"odinger operators of finite range associated with symbolic dynamical systems $(\as^G,G,\alpha)$ defined in Section~\ref{Chap2-Ssect-SpaceSubshifts}. These operator families $H_X:=(H_x)_{x\in X}$ on $\ell^2(G)$ include those Schr\"odinger and Jacobi operators that are usually considered in the community.

\medskip

{\sc Kellendonk} and {\sc Putnam} introduced the notion of {\em pattern equivariant functions} in \cite{KePu00,Kel03}. It provides an intuitive description of functions depending on the underlying geometry. Schr\"odinger operators with a potential coming from a pattern equivariant function lead to interesting spectral properties. In particular, it implies a strong tunneling of the particles, i.e., singular continuous spectrum of the associated Schr\"odinger operator, c.f. \cite{BIST89,BeBoGh91,BoGh93,Len02,BePo13}. Following \cite{Kel03} the notion of pattern equivariant functions is defined for functions on $\as^G$.

\begin{definition}[Pattern equivariant function, \cite{KePu00,Kel03}]
\label{Chap2-Def-PattEqFunc}
Let $G$ be a countable, discrete group and $\as$ be an alphabet. For a subshift $\Xi\in\SG\big(\as^G\big)$, a function $p:\Xi\to\CM$ is called {\em pattern equivariant function} if there exists a compact neighborhood of the neutral element $e\in G$ such that $p(\xi)=p(\eta)$ holds for all $\xi,\eta\in\Xi$ with $\xi|_K=\eta|_K$.
\end{definition}

In the following, let $G$ be a countable, discrete group and $\as$ be an alphabet. Recall the definition of a base for the topology on $\as^G$ introduced in Section~\ref{Chap2-Sect-SymbDynSyst}.

\begin{proposition}
\label{Chap2-Prop-CharPattEqFunc}
Let $\Xi\in\SG\big(\as^G\big)$ be a subshift. For every function $p:\Xi\to\CM$, the following assertions are equivalent. 
\begin{itemize}
\item[(i)] The function $p:\Xi\to\CM$ is pattern equivariant. 
\item[(ii)] There exist an $N\in\NM$, coefficients $p_j\in\CM$, compact sets $K_j\in\ks(G)$ and $[u_j]\in\as^{[K_j]}$ for $1\leq j\leq N$ such that 
$$
p=\sum\limits_{j=1}^N p_j\cdot\chi_{\os(K_j,[u_j])}.
$$
\item[(iii)] The function $p:\Xi\to\CM$ is continuous and takes finitely many values.
\end{itemize}
In particular, each pattern equivariant function $p:\Xi\to\CM$ on $\Xi$ is extendable to a pattern equivariant function on $\as^G$. Furthermore, the finite linear combination and multiplication of pattern equivariant functions is a pattern equivariant function, i.e., the set of pattern equivariant functions is closed with respect to algebraic combinations.
\end{proposition}

\begin{proof}
Condition (ii) implies that each pattern equivariant function $p:\as^G\to\CM$ is extendable to a pattern equivariant function on $\as^G$. Furthermore, the product of two characteristic functions $\chi_{\os(K,[u])}\cdot \chi_{\os(F,[v])}$ is either identically zero or $\chi_{\os(K\cup F,[w])}$ for a $[w]\in\as^{[K\cup F]}$. Thus, (ii) implies that the algebraic combination of pattern equivariant functions is again a pattern equivariant function.

\vspace{.1cm}

(i)$\Rightarrow$(ii): Let $p:\as^G\to\CM$ be a pattern equivariant function with associated compact neighborhood $K$ of $e\in G$. Set $N$ to be the number of elements of $U:=\{\xi|_K\;|\; \xi\in\Xi\}\subseteq\as^K$ which is finite. Choose a labeling of $U=\{u_1,\ldots,u_N\}$. Then for $1\leq j\leq N$ define $p_j:=p(\xi)$ for a fixed $\xi\in\Xi$ with $\xi|_K=u_j$. Since $p$ is a pattern equivariant function, the coefficient $p_j$ is independent of the choice of $\xi\in\Xi$. By construction, the equation $p=\sum_{j=1}^N p_j\cdot\chi_{\os(K_j,[u_j])}$ holds on $\Xi$.

\vspace{.1cm}

(ii)$\Rightarrow$(iii): Note that the sets $\os(K,[u])$ are compact and open for $K\in\ks(G)$ and $[u]\in\as^{[K]}$. Thus, the characteristic function $\chi_{\os(K,[u])}$ defines a continuous function. By definition, each function defined in (ii) takes only finitely many values, namely the values $p_j$ for $1\leq j\leq N$.

\vspace{.1cm}

(iii)$\Rightarrow$(i): If $p$ takes finitely many values $p_1,\ldots, p_N\in\CM$ then these values are separated by open sets. The continuity of $p$ then implies that, for each $p_j$, there exists an open neighborhood $\os(K_j,[u_j])$ with $K_j\in\ks(G)$ and a $[u_j]\in\as^{[K_j]}$ such that $p(\xi)=p_j$ holds for all $\xi\in\os(K_j,[u_j])$. Note that it is used that these sets form a basis for the topology on $\as^G$. Then $K:=\bigcup_{j=1}^N K_j\cup\{e\}$ defines a compact neighborhood of $\{e\}$ and, by construction, the equation $p(\xi)=p(\eta)$ holds for all $\xi,\eta\in\Xi$ with $\xi|_K=\eta|_K$.
\end{proof}

\medskip

Following Section~\ref{Chap2-Sect-TransformationGroupGroupoid}, a groupoid $\Gamma(\Xi):=\Xi\rtimes_\alpha G$ is associated with each subshift $\Xi\in\SG\big(\as^G\big)$. Let $\PG^\ast(\Gamma(\Xi))$ be the set of functions $\fz\in\Cc_c\big(\Gamma^{(1)}(\Xi)\big)$ such that, for each $g\in G$, the map $\Xi\ni\xi\mapsto\fz(\xi|g)\in\CM$ is a pattern equivariant function.

\begin{lemma}
\label{Chap2-Lem-PaEqAlgebra}
The set $\PG^\ast(\Gamma(\Xi))$ is a $\ast$-subalgebra of $\Cc_c\big(\Gamma^{(1)}(\Xi)\big)$.
\end{lemma}

\begin{proof}
It suffices to check that the set $\PG^\ast(\Gamma(\Xi))$ is invariant with respect to the convolution and involution defined in Proposition~\ref{Chap2-Prop-InvolutiveGroupoidAlgebra}. Compact subsets of $G$ are finite as $G$ is discrete. Furthermore, by Proposition~\ref{Chap2-Prop-CharPattEqFunc} and $\PG^\ast(\Gamma(\Xi))\subseteq\Cc_c\big(\Gamma^{(1)}(\Xi)\big)$ it follows that every $\fz\in\Cc_c\big(\Gamma^{(1)}(\Xi)\big)$ is a finite linear combination of $\delta_g \cdot \chi_{\os(K,[u])}$ for a $g\in G$, $K\in\ks(G)$ and $[u]\in\as^{[K]}$. Thus, it suffices to consider functions of the form $\delta_g \cdot \chi_{\os(K,[u])}$. A short computation leads to
\begin{align*}
\delta_g \cdot \chi_{\os(K,[u])}\star \delta_h \cdot \chi_{\os(F,[v])} \; 
	&= \; \delta_{gh} \cdot \chi_{\os(K,[u])} \cdot \chi_{\os(gF,[v])}\\
\left(\delta_g \cdot \chi_{\os(K,[u])}\right)^\ast \;
	&= \; \delta_{g^{-1}} \cdot \chi_{\os(g^{-1}K,[u])}
\end{align*}
where the function on the right hand side are elements of $\PG^\ast(\Gamma(\Xi))$.
\end{proof}

\medskip

The previous assertion justifies the following definition.

\begin{definition}[Pattern equivariant algebra]
\label{Chap2-Def-PaEqAlg}
Let $\Xi\in\SG\big(\as^G\big)$ be a subshift. A function $\fz\in\Cc_c\big(\Gamma^{(1)}(\Xi)\big)$ is called {\em pattern equivariant} whenever the maps $\Xi\ni\xi\mapsto\fz(\xi|g)\in\CM$ for all $g\in G$ are pattern equivariant. The $\ast$-algebra \gls{PGXi} is called the {\em pattern equivariant algebra}.
\end{definition}

Note that Proposition~\ref{Chap2-Prop-CharPattEqFunc} (iii) guarantees that $\PG^\ast(\Gamma(\Xi))$ is non-empty.

\begin{theorem}
\label{Chap2-Theo-PaEqAlgebraGenerator}
Let $\Xi\in\SG\big(\as^G\big)$ be a subshift. Then the set
$$
\AG(G,\as)\; 
	:= \; \left\{
		\sz_{g,a}:=\delta_g \cdot \chi_{\os(\{e\},[a])} \;|\; g\in G,\; a\in\as 
	\right\}
	\subsetneq\PG^\ast(\Gamma(\Xi))
$$
generates the pattern equivariant algebra $\PG^\ast(\Gamma(\Xi))$. If $G$ is finitely generated by the finite set $\{g_1,\ldots,g_n\}\subseteq G$ then
$$
\AG_{fin}(G,\as) \;
	:= \; \left\{
		\sz_{g_j,a}:=\delta_{g_j} \cdot \chi_{\os(\{e\},[a])} \;|\; 1\leq j\leq n,\; a\in\as 
	\right\}
	\subsetneq\PG^\ast(\Gamma(\Xi))
$$
generates the pattern equivariant algebra $\PG^\ast(\Gamma(\Xi))$. In particular, $\PG^\ast(\Gamma(\Xi))$ is finitely generated.
\end{theorem}

\begin{proof}
For $g\in G$ and $a\in\as$, a short computation leads to
$$
\delta_g \cdot \chi_{\os(\{e\},[a])}\star \delta_h \cdot \chi_{\os(\{e\},[b])} \; 
	= \; \delta_{gh} \cdot \chi_{\os(\{e\},[a])} \cdot \chi_{\os(\{g\},[b])}
		= \; \delta_{gh} \cdot \chi_{\os(\{e,g\},[u])}\\
$$
where $u:\{e,g\}\to\as$ is defined by $e\mapsto a, g\mapsto b$. Since the equation $\sum_{a\in\as} \sz_{g,a} = \delta_g$ holds and compact subsets of $G$ are finite, it follows that every function $\chi_{\os(K,[v])}$ is a finite algebraic combination of elements of $\AG(G,\as)$. Proposition~\ref{Chap2-Prop-CharPattEqFunc} (ii) concludes the proof of the first assertion.

\vspace{.1cm}

That $\AG_{fin}(G,\as)$ generates $\PG^\ast(\Gamma(\Xi))$ follows by using the observation
$$
\left(\delta_g \cdot \chi_{\os(\{e\},[a])}\right)^\ast \;
	= \; \delta_{g^{-1}} \cdot \chi_{\{\os(\{g^{-1}\},[a])},\qquad g\in G,\; a\in\as,
$$
and the previous considerations.
\end{proof}

\medskip

Following the more general assertions \cite[Section~5]{KePu00}, \cite[Section~4]{LeSt03-Algebras}, \cite[Proposition~4.13]{Whittaker05}, \cite[Corollary~4.6.4, Corollary~4.6.5]{StarlingThesis12} and \cite[Lemma~4.2, Proposition~4.3]{Sta14}, the pattern equivariant algebra is a dense subalgebra of the reduced and full $C^\ast$-algebra associated with $\Gamma(\Xi)$. Furthermore, whenever the group $G$ is finitely generated, the pattern equivariant algebra itself is finitely generated.

\medskip

For the following theorem, recall the notion of the inductive limit topology defined on $\Cc_c\big(\Gamma^{(1)}(\Xi)\big)$ for a subshift $\Xi\in\SG\big(\as^G\big)$. Since topological groupoids are second-countable, the inductive limit topology is described in terms of sequences: 
A sequence $(\fz_n)_{n\in\NM}\subseteq\Cc_c\big(\Gamma^{(1)}(\Xi)\big)$ converges in the inductive limit topology to $\fz\in\Cc_c\big(\Gamma^{(1)}(\Xi)\big)$ if there exist an $n_0\in\NM$ and a compact $K\in\ks(G)$ such that $\supp(\fz_n)\subseteq K$ for $n\geq n_0$ and $(f_n)_{n\geq n_0}$ converges uniformly to $\fz$ on $K$.

\begin{theorem}[\cite{KePu00}]
\label{Chap2-Theo-PaEqDense}
Let $\Xi\in\SG\big(\as^G\big)$ be a subshift. Then the pattern equivariant algebra $\PG^\ast(\Gamma(\Xi))$ is a dense subalgebra of $\CG^\ast_{red}(\Gamma(\Xi))$ and $\CG^\ast_{full}(\Gamma(\Xi))$. In particular, if $\Gamma(\Xi)$ is topologically amenable the $C^\ast$-algebra $\CG^\ast(\Gamma(\Xi))$ is the closure of pattern equivariant algebra $\PG^\ast(\Gamma(\Xi))$.
\end{theorem}

\begin{proof}
The set of pattern equivariant functions $p:\Xi\to\CM$ gets an algebra if equipped with the multiplication
$$
\chi_{\os(K,[u])}\cdot \chi_{\os(F,[v])} \;
	= \; \begin{cases}
		\chi_{\os(K\cup F,[w])}\quad &,\; u|_{K\cap F}=v|_{K\cap F} \, ,\\
		0\quad &,\; \text{otherwise} \, ,
	\end{cases}
$$
where $w\in\as^{K\cup F}$ is defined by
$$
w(g)\; 
	:= \; \begin{cases}
		u(g)\quad &,\; g\in K \, ,\\
		v(g)\quad &,\; g\in F\setminus K \, ,
	\end{cases} \in\as^{K\cup F}
$$
for $F,K\in\ks(G)$, $u\in\as^K$ and $v\in\as^F$. Since the sets $\os(K,[u])$ for $K\subseteq G$ compact and $u\in\as^K$ form a base of the topology on $\Xi$, the pattern equivariant functions are dense with respect to the sup-norm by the Theorem of Stone-Weierstrass.

\vspace{.1cm}

For $\fz\in\Cc_c\big(\Gamma^{(1)}(\Xi)\big)$, there is a compact (finite) subset $K\subseteq G$ such that $\supp(\fz)\subseteq \Xi\times K$. As $G$ is discrete, the equation $\fz=\sum_{g\in K} \fz(\cdot|g)\cdot \delta_g$ is deduced. For $g\in K$, there exists a sequence of pattern equivariant functions $\fz^g_n:\Xi\to\CM,\; n\in\NM,$ converging uniformly to $\fz(\cdot|g)$. By the finiteness of $K$, the pattern equivariant functions $\fz_n:=\sum_{g\in K} \fz^g_n\cdot \delta_g\,,\; n\in\NM\,,$ converge uniformly to $\fz$. Furthermore, it immediately follows by definition that $(\fz_n)_{n\in\NM}$ also converges in the inductive limit topology to $\fz$. Thus, $\PG^\ast(\Gamma(\Xi))$ is dense in $\Cc_c\big(\Gamma^{(1)}(\Xi)\big)$ with respect to the inductive limit topology.

\vspace{.1cm}

According to \cite[Proposition~II.1.4]{Renault80}, the topology induced by the $I$-norm $\|\cdot\|_I$ is coarser than the inductive limit topology. Hence, $\PG^\ast(\Gamma(\Xi))$ is dense in $\Cc_c\big(\Gamma^{(1)}(\Xi)\big)$ with respect to the I-norm defined in Definition~\ref{Chap2-Def-Inorm}. Then Proposition~\ref{Chap2-Prop-LeftRegulRepres} leads to the denseness of $\PG^\ast(\Gamma(\Xi))$ in the $C^\ast$-algebra $\CG^\ast_{red}(\Gamma(\Xi))$ as $\|\fz\|_{red}\leq\|\fz\|_I$. Since $\|\fz\|_{full}\leq\|\fz\|_I$ holds, the pattern equivariant algebra $\PG^\ast(\Gamma(\Xi))$ is also dense in $\CG^\ast_{full}(\Gamma(\Xi))$.
\end{proof}

\medskip

By the previous considerations, the class of pattern equivariant functions $\PG^\ast(\Gamma(\Xi))$ is a $\ast$-subalgebra of $\Cc_c\big(\Gamma^{(1)}(\Xi)\big)$. According to Proposition~\ref{Chap2-Prop-LeftRegReprTransGrouGrou}, the left-regular representation maps an element of $\Cc_c\big(\Gamma^{(1)}(\Xi)\big)$ to a family of operator on $\ell^2(G)$. In the fol\-low\-ing, a class of operators on $\ell^2(G)$ is introduced. These operators are typically considered in the community of Schr\"odinger operators and so they are of particular interest. It turns out that, for a subshift $\Xi$, these operators are exactly those one that arise by the left-regular representation of self-adjoint elements of $\PG^\ast(\Gamma(\Xi))$, c.f. Theorem~\ref{Chap2-Theo-PESchrOpFinRang}. For the following, note that compact sets of $G$ are characterized by finite subsets as $G$ is equipped with the discrete topology.
 
\begin{lemma}
\label{Chap2-Lem-PatEqSchrOp}
Consider a finite set $K\in\ks(G)$ and pattern equivariant functions $p_h:\as^G\to\CM,\; h\in K,$ and $p_e:\as^G\to\RM$. For $\xi\in\as^G$, define the operator $H_\xi:\ell^2(G)\to\ell^2(G)$ by
$$
(H_\xi\psi)(g) \; 
	:= \; \left( 
			\sum\limits_{h\in K} p_h\big(\alpha_{g^{-1}}(\xi)\big) \cdot \psi(g\,h^{-1}) + \overline{p_h\big(\alpha_{(gh)^{-1}}(\xi)\big)} \cdot \psi(gh)
		\right)
		+ p_e\big(\alpha_{g^{-1}}(\xi)\big) \cdot \psi(g)
$$
where $\psi\in\ell^2(G)$ and $g\in G$. Then $H_\xi$ is a bounded self-adjoint linear operator for $\xi\in\as^G$.
\end{lemma}

\begin{proof}
Let $\xi\in\as^G$. It immediately follows that $H_\xi$ is a linear operator. Since the maps $p_h,\; h\in K\cup\{e\},$ are pattern equivariant, they are uniformly bounded by a constant $C>0$, i.e., $\sup_{\xi\in\as^G}\sup_{h\in K\cup\{e\}} |p_h(\xi)|<C<\infty$. Hence, $H_\xi$ defines a bounded operator. Let $\psi,\varphi\in\Cc_c(G)$ where $\Cc_c(G)$ is a dense subset of $\ell^2(G)$. Then the equations
\begin{align*}
\langle H_\xi\psi,\varphi\rangle \;
	&= \; \sum_{g\in G} \overline{\big(H_\xi\psi\big)(g)} \cdot \varphi(g) \\
	&= \; \sum_{h\in K} \left(
		\sum_{g\in G} \overline{p_h\big(\alpha_{g^{-1}}(\xi)\big)} \cdot \varphi(g) \cdot \psi(gh^{-1})
		+ \sum_{g\in G} p_h\big(\alpha_{(gh)^{-1}}(\xi)\big) \cdot \varphi(g) \cdot \psi(gh)\right.\\
	&\qquad\qquad + \left.
		\sum_{g\in G} p_e\big(\alpha_{g^{-1}}(\xi)\big) \cdot \psi(g)\cdot\varphi(g)
		\right)\\
	&= \; \sum_{g\in G} \left(\left(
			\sum_{h\in K} \overline{p_h\big(\alpha_{(gh)^{-1}}(\xi)\big)} \cdot \varphi(gh) 
			+ p_h\big(\alpha_{g^{-1}}(\xi)\big) \cdot \varphi(gh^{-1})
			\right)\right.\\
	&\qquad\qquad + p_e\big(\alpha_{g^{-1}}(\xi)\big) \cdot \varphi(g)
			\bigg)\cdot \psi(g)\\
	&= \; \langle \psi,H_\xi\varphi\rangle
\end{align*}
hold. Consequently, $H_\xi$ defines a self-adjoint operator for $\xi\in\Xi$.
\end{proof}

\medskip

Note that without loss of generality the compact set $K$ can be chosen such that $h^{-1}\not\in K$ if $h\in K$. This holds since the sum of pattern equivariant functions is pattern equivariant, c.f. Proposition~\ref{Chap2-Prop-CharPattEqFunc}.

\newpage
\begin{definition}[pattern equivariant Schr\"odinger operator]
\label{Chap2-Def-SchrOp-l2(G)}
Let $\Xi\in\SG\big(\as^G\big)$ be a subshift and $K\in\ks(G)$ with pattern equivariant functions $p_h:\as^G\to\CM,\; h\in K,$ and $p_e:\as^G\to\RM$. Then the family of bounded self-adjoint linear operators $H_\xi:\ell^2(G)\to\ell^2(G),\; \xi\in\Xi\, ,$ defined in Lemma~\ref{Chap2-Lem-PatEqSchrOp} is called {\em pattern equivariant Schr\"odinger operator}.
\end{definition}

\begin{remark}
\label{Chap2-Rem-ExtPatEqFunct}
According to Proposition~\ref{Chap2-Prop-CharPattEqFunc} every pattern equivariant function $p:\Xi\to\CM$ is extendable to a pattern equivariant function $\tilde{p}:\as^G\to\CM$ such that $\tilde{p}|_\Xi=p$. Thus, there is no loss of generality if we assume directly that the pattern equivariant function is defined on $\as^G$ instead of $\Xi$. For the further applications, c.f. Section~\ref{Chap4-Sect-ContSpecGenSchrOp} and Remark~\ref{Chap4-Rem-ChoicePEFunctOnAsG}, it is useful to define the functions directly on $\as^G$ so that there is no need to handle with extensions. Note that all the results are independent of the choice of the extension, c.f. Remark~\ref{Chap4-Rem-ChoicePEFunctOnAsG}~(ii).
\end{remark}

The compact (finite) set $K\cup K^{-1}\cup\{e\}\subseteq G$ of a pattern equivariant Schr\"odinger operator $H_\Xi$ is interpreted as the range of $H_\Xi$. The following statement provides the link between the abstract $C^\ast$-algebra of a groupoid and the pattern equivariant Schr\"odinger operators on $\ell^2(G)$. It is shown that a pattern equivariant Schr\"odinger operator is a generalized Schr\"odinger operator of finite range.

\begin{theorem}
\label{Chap2-Theo-PESchrOpFinRang}
Let $\Xi\in\SG\big(\as^G\big)$ be a subshift.
\begin{itemize}
\item[(a)] Consider a self-adjoint element $\fz\in\PG^\ast(\Gamma(\Xi))$. Then $H_\Xi:=(H_\xi)_{\xi\in\Xi}$ defined by $H_\xi:=\pi^\xi(\fz)$ is a pattern equivariant Schr\"odinger operator.
\item[(b)] Let $H_\Xi:=(H_\xi)_{\xi\in\Xi}$ be a pattern equivariant Schr\"odinger operator on $\ell^2(G)$ with finite $K\in\ks(G)$ and pattern equivariant functions $p_h:\as^G\to\CM,\; h\in K\, ,$ and $p_e:\as^G\to\RM$. Then there exists a self-adjoint $\fz\in\PG^\ast(\Gamma(\Xi))$ such that $H_\xi=\pi^\xi(\fz)$ for all $\xi\in\Xi$.
\end{itemize}
In particular, $H_\Xi$ is a generalized Schr\"odinger operator of finite range and the equations $\sigma(\fz)=\overline{\bigcup_{\xi\in\Xi} \sigma(H_\xi)} = \sigma(H_\Xi)$ hold. If $G$ is additionally amenable, then $\bigcup_{\xi\in\Xi} \sigma(H_\xi) = \sigma(H_\Xi)$.
\end{theorem}

\begin{proof}
As $\PG^\ast(\Gamma(\Xi))\subseteq\Cc_c\big(\Gamma^{(1)}(\Xi)\big)$, the family of $H_\Xi$ is a generalized Schr\"odinger operator of finite range, c.f. Definition~\ref{Chap2-Def-SchrodingerOperator}. Then the equations $\sigma(\fz)=\bigcup_{\xi\in\Xi} \sigma(H_\xi)= \sigma(H_\Xi)$ immediately follow by \cite[Proposition~II.1.11]{Renault80} and \cite[Proposition~2.5]{NiPr15}. Whenever $G$ is additionally amenable, the groupoid $\Gamma(\Xi):=\Xi\rtimes_\alpha G$ is topologically amenable by Proposition~\ref{Chap2-Prop-TransGroupAmen}. Thus, $\bigcup_{\xi\in\Xi} \sigma(H_\xi)= \sigma(H_\Xi)$ follows by Proposition~\ref{Chap2-Prop-CovFamOp-Spect}.

\vspace{.1cm}

(a): Let $\xi\in\Xi$ and $\fz\in\PG^\ast(\Gamma(\Xi))$ be self-adjoint. The group $G$ is discrete and so the function $\fz$ is equal to $\sum_{h\in G} \delta_h \cdot \fz(\cdot| h)$. Since $\fz$ has compact support, there exists a $K\in \ks(G)$ such that 
\begin{itemize}
\item[(1)] $h^{-1}\not\in K$ for all $h\in K$,
\item[(2)] $\fz(\cdot|g)\equiv 0$ if $g\in G\setminus (K\cup\{e\})$,
\item[(3)] $\fz(\cdot|h)\not\equiv 0$ if $h\in K$.
\end{itemize}
Note that the set $K$ is not unique as, for instance, $\tilde{K}:= K\setminus\{h\}\cup\{h^{-1}\}$ fulfills the same properties for a $h\in K$. Since $\fz$ is self-adjoint, the identities 
$$
\fz(\xi|h)\; 
	= \; \fz^\ast(\xi|h)\;
		= \; \overline{\fz\big(\alpha_{h^{-1}}(\xi)|h^{-1} \big)}
			,\qquad \xi\in\as^G\, ,\; h\in K\cup\{e\}\, ,
$$
hold. Hence, $\fz(\cdot|h)\not\equiv 0$ is equivalent to $\fz(\cdot|h^{-1})\not\equiv 0$ for $h\in G$. Thus, $\fz(\cdot|e)$ is a real-valued function defined on $\as^G$ which is possibly equal to zero. Then $\fz$ is equal to 
$$
\sum_{g\in K} \delta_h \cdot \fz(\cdot| h) + \delta_{h^{-1}}\cdot\fz(\cdot|h^{-1}) + \delta_e\cdot \fz(\cdot|e) \;
	= \; \sum_{g\in K} \delta_h \cdot \overline{\fz\big(\alpha_{h^{-1}}(\cdot)|h\big)} + \delta_{h^{-1}}\cdot\fz(\cdot|h^{-1})  + \delta_e\cdot \fz(\cdot|e)
$$
where the second equation follows by the previous considerations. Note that the representation of $\fz$ is independent of the choice of the set $K$. For $h\in K$, $g\in G$ and $\psi\in\ell^2(G)$, a short computation leads to
\begin{align*}
\left( \pi^\xi\Big(\delta_{h^{-1}}\cdot\fz(\cdot|h^{-1})\Big) \psi \right)(g) \; 
	&= \; \sum\limits_{\tilde{h}\in G} 
		\delta_{h^{-1}} \left(g^{-1}\tilde{h}\right) \cdot \fz\big(\alpha_{g^{-1}}(\xi)|h^{-1}\big) \cdot \psi(\tilde{h})\\
	&= \; \fz\big(\alpha_{g^{-1}}(\xi)|h^{-1}\big) \cdot \psi\left( g h^{-1} \right)\, ,\\
\left( \pi^\xi\Big( \delta_h \cdot \overline{\fz\big(\alpha_{h^{-1}}(\cdot)|h\big)} \Big) \psi \right)(g) \; 
	&= \; \overline{\fz\big( \alpha_{(gh)^{-1}}(\xi)|h^{-1} \big)} \cdot \psi(gh) \, ,\\
\left( \pi^\xi\Big( \delta_e \cdot \fz(\cdot|e) \Big) \psi \right)(g) \; 
	&= \; \fz\big(\alpha_{g^{-1}}(\xi)|e\big) \cdot \psi\left( g \right) \, .
\end{align*}
For $h\in K\cup \{e\}$, define $p_h:\as^G\to\CM$ by $p_h(\cdot):=\fz(\cdot|h^{-1})$. Since $\fz\in\PG^\ast(\Gamma(\Xi))$, the functions $p_h\,,\; h\in K\cup \{e\}\,,$ are pattern equivariant. Using the linearity of the left-regular representation, this implies the equality
$$
(H_\xi\psi)(g) \; 
	= \; \left( 
			\sum\limits_{h\in K} p_h\big(\alpha_{g^{-1}}(\xi)\big) \cdot \psi(gh^{-1}) + \overline{p_h\big(\alpha_{(gh)^{-1}}(\xi)\big)} \cdot \psi(gh)
		\right)
		+ p_e\big(\alpha_{g^{-1}}(\xi)\big) \cdot \psi(g) \, ,
$$
for $\psi\in\ell^2(G)$ and $g\in G$. Thus, $H_\Xi:=(H_\xi)_{\xi\in\Xi}$ defined by $H_\xi:=\pi^\xi(\fz)$ is a pattern equivariant Schr\"odinger operator.

\vspace{.1cm}

(b): Let $K\in\ks(G)$ with pattern equivariant functions $p_h:\as^G\to\CM,\; h\in K,$ and $p_e:\as^G\to\RM$. Define the functions $\fz_h\in\PG^\ast(\Gamma(\Xi)),\; h\in K\cup\{e\},$ by
$$
\fz_h(\xi|g) \;
	:= \; \delta_{h^{-1}}(g)\cdot p_h(\xi) \, ,
	\qquad (\xi|g)\in\Gamma^{(1)}(\Xi) \, .
$$
Then
$$
\fz \;
	:= \; \left(
		\sum\limits_{h\in K} \fz_h + \fz_h^\ast
	\right) 
	+ \fz_e
$$
is an element of the pattern equivariant algebra $\PG^\ast(\Gamma(\Xi)) \subseteq \Cc_c\big(\Gamma^{(1)}(\Xi)\big)$. Since $p_e$ is real-valued $\fz$ defines a self-adjoint element. Then, for $h\in K\cup\{e\}$, $g\in G$ and $\psi\in\ell^2(G)$, a short computation leads to
\begin{align*}
\left( \pi^\xi(\fz_h) \psi \right)(g) \; 
	&= \; \sum\limits_{\tilde{h}\in G} 
		\delta_{h^{-1}} \left(g^{-1}\tilde{h}\right) \cdot p_h\big(\alpha_{g^{-1}}(\xi)\big) \cdot \psi(\tilde{h})\\
		&= \; p_h\big(\alpha_{g^{-1}}(\xi)\big) \cdot \psi\left( gh^{-1} \right)\, ,\\
\left( \pi^\xi(\fz_h^\ast) \psi \right)(g) \; 
	&= \; \sum\limits_{\tilde{h}\in G} 
		\delta_{h^{-1}} \left(\tilde{h}^{-1} g\right) \cdot \overline{p_h\big(\alpha_{\tilde{h}^{-1}}(\xi)\big)} \cdot \psi(\tilde{h})\\
	&= \; \overline{p_h\big(\alpha_{(gh)^{-1}}(\xi)\big)} \cdot \psi(gh) \, .
\end{align*}
Using the linearity of the left-regular representation $\pi^\xi$, proven in Proposition~\ref{Chap2-Prop-LeftRegulRepres}, the identity $\pi^\xi(\fz)=H_\xi$ follows for $\xi\in\Xi$.
\end{proof}

\begin{remark}
\label{Chap2-Rem-PESchrOpFinRang}
Clearly, the Schr\"odinger operators and Jacobi operators over the group $G=\ZM$, which are typically studied in the literature, are specific pattern equivariant Schr\"odinger operators, c.f. Section~\ref{Chap2-Sect-ExampleAsZMSchrOp}.
\end{remark}

It is shown in Proposition~\ref{Chap4-Prop-MonotonConvSpectr} below that the spectrum of a generalized Schr\"odinger operator associated with a subshift of finite type is of particular interest. Under certain assumptions, the spectrum is determined by the spectra of the periodic elements of this subshift which is proven in the following proposition. The proof is mainly based on the strong convergence of the generalized Schr\"odinger operators, c.f. Proposition~\ref{Chap2-Prop-CovFamOp-Spect}~(c). 

\medskip

Recall the notion of patterns with support $[K]$ for a compact $K\subseteq G$ of a discrete, countable group $G$, c.f. Definition~\ref{Chap2-Def-Pattern}. Then 
$$
\Xi(K,U) \; 
	:= \; \big\{ 
		\xi\in\as^G \;\big|\; \ws(\xi)\cap\as^{[K]}\subseteq U
	\big\} \;
	\in \SG\big(\as^G\big)
$$
is a {\em subshift of finite type}. Its strongly periodic elements are denoted by $\Xi_{per}\subseteq\Xi(K,U)$, i.e., $\Xi_{per}$ is the set of all $\xi\in\Xi(K,U)$ such that $\Orb(\xi)$ is finite, c.f. Proposition~\ref{Chap2-Prop-AperSubs}. According to Corollary~\ref{Chap2-Cor-DictOrbitSubshift}, the inclusion $\overline{\Orb(\xi)}\subseteq\Xi(K,U)$ holds for every $\xi\in\Xi(K,U)$. A subshift of finite type $\Xi(K,U)$ has a {\em dense subset of strongly periodic elements} if $\Xi_{per}\subseteq\Xi$ is dense.

\begin{proposition}
\label{Chap2-Prop-SubsFinTypPerSpectr}
Let $\as$ be an alphabet and $G$ be a discrete, countable group. Consider a compact set $K\in\ks(G)$ and a subset $U\subseteq\as^{[K]}$ of pattern with support $[K]$ satisfying that the associated subshift of finite type $\Xi:=\Xi(K,U)$ has a dense set of strongly periodic elements. Then, for every generalized Schr\"odinger operator $H_\Xi:=\big(H_\xi\big)_{\xi\in\Xi}$, the equality
$$
\sigma\big( H_\Xi\big) \;
	= \; \overline{\bigcup_{\xi\in\Xi_{per}} \sigma\big(H_\xi\big)}
$$
holds.
\end{proposition}

\begin{proof}
Let $K\in\ks(G)$ and $U\subseteq\as^{[K]}$ be such that $\Xi_{per}\subseteq\Xi$ is dense in the associated subshift of finite type $\Xi$. Thus, for each $\xi\in\Xi$, there exists a sequence $\xi_n\in\Xi_{per}\,,\; n\in\NM\,,$ tending to $\xi$. Consider a generalized Schr\"odinger operator $H_\Xi$ defined in Definition~\ref{Chap2-Def-SchrOp-l2(G)} which is self-adjoint by definition. According to Proposition~\ref{Chap2-Prop-CovFamOp-Spect}~(c), the operator sequence $H_{\xi_n}:\ell^2(G)\to\ell^2(G)$ converges strongly to $H_\xi:\ell^2(G)\to\ell^2(G)$. Consequently, the inclusion
$$
\sigma\big(H_\xi\big) \;
	\subseteq\; 
	\bigcap_{n=1}^\infty\overline{\left(
			\bigcup_{m=n}^{\infty}\sigma\big(H_{\xi_n}\big)
		\right)}
$$
is derived by the strong convergence, c.f. Section~\ref{Chap3-Sect-OpTop}. The set on the right hand side is contained in $\overline{\bigcup_{\xi\in\Xi_{per}} \sigma\big(H_\xi\big)}$ by construction. Hence, the inclusion
$$
\sigma\big( H_\Xi\big) \;
	= \; \overline{\bigcup_{\xi\in\Xi} \sigma\big( H_\xi \big)} \;
	\overset{\text{Prop.~\ref{Chap2-Prop-CovFamOp-Spect}}}{=} \; 
		 \bigcup_{\xi\in\Xi} \sigma\big( H_\xi \big)
	\subseteq \; \overline{\bigcup_{\xi\in\Xi_{per}} \sigma\big(H_\xi\big)}
$$
follows. The converse inclusion is deduced by $\Xi_{per}\subseteq\Xi$.
\end{proof}

\begin{remark}
\label{Chap2-Rem-SubsFinTypPerSpectr}
(i) The set $\bigcap_{n=1}^\infty\overline{\left( \bigcup_{m=n}^{\infty}\sigma\big(H_{\Xi_n}\big) \right)}$ is often denoted by $\lim_{n\to\infty}\sigma\big( H_{\xi_n}\big)$ in the literature. This notation is misleading since it does not agree with the limit of the sets $\big(\sigma\big( H_{\xi_n}\big)\big)_{n\in\NM}$ with respect to the Vietoris-topology on the closed subsets $\ks(\RM)$ of $\RM$ (respectively the convergence with respect to the Hausdorff metric on $\ks(\RM)$). More precisely, this sequence of sets does not need to converge in the Vietoris-topology in general. In general, this set is only the limit superior $\limsup_{n\to\infty}\sigma(H_{\xi_n})$.

\vspace{.1cm}

(ii) Let $\Xi(K,U)\in\SG\big(\as^G\big)$ be a subshift of finite type arising $K\in\ks(G)$ and $U:=\ws(\Xi')\cap\as^{[K]}$ for a periodically approximable subshift $\Xi'\in\SG\big(\as^G\big)$. If the group $G$ is residually finite and $\Xi(K,U)$ is strongly irreducible, then $\Xi_{per}\subseteq\Xi(K,U)$ is dense. This is deduced by \cite[Theorem~1.1]{SiCo12} and the fact that $\Xi(K,U)$ contains at least one periodic configuration as $\Xi'$ is periodically approximable. If $G=\ZM^d$, the denseness follows if $\Xi(K,U)$ is finite-orbit square-mixing, i.e., it is square mixing and contains an element with finite orbit (strongly periodic element), c.f. \cite[Lemma~5.4]{Lig03}. The existence of an element with finite orbit in $\Xi(K,U)$ is derived from the fact that $\Xi'$ is periodically approximable and $\Xi(K,U)$ is defined via $\Xi'$, c.f. Proposition~\ref{Chap2-Prop-AperSubs}.

\vspace{.1cm}

(iii) Proposition~\ref{Chap2-Prop-SubsFinTypPerSpectr} does not imply that the subshift $\overline{\Orb(\xi)}\subseteq\Xi(K,U)$ is periodically approximable for every $\xi\in\Xi(K,U)$, c.f. Example~\ref{Chap2-Ex-OnlyLowerSemiContSpectr}. It only states that there exists a sequence $(\Xi_n)_{n\in\NM}$ of strongly periodic subshifts so that $\sigma(H_\Xi)\subseteq\limsup_{n\to\infty}\sigma\big(H_{\Xi_n}\big) = \bigcap_{n=1}^\infty\overline{\big(\bigcup_{m=n}^{\infty}\sigma\big(H_{\Xi_n}\big))}$. 
\end{remark}

\begin{example}
\label{Chap2-Ex-OnlyLowerSemiContSpectr}
Consider the two-sided infinite word $\xi = \ldots aaab|abbb\ldots \in\as^{\ZM}$ defined in Example~\ref{Chap5-Ex-deBruijnEventuallyNotStrongConn} over the alphabet $\as:=\{a,b\}$. Consider the subshift of finite type $\Xi(K,U):=\Xi(\{0,1\},\{ab,ba,aa,bb\})$. Thus, $\eta\in\as^\ZM$ is contained in this subshift of finite type if and only if all subwords of $\eta$ up to length $2$ are elements of $\{ab,ba,aa,bb\}$. It exists a periodic element in $\Xi(K,U)$. For instance, $\eta(j):=a\,,\; j\in\ZM\,,$ is an element of $\Xi(K,U)$. Thus, the set of periodic elements of $\Xi(K,U)$ is dense in $\Xi(K,U)$ by \cite[Theorem~1.1]{SiCo12} as $\Xi(K,U)$ is strongly irreducible. Furthermore, $\xi$ is an element of $\Xi(K,U)$. On the other hand, the associated subshift $\overline{\Orb(\xi)}$ is not periodically approximable by Theorem~\ref{Chap5-Theo-ExPerAppr} since the de Bruijn graph $\gs_2$ of order $2$ associated with $\xi$ is not strongly connected, c.f. Example~\ref{Chap5-Ex-deBruijnEventuallyNotStrongConn}. On the other hand, $\eta_n:=\big(\xi|_{[-n,-1]}\big|\,\xi|_{[0,n]}\big)^\infty\,,\; n\in\NM\,,$ are periodic elements of $\as^\ZM$ such that $\lim_{n\to\infty}\eta_n=\xi$. Then the inclusion $\sigma(H_\xi)\subseteq\limsup_{n\to\infty}\sigma(H_{\eta_n})$ follows for any generalized Schr\"odinger operator $H$ by the strong convergence of the operators, c.f. Proposition~\ref{Chap2-Prop-CovFamOp-Spect}.
\end{example}

\cleardoublepage

\chapter{A tool to prove the continuous behavior of the spectra and its application to generalized Schr\"odinger operators}
\chaptermark{Continuous behavior of the spectra: tool and Schr\"odinger operators}
\label{Chap4-ToolContBehavSpectr}
\stepcounter{section}
\setcounter{section}{0}

The characterization of the continuous behavior of the spectra of generalized (pattern equivariant) Schr\"odinger operators by the continuous variation of the underlying dynami\-cal systems is proven in this chapter, c.f. Theorem~\ref{Chap4-Theo-CharDynSystConvSpectr} and Theorem~\ref{Chap4-Theo-CharSubshiftConvSpectr}. This new approach is based on joint work with {\sc J. Bellissard} and {\sc G. de Nittis} \cite{BeBeNi16}. The proof relies on the investigation of the universal dynamical system associated with a dynamical system $(X,G,\alpha)$. With this at hand, periodic approximations of these operators are defined by defining periodic approximations of the underlying dynamical systems so that the spectra of the operators converge in the Hausdorff metric and the operators converge in the strong operator topology, c.f. Theorem~\ref{Chap4-Theo-PeriodicApproximations}.

\medskip

\begin{figure}[htb]
\centering
\includegraphics[scale=0.93]{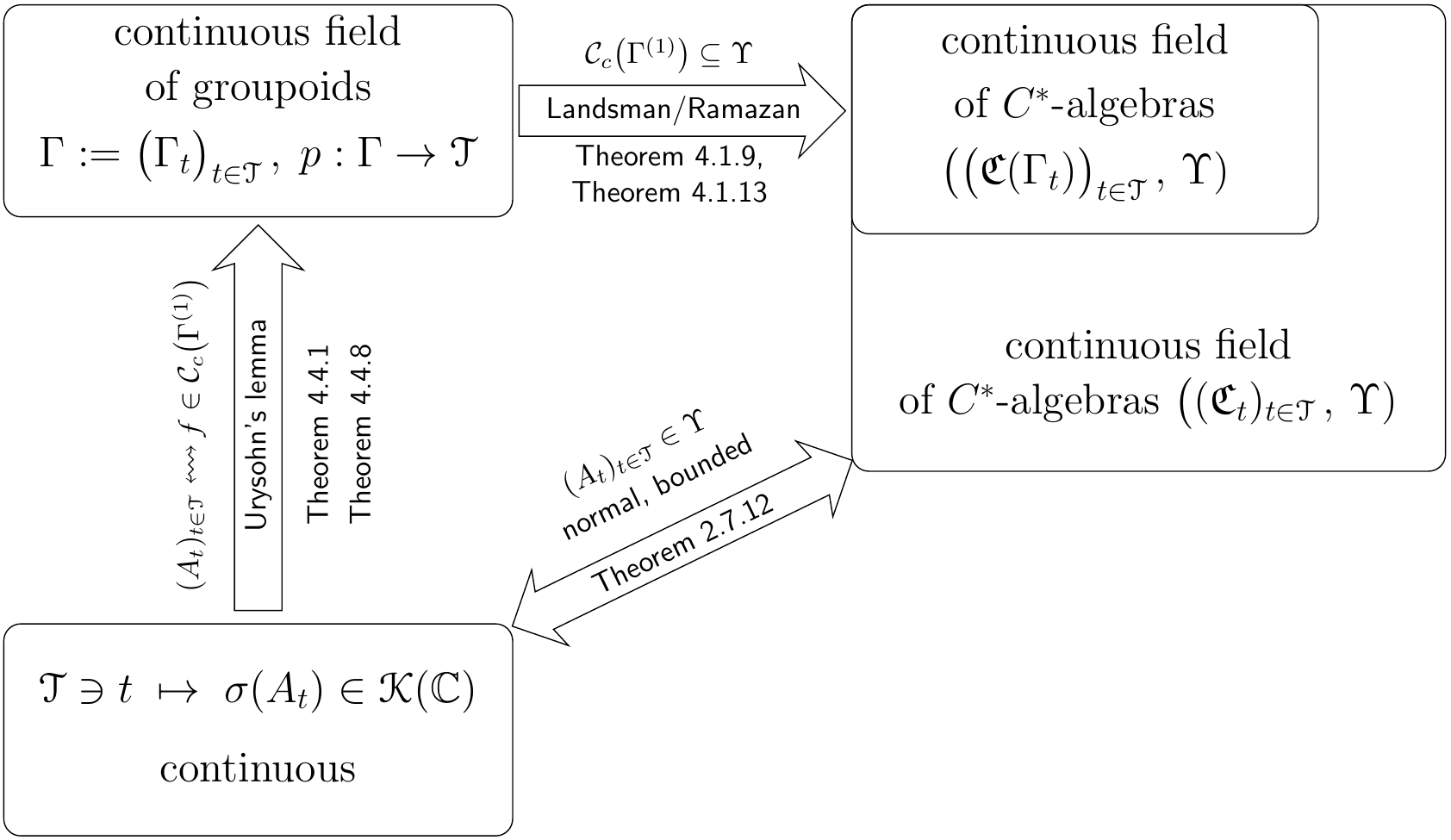}
\caption{The connections of the different concepts.}
\label{Chap4-Fig-Concept}
\end{figure}

The chapter is organized as follows. The notion of continuous fields of groupoids and its connection to continuous fields of $C^\ast$-algebras goes back to {\sc Landsman} and {\sc Ramazan} \cite{LaRa99}. In combination with Theorem~\ref{Chap3-Theo-PContEquivContSpectNormal}, this yields a reasonable tool to prove the continuous behavior of the spectra also in different situations than Schr\"odinger operators. A detailed elaboration of this concept is provided in Section~\ref{Chap4-Sect-ContFieldGroupCAlg}. Section~\ref{Chap4-Sect-GroupoidIsomorphism} provides the standard notion of measurable preserving groupoid isomorphisms and their connections to isomorphisms of $C^\ast$-algebras that are used in the following. Then the new concept of universal dynamical system $(\Xun,G,\alpha)$ associated with $(X,G,\alpha)$ are established in Section~\ref{Chap4-Sect-UnivGroupDynSyst}. There it is also proven that the universal groupoid $\Gun:=\Xun\rtimes_\alpha G$ naturally defines a continuous field of groupoids over the space $\SG(X)$ of dynamical subsystems of $(X,G,\alpha)$. With this at hand it is proven that, for all generalized Schr\"odinger opera\-tors, the spectra converge in $\ks(\CM)$ if and only if the associated dynamical subsystems converge in $\SG(X)$, c.f. Section~\ref{Chap4-Sect-ContSpecGenSchrOp}. One implication follows by using the concept of universal dynamical system and the result of \cite{LaRa99} while the converse is based on the fact that the operators are represented as continuous functions on the transformation group groupoid, c.f. Figure~\ref{Chap4-Fig-Concept}. In case of symbolic dynamical systems, this result is reduced to specific pattern equivariant Schr\"odinger operators. Finally, this new approach delivers a tool to approximate Schr\"odinger operators by periodic ones, c.f. Section~\ref{Chap4-Sect-ApprPerAppr}. While in the case of dynamical systems only the existence of a periodic approximation is verified, the approximations are given explicitly for pattern equivariant Schr\"odinger operators on symbolic dynamical systems.

\section{Continuous fields of groupoids and \texorpdfstring{$C^\ast$-algebras}{C-algebras}}
\label{Chap4-Sect-ContFieldGroupCAlg}

This section provides the link between continuous fields of $C^\ast$-algebras and continuous fields of groupoids which goes back to {\sc Landsman} and {\sc Ramazan} \cite{LaRa99}. Their result is a generalization of the group case provided by {\sc Rieffel} \cite{Rie89} and they use some general concepts of $C^\ast$-algebras proven by {\sc Blanchard} \cite{Bla96}. 

\medskip

Recall that a map between topological spaces is called {\em open} whenever the image of open sets is open. The following definition was first given in \cite{LaRa99}.

\begin{definition}[Continuous fields of groupoids, \cite{LaRa99}]
\label{Chap4-Def-ContFieldGroupoid}
Let $\Gamma$ be a topological groupoid. A triple $(\Gamma,\ts,p)$ is called {\em continuous field of groupoids} if $\ts$ is a topological Hausdorff space and $p: \Gamma^{(1)} \to \ts$ is a continuous, open and surjective map satisfying $p = p|_{\Gamma^{(0)}} \circ r = p|_{\Gamma^{(0)}} \circ s$.
\end{definition} 

This definition entails more constraints on the topology of the space $\ts$.

\begin{lemma}
\label{Chap4-Lem-TNormalContFieldGroup}
Let $(\Gamma,\ts,p)$ be a continuous field of groupoids. Then $\ts$ is a second countable, locally compact, Hausdorff space, i.e., $\ts$ is a normal space. If, additionally, the unit space $\Gamma^{(0)}$ is compact, then $\ts$ is a compact space
\end{lemma}

\begin{proof}
By definition of a continuous field of groupoids, $\ts$ is Hausdorff. Proposition~\ref{App1-Prop-ContOpenSurjMapPresSecCountable} implies that $\ts$ is also second countable and locally compact since $p:\Gamma^{(1)}\to\ts$ is a continuous, open, surjective map and $\Gamma^{(1)}$ is a second countable, locally compact space.

\vspace{.1cm}

Suppose $\Gamma^{(0)}$ is compact. The map $p$ is surjective by definition. Thus, the image $p(\Gamma^{(0)})$ is equal to $\ts$ since, additionally, the equations $p = p|_{\Gamma^{(0)}} \circ r$ and $\Gamma^{(0)}=r(\Gamma^{(1)})$ hold. Hence, the continuity of $p$ implies that the space $\ts$ is compact as it is the continuous image of a compact set. 
\end{proof}

\medskip

The name continuous field of groupoids for a triple $(\Gamma,\ts,p)$ is due to the fact that $\Gamma$ fibers over the space $\ts$ with respect to the map $p$ which is presented in the following lemma.

\begin{lemma}
\label{Chap4-Lem-FiberClosGroup}
Let $(\Gamma,\ts,p)$ be a continuous field of groupoids. For $t\in\ts$, the preimage $\Gamma_t^{(1)}:=p^{-1}(\{t\})$ is a closed, subgroupoid of $\Gamma^{(1)}$. More precisely, $\Gamma_t$ is a topological groupoid such that the $r$-fibers satisfy $\Gamma^x=\Gamma^x_t$ for $x\in\Gamma^{(0)}_t$. If, additionally, $\Gamma^{(0)}$ is compact, then $\Gamma_t$ has also compact unit space $\Gamma_t^{(0)}$ for each $t\in\ts$.
\end{lemma}

\begin{proof}
The proof is organized as follows.
\begin{itemize}
\item[(i)] The set $\Gamma_t$ is a topological groupoid with $\Gamma^x=\Gamma^x_t$ for $x\in\Gamma^{(0)}_t$.
\item[(ii)] The unit space $\Gamma^{(0)}_t$ is compact for each $t\in\ts$ if $\Gamma^{(0)}$ is compact.
\end{itemize}

(i): Since $\ts$ is a Hausdorff space, a singleton $\{t\}$ for $t\in\ts$ is closed. Then the continuity of $p$ implies that $\Gamma_t$ is a closed subset of $\Gamma^{(1)}$. Let $\gamma\in\Gamma^{(1)}$. Due to $p = p|_{\Gamma^{(0)}} \circ r = p|_{\Gamma^{(0)}} \circ s$, the relation $\gamma\in\Gamma^{(1)}_t$ holds if and only if $p(r(\gamma))=t$ if and only if $p(s(\gamma))=t$. Let $\gamma,\eta\in\Gamma_t$ be such that $s(\gamma)=r(\eta)$. According to Proposition~\ref{Chap2-Prop-BasGroupoid}~(b) and (d), the identities $p(r(\gamma^{-1}))=p(s(\gamma))=t$ and $p(r(\gamma\eta))=p(r(\gamma))=t$ are concluded. Thus, the inverse $\gamma^{-1}$ and the composition $\gamma\eta$ are contained in $\Gamma^{(1)}_t$, namely $\Gamma_t$ is a subgroupoid of $\Gamma$. Since $p = p|_{\Gamma^{(0)}} \circ r$ holds, the $r$-fiber $\Gamma^{x}$ in $\Gamma$ is equal to the $r$-fiber $\Gamma^{x}_t$ in $\Gamma_t$ for $x\in\Gamma_t^{(0)}=\Gamma^{(0)}\cap\Gamma^{(1)}_t$.

\vspace{.1cm}

(ii): The unit space $\Gamma^{(0)}_t$ is equal to the intersection $\Gamma^{(0)}\cap\Gamma^{(1)}_t$. Thus, the groupoid $\Gamma_t$ has a compact unit space since $\Gamma^{(0)}$ is compact and $\Gamma^{(1)}_t$ closed.
\end{proof}

\medskip

For a continuous field of groupoids $(\Gamma,\ts,p)$, the groupoid $\Gamma$ fibers with respect to the map $p$ and the topological space $\ts$, c.f. Lemma~\ref{Chap4-Lem-FiberClosGroup}. Then, for further applications, the exist\-ence of extensions onto $\Gamma^{(1)}$ of continuous maps $\fz:F\to\CM$ supported on closed subsets $F\subseteq\Gamma^{(1)}$ is useful. This existence is a direct consequence of the fact that a topological groupoid is a normal space and Tietze's Extension Theorem, c.f. Theorem~\ref{App1-Theo-TietzesThm}. 

\begin{lemma}
\label{Chap4-Lem-CompContExtGroupoid}
Let $\Gamma$ be a topological groupoid. Then, for every closed subset $F\subseteq\Gamma^{(1)}$ and each continuous function $\fz:F\to\CM$, there exists a continuous map $\gz:\Gamma^{(1)}\to\CM$ such that $\gz|_F=\fz$. If, additionally, $F$ is compact such a $\gz$ exists in $\Cc_c\big(\Gamma^{(1)}\big)$.
\end{lemma}

\begin{proof}
The topological space $\Gamma^{(1)}$ is locally compact, Hausdorff and second countable. Thus, $\Gamma^{(1)}$ is a normal space, c.f. Proposition~\ref{App1-Prop-LocCompHausSecCountImplNormal}. Since $\CM$ is homeomorphic to $\RM^2$, Tietze's Extension Theorem applies to the map $\fz:F\to\CM$, c.f. \cite[Satz I.8.5.4]{Schubert75}. Precisely, there exists a continuous function $\gz:\Gamma^{(1)}\to\CM$ such that $\gz(\gamma)=\fz(\gamma)$ for all $\gamma\in F$. As $\fz$ has compact support, the map $\gz$ can be chosen with compact support as well, c.f. Theorem~\ref{App1-Theo-TietzesThm}.
\end{proof}

\medskip

A subset of $\Delta\subseteq\Gamma^{(1)}$ is called a {\em subgroupoid} if it is closed under inversion and composition, i.e., $\gamma,\eta\in\Delta$ with $s(\gamma)=r(\eta)$ imply $\gamma^{-1},\gamma\eta\in\Delta$. A subgroupoid $\Delta$ itself defines a topological groupoid with the induced topology, composition and inverse by $\Gamma$.

\begin{lemma}
\label{Chap4-Lem-FiberClosGroupHaarSyst}
Let $(\Gamma,\ts,p)$ be a continuous field of groupoids and $(\mu^x)_{x\in\Gamma^{(0)}}$ be a left-continuous Haar system of $\Gamma$. For $t\in\ts$, $(\mu^x)_{x\in\Gamma^{(0)}_t}$ defines a left-continuous Haar system of $\Gamma_t$.
\end{lemma}

\begin{proof}
Let $t\in\ts$ and $x\in\Gamma^{(0)}_t$. According to Definition~\ref{Chap2-Def-HaarSystem}, $(\mu^x)_{x\in\Gamma^{(0)}_t}$ is a left-continuous Haar system of $\Gamma_t$ if it satisfies \nameref{(H1)}-\nameref{(H3)}. These conditions are proven in the following.

\vspace{.1cm}

Since $\Gamma^x=\Gamma^x_t$, the support $\supp(\mu^x)$ is equal to $\Gamma^x_t$. Thus, \nameref{(H1)} follows for the measures  $(\mu^x)_{x\in\Gamma^{(0)}_t}$ on the groupoid $\Gamma_t$.

\vspace{.1cm}

Let $\fz\in\Cc_c\big(\Gamma^{(1)}_t\big)$. Since $\Gamma^{(1)}_t$ is a closed subset of $\Gamma^{(1)}$, the support of $\fz$ is a closed subset of $\Gamma^{(1)}$. According to Lemma~\ref{Chap4-Lem-CompContExtGroupoid}, there exists a $\gz\in\Cc_c\big(\Gamma^{(1)}\big)$ such that $\gz|_{\Gamma^{(1)}_t}=\fz$. 

\vspace{.1cm}

The equation $\mu^x(\fz)=\mu^x(\gz)$ holds for $x\in\Gamma^{(0)}_t$ since the measure $\mu^x$ is supported on $\Gamma^{x}=\Gamma^{x}_t$. Thus, \nameref{(H2)} is derived by using that $(\mu^x)_{x\in\Gamma^{(0)}}$ is a left-continuous Haar system of $\Gamma$.

\vspace{.1cm}

For $\gamma\in\Gamma^{(1)}$, the equation
$$
\int_{\Gamma^{s(\gamma)}} 
	\gz(\gamma\cdot\varrho)
\; d\mu^{s(\gamma)}(\varrho) \; 
		= \; \int_{\Gamma^{r(\gamma)}} 
			\gz(\varrho)
		\; d\mu^{r(\gamma)}(\varrho)
$$
is concluded as \nameref{(H3)} is valid for $\Gamma$ with Haar system $(\mu^x)_{x\in\Gamma^{(0)}}$. If $\gamma\in\Gamma^{(1)}_t$, the measure $\mu^{s(\gamma)}$ is supported on $\Gamma^{s(\gamma)}_t$ and $\mu^{r(\gamma)}$ is supported on $\Gamma^{r(\gamma)}_t$. Hence, for $\gamma\in\Gamma^{(1)}_t$, the identities
\begin{align*}
\int_{\Gamma^{s(\gamma)}} 
	\fz(\gamma\cdot\varrho)
\; d\mu^{s(\gamma)}(\varrho) \;
	&= \; \int_{\Gamma^{s(\gamma)}} 
		\gz(\gamma\cdot\varrho)
	\; d\mu^{s(\gamma)}(\varrho) \, , \\
\; \int_{\Gamma^{r(\gamma)}} 
	\fz(\varrho)
\; d\mu^{r(\gamma)}(\varrho)		
	&= \; \int_{\Gamma^{r(\gamma)}} 
		\gz(\varrho)
	\; d\mu^{r(\gamma)}(\varrho)
\end{align*}
are deduced since $\gz|_{\Gamma^{(1)}_t}=\fz$. Consequently, \nameref{(H3)} holds for the family of Borel measures $(\mu^x)_{x\in\Gamma^{(0)}_t}$ on the groupoid $\Gamma_t$. 
\end{proof}

\begin{definition}[Fiber groupoid]
\label{Chap4-Def-FiberFieldGroupoid}
Let $(\Gamma,\ts,p)$ be a continuous field of groupoids and $(\mu^x)_{x\in\Gamma^{(0)}}$ is a left-continuous Haar system of $\Gamma$. Then the closed subgroupoid $\Gamma_t:=p^{-1}(\{t\})$ is called {\em fiber groupoid} for $t\in\ts$ with induced topology, composition, inversion and left-continuous Haar system $(\mu^x)_{x\in\Gamma^{(0)}_t}$.
\end{definition}

Recall that a topological groupoid is called \'etale if the range map of the groupoid is a local homeomorphism, c.f. Definition~\ref{Chap2-Def-Etale}.

\begin{lemma}
\label{Chap4-Lem-FiberEtale}
Let $(\Gamma,\ts,p)$ be a continuous field of groupoids such that $\Gamma$ is \'etale. Then the fiber groupoid $\Gamma_t$ is \'etale for every $t\in\ts$. 
\end{lemma}

\begin{proof}
Let $r:\Gamma^{(1)}\to\Gamma^{(1)}$ be the range map of $\Gamma$ and $t\in\ts$. For $\gamma\in\Gamma^{(1)}_t\subseteq\Gamma^{(1)}$, there exists an open neighborhood $\tilde{U}$ in $\Gamma^{(1)}$ such that $r|_{\tilde{U}}$ is a homeomorphism. The space $\Gamma^{(1)}_t$ is equipped with the subspace topology, i.e., $V\subseteq\Gamma^{(1)}_t$ is open if and only if there is an open set $U\subseteq\Gamma^{(1)}$ such that $V=U\cap\Gamma^{(1)}_t$. Consequently, the restriction of the range $r|_{\tilde{U}\cap\Gamma_t^{(1)}}$ is a homeomorphism on $\Gamma^{(1)}_t$. Thus, $V:=\tilde{U}\cap\Gamma^{(1)}_t$ is an open neighborhood of $\gamma$ in the fiber groupoid $\Gamma^{(1)}_t$ and the range map $r$ of $\Gamma^{(1)}_t$ restricted to $V$ is a homeomorphism. Hence, $\Gamma^{(1)}_t$ is \'etale.
\end{proof}

\medskip

Recall that a groupoid with left-continuous Haar system induces the reduced and full $C^\ast$-algebra introduced in Section~\ref{Chap2-Sect-GroupoidCalgebras}. These algebras were constructed by the compactly supported, continuous functions $\Cc_c\big(\Gamma^{(1)}\big)$ of the groupoid $\Gamma$. According to Proposition~\ref{Chap2-Prop-Red=FullCalgebra}, these $C^\ast$-algebras agree whenever the groupoid is topologically amenable (Definition~\ref{Chap2-Def-AmenableGroupoid}). This is crucial for the further considerations. The property of topologically amenable transfers to subgroupoids, c.f. \cite[Proposition~5.1.1]{AnantharamanRenault00}.

\begin{proposition}[\cite{AnantharamanRenault00}]
\label{Chap4-Prop-ContFieldGroupoidAmenable}
Let $(\Gamma,\ts,p)$ be a continuous field of locally compact group\-oids with continuous Haar system $(\mu^x)_{x\in\Gamma^{(0)}}$. If $\Gamma$ is topologically amenable, then the fiber groupoids $\Gamma_t,\; t\in\ts,$ are topologically amenable.
\end{proposition}

\begin{proof}
According to \cite[Proposition~5.1.1]{AnantharamanRenault00}, a closed subgroupoid of a groupoid $\Gamma$ is topologically amenable if $\Gamma$ is so. Lemma~\ref{Chap4-Lem-FiberEtale} asserts that $\Gamma_t$ is a closed subgroupoid of $\Gamma$ for $t\in\ts$. Thus, $\Gamma_t$ is topologically amenable for $t\in\ts$.
\end{proof}

\medskip

Let $(\Gamma,\ts,p)$ be a continuous field of groupoids and $(\mu^x)_{x\in\Gamma^{(0)}}$ a left-continuous Haar system of $\Gamma$. Then $\CG^\ast_{red}(\Gamma_t)$ is the reduced $C^\ast$-algebra of the topological groupoid $\Gamma_t$ with left-continuous Haar-system $\mu^x,\; x\in\Gamma^{(0)}_t\,,\;$ for $t\in\ts$. For $\fz\in\Cc_c\big(\Gamma^{(1)}\big)$, the restriction $\fz_t:=\fz|_{\Gamma_t}$ of $\fz$ to the closed subset $\Gamma_t$ is an element of $\Cc_c\big(\Gamma_t^{(1)}\big)$. Hence, $\fz_t\in\CG^\ast_{red}(\Gamma_t)$ and $\fz_t\in\CG^\ast_{full}(\Gamma_t)$ follow. Then $\|\fz_t\|_{red}$ is the reduced-norm of $\fz_t$ in the $C^\ast$-algebra $\CG^\ast_{red}(\Gamma_t)$ and $\|\fz_t\|_{full}$ is the full-norm of $\fz_t$ in the $C^\ast$-algebra $\CG^\ast_{full}(\Gamma_t)$. Altogether, an element $\fz\in\Cc_c\big(\Gamma^{(1)}\big)$ is represented by the vector field $(\fz_t)_{t\in\ts}\in\prod_{t\in\ts}\Cc_c\big(\Gamma^{(1)}_t\big)$. As long as there is no confusion in the notation, the representation of $\Cc_c\big(\Gamma^{(1)}\big)$ as a subset of $\prod_{t\in\ts} \Cc_c\big(\Gamma^{(1)}_t\big)$ is also denoted by $\Cc_c\big(\Gamma^{(1)}\big)$.

\medskip

The following theorem, proved in \cite{LaRa99}, provides the key information to obtain a continuous field of unital $C^\ast$-algebras from a continuous field of groupoids. It is the generalization of \cite[Theorem 2.5]{Rie89} from the category of locally compact groups to the category of locally compact groupoids. The crucial ingredient for both results is the proof of the upper semi-continuity of the norms. For groupoids, this is proved with a general argument of $C^\ast$-algebras provided by {\sc Blanchard} \cite{Bla96}.

\medskip

In order to formulate the result, recall the notion of semi-continuity. Let $\ts$ be a topological spaces then a map $\Phi:\ts\to\RM$ is {\em upper (respectively, lower) semi-continuous at $t_0\in\ts$} if, for every $\varepsilon>0$, there exists an open neighborhood $U_0$ of $t_0$ such that $\Phi(t)\leq \Phi(t_0)+\varepsilon$ (respectively, $\Phi(t)\geq \Phi(t_0)-\varepsilon$) for all $t\in U_0$.

\begin{theorem}[{\cite[Theorem 5.5]{LaRa99}}]
\label{Chap4-Theo-LanRamContFieldGroupoid}
Let $(\Gamma,\ts,p)$ be a continuous field of groupoids and $(\mu^x)_{x\in\Gamma^{(0)}}$ is a left-continuous Haar system of $\Gamma$. For every $\fz\in\Cc_c\big(\Gamma^{(1)}\big)$, let $\fz_t\in\Cc_c\big(\Gamma^{(1)}_t\big)$ be the restriction of $\fz$ to the fiber groupoid $\Gamma_t$ for $t\in\ts$. Then the following assertions hold.
\begin{itemize}
\item[(a)] The map $\ts\ni t\mapsto\|\fz_t\|_{full}\in[0,\infty)$ is upper semi-continuous.
\item[(b)] The map $\ts\ni t\mapsto\|\fz_t\|_{red}\in[0,\infty)$ is lower semi-continuous.
\end{itemize}
\end{theorem}

\begin{remark}
\label{Chap4-Rem-LanRamContFieldGroupoid}
In order to show the upper semi-continuity of the full-norm, {\sc Landsman} and {\sc Ramazan} use that, for a compact subset $K\subseteq \ts$ and a $t\in\ts\setminus K$, there exists an open neighborhood $V$ of $K$ such that $t\not\in V$ and there is an $f\in\Cc_c(\ts)$ with $\supp(f)\subseteq V$ and $f|_K\equiv 1$. This follows from the Lemma of Urysohn, c.f. Proposition~\ref{App1-Prop-LemmaUrysohn}. The Lemma of Urysohn is applicable as $\ts$ is automatically a normal space even if in the Definition~\ref{Chap4-Def-ContFieldGroupoid} of continuous fields of groupoids it is only required that $\ts$ is Hausdorff, c.f. Lemma~\ref{Chap4-Lem-TNormalContFieldGroup} and Theorem~\ref{App1-Theo-TietzesThm}.
\end{remark}

For convenience of the reader, the main steps of the proof are presented here.

\medskip

\begin{proof} 
Let $\fz\in\Cc_c\big(\Gamma^{(1)}\big)$. For $t\in\ts$, the fiber groupoid $p^{-1}(\{t\})$ is denoted by $\Gamma_t$ with unit space $\Gamma^{(0)}_t:=p^{-1}(\{t\})\cap\Gamma^{(0)}$. The fiber groupoid $\Gamma_t$ is closed as a subset of $\Gamma$, c.f. Lemma~\ref{Chap4-Lem-FiberClosGroup}.

\vspace{.1cm}

(a): The vector space of continuous functions $\psi:\ts\to\CM$ that vanish at infinity is denoted by $\Cc_0(\ts)$. Then $\Cc_0(\ts)$ is a $C^\ast$-algebra equipped with the pointwise multiplication and the uniform norm $\|\psi\|_\infty:=\sup_{t\in\ts}|\psi(t)|$, c.f. Example~\ref{Chap3-Ex-ContFunctCalgebra}. For $t\in\ts$, define the subalgebra $\Cc_t(\ts):=\big\{\psi\in\Cc_0(\ts)\;|\; \psi(t)=0\big\}$. This subalgebra is a closed, two-sided ideal of $\Cc_0(\ts)$. Then $\Cc_t(\ts)\cdot\CG^\ast_{full}(\Gamma)$ defines also a two-sided ideal of the full $C^\ast$-algebra $\CG^\ast_{full}(\Gamma)$ according to \cite[Corollary~1.9]{Bla96}. Consider $\rho_t:\CG^\ast_{full}(\Gamma)\to\CG^\ast_{full}(\Gamma)\big/\big(\Cc_t(\ts)\cdot\CG^\ast_{full}(\Gamma)\big)$ the corresponding quotient map. Then the equality 
$$
\|\rho_t(\fz)\| \;
	= \; \inf\Big\{
		\big\|\big(1-\psi+\psi(t)\big)\cdot\fz\big\|\; \Big|\; \psi\in\Cc_0(\ts)
	\Big\}\,,
	\qquad
	t\in\ts\,,
$$
is deduced by \cite[Lemma~1.10]{Bla96}. The map $\ts\ni t\mapsto\big\|\big(1-\psi+\psi(t)\big)\cdot\fz\big\|\in\RM$ is continuous for each $\psi\in\Cc_0(\ts)$. Hence, $\ts\ni t\mapsto\|\rho_t(\fz)\|\in\RM$ is upper semi-continuous since it is the infimum over continuous functions. With this at hand, it suffices to prove the identity $\|\rho_t(\fz)\|=\|\fz_t\|$ for $t\in\ts$. This is done by showing the existence of a $\ast$-isomorphism between the $C^\ast$-algebras $\CG^\ast_{full}(\Gamma_t)$ and $\CG^\ast_{full}(\Gamma)\big/\big(\Cc_t(\ts)\cdot\CG^\ast_{full}(\Gamma)\big)$ since each $\ast$-isomorphism preserves the norms, see e.g. \cite[Theorem~3.1.4]{Murphy90}.

\vspace{.1cm} 

In order to do so, a short exact sequence of $C^\ast$-algebras is defined. Let $t\in\ts$. The set $U:=\Gamma^{(0)}\setminus\Gamma^{(0)}_t$ induces a subgroupoid  $\Gamma^{(1)}_U:= p^{-1}(U)$ of $\Gamma^{(1)}$ which is equal to $\Gamma^{(1)}\setminus\Gamma^{(1)}_t$. A function $\fz\in\Cc_c\big(\Gamma^{(1)}_U\big)$ is extended to a function $\widehat{\fz}\in\Cc_c\big(\Gamma^{(1)}\big)$ by setting $\widehat{\fz}$ equal to $0$ on $\Gamma^{(1)}_t$. This extends to an injective $\ast$-homomorphism $e:\CG^\ast_{full}(\Gamma_U)\to\CG^{\ast}_{full}(\Gamma)$ by taking the norm closure. Let $\imath:\Gamma_t\to\Gamma^{(1)}$ be the embedding of the groupoid $\Gamma_t$ for $t\in\ts$. This induces a $\ast$-homomorphism $\imath^\ast:\CG^\ast_{full}(\Gamma)\to\CG^\ast_{full}(\Gamma_t)$ defined via the norm closure of the map $\fz\in\Cc_c\big(\Gamma^{(1)}\big)\mapsto\fz\circ \imath\in\Cc_c\big(\Gamma^{(1)}_t\big)$. The notation $\imath^\ast$ is chosen as it is the pullback of the embedding $\imath$. According to \cite{Tor85,HiSk87} (see also \cite[Proposition~5.1]{LaRa99}),  
$$
0\longrightarrow \CG^\ast_{full}(\Gamma_U)\overset{e}{\longrightarrow}\CG^\ast_{full}(\Gamma)\overset{\imath^\ast}{\longrightarrow}\CG^\ast_{full}(\Gamma_t)\longrightarrow 0
$$
is a short exact sequence, i.e., the image of $e$ is equal to the kernel of the map $\imath^\ast$. Thus, the exactness implies that
$\CG^\ast_{full}(\Gamma)\big/\big(\Im(e)\big)$ is isomorphic to $\CG^\ast_{full}(\Gamma_t)$ where $\Im(e)$ denotes the image of the $\ast$-homomorphism $e$.

\vspace{.1cm}

In view of that it suffices to prove the equality $e\big(\CG^\ast_{full}(\Gamma_U^{(1)})\big)=\Cc_t(\ts)\cdot \CG^\ast_{full}(\Gamma^{(1)})$ since this leads to $\|\rho_t(\fz)\|=\|\fz_t\|$ for $t\in\ts$ by the previous considerations. Taking Remark~\ref{Chap4-Rem-LanRamContFieldGroupoid} into account, the equality $e\big(\Cc_c(\Gamma_U^{(1)})\big)=\Cc_t(\ts)\cdot \Cc_c(\Gamma^{(1)})$ is derived. Consequently, the desired identity $e\big(\CG^\ast_{full}(\Gamma_U^{(1)})\big)=\Cc_t(\ts)\cdot \CG^\ast_{full}(\Gamma^{(1)})$ follows by taking the norm closure.

\vspace{.1cm}

(b): Denote by $\gz|_{\Gamma^x}$ the restriction of the function $\gz\in\Cc_c\big(\Gamma^{(1)}\big)$ to the $r$-fiber $\Gamma^x$. Then $\gz|_{\Gamma^x}$ is an element of the Hilbert space $L^2(\Gamma^x,\mu^x)$. Furthermore, each continuous, compactly supported function $L^2(\Gamma^x,\mu^x)$ can be extended to a $\gz\in\Cc_c\big(\Gamma^{(1)}\big)$ by Lemma~\ref{Chap4-Lem-CompContExtGroupoid} since $\Gamma^x:=r^{-1}(\{x\})$ is a closed subset of $\Gamma^{(1)}$. Consequently, the set $\big\{\gz|_{\Gamma^x}\;\big|\; \gz\in\Cc_c\big(\Gamma^{(1)}\big)\big\}$ is equal to $\Cc_c\big(\Gamma^x\big)$ implying that it is a dense subset of $L^2(\Gamma^x,\mu^x)$ with respect to the $L^2$-norm. Define $\|\gz\|_{2,\infty}:=\sup_{x\in\Gamma^{(0)}}\|\gz\|_{L^2(\Gamma^x,\mu^x)}$ for $\gz\in\Cc_c\big(\Gamma^{(1)}\big)$. The equality
$$
\|\pi^x(\fz)\| \;
	= \; \sup\Big\{
		\Big\langle\,\gz|_{\Gamma^x}\; ,\; \pi^x(\fz)\hz|_{\Gamma^x}\Big\rangle \; \Big|\;
			\gz,\hz\in\Cc_c\big(\Gamma^{(1)}\big) \text{ with }\|\gz\|_{2,\infty},\|\hz\|_{2,\infty}\leq 1
	\Big\}
$$
is derived by a standard argument of operator theory and the previous considerations. Note that $\pi^x(\fz)$ is a bounded operator on the Hilbert space $L^2(\Gamma^x,\mu^x)$. For each $\gz,\hz\in\Cc_c\big(\Gamma^{(1)}\big)$, the map $x\in\Gamma^{(0)}\mapsto \big\langle\gz|_{\Gamma^x},\pi^x(\fz)\hz|_{\Gamma^x}\big\rangle\in\CM$ is continuous since all the functions are continuous with compact support.

\vspace{.1cm}

The previous considerations imply that $\Gamma^{(0)}\ni x\mapsto \|\pi^x(\fz)\|$ is lower semi-continuous as it is the supremum of continuous functions. Define $\Psi:\ts\to\RM$ by
$$
\Psi(t) \;
	= \; \|\fz_t\|_{red} \; 
		= \; \sup\big\{ 
			\|\pi^x(\fz_t)\| \;\big|\; 
				x\in\Gamma^{(0)}\text{ with } p(x)=t 
		\big\}\,,
		\qquad
		t\in\ts\,.
$$
The map $q:\Gamma^{(0)}\to\ts$ defined by $q:=p|_{\Gamma^{(0)}}$ is open and surjective, c.f. proof of Corollary~\ref{Chap4-Cor-FieldGroupContCAlgebras}. Since $\Psi$ is represented by the supremum over a lower semi-continuous map and $q$ is open and surjective, a standard topological argument leads to the lower semi-continuity of $\Psi$, c.f. \cite[Lemma~5.4]{LaRa99}.
\end{proof}

\begin{remark}
\label{Chap4-Rem-LRUseAmenable}
The proof of claim (a) can be also made with the reduced $C^\ast$-algebra whenever the short sequence is exact for the reduced groupoid $C^\ast$-algebras. This can be the case even if the groupoid is not topologically amenable. Then the map $\ts\ni t\mapsto\|\fz_t\|_{red}\in\RM$ is continuous by combining (a) and (b) of Theorem~\ref{Chap4-Theo-LanRamContFieldGroupoid}. This would be sufficient for our purposes as the generalized Schr\"odinger operator considered in Section~\ref{Chap2-Sect-SchrOp} are defined via the left-regular representation on the reduced $C^\ast$-algebra. On the other hand, there are groupoids so that this short sequence for the reduced $C^\ast$-algebras is not exact, c.f. \cite[Remark~4.10]{Ren91}. 
\end{remark}

\begin{remark}
\label{Chap4-Rem-ExactSequence}
The proof of the exactness of the short sequence
$$
0\longrightarrow \CG^\ast_{full}(\Gamma_U)\overset{e}{\longrightarrow}\CG^\ast_{full}(\Gamma)\overset{\imath^\ast}{\longrightarrow}\CG^\ast_{full}(\Gamma_t)\longrightarrow 0
$$
was provided by {\sc Torpe} \cite{Tor85} for $C^\ast$-algebras associated with foliations. Later {\sc Hilsum} and {\sc Skandalis} \cite{HiSk87} proved the exactness for general topological groupoids with left-continuous Haar system. The result can also be found in \cite{Renault80,Ren87,Ramazan98,Renault09}.
\end{remark}

Let $(\Gamma,\ts,p)$ be a continuous field of groupoids such that the fiber groupoids $\Gamma_t\, ,\; t\in\ts\,,$ are topologically amenable. Then the reduced and full $C^\ast$-algebra of $\Gamma_t$ coincide and we denote these $C^\ast$-algebras by $\CG^{\ast}(\Gamma_t)$ for $t\in\ts$, c.f. Definition~\ref{Chap2-Def-GroupoidCalgebra}. With this at hand, Theorem~\ref{Chap4-Theo-LanRamContFieldGroupoid} provides a recipe for building a continuous field of unital $C^*$-algebras from a continuous field of groupoids, c.f. \cite[Section~5]{LaRa99}.

\begin{theorem}[\cite{LaRa99}, continuous fields of groupoids and $C^\ast$-algebras]
\label{Chap4-Theo-FieldGroupContCAlgebras}
Let $(\Gamma,\ts,p)$ be a continuous field of groupoids such that $\Gamma$ is \'etale with compact unit space and left-continuous Haar system $\mu:=(\mu^x)_{x\in\Gamma^{(0)}}$. Suppose that the fiber groupoid $\Gamma_t$ is topologically amenable for each $t\in\ts$. Then there exists a unique $\Upsilon\subseteq\prod_{t\in\ts} \CG^{\ast}(\Gamma_t)$ containing $\Cc_c\big(\Gamma^{(1)}\big)$ such that $(\CG^{\ast}(\Gamma_t)_{t\in\ts},\Upsilon)$ is a continuous field of unital $C^\ast$-algebras. In particular, this statement holds if $\Gamma$ itself is topologically amenable.
\end{theorem}

\begin{proof}
The fiber groupoids $\Gamma_t\, , \; t\in\ts\, ,$ are topologically amenable. De\-fine $\Lambda:=\Cc_c\big(\Gamma^{(1)}\big)$. Then every $\fz\in\Cc_c\big(\Gamma^{(1)}\big)$ is represented by a vector field $(\fz_t)_{t\in\ts}\in\prod_{t\in\ts}\CG^\ast(\Gamma_t)$ where $\fz_t:=\fz|_{\Gamma_t}$. Consequently, $\Lambda$ is represented as a subset of $\prod_{t\in\ts}\CG^\ast(\Gamma_t)$. According to Proposition~\ref{Chap3-Prop-GenerFamiContField} it suffices to show that $\Lambda$ is a generating family, i.e., $\Lambda$ satisfies \nameref{(CFC1)}-\nameref{(CFC3)}. The in particular part follows then by Proposition~\ref{Chap4-Prop-ContFieldGroupoidAmenable}.

\vspace{.1cm}

\nameref{(CFC1)}: According to Proposition~\ref{Chap2-Prop-InvolutiveGroupoidAlgebra}, $\Lambda$ is a complex linear space that is closed under multiplication and involution. Since $\Gamma$ is \'etale, the unit space $\Gamma^{(0)}$ is open, c.f. Lemma~\ref{Chap2-Lem-CharUnitSpaceOpen} and Remark~\ref{Chap2-Rem-CharUnitSpaceOpen}. Moreover, by assumption, the unit space $\Gamma^{(0)}$ is compact. Hence, the unit $\mathpzc{1}$ of the reduced $C^\ast$-algebra $\CG_{red}^\ast(\Gamma)$ is the characteristic function $\chi_{\Gamma^{(0)}}\in\Cc_c\big(\Gamma^{(1)}\big)$, c.f. Theorem~\ref{Chap2-Theo-UnitGroupoidCalgebra}. Then $\mathpzc{1}_t=\chi_{\Gamma^{(0)}\cap\Gamma^{(1)}_t}\in\Cc_c\big(\Gamma^{(1)}_t\big)$ is the unit of $\CG^\ast(\Gamma_t)$. Thus, $(\mathpzc{1}_t)_{t\in\ts}$ is an element of $\Lambda$.

\vspace{.1cm}

\nameref{(CFC2)}: Let $t_0\in\ts$ and $\fz\in\Cc_c\big(\Gamma_t^{(1)}\big)$. Then $\Gamma^{(1)}_t$ is a closed subset of $\Gamma^{(1)}$, c.f. Lemma~\ref{Chap4-Lem-FiberClosGroup}. Thus, there exists a $\gz\in\Cc_c\big(\Gamma^{(1)}\big)$ such that $\gz(\gamma)=\fz(\gamma)$ for $\gamma\in\Gamma^{(1)}_t$, c.f. Lemma~\ref{Chap4-Lem-CompContExtGroupoid}. Consequently, $\Cc_c\big(\Gamma^{(1)}_t\big)$ is a subset of $\Lambda_t:=\{\fz_t\;|\; \fz\in\Lambda\}\subseteq\CG^\ast(\Gamma_t)$. Due to definition of $\CG^\ast(\Gamma_t)$, the set $\Cc_c\big(\Gamma^{(1)}_t\big)$ is a dense subset of $\CG^\ast(\Gamma_t)$. Hence, $\Lambda_t$ is dense in $\CG^\ast(\Gamma_t)$.

\vspace{.1cm}

\nameref{(CFC3)}: Proposition~\ref{Chap2-Prop-Red=FullCalgebra} implies that, for $\fz_t\in\Cc_c\big(\Gamma_t^{(1)}\big)$, the equation $\|\fz_t\|_{red}=\|\fz_t\|_{full}$ holds. Hence, Theorem~\ref{Chap4-Theo-LanRamContFieldGroupoid} leads to the continuity of $\ts\ni t\mapsto\|\fz_t\|\in\RM$ for $\fz\in\Cc_c\big(\Gamma^{(1)}\big)$ where $\|\fz_t\|:=\|\fz_t\|_{red}=\|\fz_t\|_{full}$ is the $C^\ast$-norm of $\CG^\ast(\Gamma_t)$, c.f. Definition~\ref{Chap2-Def-GroupoidCalgebra}. 
\end{proof}

\begin{remark}
\label{Chap4-Rem-TopoAmenableWeaken}
(i) The amenability assumption in Theorem~\ref{Chap4-Theo-FieldGroupContCAlgebras} only serves the purpose to ensure that the reduced and the full $C^\ast$-algebra of $\Gamma_t$ agree for $t\in\ts$. The assertion is also valid if the condition topologically amenable is dropped and only the reduced and full $C^\ast$-algebras associated with all fiber groupoids coincide. It was shown in \cite{Wil15} that the property of a groupoid $\Gamma$ being topologically amenable is not equivalent to the fact that the groupoid $C^\ast$-algebras $\CG^\ast_{red}(\Gamma)$ and $\CG^\ast_{full}(\Gamma)$ are isomorphic. More precisely, there are \'etale, topological groupoids with compact unit space being not topologically amenable while the $C^\ast$-algebras $\CG^\ast_{red}(\Gamma)$ and $\CG^\ast_{full}(\Gamma)$ are isomorphic. This is different to the group case where amenability is equivalent to the desired existence of an isomorphism, c.f. \cite{Hul66}. The reader is referred to \cite{Ana16} and references therein for recent results about the question when the reduced and the full $C^\ast$-algebras are isomorphic.

\vspace{.1cm}

(ii) The compactness of $\Gamma^{(0)}$ is required so that the associated $C^\ast$-algebras have a unit. In view of Example~\ref{Chap3-Ex-ConstTermPolNecess}, this is crucial for our purposes. More precisely, the continuous behavior of the spectra holds only after excluding $0\in\RM$ if the vector field of units is not a continuous section of the continuous field of $C^\ast$-algebras, c.f. Remark~\ref{Chap3-Rem-ContSpectrNonUnitalAlg}. On the other hand, the compactness of $\Gamma^{(0)}$ implies also the compactness of $\ts$ which seems to be restrictive. Since the space $\ts$ is locally compact (Lemma~\ref{Chap4-Lem-TNormalContFieldGroup}), the space can always be restricted to a compact neighborhood of a $t_0\in\ts$. In view of that the requirement of $\ts$ being compact is not too limiting.
\end{remark}

The compactness of $\Gamma^{(0)}$ can be relaxed as shown by the following corollary. Recall that a map between topological spaces is called closed whenever it maps closed sets to closed sets. 

\begin{corollary}
\label{Chap4-Cor-FieldGroupContCAlgebras}
Let $(\Gamma,\ts,p)$ be a continuous field of groupoids such that $\Gamma$ is \'etale with left-continuous Haar system $\mu:=(\mu^x)_{x\in\Gamma^{(0)}}$ and the map $p:\Gamma^{(1)}\to\ts$ is closed. Suppose that the fiber groupoids $\Gamma_t\,,\; t\in\ts\,,$ have compact unit space and they are topologically amenable. Then there exists a unique $\Upsilon\subseteq\prod_{t\in\ts} \CG^{\ast}(\Gamma_t)$ containing $\Cc_c\big(\Gamma^{(1)}\big)$ such that $(\CG^{\ast}(\Gamma_t)_{t\in\ts},\Upsilon)$ is a continuous field of unital $C^\ast$-algebras. 
\end{corollary}

\begin{proof}
According to Theorem~\ref{Chap4-Theo-FieldGroupContCAlgebras} and Remark~\ref{Chap4-Rem-TopoAmenableWeaken}~(ii), it suffices to show that, for every $t_0\in\ts$, there exist an open neighborhood $U_0\subseteq\ts$ of $t_0$ and an $\fz\in\Cc_c\big(\Gamma^{(1)}\big)$ such that $\fz_t=\mathpzc{1}_t$ for $t\in U_0$ where $\mathpzc{1}_t$ is the unit of $\CG_t$, namely the second part of \nameref{(CFC1)} has to be verified.

\vspace{.1cm}

Let $t_0\in\ts$. According to Lemma~\ref{Chap4-Lem-TNormalContFieldGroup}, $\ts$ is a normal space. Hence, there are open sets $U_0,U_1\subseteq\ts$ and a compact set $K_0\in\ts$ such that $t_0\in U_0\subseteq K_0\subseteq U_1$. Furthermore, the inclusion
$$
F \;
	:= \; p^{-1}(K_0)\cap\Gamma^{(0)} \;
	\subseteq \; p^{-1}(U_1)\cap\Gamma^{(0)} \;
	=: \; V
$$
follows by assumption $K_0\subseteq U_1$. The unit space $\Gamma^{(0)}$ is open as $\Gamma$ is \'etale, c.f. Lemma~\ref{Chap2-Lem-CharUnitSpaceOpen} and Remark~\ref{Chap2-Rem-CharUnitSpaceOpen}. Since $p$ is, additionally, continuous, the set $V$ is open as an intersection of open sets. Thus, its complement $\Gamma^{(1)}\setminus V$ is closed. 

\vspace{.1cm}

With this at hand, the set $F$ is compact by the following reason. Denote by $q:\Gamma^{(0)}\to\ts$ the restriction of the map $p$ to the unit space $\Gamma^{(0)}$. Since $p$ is continuous, $q$ is continuous as well. The surjectivity of $p$ and the identity $p|_{\Gamma^{(0)}}\circ r =p$ imply that $q$ is surjective. The set $\Gamma^{(0)}$ is open as $\Gamma$ is \'etale. Proposition~\ref{Chap2-Prop-BasGroupoid}~(g) implies that $\Gamma^{(0)}$ is also closed. Thus, the restriction $q$ of the map $p$ to $\Gamma^{(0)}$ is an open and closed map since $p$ is so. Altogether, the map $q:\Gamma^{(0)}\to\ts$ is an open, closed, continuous and surjective map and the preimage $q^{-1}(\{t\})=\Gamma^{(0)}_t$ is compact for each $t\in\ts$ by assumption. Consequently, Proposition~\ref{App1-Prop-ContClosMapIsProper} leads to the compactness of $q^{-1}(K_0)=F$ as $K_0\subseteq\ts$ is compact.

\vspace{.1cm}

The space $\Gamma^{(1)}$ is normal since it is second countable, locally compact and Hausdorff, c.f. Proposition~\ref{App1-Prop-LocCompHausSecCountImplNormal}. Hence, the Lemma of Urysohn applies, c.f. Proposition~\ref{App1-Prop-LemmaUrysohn}. Thus, there exists a continuous function $\fz:\Gamma^{(1)}\to [0,1]$ such that $\fz|_{F}\equiv 1$ and $\fz|_{\Gamma^{(1)}\setminus V}\equiv 0$. According to Proposition~\ref{App1-Prop-LemmaUrysohn}, there is no loss of generality in assuming that the map $\fz$ has compact support as $F$ is compact. Consequently, $\fz$ is an element of $\Cc_c\big(\Gamma^{(1)}\big)$. 

\vspace{.1cm}

By construction, $\fz$ is identically to $1$ on $q^{-1}(\{t\})=\Gamma^{(0)}_t$ for each $t\in U_0$. Thus, $\fz_t$ is the unit $\mathpzc{1}_t$ of $\CG^\ast(\Gamma_t)$ for $t\in U_0$, c.f. Theorem~\ref{Chap2-Theo-UnitGroupoidCalgebra}. Since $\fz\in\Cc_c\big(\Gamma^{(1)}\big)$, the second part of \nameref{(CFC1)} is satisfied leading to the desired result by Theorem~\ref{Chap4-Theo-FieldGroupContCAlgebras}.
\end{proof}

\begin{remark}
\label{Chap4-Rem-FieldGroupContCAlgebras}
Note that the condition that $\Gamma_t$ has compact unit space cannot be relaxed since otherwise $\CG^\ast(\Gamma_t)$ does not have a unit in general, c.f. Theorem~\ref{Chap2-Theo-UnitGroupoidCalgebra} as well as discussion in Section~\ref{Chap8-Sect-Groupoids}.
\end{remark}

Let $(\Gamma,\ts,p)$ be a continuous field of groupoids where $\Gamma$ is an \'etale groupoid with left-continuous Haar system. According to Proposition~\ref{Chap2-Prop-C0functionRedCAlgebra}, every element $\fz\in\CG^\ast_{red}(\Gamma)$ is uniquely identified with a continuous function on $\Gamma^{(1)}$ that vanishes at infinity. Thus, the restriction $\fz|_{\Delta^{(1)}}$ is also a continuous function on $\Delta^{(1)}$ that vanishes at infinity if $\Delta$ is a closed subgroupoid of $\Gamma$. Furthermore, $\fz|_{\Delta^{(1)}}$ defines an element of $\CG^\ast_{red}(\Delta)$ as the reduced norm of the groupoid $\Delta$ is bounded by the reduced-norm of $\Gamma$. Hence, each $\fz\in\CG^\ast_{red}(\Gamma)$ defines a vector field $(\fz_t)_{t\in\ts}$ in $\prod_{t\in\ts}\CG^\ast_{red}(\Gamma_t)$ where $\fz_t:=\fz|_{\Gamma^{(1)}_t}$ for $t\in\ts$.

\begin{proposition}
\label{Chap4-Prop-CAlgGUnivContVectField}
Let $(\Gamma,\ts,p)$ be a continuous field of groupoids such that $\Gamma$ is \'etale with compact unit space and left-continuous Haar system $\mu:=(\mu^x)_{x\in\Gamma^{(0)}}$. Suppose that the fiber groupoids $\Gamma_t\,,\; t\in\ts\,,$ are topologically amenable. Consider the unique continuous field of $C^\ast$-algebras $(\CG^{\ast}(\Gamma_t)_{t\in\ts},\Upsilon)$ with $\Cc_c\big(\Gamma^{(1)}\big)\subseteq\Upsilon$ which was defined in Theorem~\ref{Chap4-Theo-FieldGroupContCAlgebras}. Then $\CG^\ast_{red}(\Gamma)$ is contained in $\Upsilon$, i.e., each $\fz\in\CG^\ast_{red}(\Gamma)$ defines a continuous vector field $(\fz_t)_{t\in\ts}\in\prod_{t\in\ts}\CG^\ast(\Gamma_t)$ where $\fz_t:=\fz|_{\Gamma_t}$ is the restriction of the function $\fz$ to $\Gamma^{(1)}_t$ for $t\in\ts$.
\end{proposition}

\begin{proof}
For the convenience of the reader, the norms are indexed in this proof by the corresponding $C^\ast$-algebra in which the norm is taken. Let $\fz\in\CG^\ast_{red}(\Gamma)$. By definition of the reduced $C^\ast$-algebra, there is a sequence $\fz_n\in\Cc_c\big(\Gamma^{(1)}\big)$ such that $\lim_{n\to\infty}\|\fz-\fz_n\|_{\CG^\ast_{red}(\Gamma)}=0$. Thus, $(\fz_n)_{n\in\NM}$ is a Cauchy sequence in $\CG^\ast_{red}(\Gamma)$. Let $\varepsilon>0$. Then there exists an $n(\varepsilon)\in\NM$ such that $\|\fz_n-\fz_m\|_{\CG^\ast_{red}(\Gamma)}<\varepsilon$ for $n,m\geq n(\varepsilon)$. Additionally, the restriction $\fz_n|_{\Gamma_t^{(1)}}$ is an element of $\Cc_c\big(\Gamma^{(1)}_t\big)$ for $t\in\ts$. Since $\Gamma^{(0)}_t\subseteq\Gamma^{(0)}$, the sequence $\Big(\fz_n|_{\Gamma_t^{(1)}}\Big)_{n\in\NM}$ in $\CG^\ast_{red}(\Gamma_t)$ satisfies 
$$
\Big\|\fz_n|_{\Gamma_t^{(1)}}-\fz_m|_{\Gamma_t^{(1)}}\Big\|_{\CG^\ast_{red}(\Gamma_t)} \;
	\leq \; \big\|\fz_n-\fz_m\big\|_{\CG^\ast_{red}(\Gamma)} \;
	<\; \varepsilon\,,
	\qquad n,m\geq n(\varepsilon)\, .
$$
Thus, it is a Cauchy-sequence in $\CG^\ast_{red}(\Gamma_t)$ and its limit is denoted by $\fz_t$. Furthermore, it converges uniformly in $t\in\ts$ since the estimate in $\varepsilon$ is independent of $t\in\ts$. By construction, $\fz_t$ is exactly the restriction $\fz|_{\Gamma_t}$. For each $n\in\NM$, $\fz_n\in\Cc_c\big(\Gamma^{(1)}_t\big)$ is a continuous vector field. Hence, \nameref{(CFC4)} implies that $(\fz_t)_{t\in\ts}$ is also a continuous vector field as 
$$
\Big\|\fz_t-\fz_n|_{\Gamma^{(1)}_t}\Big\|_{\CG^\ast_{red}(\Gamma_t)} \;
	< \; \varepsilon\,,
	\qquad t\in\ts\,,\; \varepsilon>0\,,
$$ 
hold by the previous considerations for $n\geq n(\varepsilon)$.
\end{proof}

\begin{remark}
\label{Chap4-Rem-CAlgGUnivContVectField}
The statement of Proposition~\ref{Chap4-Prop-CAlgGUnivContVectField} extends to the setting of Corollary~\ref{Chap4-Cor-FieldGroupContCAlgebras}. More specifically, if $\Gamma^{(0)}$ is not compact but the sets $\Gamma^{(0)}_t\,,\; t\in\ts\,,$ are still compact and $p:\Gamma^{(1)}\to\ts$ is a closed map, then the elements of $\CG^\ast_{red}(\Gamma)$ are represented by a continuous vector field in $\prod_{t\in\ts}\CG^\ast(\Gamma_t)$.
\end{remark}

\section{Groupoid isomorphisms and their \texorpdfstring{$C^\ast$-algebras}{C-algebras}}
\label{Chap4-Sect-GroupoidIsomorphism}

In the following, the notion of a measure preserving groupoid isomorphism is intro\-duced. Such an isomorphism naturally induces an isometric, $\ast$-isomorphism between the associated reduced $C^\ast$-algebras which is used in Section~\ref{Chap4-Sect-UnivGroupDynSyst}.

\begin{definition}[Groupoid homomorphism / isomorphism]
\label{Chap4-Def-GroupoidIsomorphism}
Let $\Gamma_1,\Gamma_2$ be two topological groupoids. A continuous map $\Phi:\Gamma^{(1)}_1\to\Gamma^{(1)}_2$ is called a {\em (groupoid) homomorphism} if $(\gamma,\eta)\in\Gamma^{(2)}_1$ implies that $(\Phi(\gamma),\Phi(\eta))\in\Gamma^{(2)}_2$ and $\Phi(\gamma\circ\eta)=\Phi(\gamma)\circ\Phi(\eta)$. The map $\Phi$ is called a {\em (groupoid) isomorphism} if, additionally, $\Phi$ is a homeomorphism.
\end{definition}

The following proposition collects some basic properties for a groupoid isomorphism, see e.g. \cite[Proposition~1.7]{GoehleThesis09}.

\begin{proposition}[\cite{GoehleThesis09}]
\label{Chap4-Prop-PropGroupoidIsomorphism}
Let $\Gamma_1,\Gamma_2$ be two topological groupoids with groupoid isomorphism $\Phi:\Gamma_1^{(1)}\to\Gamma_2^{(1)}$. Then the following assertions hold.
\begin{itemize}
\item[(a)] The tuple $(\gamma,\eta)$ is an element of $\Gamma^{(2)}_1$ if and only if $(\Phi(\gamma),\Phi(\eta))\in\Gamma^{(2)}_2$.
\item[(b)] The relation $x\in\Gamma^{(0)}_1$ holds if and only if $\Phi(x)\in\Gamma^{(0)}_2$.
\item[(c)] The map $\Phi$ commutes with the inverse of the groupoids, i.e., $\Phi(\gamma^{-1})=\Phi(\gamma)^{-1}$ for all $\gamma\in\Gamma^{(1)}_1$.
\item[(d)] For each $\gamma\in\Gamma^{(1)}_1$, the equalities $r\big(\Phi(\gamma)\big)=\Phi\big(r(\gamma)\big)$ and $s\big(\Phi(\gamma)\big)=\Phi\big(s(\gamma)\big)$ are valid. In particular, the identity $\Phi(\Gamma_1^x)=\Gamma^{\Phi(x)}_2$ holds for the $r$-fiber $\Gamma_1^x$ associated with $x\in\Gamma^{(0)}_1$.
\end{itemize}
\end{proposition}

\begin{proof}
Since $\Phi$ is bijective, it has an inverse. Then it is not difficult to check that the inverse $\Phi^{-1}:\Gamma^{(1)}_2\to\Gamma^{(1)}_1$ is a groupoid homomorphism as well.

\vspace{.1cm}

(a): By definition, $(\gamma,\eta)\in\Gamma^{(2)}_1$ implies that $(\Phi(\gamma),\Phi(\eta))\in\Gamma^{(2)}_2$. The converse follows as $\Phi^{-1}$ is a groupoid homomorphism as well.

\vspace{.1cm}

(b): It suffices to show that $x\in\Gamma^{(0)}_1$ implies $\Phi(x)\in\Gamma^{(0)}_2$ since the converse follows by passing to groupoid homomorphism $\Phi^{-1}$. Let $x\in\Gamma^{(0)}_1$. Then the equalities $\Phi(x) = \Phi(x x) = \Phi(x)\Phi(x)$ is derived since $x$ is a unit. Hence, $\Phi(x)\Phi(x)^{-1} = \Phi(x)$ is deduced. Thus, $\Phi(x)\in r\big(\Gamma_2^{(1)}\big)=\Gamma^{(0)}_2$.

\vspace{.1cm}

(c): Let $\gamma\in\Gamma^{(1)}_1$. Then the relations $\big(\gamma, \gamma^{-1}\big)\in\Gamma^{(2)}_1$ and $\big(\Phi(\gamma)^{-1},\Phi(\gamma)\big)\in\Gamma^{(2)}_2$ are valid by \nameref{(G2)}. Consequently, the equations
$$
\Phi(\gamma)^{-1}\Phi\big(\gamma\gamma^{-1}\big) \; 
	= \; \Phi(\gamma)^{-1}\Phi(\gamma)\Phi\big(\gamma^{-1}\big) \;
		= \; \Phi\big(\gamma^{-1}\big)
$$
are derived. According to (b), $\Phi\big(\gamma\gamma^{-1}\big)$ is a unit of the groupoid $\Gamma_2$. Thus, the left hand side is equal to $\Phi(\gamma)^{-1}$ by \nameref{(G3)} implying the desired identity.

\vspace{.1cm}

(d): The equalities $r\big(\Phi(\gamma)\big)=\Phi\big(r(\gamma)\big)$ and $s\big(\Phi(\gamma)\big)=\Phi\big(s(\gamma)\big)$ for $\gamma\in\Gamma^{(1)}_1$ immediately follow by (iii) and the definition of the range and the source map. Consequently,
$$
\Phi(\Gamma_1^x) \;
	= \; \Phi\Big( \Big\{ \gamma\in\Gamma^{(1)}_1 \;\Big|\; r(\gamma)=x \Big\} \Big) \;
		= \; \Big\{ \Phi(\gamma)\in\Gamma^{(1)}_2 \;\Big|\; \Phi(r(\gamma))=\Phi(x) \Big\} 
			= \;\Gamma^{\Phi(x)}_2
$$
is deduced from the bijectivity of $\Phi$.
\end{proof}

\medskip

Since every groupoid isomorphism $\Phi$ is by definition a homeomorphism, it is also Borel-measurable. Thus, $\Phi(F)$ is measurable if $F\subseteq\Gamma^{(1)}_1$ is Borel-measurable.

\begin{definition}[Measure preserving groupoid isomorphism]
\label{Chap4-Def-MeasPreservGroupIsomorph}
For $i=1,2$, let $\Gamma_i$ be a topological groupoid with left-continuous Haar system $\mu_i^x\, ,\; x\in\Gamma^{(0)}_i$. A groupoid isomorphism $\Phi:\Gamma^{(1)}_1\to\Gamma^{(1)}_2$ is called {\em measure preserving} if $\mu_1^x(F)=\mu_2^{\Phi(x)}(\Phi(F))$ holds for all $x\in\Gamma^{(0)}_1$ and each Borel measurable set $F\subseteq\Gamma^{(1)}_1$.
\end{definition}

The following assertion is a basic property of a measure preserving groupoid isomorphism that follows by Lebesgue monotone convergence theorem.

\begin{lemma}
\label{Chap4-Lem-IdentIntegrMeasPresIsom}
For $i=1,2$, let $\Gamma_i$ be a topological groupoid with left-continuous Haar system $\mu_i^x\, ,\; x\in\Gamma^{(0)}_i$. Consider a measure preserving groupoid isomorphism $\Phi:\Gamma^{(1)}_1\to\Gamma^{(1)}_2$. Then
$$
\int_{\Gamma^x_1} \fz(\rho) \; d\mu_1^x(\rho) \; 
	= \; \int_{\Gamma^{\Phi(x)}_2} \fz\circ\Phi^{-1}(\zeta) \; d\mu_2^{\Phi(x)}(\zeta)
$$
holds for all $x\in\Gamma^{(0)}_1$ and $\fz\in\Cc_c\big(\Gamma^{(1)}_1\big)$.
\end{lemma}

\begin{proof}
According to Proposition~\ref{Chap4-Prop-PropGroupoidIsomorphism}~(iv), $\Gamma^x_1$ is mapped bijectively onto $\Gamma^{\Phi(x)}_2$. Since $\Phi$ is measure preserving it immediately follows that the desired identity holds if $\fz$ is a simple function on the Borel-measure space $\big(\Gamma^{(1)}_1,\mu_1^x\big)$ for $x\in\Gamma^{(0)}_1$. Then, for each non-negative $\fz\in\Cc_c\big(\Gamma^{(1)}_1\big)$, there exists a non-decreasing sequence of simple functions $(\fz_n)_{n\in\NM}$ such that it converges pointwise to $\fz$. According to Lebesgue monotone convergence theorem, the limit $\lim_{n\to\infty}\int\fz_n\; d\mu_1^x$ coincides with $\int\fz\; d\mu_1^x$. Then $\big(\fz_n\circ\Phi^{-1}\big)_{n\in\NM}$ converges pointwise to $\fz\circ\Phi^{-1}$ implying 
$$
\lim_{n\to\infty}\int_{\Gamma^{\Phi(x)}_2}\fz_n\circ\Phi^{-1}(\zeta)\; d\mu_2^{\Phi(x)}(\zeta) \; 
	= \; \int_{\Gamma^{\Phi(x)}_2}\fz\circ\Phi^{-1}(\zeta)\; d\mu_2^{\Phi(x)}(\zeta)\,,
	\qquad
	x\in\Gamma^{(0)}_1\,.
$$
This leads to
$$
\int_{\Gamma^x_1} \fz(\rho) \; d\mu_1^x(\rho) \; 
	= \; \int_{\Gamma^{\Phi(x)}_2} \fz\circ\Phi^{-1}(\zeta) \; d\mu_2^{\Phi(x)}(\zeta)\,,
	\qquad
	x\in\Gamma^{(0)}_1\,.
$$
By taking the positive and the negative part of $\fz\in\Cc_c\big(\Gamma^{(1)}_1\big)$, the desired assertion follows for $\fz\in\Cc_c\big(\Gamma^{(1)}_1\big)$ since the integral is linear.
\end{proof}

\medskip

For a topological groupoid with left-continuous Haar system, the reduced $C^\ast$-algebra is defined in Section~\ref{Chap2-Sect-GroupoidCalgebras}. The notion of a measure preserving groupoid isomorphism leads to an isometric, $\ast$-isomorphism of the reduced $C^\ast$-algebras.

\begin{definition}
\label{Chap4-Def-CAlgIsometricIsomorphism}
Let $\CG_1,\CG_2$ be $C^\ast$-algebras. A map $\Theta:\CG_1\to\CG_2$ is called a {\em $\ast$-homo\-morphism} if $\Theta(\fz^\ast)=\Theta(\fz)^\ast$ and 
$$
\Theta(\fz+\lambda\cdot\gz) \;
	= \; \Theta(\fz)+\lambda\cdot \Theta(\gz)\,,
		\qquad
			\Theta(\fz\star\gz)=\Theta(\fz)\star\Theta(\gz)\,,
$$
hold for all $\fz,\gz\in\CG_1$ and $\lambda\in\CM$. The map $\Theta$ is called {\em isometric} if, additionally, $\|\fz\|=\|\Theta(\fz)\|$ is valid for all $\fz\in\CG_1$. If a $\ast$-homomorphism $\Theta:\CG_1\to\CG_2$ is bijective, $\Theta$ is called a {\em $\ast$-isomorphism}.
\end{definition}

It is well-known that a $\ast$-homomorphism between $C^\ast$-algebras is always norm-decreasing, see e.g. \cite[Theorem~2.1.7]{Murphy90}. Furthermore, if $\Theta$ is an injective $\ast$-homomorphism it follows that $\Theta$ is also isometric, see e.g. \cite[Theorem~3.1.4]{Murphy90}. Thus, a $\ast$-isomorphism of $C^\ast$-algebras is automatically isometric.

\medskip

A $\ast$-homomorphism between unital $C^\ast$-algebras is called {\em unital} if it maps the unit to the unit. In general, it is not clear that a $\ast$-homomorphism is unital since the map $\Theta:\CG_1\to\CG_2$ that sends every element of $\CG_1$ to $0\in\CG_2$ is a $\ast$-homomorphism. On the other hand, if $\Theta$ is a $\ast$-isomorphism, it is automatically unital. 

\begin{lemma}
\label{Chap4-Lem-CAlgIsomorphismUnital}
Let $\CG_1$ and $\CG_2$ be unital $C^\ast$-algebras with unit $\mathpzc{1}_1\in\CG_1$ and $\mathpzc{1}_2\in\CG_2$. Suppose there is a $\ast$-isomorphism $\Theta:\CG_1\to\CG_2$. Then the equation $\Theta(\mathpzc{1}_1)=\mathpzc{1}_2$ holds, i.e., $\Theta$ is unital.
\end{lemma}

\begin{proof}
Since $\Theta$ is surjective, there exists an $\fz\in\CG_1$ such that $\Theta(\fz)=\mathpzc{1}_2$. Then the equations
$$
\mathpzc{1}_2 \; 
	= \; \Theta(\fz) \;
		= \; \Theta(\fz \star \mathpzc{1}_1) \; 
			= \; \Theta(\fz)\star \Theta(\mathpzc{1}_1) \; 
				= \; \mathpzc{1}_2 \star \Theta(\mathpzc{1}_1) \; 
					= \; \Theta(\mathpzc{1}_1)
$$
hold leading to the desired result.
\end{proof}

\begin{remark}
\label{Chap4-Rem-IsomorphSpectr}
Let $\Theta:\CG_1\to\CG_2$ be a $\ast$-isomorphism. Then $\fz\in\CG_1$ is invertible if and only if $\Theta(\fz)$ is invertible. Thus, the spectrum $\sigma(\fz)$ coincides with the spectrum $\sigma(\Theta(\fz))$. Thus, $\Theta$ preserves the spectrum.
\end{remark}

The following well-known assertion is used in Section~\ref{Chap4-Sect-UnivGroupDynSyst}, c.f. \cite[Proposition~4.23]{GoehleThesis09}.

\begin{proposition}
\label{Chap4-Prop-CAlgIsoGroupIso}
For $i=1,2$, let $\Gamma_i$ be a topological groupoid with left-continuous Haar system $\mu_i^x\, ,\; x\in\Gamma^{(0)}_i$. If $\Phi:\Gamma^{(1)}_1\to\Gamma^{(1)}_2$ is a measure preserving groupoid isomorphism, then
$$
\Theta:\Cc_c\Big(\Gamma_1^{(1)}\Big)\to\Cc_c\Big(\Gamma_1^{(1)}\Big)\,,\; 
	f\mapsto f\circ\Phi^{-1}\,,
$$
extends uniquely to an isometric, $\ast$-isomorphism $\Theta:\CG_{red}^\ast(\Gamma_1)\to\CG_{red}^\ast(\Gamma_2)$ such that $\Theta\big(\Cc_c\big(\Gamma_1^{(1)}\big)\big)=\Cc_c\big(\Gamma_2^{(1)}\big)$.
\end{proposition}

\begin{proof} 
It suffices to show that $\Theta:\Cc_c\big(\Gamma_1^{(1)}\big)\to\Cc_c\big(\Gamma_2^{(1)}\big)$ is a bijective, $\ast$-homomorphism satisfying $\|\fz\|=\|\Theta(\fz)\|$ for all $\fz\in\Cc_c\big(\Gamma_1^{(1)}\big)$. Then $\Theta$ extends uniquely since $\Cc_c\big(\Gamma^{(1)}_1\big)\subseteq\CG^\ast(\Gamma_1)$ and $\Cc_c\big(\Gamma^{(1)}_2\big)\subseteq\CG^\ast(\Gamma_2)$ are dense $\ast$-subalgebras.

\vspace{.1cm}

Define $\Theta:\Cc_c\big(\Gamma_1^{(1)}\big)\to\Cc_c\big(\Gamma_2^{(1)}\big)$ by $\Theta(\fz):=\fz\circ\Phi^{-1}$ for $\fz\in\Cc_c\big(\Gamma_1^{(1)}\big)$. Since $\Phi$ is a homeomorphism, the map $\Theta$ is bijective as well. Let $\fz,\gz\in\Cc_c\big(\Gamma_1^{(1)}\big)$ and $\lambda\in\CM$. Then, the equations 
$$
\Theta(\fz+\lambda\cdot\gz) \;
	= \; (\fz+\lambda\cdot\gz)\circ\Phi^{-1} \;
	= \; \Theta(\fz) + \lambda\cdot \Theta(\gz)
$$
are derived. Thus, $\Theta$ is linear. For $\gamma\in\Gamma_2^{(1)}$, the identities
\begin{align*}
\Theta(\fz\star\gz)(\gamma) \; &\overset{\quad\quad\quad\quad\;\;\;\;}{=} \; \big(\fz\star\gz\big)(\Phi^{-1}(\gamma))\\
	&\overset{\text{Prop.~\ref{Chap4-Prop-PropGroupoidIsomorphism}~(iv)}}{=} \; 
		\int_{\Gamma_1^{\Phi^{-1}(r(\gamma))}} 
			\fz(\rho)\cdot 
			\gz\big( \rho^{-1}\circ \Phi^{-1}(\gamma) \big) 
		\; d\mu_1^{\Phi^{-1}(r(\gamma))}(\rho)\\
	&\overset{\;\;\;\;\text{Lem.~\ref{Chap4-Lem-IdentIntegrMeasPresIsom}}\;\;\;\;}{=} \; 
		\int_{\Gamma_2^{r(\gamma)}} 
			\big(\fz\circ\Phi^{-1}\big)(\zeta) \cdot 
			\big(\gz\circ\Phi^{-1}\big)\big( \zeta^{-1} \gamma \big) 
		\; d\mu_2^{r(\gamma)}(\zeta)\\
	&\overset{\quad\quad\quad\quad\;\;\;\;}{=} \; 
		\big(\Theta(\fz)\star\Theta(\gz)\big)(\gamma)
\end{align*}
and
\begin{align*}
\Theta(\fz)^\ast(\gamma)\; & \overset{\quad\quad\quad\quad\;\;\;\;}{=} \; (\fz\circ\Phi^{-1})^\ast(\gamma) \\
	& \overset{\quad\quad\quad\quad\;\;\;\;}{=} \; 
		\overline{\fz\big(\Phi^{-1}(\gamma)^{-1}\big)} \\
	& \overset{\quad\quad\quad\quad\;\;\;\;}{=} \; 
		\Theta(\fz^\ast)(\gamma)
\end{align*}
follow by using that $\Phi^{-1}$ is a groupoid homomorphism. Hence, $\Theta$ defines a $\ast$-homo\-morphism. 

\vspace{.1cm}

According to Proposition~\ref{Chap4-Prop-PropGroupoidIsomorphism}~(iv) and the fact that $\Phi$ is measure preserving, the map 
$$
U_x:L^2\left(\Gamma_1^x,\mu_1^x\right)\to L^2\left(\Gamma_2^{\Phi(x)},\mu_2^{\Phi(x)}\right)\, ,
	\qquad \psi\mapsto\psi\circ\Phi^{-1}\, ,
$$ 
is an isomorphism for each $x\in\Gamma^{(0)}_1$. For $\fz\in\Cc_c\big(\Gamma_1^{(1)}\big),\; \gamma\in\Gamma^{x}_1$ and $\psi\in L^2\big(\Gamma_1^x,\mu_1^x\big)$, the equalities
\begin{align*}
\left(\pi^{x}(\fz)\, \psi\right)(\gamma)\; 
	&\overset{\quad\quad\quad\;\,}{=} \; 
		\int_{\Gamma_1^x} 
			\fz\big(\gamma^{-1}\rho\big) \cdot \psi(\rho) 
		\; d\mu_1^x(\rho)\\
	&\overset{\text{Lem.~\ref{Chap4-Lem-IdentIntegrMeasPresIsom}}}{=} \; 
		\int_{\Gamma_2^{\Phi(x)}} 
			\big(\fz\circ\Phi^{-1}\big)\big(\Phi(\gamma)^{-1}\zeta\big) 
			\cdot \big(\psi\circ\Phi^{-1}\big)(\zeta) 
		\; d\mu_2^{\Phi(x)}(\zeta)\\
	&\overset{\quad\quad\quad\;\,}{=} \; 
		\Big(\pi^{\Phi(x)}\big(\Theta(\fz)\big) \, U_{x}(\psi)\Big)(\Phi(\gamma))
\end{align*}
are deduced. Consequently, the equation $\|\pi^{x}(\fz)\psi\| = \|\pi^{\Phi(x)}(\Theta(\fz))U_{x}(\psi)\|$ is derived since $\Phi$ is measure preserving and bijective. This leads to 
\begin{align*}
\|\fz\|_{red} \; 
	&= \; \sup_{x\in\Gamma_1^{(0)}} \; \sup\Big\{ \|\pi^{x}(\fz)\psi\| \;\Big|\; \psi\in L^2\big(\Gamma_1^x,\mu_1^x\big) \text{ such that } \|\psi\|\leq 1 \Big\}\\
	&= \; \sup_{y\in\Gamma_2^{(0)}} \; \sup\Big\{ \|\pi^{y}(\Theta(\fz))\varphi\| \;\Big|\; \varphi\in L^2\big(\Gamma_2^y,\mu_2^y\big) \text{ such that } \|\varphi\|\leq 1 \Big\}\\
	&= \; \|\Theta(\fz)\|_{red}\, .
\end{align*}
as $U_x$ is an isomorphism and $\Phi$ bijectively maps $\Gamma^{(0)}_1$ onto $\Gamma^{(0)}_2$. Altogether, the previous considerations imply that $\Theta:\Cc_c\big(\Gamma_1^{(1)}\big)\to\Cc_c\big(\Gamma_2^{(1)}\big)$ is an isometric, bijective $\ast$-homomor\-phism. By density, $\Theta$ extends uniquely to an isometric, $\ast$-isomorphism from $\CG^\ast_{red}(\Gamma_1)$ onto $\CG^\ast_{red}(\Gamma_2)$.
\end{proof}

\medskip

In Theorem~\ref{Chap4-Theo-FieldGroupContCAlgebras}, the condition that a groupoid is topologically amenable is essential. This property is preserved by groupoid isomorphisms.

\begin{corollary}
\label{Chap4-Cor-GroupIsomTopAmen}
For $i=1,2$, let $\Gamma_i$ be a topological groupoid with left-continuous Haar system $\mu_i^x\, ,\; x\in\Gamma^{(0)}_i$. Suppose there is a measure preserving groupoid isomorphism $\Phi:\Gamma^{(1)}_1\to\Gamma^{(1)}_2$. Then $\Gamma_1$ is topologically amenable if and only if $\Gamma_2$ is topologically amenable.
\end{corollary}

\begin{proof}
Since $\Phi^{-1}:\Gamma^{(1)}_2\to\Gamma^{(1)}_1$ defines also a measure preserving isomorphism, it suffices to show that $\Gamma_2$ is topologically amenable if $\Gamma_1$ is so. Suppose that $\Gamma_1$ is topologically amenable with the left-continuous Haar system $\mu_1^x\,,\; x\in\Gamma^{(0)}_1$. Then there exists a sequence of non-negative functions $\{\fz_n\}_{n\in\NM}\subseteq\Cc_c\big(\Gamma^{(1)}_1\big)$ such that \nameref{(A1)} and \nameref{(A2)} hold. Let $\Theta:\CG^\ast_{red}(\Gamma_1)\to\CG^\ast_{red}(\Gamma_2)$ be the $\ast$-isomorphism of Proposition~\ref{Chap4-Prop-CAlgIsoGroupIso}. Then $\{\Theta(\fz_n)\}_{n\in\NM}\subseteq\Cc_c\big(\Gamma^{(1)}_2\big)$ is a sequence of non-negative functions. It suffices to verify \nameref{(A1)} and \nameref{(A2)} for this sequence.

\vspace{.1cm}

Let $x\in\Gamma^{(0)}_2$. Then the equations
$$
1 \; 
	= \; \int_{\Gamma^{\Phi^{-1}(x)}_1} \fz_n(\rho)\; d\mu_1^{\Phi^{-1}(x)}(\rho) \;
	\overset{\text{Lem.~\ref{Chap4-Lem-IdentIntegrMeasPresIsom}}}{=} \; 
		\int_{\Gamma^{x}_2} 
			\underbrace{\big(\fz_n\circ\Phi^{-1}\big)(\zeta)}_{=\Theta(\fz_n)(\zeta)}
		\; d\mu_2^{x}(\zeta) 
$$
are derived, namely $\{\Theta(\fz_n)\}_{n\in\NM}$ satisfies \nameref{(A1)}. Furthermore, Lemma~\ref{Chap4-Lem-IdentIntegrMeasPresIsom} leads to
\begin{align*}
m_n(\gamma) \; 
	&:= \; \int_{\Gamma_2^{r(\gamma)}} 
			\big| \Theta(\fz_n)\big(\gamma^{-1}\cdot\zeta\big)- \Theta(\fz_n)(\zeta)\big|
		\; d\mu_2^{r(\gamma)}(\zeta)\\
	&= \; \int_{\Gamma_1^{r(\Phi^{-1}(\gamma))}} 
			\Big|\fz_n\Big(\Phi^{-1}(\gamma)^{-1} \rho\Big)- \fz_n(\rho)\Big|
		\; d\mu_1^{r(\Phi^{-1}(\gamma))}(\rho)\\
	&=: \; \widehat{m}_n\big(\Phi^{-1}(\gamma)\big)
\end{align*}
for $n\in\NM$ and $\gamma\in\Gamma^{(1)}_2$. Since $\Phi$ is a homeomorphism it maps compact sets to compact sets. Consequently, $(m_n)_{n\in\NM}$ converges uniformly on compact sets to zero since $(\widehat{m}_n)_{n\in\NM}$ does so by assumption. Thus, $\{\Theta(\fz_n)\}_{n\in\NM}$ satisfies \nameref{(A2)} implying that $\Gamma_2$ is topologically amenable.
\end{proof}

\section{The universal dynamical system (groupoid) associated with a dynamical system}
\label{Chap4-Sect-UnivGroupDynSyst}

For a topological dynamical system $(X,G,\alpha)$, a new dynamical system is constructed by using the space $\SG(X)$. The so called universal dynamical system $(\Xun,G,\alpha)$ naturally induces a continuous field of groupoids where each fiber groupoid is isomorphic to a specific transformation group groupoid $Y\rtimes_\alpha G$ of a dynamical subsystem $Y\in\SG(X)$. This leads to a construction of a continuous field of unital $C^\ast$-algebras where specific families of generalized Schr\"odinger operators indexed by the topological space $\ts:=\SG(X)$ define continuous vector fields. Then Theorem~\ref{Chap3-Theo-ContFieldCALgContSpectr} implies the continuous behavior of the spectra of these operators in the parameter $Y\in\SG(X)$. The construction of the universal dynamical system and its continuous structure is based on joint work with {\sc J. Bellissard} and {\sc G. de Nittis} \cite{BeBeNi16}.

\medskip

Let $(X,G,\alpha)$ be a topological dynamical system. Since $X$ is compact, the space $\ks(X)$ of compact subsets of $X$ and the space $\cs(X)$ of closed subsets of $X$ agree. According to Definition~\ref{Chap2-Def-SpaDynSyst}, the space of dynamical systems $\SG(X)$ is defined by the set of closed and $G$-invariant subsets of $X$. This space is equipped with the subspace topology of $\ks(X)$, namely the Hausdorff-topology, c.f. Definition~\ref{Chap2-Def-SpaDynSyst}.

\begin{lemma}[\cite{BeBeNi16}]
\label{Chap4-Lem-XunivClosed}
Let $(X,G,\alpha)$ be a topological dynamical system with the associated space of dynamical subsystems $\SG(X)$. The space $\SG(X)\times X$ is equipped with the product topology. Then the set 
$$
\text{\gls{Xun}} \;
	:= \; \big\{(Y,y)\in \SG(X)\times X\;\big|\; y\in Y\big\}
$$
is a closed subset of $\SG(X)\times X$. In particular, the space $\Xun$ equipped with the induced topology of $\SG(X)\times X$ is second countable, compact and Hausdorff.
\end{lemma}

\begin{proof}
Let $(Y_n,y_n)\in\Xun,\; n\in\NM,$ be a convergent sequence to $(Y,y)\in\SG(X)\times X$. In detail, $(Y_n)_{n\in\NM}$ converges in the Hausdorff-topology to $Y\in\SG(X)$ and $y_n\in Y_n,\; n\in\NM$ converges in $X$ to $y$. It suffices to prove that $y\in Y$. This is proved in the claim (i) of the proof of Proposition~\ref{Chap2-Prop-SpaDynSyst}. Thus, $\Xun$ is a closed subset of $\SG(X)\times X$.

\vspace{.1cm}

The space $\SG(X)\times X$ is second countable, compact and Hausdorff, c.f. Definition~\ref{Chap2-Def-DynSyst} and Proposition~\ref{Chap2-Prop-SpaDynSyst}. Consequently, $\Xun$ is also a second countable, compact, Hausdorff space with respect to the induced topology of $\SG(X)\times X$. 
\end{proof}

\medskip

Since the group $G$ acts on each $Y\in\SG(X)$, it naturally acts on the so called {\em universal space $\Xun$}. 

\begin{lemma}[\cite{BeBeNi16}]
\label{Chap4-Lem-XunivDynSyst}
Let $(X,G,\alpha)$ be a topological dynamical system with the topological space $\Xun$. Then the map $\beta:G\times \Xun\to \Xun$ defined by
$$
\beta_g(Y,y):=(Y,\alpha_g(y))
$$
is a continuous action of $G$ on $\Xun$. In particular, $(\Xun,G,\beta)$ is a topological dynamical system.
\end{lemma}

\begin{proof}
Denote by $id:\SG(X)\to\SG(X)$ the identity map on $\SG(X)$. Then $\beta$ is equal to $id\times\alpha$. Since both maps are continuous actions of $G$, the map $\beta$ is a continuous action of $G$ as well, c.f. Definition~\ref{Chap2-Def-DynSyst}. Due to Lemma~\ref{Chap4-Lem-XunivClosed}, $\Xun$ is a second-countable, compact, Hausdorff space. The pair $(X,x)$ is an element of $\Xun$ for a $x\in X$. Hence, $\Xun$ is non-empty. Thus, $(\Xun,G,\beta)$ is a topological dynamical system in terms of Definition~\ref{Chap2-Def-DynSyst}.
\end{proof}

\begin{definition}[Universal dynamical system, \cite{BeBeNi16}]
Let $(X,G,\alpha)$ be a topological dynamical system. Then $(\Xun,G,\beta)$ is called the {\em universal dynamical system of $X$} associated with $(X,G,\alpha)$ and $\Xun$ is called the {\em universal space of $X$}.
\end{definition}

For symbolic dynamical systems, the universal space is totally disconnected.

\begin{corollary}[\cite{BeBeNi16}]
\label{Chap4-Cor-XunivTotDisc}
Let $\as$ be an alphabet and $G$ a discrete, countable group. Then the universal space $\asGun$ is totally disconnected.
\end{corollary}

\begin{proof}
According to Proposition~\ref{Chap2-Prop-FullShift} and Corollary~\ref{Chap2-Cor-SGTotaDisco} the space $\SG\big(\as^G\big)\times \as^G$ is totally disconnected. Consequently, the universal space $\asGun$ is totally disconnected as a closed subspace of $\SG\big(\as^G\big)\times \as^G$.
\end{proof}

\medskip

According to Definition~\ref{Chap2-Def-TransformationGroupGroupoid}, a topological dynamical system naturally induces the so called transformation group groupoid. 

\begin{definition}[Universal groupoid, \cite{BeBeNi16}]
\label{Chap4-Def-UnivGroupoid}
Let $(X,G,\alpha)$ be a topological dynami\-cal system. Then the transformation group groupoid $\text{\gls{Gun}}:=\Gun(X):=\Xun\rtimes_{\beta} G$ is called the {\em universal groupoid} associated with $(X,G,\alpha)$.
\end{definition}

According to Proposition~\ref{Chap2-Prop-TransGrGroupProp}, the unit space $\Gun^{(0)}$ is equal to the topological space $\Xun\times\{e\}\subsetneq\Gun^{(1)}$ where $e\in G$ is the unit. Thus, it immediately follows that the unit space $\Gun^{(0)}$ is homeomorphic to $\Xun$. For the sake of convenience, the unit space $\Gun^{(0)}$ is identified with the space $\Xun$ as long as there is no confusion. Denote by $\lambda$ the counting measure on a discrete, countable group $G$, namely $\lambda$ is the left-invariant Haar measure of $G$.

\begin{proposition}[\cite{BeBeNi16}]
\label{Chap4-Prop-PropUnivGroup}
Let $(X,G,\alpha)$ be a topological dynamical system where $G$ is a discrete, countable group with counting measure $\lambda$. Then $\Gun:=\Gun(X)$ is an \'etale groupoid with compact unit space $\Gamma^{(0)}=\Xun\times\{e\}$ and left-continuous Haar system $\nu^{(Y,y|e)}:=\delta_Y\times\delta_y\times\lambda,\; (Y,y|e)\in \Gamma^{(0)}$. The groupoid $\Gun$ is topologically amenable whenever $G$ is an amenable group. Furthermore, the equations $r(Y,y|g)=(Y,y|e)$ and $s(Y,y|g)=(Y,\alpha_{g^{-1}}(y)|e)$ are satisfied for $(Y,y|g)\in\Gun^{(1)}$.
\end{proposition}

\begin{proof}
Proposition~\ref{Chap2-Prop-Etale}~(b) implies that $\Gun$ is an \'etale groupoid. Due to Proposition~\ref{Chap2-Prop-TransGrGroupProp}, the unit space is equal to the compact set $\Xun\times\{e\}$ and so it is homeomorphic with $\Xun$, c.f. Proposition~\ref{Chap4-Lem-XunivClosed}. Then Proposition~\ref{Chap2-Prop-TransGroupAmen} implies that the measures 
$$
\nu^{(Y,y|e)} \; 
	:= \; \nu^{(Y,y)} \;
		:= \; \delta_Y\times\delta_y\times\lambda\; ,
			\qquad (Y,y|e)\in \Gun^{(0)}\; = \; \Xun\times\{e\}\,,
$$
define a left-continuous Haar system of the groupoid $\Gun$. Furthermore, $\Gun$ is topologically amenable whenever $G$ is amenable by Proposition~\ref{Chap2-Prop-TransGroupAmen}. The identities for the range and source map are proven in Proposition~\ref{Chap2-Prop-TransGrGroupProp}.
\end{proof}

\medskip

It turns out that the universal groupoid induces a continuous field of groupoids. The fiber groupoids are isomorphic to $Y\rtimes_\alpha G\,,\; Y\in\SG(X)$. 

\medskip

For simplification of the notation, the universal groupoid $\Gun(X)$ is denoted by $\Gamma$ whenever there is no confusion.

\begin{proposition}[\cite{BeBeNi16}]
\label{Chap4-Prop-UnivGroupContField}
Let $(X,G,\alpha)$ be a topological dynamical system where $G$ is a discrete, countable group with Haar measure $\lambda$ and $\Gamma:=\Gun:=\Gun(X)$ is the associated universal groupoid. Define the map $p:\Gamma^{(1)}\to\SG(X)$ by
$$
p(Y,y|g) \; 
	:= \; Y,\qquad (Y,y|g)\in\Gamma^{(1)} \, .
$$
Then $(\Gamma,\SG(X),p)$ is a continuous field of groupoids. Furthermore, for $Y\in\SG(X)$, the map $\Phi_Y:\Gamma_Y^{(1)}\to\Gamma^{(1)}(Y),\; (Y,y|g)\mapsto(y|g)\, ,$ defines a measure preserving isomorphism where $\Gamma(Y):=Y\rtimes_\alpha G$.
\end{proposition}

\begin{proof}
Proposition~\ref{Chap4-Prop-PropUnivGroup} implies that the universal groupoid $\Gamma$ is a topological groupoid with compact unit space and left-continuous Haar system 
$$
\nu^{(Y,y|e)} \; 
	:= \; \nu^{(Y,y)} \;
		:= \; \delta_Y\times\delta_y\times\lambda\; ,
			\qquad (Y,y|e)\in \Gamma^{(0)}\; = \; \Xun\times\{e\}\, .
$$
Furthermore, the range and the source map of $\Gamma$ leave the first component of an element $(Y,y|g)\in\Gamma^{(1)}$ invariant. Thus, the equations $p=p|_{\Gamma^{(0)}}\circ r=p|_{\Gamma^{(0)}}\circ s$ follow. The proof is organized as follows.
\begin{itemize}
\item[(i)] The map $p$ is surjective.
\item[(ii)] The map $p$ is open. 
\item[(iii)] The map $p$ is continuous. 
\item[(iv)] The map $\Phi_Y:\Gamma_Y^{(1)}\to\Gamma^{(1)}(Y)$ defines a groupoid isomorphism for all $Y\in\SG(X)$.
\item[(v)]  The map $\Phi_Y:\Gamma_Y^{(1)}\to\Gamma^{(1)}(Y)$ is measure preserving for all $Y\in\SG(X)$.
\end{itemize}

Assertion (i)-(iii) imply that $(\Gamma,\SG(X),p)$ is a continuous field of groupoids while (iv) and (v) yield that there is a measure preserving groupoid isomorphism between the fiber groupoid $\Gamma_Y$ and the transformation group groupoid $\Gamma(Y):=Y\rtimes_\alpha G$ for $Y\in\SG(X)$.

\vspace{.1cm}

(i): Let $Y\in\SG(X)$. Since $Y$ is a dynamical subsystem, it is non-empty and so there is a $y\in Y$. Then $(Y,y|e)$ is an element of $\Gamma^{(1)}$ where $e\in G$ is the neutral element. Hence, $p(Y,y,e)=Y$ holds. Hence, $p$ is surjective.

\vspace{.1cm}

(ii): It suffices to show that $p(\vs)$ is open for an element $\vs$ of the base for the topology of $\Gamma^{(1)}$. Consider $\vs:=(\us\times\os\times \{g\})\cap(\Xun\times G)$ where $\us:=\us(F,\os_1,\ldots,\os_n)\subseteq\SG(X)$ and $\os\subseteq X$ are open subsets and $g\in G$. Then the equation
$$
p(\vs) \;
	= \; \left\{ Y\in\SG(X)\;\left|\; 
		\begin{array}{c}
			\text{there exists a } y\in Y \text{ such}\\ 
			\text{that } (Y,y)\in \us\times\os
		\end{array}
		\right.
	 \right\}
$$
follows. It is immediate to check that $Y\in p(\vs)$ holds if and only if $Y\in\us$ and the intersection $Y\cap\os$ is non-empty. Thus, the identity $p\big(\vs\big) = \us\big(F,\{\os_1,\ldots,\os_n,\os\}\big)$ is derived where the set on the right hand side is open in $\SG(X)$. 

\vspace{.1cm}

(iii): Consider an open set $\us\subseteq\SG(X)$. Then the equation 
$$
p^{-1}(\us) \; 
	= \; (\us\times X\times G)\cap \Gamma^{(1)}
$$
follows immediately. The set on the right hand side is open in $\Gamma^{(1)}$ implying that $p$ is continuous.

\vspace{.1cm}

(iv): Let $Y\in\SG(X)$. Due to definition of $\Gamma_Y$ the equation $\Gamma^{(1)}_Y = \{Y\}\times Y\times G$ holds. The first component only plays the role of a label. Then the fact that the groupoid $\Gamma_Y$ arises as a subgroupoid of the transformation group groupoid $\Gamma:=\Xun\rtimes_\beta G$ leads to the desired result. 

\vspace{.1cm}

In detail, the groupoid $\Gamma_Y$ is equipped with the subspace topology of 
$$
\Gamma^{(1)} \;
	:= \; \Xun\times G\subsetneq \SG(X)\times X\times G\, .
$$ 
Thus, a base of the topology of $\Gamma^{(1)}_Y$ is given by the sets $\{Y\}\times\os\times U$ where $\os\subseteq Y$ and $U\subseteq G$ are open. Since $\Gamma^{(1)}(Y):=Y\times G$ is equipped with the product topology, the map $\Phi:\Gamma^{(1)}_Y\to\Gamma^{(1)}(Y),\; (Y,y|g)\in\Gamma^{(1)}_Y\mapsto (y|g)$ defines a homeomorphism.

\vspace{.1cm}

The equivalences
\begin{align*}
\big((Y,y|g),(Y,z|h)\big)\in\Gamma^{(2)}_Y 
	\;&\Longleftrightarrow\; \beta_{g^{-1}}(Y,y)=(Y,z)\\
	\;&\Longleftrightarrow\; \alpha_{g^{-1}}(y)=z\\
	\;&\Longleftrightarrow\; \big((y|g),(z|g)\big)\in\Gamma^{(2)}(Y)
\end{align*}
is derived by \nameref{(TG2)}. Then a short computation leads to
$$
\Phi\Big((Y,y|g) \circ (Y,z|h) \Big) \; 
	= \; (y|gh)
		= \; \Phi(Y,y|g) \circ \Phi(Y,z|h)\, ,
			\quad \big((Y,y|g),(Y,z|h)\big)\in\Gamma^{(2)}_Y\, .
$$
Thus, $\Phi$ is a groupoid isomorphism in terms of Definition~\ref{Chap4-Def-GroupoidIsomorphism}. 

\vspace{.1cm}

(v): The unit space of $\Gamma^{(0)}_Y$ is equal to $\{Y\}\times Y\times\{e\}$ and its left-continuous Haar system is the family $\nu^{(Y,y|e)}:=\delta_Y\times\delta_y\times\lambda\,,\; (Y,y|e)\in\Gamma^{(0)}_Y$ of Borel measures on $\Gamma^{(1)}_Y$. Furthermore, the left-continuous Haar system on $\Gamma(Y)$ is defined by the family of Borel measures $\mu^y:=\mu^{(y|e)}:=\delta_y\times\lambda$ for $(y|e)\in\Gamma^{(0)}(Y)=Y\times\{e\}$, c.f. Definition~\ref{Chap2-Prop-TransGroupAmen}. Then, for a Borel measurable set $\{Y\}\times F\subseteq\Gamma^{(1)}_Y$ and $(Y,y|e)\in\Gamma^{(0)}_Y$, the equations
$$
\nu^{(Y,y|e)}\big(\{Y\}\times F\big) \; 
	= \; \big(\delta_Y\times\delta_y\times\lambda \big)(\{Y\}\times F) \;
	= \; \big(\delta_y\times\lambda\big)(F) \;
	= \; \mu^{\Phi(Y,y|e)}\big(\Phi(\{Y\}\times F)\big)
$$
are derived implying that $\Phi$ is measure preserving.
\end{proof}

\medskip

Combining Proposition~\ref{Chap4-Prop-CAlgIsoGroupIso} and Proposition~\ref{Chap4-Prop-UnivGroupContField}, the following assertion is immediately deduced.

\begin{corollary}[\cite{BeBeNi16}]
\label{Chap4-Cor-IsometrIsomorRedCAlg}
Let $(X,G,\alpha)$ be a topological dynamical system where $G$ is a discrete, countable group with Haar measure $\lambda$. Then, for each $Y\in\SG(X)$, there exists an isometric, $\ast$-isomorphism $\Theta_Y:=\Theta:\CG^\ast_{red}(\Gamma_Y)\to\CG^\ast_{red}(\Gamma(Y))$ mapping $\Cc_c\big(\Gamma_Y^{(1)}\big)$ onto $\Cc_c\big(\Gamma^{(1)}(Y)\big)$.
\end{corollary}

\begin{proof}
According to Proposition~\ref{Chap4-Prop-UnivGroupContField}, there exists a measure preserving groupoid isomorphism between the groupoids $\Gamma_Y$ and $\Gamma(Y)$ for each $Y\in\SG(X)$. Thus, Proposition~\ref{Chap4-Prop-CAlgIsoGroupIso} leads to the desired result.
\end{proof}

\medskip

Theorem~\ref{Chap4-Theo-FieldGroupContCAlgebras}, Proposition~\ref{Chap4-Prop-UnivGroupContField} and Corollary~\ref{Chap4-Cor-IsometrIsomorRedCAlg} imply that the field of $C^\ast$-algebras $\CG^\ast(Y\rtimes_\alpha G)\, ,\; Y\in\SG(X)\, ,$ has a continuous structure whenever the associated groupoids are topologically amenable.

\begin{theorem}[\cite{BeBeNi16}]
\label{Chap4-Theo-TransGroupContFieldCAlg}
Let $(X,G,\alpha)$ be a topological dynamical system where $G$ is an amenable, discrete, countable group with counting measure $\lambda$. Then
$$
\Lambda \; 
	:= \; \Bigg\{
		\Big(\Theta_Y(\fz_Y)\Big)_{Y\in\SG(X)} \in\prod_{Y\in\SG(X)} \Cc_c\big(\Gamma^{(1)}(Y)\big) 
		\; \Bigg|\;  
		\fz\in\Cc_c\big(\Gun^{(1)}\big)
	\Bigg\}
$$
is a generating family of $\prod_{Y\in\SG(X)} \CG^\ast(\Gamma(Y))$. In particular, there exists a unique $\Upsilon\subseteq\prod_{Y\in\SG(X)} \CG^\ast(\Gamma(Y))$ such that $\Lambda\subseteq\Upsilon$ and $\big((\CG^\ast(\Gamma(Y)))_{Y\in\SG(X)}, \Upsilon \big)$ defines a continuous field of unital $C^\ast$-algebras.
\end{theorem}

\begin{proof}
The subset 
$$
\tilde{\Lambda} \;
	:= \;  \big\{(\gz_Y)_{Y\in\SG(X)}\; \big|\; \gz\in\Cc_c\big(\Gun^{(1)}\big) \big\} \;
		\subseteq \; \prod_{Y\in\SG(X)} \CG^\ast(\Gamma_Y)
$$
satisfies \nameref{(CFC1)}-\nameref{(CFC3)}. Hence, $\tilde{\Lambda}$ is a generating family , c.f. Proposition~\ref{Chap4-Prop-UnivGroupContField} and Theorem~\ref{Chap4-Theo-FieldGroupContCAlgebras}. For $\gz\in\Cc_c\big(\Gun^{(1)}\big)$, the restriction $\gz_Y$ to $\Gamma_Y^{(1)}$ is an element of $\Cc_c\big(\Gamma^{(1)}_Y\big)$. Furthermore, $\Theta_Y:\CG^\ast(\Gamma_Y)\to\CG^\ast(\Gamma(Y))$ maps $\Cc_c\big(\Gamma_Y^{(1)}\big)$ onto $\Cc_c\big(\Gamma^{(1)}(Y)\big)$. Thus, every vector field $(\fz_Y)_{Y\in\SG(X)}\in\Lambda$ satisfies $\fz_Y\in\Cc_c\big(\Gamma^{(1)}(Y)\big)$ for all $Y\in\SG(X)$. Since $\tilde{\Lambda}$ is a generating family and $\Theta_Y:\CG^\ast(\Gamma_Y)\to\CG^\ast(\Gamma(Y))\,,\; Y\in\SG(X)\,,$ are $\ast$-isomorphisms the desired result is derived.

\vspace{.1cm}

In detail, $\tilde{\Lambda}$ is a complex linear subspace that is closed under multiplication and involution. Thus, $\Lambda$ is so as well since $\Theta$ is a $\ast$-homomorphism, c.f. Corollary~\ref{Chap4-Cor-IsometrIsomorRedCAlg}. Then Lemma~\ref{Chap4-Lem-CAlgIsomorphismUnital} implies that $\Theta_Y$ maps the unit of $\CG^\ast(\Gamma_Y)$ to the unit of $\CG^\ast(\Gamma(Y))$. Hence, $\Lambda$ satisfies \nameref{(CFC1)}. Let $Y_0\in\SG(X)$. Due to Theorem~\ref{Chap4-Theo-FieldGroupContCAlgebras}, the set $\tilde{\Lambda}_{Y_0}:=\{\gz_{Y_0}\;|\; (\gz_Y)_{Y\in\SG(X)}\in \tilde{\Lambda}\}$ contains the set $\Cc_c\big(\Gamma^{(1)}_Y\big)$. Since $\Theta$ maps $\Cc_c\big(\Gamma^{(1)}_Y\big)$ onto $\Cc_c\big(\Gamma^{(1)}(Y)\big)$, this implies that 
$$
\Cc_c\big(\Gamma^{(1)}(Y_0)\big)\; 
	\subseteq \; \Lambda_{Y_0}\;
		:= \; \{\fz_{Y_0}\;|\; (\fz_Y)_{Y\in\SG(X)}\in \Lambda\}\, .
$$
Consequently, $\Lambda_{Y_0}$ is a dense subset of $\CG^\ast(\Gamma(Y))$. Thus, \nameref{(CFC2)} holds for the set $\Lambda$. Condition~\nameref{(CFC3)} of $\tilde{\Lambda}$ yields that $\SG(X)\ni Y\mapsto\|\gz_Y\|$ is continuous for every vector field $(\gz_Y)_{Y\in\SG(X)}\in\tilde{\Lambda}$. Then, for $(\fz_Y)_{Y\in\SG(X)}\in\Lambda$, the vector field $\big(\Theta_Y^{-1}(\fz_Y)\big)_{Y\in\SG(X)}$ is an element of $\tilde{\Lambda}$ and so the map
$$
\SG(X)\ni Y\;
	\mapsto\; \|\Theta_Y^{-1}(\fz_Y)\|\in[0,\infty)
$$
is continuous. Hence, the map $\SG(X)\ni Y\mapsto\|\fz_Y\|\in[0,\infty)$ is continuous as well for the vector field $(\fz_Y)_{Y\in\SG(X)}\in\Lambda$ since $\|\fz_Y\|=\|\Theta_Y^{-1}(\fz_Y)\|$ holds for $Y\in\SG(X)$, c.f. Corollary~\ref{Chap4-Cor-IsometrIsomorRedCAlg}. Consequently, $\Lambda$ satisfies also \nameref{(CFC3)}. Altogether, $\Lambda$ is a generating family of $\prod_{Y\in\SG(X)} \CG^\ast(\Gamma(Y))$. Then Proposition~\ref{Chap3-Prop-GenerFamiContField} leads to the desired result.
\end{proof}

\begin{remark}
\label{Chap4-Rem-TransGroupContFieldCAlg}
(i) Note that the condition that the group $G$ is amenable can be relaxed, c.f. Remark~\ref{Chap4-Rem-LRUseAmenable}. For instance, it suffices that the reduced and the full $C^\ast$-algebras of the fiber groupoids coincide, c.f. Remark~\ref{Chap4-Rem-TopoAmenableWeaken}~(i).

\vspace{.1cm}

(ii) For $Y_0\in\SG(X)$, the equality 
$$
\Cc_c\big(\Gamma^{(1)}(Y_0)\big)\; 
	= \; \Lambda_{Y_0}\;
		:= \; \{\fz_{Y_0}\;|\; (\fz_Y)_{Y\in\SG(X)}\in \Lambda\}\, .
$$
is actually valid by Lemma~\ref{Chap4-Lem-CompContExtGroupoid}, c.f. Proposition~\ref{Chap4-Prop-ExContSectContSpectr}~(a) below.
\end{remark}

According to Section~\ref{Chap3-Sect-ContFieldCAlg}, the spectra of a normal continuous vector field vary continuously in the parameter. In combination with the notion of the universal groupoid and its properties, this leads to the following fundamental result of this section.

\begin{corollary}[\cite{BeBeNi16}]
\label{Chap4-Cor-ContSpectrContSectUnivGroup}
Let $(X,G,\alpha)$ be a topological dynamical system where $G$ is an amen\-able, discrete, countable group with counting measure $\lambda$. Consider $(\fz_Y)_{Y\in\SG(X)}$ a normal continuous vector field of $\big((\CG^\ast(\Gamma(Y)))_{Y\in\SG(X)}, \Upsilon \big)$. Then the map 
$$
\SG(X)\ni Y\;
	\mapsto\;\sigma(\fz_Y)\in\ks(\CM)
$$
is continuous with respect to the Hausdorff metric.
\end{corollary}

\begin{proof}
This follows by Theorem~\ref{Chap3-Theo-ContFieldCALgContSpectr} and Theorem~\ref{Chap4-Theo-TransGroupContFieldCAlg}.
\end{proof}

\begin{remark}
\label{Chap4-Rem-ContSpectrContSectUnivGroup}
Note that each normal element of $\CG_{red}^\ast(\Gun)$ defines a normal continuous vector field by Proposition~\ref{Chap4-Prop-CAlgGUnivContVectField}.
\end{remark}

Let $(X,G,\alpha)$ be a topological dynamical system where $G$ is an amenable, discrete and countable group. For $Y_0\in\SG(X)$, a self-adjoint element $\fz$ of $\Cc_c(\Gamma(Y_0))$ defines a generalized Schr\"odinger operator of finite range. The section is finished by showing the existence of a continuous vector field $\gz:=(\gz_Y)_{Y\in\SG(X)}$ satisfying $\gz_{Y_0}=\fz$. Theorem~\ref{Chap3-Theo-ContFieldCALgContSpectr} leads to the continuity of the map $\SG(X)\ni Y\mapsto\sigma(\gz_Y)\in\ks(\CM)$ whenever $\gz$ is normal or self-adjoint.

\begin{proposition}[\cite{BeBeNi16}]
\label{Chap4-Prop-ExContSectContSpectr}
Let $(X,G,\alpha)$ be a topological dynamical system where $G$ is an amen\-able, discrete, countable group with counting measure $\lambda$. Consider the continuous field of unital $C^\ast$-algebras $\big((\CG^\ast(\Gamma(Y)))_{Y\in\SG(X)}, \Upsilon \big)$ defined in Theorem~\ref{Chap4-Theo-TransGroupContFieldCAlg} with generating family $\Lambda$. The following assertions hold for $Y_0\in\SG(X)$ and $\fz\in\Cc_c\big(\Gamma^{(1)}(Y_0)\big)$.
\begin{itemize}
\item[(a)] There exists a continuous vector field $(\gz_Y)_{Y\in\SG(X)}\in\Lambda$ such that $\gz_{Y_0}=\fz$. 
\item[(b)] If, additionally, $\fz$ is self-adjoint, there exists a self-adjoint continuous vector field $(\gz_Y)_{Y\in\SG(X)}\in\Lambda$ such that $\gz_{Y_0}=\fz$. 
\end{itemize}
\end{proposition}

\begin{proof}
(a): Let $\fz\in\Cc_c\big(\Gamma^{(1)}(Y_0)\big)$ for $Y_0\in\SG(X)$. According to Corollary~\ref{Chap4-Cor-IsometrIsomorRedCAlg}, $\Theta_{Y_0}^{-1}(\fz)$ is a continuous function defined on $\Gamma^{(1)}_Y$ with compact support. Thus, Lemma~\ref{Chap4-Lem-CompContExtGroupoid} implies the existence of a $\hz\in\Cc_c\big(\Gun^{(1)}\big)$ such that the restriction $\hz_{Y_0}$ of $\hz$ to $\Gamma^{(1)}_{Y_0}$ coincides with $\Theta_{Y_0}^{-1}(\fz)$. Define $\gz_Y:=\Theta_Y(\hz_Y)$ for $Y\in\SG(X)$. Then $\gz_{Y_0}=\fz$ holds and $(\gz_Y)_{Y\in\SG(X)}\in\Lambda$ is a continuous vector field, c.f. Theorem~\ref{Chap4-Theo-TransGroupContFieldCAlg}.

\vspace{.1cm}

(b): Let $Y_0\in\SG(X)$ and consider a self-adjoint $\fz\in\Cc_c\big(\Gamma^{(1)}(Y_0)\big)$. For $Y\in\SG(X)$, the unit of the $C^\ast$-algebra $\CG^\ast(\Gamma(Y))$ is denoted by $\mathpzc{1}_Y$. Since $\fz$ is self-adjoint, the spectrum $\sigma(\fz)$ is contained in $\RM$ and it is compact. Thus, there exists an $m\geq 0$ such that $\sigma(\fz+m\cdot \mathpzc{1}_{Y_0})\subseteq[0,\infty)$, i.e., $\fz+m\cdot\mathpzc{1}$ is a non-negative element of the $C^\ast$-algebra $\CG^\ast(\Gamma(Y_0))$, c.f. \cite[Section~2.2]{Murphy90}. For instance, choose $m$ to be the norm $\|\fz\|$. Recall that the unit $\mathpzc{1}_{Y_0}$ of $\CG^\ast(\Gamma(Y_0))$ is the characteristic function $\chi_{\Gamma^{(0)}(Y_0)}\in\Cc_c\big(\Gamma^{(1)}(Y_0)\big)$, c.f. Theorem~\ref{Chap2-Theo-UnitGroupoidCalgebra}. Then $\fz+m\cdot \mathpzc{1}_{Y_0}$ is also an element of $\Cc_c\big(\Gamma^{(1)}(Y_0)\big)$, c.f. Proposition~\ref{Chap2-Prop-InvolutiveGroupoidAlgebra}. Furthermore, $\fz+m\cdot \mathpzc{1}_{Y_0}$ is self-adjoint since $\fz$ and $\mathpzc{1}_{Y_0}$ are self-adjoint and $m\in\RM$. By (a), there exists a continuous vector field $(\hz_Y)_{Y\in\SG(X)}\in\Lambda$ such that $\hz_{Y_0}=\fz+m\cdot \mathpzc{1}_{Y_0}$. For $Y\in\SG(X)$, the convolution $\hz_Y^\ast\star\hz_Y$ satisfies
$$
\big(\hz_Y^\ast\star\hz_Y\big)^\ast \;
	= \; \hz_Y^\ast\star \big(\hz_Y^\ast\big)^\ast \;
		=\; \hz_Y^\ast\star\hz_Y\, .
$$
Hence, the vector field $(\hz_Y^\ast\star\hz_Y)_{Y\in\SG(X)}$ is self-adjoint and by construction it is non-negative, c.f. \cite[Theorem~2.2.4]{Murphy90}. Thus, there exists a unique non-negative self-adjoint vector field $(\hz'_Y)_{Y\in\SG(X)}$ such that $\hz'_Y\star\hz'_Y=\hz_Y^\ast\star\hz_Y$ for all $Y\in\SG(X)$. In detail, the square root $\phi:\CM\to\CM\,,\; z\mapsto\sqrt{|z|}$ is a continuous map. Then $\hz'_Y$ is equal to $\phi(\hz_Y^\ast\star\hz_Y)$ for $Y\in\SG(X)$ by the functional calculus. Thus, $\hz'_{Y_0}$ is equal to $\fz+m\cdot \mathpzc{1}_{Y_0}$ as $\hz_{Y_0}=\fz+m\cdot \mathpzc{1}_{Y_0}$ is self-adjoint and non-negative. 

\vspace{.1cm}

Since $\Lambda$ is a generating family and $\hz_Y\in\Cc_c\big(\Gamma^{(1)}(Y)\big)$, the vector field $(\hz_Y^\ast\star\hz_Y)_{Y\in\SG(X)}$ is an element of $\Lambda$. Together with Proposition~\ref{Chap3-Prop-ContFuncContVectFie}, this implies that $(\hz'_Y)_{Y\in\SG(X)}$ is a continuous vector field as $\phi$ is a continuous function and $h'_Y=\phi\big(\hz_Y^\ast\star\hz_Y\big)$.

\vspace{.1cm}

Now, define the self-adjoint vector field $\gz:=(\gz_Y)_{Y\in\SG(X)}$ by $\gz_Y:=\hz'_Y-m\cdot\mathpzc{1}_Y$ for $Y\in\SG(X)$. According to Lemma~\ref{Chap3-Lem-PolContVectField}, $\gz$ is an element of $\Lambda$ and it satisfies $\gz_{Y_0}=\fz$. 
\end{proof}

\medskip

For a symbolic dynamical system $(\as^G,G,\alpha)$, the assertion of Proposition~\ref{Chap4-Prop-ExContSectContSpectr} is made more explicitly for pattern equivariant Schr\"odinger operators associated with subshifts $\Xi\in\SG\big(\as^G\big)$. These operators are defined in Definition~\ref{Chap2-Def-SchrOp-l2(G)} and they were studied intensively in the literature, c.f. Section~\ref{Chap2-Sect-ExampleAsZMSchrOp} and Section~\ref{Chap2-Sect-SchrOp}. It is proven in Theorem~\ref{Chap2-Theo-PESchrOpFinRang} that the pattern equivariant Schr\"odinger operators arise exactly by the self-adjoint elements of $\PG^\ast(\Gamma(\Xi))$ which is a dense $\ast$-subalgebra of $\CG^\ast_{red}(\Gamma(\Xi))$ and $\CG^\ast_{full}(\Gamma(\Xi))$ by Theorem~\ref{Chap2-Theo-PaEqDense}.

\begin{corollary}[\cite{BeBeNi16}]
\label{Chap4-Cor-ExContSectContSpectrSymbDynSyst}
Let $\as$ be a finite alphabet and $G$ be an amenable, discrete, countable group with counting measure $\lambda$. Let $K\in\ks(G)$ be a finite set with pattern equivariant functions $p_h:\as^G\to\CM\,,\; h\in K\,,$ and $p_e:\as^G\to\RM$. Consider the family of bounded self-adjoint $H_\xi:\ell^2(G)\to\ell^2(G)\,,\; \xi\in\Xi\,,$ defined by
$$
(H_\xi\psi)(g) \; 
	:= \; \left( 
			\sum\limits_{h\in K} p_h\big(\alpha_{g^{-1}}(\xi)\big) \cdot \psi(g\,h^{-1}) + \overline{p_h\big(\alpha_{(gh)^{-1}}(\xi)\big)} \cdot \psi(gh)
		\right)
		+ p_e\big(\alpha_{g^{-1}}(\xi)\big) \cdot \psi(g)
$$
where $\psi\in\ell^2(G)$ and $g\in G$. For $\Xi\in\SG\big(\as^G\big)$, $H_\Xi:=\big(H_\xi\big)_{\xi\in\Xi}$ is the associated pattern equivariant Schr\"odinger operator with $\Xi$. 

\vspace{.1cm}

Then the map 
$$
\SG\big(\as^G\big)\ni \Xi \;
	\mapsto\;\sigma(H_\Xi)\in\ks(\RM)
$$ 
is continuous with respect to the Hausdorff metric on $\ks(\RM)$.
\end{corollary}

\begin{proof}
Consider the continuous field of unital $C^\ast$-algebras $\big(\big(\CG^\ast(\Gamma(\Xi))\big)_{\Xi\in\SG(\as^G)}, \Upsilon \big)$ defined in Theorem~\ref{Chap4-Theo-TransGroupContFieldCAlg} with generating family $\Lambda:=\Cc_c\big(\Gun^{(1)}(\as^G)\big)$. It suffices to prove that there exists a continuous vector field $\fz\in\PG^\ast\big(\Gun\big(\as^G\big)\big)$ such that $\pi^\xi(\fz)=H_\xi$ for $\xi\in\as^G$. In this case the continuity of the spectra follows by Corollary~\ref{Chap4-Cor-ContSpectrContSectUnivGroup}.

\vspace{.1cm}

This proof for the existence of such an $\fz$ follows the lines of Theorem~\ref{Chap2-Theo-PESchrOpFinRang}~(b). Define $\fz:\Gun^{(1)}(\as^G) \to\CM$ by
$$
\fz \;
	:= \; \left( \sum_{h\in K} \fz_h + \fz_h^\ast\right) +  \fz_e
$$
where the maps $\fz_h\in\Cc_c\big(\Gun^{(1)}(\as^G)\big)$ for $h\in K\cup\{e\}$ are defined by
$$
\fz_h(\Xi,\xi|g) \; 
	:= \; \delta_{h^{-1}}(g)\cdot p_h(\xi)
		\, , \qquad (\Xi,\xi|g)\in\Gun^{(1)}(\as^G)\, .
$$
Since $K$ is finite, the set $\Xun\times K\cup K^{-1}\subseteq\Gun^{(1)}$ is compact, c.f. Lemma~\ref{Chap4-Lem-XunivClosed}. Thus, the map $\fz$ is compactly supported as $\supp(\fz)\subseteq \Xun\times K\cup K^{-1}$. It is straight forward to check that $\fz$ is also continuous. According to Proposition~\ref{Chap2-Prop-CharPattEqFunc}, the finite sum of pattern equivariant functions is again pattern equivariant. Hence, $\fz$ is an element of $\PG^\ast\big(\Gun\big(\as^G\big)\big)\subseteq\Cc_c\big(\Gun^{(1)}\big)$ as the maps $p_h\,,\; h\in K\cup\{e\}\,,$ are pattern equivariant. Since $p_e$ is real-valued, $\fz$ is self-adjoint. Along the lines of Theorem~\ref{Chap2-Theo-PESchrOpFinRang}~(b), it follows that $\pi^{(\Xi,\xi)}(\fz)=H_\xi$ for all $(\Xi,\xi)\in \Xun$. Then Proposition~\ref{Chap2-Prop-CovFamOp-Spect}~(a) implies 
$$
\sigma(\fz_{\,\Xi}) \; 
	= \; \sigma(H_\Xi) \;
	= \; \bigcup_{\xi\in\Xi} \sigma\big(\pi^{(\Xi,\xi)}(\fz)\big) 		
$$
for $\Xi\in\SG\big(\as^G\big)$. Since $(\fz_{\,\Xi})_{\Xi\in\SG(\as^G)}$ is a self-adjoint continuous vector field, the map 
$$
\SG\big(\as^G\big)\ni \Xi\longmapsto\sigma(\fz_{\,\Xi}) \; 
	= \; \sigma(H_\Xi)\in\ks(\RM)
$$
is continuous with respect to the Hausdorff metric on $\ks(\RM)$ by Corollary~\ref{Chap4-Cor-ContSpectrContSectUnivGroup}.
\end{proof}

\begin{remark}
\label{Chap4-Rem-ExContSectContSpectrSymbDynSyst}
Clearly, the assertion of Corollary~\ref{Chap4-Cor-ExContSectContSpectrSymbDynSyst} holds also for all normal elements of the pattern equivariant algebra $\PG^\ast\big(\Gun\big(\as^G\big)\big)$. 
\end{remark}

\section[Continuous behavior of the spectra for generalized Schr\"o\-dinger operators]{Characterization of the continuous behavior of the spectra for generalized Schr\"odinger operators}
\label{Chap4-Sect-ContSpecGenSchrOp}

In the previous section, it was shown that the $C^\ast$-algebras associated with every $Y\in\SG(X)$ naturally define a continuous field leading to the continuity of the spectra of generalized Schr\"odinger operators. The continuity of the spectra was obtained by using the underlying groupoid structure. This is used in the following to characterize the convergence of the spectra for all generalized Schr\"odinger operators by the convergence of the associated dynamical subsystems. In this context, it also turns out that monotone convergence of the dynamical subsystems leads to the monotone convergence of the spectra. Furthermore, the strong continuity of the operators is verified. Finally, Example~\ref{Chap4-Ex-KeineNormKonv} shows that these operators do not converge in operator norm, in general.

\medskip

This section is divided up into two parts. First the case of a dynamical system $(X,G,\alpha)$ is treated. Secondly, the case of symbolic dynamical systems $(\as^G,G,\alpha)$ is discussed which is more explicit than the general case of dynamical systems.

\medskip

According to Section~\ref{Chap4-Sect-UnivGroupDynSyst}, the universal groupoid associated with a dynamical system defines a continuous field of groupoids. Due to \cite{LaRa99}, this yields to the continuity of the associated field of $C^\ast$-algebra. Together with Theorem~\ref{Chap3-Theo-ContFieldCALgContSpectr}, the continuity of the spectra corresponding to a continuous vector field is derived, c.f. Corollary~\ref{Chap4-Cor-ContSpectrContSectUnivGroup}. In this section, this chain of assertions is closed in spirit of Figure~\ref{Chap4-Fig-Concept}. More precisely, the convergence of the spectra of all generalized Schr\"odinger operators delivers the convergence of the corresponding dynamical subsystems $(X_n)_{n\in\NM}$ to $X_\infty\in\SG(X)$. The main idea is based on Urysohn's Lemma meaning that the topology can be described by continuous functions. Since elements of the groupoid $C^\ast$-algebra are defined via continuous functions on the groupoid, the convergence of the dynamical subsystems is deduced.

\medskip

In general, this can be interpreted as follows. The map $\imath:\oNM\to\SG(X)\,,\; n\mapsto X_n\,,$ induces a new groupoid $\imath^\ast\Gamma$ which is called the {\em pullback groupoid} where $\oNM:=\NM\cup\{\infty\}$ is the one point compactification of $\NM$. The groupoid $\imath^\ast\Gamma$ itself defines a continuous field of groupoids with $\imath^\ast p:\imath^\ast\Gamma\to\oNM$. Altogether, the following assertions are equivalent.
\begin{itemize}
\item[(i)] The sequence of dynamical subsystems $(X_n)_{n\in\NM}$ converges to $X_\infty\in\SG(X)$. 
\item[(ii)] The triple $(\imath^\ast\Gamma,\oNM,\imath^\ast p)$ is a continuous field of groupoids.
\item[(iii)] The field $(\CG^\ast(\imath^\ast\Gamma_n)_{n\in\oNM}$ of $C^\ast$-algebras is equipped with a continuous structure.
\item[(iv)] The map $n\in\oNM\mapsto \sigma\big(\gz_{X_n}\big)\in\ks(\CM)$ is continuous for all generalized Schr\"odinger operators arising by $\gz\in\Cc_c(X\rtimes_\alpha G)$.
\end{itemize}
Only the equivalence of (i) and (iv) is proven here. The construction of the continuous field of groupoids $(\imath^\ast\Gamma,\oNM,\imath^\ast p)$ follows the lines of the definition of the universal groupoid by adding a label $n\in\oNM$. Then the general theory of Section~\ref{Chap4-Sect-ContFieldGroupCAlg} and Theorem~\ref{Chap3-Theo-ContFieldCALgContSpectr} leads to the desired equivalences.

\subsection{General dynamical systems}
\label{Chap4-Ssect-GenDynSyst}

Let $(X,G,\alpha)$ be a topological dynamical system where $G$ is an amenable, discrete, countable group with counting measure $\lambda$. For simplification of the notation, the universal groupoid $\Gun$ is denoted by $\Gamma$ in this section whenever there is no confusion.

\begin{theorem}[Continuous behavior of the spectra and dynamical systems, \cite{BeBeNi16}]
\label{Chap4-Theo-CharDynSystConvSpectr}
Consider a topological dynamical system $(X,G,\alpha)$ where $G$ is an amenable, discrete, countable group with counting measure $\lambda$. Furthermore, let $\Gamma:=\Gun(X)$ be the universal groupoid associated with $(X,G,\alpha)$. Consider a dynamical subsystem $X_n\in\SG(X)$ for each $n\in\oNM=\NM\cup\{\infty\}$. Then the following assertions are equivalent.
\begin{itemize}
\item[(i)] The sequence $(X_n)_{n\in\NM}$ converges to $X_\infty\in\SG(X)$ in the Hausdorff-topology of $\SG(X)$.
\item[(ii)] For all self-adjoint (normal) $\fz\in\Cc_c\big(\Gun^{(1)}\big)$, the map 
$$
\Sigma_\fz:\oNM\to \ks(\CM)\,,\;
	n\;\mapsto\; \sigma\big(\fz_{X_n}\big)\,,
$$ 
is continuous with respect to the Hausdorff metric $\ks(\CM)$, i.e., $\lim_{n\to\infty}\sigma(\fz_{X_n})=\sigma(\fz_{X_\infty})$ holds.
\end{itemize}
\end{theorem}

\begin{proof}
(i)$\Rightarrow$(ii): Let $\fz\in\Cc_c\big(\Gun^{(1)}\big)$ be self-adjoint or normal. Consider the associated spectral map $\Sigma_\fz:\SG(X)\to\ks(\CM)\,,\; Y\mapsto\sigma(\fz_Y)$. Due to Corollary~\ref{Chap4-Cor-ContSpectrContSectUnivGroup} the map $\Sigma_\fz$ is continuous. The map $\imath:\oNM\to\SG(X)\,,\; \imath(n):=X_n\,,\; n\in\oNM\,,$ is also continuous by (i). Thus, $\oNM\ni n\mapsto\sigma(\fz_{X_n})\in\ks(\CM)$ is continuous since it is equal to the composition $\Sigma_\fz\circ \imath$ of continuous maps.

\vspace{.1cm}

(ii)$\Rightarrow$(i): Assume that $(X_n)_{n\in\NM}$ does not converge to $X_\infty\in\SG(X)$. Since $\SG(X)$ is compact, there is a subsequence $(X_{n_k})_{k\in\NM}$ such that $\lim_{k\to\infty} X_{n_k}=:Y$ exists and $Y\neq X_\infty$. There are two cases, namely (a) $X_\infty\setminus Y\neq\emptyset$ and (b) $Y\setminus X_\infty\neq\emptyset$. Both cases lead as follows to a contradiction. 

\vspace{.1cm}

(a) If $X_\infty\setminus Y$ is non-empty, consider an $x\in X_\infty\setminus Y$. The space $X$ is a second countable, locally compact and Hausdorff, c.f. Definition~\ref{Chap2-Def-DynSyst}. Thus, $X$ is normal and so the Lemma of Urysohn applies. Consequently, there exists a continuous function $\phi:X\to[0,1]$ with compact support such that $\phi(x)=1$ and $\phi(z)=0$ for all $z\in Y$, c.f. Proposition~\ref{App1-Prop-LocCompHausSecCountImplNormal} and Proposition~\ref{App1-Prop-LemmaUrysohn}. Then $\fz:\Gun^{(1)}\to\RM\,,\; (Z,z|g)\mapsto\phi(z)\cdot\delta_e(g)$ defines a continuous map with compact support. Furthermore, $\fz\in\Cc_c\big(\Gun^{(1)}\big)$ is self-adjoint since it is supported on $\supp(\phi)\times\{e\}$ and $\phi$ is real-valued. Set $z_\infty:=(X_\infty,x)\in \{X_\infty\}\times X_\infty$. Then the estimates
\begin{align*}
\|\fz_{X_\infty}\|
	&\overset{\qquad\quad\;\;}{:=} \; \sup\Big\{ \|\pi^z(\fz)\| \;\big|\; z \, :=\, (X_\infty,y)\in \{X_\infty\}\times X_\infty \Big\} \\
		&\overset{\qquad\quad\;\;}{\geq} \; \sup\big\{ \left\|\pi^{z_\infty}(\fz)\psi\right\|_{\ell^2(G)} \; \big|\; \psi\in\ell^2(G) \text{ with } \|\psi\|=1 \big\}\\
			&\overset{\qquad\quad\;\;}{\geq} \; \left\|\pi^{z_\infty}(\fz) \delta_e\right\|_{\ell^2(G)}\\
				&\overset{\qquad\quad\;\;}{\geq} \; \left|\big(\pi^{z_\infty}(\fz) \delta_e\big)(e)\right|\\
					&\overset{\text{Prop.~\ref{Chap2-Prop-LeftRegReprTransGrouGrou}}}{=} \; |\phi(x)|\\
						&\overset{\qquad\quad\;\;}{=} \; 1
\end{align*}
are derived. Note that $e\in G$ is the neutral element and the Kronecker delta function $\delta_e$ is an element of $\ell^2(G)$ with $\ell^2$-norm equals to $1$. Since $\fz_{X_\infty}$ is self-adjoint, the spectral radius agrees with the norm $\|\fz_{X_\infty}\|$. Consequently, there exists a $\lambda\in\sigma\big(\fz_{X_\infty}\big)$ with $|\lambda|\geq 1$. On the other hand, by convergence of the sequence $(X_{n_k})_{k\in\NM}$ to $Y\in\SG(X)$ with respect to the Hausdorff-topology, there is a $k_0\in\NM$ such that 
$$
X_{n_k}\in
	\us(\supp(\phi),\{X\}) \;
		= \; \big\{Z\in\SG(X)\;\big|\; Z\cap\supp(\phi)=\emptyset\,,\; Z\cap X\neq\emptyset\big\}\,,
			\quad k\geq k_0\, .
$$
Note that $\us(\supp(\phi),\{X\})$ is an open neighborhood of $Y$ in the Hausdorff-topology of $\SG(X)$ as $\supp(\phi)\cap Y=\emptyset$ by construction. Thus, $\fz_{X_{n_k}}\equiv 0$ holds for all $k\geq k_0$ since the intersection $X_{n_k}\cap\supp(\phi)$ is empty for $k\geq k_0$. This implies that $\sigma\big(\fz_{X_{n_k}}\big)=\{0\}$ for each $k\geq k_0$. Altogether, the inequality
$$
d_H\left(\sigma\big(\fz_{X_\infty}\big),\sigma\left(\fz_{X_{n_k}}\right)\right) \;
	 \geq \; 1\,,
	 \qquad k\geq k_0 \,,
$$
is derived. Hence, the map $\oNM\ni n\mapsto \sigma(\fz_{X_n})\in\ks(\CM)$ is not continuous with respect to the Hausdorff metric. This is a contradiction as by (ii) every self-adjoint $\fz\in\Cc_c\big(\Gun^{(1)}(X)\big)$ satisfies that $\oNM\ni n\mapsto \sigma(\fz_{X_n})\in\ks(\CM)$ is continuous.

\vspace{.1cm}

(b) If $Y\setminus X_\infty$ is non-empty, consider an $x\in Y\setminus X_\infty$. Like in the case (a), define a self-adjoint element $\fz\in\Cc_c\big(\Gun^{(1)}(X)\big)$ by $\fz(Z,z|g)=\phi(z)\cdot\delta_e(g)$ for $(Z,z|g)\in\Gun(X)$ where $\phi:X\to[0,1]$ is a continuous map satisfying $\supp(\phi)\cap X_\infty = \emptyset$ and $\phi(x)=1$, c.f. Proposition~\ref{App1-Prop-LocCompHausSecCountImplNormal} and Proposition~\ref{App1-Prop-LemmaUrysohn}. According to claim (ii) of Proposition~\ref{Chap2-Prop-SpaDynSyst}, there exists an $x_{n_k}\in X_{n_k}$ for each $k\in\NM$ satisfying $\lim_{k\to\infty}x_{n_k}=x$. By continuity of $\phi$ and $\phi(x)=1$, there is a $k_0\in\NM$ such that $\phi(x_{n_k})\geq 1/2$ for all $k\geq k_0$. Similar to the considerations in (a), the estimates $\|\fz_{X_{n_k}}\| \geq \big|\phi\big(x_{n_k}\big)\big| \geq 1/2$ is deduced for $k\geq k_0$. Since $\fz_{X_n}$ is self-adjoint for every $n\in\NM$, the spectral radius agrees with the norm $\|\fz_{X_n}\|$. Consequently, there exists a $\lambda_k\in\sigma(X_{n_k})$ with $|\lambda_k|\geq 1/2$ for all $k\geq k_0$. Furthermore, $\fz_{X_\infty}\equiv 0$ holds as $\supp(\phi)\cap X_\infty=\emptyset$. This implies $\sigma(\fz_{X_\infty})=\{0\}$. Hence, the Hausdorff distance of the associated spectra satisfies
$$
d_H\left(\sigma\big(\fz_{X_\infty}\big),\sigma\left(\fz_{X_{n_k}}\right)\right) \; 
	\geq \; \frac{1}{2}\,,
		\qquad k\geq k_0 \, .
$$
Thus, the map $\oNM\ni n\mapsto \sigma(\fz_{X_n})\in\ks(\CM)$ is not continuous with respect to the Hausdorff metric. Similar to (a), this is a contradiction to (ii).
\end{proof}

\begin{remark}
\label{Chap4-Rem-CharDynSystConvSpectr}
The condition that $G$ is amenable is not necessary for the assertion in Theorem~\ref{Chap4-Theo-CharDynSystConvSpectr}. In detail, it suffices that $\Gamma(Y)$ is topologically amenable for each $Y\in\SG(X)$. Actually, it suffices that the associated reduced and full groupoid $C^\ast$-algebras of $\Gamma(Y)$ coincide for each $Y\in\SG(X)$. According to Remark~\ref{Chap4-Rem-TopoAmenableWeaken}~(i), there exist topological groupoids $\Gamma$ such that $\CG^\ast_{red}(\Gamma)=\CG^\ast_{full}(\Gamma)$ while $\Gamma$ is not topologically amenable. 

\vspace{.1cm}

According to Proposition~\ref{Chap4-Prop-CAlgGUnivContVectField}, one can replace $\fz\in\Cc_c\big(\Gun^{(1)}\big)$ by all self-adjoint (normal) $\fz\in\CG^\ast_{red}(\Gun)$ in Theorem~\ref{Chap4-Theo-CharDynSystConvSpectr}~(ii).
\end{remark}

The main difference between the following corollary and Theorem~\ref{Chap4-Theo-CharDynSystConvSpectr} is that all the dynamical subsystems are assumed to be minimal.

\begin{corollary}[Continuous behavior of the spectra and minimal dynamical subsystems, \cite{BeBeNi16}]
\label{Chap4-Cor-CharDynSystConvSpect}
Let $(X,G,\alpha)$ be a topological dynamical system where $G$ is an amenable, discrete, countable group with counting measure $\lambda$. The associated universal groupoid is denoted by $\Gamma:=\Gun$. Consider a minimal dynamical subsystem $X_n\in\SG(X)$ for each $n\in\oNM=\NM\cup\{\infty\}$. Then the following assertions are equivalent.
\begin{itemize}
\item[(i)] The sequence $(X_n)_{n\in\NM}$ converges to $X_\infty\in\SG(X)$ in the Hausdorff-topology of $\SG(X)$.
\item[(ii)] For all self-adjoint (normal) $\fz\in\Cc_c\big(\Gun^{(1)}\big)$ and all $x_n\in X_n\,,\; n\in\oNM$, the map 
$$
\Sigma:\oNM\to \ks(\CM)\,,\;
	n\; \mapsto \; \sigma\big(\pi^{x_n}(\fz)\big) \;
		= \; \sigma_{ess}\big(\pi^{x_n}(\fz)\big)
$$
is continuous with respect to the Hausdorff metric on $\ks(\CM)$, i.e.,
$$
\lim_{n\to\infty}\sigma\big(\pi^{x_n}(\fz)\big)\; 
	=\; \lim_{n\to\infty} \sigma_{ess}\big(\pi^{x_n}(\fz)\big)
	=\; \sigma\big(\pi^{x_\infty}(\fz)\big)
	=\; \sigma_{ess}\big(\pi^{x_\infty}(\fz)\big)
$$
holds.
\end{itemize}
\end{corollary}

\begin{proof}
The equivalence is a consequence of Theorem~\ref{Chap4-Theo-CharDynSystConvSpectr} and Theorem~\ref{Chap2-Theo-ConstSpectrMinimal}. The desired minimality of $X_n$ for $n\in\oNM$ implies that $\sigma(\fz_{X_n})=\sigma\big(\pi^{x_n}(\fz)\big)=\sigma_{ess}\big(\pi^{x_n}(\fz)\big)$ for all $n\in\oNM$. Then (ii) asserts $\lim_{n\to\infty}\sigma\big(\pi^{x_n}(\fz)\big)=\sigma\big(\pi^{x_\infty}(\fz)\big)=\sigma_{ess}\big(\pi^{x_\infty}(\fz)\big)$.
\end{proof}

\medskip

The previous assertions provide a characterization of the \pt-continuity for generalized Schr\"odinger operators, namely the operators are continuous with respect to the \pt-topology on $\Ll(\ell^2(G))$. In Section~\ref{Chap3-Sect-OpTop}, the relation to other topologies on $\Ll(\ell^2(G))$ was analyzed. It was shown that there is no relations between the strong operator topology and the \pt-topology, c.f. Proposition~\ref{Chap3-Prop-P2NotImplStrong} and Proposition~\ref{Chap3-Prop-StrongNotImplP2}. Due to Propo\-sition~\ref{Chap2-Prop-CovFamOp-Spect}~(c), the strong continuity holds also for generalized Schr\"odinger operators.

\begin{proposition}[Strong convergence and convergence of the spectra, \cite{BeBeNi16}]
\label{Chap4-Prop-StrongConv+ConvSpectr}
Let $(X,G,\alpha)$ be a topological dynamical system where $G$ is an amenable, discrete, countable group with counting measure $\lambda$. Furthermore, let $\Gamma:=\Gun$ be the universal groupoid associated with $(X,G,\alpha)$. Consider a sequence of minimal dynamical subsystem $X_n\in\SG(X)\,,\; n\in\NM\,$ tending to the minimal $X_{\infty}\in\SG(X)$. Then the following assertions hold for each self-adjoint (normal) $\fz\in\CG^\ast(\Gun)$ and all sequence $x_n\in X_n\,,\; n\in\NM\,,$ with $\lim_{n\to\infty}x_n=:x_\infty\in X_\infty$.
\begin{itemize}
\item[(a)] The map $\oNM\ni n\mapsto \pi^{x_n}(\fz)\in\Ll(\ell^2(G))$ is continuous with respect to the strong operator topology on $\Ll(\ell^2(G))$.
\item[(b)] The equality $\lim_{n\to\infty}\sigma\big(\pi^{x_n}(\fz)\big) = \sigma\big(\pi^{x_\infty}(\fz)\big)=\sigma_{ess}\big(\pi^{x_\infty}(\fz)\big)\in\ks(\CM)$ holds where the limit is taken with respect to the Hausdorff metric on $\ks(\CM)$.
\end{itemize}
\end{proposition}

\begin{proof}
According to Proposition~\ref{Chap4-Prop-PropUnivGroup}, the universal groupoid $\Gamma$ is topologically amenable and \'etale. With this at hand, Proposition~\ref{Chap2-Prop-CovFamOp-Spect}~(c) leads to the statement (a) since $\fz$ is a self-adjoint (normal) element of the groupoid $C^\ast$-algebra $\CG^\ast(\Gun)$. Finally, assertion (b) follows by Corollary~\ref{Chap4-Cor-CharDynSystConvSpect}.
\end{proof}

\begin{remark}
\label{Chap4-Rem-GeneralSchrOpStrongCont}
The strong operator continuity of the generalized Schr\"odinger operators does not rely on the continuous behavior of the spectra. It is obtained since the operators are elements of a groupoid $C^\ast$-algebra and so they arise as continuous functions on a topological space (groupoid).
\end{remark}

\begin{proposition}[Monotone convergence, \cite{BeBeNi16}]
\label{Chap4-Prop-MonotonConvSpectr}
Let $(X,G,\alpha)$ be a topological dynamical system where $G$ is an amen\-able, countable, discrete group. Consider a convergent sequence $(Y_n)_{n\in\NM}$ with $\lim_{n\to\infty}Y_n=:Y_\infty\in\SG(X)$ such that $Y_{n+1}\subseteq Y_n$ holds for all $n\in\NM$. Let $\fz\in\Cc_c\big(\Gun^{(1)}\big)$ be a normal element and consider the generalized Schr\"odinger operators $H_{Y_n}:=(H_x)_{x\in Y_n}$ where $H_x:=\pi^x(\fz):\ell^2(G)\to\ell^2(G)$ for $n\in\oNM$. Then the inclusions $\sigma\big(H_{Y_{n+1}}\big)\subseteq\sigma\big(H_{Y_{n}}\big)\,,\; n\in\NM\,,$ hold and
$$
\sigma\big(H_{Y_{\infty}}\big) \; 
	= \; \bigcap_{n\in\NM}\sigma\big(H_{Y_{n}}\big)\ .
$$
\end{proposition}

\begin{proof}
For $n\in\NM$, the equality $\sigma\big(H_{Y_n}\big)=\bigcup_{y\in Y_n} \sigma(H_y)$ holds by Proposition~\ref{Chap2-Prop-CovFamOp-Spect}~(a). Thus, the inclusion $\sigma\big(H_{Y_{n+1}}\big)\subseteq\sigma\big(H_{Y_{n}}\big)$ is derived by $Y_{n+1}\subseteq Y_n$. With this at hand, the identity $\bigcap_{n=1}^m\sigma\big(H_{Y_{n}}\big) = \sigma\big(H_{Y_{m}}\big)$ is deduced. Since $\lim_{m\to\infty}Y_m=Y_\infty$, the convergence $\lim_{m\to\infty}\sigma\big(H_{Y_{m}}\big)=\sigma\big(H_{Y_{\infty}}\big)$ follows by Theorem~\ref{Chap4-Theo-CharDynSystConvSpectr} leading to the desired result.
\end{proof}

\begin{remark}
\label{Chap4-Rem-MonotonConvSpectr}
Consider a symbolic dynamical system $(\as^G,G,\alpha)$ and an exhausting sequence $(K_n)_{n\in\NM}$ of $G$, c.f. Definition~\ref{Chap2-Def-ExhaustSeq}. Let $\Xi\in\SG\big(\as^G\big)$ be a subshift and define the subshifts of finite type
$$
\Xi_n \;
	:= \; \big\{
		\eta\in\as^G\;\big|\; 
		\ws(\eta)\cap\as^{[K_n]}\subseteq \ws(\Xi)\cap\as^{[K_n]} 
	\big\}\,,
	\qquad n\in\NM.
$$
Then $(\Xi_n)_{n\in\NM}$ converges to $\Xi$ in $\SG\big(\as^G\big)$ by Theorem~\ref{Chap2-Theo-Shift+DictSpace}. Furthermore, the inclusion $\Xi_{n+1}\subseteq\Xi_n$ is satisfied for $n\in\NM$ and so the requirements of Proposition~\ref{Chap4-Prop-MonotonConvSpectr} are fulfilled. Suppose in addition that the subshift of finite type $\Xi_n$ has a dense subset of strongly periodic elements for each $n\in\NM$. Then the spectrum of a Schr\"odinger operator $H_{\Xi_n}=\big(H_\xi\big)_{\xi\in\Xi_n}$ is already determined by the spectra $\sigma(H_\xi)\,,\;\xi\in\Xi_{per}\,,$ where $\Xi_{per}\subseteq\Xi_n$ is the set of strongly periodic elements of $\Xi_n$, c.f. Proposition~\ref{Chap2-Prop-SubsFinTypPerSpectr}. Under suitable conditions on the group $G$, the subshift $\Xi_n$ has a dense subset of strongly periodic elements if $\Xi$ is periodically approximable, c.f. Remark~\ref{Chap2-Rem-SubsFinTypPerSpectr}~(ii) for details.
\end{remark}

\subsection{Symbolic dynamical systems}
\label{Chap4-Ssect-SymbDynSyst}

Looking at the proof of the implication (ii)$\Rightarrow$(i) in Theorem~\ref{Chap4-Theo-CharDynSystConvSpectr}, we see that only specific generalized Schr\"odinger operators are necessary to guarantee the convergence of the dynamical subsystems. Specifically, it is sufficient to consider those operators defined by continuous functions on the universal groupoid $\Gun$ that describe the topology on $\Gun$. Having this in mind, assertion (ii) in Theorem~\ref{Chap4-Theo-CharDynSystConvSpectr} can be relaxed to the class of pattern equivariant Schr\"odinger operators for symbolic dynamical systems. Therefore, its Definition~\ref{Chap2-Def-SchrOp-l2(G)} is recalled in the following.

\medskip

Let $K\in\ks(G)$ be a compact set with pattern equivariant functions $p_h:\as^G\to\CM,\; h\in K,$ and $p_e:\as^G\to\RM$. Define, for every subshift $\Xi\in\SG\big(\as^G\big)$, the associated family of bounded self-adjoint operators $H_\xi:\ell^2(G)\to\ell^2(G)\,,\; \xi\in\Xi\,,$ by
$$
(H_\xi\psi)(g) \; 
	:= \; \left( 
			\sum\limits_{h\in K} p_h\big(\alpha_{g^{-1}}(\xi)\big) \cdot \psi(g\,h^{-1}) + \overline{p_h\big(\alpha_{(gh)^{-1}}(\xi)\big)} \cdot \psi(gh)
		\right)
		+ p_e\big(\alpha_{g^{-1}}(\xi)\big) \cdot \psi(g)
$$
for $\psi\in\ell^2(G)$ and $g\in G$. Then the family $H_\Xi:=(H_\xi)_{\xi\in\Xi}$ is called a pattern equivariant Schr\"odinger operator associated with $\Xi\in\SG(X)$ for a given compact set $K\in\ks(G)$ and the pattern equivariant functions $p_h:\as^G\to\CM,\; h\in K,$ and $p_e:\as^G\to\RM$.

\medskip

Recall that $(\as^G,G,\alpha)$ defines a topological dynamical system for a finite alphabet $\as$ and a discrete, countable group $G$, c.f. Section~\ref{Chap2-Ssect-SpaceSubshifts} and Proposition~\ref{Chap2-Prop-FullShift}.

\begin{theorem}[Continuous behavior of the spectra and subshifts, \cite{BeBeNi16}]
\label{Chap4-Theo-CharSubshiftConvSpectr}
Let $\as$ be a finite alphabet and $G$ a discrete, countable group such that the topological groupoid $\Gamma:=\as^G\rtimes_\alpha G$ is topologically amenable (say e.g. $G$ is an amenable group). Consider a subshift $\Xi_n\in\SG\big(\as^G\big)$ for each $n\in\oNM=\NM\cup\{\infty\}$. Then the following assertions are equivalent.
\begin{itemize}
\item[(i)] The sequence $(\Xi_n)_{n\in\NM}$ converges to $\Xi_\infty\in\SG\big(\as^G\big)$ in the Hausdorff-topology of the space $\SG\big(\as^G\big)$.
\item[(ii)] The sequence of dictionaries $(\ws(\Xi_n))_{n\in\NM}$ converges to $\ws(\Xi_\infty)\in\DG$ in the local pattern topology on $\DG$.
\item[(iii)] For every $K\in\ks(G)$ and all pattern equivariant functions $p_h:\as^G\to\CM,\; h\in K,$ and $p_e:\as^G\to\RM$, the map 
$$
\Sigma:\oNM\to \ks(\RM)\,,\;
	n\;\mapsto\; \sigma(H_{\Xi_n})\,,
$$ 
is continuous with respect to the Hausdorff metric on $\ks(\RM)$ where $H_{\Xi_n}$ is the pattern equivariant Schr\"odinger operator associated with $\Xi_n$ for $n\in\oNM$.
\item[(iv)] For all pattern equivariant functions $p:\as^G\to\RM$, the map $\oNM\ni n\mapsto\sigma(V_{\Xi_n})\in\ks(\RM)$ is continuous with respect to the Hausdorff metric on $\ks(\RM)$ where $V_{\Xi_n}:=(V_\xi)_{\xi\in\Xi_n}$ for $n\in\oNM$ are defined by
$$
(V_\xi\psi)(g) \; 
	= \; p\left(\alpha_{g^{-1}}(\xi)\right) \cdot \psi(g)\, ,
		\qquad\, \xi\in\as^G\, .
$$
\end{itemize}
\end{theorem}

\begin{proof}
The equivalence of (i) and (ii) is proven in Theorem~\ref{Chap2-Theo-Shift+DictSpace}. 

\vspace{.1cm}

(i)$\Rightarrow$(iii): Let $K\in\ks(G)$ be finite and $p_h:\as^G\to\CM,\; h\in K,$ and $p_e:\as^G\to\RM$ be pattern equivariant functions. For $\Xi\in\SG\big(\as^G\big)$, consider the associated pattern equivariant Schr\"odinger operator $H_\Xi$. Corollary~\ref{Chap4-Cor-ExContSectContSpectrSymbDynSyst} implies the continuity of the map $\SG(X)\ni \Xi\mapsto\sigma(H_\Xi)\in\ks(\RM)$ with respect to the Hausdorff metric. Furthermore, the map $\oNM\ni n\mapsto \Xi_n\in\SG\big(\as^G\big)$ is continuous by (i). Hence, $\Sigma:\oNM\to \ks(\RM)\,,\; n\mapsto \sigma(H_{\Xi_n})\,,$ is continuous as a composition of continuous maps proving assertion (iii).

\vspace{.1cm}

(iii)$\Rightarrow$(iv): Let $p:\as^G\to\RM$ be a pattern equivariant function. Set $K:=\emptyset$ and $p_e:=p$. Then the associated pattern equivariant Schr\"odinger operator $H_{\Xi_n}$ is equal to $V_{\Xi_n}$ for all $n\in\oNM$. Then (iii) implies the desired continuity $\oNM\ni n\mapsto\sigma(V_{\Xi_n})\in\ks(\RM)$.

\vspace{.1cm}

(iv)$\Rightarrow$(i): This proof follows the lines of the implication (ii)$\Rightarrow$(i) in Theorem~\ref{Chap4-Theo-CharDynSystConvSpectr}. In detail, assume that $(\Xi_n)_{n\in\NM}$ does not converge to $\Xi_\infty\in\SG\big(\as^G\big)$. By compactness of $\SG\big(\as^G\big)$, there exists a convergent subsequence $(\Xi_{n_k})_{k\in\NM}$ such that $Y:=\lim_{k\to\infty}\Xi_{n_k}\neq\Xi_\infty$. Then (a) $\Xi_\infty\setminus Y$ is non-empty or (b) $Y\setminus\Xi_\infty$ is non-empty. As in  Theorem~\ref{Chap4-Theo-CharDynSystConvSpectr} both cases lead to a contradiction. Only the case (a) is presented here. The case (b) is similarly treated.

\vspace{.1cm}

(a) If $\Xi_\infty\setminus Y$ is non-empty, consider an $\eta\in \Xi_\infty\setminus Y$. A base for the topology of $\as^G$ is given by the clopen sets of the form $\os(K,[u]):=\{\xi\in\as^G\;|\; \xi|_K\in[u]\}$ where $K\in\ks(G)$ and $[u]\in\as^{[K]}$, c.f. Section~\ref{Chap2-Ssect-SpaceSubshifts}. Recall that $[K]$ and $[u]$ for $K\in\ks(G)$ and $u:K\to\as$ denote the corresponding equivalence classes with respect to the translation by $G$. The topological space $\as^G$ is second countable, locally compact and Hausdorff, c.f. Proposition~\ref{Chap2-Prop-FullShift}. Thus, $\as^G$ is normal and so there exists an open neighborhood $\os:=\os(K,[u])$ of $\eta$ with $K\in\ks(G)$ and $[u]\in\as^{[K]}$ such that $\os\cap Y=\emptyset$, c.f. Definition~\ref{App1-Def-NormalSpace}. The set $\os$ is compact and open implying that the characteristic function $p:\as^G\to\{0,1\}\,,\; p(\xi):=\chi_\os(\xi)\,,$ of the set $\os$ is continuous with compact support. According to Proposition~\ref{Chap2-Prop-CharPattEqFunc}~(ii), the function $p$ is pattern equivariant. Now, consider the families of operators $V_{\Xi_n}$ for $n\in\oNM$ associated with the previously defined pattern equivariant function $p$.

\vspace{.1cm}

The convergence $\lim_{k\to\infty}\Xi_{n_k}=Y$ and the fact that $\os\cap Y=\emptyset$ lead to the existence of a $k_0\in\NM$ such that $\Xi_{n_k}\in\us(\os,\{\as^G\})$ for $k\geq k_0$. Here $\us(\os,\{\as^G\})$ is an open neighborhood of $Y\in\SG\big(\as^G\big)$ defined by
$$
\big\{\Xi\in\SG\big(\as^G\big)\;\big|\; \Xi\cap\os=\emptyset\,,\; \Xi\cap\as^G\neq\emptyset\big\} \, .
$$
Thus, the spectra $\sigma(V_{\Xi_{n_k}})\,,\; k\geq k_0\,,$ are equal to $\{0\}$ as the function $p$ vanishes on $\Xi_{n_k}$ for all $k\geq k_0$. Since $\eta\in\os$, the inequalities $\|V_\eta\delta_e\|\geq|p(\eta)|=1$ hold where $e\in G$ is the neutral element and $\delta_e\in\ell^2(G)$ is the Kronecker delta function with $\ell^2$-norm equal to $1$. Consequently, the estimate $\|V_{\Xi_\infty}\|\geq\|V_\eta\|\geq 1$ follows. The family of operators $V_{\Xi_\infty}$ is self-adjoint and so the spectral radius agrees with its norm $\|V_{\Xi_\infty}\|$. Thus, there exists a $\lambda\in\sigma(V_{\Xi_\infty})$ such that $|\lambda|\geq 1$. Hence, the Hausdorff distance of the spectra satisfies $d_H(\sigma(V_{\Xi_\infty}),\sigma(V_{\Xi_{n_k}}))\geq 1$ whenever $k\geq k_0$. This implies that the map $\oNM\ni n\mapsto \sigma(V_{\Xi_n})\in\ks(\RM)$ is not continuous with respect to the Hausdorff metric being a contradiction to (iv).
\end{proof}

\begin{remark}
\label{Chap4-Rem-ChoicePEFunctOnAsG}
(i) Note that instead of considering each $K\in\ks(G)$ in the assertion of Theorem~\ref{Chap4-Theo-CharSubshiftConvSpectr}~(iii), it suffices to consider any singletons $K\in\ks(G)$, c.f. Theorem~\ref{Chap4-Theo-CharSubshiftConvSpectr}~(iv). 

\vspace{.1cm}

(ii) In Definition~\ref{Chap2-Def-SchrOp-l2(G)} of generalized Schr\"odinger operators, it is assumed that the pattern equivariant functions are defined on $\as^G$ and not only on $\Xi\in\SG\big(\as^G\big)$. Proposition~\ref{Chap2-Prop-CharPattEqFunc} states that every pattern equivariant function is extendable to a pattern equivariant function on $\as^G$, c.f. Remark~\ref{Chap2-Rem-ExtPatEqFunct} as well. The statement of Theorem~\ref{Chap4-Theo-CharSubshiftConvSpectr} is independent of the choice of the extension for the following reason. Consider a pattern equivariant Schr\"odinger operator with a finite number of pattern equivariant functions. Then there exists an $N\in\NM$ such that, for all of these pattern equivariant functions, the value at $\xi\in\Xi_\infty$ only depend on $\xi|_{[-N,N]}$, c.f. Definition~\ref{Chap2-Def-PattEqFunc}. Thus, if $(\Xi_n)_{n\in\NM}$ is a sequence of subshifts tending to $\Xi_\infty$, then there exists an $n_0\in\NM$ such that $\ws(\Xi_n)\cap\as^{2N+1}=\ws(\Xi_\infty)\cap\as^{2N+1}$ for all $n\geq n_0$, c.f. Theorem~\ref{Chap2-Theo-Shift+DictSpace}. Consequently, the operators $H_{\Xi_n}$ are independent of the choice of extensions of the pattern equivariant functions on $\as^G$ for $n\geq n_0$. This justifies that the pattern equivariant Schr\"odinger operators are directly defined for pattern equivariant function on the whole space $\as^G$.
\end{remark}

In the following the special case $G=\ZM$ is considered. The families of operators in (iii) are the so called Jacobi operators whereas the families of operators in (iv) are called Schr\"odinger operators, c.f. Section~\ref{Chap2-Sect-ExampleAsZMSchrOp}. The case $G=\ZM$ is presented here as it is the most studied case in the literature. Furthermore, the assertion is used in Chapter~\ref{Chap7-Examples} for the one-dimensional examples.

\begin{corollary}[Continuous behavior of the spectra and subshifts, \cite{BeBeNi16}]
\label{Chap4-Cor-CharSubshiftConvSpectrZM}
Let $\as$ be a finite alphabet. Consider a subshift $\Xi_n\in\SZ\big(\as^\ZM\big)$ for each $n\in\oNM=\NM\cup\{\infty\}$. Then the following assertions are equivalent.
\begin{itemize}
\item[(i)] The sequence $(\Xi_n)_{n\in\NM}$ converges to $\Xi_\infty$ in the Hausdorff-topology of $\SZ\big(\as^\ZM\big)$.
\item[(ii)] The sequence of dictionaries $(\ws(\Xi_n))_{n\in\NM}$ converges to $\ws(\Xi_\infty)\in\mathfrak{D}_{\ZM}(\as)$ in the local pattern topology on $\mathfrak{D}_{\ZM}(\as)$.
\item[(iii)] Let $p:\as^\ZM\to\CM$ and $q:\as^\ZM\to\RM$ be pattern equivariant functions. For $n\in\oNM$, consider the associated family $J_{\Xi_n}:=(J_{\xi})_{\xi\in\Xi_n}$ of bounded self-adjoint operators $J_\xi:\ell^2(\ZM)\to\ell^2(\ZM)$ defined by
$$
\!\!\!
(J_\xi\psi)(m) 
	:= \left( 
			p\big(\alpha_{-m}(\xi)\big) \cdot \psi(m-1) + \overline{p\big(\alpha_{-(m+1)}(\xi)\big)} \cdot \psi(m+1)
		\right)
		+ q\big(\alpha_{-m}(\xi)\big) \cdot \psi(m)
$$
for $\psi\in\ell^2(\ZM)$ and $m\in \ZM$. Then the equation $\lim_{n\to\infty}\sigma(J_{\Xi_n})=\sigma(J_{\Xi_\infty})$ holds where the limit is taken with respect to the Hausdorff metric on $\ks(\RM)$.
\item[(iv)] Let $p:\as^\ZM\to\RM$ be a pattern equivariant functions. For $n\in\oNM$, consider the associated family $S_{\Xi_n}:=(S_\xi)_{\xi\in\Xi_n}$ of bounded self-adjoint operators $S_\xi:\ell^2(\ZM)\to\ell^2(\ZM)$ defined by
$$
\!\!\!
(S_\xi\psi)(m) 
	:= \psi(m-1) + \psi(m+1)
		+ p\big(\alpha_{-m}(\xi)\big) \cdot \psi(m)
$$
for $\psi\in\ell^2(\ZM)$ and $m\in \ZM$. Then the equation $\lim_{n\to\infty}\sigma(S_{\Xi_n})=\sigma(S_{\Xi_\infty})$ holds where the limit is taken with respect to the Hausdorff metric on $\ks(\RM)$.
\end{itemize}
\end{corollary}

\begin{proof}
The sets $F_k:=\{-k,-k+1,\ldots,k\}\,,\; k\in\NM\,,$ define a weak-F\o lner sequence for $\ZM$, see page~\pageref{Page-weak-Folner} for the definition of weak-F\o lner sequence. Thus, $\ZM$ is a finitely generated, discrete, amenable group. Hence, $\Gamma(\Xi):=\Xi\rtimes_\alpha\ZM$ is topologically amenable, c.f. Proposition~\ref{Chap2-Prop-TransGroupAmen}. The equivalence of (i) and (ii) is proven in Theorem~\ref{Chap2-Theo-Shift+DictSpace}. Then the implication (i)$\Rightarrow$(iii) follows by Theorem~\ref{Chap4-Theo-CharSubshiftConvSpectr}. By setting $p\equiv 1$, the implication (iii)$\Rightarrow$(iv) is deduced. Finally, the implication (iv)$\Rightarrow$(i) is proven in analogy to the proof of (iv)$\Rightarrow$(i) in Theorem~\ref{Chap4-Theo-CharSubshiftConvSpectr} as follows.

\vspace{.1cm}

(iv)$\Rightarrow$(i): Assume that $(\Xi_n)_{n\in\NM}$ does not converge to $\Xi_\infty\in\SZ\big(\as^\ZM\big)$. By compactness of $\SZ\big(\as^\ZM\big)$, there exists a convergent subsequence $(\Xi_{n_k})_{k\in\NM}$ such that $Y:=\lim_{k\to\infty}\Xi_{n_k}\neq\Xi_\infty$. Then (a) $\Xi_\infty\setminus Y$ is non-empty or (b) $Y\setminus\Xi_\infty$ is non-empty. Like in Theorem~\ref{Chap4-Theo-CharSubshiftConvSpectr}, both cases lead to a contradiction to (iv). Since the proofs follow the lines of the proof of (iv)$\Rightarrow$(i) in Theorem~\ref{Chap4-Theo-CharSubshiftConvSpectr}, the main steps for the case (a) are only presented.

\vspace{.1cm}

(a) Let $\eta\in\Xi_\infty\setminus Y$. Define $p:\as^\ZM\to\RM$ by $p(\xi):=\sqrt{7}\cdot\chi_\os(\xi)$ where $\os$ is a compact and open neighborhood of $\eta$ such that $\os\cap Y=\emptyset$, c.f. proof of Theorem~\ref{Chap4-Theo-CharSubshiftConvSpectr}. Hence, $p$ is a pattern equivariant function. Since $Y=\lim_{k\to\infty}\Xi_{n_k}$ and $\os\cap Y=\emptyset$, there exists a $k_0\in\NM$ such that $\Xi_{n_k}\cap\os=\emptyset$ for $k\geq k_0$. Consider the associated families of Schr\"odinger operators $S_{\Xi_n}:=\big(S_\xi:\ell^2(\ZM)\to\ell^2(\ZM)\big)_{\xi\in\Xi_n}\,,\; n\in\oNM$ with the previously defined $p:=\sqrt{7}\cdot\chi_\os$. For $k\geq k_0$, the pattern equivariant function $p$ vanishes on $\Xi_{n_k}$ implying 
$$
\!\!\!
(S_\xi\psi)(m) 
	:= \psi(m-1) + \psi(m+1)\,,
	\qquad \psi\in\ell^2(\ZM)\,, \; m\in\ZM\,, \; \xi\in\Xi_{n_k}\,.
$$
Thus, $S_{\Xi_{n_k}}$ satisfies $\|S_{\Xi_{n_k}}\|\leq 2$ for $k\geq k_0$. Note that it is well-known that actually the equation $\|S_{\Xi_{n_k}}\|= 2$ holds for $k\geq k_0$. On the other hand, $\eta$ is an element of $\os$ leading to $p(\eta)=\sqrt{7}$. The neutral element of the group $\ZM$ is $0$ and $\delta_0\in\ell^2(\ZM)$ denotes the Kronecker delta function at $0\in\ZM$. Then a short calculation leads to
\begin{align*}
\|S_\eta\delta_0\|^2 \;
	&= \; \sum_{j\in\ZM} \big| \big(S_\eta\delta_0 \big)(j)\big|^2\\
	&= \; \sum_{j\in\ZM} \big|
					\delta_0(j-1) +\delta_0(j+1) + p(\alpha_{-j}(\eta))\cdot \delta_0(j)
				\big|^2\\
	&= \; 1 + 1 + |p(\eta)|^2\\
	&= \; 9\, .
\end{align*}
Since $\delta_0$ has $\ell^2$-norm equal to $1$, the last equalities imply $\|S_{\Xi_\infty}\|\geq\|S_\eta\|\geq 3$. Hence, there exists a $\lambda\in\sigma(S_{\Xi_\infty})$ with $|\lambda|\geq 3$ since $S_\eta$ is self-adjoint. Consequently, the estimate $d_H(\sigma(S_{\Xi_\infty}),\sigma(S_{\Xi_{n_k}}))\geq 1$ holds for all $k\geq k_0$ by the previous considerations. Thus, $\big(\sigma(S_{\Xi_n})\big)_{n\in\NM}$ does not tend to $\sigma(S_{\Xi_\infty})$ contradicting (iv).
\end{proof}

\medskip

Assertion (iii) and (iv) of Corollary~\ref{Chap4-Cor-CharSubshiftConvSpectrZM} meet that the spectra of \underline{all} Jacobi or Schr\"o\-dinger operators behave continuously. In particular, it can happen for specific pattern equivariant functions $q:\as^\ZM\to\CM$ and $p:\as^\ZM\to\RM$ that the equation $\lim_{n\to\infty}\sigma\big(J_{\Xi_n}\big)=\sigma\big(J_{\Xi_\infty}\big)$ holds while $(\Xi_n)_{n\in\NM}$ does not converge to $\Xi_\infty$. For instance, set $q\equiv 1$ and $p\equiv 0$. Then, for each family of subshifts $\Xi_n\,,\; n\in\oNM\,,$ the equality $\lim_{n\to\infty}\sigma\big(J_{\Xi_n}\big)=\sigma\big(J_{\Xi_\infty}\big)$ follows since the pattern equivariant functions do not distinguish any pattern. Actually, the equation $\sigma\big(J_{\Xi}\big)=\sigma\big(J_{\Xi'}\big)$ is valid for all subshifts $\Xi,\Xi'\in\SZ\big(\as^\ZM\big)$.

\begin{remark}
\label{Chap4-Rem-StrongConvSpectrConvGenSchrOp}
For a subshift $\Xi\in\SG\big(\as^G\big)$, each pattern equivariant Schr\"odinger operator $H_\Xi:=(H_\xi)_{\xi\in\Xi}$ is represented by a bounded self-adjoint operator $\bigoplus_{\xi\in\Xi}H_\xi$ on the Hilbert space $\bigoplus_{\xi\in\Xi}\ell^2(G)$. Theorem~\ref{Chap4-Theo-CharSubshiftConvSpectr} and its corollaries provide a characterization of the \pt-con\-tinuity of all generalized Schr\"odinger operators. Proposition~\ref{Chap3-Prop-P2NotImplStrong} and Proposition~\ref{Chap3-Prop-StrongNotImplP2} assert that there is no relation between the \pt-topology and the strong operator topology in general. By the previous considerations and Proposition~\ref{Chap2-Prop-CovFamOp-Spect}~(c), generalized Schr\"odinger operators  can be continuous with respect to both topologies. 

\vspace{.1cm}

More specifically, consider $(\as^G,G,\alpha)$ a symbolic dynamical system satisfying that the transformation group groupoid $\as^G\rtimes_\alpha G$ is topologically amenable. For each $\eta\in\as^G$ a pattern equivariant Schr\"odinger operator is defined by a fixed finite $K\in\ks(G)$ with pattern equivariant functions $p_h:\as^G\to\CM\,,\; h\in K\,,$ and $p_e:\as^G\to\RM$. Suppose that a sequence $(\xi_n)_{n\in\NM}$ in $\as^G$ converges to $\xi_\infty\in\as^G$  such that the dictionaries $\big(\ws(\xi_n)\big)_{n\in\NM}$ converge to $\ws(\xi_\infty)$ in the local pattern topology. Then the following assertions hold.
\begin{description}
\item[(a)] The pattern equivariant Schr\"odinger operators $\big(H_{\xi_n}\big)_{n\in\NM}$ converge strongly to $H_{\xi_\infty}$.
\item[(b)] The equation $\lim_{n\to\infty}\sigma(H_{\Xi_n})=\sigma(H_{\Xi_\infty})$ holds for $\Xi_n:=\overline{\Orb(\xi_n)}\,,\; n\in\oNM\,,$ where the limit is taken with respect to the Hausdorff metric on $\ks(\RM)$.
\end{description}
Assertion (a) follows by Proposition~\ref{Chap2-Prop-CovFamOp-Spect}~(c) while statement (b) is derived from Theorem~\ref{Chap4-Theo-CharSubshiftConvSpectr}. If, additionally, all the subshift $\Xi_n\,,\; n\in\oNM\,,$ are minimal, then 
$$
\lim_{n\to\infty}\sigma(H_{\xi_n}) \; 
	= \; \lim_{n\to\infty}\sigma_{ess}(H_{\xi_n}) \;
		= \; \sigma(H_{\xi_\infty}) \;
			= \; \sigma_{ess}(H_{\xi_\infty})
$$
is valid by Theorem~\ref{Chap2-Theo-ConstSpectrMinimal}.
\end{remark}

The following example shows that the convergence of $(\xi_n)_{n\in\NM}$ to $\xi_\infty$ is not sufficient for the convergence of the associated dictionaries. On the other hand, the convergence of the dictionaries is needed in general to get the convergence of the spectra, c.f. Theorem~\ref{Chap4-Theo-CharSubshiftConvSpectr}. Conversely, the convergence of the dictionaries $\lim_{n\to\infty}\ws(\xi_n)=\ws(\xi_\infty)$ does not imply $\lim_{n\to\infty}\xi_n=\xi_\infty$ since $\ws\big(\alpha_g(\xi)\big)=\ws(\xi)$ holds, c.f. Lemma~\ref{Chap2-Lem-xiDict}.

\begin{example}
\label{Chap4-Ex-StrongConvSpectrConvGenSchrOp}
Let $G=\ZM$. For two finite words $u\in\as^n$ and $v\in\as^m$, the periodic extension $(u|v)^\infty$ of $u|v$ to a two-sided infinite word is defined by the unique $\xi\in\as^\ZM$ satisfying $\xi|_{[-n+k(n+m),(m-1)+k(n+m)]}=uv$ for all $k\in\ZM$, c.f. Section~\ref{Chap2-Sect-ExampleAsZMSchrOp} and Section~\ref{Chap5-Sect-SymbDynSystZM}. Let $\as:=\{a,b\}$ be an alphabet with two elements. Define $\xi_\infty:=(a|b)^\infty$ and $\xi_n:=\big(b(ab)^na|b(ab)^n\big)$ for $n\in\NM$, i.e.,
\begin{align*}
\xi_\infty \;
	&= \; \ldots b\, a\, b\, a\, b\, a\, b\, a\, b\, a\, |\, b\, a\, b\, a\, b\, a\, b\, a\, b\, a\, \ldots \, ,\\
\xi_n \;
	&= \; \ldots b\, b\, \underbrace{a\, b\, \ldots\, a\, b}_{n-\text{times}}\, a\, |\, b\, \underbrace{a\, b\, \ldots a\, b}_{n-\text{times}}\, b\, a\,\ldots \, .
\end{align*}
By construction, the restrictions $\xi_\infty|_{[-2n-1,2n+1]}$ and $\xi_n|_{[-2n-1,2n+1]}$ coincide for $n\in\NM$. Thus, the sequence $(\xi_n)_{n\in\NM}$ converges to $\xi_\infty$. On the other hand, the word $bb$ is an element of $\ws(\xi_n)$ for all $n\in\NM$ while $bb\not\in\ws(\xi_\infty)$. Hence, the dictionaries $(\ws(\xi_n))_{n\in\NM}$ do not converge to $\ws(\xi_\infty)$ in the local pattern topology.
\end{example}

Proposition~\ref{Chap3-Prop-NormImplP2} asserts that the convergence in operator norm is stronger than the convergence with respect to the \pt-topology. The following example shows that the \pt-convergence and the strong operator convergence together do not imply the opera\-tor norm convergence. Specifically, it shows that the pattern equivariant Schr\"odinger operators do not converge in operator norm.

\begin{example}
\label{Chap4-Ex-KeineNormKonv}
Let $\as:=\{a_1,\ldots,a_N\}$ be a finite alphabet and $G$ a discrete, countable group. Suppose $(\xi_n)_{n\in\NM}$ does converge to $\xi_\infty\in\as^G$ such that $\xi_n\neq_\infty$ for all $n\in\NM$. Let $K\in\ks(G)$ be compact and define $p_h\equiv 1$ for $h\in K$. Define $p_e:\as^G\to\RM$ by 
$$
p_e \;
	= \; \sum_{j=1}^N j \cdot \chi_{\os(\{e\},[a_j])} \,,
$$
i.e., $p_e(\eta)=j$ if $\eta(e)=a_j\in\as$. Clearly, $p_h$ and $p_e$ are pattern equivariant functions.  According to Definition~\ref{Chap2-Def-SchrOp-l2(G)}, define the bounded self-adjoint operator $H_\eta:\ell^2(G)\to\ell^2(G)$ for $\eta\in\as^G$ associated with $K$ and $p_h\,,\; h\in K\cup\{e\}$. More precisely, $H_\eta$ is defined by
$$
(H_\xi\psi)(g) \; 
	:= \; \left( 
			\sum\limits_{h\in K} \psi(g\,h^{-1}) + \psi(gh)
		\right)
		+ p_e\left(\alpha_{g^{-1}}(\xi)\right) \cdot \psi(g) \,,
	\qquad \psi\in\ell^2(G)\,,\; g\in G\, .
$$
For each $n\in\NM$, consider a $g_n\in G$ such that $\xi_n(g_n)\neq\xi_\infty(g_n)$. Such a $g_n$ exists for each $n\in\NM$ since $\xi_n\neq\xi_\infty$ for all $n\in\NM$. Now, the shift $\alpha_g(\eta)$ is an element of $\os(\{e\},[a_j])$ if and only if $a_j=\alpha_g(\eta)(e)=\eta(g^{-1})$ for $\eta\in\as^G$, $g\in G$ and $j\in\{1,\ldots,N\}$. Thus, the estimate
$$
\Big|p_e\Big(\alpha_{g_n^{-1}}(\xi_n)\Big) - p_e\Big(\alpha_{g_n^{-1}}(\xi_\infty)\Big)\Big|\;
	\geq\; 1
$$
is derived. Consider the Kronecker delta function $\delta_{g_n}\in\ell^2(G)$ for $n\in\NM$ with $\ell^2$-norm $\|\delta_{g_n}\|_{\ell^2(G)}=1$. Then a short computation leads to 
$$
\Big\| \big(H_{\xi_n} - H_{\xi_\infty}\big)\Big\| \;
	\geq\; \Big\| \big(H_{\xi_n} - H_{\xi_\infty}\big) \delta_{g_n} \Big\| \;
	= \; \Big|p_e\Big(\alpha_{g_n^{-1}}(\xi_n)\Big) - p_e\Big(\alpha_{g_n^{-1}}(\xi_\infty)\Big)\Big| \;
	\geq\; 1
$$
for all $n\in\NM$. Consequently, $(H_{\xi_n})_{n\in\NM}$ does not converge to $H_{\xi_\infty}$ in operator norm. On the other hand, the operators converge in the \pt-topology and the strong operator topology if $\lim_{n\to\infty}\ws(\xi_n)=\ws(\xi_\infty)$, $\lim_{n\to\infty}\xi_n=\xi_\infty$ and the subshifts $\overline{\Orb(\xi_n)}\,,\; n\in\oNM\,,$ are minimal, c.f. Remark~\ref{Chap4-Rem-StrongConvSpectrConvGenSchrOp}.
\end{example}

It is worth noticing that the norm convergence of pattern equivariant Schr\"odinger operators cannot be expected in general by following the lines of Example~\ref{Chap4-Ex-KeineNormKonv}. Specifically, let $(\xi_n)_{n\in\NM}$ be a sequence in $\as^G$ converging to $\xi$ such that $\xi_n\neq\xi$. Then there exists a $g_n\in G$ such that $\xi_n(g_n)\neq\xi(g)$ for $n\in\NM$. Consider a pattern equivariant Schr\"odinger operator $H_\xi:\ell^2(G)\to\ell^2(G)\,,\; \xi\in\as^G$. In the general case, there will exists a constant $C>0$ corresponding to one pattern equivariant function $p:\as^G\to\CM$ associated with the operator $H$ such that $\big|p\big(\alpha_{g_n^{-1}}(\xi_n)\big)-p\big(\alpha_{g_n^{-1}}(\xi)\big)\big|\geq C\,,\; n\in\NM$. This implies $\big\|\big(H_{\xi_n} - H_\xi\big)\psi_n\big\|\geq C$ for a suitable functions $\psi_n\in\ell^2(G)$ with $\|\psi_n\|=1$ and $n\in\NM$.

\section{Approximation by periodic Schr\"odinger operators}
\label{Chap4-Sect-ApprPerAppr}

In this section the periodic approximation of generalized Schr\"odinger operators is presented by applying the result of Section~\ref{Chap4-Sect-ContSpecGenSchrOp}.

\medskip

According to Theorem~\ref{Chap4-Theo-CharDynSystConvSpectr}, the convergence of the spectra of generalized Schr\"odinger operators is characterized by the convergence of the underlying dynamical systems. Consequently, these operators can be approximated by periodic operators if the dynamical systems is periodically approximable.

\medskip

Periodic approximations are of particular interest since the spectra of these operators can be computed with the so called Floquet-Bloch theory. The Floquet-Bloch theory has its roots in the work of {\sc Floquet} \cite{Flo1883} in 1883 where he studied ordinary differential equations with periodic coefficients. He understands the question of stability very close to the modern concept of band theory. {\sc Bloch} \cite{BlochThesis29} rediscovered these concepts in a different context in his PhD thesis published 1929. He used the theory of Hilbert spaces and self-adjoint operators which applies independent of the dimension whereas {\sc Floquet} used the transfer matrix formulation which is only valid in one-dimension. The work of {\sc Bloch} leads to the point of view of band theory which constitutes the fundamental aspect of Solid State Physics. In 1937, {\sc Wannier} introduced in his thesis \cite{WannierThesis37} a generalization of the Fourier transform called nowadays the {\em Wannier transform}. There he used the translation invariance of the underlying system on which the Schr\"odinger operator is defined. By using that the orbit is finite and the stabilizer group of the system is abelian and cocompact, the operator is unitarily represented as a direct integral of finite dimensional Hilbert space over a compact group.

\medskip

Let $(X,G,\alpha)$ be a dynamical system. Recall that a dynamical subsystem $Y\in\SG(X)$ is called strongly periodic if $Y$ is minimal and each $y\in Y$ is strongly periodic, c.f. Definition~\ref{Chap2-Def-AperSubs}. Then $\SP(X)$ denotes the set of all strongly periodic dynamical subsystems of $X$ and $\PA(X):=\overline{\SP(X)}$ is the set of periodically approximable subshifts. A generalized Schr\"odinger operator $H_Y:=(H_y)_{y\in Y}$ associated with a dynamical system $(Y,G,\alpha)$ is called periodic whenever $Y$ is a strongly periodic dynamical system, c.f. Definition~\ref{Chap2-Def-SchrodingerOperator}. In this case, $H_y:\ell^2(G)\to\ell^2(G)$ is also called a periodic generalized Schr\"odinger operator for $y\in Y$.

\medskip

\begin{theorem}[Approximation by strongly periodic systems, \cite{BeBeNi16}]
\label{Chap4-Theo-PeriodicApproximations}
Let $(X,G,\alpha)$ be a topological dynamical system where $G$ is an amenable, countable, discrete group. Consider a periodically approximable subshift $Y\in\SG(X)$. Then, for every self-adjoint generalized Schr\"odinger operator $H_Y:=(H_y)_{y\in Y}$ of finite range and a $y_\infty\in Y$, there exists a sequence of periodic generalized Schr\"odinger operators $H_n:\ell^2(G)\to\ell^2(G)$ of finite range such that the following holds.
\begin{itemize}
\item[(a)] The sequence of operators $(H_n)_{n\in\NM}$ tends strongly to $H_{y_\infty}$.
\item[(b)] The associated spectra converge, i.e., the equalities
$$
\lim_{n\to\infty}\sigma(H_n) \;
	= \; \bigcup_{y\in Y} \sigma\big(H_y\big) \;
	= \; \sigma(H_Y)
$$
hold where the limit is taken with respect to the Hausdorff metric on $\ks(\RM)$.
\end{itemize}
If, additionally, $Y$ is minimal, then the equations $\lim\limits_{n\to\infty}\sigma(H_n)=\sigma(H_{y_\infty})=\sigma_{ess}(H_{y_\infty})$ hold.
\end{theorem}

\begin{proof}
According to the definition of a generalized Schr\"odinger operator of finite range, there is an $\fz\in\Cc_c\big(\Gamma^{(1)}(Y)\big)$ such that $\pi^y(\fz)=H_y$ for each $y\in Y$, c.f. Definition~\ref{Chap2-Def-SchrOp-l2(G)}. Let $\Gamma:=\Gun$ be the universal groupoid associated with $(X,G,\alpha)$. Proposition~\ref{Chap4-Prop-ExContSectContSpectr}~(b) implies the existence of a self-adjoint $\gz\in\Cc_c\big(\Gun^{(1)}\big)$ such that $\Theta_Y(\gz_Y)=\fz$ while the vector field $\big(\Theta_Z(\gz_Z)\big)_{Z\in\SG(X)}$ is self-adjoint and continuous in the continuous field of $C^\ast$-algebras $\big((\CG^\ast(\Gamma(Y)))_{Y\in\SG(X)}, \Upsilon \big)$ defined in Theorem~\ref{Chap4-Theo-TransGroupContFieldCAlg}. Here $\Theta_Y:\CG^\ast(\Gamma_Y)\to\CG^\ast(\Gamma(Y))$ denotes the $\ast$-isomorphism for $Y\in\SG(X)$, c.f. Corollary~\ref{Chap4-Cor-IsometrIsomorRedCAlg}. By construction, the generalized Schr\"odinger operator $H_Y$ is exactly the family of operators $\big(\pi^y\big(\Theta_Y(\gz_Y)\big)\big)_{y\in Y}$. Thus, the spectrum $\sigma(H_Y)$ coincides with $\sigma(\gz_Y)$. Now, define the generalized Schr\"odinger operator  $H_{Z}:=(H_x)_{x\in Z}$ of finite range by $H_x:=\pi^x\big(\Theta_Z(\gz_Z)\big)\,,\; x\in Z\,,$ for each $Z\in\SG(X)$ with $Z\neq Y$. 

\vspace{.1cm}

Since $Y$ is periodically approximable, there is a sequence of strongly periodic subshifts $(Y_n)_{n\in\NM}$ converging to $Y$ in $\SG(X)$. Then Theorem~\ref{Chap4-Theo-CharDynSystConvSpectr} leads to $\lim_{n\to\infty}\sigma(H_{Y_n})=\sigma(H_Y)$ since $\sigma(H_Z)=\sigma(\gz_Z)$ for all $Z\in\SG(X)$. 

\vspace{.1cm}

Let $y_\infty\in Y$. According to the claim (ii) in the proof of Proposition~\ref{Chap2-Prop-SpaDynSyst}, there is a $y_n\in Y_n$ for each $n\in\NM$ such that $\lim_{n\to\infty}y_n=y_\infty$. Hence, Proposition~\ref{Chap4-Prop-StrongConv+ConvSpectr} implies that the operators $H_n:=H_{y_n}\,,\; n\in\NM\,$ converge to $H_{y_\infty}$ in the strong operator topology on $\Ll(\ell^2(G))$. Since $Y_n\,,\; n\in\NM\,$ are strongly periodic, these dynamical subsystems are minimal. Consequently, the equality $\sigma(H_{y_n})=\sigma(H_{Y_n})$ follows for each $n\in\NM$ by Theorem~\ref{Chap2-Theo-ConstSpectrMinimal}. Together with the previous considerations, this implies $\lim_{n\to\infty}\sigma(H_n) = \sigma(H_Y)$ and $\sigma(H_Y)=\sigma(H_{y_\infty})=\sigma_{ess}(H_{y_\infty})$ if $Y$ is minimal.
\end{proof}

\begin{remark}
\label{Chap4-Rem-PeriodicApproximations}
(i) For a general topological dynamical system $(X,G,\alpha)$, the operators $(H_n)_{n\in\NM}$ are not given in a constructive way since the operators are defined by extensions of continuous function on the groupoids, c.f. Proposition~\ref{Chap4-Prop-ExContSectContSpectr}.

\vspace{.1cm}

(ii) Suppose $X=\as^G$ for a finite alphabet $\as$. Then the operators $(H_n)_{n\in\NM}$ and $(H_{Y_n})_{n\in\NM}$ are defined explicitly if $H_Y$ is a pattern equivariant Schr\"odinger operator, c.f. Theorem~\ref{Chap4-Theo-CharSubshiftConvSpectr} and its corollaries. 

\vspace{.1cm}

More precisely, let $Y\in\SG\big(\as^G\big)$ be periodically approximable with a sequence $(Y_n)_{n\in\NM}$ of strongly periodic subshifts converging to $Y$. Consider a pattern equivariant Schr\"odinger operator $H_Y:=(H_y)_{y\in Y}$ with compact $K\in\ks(G)$ and pattern equivariant functions $p_h:\as^G\to\CM\,,\; h\in K\,,$ and $p_e:\as^G\to\RM$. Then assertion (a) and (b) of Theorem~\ref{Chap4-Theo-PeriodicApproximations} hold where $H_n:=H_{y_n}:\ell^2(G)\to\ell^2(G)$ is given by $y_n\in Y_n\,,\; n\in\NM\,,$ converging to $y_\infty\in Y$ and the pattern equivariant Schr\"odinger operator with the same compact $K\in\ks(G)$ and pattern equivariant functions $p_h\,,\; h\in K\cup\{e\}$.
\end{remark}

\cleardoublepage


\chapter{The one-dimensional case}
\label{Chap5-OneDimCase}
\stepcounter{section}
\setcounter{section}{0}

Due to Theorem~\ref{Chap4-Theo-PeriodicApproximations}, generalized Schr\"odinger operators are approximable by periodic operators if the corresponding dynamical system is periodically approximable. Having this in mind, the set of periodically approximable subshifts of $(\as^\ZM,\ZM,\alpha)$ is characterized in this chapter, c.f. Theorem~\ref{Chap5-Theo-ExPerAppr}. With this at hand, it turns out that each minimal subshift of $(\as^\ZM,\ZM,\alpha)$ is periodically approximable, c.f. Proposition~\ref{Chap5-Prop-MinAllPath}. The investigated theory of this chapter is also a first step of analyzing the topological properties of periodically approximable subshifts as a subset of $\SG(X)$, c.f. discussion in Section~\ref{Chap8-Sect-PerApprox}.

\medskip

Specifically, we show that the property of a subshift being periodically approximable is equivalent to the fact that the associated de Bruijn graphs are strongly connected, c.f. Theorem~\ref{Chap5-Theo-ExPerAppr}. The proof uses the one-dimensional nature of the group $\ZM$ and the fact that subshifts can be identified with its dictionaries, c.f. Theorem~\ref{Chap2-Theo-Shift+DictSpace}. In detail, the so called de Bruijn graphs are defined via a dictionary $\ws$ which encode the local properties of the subshift $\Xi(\ws)\in\SZ\big(\as^\ZM\big)$. Then periodic approximations are defined via global paths in the de Bruijn graphs in Section~\ref{Chap5-Sect-PerApprSubs}. With this at hand, sufficient conditions are proven for a subshift $\Xi\in\SZ\big(\as^\ZM\big)$ to be periodically approximable, c.f. Proposition~\ref{Chap5-Prop-SuffCondPerAppr} and Corollary~\ref{Chap5-Cor-SuffCondPerAppr}. 

\medskip

The property of a dynamical subsystem being periodically approximable is a property of the Hausdorff-topology on $\SG(X)$ for a dynamical system $(X,G,\alpha)$. More precisely, the intersections $U_n(Y)\cap\SP(X)\,,\; n\in\NM\,,$ are non-empty for a neighborhood base $U_n(Y)\,,\; n\in\NM\,,$ of $Y\in\SG(X)$. The strongly periodic subshifts in $\SZ\big(\as^\ZM\big)$ are isolated in the Hausdorff-topology, c.f. Proposition~\ref{Chap5-Prop-PerPointIsolated}. Moreover, there exist also aperiodic subshifts that are isolated in the Hausdorff-topology, c.f. Example~\ref{Chap5-Ex-DeBruijnNotStrongConnect}.

\medskip

\label{Page-GeneralExStrongPerSubshFinitType} Note that a subshift in a symbolic dynamical system $(\as^G,G,\alpha)$ is (finitely) periodically approxi\-mable if and only if all the associated subshifts of finite type 
$$
\Xi_n \; 
	:= \; \big\{ 
		\eta\in\as^G \;\big|\; \ws(\eta)\cap\as^{[K_n]}\subseteq\ws(\Xi)\cap\as^{[K_n]}	
	\big\}\,,
	\qquad
	n\in\NM\,,
$$
contain sufficiently many strongly periodic elements since $\Xi=\bigcap_{n\in\NM}\Xi_n$ holds by Theorem~\ref{Chap2-Theo-Shift+DictSpace} where $(K_n)_{n\in\NM}$ is an exhausting sequence of $G$. Here finitely periodically approximable means that there exists a sequence $(\Xi_n)_{n\in\NM}$ converging to $\Xi$ in $\SZ\big(\as^\ZM\big)$ such that $\Xi_n$ is the finite union of strongly periodic subshifts. The existence of a (dense) subset of strongly periodic elements in a subshift of finite type is a current question in the field of dynamical systems \cite{Pia08,Hoc09,SiCo12,Coh14,CaPe15} which is also related to the so called Wang tiles \cite{Wan61,Ber66,CuKa95}, c.f. discussion on page~\pageref{Page-ExStrongPerSubshFinitType}.

\medskip

\begin{figure}[htb]
\centering
\includegraphics[scale=0.85]{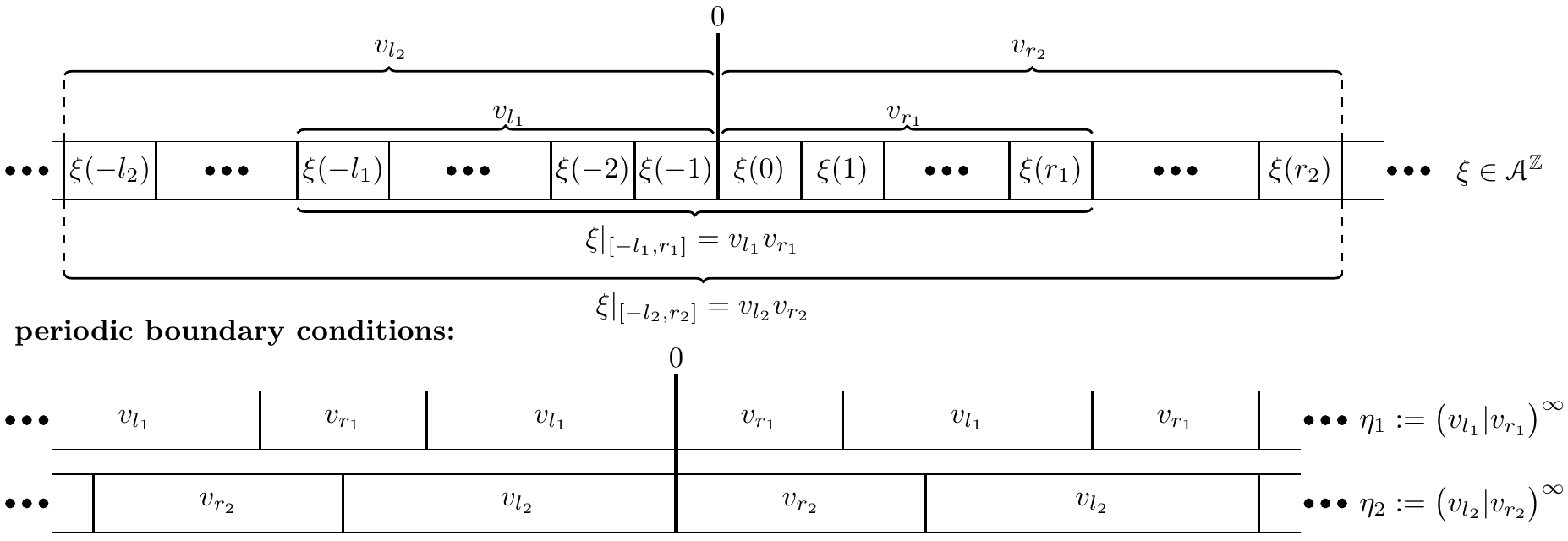}
\caption{Restrict $\xi\in\as^\ZM$ to $[-l_k,r_k]\subseteq\ZM\,,\; k\in\NM\,,$ around the origin and take periodic boundary conditions. This defines a sequence of strongly periodic $(\eta_k)_{k\in\NM}$ converging to $\xi$. If $l_k,r_k\in\NM$ are chosen correctly for $k\in\NM$, this implies the convergence of $(\Orb(\eta_k))_{k\in\NM}$ to $\overline{\Orb(\xi)}$ in $\SZ\big(\as^\ZM\big)$, c.f. assertion \nameref{(5.IV)}.}
\label{Chap5-Fig-StrategyPerAppr}
\end{figure}

Note that the notions of weakly periodic and strongly periodic elements or dynamical subsystems coincide for a given dynamical system $(X,\ZM,\alpha)$ over the group $\ZM$, c.f. discussion before Definition~\ref{Chap2-Def-AperSubs}. 

\medskip

The aim of this chapter is to find sufficient conditions for a subshift $\Xi\in\SZ\big(\as^\ZM\big)$ to be periodically approximable. The strategy is formulated as follows. 

\begin{description}
\item[(5.I)\label{(5.I)}] First, suppose $\Xi$ is topologically transitive, i.e., there exists a $\xi\in\as^\ZM$ such that $\Xi=\overline{\Orb(\xi)}$. Thanks to Corollary~\ref{Chap2-Cor-DictOrbitSubshift}, the equation $\ws(\Xi)=\ws(\xi)$ is derived. With this at hand, we use in the following that the convergence of subshifts is characterized by the convergence of the associated dictionaries with respect to the local pattern topology, c.f. Theorem~\ref{Chap2-Theo-Shift+DictSpace}. 
\item[(5.II)\label{(5.II)}] In the case of a topologically transitive subshift, it is natural to define a periodic approximation of $\Xi$ via a periodic approximation of $\xi$, c.f. Figure~\ref{Chap5-Fig-StrategyPerAppr}. Let $(l_k)_{k\in\NM}$ and $(r_k)_{k\in\NM}$ be sequences in $\NM$ that converge to infinity. For $k\in\NM$, consider the restriction $\xi_{[-l_k,r_k]}$ of $\xi\in\as^\ZM$ to $[-l_k,r_k]$ around the origin. Take periodic boundary condition for $\xi_{[-l_k,r_k]}$, i.e., define $\eta_k:=\big(\xi_{[-l_k,-1]}|\xi|_{[0,r_k]}\big)^{\infty}\,,\;k\in\NM$. Then $(\eta_k)_{k\in\NM}$ defines a sequence of strongly periodic element in $\as^\ZM$ that converges to $\xi$ in $\as^\ZM$. In view of \nameref{(5.I)}, it is necessary to choose $(\eta_k)_{k\in\NM}$ such that the dictionaries $\big(\ws(\eta_k)\big)_{k\in\NM}$ converge to $\ws(\xi)$ with respect to the local pattern topology. Due to Corollary~\ref{Chap2-Cor-DictExhSeq}, this is equivalent to the existence of a $k_0\in\NM$ for each $m\in\NM$ satisfying 
$$
\ws(\eta_k)\cap\as^m \;
	=\; \ws(\xi)\cap\as^m\,,
	\qquad k\geq k_0\,.
$$
\item[(5.III)\label{(5.III)}] Let $m\in\NM$. Since the set $\as^m$ is finite, there exists an $N\in\NM$ such that $\ws(\xi)\cap\as^m\subseteq\Sub\big[\xi|_{[-N,N]}\big]$ (the subwords of $\xi|_{[-N,N]}$). Thus, the subwords $\ws(\xi)\cap\as^m$ occur in $\xi$ with finite distance $N\in\NM$ from the origin. Consequently, there exists a $k_0\in\NM$ such that $\xi|_{[-N,N]}\in\Sub\big[\xi|_{[-l_k,r_k]}\big]\subseteq\ws(\eta_k)$ for all $k\geq k_0$. Hence, the inclusion $\ws(\eta_k)\cap\as^m \supseteq \ws(\xi)\cap\as^m$ is derived for $k\geq k_0$. This provides the first step of the convergence of the dictionaries.
\item[(5.IV)\label{(5.IV)}] To deduce the converse inclusion $\ws(\eta_k)\cap\as^m \subseteq \ws(\xi)\cap\as^m$ for $k$ sufficiently large, the coefficients $l_k, r_k\in\NM\,,\; k\in\NM\,,$ have to be chosen with much more care. More precisely, $[-l_k,r_k]\,,\; k\in\NM\,,$ need to satisfy that no defects are created at the boundaries by taking periodic boundary conditions. Here a defect is a pattern of certain length $m$ created at the boundary that does not occur in the original system $\ws(\xi)$, c.f. Example~\ref{Chap4-Ex-StrongConvSpectrConvGenSchrOp} and Example~\ref{Chap5-Ex-PerApprCutNotConv} below. In order to avoid such defects, the coefficients $l_k, r_k\in\NM\,,\; k\in\NM\,,$ are defined via closed paths in the de Bruijn graphs $\gs_k$ of order $k\in\NM$ associated with $\ws(\xi)$.
\end{description}

\begin{remark}
\label{Chap5-Rem-GeneralStrategy}
If a sequence $\Xi_k:=\Orb(\eta_k)\,,\; k\in\NM\,,$ of strongly periodic subshifts is defined according to step \nameref{(5.II)} without any further assumptions, then $(\Xi_k)_{k\in\NM}$ does not converge to $\Xi$ in general, c.f. Example~\ref{Chap5-Ex-PerApprCutNotConv}. More precisely, only the inclusion $\Xi\subseteq\limsup_{k\to\infty}\Xi_k=\bigcap_{k\in\NM}\overline{\big(\bigcup_{m=k}^\infty\Xi_m\big)}$  follows thanks to the argument provided in \nameref{(5.III)}. The inclusion in \nameref{(5.III)} is just a semi-continuity that leads to a semi-continuity of the subshifts. For an associated generalized Schr\"odinger operators $H$, this only leads to a semi-continuity of the spectra. Specifically, $\sigma(H_\Xi)$ is contained in $\limsup_{k\to\infty}\sigma(H_{\Xi_k}) = \bigcap_{k\in\NM}\overline{\big(\bigcup_{m=k}^\infty \sigma(H_{\Xi_m})\big)}$, c.f. Remark~\ref{Chap3-Rem-PContEquivContSpectNormal}.
\end{remark}

The following example shows that, if the sequences $(l_k)_{k\in\NM}$ and $(r_k)_{k\in\NM}$ in \nameref{(5.II)} are not carefully chosen, then the dictionaries $\big(\ws(\eta_k)\big)_{k\in\NM}$ do not converge to $\ws(\xi)$ in the local pattern topology.

\begin{example}
\label{Chap5-Ex-PerApprCutNotConv}
Consider the alphabet $\as:=\{a,b\}$ with two letters and the Fibonacci subshift $\Xi_S\in\SZ\big(\as^\ZM\big)$ is defined via the primitive substitution $S:\as\to\as^+\,,\; a\mapsto ab\,,\; b\mapsto a\,,$ c.f. Example~\ref{Chap2-Ex-NonPeriodic-Fibonacci}. Then there exists a $2$-periodic $\xi\in\as^\ZM$ with respect to the substitution $S$, i.e., $S^2(\xi)=\xi$, such that $\Xi_S=\overline{\Orb(\xi)}$ and $\ws(\Xi)=\ws(\xi)$, c.f. Section~\ref{Chap6-HigherDimPerAppr} for more details. According to Proposition~\ref{Chap6-Prop-PrimSubstSubshZdMin}, the subshift $\Xi_S$ is minimal as $S$ is a primitive substitution. Thus, the letter $b$ occurs infinitely often in $\xi$ to the left and the right with bounded gaps. More precisely, for $k\in\NM$, there are $l_k, r_k\in\NM$ satisfying $\xi(-l_k)=b=\xi(r_k)$ and $\lim_{k\to\infty}l_k=\lim_{k\to\infty} r_k=\infty$. Consequently, the strongly periodic word $\eta_k:=(v_{l_k}|v_{r_k})^\infty\in\as^\ZM$ satisfies $bb\in\ws(\eta_k)$ for each $k\in\NM$ where $v_{l_k}:=\xi_{[-l_k,-1]}$ and $v_{r_k}:=\xi_{[0,r_k]}$. On the other hand, $bb$ is not an element of $\ws(\Xi_S)=\ws(\xi)$. Specifically, $\ws(\eta_k)\cap\as^2\not\subseteq\ws(\xi)\cap\as^2$ follows for $k\in\NM$ since $bb\in\ws(\eta_k)$ holds while $bb\not\in\ws(\xi)$. Corollary~\ref{Chap2-Cor-DictExhSeq} asserts that a sequence $(\ws_k)_{k\in\NM}$ of dictionary converges to a dictionary $\ws$ if and only if $\ws_k\cap\as^{[K_n]}=\ws\cap\as^{[K_n]}$ for $k$ big enough where $(K_n)_{n\in\NM}$ is an exhausting sequence of $G=\ZM$. For $\ZM$, such an exhausting sequence can be chosen such that $\as^{[K_2]}=\as^2$, c.f. Section~\ref{Chap5-Sect-SymbDynSystZM}. Thus, the dictionaries $\big(\ws(\eta_k)\big)_{k\in\NM}$ do not converge to $\ws(\xi)$ by the previous considerations. Hence, the sequence of subshifts $\Xi_k:=\overline{\Orb(\eta_k)}\,,\; k\in\NM\,,$ does not converge to $\Xi_S$ by Theorem~\ref{Chap2-Theo-Shift+DictSpace}.
\end{example}

In accordance with \nameref{(5.IV)} and Example~\ref{Chap5-Ex-PerApprCutNotConv}, the main task is to show that no new words (patterns) are created up to a certain length by taking periodic boundary conditions.

\medskip

Note that there is a unique identification between subshifts and dictionaries, c.f. Theorem~\ref{Chap2-Theo-Shift+DictSpace}. Additionally, a dictionary in $\DZ$ has a unique correspondence to a sequence of de Bruijn graphs $(\gs_k)_{k\in\NM}$, c.f. Definition~\ref{Chap5-Def-DeBruijn} below. Conversely, each sequence of compatible de Bruijn $(\gs_k)_{k\in\NM}$ graphs defines a unique dictionary, c.f. Remark~\ref{Chap5-Rem-SpaceBruijnGraphs}. This fact is used in this chapter to characterize the periodically approximable subshifts of $\as^{\ZM}$, c.f. Theorem~\ref{Chap5-Theo-ExPerAppr}.

\medskip

Throughout this chapter, it is always assumed that the alphabet contains at least two letters since $\as^\ZM$ only contains one element if $\sharp\as=1$. Since the notations simplifies for subshifts over the group $\ZM$, a short reminder of the framework for the case $G=\ZM$ is provided in Section~\ref{Chap5-Sect-SymbDynSystZM}.

\section{Symbolic dynamical systems over the group \texorpdfstring{$\ZM$}{ZM}}
\label{Chap5-Sect-SymbDynSystZM}

For convenience of the reader, basic definitions and notations for symbolic dynamical systems over the group $\ZM$ are recalled in the following, c.f. also Section~\ref{Chap2-Sect-ExampleAsZMSchrOp}. An {\em alphabet} is a {\em finite} set $\as$. A {\em word} with letters in $\as$ is an element of the Cartesian product $\as^n$. For the sake of convenience, a word $u\in\as^n$ is represented as a concatenation of letters $u=a_1 a_2\ldots a_n$ with $a_k\in\as\,,\; 1\leq k\leq n$. The number of letters $|u|:=n$ is called its {\em word length}. The {\em empty word} $\epsilon$ has by definition zero length. Given two words $u=a_1\ldots a_n\in\as^n$ and $v=b_1\ldots b_m\in\as^m$, the {\em concatenation} $uv=a_1\ldots a_n b_1\ldots b_m$ is an element of $\as^{n+m}$ for $n,m\in\NM$. Then the concatenation map $(u,v)\mapsto uv$ is associative, but in general not commutative. The set $\as^\ast:=\bigcup_{n\in\NM}\as^n\cup\{\epsilon\}$ is called the {\em free monoid} on the set $\as$ while $\as^+:=\bigcup_{n\in\NM}\as^n$ is called the {\em free semigroup} on $\as$. For a word $v\in\as^+$ and $j\in\NM$, the {\em $j$-th time concatenation} of the word $v$ is denoted by $v^j:=vv\ldots v\in\as^{j|v|}$. A word $v$ is called a {\em subword} of $u$ if there exist two other words $w_1,w_2\in\as^\ast$ (possibly the empty words) such that $u=w_1vw_2$. This implies $|v|\leq |u|$. In the literature, a subword of $v$ is also called a {\em factor} of $v$. Let $u,v\in\as^\ast$. Then $u$ is called a {\em prefix of $v$} (respectively, a {\em suffix of $v$}) if there exists a finite word $w\in\as^\ast$ such that $uw=v$ (respectively, $wu=v$) holds. 

\medskip

Recall the notion of a pattern introduced in Definition~\ref{Chap2-Def-Pattern}. The compact subsets $\ks(\ZM)$ of $\ZM$ are given by the finite subsets of $\ZM$. Two sets $K,F\in\ks(\ZM)$ are called {\em $\ZM$-equivalent} $(K\sim_\ZM F)$ if there is an $n\in\ZM$ such that $K+n=F$. In analogy to Subsection~\ref{Chap2-Ssect-SpaceDictionaries}, the quotient $\ks(\ZM)/\sim_\ZM$ is defined. Consider the subset 
$$
C \;
	:= \; \big\{
		\{1,\ldots,n\} \;\big|\; n\in\NM 
	\big\} \;
	\subsetneq \ks(\ZM)	\,.
$$
A subset $K\subseteq\ZM$ is called {\em connected} if $K$ is $\ZM$-equivalent to an element of $C$. For instance, $\{3,4,5\}$ is a connected set whereas $\{1,2,5\}\subseteq\ZM$ is not connected. The elements of $C$ are of particular interest as an exhausting sequence $(K_n)_{n\in\NM}$ of $\ZM$ exists satisfying $K_n\sim_\ZM\{1,\ldots,n\}$ for each $n\in\NM$, see below. Then the set of patterns $\as^{[K_n]}$ with support $[K_n]=[\{1,\ldots,n\}]$ is exactly the Cartesian product $\as^n$. Hence, the free semigroup $\as^+$ on $\as$ is equal to $\big\{[u]\in\as^{[K]}\;|\; K\in C\big\}$. In view of that a finite word is a pattern according to Definition~\ref{Chap2-Def-Pattern}. For the sake of simplicity, the notation $u\in\as^n$ is used instead of $[u]\in\as^{[K_n]}$. Then the notions of subpatterns and extensions transfer to words. Specifically, $\Sub[u]$ is the set of all subwords of $u\in\as^n$ and $u$ is an extension of every element of $\Sub[u]$.

\medskip

The finite set $\as$ carries the discrete topology. Then the space $\as^\ZM:=\prod_{j\in\ZM}\, \as$ endowed with the product topology which is called the {\em space of two-sided infinite sequences} over the alphabet $\as$. Due to the finiteness of $\as$ the space $\as^\ZM$ defines a compact, second countable, Hausdorff space, c.f. \cite{Tyc30,Cec37} or Proposition~\ref{Chap2-Prop-FullShift}. Then the continuous action of $\ZM$ on the compact space $\as^\ZM$ is given by $\alpha:\ZM\times\as^\ZM\to\as^\ZM\,,\; \alpha_n(\xi):= \xi(\cdot - n)$. Furthermore, a dynamical subsystem of $(\as^\ZM,\ZM,\alpha)$, i.e., a closed, $\alpha$-invariant subset $\Xi\subseteq\as^\ZM$, is called a subshift, c.f. Section~\ref{Chap2-Sect-SymbDynSyst}. 

\medskip

An element in $\as^\ZM$, i.e., a {\em two-sided infinite word} $\xi$, is a map $\xi:\ZM\to\as$. For each $\xi\in\as^\ZM$ and $l,r\in\ZM$ with $l\leq r$, the restriction $\xi|_{[l,r]}$ denotes the {\em subword} $\xi(l)\xi(l+1)\ldots \xi(r-1)\xi(r)$ of length $r-l+1$. 

\medskip

Let $u=a_1\ldots a_n\in\as^n$ and $v=b_1\ldots b_m\in\as^m$. Denote by $(u|v)^\infty$ the periodic extension of $u|v$ to a two-sided infinite word $\xi\in\as^\ZM$ such that $\xi|_{[-n+k(n+m),(m-1)+k(n+m)]}=uv$ for all $k\in\ZM$, i.e.,
$$
\xi\; 
	:= \; \ldots a_{n}\, b_1\ldots b_{m}\, a_1\ldots a_{n}\, | \, b_1\ldots b_{m}\, a_1\ldots a_{n}\, b_1\ldots
	\; \in\as^\ZM  \, .
$$
In this notation, it is allowed that one of the words $u,v$ is the empty word $\epsilon$. Then the notation $v^\infty$ is used instead of $(\epsilon|v)^\infty = (v|\epsilon)^\infty$. The periodic extensions of words $u,v\in\as^+$ are strongly periodic with respect to the $\ZM$-action, i.e., there exists a $p\in\NM$ such that $\alpha_p(\xi)=\xi$ with $p\leq |v|+|u|$.

\medskip

With the previous considerations at hand, the convergence of dictionaries $\DZ$ is reformulated as follows. Let $K_n:=\{-l_n,\ldots,r_n\}\,,\; n\in\NM\,,$ be a sequence of compact subsets of $\ZM$ where $l_n,r_n\in\NM_0$ satisfies $l_n+r_n+1=n$ for $n\in\NM$ and $\lim_{n\to\infty} l_n=\lim_{n\to\infty} r_n=\infty$. It immediately follows that $(K_n)_{n\in\NM}$ defines an exhausting sequence of $\ZM$. Furthermore, $K_n$ is $\ZM$-equivalent to $\{1,\ldots,n\}$. According to Remark~\ref{Chap2-Rem-DictExhSeq}, every dictionary $\ws\in\DZ$ is uniquely defined by the sets $\ws\cap\as^{[K_n]}=\ws\cap\as^n\,,\; n\in\NM$. Thus, the dictionary $\ws(\xi)$ associated with $\xi\in\as^\ZM$ is identified with the set of all subwords of $\xi$. 

\medskip

A sequence of dictionaries $(\ws_k)_{k\in\NM}$ converges to $\ws\in\DZ$ if and only if, for every $n\in\NM$, there is a $k_0\in\NM$ such that
$$
\ws_k\cap\as^n \;
	= \; \ws\cap\as^n\,,
	\qquad 
	k\geq k_0\,.
$$
Let $U$ be a finite subset of $\as^k$ for a $k\in\NM$. Then a {\em subshift of finite type $\Xi(U)$} associated with $U$ is defined by
$$
\Xi(U) \;
	:=\; \big\{ 
		\eta\in\as^\ZM \;\big|\; \ws(\eta)\cap\as^k \subseteq U
	\big\} \,,
$$
c.f. Remark~\ref{Chap4-Rem-MonotonConvSpectr}. In general, the set $\Xi(U)$ can be empty. If, for instance, $\as:=\{a,b\}$ and $U:=\{ab\}\subseteq\as^2$, then the subshift of finite type $\Xi(U)$ is empty since every $\xi\in\as^\ZM$ that has $ab$ as a subword also contains a subword that begins with the letter $b$.

\medskip

Let $\Xi\in\SZ\big(\as^\ZM\big)$ be a subshift. Then the sequence
$$
\Xi_n \; 
	:= \; \big\{ 
		\eta\in\as^\ZM \;\big|\; \ws(\eta)\cap\as^n\subseteq\ws(\Xi)\cap\as^n	
	\big\}\,,
	\qquad
	n\in\NM\,,
$$
of subshifts of finite type monotonically converges to $\Xi$ by Theorem~\ref{Chap2-Theo-Shift+DictSpace}, i.e., $\Xi_{n+1}\subseteq\Xi_n$ holds for $n\in\NM$ and $\bigcap_{n\in\NM}\Xi_n=\Xi$, c.f. Remark~\ref{Chap4-Rem-MonotonConvSpectr}.

\medskip

\label{Page-ExStrongPerSubshFinitType}
According to the definition of the local pattern topology and Theorem~\ref{Chap2-Theo-Shift+DictSpace}, the subshift $\Xi$ is (finitely) periodically approximable if and only if all the subshifts of finite type $(\Xi_n)_{n\in\NM}$ contain sufficiently many strongly periodic elements. The question whether a subshift of finite type contains strongly periodic elements is a current question in the field of dynamical systems, c.f. \cite{Pia08,Hoc09,Fio09,SiCo12,Coh14,CaPe15}. It is known \cite{Fio09} that there are sufficient conditions for one-dimensional subshifts of finite type to contain strongly periodic el\-ements. The proof relies on the de Bruijn graphs associated with a subshift of finite type. More precisely, closed paths give rise to strongly periodic elements in the subshift of finite type. 

\medskip

A similar strategy is followed here but with a different perspective. Specifically, \cite{Fio09} studies a fixed subshift of finite type whereas the whole topological space $\SZ\big(\as^\ZM\big)$ is analyzed in this thesis. More precisely, a subshift $\Xi\in\SZ\big(\as^\ZM\big)$ is considered which is not a subshift of finite type, in general. Then a sequence of subshifts of finite type $(\Xi_n)_{n\in\NM}$ is defined via $U_n:=\ws(\Xi)\cap\as^n\,,\; n\in\NM$, c.f. the previous considerations. Then strongly periodic elements $\eta_n\in\Xi_n$ with $\ws(\eta_n)\cap\as^n=U_n$ lead to strongly periodic subshifts $\Orb(\eta_n)\subseteq\Xi_n$ converging to $\Xi$.

\medskip

In higher dimensions, the questions of strongly periodic elements in a subshift of finite type is more difficult, c.f. \cite{Hoc09}. Furthermore, whenever a subshift of finite type contains one strongly periodic element then it contains already many strongly periodic elements under suitable assumptions. Specifically, the set of strongly periodic elements of the subshift of finite type is a dense subset, c.f. \cite{Lig03,SiCo12} and Remark~\ref{Chap2-Rem-SubsFinTypPerSpectr}.

\medskip

The existence of strongly periodic elements in a subshift of finite type corresponds also to the existence of so called Wang tiles. More precisely, {\sc Wang} \cite{Wan61} was the first who asked the question whether, for a given finite set of tiles, there is a periodic tessellation of the plane with these tiles. {\sc Berger} \cite{Ber66}, a student of {\sc Wang}, provides an example of tiles that do not tile the plane in a periodic way. Such tilings are called nowadays {\em Wang tiles}. Even in the case of a lattice $\ZM^d$ for $d\geq 3$, there exist tilings, i.e., $\xi:\ZM^d\to\as$, that do not allow any periodic tiling, c.f. \cite{CuKa95}. This shows that there exist subshifts that are not periodically approximable.

\begin{remark}
\label{Chap5-Rem-GAP}
The original starting point original of \cite{BeBeNi16} was the work of {\sc Anderson} and {\sc Putnam} \cite{AnPu98} who defined CW-complexes associated with tilings arising by a substitution that are called {\em Anderson-Putnam complexes} nowadays. They show that the original dynamical system is recovered by an inverse limit which was used to compute the cohomology and K-theory of the tiling space. This approach was then extended by {\sc Bellissard}, {\sc Benedetti} and {\sc Gambaudo} to general finite pattern tilings, c.f. \cite{BeGa03,BeBeGa06}.  {\sc G\"ahler} \cite{Gah02} proposed a different approach by increasing the collar instead of the tiles. This approach was taken up in the works of {\sc Sadun} \cite{Sad03,Sad08}. It turns out that the CW-complexes for subshifts of $\as^\ZM$ are essentially isomorphic as graphs to the de Bruijn graphs. For minimal Cantor sets, the work \cite{GaMa06} by {\sc Gambaudo} and {\sc Martens} also provides a sequence of graphs that describe a dynamical system by taking a projective limit. Also these graphs coincide with the de Bruijn graphs if the minimal dynamical is a subshift over an alphabet, c.f. the following discussion.

\vspace{.1cm}

The Anderson-Putnam complexes were used to calculate the \v{C}hech cohomology \cite{AnPu98,GaHuKe05,GaHuKe08,Sad08} and the Gap Labeling Theorem \cite{Bel92,BeOy02,BeOy07,KaPu03,BeBeGa06}. These complexes are also important for the existence of strongly periodic approximations of subshifts in higher dimensions, c.f. Remark~\ref{Chap6-Rem-StrategySubstitution} and the discussion in Section~\ref{Chap8-Sect-APComplex}.
\end{remark}

\section{De Bruijn graphs}
\label{Chap5-Sect-DBGraphs}

This section deals with the so called de Bruijn graphs associated with a dictionary $\ws\in\DZ$. These graphs describe the local structures of the associated dynamical system $\Xi(\ws)$. Furthermore, these graphs allow one to define sequences of strongly periodic subshifts such that they converge in $\SZ\big(\as^\ZM\big)$ to $\Xi$.

\medskip

\begin{figure}[htb]
\centering
\includegraphics[scale=0.95]{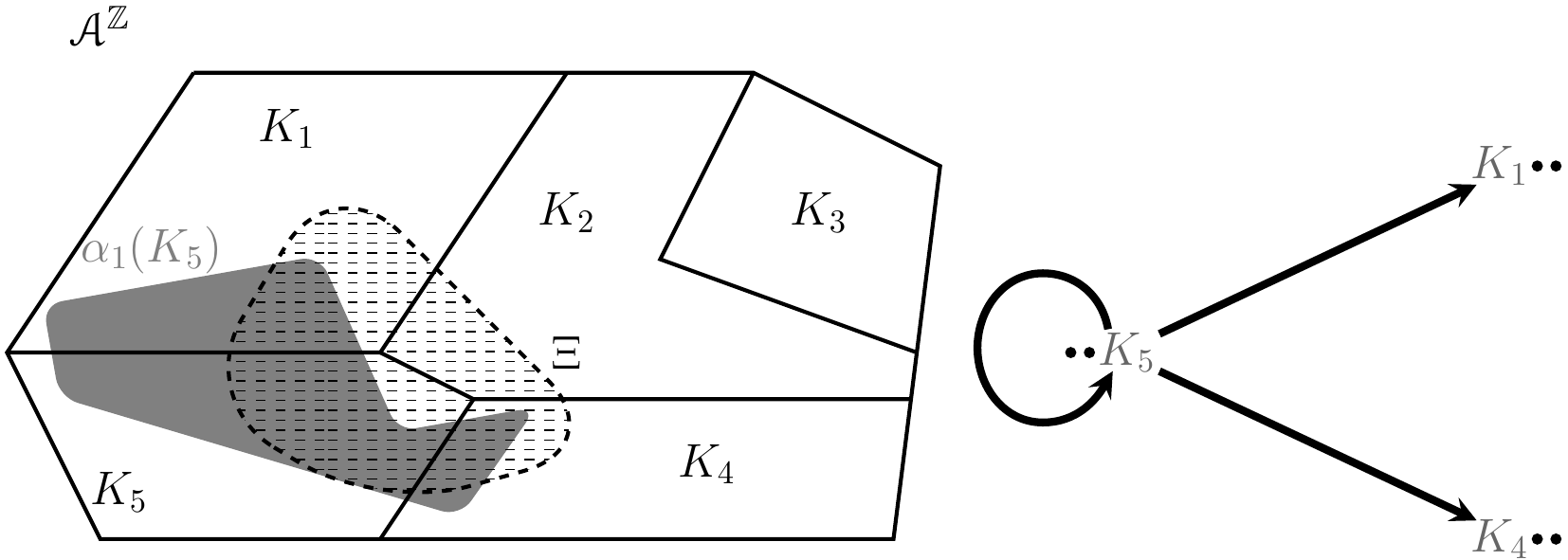}
\caption{Sketch of a partition of $\as^\ZM$ and the shift $\alpha_1(K_5)$. The vertex set of $\gs_5$ associated with $\Xi$ is given by $\{ K_1,\, K_2,\, K_4,\, K_5 \}$. Since $\big(K_1\cap\alpha_1(K_5)\big)\cap\Xi\neq\emptyset$ and $\big(K_4\cap\alpha_1(K_5)\big)\cap\Xi\neq\emptyset$, there is an edge joining $K_5$ and $K_1$ and an edge joining $K_5$ and $K_4$. Additionally, there is an edge linking $K_5$ with itself. Thus, $K_5$ is a branching vertex of $\gs_5$, c.f. Subsection~\ref{Chap5-Ssect-GraphRemind}.}
\label{Chap5-Fig-Partition}
\end{figure}

The main idea for the construction of the de Bruijn graph is described next. Later, an equivalent definition of the de Bruijn graphs associated with a dictionary is provided since a dictionary is uniquely identified with a subshift by Theorem~\ref{Chap2-Theo-Shift+DictSpace}. Let $\alpha:\ZM\times\as^\ZM\to\as^\ZM\,,\; \alpha_n(\xi):=\xi(\cdot-n)\,,$ be the continuous action of $\ZM$ on $\as^\ZM$. The topological space $\as^\ZM$ is a compact space. Let $(K_n)_{n=1}^N$ be a finite partition of $\as^\ZM$, i.e., $K_n\subseteq\as^\ZM\,,\; 1\leq n\leq N\,,$ are clopen and pairwise disjoint sets satisfying $\bigsqcup_{n=1}^N K_n =\as^\ZM$. Then a directed graph $\gs_N$ associated with a subshift $\Xi\in\SZ\big(\as^\ZM\big)$ is constructed as follows: 
\begin{description}
\item[(i)] The vertex set $\vs_N$ is given by $\big\{K_j\;\big|\; K_j\cap\Xi\neq\emptyset\,,\; 1\leq j\leq N \big\}$. 
\item[(ii)] There is an edge joining $K_i$ to $K_j$ if the intersection $\big(\alpha_1(K_i)\cap K_j\big)\cap\Xi$ is non-empty, c.f. Figure~\ref{Chap5-Fig-Partition}.
\end{description}

Closed paths that visit all vertices in such a graph give rise to strongly periodic elements $\eta$ of $\as^\ZM$ such that $\eta$ shares the same words with $\Xi$ up to a certain length. The element $\eta$ is defined via the associated dictionary. In view of Section~\ref{Chap5-Sect-SymbDynSystZM}, this implies that the subshift $\Orb(\eta)$ is close to $\Xi$ in the Hausdorff-topology of $\SZ\big(\as^\ZM\big)$. The correct generalization of these graphs is the Anderson-Putnam complex \cite{AnPu98,BeGa03,BeBeGa06} for the group $\ZM^d$ or even more general for a tiling in $\RM^d$, c.f. Remark~\ref{Chap5-Rem-GAP}, Remark~\ref{Chap6-Rem-StrategySubstitution} and Section~\ref{Chap8-Sect-APComplex}. In this case, a torus in the Anderson-Putnam complex is associated with a strongly periodic subshift close to the original system.

\subsection{Graphs a reminder}
\label{Chap5-Ssect-GraphRemind}

This section recalls some basic concepts of directed graphs, c.f. \cite{Bol98}. In particular, the notion of strongly connected graphs is discussed which plays an important role for the existence of strongly periodic elements, c.f. Theorem~\ref{Chap5-Theo-ExPerAppr}. 

\medskip

\begin{figure}[htb]
\centering
\includegraphics[scale=1.02]{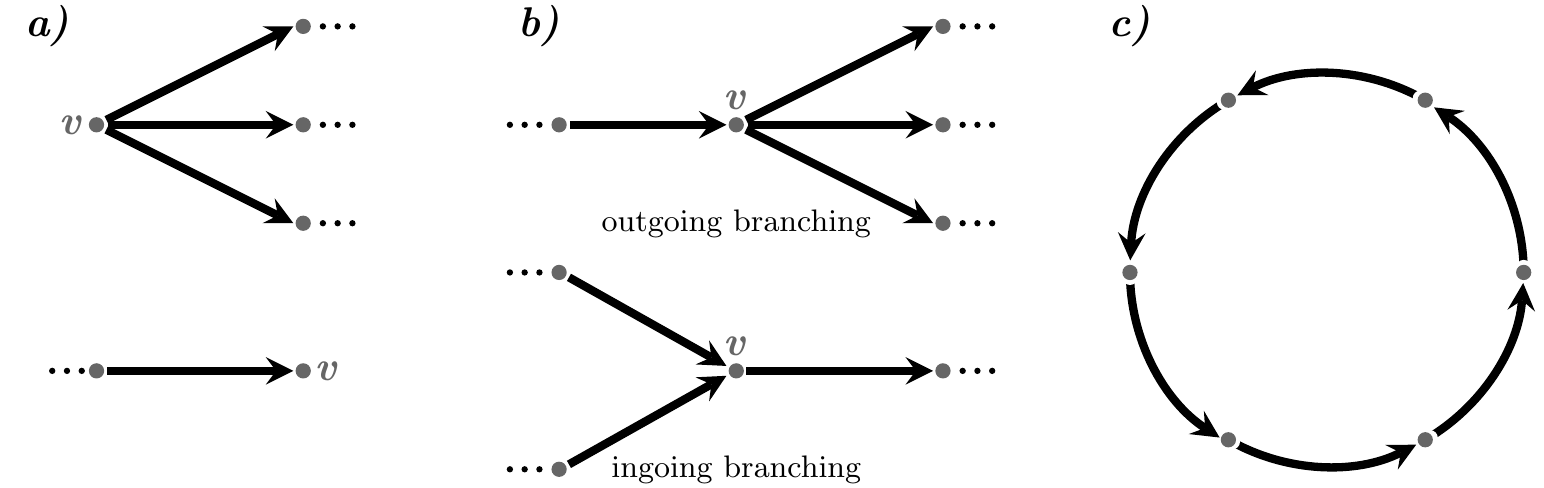}
\caption{a) dandling vertices $v\,,$ b) branching vertices $v\,,$ c) a directed graph where each vertex is not dandling and not branching.}
\label{Chap5-Fig-DandlingBranching}
\end{figure}

A {\em directed graph} is a family $\gs=(\vs,\es,\partial_0,\partial_1)$ where $\vs$ and $\es$ are discrete countable sets and the maps $\partial_i:\es\to\vs\,,\; i=0,1\,,$ are called {\em boundaries}. An element of $\vs$ is called a {\em vertex} while an element of $\es$ is called an {\em edge}. An edge $e\in\es$ is an {\em arrow} joining its {\em origin} $\partial_0\,e$ to its {\em end} $\partial_1\,e$. An edge $e$ is called {\em outgoing} from the vertex $v$ if $\partial_0\,e=v$ while the edge $e$ is called {\em incoming} if $\partial_1\,e=v$. The {\em vertex degree} ${\rm deg}(v)$ is defined by the sum of the number of incoming and outgoing edges from the vertex $v\in\vs$. A vertex is called {\em dandling} if it has either only incoming or only outgoing edges. If $v\in\vs$ is not dandling and ${\rm deg}(v)>2$, the vertex $v$ is called {\em branching}, c.f. Figure~\ref{Chap5-Fig-DandlingBranching}.

\medskip

As previously mentioned, directed graphs are associated with a dynamical system. Then paths in a directed graph are interpreted as a flow along the dynamics. Thus, the structure of the graphs encode dynamical properties.

\medskip

A graph is called {\em finite} if the vertex set and the edge set are finite. Clearly, a finite graph is a circle if all vertices are not dandling and not branching, c.f. Figure~\ref{Chap5-Fig-DandlingBranching}~c).

\medskip

A graph is called {\em simple} whenever
\begin{description}
\item[(i)] the source $\partial_0(e)$ and the range $\partial_1(e)$ are not equal for all edges $e\in\es$,
\item[(ii)] for every pair of different vertices $u,v\in\vs$, there exists at most one edge $e$ with origin $u$ and end $v$, i.e., edges joining $u$ to $v$ are unique if there is an edge. 
\end{description}

\begin{figure}[htb]
\centering
\includegraphics[scale=1]{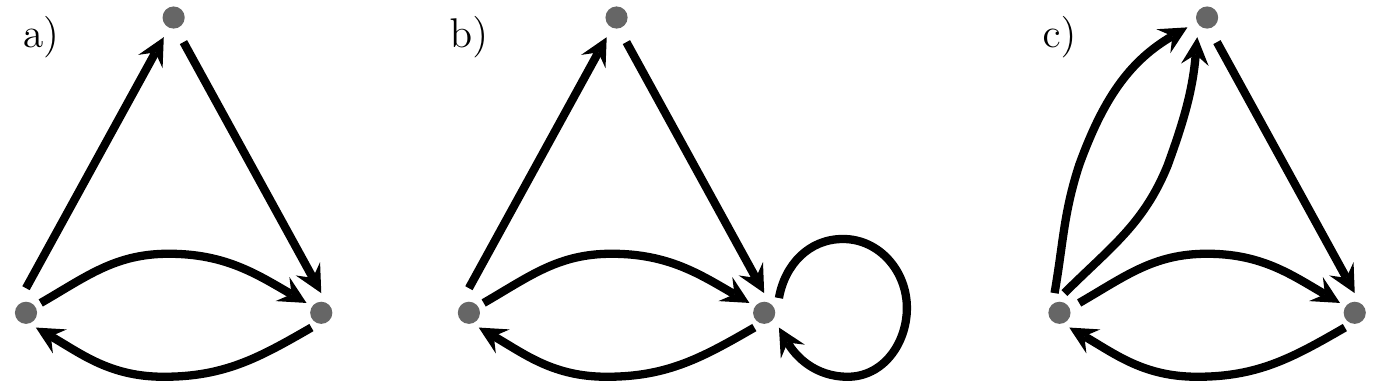}
\caption{a) a simple graph, b) a semi-simple graph that is not simple, c) a graph that is not semi-simple.}
\label{Chap5-Fig-GraphSimpleSemi}
\end{figure}

A graph is called  {\em semi-simple} if edges linking one vertex to itself are not excluded, i.e., the graph only satisfies condition (ii), c.f. Figure~\ref{Chap5-Fig-GraphSimpleSemi}~b).

\medskip

A {\em path} is a finite chain $\wp=(e_1,e_2,\ldots, e_k)$ of edges such that $\partial_1(e_j)=\partial_0(e_{j+1})$ for $1\leq j\leq k-1$. The number of edges $|\wp|:=k$ in a path $\wp$ is called the {\em length of $\wp$}. A path $(e_1,\ldots,e_k)$ is called {\em closed path} whenever $\partial_0 e_1=\partial_1 e_{k}$. A closed path $(e_1,\ldots,e_k)$ is called {\em global} in $\gs=(\vs,\es,\partial_0,\partial_1)$ if $\{\partial_0(e_j)\;|\; 1\leq j\leq k\}$ is equal to $\vs$, namely the closed path visits all vertices of the graph.

\medskip

\begin{figure}[htb]
\centering
\includegraphics[scale=1]{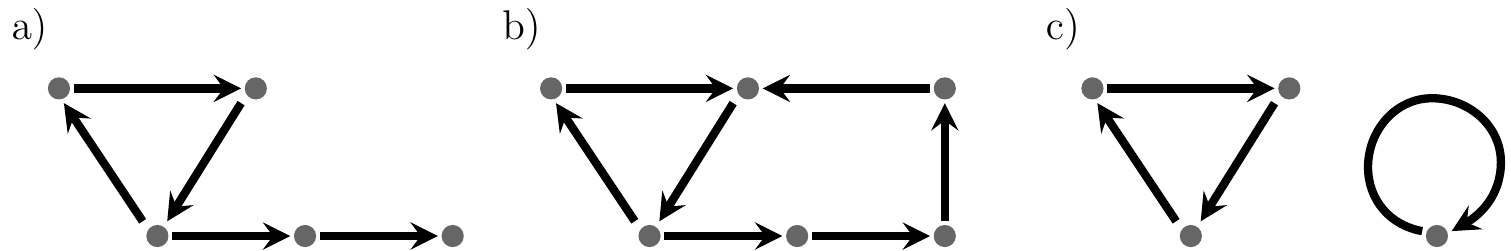}
\caption{a) a connected graph that is not strongly connected, b) a strongly connected graph, c) a graph that is not connected.}
\label{Chap5-Fig-StronglyConnectedGraph}
\end{figure}

A graph $\gs$ is called {\em connected} if there is at least one path connecting $u$ and $v$ for every pair of distinct vertices $u,v\in\vs$, i.e., there exists a path $(e_1,\ldots,e_k)$ such that $\partial_0(e_1)=u,\,\partial_1(e_k)=v$ or $\partial_0(e_1)=v,\,\partial_1(e_k)=u$, c.f. Figure~\ref{Chap5-Fig-StronglyConnectedGraph}. Connected graphs are not sufficient for our purposes. More precisely, the graphs need to be strongly connected meaning that there exist paths connecting $u$ to $v$ and $v$ to $u$ for each pair $u, v$ of vertices.

\begin{definition}[Strongly connected]
\label{Chap5-Def-StronglyConnectedGraph}
A directed graph $\gs:=(\vs,\es,\partial_0,\partial_1)$ is called {\em strongly connected} if, for every given pair of vertices $u,v\in\vs$, there are paths $(e_1,\ldots,e_k)$ and $(f_1,\ldots,f_l)$ such that $\partial_0(e_1)=u=\partial_1(f_l)$ and $\partial_1(e_k)=v=\partial_0(f_1)$.
\end{definition}

A strongly connected graph is connected by definition whereas the converse is false in general, c.f. Figure~\ref{Chap5-Fig-StronglyConnectedGraph}. The section is finished with a characterization for a graph being strongly connected. This characterization is crucial for the further purposes.

\begin{proposition}
\label{Chap5-Prop-GlobalPath}
Let $\gs=(\vs,\es,\partial_0,\partial_1)$ be a finite directed graph. Then there exists a global path $\wp=(e_1,e_2,\ldots, e_k)$ if and only if $\gs$ is strongly connected.
\end{proposition}

\begin{proof}
This immediately follows by definition of a strongly connected graph and the fact that $(e_j,\ldots, e_k,e_1,\ldots e_{j-1})$ defines also a closed path in $\gs$ for any $1\leq j\leq k$.
\end{proof}

\subsection{De Bruijn graphs and the relation to periodicity and aperiodicity}
\label{Chap5-Ssect-PropDeBruijn}

In 1894, {\sc Flye} introduced a graph representing possible continuation of words of finite lengths in \cite{Fly1894}. {\sc De Bruijn} \cite{Bru46} and {\sc Good} \cite{Goo46} specified the construction of these graphs in 1946, independently. However, the community decided to call these graphs {\em de Bruijn graphs}. {\sc Rauzy} \cite{Rau83} provided the first use of the de Bruijn graphs where he analyzed the combinatorial structure of minimal sequences with three letters such that the subword complexity $p(n)$ is equal to $2n+1$ for $n\in\NM$. In view of that these graphs are also called {\em Rauzy graphs} elsewhere \cite{Cas97,Jul10}. As mentioned in Remark~\ref{Chap5-Rem-GAP}, these graphs are isomorphic (in terms of graph isomorphisms) to the one-dimensional version of the Anderson-Putnam complexes defined in \cite{AnPu98}, see \cite{Gah02,Sad03,BeGa03,BeBeGa06,Sad08} as well. Thus, the Anderson-Putnam complex is the correct generalization for subshifts defined on $\ZM^d$ or for tiling spaces in $\RM^d$, c.f. Section~\ref{Chap8-Sect-APComplex}.

\medskip

These graph are also called {\em Nth higher edge graph} in \cite[Definition~2.3.4]{LindMarcus95} for $N\in\NM$. There it is proven that these graphs associated with a sofic subshift are strongly connected if and only if the subshift is irreducible, c.f. \cite[Lemma~3.3.10, Proposition~2.2.14]{LindMarcus95}. In \cite{More05}, the existence of the so called de Bruijn sequences is studied for a given set of words of the same length. A {\em de Bruijn sequence} is a closed path in the de Bruijn graph such that every vertex is visited exactly once. The characterization of a strongly connected de Bruijn graph by the fact that the subshift is irreducible is also proven there \cite[Lemma~9]{More05}. Note that in \cite{More05} a different notion of a dictionary is used than here. More precisely, a dictionary is a set of words of a fixed length $n\in\NM$. However, this result directly transfers to our notion of a dictionary in terms of Definition~\ref{Chap2-Def-Dictionary}. The concept of strongly connected graphs is used to construct strongly periodic elements (Definition~\ref{Chap5-Def-AssPerWor}) that are close in the Hausdorff-topology to a subshift $\Xi\in\SZ\big(\as^\ZM\big)$. The idea of using de Bruijn graphs to construct strongly periodic two-sided infinite words is not new, c.f. \cite{Fio09} and the discussion in Section~\ref{Chap5-Sect-SymbDynSystZM}.

\medskip

In the following the de Bruijn graphs arising by a subshift are studied. The behavior of closed paths as well as the existence of global paths is presented, i.e., sufficient conditions so that the de Bruijn graphs are strongly connected.

\medskip

The de Bruijn graphs are defined for a fixed dictionary having in mind that every dictionary is uniquely associated with a subshift by Theorem~\ref{Chap2-Theo-Shift+DictSpace}.

\begin{definition}
\label{Chap5-Def-DeBruijn}
Let $\ws\in\DZ$ be a dictionary. For $k\in\NM$, define the vertex set $\vs_k:=\ws\cap\as^k$ and the edge set $\es_k:=\ws\cap\as^{k+1}$. The boundary maps $\partial_0\,,\,\partial_1:\es_k\to \vs_k$ are defined by
$$\partial_0(a_0a_1\ldots a_k):=
   a_0a_1\ldots a_{k-1}\,,
    \hspace{2cm}
     \partial_1(a_0a_1\ldots a_k):=
      a_1a_2\ldots a_k\,.
$$
The oriented graph $\gs_k:=\gs_k^\ws:=(\vs_k,\es_k,\partial_0,\partial_1)$ is called the {\em de Bruijn graph of order $k$ associated with the dictionary $\ws$}.
\end{definition}

Clearly, the de Bruijn graph $\gs_k^\ws$ of order $k\in\NM$ associated with $\ws\in\DZ$ only depends on the set $\ws\cap\as^{k+1}$. Thus, the de Bruijn graphs $\gs_k^\ws$ and $\gs_k^{\ws'}$ agree if the equality $\ws\cap\as^{k+1}=\ws'\cap\as^{k+1}$ holds for the dictionaries $\ws,\ws'\in\DZ$.

\begin{figure}[htb]
\centering
\includegraphics[scale=0.933]{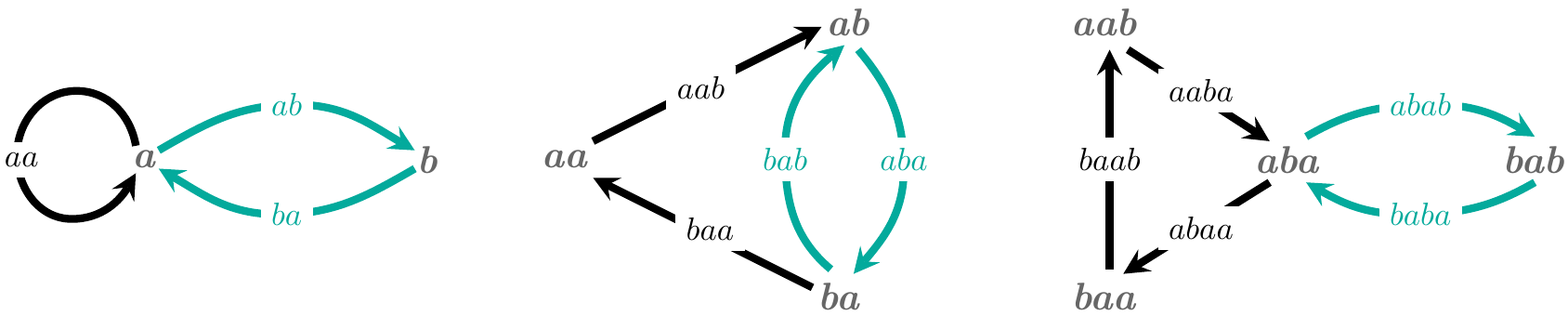}
\caption{The de Bruijn graphs of order $1$, $2$ and $3$ associated with the Fibonacci sequence.}
\label{Chap5-Fig-Fibonacci}
\end{figure}

\begin{example}[Fibonacci sequence]
\label{Chap5-Ex-Fibonacci}
Recall Example~\ref{Chap2-Ex-NonPeriodic-Fibonacci} of an aperiodic subshift $\Xi_S$. Specifically, the alphabet $\as$ consists of two letters $a$ and $b$ and the subshift $\Xi_S$ was defined via the primitive substitution $S:\as\to\as^+\,,\; a\to ab\,,\,b\to a$, c.f. Definition~\ref{Chap6-Def-BlockPrimSubst}. Thus, the associated subshift $\Xi_S$ is minimal by Proposition~\ref{Chap6-Prop-PrimSubstSubshZdMin}. This subshift is called the {\em Fibonacci subshift}. Denote by $\ws:=\ws(\Xi_S)$ the associated dictionary with the subshift $\Xi_S$. The corresponding subword complexity $\comp:\NM\to\NM\,,\; \comp(k):=\sharp\big(\ws\cap\as^k\big)\,,$ is given by $\comp(k)=k+1\,,\; k\in\NM$, see e.g. \cite[Corollary~5.4.10]{Fogg02}. Thus, the associated de Bruijn graph of order $k\in\NM$ has exactly $k+1$ vertices and $k+2$ edges. The de Bruijn graphs of order $1$, $2$ and $3$ are given in Figure~\ref{Chap5-Fig-Fibonacci}.
\end{example}

\begin{remark}
\label{Chap5-Rem-ConstraintBruijnGrWords}
It is important to notice that the boundary maps $\partial_0,\partial_1$ provide a constraint on the neighboring vertices and edges. In fact, for edges $e:=a_0a_1\ldots a_k,\, f:=b_0b_1\ldots b_k\in\es_k$, the condition $\partial_1(e)=\partial_0(f)$ implies  that $a_{i+1}=b_i$ for all $i=0,\ldots, k-1$. Thus, an edge $e\in\es_k$ encodes the possible continuation of the word $\partial_0(e)\in\vs_k$ to the right or the possible continuation of the word $\partial_1(e)\in\vs_k$ to the left. Consequently, if there is exactly one outgoing (respectively, ingoing) edge $e\in\es_k$ of a vertex $v\in\vs_k$, then the word $v$ of length $k$ has a unique continuation to the right (respectively, to the left). This observation emphasizes the importance of the branching points. 

\vspace{.1cm}

For instance, in the case of an outgoing branching vertex, i.e., $v\in\vs_k$ has at least two outgoing edges, the word $v$ does not have a unique continuation to the right and similarly for an ingoing branching vertex. This is observed for non-periodic, infinite words. In particular, it is impossible to reconstruct uniquely a non-periodic word $\xi\in\as^\ZM$ out of the knowledge of a finite subwords of $\xi$. In contrast, strongly periodic elements $\xi\in\as^\ZM$ are determined by the finite word $\xi|_{[0,p-1]}$ where $p\in\NM$ is the period of $\xi$, c.f. Lemma~\ref{Chap6-Lem-ReprStrongPer}. In view of that Proposition~\ref{Chap5-Prop-ClosedPathBehaviour} asserts that the de Bruijn graphs associated with non-periodic systems have branching vertices whereas the de Bruijn graphs $\gs_k$ of order $k$ associated with a strongly periodic subshift is a circles for $k\in\NM$ sufficiently large, c.f. Proposition~\ref{Chap5-Prop-ClosedPathBehaviour}~(d).

\vspace{.1cm}

The behavior in $k\in\NM$ of the number of branching vertices in the de Bruijn graphs $(\gs_k^\ws)_{k\in\NM}$ depends on the subword complexity function $\comp:\NM\to\NM$ of the dictionary $\ws\in\DZ$ which is studied in Section~\ref{Chap5-Ssect-SubwordComplexity}.
\end{remark}

Having this in mind the following key observation follows.

\begin{lemma}[\cite{BeBeNi16}]
\label{Chap5-Lem-WordExPath}
Let $\as$ be an alphabet, $\ws\in\DZ$ be a dictionary and $k\in\NM$. Every word $w:=w_1\ldots w_m\in\ws\cap\as^m$ for $m\geq k+1$ induces a path $(e_1,\ldots,e_{m-k})$ in the de Bruijn graph $\gs_k$ of order $k$ associated with $\ws$ such that $\partial_0(e_1)=w_1\ldots w_k$ and $\partial_1(e_{m-k})=w_{m-k+1}\ldots w_m$.
\end{lemma}

\begin{proof}
Define $e_j:=w_j\ldots w_{j+k}\in\as^{k+1}$ for $1\leq j\leq m-k$ which are subwords of $w$. Hence, $e_j\in\Sub[w]\subseteq\ws$ follows for $1\leq j\leq m-k$ by \nameref{(D1)}. Then the equalities $\partial_1(e_j)=\partial_0(e_{j+1})\,,\; 1\leq j\leq m-k-1\,,$ are derived by construction. Altogether, $(e_1,\ldots, e_{m-k})$ is a path in the de Bruijn graph $\gs_k$ of order $k$ with the desired properties.
\end{proof}

\medskip

Let $u,v\in\as^+$ be finite words over the alphabet $\as$. Then $u$ is called a {\em prefix of $v$} (respectively, a {\em suffix of $v$}) if there exists a finite word $w\in\as^\ast$ such that $uw=v$ (respectively, $wu=v$) holds. Here $\epsilon$ denotes the empty word and so $v$ is a prefix and a suffix of itself.

\newpage
\begin{proposition}[\cite{BeBeNi16}]
\label{Chap5-Prop-GAPBasProper}
Let $\ws\in\DZ$ be a dictionary. For each $k\in\NM$, the de Bruijn graph $\gs_k$ associated with $\ws$ is semi-simple with no dandling vertices. If, additionally, $\Xi(\ws)$ is topologically transitive, i.e., there exists a $\xi\in\as^\ZM$ with $\ws(\xi)=\ws$, then $\gs_k$ is connected for each $k\in\NM$.
\end{proposition}

\begin{proof}
Let $\ws$ be a dictionary with associated de Bruijn graphs $(\gs_k)_{k\in\NM}$. The proof is organized as follows: It is proven that the graph $\gs_k$ for fixed $k\in\NM$ is (i) semi-simple, (ii) has no dandling vertices and (iii) is connected if $\ws=\ws(\xi)$ for $\xi\in\as^\ZM$.

\vspace{.1cm}

(i): Let $u\neq v$ be vertices in $\vs_k$ with an edge $e=a_0\ldots a_{k}\in\es_k$ connecting them, i.e., $u=a_0\ldots a_{k-1}=\partial_0(e)$ and $v=a_1\ldots a_k=\partial_1(e)$. Thus, every other edge linking $u$ to $v$ must be equal to $e$ meaning that $\gs_k$ is semi-simple. 

\vspace{.1cm}

(ii): For each $u\in\vs_n$, there exist an extension $e\in\ws\cap\as^{n+1}$ to the left and $f\in\ws\cap\as^{n+1}$ to right by \nameref{(D2)}, i.e., $u$ is a suffix of $e$ and a prefix of $f$. Hence, $e$ and $f$ are edges of $\gs_k$ satisfying $\partial_1(e)=u=\partial_0(f)$. Thus, $\gs_k$ has no dandling vertices. 

\vspace{.1cm}

(iii): Let $u, v\in\vs_k=\ws(\xi)\cap\as^k$ be two distinct vertices. Then there is an $n_u,n_v\in\ZM$ such that $\xi|_{[n_u,n_u+k-1]}=u$ and $\xi|_{[n_v,n_v+k-1]}=v$. Without loss of generality, suppose that $n_u\leq n_v$. Then $u$ is a prefix and $v$ is a suffix of the word $w:=\xi|_{[n_u,n_v+k-1]}\in\ws(\xi)$. Lemma~\ref{Chap5-Lem-WordExPath} implies that there is a path joining $u$ to $v$ in $\gs_k$. Thus, $\gs_k$ is connected.
\end{proof}

\medskip

Note that there exist dictionaries so that the associated de Bruijn graphs $(\gs_k)_{k\in\NM}$ are not connected. For instance, consider $\Xi=\{\xi,\eta\}\in\SZ\big(\as^\ZM\big)$ defined over the alphabet $\as:=\{a,b\}$ where $\xi:=a^\infty$ and $\eta:=b^\infty$. Then all the de Bruijn graphs associated with $\ws(\Xi)$ are not connected. The following example shows that de Bruijn graphs are also not strongly connected in general.

\begin{figure}[htb]
\centering
\includegraphics[scale=0.67]{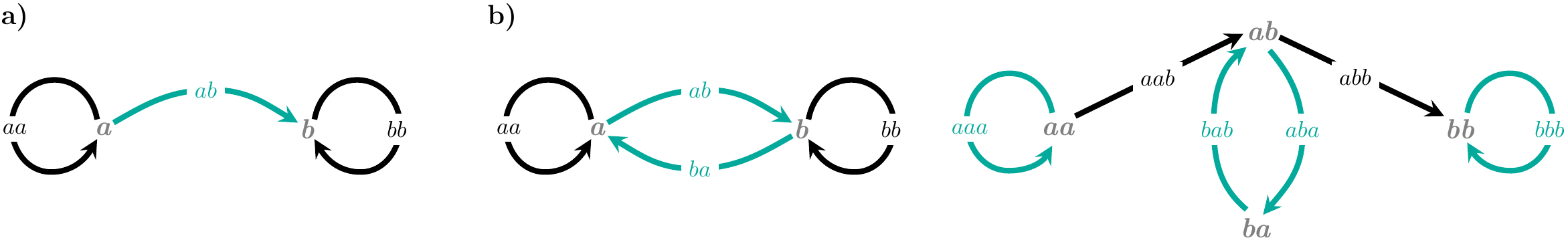}
\caption{a) The de Bruijn graph of order $1$ corresponding to Example~\ref{Chap5-Ex-DeBruijnNotStrongConnect}, b) The de Bruijn graphs of order $1$ and $2$ corresponding to  Example~\ref{Chap5-Ex-deBruijnEventuallyNotStrongConn}.}
\label{Chap5-Fig-BruijnNotStrongConn}
\end{figure}

\begin{example}
\label{Chap5-Ex-DeBruijnNotStrongConnect}
Let $\as:=\{a,b\}$ be the alphabet and $\xi\in\as^\ZM$ be defined by
$$
\xi \; 
	= \; \begin{cases}
		a\,, \qquad & j\leq -1\,,\\
		b\,, \qquad & j\geq  0\,.
	\end{cases}
$$
Then all the de Bruijn graphs are only connected but not strongly connected. This follows by Proposition~\ref{Chap5-Prop-HereditaryBruijn} and the fact that the de Bruijn graph of order $1$ is not strongly connected, c.f. Figure~\ref{Chap5-Fig-BruijnNotStrongConn}~a). Furthermore, the subshift $\Xi:=\overline{\Orb(\xi)}$ is isolated in the Hausdorff-topology of $\SZ\big(\as^\ZM\big)$ by using Theorem~\ref{Chap2-Theo-Shift+DictSpace} since the open set $\{\ws\in\DZ\;|\; \ws\cap\as^2=\ws(\Xi)\cap\as^2\}$ in the local pattern topology only contains $\ws(\Xi)$. Hence, $\Xi$ is also not periodically approximable as it is isolated in the Hausdorff topology of $\SZ\big(\as^\ZM\big)$.
\end{example}

According to Theorem~\ref{Chap5-Theo-ExPerAppr}, strongly connected de Bruijn graphs associated with a subshift $\Xi\in\SZ\big(\as^\ZM\big)$ give rise to strongly periodic subshifts close to $\Xi$ in the Hausdorff-topology of $\SZ\big(\as^\ZM\big)$. In the following a hereditary property is studied.

\begin{proposition}[\cite{BeBeNi16}]
\label{Chap5-Prop-HereditaryBruijn}
Let $\ws\in \DZ$ be a dictionary and $\gs_k$ be the related de Bruijn graph of order $k\in\NM$. If, for a given $k_0\in\NM$, the de Bruijn graphs $\gs_{k_0}$ is strongly connected, then all de Bruijn graphs $\gs_{k}$ of order $k<k_0$ are strongly connected as well.
\end{proposition}

\begin{proof}
Let $k_0\in\NM$ be such that $\gs_{k_0}$ is strongly connected. Consider a $k<k_0$ and $u,v\in\vs_k=\ws\cap\as^k$. Since $\ws$ satisfies \nameref{(D2)}, there are two words $\tilde{u}=a_1\ldots a_{k_0}\in\vs_{k_0}$ and $\tilde{v}=b_1\ldots b_{k_0}\in\vs_{k_0}$ contained in $\ws$ such that $a_1\ldots a_k=u$ and $b_{k_0-k+1}\ldots b_{k_0}=v$. Thus, $u$ is a prefix of $\tilde{u}$ and $v$ is a suffix of $\tilde{v}$. Since $\tilde{u},\tilde{v}\in\vs_{k_0}$, there exists by hypothesis a path $(\tilde{e}_1,\ldots, \tilde{e}_l)$ in $\gs_{k_0}$ such that $\partial_0(\tilde{e}_1)=\tilde{u}$ and $\partial_1(\tilde{e}_l)=\tilde{v}$. Recall that an edge $e\in\es_{k_0}$ is a word $e(0)e(1)\ldots e(k_0)$ of length $k_0+1$. This path induces a path joining $u$ to $v$ in $\gs_k$ via 
\begin{gather*}
\begin{aligned}
e_j \;
	&:= \; \tilde{e}_1(j)\ldots \tilde{e}_1(j+k)\,,
	\qquad 
	&0\leq j \leq k_0-k\,,\\
e_{k_0-k+i} \; 
	&:= \; \tilde{e}_{i}(k_0-k)\ldots \tilde{e}_{i}(k_0)\,,
	\qquad 
	&2\leq i \leq l\,.	
\end{aligned}
\end{gather*} 
Then $(e_1,\ldots, e_{k_0-k+l})$ is a path in $\gs_k$ satisfying $\partial_0(e_1)=u$ and $\partial_1(e_{k_0-k+l})=v$.

\vspace{.1cm}

Similarly, there are $u',v'\in\vs_{k_0}$ such that $u$ is a suffix of $u'$ and $v$ is a prefix of $v'$. Since $\gs_{k_0}$ is strongly connected, there exists a path $(f'_1,\ldots,f'_m)$ in $\gs_{k_0}$ satisfying $\partial_0(f'_1)=v'$ and $\partial_1(f'_m)=u'$. In analogy to the previous considerations, a path $(f_1,\ldots,f_{k_0-k+m})$ in $\gs_k$ is constructed such that $\partial_0(f_1)=v$ and $\partial_1(f_{k_0-k+m})=u$. Altogether, the de Bruijn graph $\gs_k$ is strongly connected for $k<k_0$.
\end{proof}

\medskip

There exist subshifts such that the de Bruijn graphs of order $k_0\in\NM$ is strongly connected while the de Bruijn of order $k_0+1$ is not strongly connected, c.f. the following example.

\begin{example}
\label{Chap5-Ex-deBruijnEventuallyNotStrongConn}
Let $\as:=\{a,b\}$ be the alphabet and consider $\xi := \ldots aaab|abbb\ldots\in\as^\ZM\,,$ i.e.,
$$
\xi(j) \;
	:=\; \begin{cases}
		a\,,\quad & j=0 \text{ or } j\leq -2\,,\\
		b\,,\quad & j=-1 \text{ or } j\geq 1\,,
	\end{cases}
	\qquad j\in\ZM\,.
$$
Define the subshift $\Xi:=\overline{\Orb(\xi)}$ which only contains the orbit $\Orb(\xi)$ and the strongly periodic word $a^\infty\in\as^\ZM$ and $b^\infty\in\as^\ZM$. Then the de Bruijn graph of order $1$ is strongly connected while the de Bruijn graphs of order $k\geq 2$ are connected but not strongly connected by Proposition~\ref{Chap5-Prop-GAPBasProper} and Proposition~\ref{Chap5-Prop-HereditaryBruijn}, c.f. Figure~\ref{Chap5-Fig-BruijnNotStrongConn}~b). This also implies that $\Xi$ is not periodically approximable by Theorem~\ref{Chap5-Theo-ExPerAppr} below.
\end{example}

Recall the notions of strongly periodic and completely aperiodic dynamical system, c.f. Definition~\ref{Chap2-Def-AperSubs}. The growth of the length of closed paths in the de Bruijn graphs is an indicator of the aperiodicity of the system as showen in the following result. 

\begin{proposition}[\cite{BeBeNi16}]
\label{Chap5-Prop-ClosedPathBehaviour}
Let $\Xi\in\SG\big(\as^G\big)$ be a subshift with dictionary $\ws(\Xi)$ and associated sequence of de Bruijn graphs $(\gs_k)_{k\in\NM}$.
\begin{itemize}
\item[(a)] If $\Xi$ is completely aperiodic, then, for each $m\in\NM$, there is a $k_0:=k_0(m)\in\NM$ such that the de Bruijn graphs $\gs_k$ of order $k\geq k_0$ have no closed path of length smaller than or equal to $m$. 
\item[(b)] If $\Xi$ is completely aperiodic, then the de Bruijn graphs $(\gs_k)_{k\in\NM}$ are eventually simple, i.e., there exists a $k_0\in\NM$ such that $\gs_k$ is simple for $k\geq k_0$.
\item[(c)] If $\Xi$ is not completely aperiodic, i.e., there is an $\eta\in\Xi$ with period $p\in\NM$, then there exists a closed path of length $p$ in the de Bruijn graph $\gs_k$ of order $k$ for each $k\in\NM$.
\item[(d)] The subshift $\Xi=\Orb(\eta)$ is strongly periodic with $\eta\in\as^\ZM$ satisfying $\alpha_p(\eta)=\eta$ for $p\in\NM$ if and only if the de Bruijn graph $\gs_p$ of order $p$ associated with $\ws(\Xi)$ is connected with at most $p$ vertices such that $\gs_p$ has no branching vertices. In this case, $\gs_k$ is connected and has no branching and dandling vertices for $k\geq p$.
\end{itemize}
\end{proposition}

\begin{proof}
First, observe the following. Let $\Xi\in\SZ\big(\as^\ZM\big)$ be a subshift, $p\in\NM$ and $\gs_k$ be the de Bruijn graph of order $k> p$. Suppose that $\gs_k$ has a closed path $\wp=(e_1,\ldots,e_p)$ with $e_1=a_0a_1\ldots a_{k}$. Due to the constraints given by the boundary maps $\partial_0(e_1)=\partial_1(e_p)$ and $\partial_1(e_i)=\partial_0(e_{i+1})\,, 1\leq i\leq p-1\,,$ the equalities $a_i = a_{i+p}\,,\; 0\leq i\leq k-p\,,$ are derived. These constraints imply that $e_1\in\es_k=\ws(\Xi)\cap\as^{k+1}$ is represented by a word $v_ku$ for $v_k:=(a_0a_1\ldots a_{p-1})^{j_k}$ where $j_k\in\NM$ is the biggest natural number such that $j_k\cdot p\leq k$ and $u$ a word of length smaller than $p$. Note that $u$ is actually a prefix of $v_k$ since $a_i = a_{i+p}\,,\; 0\leq i\leq k-p$. Due to \nameref{(D1)}, all the subwords of $e_1=a_0\ldots a_{k}$ are contained in $\ws(\Xi)$. Thus, $v_k$ is an element of the dictionary $\ws(\Xi)$. Hence, $\ws(\Xi)$ contains a word that is the $j_k$-th concatenation of a word of length $p$ with $j_k:=\max\{ j\in\NM\;|\; j\cdot p\leq k \}$.

\vspace{.1cm}

(a): Let $\Xi$ be a completely aperiodic subshift and $p\in\NM$. Consider the set $I\subseteq\NM$ of natural numbers $k\in\NM$ such that $\gs_k$ has a closed path of length $p$. Assume that assertion (a) does not hold, i.e., the set $I$ is infinite. By the previous considerations, there is for each such $k\in I$ a word $v_k\in\ws(\Xi)$ satisfying that $v_k$ is the $j_k$-th time concatenation of a word of length $p$. The length $|v_k|$ is equal to $j_k\cdot p$. Consequently, there is a sequence of words $v_k\in\ws(\Xi)\,,\; k\in I\,,$ such that $|v_k|=j_k\cdot p$, $\lim_{k\to\infty}j_k=\infty$ and $v_k$ is the $j_k$-th time concatenation of a word of length $p$. 

\vspace{.1cm}

Since $\as^p$ is finite, there is a subsequence of words $v_{k_n}\in\ws(\Xi)$ that are the $j_{k_n}$-th time concatenation of a fixed word $u\in\as^p$ and $\lim_{n\to\infty}|v_{k_n}|=\infty$. Define 
\begin{align*}
l_n \; 
	&:= \; \left\lfloor \frac{j_{k_n}}{2} \right\rfloor \;
	:= \; \max\left\{ 
		i\in\NM\;\Big|\; i\leq \frac{j_{k_n}}{2} 
	\right\}\,,\\
r_n\; 
	&:= \; \left\lceil \frac{j_{k_n}}{2} \right\rceil \;
	:= \; \min\left\{ 
		i\in\NM\;\Big|\; i\geq \frac{j_{k_n}}{2} 
	\right\}\,,
\end{align*}
for $n\in\NM$. Then $l_n+r_n=j_{k_n}$ and $\lim_{n\to\infty}l_n=\lim_{n\to\infty}r_n=\infty$ hold. Thus, the compact sets $K_n:=\{-l_n,\ldots,r_n\}\,,\; n\in\NM\,,$ define an exhausting sequence for $\ZM$. For $n\in\NM$, consider an $\eta_n\in\Xi$ satisfying $\eta_n|_{[-l_n\cdot p,r_n\cdot p]}=v_{k_n}$ which exists by Lemma~\ref{Chap2-Lem-ExInfWord}. Since $\eta_{n+1}|_{[-l_n,r_n]}  = \eta_n|_{[-l_n,r_n]}\,,\; n\in\NM\,,$ hold, Lemma~\ref{Chap2-Lem-XiDictClosed} yields that $\eta:=\lim_{n\to\infty}\eta_n$ exists and $\eta$ is an element of $\Xi$. The limit $\eta$ is equal to the strongly periodic word $u^\infty\in\as^\ZM$ with period $p$. This contradicts that the subshift $\Xi$ is completely aperiodic.

\vspace{.1cm}

(b): According to (a), choose $k_0\in\NM$ such that all graphs $\gs_k,\; k\geq k_0\,,$ do not contain a closed path of length $1$, namely a loop. Thus, the de Bruijn graph $\gs_k$ is simple for all $k\geq k_0$ by using Proposition~\ref{Chap5-Prop-GAPBasProper}.

\vspace{.1cm}

(c): Let $\eta\in\Xi$ be strongly periodic with period $p$. Then the dictionary $\ws(\eta)$ of $\eta$ is a subset of $\ws(\Xi)$. Fix $k\in\NM$ and define $e_j:=\eta|_{[j,j+k]}\in\ws(\eta)\cap\as^{k+1}\subseteq\es_k$ for $1\leq j\leq p$. The chain $(e_1,\ldots, e_p)$ of edges defines a path in $\gs_k$. Due to the $p$-periodicity of $\eta$, it follows that $\partial_0(e_1)=\partial_1(e_p)$. Hence, $(e_1,\ldots,e_p)$ is a closed path of length $p$ in $\gs_k$.

\vspace{.1cm}

(d): Suppose $\Xi$ is strongly periodic such that $\Xi=\Orb(\eta)$ for $\eta\in\as^\ZM$ with $\alpha_p(\eta)=\eta$ for $p\in\NM$. Due to Proposition~\ref{Chap5-Prop-GAPBasProper}, the graph $\gs_k$ does not have any dandling vertex and is connected for all $k\in\NM$. Let $p':=\min\{k\in\NM\;|\; \alpha_k(\eta)=\eta\}$ which is smaller than or equal to $p$. Assertion (iii) implies that $\gs_k$ has a closed path of length $p'$ that visits all vertices for each $k\in\NM$, i.e., the graph $\gs_k$ has at most $p'$ vertices. For $k\geq p'$, $\gs_k$ has exactly $p'$ vertices since otherwise there exists a $1\leq j<p'$ satisfying $\alpha_j(\eta)=\eta$ which is not possible by definition of $p'$. If $k\geq p'$, then a word $v\in\ws(\Xi)\cap\as^k$ has a unique continuation to the left and a unique continuation to the right as a subword of $\eta$ since $\alpha_{p'}(\eta)=\eta$. Altogether, $\gs_k$ is connected without branching vertices and $\sharp\vs_k=p'\leq p$ for $k\geq p'$.

\vspace{.1cm}

Suppose that the de Bruijn graph $\gs_p$ associated with $\ws(\Xi)$ is connected without branching vertices and $\sharp\vs_p=:p'\leq p$. Then $\gs_p$ is a circle consisting of $p'$ vertices since $\gs_p$ has no dandling vertices by Proposition~\ref{Chap5-Prop-GAPBasProper}. Let $v:=v_1\ldots v_p\in\ws(\Xi)\cap\as^p$. Due to the constraints given by the boundary maps, $\vs_p$ is given by $\{u\in\as^p \;|\; \exists 1\leq j\leq p \text{ s.t. } u=v_j\ldots v_pv_1\ldots v_{j-1}\}$ since $\vs_p$ has at most $p'\leq p$ vertices. Thus, the periodic element $\eta:=v^\infty$ satisfies $\alpha_p(\eta)=\eta$, $\ws(\eta)\cap\as^p=\ws(\Xi)\cap\as^p$ and $\ws(\eta)\subseteq\ws(\Xi)$. Let $u_1\ldots u_k\in\ws(\Xi)\cap\as^k$ for $k>p$. By the previous considerations, there is an $l\in\NM$ such that $u_1\ldots u_p= \eta|_{[l,l+p-1}]$ follows. Since $\gs_p$ is a circle, this word has a unique continuation to the right in the dictionary $\ws(\Xi)$ given by $\eta|_{[l,l+p}]=u_1\ldots u_{p+1}$. Hence, the equality $\eta|_{[l,l+k-1}]=u$ is derived by iteration leading to $u\in\ws(\eta)$. Consequently, $\ws(\eta)=\ws(\Xi)$ is deduced. Then Theorem~\ref{Chap2-Theo-Shift+DictSpace} implies $\Xi=\Orb(\eta)$.
\end{proof}

\medskip

Note that there exist strongly periodic elements $\eta\in\as^\ZM$ with $\alpha_p(\eta)=\eta$ such that the de Bruijn graphs of order $k< p$ associated with $\ws(\eta)$ have branching vertices, c.f. Example~\ref{Chap5-Ex-NoBranch}. 

\begin{example}
\label{Chap5-Ex-NoBranch}
Let $\as:=\{a,b\}$ be the alphabet and consider the strongly periodic two-sided infinite word $\eta=(baab)^\infty\in\as^\ZM$. Then, $\eta$ has period $4$ and the de Bruijn graph $\gs_1$ of order $1$ has a branching vertex $v:=a\in\vs_1$ since $ab, aa\in\ws(\eta)$.
\end{example} 

The following proposition asserts that every strongly periodic element of $\SZ\big(\as^\ZM\big)$ is isolated in the Hausdorff-topology. Following the lines of the proof, this assertion also holds in $\ZM^d$. Specifically, the strongly periodic subshifts in $\SZd\big(\as^{\ZM^d}\big)$ are isolated in the Hausdorff-topology. 

\begin{proposition}[\cite{BeBeNi16}]
\label{Chap5-Prop-PerPointIsolated}
Let $\as$ be an alphabet. Then every strongly periodic subshift $\Xi\in\SZ\big(\as^\ZM\big)$ is isolated in the Hausdorff-topology of $\SZ\big(\as^\ZM\big)$.
\end{proposition}

\begin{proof}
Let $\Xi=\Orb(\xi)\in\SZ\big(\as^\ZM\big)$ be a strongly periodic subshift such that $\alpha_p(\xi)=\xi$ for a $p\in\NM$. The set $\vs:=\{\ws\in\DZ\;|\; \ws\cap\as^{p+1}=\ws(\Xi)\cap\as^{p+1}\}$ defines an open neighborhood in the local pattern topology of $\ws(\Xi)$. Due to Theorem~\ref{Chap2-Theo-Shift+DictSpace}, it suffices to show that $\ws(\Xi)$ is isolated in the local pattern topology. Hence, $\vs$ only contains the dictionary $\ws(\Xi)$.

\vspace{.1cm}

Let $\ws\in\vs$. According to Proposition~\ref{Chap5-Prop-ClosedPathBehaviour}~(d), every word of $\ws(\Xi)\cap\as^p$ has a unique continuation to the left and to the right in $\ws(\Xi)$. Thus, $\gs_p$ has no branching vertices. Then the de Bruijn graph $\gs_p^{\ws}$ is equal to $\gs_p$ since it only depends on the words up to length $p+1$ and $\ws\in\vs$. By assumption, $\gs_p$ is a circle with at most $p$ vertices. Consequently, $\gs_p^\ws$ is also a circle with at most $p$ vertices, c.f. Proposition~\ref{Chap5-Prop-ClosedPathBehaviour}~(d). Thus, the subshift $\Phi^{-1}(\ws)$ is strongly periodic with an $\eta\in\as^\ZM$ satisfying $\Phi^{-1}(\ws)=\Orb(\eta)$ and $\alpha_p(\eta)=\eta$ by Proposition~\ref{Chap5-Prop-ClosedPathBehaviour}~(d). Furthermore, $\ws=\ws(\eta)$ holds by Corollary~\ref{Chap2-Cor-DictOrbitSubshift}. Since $\ws\in\vs$, the equation $\ws(\eta)\cap\as^{p+1}=\ws(\xi)\cap\as^{p+1}$ is derived. Consequently, $\ws(\eta)=\ws(\xi)$ follows as $\alpha_p(\eta)=\eta$ and $\alpha_p(\xi)=\xi$. Thus, $\ws=\ws(\Xi)$ is deduced by Corollary~\ref{Chap2-Cor-DictOrbitSubshift}.
\end{proof}

\begin{remark}
\label{Chap5-Rem-SpaceBruijnGraphs}
The set of de Bruijn graphs $\bs(\as)$ over the alphabet $\as$ can be defined. This set is natural equipped with a topology so that $\bs(\as)$ is homeomorphic to $\DZ$ and to $\SZ\big(\as^\ZM\big)$ by using Theorem~\ref{Chap2-Theo-Shift+DictSpace}. The main steps for this assertion are described in the following. It is left to the reader to fill the details.

\vspace{.1cm}

Let $\as$ be an alphabet and $k\in\NM$. Let $\partial_0:\as^{k+1}\to\as^k\,,\; a_0\ldots a_k\mapsto a_0\ldots a_{k-1}\,,$ and $\partial_1:\as^{k+1}\to\as^k\,,\; a_0\ldots a_k\mapsto a_1\ldots a_k\,,$ be the boundary maps defined in Definition~\ref{Chap5-Def-DeBruijn}. Consider a pair of subsets $\vs_k\subseteq\as^k$ and $\es_k\subseteq\as^{k+1}$ satisfying $\partial_0(\es_k)\cup\partial_1(\es_k)=\vs_k$. Define for each such pair of subsets the graph $(\vs_k,\es_k,\partial_0,\partial_1)$ like in Definition~\ref{Chap5-Def-DeBruijn}. Then $\bs_k(\as)$ denotes the set of all these graphs for fixed $k\in\NM$ that have no dandling vertices. Since $\as$ is finite, the set $\bs_k(\as)$ is finite as well and so equipped with the discrete topology. Note that these graphs are not necessarily connected. A graph map is a tuple of maps mapping the vertex and the edge set of a graph to the vertex and edge set that respect the boundary maps, i.e., the maps preserve neighbors. For $k\geq 2$, the graph map $\pi_k:=\big(\pi_k^\vs,\pi_k^\es\big):\bs_k(\as)\to\bs_{k-1}(\as)$ is defined by
\begin{align*}
\pi_k^\vs(\vs_k) \;
	&:= \; \big\{ 
		a_0\ldots a_{k-1}\;|\; 
		\exists b\in\as \text{ such that } a_0\ldots a_{k-1}b\in\vs_k
	\big\}\,,\\
\pi_k^\es(\es_k) \;
	&:= \; \big\{ 
		a_0\ldots a_{k}\;|\; 
		\exists b\in\as \text{ such that } a_0\ldots a_{k}b\in\es_k
	\big\}\,,
\end{align*}
where $\gs_k:=(\vs_k,\es_k,\partial_0,\partial_1)\in\bs_k(\as)$. Note that there may exists more than one vertex of $\vs_k$ mapped to the same vertex in $\pi_k^\vs(\vs_k)$ and similarly for the edge set. In detail, the graph map $\pi_k$ identifies vertices and edges in the image. Since $\gs_k\in\bs_k(\as)$ has no dandling vertices and $\pi_k$ is a graph map, the image $\pi_k(\gs_k)$ also has no dandling vertices. Then the set
$$
\bs(\as)\;
	:= \; \Big\{
		(\gs_k)_{k\in\NM}\in\prod_{k\in\NM}\bs_k(\as) \;\Big|\;
		\pi_k(\gs_{k+1})=\gs_{k} \text{ for all } k\in\NM
	\Big\}\,.
$$
is called the {\em set of de Bruijn graphs} since each element of $\bs(\as)$ is identified with a sequence of de Bruijn graphs associated with a dictionary which is described in the following.

\vspace{.1cm}

In addition, the space $\bs(\as)$ is naturally equipped with the product topology. Thus, a base for the topology is given by the sets
$$
\us\left((\widehat{\gs}_k)_{k=1}^m\right) \;
	:= \; \Big\{
		(\gs_k)_{k\in\NM} \;|\; \gs_j=\widehat{\gs}_j \text{ for all } 1\leq j\leq m
	\Big\}
$$
for all $m\in\NM$ and $(\widehat{\gs}_k)_{k=1}^m\in\prod_{k=1}^m\bs_k(\as)$ satisfying $\pi_j(\gs_j)=\gs_{j-1}$ for $2\leq j\leq m$. Then $\bs(\as)$ is a second countable, compact, Hausdorff space, c.f. \cite{Tyc30,Cec37}.

\vspace{.1cm}

Now, consider the map $\Psi:\DZ\to\bs(\as)$ defined by $\ws\in\DZ\mapsto(\gs_k^\ws)_{k\in\NM}\in\bs(\as)$ where $\gs_k^\ws$ is the de Bruijn graph of order $k$ associated with $\ws$, c.f. Definition~\ref{Chap5-Def-DeBruijn}. Due to \cite[Satz~8.11]{Querenburg2001}, it suffices to verify that $\Psi$ is continuous and bijective so that $\Psi$ is a homeomorphism. Since the de Bruijn graph $\gs_k^\ws$ of order $k\in\NM$ only depends on $\ws\cap\as^{k+1}$, the map $\Psi$ is injective and continuous. For the surjectivity, let $\gs:=(\gs_k)_{k\in\NM}\in\bs(\as)$. According to Remark~\ref{Chap2-Rem-DictExhSeq} (see also Section~\ref{Chap5-Sect-SymbDynSystZM}), a dictionary $\ws\in\DZ$ is uniquely defined by the sets $\ws\cap\as^m\,,\; m\in\NM$. Consider the vertex sets $\vs_k\subseteq\as^k\,,\; k\in\NM\,,$ of the de Bruijn graphs $(\gs_k)_{k\in\NM}$. The constraint that $\gs_k$ has no dandling vertices becomes important. It implies that for each $v\in\vs_k$, there exist an $e,f\in\es_k$ with $\partial_1(e)=v=\partial_0(f)$. Thus, the word $v$ has an extension to the left and to the right. Together with the constraints provided by the graph map $\pi_k$, there exists a unique dictionary $\ws(\gs)\in\DZ$ satisfying $\ws(\gs)\cap\as^k=\vs_k$. Then $\Psi\big(\ws(\gs)\big)=\gs$ follows. Hence, $\Psi:\DZ\to\bs(\as)$ is bijective and continuous and so $\Psi$ defines a homeomorphism by \cite[Satz~8.11]{Querenburg2001}.
\end{remark}

\subsection{The subword complexity function}
\label{Chap5-Ssect-SubwordComplexity}

The branching vertices are strongly related to the non-periodicity of a subshift according to Proposition~\ref{Chap5-Prop-ClosedPathBehaviour} and Remark~\ref{Chap5-Rem-ConstraintBruijnGrWords}. This section draws the connection between the growth of the number of branching vertices and the subword complexity function. It is shown in Section~\ref{Chap5-Sect-PerApprSubs} that the subword complexity function provides also a lower bound for the length of the period of the approximating periodic subshifts defined via the de Bruijn graphs, c.f. Corollary~\ref{Chap5-Cor-PerGrowth}.

\medskip

Recall that $\sharp A$ of a finite set $A$ denotes the number of elements in $A$.

\begin{definition}[Subword complexity]
\label{Chap5-Def-SubwordComplex}
The {\em subword complexity function} $\text{\gls{comp}}:\NM\to\NM$ associated with a dictionary $\ws\in\DZ$ is defined by $\comp(k):=\sharp\big(\ws(\Xi)\cap\as^k\big)\,,\; k\in\NM$.
\end{definition}

The estimates $\comp(k)\leq \comp(k+1)\leq \sharp\as\cdot\comp(k)$ are immediately derived from the definition of the subword complexity function. If $\Xi$ is a periodic subshift with dictionary $\ws$, then the subword complexity is bounded by a constant. On the other hand, if there is a non-periodic element in $\Xi$, then $\comp$ growths at least linearly, i.e., $\comp(k)\geq k+1\,,\; k\in\NM$, c.f. \cite[Corollary, page 829]{MoHe38} and \cite{LuPl92,Lothaire02} for two sided-infinite sequences. 

\medskip

The subword complexity function was introduced by {\sc Morse} and {\sc Hedlund} \cite{MoHe38} as a tool to study symbolic dynamical systems. There the map was called {\em block growth}. The connection of the subword complexity function and the periodicity of a word is due to the seminal works \cite{MoHe38,MoHe40,CoHe73}, see also \cite{LuPl87,Lothaire02} for the case of two-sided infinite words. The name subword complexity function was later introduced by {\sc Ehrenfeucht}, {\sc Lee} and {\sc Rozenberg} \cite{EhLeRo75}. This function has been studied by many authors, see e.g. \cite{EhrRoz82,Rau83,LindMarcus95,Cas97,Lag99,Her00,Jul10}.

\begin{definition}[Branching map, \cite{BeBeNi16}]
\label{Chap5-Def-BranchPoint}
Let $\ws\in\DZ$ be a dictionary with associated sequence of de Bruijn graphs $(\gs_k)_{k\in\NM}$. The {\em branching map $\text{\gls{Nub}}:\NM\to\NM$} is defined by $\Nub(k):=\sharp\big\{ v\in\vs_k \;\big|\; \deg(v)>2\big\}$ where $\vs_k$ is the vertex set of $\gs_k$ and $\deg(v)$ is the vertex degree of $v$, i.e., the number of edges attached to the vertex $v$.
\end{definition}

The connection of the growth in $k\in\NM$ of the difference $\comp(k+1)-\comp(k)$ to $\Nub(k)$ was already studied in \cite{Jul10}. The following proposition provides a slightly different point of view.

\begin{proposition}[\cite{BeBeNi16}]
\label{Chap5-Prop-BranPoiSubwComplex}
Let $|\as|\geq 2$ and $\ws\in\DZ$ be a dictionary. Then, for $k\in\NM$, the estimates
$$
\frac{\comp(k+1)-\comp(k)}{\sharp\as-1}\leq \Nub(k)\leq 2 \big(\comp(k+1)-\comp(k)\big)
$$
hold. These estimates are optimal, i.e., there exists a dictionary such that both bounds are reached.
\end{proposition}

\begin{proof}
In the following, (i) the upper bound and (ii) the lower bound are proven. Finally, it is shown that (iii) these estimates are optimal by providing an example. Let $(\gs_k)_{k\in\NM}$ be the de Bruijn graphs associated with the dictionary $\ws$ and consider a $k\in\NM$.

\vspace{.1cm}

(i): The number of vertices of the graph $\gs_k$ is exactly the number of words of length $k\in\NM$ in the dictionary $\ws$, i.e., $\gs_k$ has $\comp(k)$ vertices. Let $u=a_1a_2\ldots a_k\in\ws\cap\as^k$ be a vertex, then it is a branching vertex if the word $u$ does not have a unique extension to the right or to the left. Consider the branching vertices $u$ having more than one outgoing edge, i.e., $u$ has not a unique extension to the right. Since $\comp(k+1)$ is the number of extensions of all words of length $k$ it follows that $\comp(k+1)-\comp(k)$ is an upper bound for the number of words in $\ws\cap\as^k$ having not a unique extension to the right. Similar an upper bound for the extension to the left is provided. Summing up, $\Nub(k)$ is bounded by the sum of the number of vertices having more than one outgoing edge and one ingoing edge. Hence, $\Nub(k)\leq 2 \cdot \big(\comp(k+1)-\comp(k)\big)$ follows for $k\in\NM$.

\vspace{.1cm}

(ii): In order to provide an upper bound on $\comp(k+1)$ in terms of $\comp(k)$ and $\Nub(k)$, it suffices to count for each word of length $k$ the number of allowed continuations to the right. This is equivalent to count the number of outgoing edges in the graph $\gs_k$. Each outgoing branching vertex has at most $\sharp\as$ outgoing edges and so 
$$
\underbrace{\comp(k)-\Nub(k)}_{
		\genfrac{}{}{0pt}{}{\text{vertices with unique}}{\text{extension to the right}}
	}
	+ \sharp\as\cdot  \underbrace{\Nub(k)}_{
		\genfrac{}{}{0pt}{}{\text{vertices without}}{\text{unique extension}}
	}\;
	\geq \;  \comp(k+1)\,.
$$
This leads to the lower bound for $\Nub(k)$.

\vspace{.1cm}

(iii): That the estimates are optimal is observed for the Fibonacci subshift over the alphabet $\as:=\{a,b\}$ with two letters. Specifically, consider the primitive substitution $S:\as\to\as^+\,,\; a\mapsto ab, b\mapsto a\,,$ defining a minimal, completely aperiodic subshift $\Xi_S$, c.f. Example~\ref{Chap2-Ex-NonPeriodic-Fibonacci} and Example~\ref{Chap5-Ex-Fibonacci}. Let $\ws:=\ws(\Xi_S)$ be the dictionary associated with the Fibonacci subshift. Then the subword complexity function of $\ws$ is given by $\comp(k)=k+1$ for $k\in\NM$, see e.g. \cite[Corollary~5.4.10]{Fogg02}. Thus, the estimates $1\leq \Nub(k)\leq 2$ hold by (i) and (ii). The de Bruijn graph $\gs_2$ of order $2$ has two branching vertices while the de Bruijn graph $\gs_3$ of order $3$ has one branching vertex, c.f. Example~\ref{Chap5-Ex-Fibonacci} and Figure~\ref{Chap5-Fig-Fibonacci}. Thus, the estimates are optimal.
\end{proof}

\medskip

With the subword complexity function at hand, the configurational entropy of a subshift is defined by $\lim_{k\to\infty}\frac{\log(\comp(k))}{k}$, see e.g. \cite[Section~4.8.2]{BaakeGrimm13}. For subshifts of finite type, the configurational entropy measures also the number of periodic elements contained in a subshift of finite type, c.f. \cite[Corollary~4.3.8]{LindMarcus95}. Due to Proposition~\ref{Chap5-Prop-BranPoiSubwComplex}, the $\limsup_{k\to\infty}\frac{\comp(k+1)-\comp(k)}{\comp(k)}$ is zero if and only if $\limsup_{k\to\infty}\frac{\Nub(k)}{\comp(k)}=0$. That the limit $\lim_{k\to\infty}\frac{\Nub(k)}{\comp(k)}$ is equal to zero means that the branching vertices become eventually negligible in comparison with the number of vertices in the de Bruijn graphs since $\comp(k)=\sharp\vs_k$ for $k\in\NM$ where $\sharp\vs_k$ is the number of vertices in $\gs_k$. Due to the connection of the branching vertices with the aperiodicity by Proposition~\ref{Chap5-Prop-ClosedPathBehaviour}, it is natural to expect that the behavior of $\Nub:\NM\to\NM$ and $\comp:\NM\to\NM$ is related to the rate of convergence for periodic approximations. For instance, a relation is provided in the next section in Corollary~\ref{Chap5-Cor-PerGrowth}.

\section{Periodic approximation of 1-dimensional Schr\"odinger op\-erators}
\label{Chap5-Sect-PerApprSubs}

This section provides a characterization of periodically approximable subshifts in terms of the fact that all the associated de Bruijn graphs are strongly connected, c.f. Theorem~\ref{Chap5-Theo-ExPerAppr}. Generalized Schr\"odinger operators associated with periodically approximable subshifts can be approximated by periodic Schr\"odinger operators by using Theorem~\ref{Chap4-Theo-PeriodicApproximations}. Theorem~\ref{Chap5-Theo-ExPerAppr} asserts not only the existence of strongly periodic subshifts that converge. It supplies a constructive procedure to define strongly periodic subshifts that converge to the periodically approximable subshift via the de Bruijn graphs.

\medskip

The connection of strongly periodic two-sided infinite words and de Bruijn graphs is based on the fact that closed paths in de Bruijn graphs define strongly periodic elements.

\begin{definition}[Associated periodic word]
\label{Chap5-Def-AssPerWor}
Let $\ws\in\DZ$ be a dictionary over the alphabet $\as$. Consider a closed path $\wp:=(e_1,\ldots, e_{l})$ in the de Bruijn graph $\gs_k$ of order $k\in\NM$. An edge $e_j\in\es_k$ is a word in $\ws$ of length $k+1$, i.e., $e_j=e_j(0)\ldots e_j(k)$. The {\em associated periodic word $\eta_\wp\in\as^\ZM$ with the closed path $\wp$} is defined by $\eta_\wp:=\big(e_1(0)e_2(0)\ldots e_l(0)\big)^\infty$. 
\end{definition}

Clearly, $\alpha_l(\eta_\wp)=\eta_\wp$ holds for $\eta_\wp\in\as^\ZM$ induced by a closed path $\wp:=(e_1,\ldots, e_{l})$. The advantage of the periodic word associated with a closed path $\wp$ in the de Bruijn graph is that no small defects are created at the boundary like in Example~\ref{Chap4-Ex-StrongConvSpectrConvGenSchrOp} and Example~\ref{Chap5-Ex-PerApprCutNotConv}. Specifically, all words (patterns) up to length $k+1$ are contained in the original dictionary $\ws$ if $\wp$ is a closed path in the de Bruijn graphs $\gs_k$ of order $k\in\NM$ associated with $\ws$.

\begin{lemma}[\cite{BeBeNi16}]
\label{Chap5-Lem-VerSetPerWor}
Let $\ws$ be a dictionary and $\wp=(e_1,\ldots, e_{l})$ be a closed path in the de Bruijn graph $\gs_k$ of order $k\in\NM$. Then the associated periodic word $\eta_\wp:=\big(e_1(0)e_2(0)\ldots e_l(0)\big)^\infty$ satisfies
\begin{align*}
\ws(\eta_\wp)\cap\as^k \;
	= \; \big\{\partial_0(e_j) \;\big|\; 1\leq j\leq l\big\} \;
	&= \; \big\{\partial_1(e_j) \;\big|\; 1\leq j\leq l\big\} \;
	&\subseteq \; \vs_k\,,\\
\ws(\eta_\wp)\cap\as^{k+1} \;
	&= \; \{e_j \;|\; 1\leq j\leq l \}\;
	&\!\!\subseteq \; \es_k\,.
\end{align*}
\end{lemma}

\begin{proof}
Let $\wp=(e_1,\ldots, e_{l})$ be a closed path in the de Bruijn graph associated with $\ws$ and $\eta=\eta_\wp$ be the associated periodic word. Since $\wp$ is a closed path, the equalities 
\begin{equation}
\label{Chap5-Eq-ClosPath} \tag{$\clubsuit$}
\partial_1(e_j) \;
	= \; \partial_0(e_{j+1})\,,
	\quad 1\leq j\leq l-1\,,
	\qquad
\partial_1(e_l) \; 
	= \; \partial_0(e_1)\,,
\end{equation}
are satisfied. Hence, the equation $\big\{\partial_0(e_j) \;\big|\; 1\leq j\leq l\big\} = \big\{\partial_1(e_j) \;\big|\; 1\leq j\leq l\big\}$ is derived. For $1\leq j\leq l$, all subwords of $e_j$ of length $k$ are given by $\partial_0(e_j)$ and $\partial_1(e_j)$. Thus, it suffices to prove $\ws(\eta)\cap\as^{k+1}=\{e_j \;|\; 1\leq j\leq l \}$ by the previous considerations.

\vspace{.1cm}

Let $1\leq j\leq l$. Any edge $e_j$ is a word $e_j(0)\ldots e_j(k)$ of length $k+1$. The equalities (\ref{Chap5-Eq-ClosPath}) imply $e_j|_{[1,k]}=e_{j+1}|_{[0,k-1]}$. Consequently, $\eta|_{[j-1,j+k-1]}=e_j$ follows where the $-1$ appears due to the convention $\eta(j-1)=e_j(0)$ for $1\leq j\leq l$ in the definition of $\eta:=\big(e_1(0)\ldots e_l(0)\big)^\infty$. Hence, $e_j\in\ws(\eta)\cap\as^{k+1}$ is deduced for each $1\leq j\leq l$. 

\vspace{.1cm}

For the converse inclusion, let $v\in\ws(\eta)\cap\as^{k+1}$. Since $\eta$ is strongly periodic with $\alpha_l(\eta)=\eta$, there exists a $1\leq j\leq l$ such that $v = \eta|_{[j-1,j+k-1]}$. By the previous considerations, the equality $\eta|_{[j-1,j+k-1]}=e_j$ holds. Thus, $v=e_j$ is deduced. Altogether, the desired equality $\{e_j \;|\; 1\leq j\leq l\} = \ws(\eta)\cap\as^{k+1}$ is concluded. 
\end{proof}

\medskip

With Lemma~\ref{Chap5-Lem-VerSetPerWor} at hand, the following assertion provides the connection between the de Bruijn graphs and periodically approximable subshifts.

\begin{theorem}[\cite{BeBeNi16}]
\label{Chap5-Theo-ExPerAppr} 
Let $\as$ be an alphabet and $\Xi\in\SZ\big(\as^\ZM\big)$ be a subshift with dictionary $\ws:=\ws(\Xi)$. Then the following assertions are equivalent.
\begin{itemize}
\item[(i)] The subshift $\Xi\in\SZ\big(\as^\ZM\big)$ is periodically approximable.
\item[(ii)] The de Bruijn graph $\gs_k$ associated with $\ws$ of order $k$ is strongly connected for each $k\in\NM$.
\end{itemize}
In particular, if one of these assertions hold, then there exists a strongly periodic word $\eta_k\in\as^\ZM$ for $k\in\NM$ defined via a global path in the de Bruijn graph $\gs_k$ such that $\ws(\eta_k)\cap\as^k=\ws(\Xi)\cap\as^k$, namely the strongly periodic subshifts $\Xi_k:=\Orb(\eta_k)\,,\; k\in\NM\,,$ converge to $\Xi$.
\end{theorem}

\begin{proof}
The proof use the fact that the de Bruijn graph of order $k\in\NM$ only depends on the words up to length $k+1$. The in particular part is proven together with the implication (ii)$\Rightarrow$(i).

\vspace{.1cm}

(i)$\Rightarrow$(ii): Let $\Xi\in\SZ\big(\as^\ZM\big)$ be a periodically approximable subshift. Consider a sequence $\Xi_k\in\SPZ(\as^\ZM),\; k\in\NM\,,$ of strongly periodic subshifts tending to $\Xi$. For $k\in\NM$, let $\eta_k\in\Xi_k$ be such that $\Xi_k=\Orb(\eta_k)$ and $l_k\in\NM$ be so that $\alpha_{l_k}(\eta)=\eta$. Then $\ws(\eta_k)=\ws(\Xi_k)$ holds by Corollary~\ref{Chap2-Cor-DictOrbitSubshift}. According to Theorem~\ref{Chap2-Theo-Shift+DictSpace}, the convergence of the subshifts is expressed in terms of the convergence of the associated dictionaries, i.e., for every $m\in\NM$, there is a $k_0\in\NM$ satisfying $\ws(\eta_k)\cap\as^m=\ws(\Xi)\cap\as^m$ for $k\geq k_0$, c.f. Section~\ref{Chap5-Sect-SymbDynSystZM} and Corollary~\ref{Chap2-Cor-DictExhSeq}. Denote by $\ws_k\in\DZ$ the dictionary $\ws(\eta_k)$ associated with $\eta_k$ for $k\in\NM$. Without loss of generality, suppose that $\ws(\eta_k)\cap\as^k=\ws(\Xi)\cap\as^k$ holds for all $k\in\NM$. Let $k\geq 2$. Then the chain $\wp_{k-1}:=(e_1,\ldots,e_{l_k})$ defined by 
$$
e_j \;
	:= \; \eta_k|_{[j-1,j+k-2]}\,,
	\qquad 1\leq j\leq l_k\,,
$$
is a path in the de Bruijn graph $\gs_{k-1}^{\ws_k}$ of order $k-1$ associated with $\ws_k$. The periodicity of $\eta_k$ implies $\eta_k(j)=\eta_k(j+l_k)$ for $j\in\ZM$. Consequently, $\wp_{k-1}$ is a closed path and, furthermore, it visits all vertices of $\gs_{k-1}^{\ws_k}$, i.e., $\ws_k\cap\as^{k-1}=\{\partial_0(e_j)\;|\; 1\leq j\leq l_k\}$. 

\vspace{.1cm}

Let $\gs_{k-1}^\ws$ be the de Bruijn graph of order $k-1$ associated with $\ws$. A de Bruijn graph of order $k-1$ only depends on the words up to length $k$. Thus, the de Bruijn graphs $\gs_{k-1}^\ws$ and $\gs_{k-1}^{\ws_k}$ coincide since $\ws_k\cap\as^k=\ws\cap\as^k$. Consequently, $\wp_{k-1}$ defines a global path in $\gs_{k-1}^\ws=\gs_{k-1}^{\ws_k}$. Hence, $\gs_{k-1}^\ws$ is strongly connected by Proposition~\ref{Chap5-Prop-GlobalPath} for $k\geq 2$. Note that $\wp_{k-1}$ visits also every edge in the graph $\gs_{k-1}^\ws$.

\vspace{.1cm}

(ii)$\Rightarrow$(i): Consider the family of de Bruijn graphs $\gs_k:=\gs_k^\ws\,, k\in\NM\,,$ associated with the dictionary $\ws$. These graphs are strongly connected by assumption. Let $k\in\NM$. Proposition~\ref{Chap5-Prop-GlobalPath} assures the existence of a global path $\wp_k:=(e_1,\ldots, e_{l_k})$ in $\gs_k$. Let $\eta_k\in\as^\ZM$ be the strongly periodic word associated with the path $\wp_k$ , c.f. Definition~\ref{Chap5-Def-AssPerWor}. Since $\wp_k$ is a global path in the de Bruijn graph $\gs_k$ of order $k$, the equation $\ws(\eta_k)\cap\as^k=\ws\cap\as^k$ is deduced by Lemma~\ref{Chap5-Lem-VerSetPerWor}. The subshift $\Xi_k:=\Orb(\eta_k)$ is strongly periodic since $\eta_k$ is a strongly periodic element of $\as^\ZM$. Thus, the equations $\ws(\Xi_k)\cap\as^k = \ws(\eta_k)\cap\as^k = \ws\cap\as^k$ are derived, c.f. Proposition~\ref{Chap2-Prop-AperSubs} and Proposition~\ref{Chap2-Prop-MinDictConst}. Hence, the sequence of dic\-tion\-aries $\big(\ws(\Xi_k)\big)_{k\in\NM}$ converges to $\ws$ in the local pattern topology. Consequently, the convergence $\lim_{k\to\infty}\Xi_k=\Xi$ follows in the Hausdorff-topology of $\SZ\big(\as^\ZM\big)$ by Theorem~\ref{Chap2-Theo-Shift+DictSpace}. Since $\Xi_k$ is strongly periodic for $k\in\NM$, the subshift $\Xi$ is periodically approximable.
\end{proof}

\medskip

According to Theorem~\ref{Chap5-Theo-ExPerAppr}, a sequence of global paths $(\wp_k)_{k\in\NM}$ in $(\gs_k)_{k\in\NM}$ defines a sequence of strongly periodic subshifts tending to $\Xi$. Furthermore, each sequence of strongly periodic subshifts that converges to $\Xi$ gives rise to a sequence of global paths in the family of de Bruijn graphs associated with $\ws(\Xi)$. In general, the sequence of strongly periodic subshifts $(\Xi_k)_{k\in\NM}$ does not converge to $\Xi$ in $\SZ\big(\as^\ZM\big)$ if $(\Xi_k)_{k\in\NM}$ is defined via a sequence $(\wp_k)_{k\in\NM}$ where $\wp_k$ is only a closed path in the de Bruijn graph $\gs_k$ that is not global, c.f. the following example.

\begin{example}
\label{Chap5-Ex-SeqClosPathNotConvSubsh}
Consider the two-sided infinite word $\xi\in\as^\ZM$ defined by
$$
\xi(j) \;
	:= \; \begin{cases}
		a\,,\quad &j\neq 0\,,\\
		b\,,\quad &j=0\,,
	\end{cases}
	\qquad
	j\in\NM\,,
$$
over the alphabet $\as:=\{a,b\}$. Then the associated subshift $\Xi:=\overline{\Orb(\xi)}=\Orb(\xi)\cup\{a^\infty\}$ is called {\em one-defect} and it has a family of strongly connected de Bruijn graphs by Corollary~\ref{Chap5-Cor-SuffCondPerAppr}, c.f. Section~\ref{Chap7-Sect-OneDefect}. A periodic approximation is defined via the two-sided infinite word $\eta_k:=v^\infty$ where $v=ba\ldots a\in\as^k$ for $k\in\NM$. The sequence of de Bruijn graphs $(\gs_k)_{k\in\NM}$ associated with $\Xi$ has a sequence of closed paths $(\wp_k)_{k\in\NM}$ where each path $\wp_k$ only consists of one edge such that the associated periodic word is given by $\eta_k:=a^\infty\in\as^\ZM$ for every $k\in\NM$. Clearly, the associated sequence of strongly periodic subshifts $\Xi_k:=\Orb(\eta_k)=\{\eta_k\}$ is constant and does not converge to $\Xi$ in the Hausdorff-topology since $b\not\in\ws(\Xi_k)$ for all $k\in\NM$. Thus, this sequence of closed paths $(\wp_k)_{k\in\NM}$ does not define a sequence of strongly periodic subshifts that converges to $\Xi$.
\end{example}

In the situation of a minimal subshift $\Xi$, each sequence of closed paths $\wp_k$ in the de Bruijn graphs $\gs_k$ gives rise to a strongly periodic sequence of subshifts converging to $\Xi$.

\begin{proposition}[\cite{BeBeNi16}]
\label{Chap5-Prop-MinAllPath}
Let $\as$ be an alphabet and $\Xi\in\SZ\big(\as^\ZM\big)$ be a minimal subshift with the associated sequence of de Bruijn graphs $(\gs_k)_{k\in\NM}$. Then $\Xi$ is periodically approximable. Furthermore, each sequence $(\wp_k)_{k\in\NM}$ of closed paths in $(\gs_k)_{k\in\NM}$, i.e., $\wp_k$ is a closed path in $\gs_k$ of order $k\in\NM$, gives rise to a sequence of strongly periodic subshifts converging to $\Xi$ with respect to the Hausdorff-topology on $\SZ\big(\as^\ZM\big)$. Specifically, the subshifts $\Xi_k:=\Orb(\eta_k)\,,\; k\in\NM\,,$ converge to $\Xi$ where $\eta_k\in\as^\ZM$ is the associated periodic word with $\wp_k$.
\end{proposition}

\begin{proof}
Let $\Xi$ be a minimal subshift over the alphabet $\as$ with associated sequence of de Bruijn graphs $(\gs_k)_{k\in\NM}$. Consider a sequence of closed paths $(\wp_k)_{k\in\NM}$ in $(\gs_k)_{k\in\NM}$ with associated sequence of strongly periodic words $(\eta_k)_{k\in\NM}$ defined in Definition~\ref{Chap5-Def-AssPerWor}. Define the strongly periodic subshift $\Xi_k:=\Orb(\eta_k)$ for $k\in\NM$. It suffices to prove that there exists a $k(m)$ for each $m\in\NM$ such that $\ws(\Xi_k)\cap\as^m=\ws(\Xi)\cap\as^m$ holds for $k\geq k(m)$. Let $m\in\NM$. Due to Corollary~\ref{Chap2-Cor-DictOrbitSubshift}, the equation $\ws(\Xi_k)=\ws(\eta_k)$ holds.

\vspace{.1cm}

According to \cite[Theorem~4.1.2]{Petersen83}, a subshift $\Xi\in\SZ\big(\as^\ZM\big)$ is minimal if and only if each word occurs with bounded gaps, i.e, there exists a $k(m)\in\NM$ for all $m\in\NM$ such that every word $u\in\ws(\Xi)$ of length $|u|\geq k(m)$ contains each word in $\ws(\Xi)\cap\as^m$. Clearly, $k(m)$ has to be greater than or equal to $m$. Let $k(m)$ be the corresponding integer associated with $m\in\NM$ and the minimal subshift $\Xi$. Then $\ws(\eta_k)\cap\as^m=\ws(\Xi)\cap\as^m$ is derived as follows for $k\geq k(m)$.

\vspace{.1cm}

Lemma~\ref{Chap5-Lem-VerSetPerWor} implies $\ws(\eta_k)\cap\as^k\subseteq \ws(\Xi)\cap\as^k$ for each $k\in\NM$. Thus, the inclusion $\ws(\eta_k)\cap\as^m\subseteq \ws(\Xi)\cap\as^m$ is deduced for $k\geq m$ by \nameref{(D2)}. Hence, this inclusion holds for $k\geq k(m)$ as $k(m)\geq m$.

\vspace{.1cm}

Let $v\in\ws(\Xi)\cap\as^m$. Then each word in $\ws(\Xi)$ of length larger than or equal to $k(m)$ contains $v$ by the previous considerations. Thus, for an edge $e$ of the closed path $\wp_k$ for $k\geq k(m)$, either $\partial_1(e)\in\vs_k$ or $\partial_0(e)\in\vs_k$ contains the word $v$. Lemma~\ref{Chap5-Lem-VerSetPerWor} yields that $\partial_0(e) $ and $\partial_1(e)$ are elements of $\ws(\eta_k)\cap\as^k$. Hence, $v\in\ws(\eta_k)\cap\as^m$ follows as $\ws(\eta_k)$ is a dictionary and so it satisfies \nameref{(D1)}. Since $v\in\ws(\Xi)\cap\as^m$ was arbitrarily chosen, the converse inclusion $\ws(\eta_k)\cap\as^m\supseteq \ws(\Xi)\cap\as^m$ follows for each $k\geq k(m)$.

\vspace{.1cm}

By the previous considerations, the sequence $(\Xi_k)_{k\in\NM}$ of strongly periodic subshifts converges to the minimal subshift $\Xi$. Thus, $\Xi$ is periodically approximable.
\end{proof}

\begin{remark}
\label{Chap5-Rem-PerGrowth}
Let $m\in\NM$. Then Proposition~\ref{Chap5-Prop-MinAllPath} asserts that, for $k\in\NM$ large enough, i.e., $k\geq k(m)$, every closed path $\wp_k$ defines a global path $\wp_m$ in the de Bruijn graph $\gs_m$ associated with $\ws(\Xi)$. Note that this path $\wp_m$ might visit some vertices and edges in $\gs_m$ more than once, in general.
\end{remark}

Let $\ws\in\DZ$ be a dictionary. Then the map $\comp:\NM\to\NM\,,\; \comp(k):=\sharp\ws(\Xi)\cap\as^k\,,$ is the subword complexity function associated with $\ws$, c.f. Subsection~\ref{Chap5-Ssect-SubwordComplexity}. It turns out that the growth of the subword complexity is a lower bound for the length of the global paths in the de Bruijn graphs.

\begin{corollary}[\cite{BeBeNi16}]
\label{Chap5-Cor-PerGrowth}
Let $\Xi\in\SZ\big(\as^\ZM\big)$ be a periodically approximable subshift with diction\-ary $\ws:=\ws(\Xi)$ and $\Xi_k:=\Orb(\eta_k)\in\SP(\as^G)\,,\; k\in\NM\,,$ be a sequence of strongly periodic subshifts tending to $\Xi$. If $\Xi$ is aperiodic, then the period $l_k$ of $\eta_k$ goes to infinity for $k\to\infty$. In particular, if $\eta_k$ arises by a global path in $\gs_k$, then the period $l_k$ is greater than or equal to $\comp(k):=\sharp\ws(\Xi)\cap\as^k$.
\end{corollary}

\begin{proof}
Without loss of generality, suppose that $\ws(\Xi_k)\cap\as^k=\ws(\Xi)\cap\as^k$ holds for $k\in\NM$. Let $k\in\NM$. Hence, every word of length $k$ in $\ws(\Xi)$ occurs as a subword of $\eta_k$ implying that its period $l_k$ is bounded from below by the number $\sharp\ws(\Xi)\cap\as^k=\compX(k)$ of words in $\ws(\Xi)$ of length $k$. Since $\ws(\xi)\subseteq\ws(\Xi)$ holds for each $\xi\in\Xi$, the inequality $\sharp\ws(\Xi)\cap\as^k\geq\sharp\ws(\xi)\cap\as^k$ is derived. According to \cite[Corollary, page 829]{MoHe38}, the limit $\lim_{k\to\infty}\sharp\ws(\xi)\cap\as^k$ is infinite if and only if $\xi\in\as^\ZM$ is not periodic, see also \cite{MoHe40,LuPl92,Lothaire02}. Altogether, $\Xi$ contains a non-periodic element $\xi\in\Xi$ by assumption. Consequently,
$$
\lim_{k\to\infty} l_k \; 
	\geq \; \lim_{k\to\infty} \sharp\big(\ws(\Xi_k)\cap\as^k\big) \;
	\geq \; \lim_{k\to\infty} \sharp\big(\ws(\xi)\cap\as^k\big) \;
	= \; \infty\
$$
is derived by the previous considerations.
\end{proof}

\begin{remark}
\label{Chap5-Rem-PeriodLength}
Let $\ws$ be a dictionary and $k\in\NM$. Then there exists a strongly periodic $\eta_k$ of period $l_k=\comp(k+1)$ if and only if the de Bruijn graph associated with $\ws$ of order $k$ is an Eulerian graph. Recall that a graph is called {\em Eulerian} if it contains a closed path visiting each edge exactly once. Such a closed path is called {\em Eulerian cycle}. The existence of such Eulerian cycles in the de Bruijn graphs is studied in \cite{More05}. There the Eulerian cycles are related to the so called de Bruijn sequences, c.f. \cite[Lemma~4]{More05}.
\end{remark}

\begin{proposition}[\cite{BeBeNi16}]
\label{Chap5-Prop-SuffCondPerAppr}
Let $\as$ be an alphabet and $\Xi\in\SZ\big(\as^\ZM\big)$ be a topologically transitive subshift, i.e.,  $\Xi=\overline{\Orb(\xi)}$ for a $\xi\in\as^\ZM$. Then $\Xi$ is periodically approximable if, for each $k\in\NM$, there exist $u_{le}$ and $u_{ri}$ in $\ws(\xi)\cap\as^k$ such that the following assertions hold.
\begin{description}
\item[(a)] The word $u_{le}$ occurs infinitely often in $\xi$ to the left, i.e., for all $l_0\in\NM$, there exists an $l\geq l_0$ satisfying $\xi|_{[-l,-l+k-1]}=u_{le}$.
\item[(b)] The word $u_{ri}$ occurs infinitely often in $\xi$ to the right, i.e.,  for all $r_0\in\NM$, there exists an $r\geq r_0$ satisfying $\xi|_{[r,r+k-1]}=u_{ri}$.
\item[(c)] There exist an $l_1,r_1\in\ZM$ such that $r_1<l_1$ and 
$$
\xi_{[l_1,l_1+k-1]} \; 
	= \; u_{le}\,,
	\qquad
\xi_{[r_1,r_1+k-1]} \;
	= \; u_{ri}\,.
$$
\end{description}
\end{proposition}

\begin{proof}
According to Theorem~\ref{Chap5-Theo-ExPerAppr}, it suffices to verify that the de Bruijn graph $\gs_k$ of order $k$ associated with $\ws(\Xi)$ is strongly connected for every $k\in\NM$. Since $\Xi=\overline{\Orb(\xi)}$, the equation $\ws(\Xi)=\ws(\xi)$ follows by Corollary~\ref{Chap2-Cor-DictOrbitSubshift}. Thus, the de Bruijn graphs $(\gs_k)_{k\in\NM}$ associated with $\ws(\Xi)$ are defined via the dictionary $\ws(\xi)$. Let $k\in\NM$ and $u_{le}\,,\; u_{ri}\in\ws(\xi)\cap\as^k$ be such that they satisfy (a),(b) and (c). The proof is organized as follows.
\begin{itemize}
\item[(i)] There exists a path $(e_1,\ldots, e_n)$ in $\gs_k$ such that $\partial_0(e_1)=u_{le}$ and $\partial_1(e_n)=u_{ri}$.
\item[(ii)] There exists a path $(f_1,\ldots, f_{n'})$ in $\gs_k$ such that $\partial_0(f_1)=u_{ri}$ and $\partial_1(f_{n'})=u_{le}$.
\item[(iii)] Let $v\in\ws(\xi)\cap\as^k$. Then there is a path
\begin{Itemize}[0.1]
\item[(iii.a)] $(\tilde{e}_1,\ldots,\tilde{e}_{l'})$ in $\gs_k$ such that $\partial_0(\tilde{e}_1)=u_{le}$ and $\partial_1(\tilde{e}_{l'})=v$.
\item[(iii.b)] $(\tilde{f}_1,\ldots,\tilde{f}_{r'})$ in $\gs_k$ such that $\partial_0(\tilde{f}_1)=v$ and $\partial_1(\tilde{f}_{r'})=u_{ri}$.
\end{Itemize}
\item[(iv)] The graph $\gs_k$ is strongly connected.
\end{itemize}

(i): Due to condition (a) and (b), there exist an $l,r\in\NM$ such that 
$$
\xi_{[-l,-l+k-1]} \;
	= \; u_{le}\,,
	\qquad
\xi_{[r,r+k-1]} \;
	= \; u_{ri}\,,
	\qquad
-l<r\,.
$$
Then $u_{le}$ is a prefix and $u_{ri}$ is a suffix of the word $w:=\xi|_{[-l,r+k-1]}\in\ws(\xi)\cap\as^{l+r+k-1}$. Thus, Lemma~\ref{Chap5-Lem-WordExPath} implies assertion (i).

\vspace{.1cm}

(ii): This follows in analogy to (i) by using (c). More precisely, $u_{ri}$ is a prefix and $u_{le}$ is a suffix of the word $w:=\xi|_{[r_1,l_1+k-1]}\in\ws(\xi)$.

\vspace{.1cm}

(iii): Let $v\in\ws(\xi)\cap\as^k$. Then there exists an $m\in\ZM$ such that $\xi_{[m,m+k-1]}=v$.
\begin{Itemize}[0.2]
\item[(iii.a):] According to (a) there exists an $l\geq |m|$ satisfying $\xi_{[-l,-l+k-1]}=u_{le}$. Thus, $u_{le}$ is a prefix and $v$ is a suffix of the word $w:=\xi|_{[-l,m+k-1]}\in\ws(\xi)$. Then Lemma~\ref{Chap5-Lem-WordExPath} leads to the desired result.
\item[(iii.b):] According to (b) there exists an $r\geq |m|$ satisfying $\xi_{[r,r+k-1]}=u_{ri}$. Thus, $v$ is a prefix and $u_{ri}$ is a suffix of the word $w:=\xi|_{[m,r+k-1]}\in\ws(\xi)$. Then Lemma~\ref{Chap5-Lem-WordExPath} leads to the desired result.
\end{Itemize}
(iv): Let $v_1\,,\, v_2\in\ws(\xi)\cap\as^k$. Consider two paths $(e_1,\ldots,e_m)$ and $(f_1,\ldots,f_{m'})$ in $\gs_k$ such that $\partial_1(e_m)=\partial_0(f_1)$. Then the chain $(e_1,\ldots, e_m,f_1,\ldots,f_{m'})$ defines also a path in the graph $\gs_k$. With this at hand, assertions (i)-(iii) imply that there are two paths $\wp_1$ and $\wp_2$ in $\gs_k$ satisfying that $\wp_1$ joins $v_1$ to $v_2$ and $\wp_2$ joins $v_2$ to $v_1$. Consequently, the de Bruijn graph $\gs_k$ associated with $\ws(\Xi)$ is strongly connected.
\end{proof}

\begin{remark}
\label{Chap5-Rem-SuffCondPerAppr}
Due to the finiteness of $\as^k$ for $k\in\NM$, it is clear that every two-sided infinite word $\xi\in\as^\ZM$ has subwords $u_{le},u_{ri}\in\ws(\xi)\cap\as^k$ satisfying (a) and (b). The essential requirement is the combination with (c). This guarantees that there exists a path joining $u_{ri}$ to $u_{le}$. Then assertion (ii) is proven by using this fact. Note that a path joining $u_{le}$ to $u_{ri}$ exists always thanks to (a) and (b).
\end{remark}

\begin{corollary}[\cite{BeBeNi16}]
\label{Chap5-Cor-SuffCondPerAppr}
Let $\as$ be an alphabet and $\Xi\in\SZ\big(\as^\ZM\big)$ be a subshift such that $\Xi=\overline{\Orb(\xi)}$ for a $\xi\in\as^\ZM$. Then $\Xi$ is periodically approximable if, for each $k\in\NM$, there exists an $u\in\ws(\xi)\cap\as^k$ such that the word $u$ occurs infinitely often in $\xi$ to the left and the right, i.e., for all $n_0\in\NM$, there exist an $l,r\geq n_0$ such that $\xi_{[-l,-l+k-1]} = u = \xi_{[r,r+k-1]}$.
\end{corollary}

\begin{proof}
This follows from Proposition~\ref{Chap5-Prop-SuffCondPerAppr}.
\end{proof}

\cleardoublepage


\chapter{Periodically approximable subshifts in \texorpdfstring{$\ZM^d$}{ZMd} arising by primitive block substitutions}
\label{Chap6-HigherDimPerAppr}
\stepcounter{section}
\setcounter{section}{0}

This chapter deals with subshifts of the symbolic dynamical system $(\as^{\ZM^d},\ZM^d,\alpha)$ defined via a primitive (block) substitution. Sufficient conditions are provided so that these subshifts are periodically approximable. In this case, an approximation of strongly periodic subshifts is defined by applying the substitution iteratively to a suitable strongly periodic element of $\as^{\ZM^d}$, c.f. Theorem~\ref{Chap6-Theo-PerApprZdPrimSubst}. As it turns out, local symmetries in the patterns of a subshift imply that the subshift is periodically approximable.

\medskip

Substitutions are a standard tool to construct subshifts in $\ZM^d$ or tilings in $\RM^d$. Most of the known and studied examples arise from primitive substitutions. Primitive substitutions are of particular interest since the associated subshifts is minimal, c.f. Proposition~\ref{Chap6-Prop-PrimSubstSubshZdMin}. The construction of a subshift associated with a substitution is standard for the group $\ZM$, see e.g. \cite{Fogg02,BaakeGrimm13}. It is straight forward to extend these notions to $\ZM^d$, see e.g. \cite{LuPl87,Sol98,Rob99,Fra05,Fra08,BaakeGrimm13}. Due to the lack of a general reference, the construction of the subshift associated with a substitution over the group $\ZM^d$ and its basic properties are presented in Section~\ref{Chap6-Sect-PrimBlockSubst}. For the further applications, it is crucial that the substitution $S$ extends to a continuous map $S:\as^{\ZM^d}\to\as^{\ZM^d}$, c.f. Lemma~\ref{Chap6-Lem-SubstitutionContinuous}. Thus, we restrict our considerations to block substitutions in this chapter whenever $d\geq 2$.

\medskip

If $d\neq 1$, sufficient conditions are given to insure that the subshift associated with a primi\-tive substitution is periodically approximable, c.f. Theorem~\ref{Chap6-Theo-PerApprZdPrimSubst} and Proposition~\ref{Chap6-Prop-SuffPerApprZ2SymmPatt}. In Section~\ref{Chap7-Sect-TableSubst} and Section~\ref{Chap7-Sect-SierpinskiCarpetSubstitution}, examples are provided satisfying this condition for $d=2$. Only the methods developed in Section~\ref{Chap2-Sect-SymbDynSyst}, Section~\ref{Chap4-Sect-ContSpecGenSchrOp} and Section~\ref{Chap5-Sect-PerApprSubs} are used to prove these convergence results. 

\medskip

Recall that if $d=1$, the class of all minimal subshift is periodically approximable, c.f. Proposition~\ref{Chap5-Prop-MinAllPath}. Actually, all subshifts satisfying conditions (a)-(c) of Proposition~\ref{Chap5-Prop-SuffCondPerAppr} are periodically approximable. A similar statement in the case $d\neq 1$ remains open, c.f. Section~\ref{Chap8-Sect-PerApprox}. The general philosophy for the connections of $d=1$ and $d\neq 1$ is the following.

\medskip

In the one-dimensional case closed paths in the associated de Bruijn graphs give rise to strongly periodic approximations, c.f. Theorem~\ref{Chap5-Theo-ExPerAppr}. If the subshift arises by a primitive substitution, then a closed path in the de Bruijn graphs defines a strongly periodic approximation by applying the substitution to the strongly periodic element associated with the closed path, c.f. Corollary~\ref{Chap6-Cor-PerApprZdPrimSubst-d=1}.

\medskip

In the higher dimensional case, the de Bruijn graphs are replaced by the so called Anderson-Putnam complex \cite{AnPu98,BeGa03,BeBeGa06}, c.f. Remark~\ref{Chap5-Rem-GAP}. The analog of a closed path is a torus in this complex. Let $\Xi_S$ be a subshift defined via a substitution $S$. Then a torus in the Anderson-Putnam complex defines a strongly periodic element $\eta\in\as^{\ZM^d}$ such that the patterns of $\eta$ up to a certain size $K\subseteq\ZM^d$ are contained in $\ws_S$. In this case, $\Xi_n:=\Orb\big(S^n(\eta)\big)\in\SG\big(\as^{\ZM^d}\big)\,,\; n\in\NM\,,$ define a sequence of strongly periodic subshifts converging to $\Xi_S\in\SG\big(\as^{\ZM^d}\big)$, c.f. Theorem~\ref{Chap6-Theo-PerApprZdPrimSubst}. The main task is to prove the existence of such a strongly periodic element, c.f. discussion in Section~\ref{Chap8-Sect-PerApprox}.

\medskip

Following the strategy \nameref{(5.I)}-\nameref{(5.IV)}, the guideline of this chapter is as follows. Consider an exhausting sequence $(K_m)_{m\in\NM}$ of $\ZM^d$. Let $\Xi_S$ be the minimal subshift with dictionary $\ws_S$ defined via a primitive block substitution $S:\as\to\PaZd$, c.f. Proposition~\ref{Chap6-Prop-PrimSubstSubshZdMin}.

\begin{description}
\item[(6.I)\label{(6.I)}] Find an $\eta\in\as^{\ZM^d}$ that is strongly periodic such that $\ws(\eta)\cap\as^{[K]}\subseteq\ws(\xi)\cap\as^{[K]}$ where $\prod_{j=1}^d\{0,1\}\subseteq K\subseteq\ZM^d$, c.f. Corollary~\ref{Chap6-Cor-ExKPeriodPoint-d=1} and Proposition~\ref{Chap6-Prop-SuffPerApprZ2SymmPatt}.
\item[(6.II)\label{(6.II)}] Prove that $S^n(\eta)\in\as^{\ZM^d}$ is strongly periodic for each $n\in\NM$, c.f. Lemma~\ref{Chap6-Lem-SnETAStronglyPeriodic}.
\item[(6.III)\label{(6.III)}] Show that, for all $m\in\NM$, there exists an $n_0\in\NM$ satisfying the inclusion $\ws\big(S^n(\eta)\big)\cap\as^{[K_m]}\subseteq\ws_S\cap\as^{[K_m]}$ for $n\geq n_0$, c.f. Lemma~\ref{Chap6-Lem-AssDictInvSubstit} and Lemma~\ref{Chap6-Lem-GrowthAcceptPatternSubst}.
\item[(6.IV)\label{(6.IV)}] For every $m\in\NM$, prove the existence of an integer $n_1\in\NM$ such that the inclusion $\ws_S\cap\as^{[K_m]}\subseteq\ws\big(S^n(\eta)\big)\cap\as^{[K_m]}$ holds for each $n\geq n_1$, c.f. Theorem~\ref{Chap6-Theo-PerApprZdPrimSubst}.
\end{description}

\begin{remark}
\label{Chap6-Rem-StrategySubstitution}
(i) The compact set $\tilde{K}:=\prod_{j=1}^d\{0,1\}$ in \nameref{(6.I)} is a $d$-dimensional discrete cube with two points in each direction. This set plays the role of a word of length $2$ if $d=1$. Even in the case $d=1$, the words of length two are necessary to define a closed path in the de Bruijn graph $\gs_1$ of order $1$ associated with a subshift. The element $\eta\in\as^{\ZM^d}$ in \nameref{(6.I)} plays a similar role like a periodic word of associated with a closed path in the de Bruijn graph for $d=1$. More precisely, $\eta$ and $K$ define a torus in the Anderson-Putnam complex associated with $\ws_S$, c.f. Remark~\ref{Chap5-Rem-GAP} and Section~\ref{Chap8-Sect-APComplex}.

\vspace{.1cm}
 
(ii) This chapter does not provide a sufficient condition in large generality for the existence of an $\eta\in\as^{\ZM^d}$ with desired properties in \nameref{(6.I)} if $d\neq 1$. Specifically, only in the case $d=2$ a class of subshifts satisfying \nameref{(6.I)} is provided, c.f. Proposition~\ref{Chap6-Prop-SuffPerApprZ2SymmPatt}. However, this result shows the essential ingredient for a subshift being periodically approximable. Specifically, Proposition~\ref{Chap6-Prop-SuffPerApprZ2SymmPatt} asserts that the existence of periodic approximations is related to a local symmetry of the patterns. Furthermore, it provides an intuition which type of properties need to be checked in higher dimensions. In the one-dimensional case, i.e., $d=1$, such an $\eta$ exists for all primitive substitution $S:\as\to\as^+$ where $\as^+\subseteq\PaZ$ is the set of finite words excluding the empty word, c.f. Corollary~\ref{Chap6-Cor-PrimSubsPerAppr} and Corollary~\ref{Chap6-Cor-PerApprZdPrimSubst-d=1}.  

\vspace{.1cm}

(iii) The primitivity is mainly used in step \nameref{(6.IV)}.

\vspace{.1cm}

(iv) If $d=1$, the notion of a block substitution is also called a substitution with constant length. This requirement is only necessary if $d\neq 1$. More precisely, all stated assertions hold for all primitive substitutions $S:\as\to\as^+$, c.f. Remark~\ref{Chap6-Rem-BlockNotNec-d=1}.
\end{remark}

The reader is referred to \cite[Chapter~5]{Queffelec87}, \cite[Section~1.2]{Fogg02} and \cite[Chapter~4]{BaakeGrimm13} for more background on substitutional dynamical systems over the group $\ZM$. Background for substitutions in $\ZM^d$ can be found in \cite{LuPl87,Sol98,Rob99,Fra05,Fra08,BaakeGrimm13}.

\section{Primitive block substitutions}
\label{Chap6-Sect-PrimBlockSubst}

In the following, the basic concepts of primitive block substitutions are introduced. It is shown that a primitive block substitution defines a dictionary $\ws_S$ and a subshift $\Xi_S$ that is minimal, c.f. Proposition~\ref{Chap6-Prop-PrimSubstSubshZdMin}. If $d=1$, this implies that $\Xi_S$ is periodically approximable, c.f. Corollary~\ref{Chap6-Cor-PerApprZdPrimSubst-d=1}.

\medskip

In the one-dimensional case, words of certain length were sufficient to describe a dictionary, c.f. Section~\ref{Chap5-Sect-SymbDynSystZM}. For the group $\ZM^d$ the analog of words are patterns defined on so called blocks. 

\begin{definition}[Block]
A compact subset $K\subseteq\ZM^d$ is called a {\em block} whenever there are $n_1,\ldots,n_d\in\NM$ such that $K=\prod_{j=1}^d \{1,\ldots, n_j\}$. Then the integers $n_1,\ldots,n_d\in\NM$ are called {\em parameters} of the block $K$.
\end{definition}

\begin{remark}
\label{Chap6-Rem-Block-Dictionary}
According to Remark~\ref{Chap2-Rem-DictExhSeq}, a dictionary in $\DZd$ is uniquely determined by the patterns with support on the sets $\{K_n\;|\; n\in\NM\}\subsetneq\ks\big(\ZM^d\big)$ where $(K_n)_{n\in\NM}$ is an exhausting sequence. Clearly, such an exhausting sequence can be defined such that $K_n$ is equivalent by $\ZM^d$-translation to a block for each $n\in\NM$. For instance, the blocks $K_m:=\prod_{j=1}^d \{-2^m,\ldots,2^m\}\subseteq\ZM^d$ for $m\in\NM$ define an exhausting sequence of $\ZM^d$ and $(2^m,\ldots,2^m)+K_m$ is a block for each $m\in\NM$.
\end{remark}

Let $K\subseteq\ZM^d$. Then $\as^K$ denotes the set of all maps $v:K\to\as$. By taking the equivalence class with respect to translation by elements of $\ZM^d$, the set of patterns $\as^{[K]}$ is defined in Definition~\ref{Chap2-Def-Pattern}. Here $[K]$ denotes the equivalence class of a compact set $K\subseteq\ZM^d$ with respect to translation by elements of $\ZM^d$. Each element of $\as^{[K]}$ is determined by one of its representative $v:K\to\as$. For sake of simplification, a pattern $[v]$ is identified in the following with a representative $v\in\as^K$ whenever there is no confusion. Then the compact set $K\subseteq\ZM^d$ is called the support of $v$. Patterns with support on a block give rise to a strongly periodic element in $\as^{\ZM^d}$.

\begin{definition}[Strongly periodic extension]
\label{Chap6-Def-StronglyPeriodExt}
Let $\as$ be an alphabet. Consider a pattern $v\in\as^K$ where $K\subseteq\ZM^d$ is a block with parameters $n_1,\ldots,n_d\in\NM$. Then the  {\em strongly periodic extension $\xi:=v^\infty\in\as^{\ZM^d}$} of $v$ is defined by $\xi|_{p_i+K}=v$ for all $i:=(i_1,\ldots, i_d)\in\ZM^d$ where $p_i:=(i_1 \cdot n_1,\ldots, i_d\cdot n_d)\in\ZM^d$.
\end{definition}

Clearly, a strongly periodic extension of a pattern $v$ is strongly periodic. Conversely, each strongly periodic $\xi\in\ZM^d$ is equal to $v^\infty$ for a suitable pattern $v\in\as^K$ and a block $K\subseteq\ZM^d$. 

\begin{lemma}
\label{Chap6-Lem-ReprStrongPer}
Let $\xi\in\as^{\ZM^d}$. Then $\xi$ is strongly periodic if and only if there exist a block $K\subseteq\ZM^d$ and a pattern $v\in\as^K$ such that $\xi=v^\infty$.
\end{lemma}

\begin{proof}
Recall that $\xi\in\as^{\ZM^d}$ is strongly periodic if the orbit $\Orb(\xi)$ is finite, c.f. Definition~\ref{Chap2-Def-AperSubs}. Consider a $\xi\in\as^{\ZM^d}$ defined by $v^\infty$ for a pattern $v\in\as^K$ supported on a block $K\subseteq\ZM^d$ with parameters $n_1,\ldots,n_d\in\NM$. Then the number of elements of $\Orb(\xi)$ is bounded from above by $\prod_{j=1}^d n_j$. Hence, $\xi$ is strongly periodic.

\vspace{.1cm}

Let $\xi\in\as^{\ZM^d}$ be strongly periodic. Then, for $1\leq j\leq d$, there exists an $n_j\in\NM$ such that $\alpha_{n_j\cdot e_j}(\xi)=\xi$ holds where $e_j\in\ZM^d$ is the vector that is $1$ at the $j$-th entry and zero everywhere else. Consequently, $\xi(i+n_j\cdot e_j)=\xi(i)$ follows for each $i\in\ZM^d$. Define $v:=\xi|_{K}$ for the block $K:=\prod_{j=1}^d\{1,\ldots, n_j\}\subseteq\ZM^d$. Then the equality $\xi=v^\infty$ is immediately derived by the previous considerations.
\end{proof}

\medskip

Like in Chapter~\ref{Chap5-OneDimCase}, it is assumed throughout the chapter that the alphabet contains at least two letters. If $\as$ only has one element there exists only one element in $\as^{\ZM^d}$ which is not interesting at all.

\begin{definition}[Primitive block substitution]
\label{Chap6-Def-BlockPrimSubst}
Let $\as$ be an alphabet with at least two letters. A map $S:\as\to\PaZd$ is called a {\em substitution} over the alphabet. A substitution $S$ is called 
\begin{itemize}
\item {\em primitive} if there exists an $l_0\in\NM$ such that every letter $b\in\as$ occurs in $S^{l_0}(a)$ for all $a,b\in\as$. Then the smallest such $l_0$ is called the {\em parameter of the primitivity of $S$};
\item a {\em block substitution} if there exists a block $K := \prod_{j=1}^d \big\{1,2\ldots,n_j\big\} \subseteq \ZM^d$ with parameters $n_1,\ldots, n_d\geq 2$ such that $S(a)\in\as^{[K]}$ for all letters $a\in\as$.
\end{itemize}
\end{definition}

\begin{remark}
\label{Chap6-Rem-BlockNotNec-d=1}
There are two reasons for us to assume that $S$ is a block substitution. 
\begin{itemize}
\item[(i)] A block substitution $S$ naturally extends to a map on $S:\PaZd\to\PaZd$ and to a map $S:\as^{\ZM^d}\to\as^{\ZM^d}$ by acting separately on each letter, c.f. Lemma~\ref{Chap6-Lem-SubstitutionContinuous}. If $S$ is not a block substitution this might fail. Specifically, a letterwise definition can create gaps in the coloring of $\ZM^d$ if the tiles do not fit together, c.f. \cite[Page~297]{Fra08}.
\item[(ii)] In the definition of a block substitution, the parameters $n_1,\ldots, n_d$ are assumed to be greater than or equal to $2$. This implies that by applying the substitution to a letter $a\in\as$ the size of $S(a)$ increases in any direction $j\in\{1,\ldots,d\}$. In the one-dimensional case ($d=1$) this means that the length of the word $S^n(a)$ tends to infinity if $n\to\infty$. 
\end{itemize}

For $d=1$, each substitution $S:\as\to\as^+$ extends to a continuous map $S:\as^\ZM\to\as^\ZM$ without assuming that the substitution has a block structure. Here $\as^+\subseteq\PaZ$ denotes the free semigroup on $\as$ introduced in Section~\ref{Chap5-Sect-SymbDynSystZM}, i.e., $\as^+:=\bigcup_{k\in\NM}\as^k$. Furthermore, if, additionally, $\as$ has at least two letters and $S:\as\to\as^+$ is a primitive substitution, then $|S^{l_0\cdot n}(a)|\geq 2^n$ follows for each letter $a\in\as$ and $n\in\NM$ where $|S^{l_0\cdot n}(a)|$ denotes the length of the word. This is derived as $S^{l_0}(a)$ contains each letter of $\as$. Consequently, $|S^{l_0}(a)|$ is greater than or equal to $2$ for each $a\in\as$. Thus, $\lim_{n\to\infty}|S^n(a)|=\infty$ is deduced implying that all primitive substitutions $S:\as\to\as^+$ satisfy (i) and (ii). Hence, every of the following assertion is valid for a primitive substitutions $S:\as\to\as^+$ that is not a block substitution if $d=1$.
\end{remark}

In the one-dimensional case, there are several examples that arise from primitive substitutions that are not block substitutions, c.f. the Fibonacci and the Silver mean substitution in Section~\ref{Chap7-Sect-SturmSeq}. Note that block substitutions for $d=1$ are also called substitutions of constant length. Examples for block substitutions over the group $\ZM$ are the Prouhet-Thue-Morse sequence (Section~\ref{Chap7-Sect-ThueMorse}), Period-Doubling sequence (Section~\ref{Chap7-Sect-PerDoubl}) and the Rudin-Shapiro sequence (Section~\ref{Chap7-Sect-GolRudSha}).

\medskip

Another example for a block substitution is the table substitution (Section~\ref{Chap7-Sect-TableSubst}) defined over the group $\ZM^2$ as follows, c.f. \cite{Fra08} and \cite[Section~4.9]{BaakeGrimm13}. Note that the table substitution was first introduced by {\sc Solomyak} \cite{Sol98}.

\begin{example}[Table substitution]
\label{Chap6-Ex-TableSubst}
Consider the alphabet $\as:=\{a,b,c,d\}$. The map $S:\as\to\textit{Pat}_{\ZM^2}(\as)$ defined by
$$
a\; \overset{S}{\longmapsto} \;
	\begin{matrix}
		b & a\\
		d & a
	\end{matrix}\,,
	\qquad
b\; \overset{S}{\longmapsto} \;
	\begin{matrix}
		a & c\\
		b & b
	\end{matrix}\,,
	\qquad
c\; \overset{S}{\longmapsto} \;
	\begin{matrix}
		c & b\\
		c & d
	\end{matrix}\,,
	\qquad
d\; \overset{S}{\longmapsto} \;
	\begin{matrix}
		d & d\\
		a & c
	\end{matrix}\,,
$$
is called the {\em table substitution}. The table substitution is a block substitution with parameters $n_1=n_2=2$. Furthermore, $S$ is primitive with $l_0=2$, i.e., the letter $y$ occurs in $S^2(x)$ for every pair of letters $x,y\in\as$. The associated subshift turns out to be aperiodic and periodically approximable, c.f. Proposition~\ref{Chap6-Prop-BasPropTableSubst} and Section~\ref{Chap7-Sect-TableSubst}.
\end{example}

A primitive block substitution defines a dictionary $\ws_S\in\DZd$ and a subshift $\Xi_S\in\SZd\big(\as^{\ZM^d}\big)$ as follows:

\begin{definition}[Dictionary/subshift associated with a primitive block substitution]
\label{Chap6-Def-AssDictSubshSubstZd}
Let $\as$ be an alphabet and $S:\as\to\PaZd$ be a primitive block substitution. The {\em dictionary $\ws_S\subseteq\PaZd$ associated with $S$} is defined by all $v\in\PaZd$ such that there is a $k_v\in\NM$ and a letter $a_v\in\as$ satisfying $v$ is a subpattern of $S^k(a)$. Furthermore, 
$$
\Xi_S \;
	:= \; \big\{
		\xi\in\as^{\ZM^d}\;|\; \ws(\xi)\subseteq\ws_S
	\big\} \;
	\subseteq\; \as^{\ZM^d}
$$
is called the {\em subshift associated with $S$}.
\end{definition}

As the names suggest, the set $\ws_S\subseteq\PaZd$ and $\Xi_S\subseteq\as^{\ZM^d}$ associated with a primitive block substitution $S:\as\to\PaZd$ define a dictionary and a subshift.

\begin{proposition}
\label{Chap6-Prop-AssDictSubshSubstZd}
Let $\as$ be an alphabet and $S:\as\to\PaZd$ be a primitive block substitution. Then the dictionary $\ws_S$ associated with $S$ defines a dictionary in terms of Definition~\ref{Chap2-Def-Dictionary}. Furthermore, $\Xi_S$ is a subshift of $\as^{\ZM^d}$, i.e., it is closed and $\ZM^d$-invariant.

\vspace{.1cm}

If $d=1$ and $S:\as\to\as^+$ is only a primitive substitution, then the set $\ws_S\subseteq\PaZ$ is a dictionary and $\Xi_S\subseteq\as^{\ZM}$ is a subshift.
\end{proposition}

\begin{figure}[htb]
\centering
\includegraphics[scale=0.93]{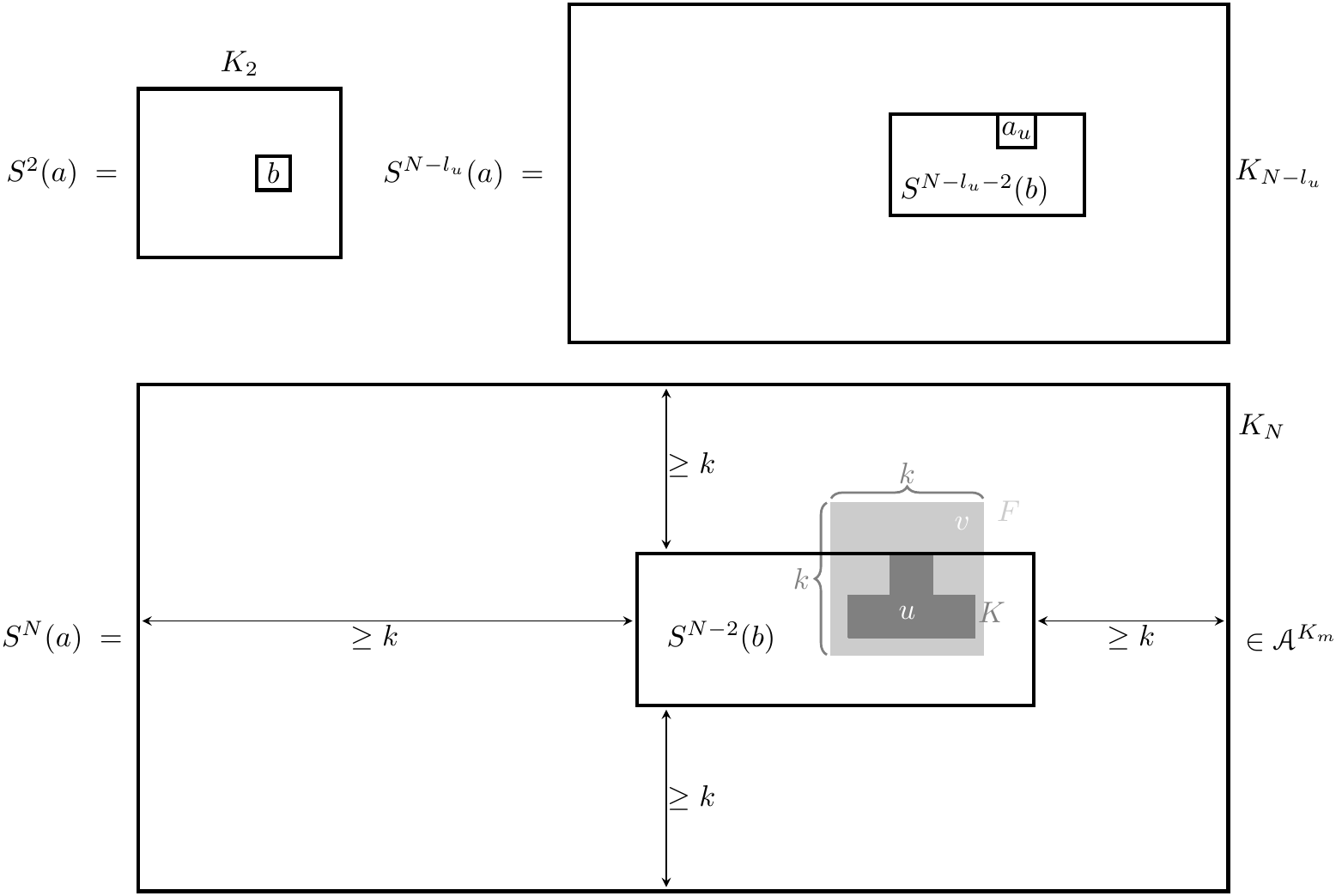}
\caption{Each rectangle represents a pattern in $\ZM^d$ over the alphabet $\as$. The number $N$ is defined by the sum $k+l_u+l_0+2$. The $2$ implies that a letter $b$ can be chosen in the interior of $S^2(a)$ since $S$ is a block substitution. Let $u\in\as^{K}$ be such that $[u]\in\ws_S$, i.e., there is an $l_u\in\NM$ and a letter $a_u\in\as$ such that $u$ occurs in $S^{l_u}(a_u)$. Furthermore, $F\in\ks(\ZM^d)$ is a block with the same edge lengths $k\in\NM$ such that $K\subseteq F$. The substitution $S$ acts letterwise. Since $N-l_u-2\geq l_0$, the pattern $S^{N-l_u-2}(b)$ contains the letter $a_u$ implying that the pattern $u$ appears in $S^{N-2}(b)$, i.e., there is a $g_u\in\ZM^d$ satisfying $[S^{N-2}(b)|_{g_u+K}]=[u]$. The blocks around the letter $b$ grow for each application of $S$ by $n_j\geq 2$ in direction $e_j\in\ZM^d$ for $1\leq j\leq d$. Thus, the distance from the boundary of $S^{N-2}(b)$ is greater than or equal to $k$ by applying at least $k$-times the substitution to $S^2(a)$. Hence, $g_u+F\subseteq K_m$ follows since $N-2\geq k$ and $b$ was chosen in the interior of $S^2(a)$.}
\label{Chap6-Fig-WSDict}
\end{figure}

\begin{proof} 
Only the proof for a primitive block substitution is presented. The case of a primitive substitution $S:\as\to\as^+$ is similarly treated by taking Remark~\ref{Chap6-Rem-BlockNotNec-d=1} into account.

\vspace{.1cm}

In the following, it is proven that (i) $\ws_S\subseteq\PaZd$ is a dictionary and (ii) $\Xi_S\subseteq\as^{\ZM^d}$ is a subshift.

\vspace{.1cm}

(i): A subset of $\PaZd$ is called a dictionary if it satisfies \nameref{(D1)} and \nameref{(D2)}, c.f. Definition~\ref{Chap2-Def-Dictionary}. Condition~\nameref{(D1)} immediately follows for $\ws_S$ by definition. Condition~\nameref{(D2)} is derived by using that $S$ is a primitive block substitution see the following discussion. For convenience of the reader, the main idea of the proof is sketched in Figure~\ref{Chap6-Fig-WSDict}.

\vspace{.1cm}

The pattern $S^m(a)$ is represented by a map $S^m(a):K_m\to\as$ for each $a\in\as$ and $m\in\NM$ with $K_m:=\prod_{j=1}^d\{1,\ldots n_j^m\}\subseteq\ZM^d$ which is possible due to the block structure of $S$. Fix $a\in\as$. Since $n_1,\ldots,n_d\geq 2$ there exists an $i:=(i_1,\ldots,i_d)\in K_2$ such that $1<i_j<n_j^2$ for each $1\leq j\leq d$. Consider the letter $b:=S^2(a)(i)$ at the $i$-th entry of the pattern $S^2(a)$. Due to construction the set $\prod_{j=1}^d\{i_j-1,i_j,i_j+1\}$ is contained in $K_m$. 

\vspace{.1cm}

Let $[u]\in\ws_S$ with representative $u:K\to\as$. Then there is a letter $a_u\in\as$ and a $l_u\in\NM$ satisfying $u$ occurs in $S^{l_u}(a_u)$. Consider a compact $F\subseteq\ZM^d$ such that $K\subseteq F$. It suffices to show the existence of a $v\in\as^{g_u+F}$ such that $[v]\in\ws_S$ and $[v|_{g_u+K}]=[u]$ for a $g_u\in\ZM^d$. Thanks to compactness of $F$, there exists a $k\in\NM$ such that $F\subseteq\prod_{j=1}^d\{1,\ldots,k\}$. Thus, there is no loss of generality in assuming by \nameref{(D1)} that $F$ is a block with same edge lengths, i.e., $F=\prod_{j=1}^d\{1,\ldots,k\}$ for a $k\in\NM$. Denote by $l_0\in\NM$ the parameter of the primitivity, i.e., $S^{l_0}(a)$ contains the letter $b$ for every pair $a,b\in\as$. Define $N:=k+l_u+l_0+2$ and the pattern $w:=S^{N}(a):K_N\to\as$. Then $[w]$ is an element of $\ws_S$ by definition of $\ws_S$. 

\vspace{.1cm}

With this at hand, it is proven in the following that there exists a $g_u\in\ZM^d$ such that $g_u+F\subseteq K_m$ and $[w|_{g_u+K}]=[u]$. This implies that $v:=w|_{g_u+F}$ satisfies $[v|_{g_u+K}]=[u]$ and $[v]\in\ws_S$ since $[w]\in\ws_S$ and $\ws_S$ satisfies \nameref{(D1)}.

\vspace{.1cm}

Recall the notation $b:=S^2(a)(i)$. Then the letter $a_u$ occurs in $S^{N-l_u-2}(b)$ since $N-l_u-2 \geq l_0$. Consequently, $u$ is a subpattern of $S^{N-2}(b)$. Thus, there exists a $g_u\in\ZM^d$ such that $S^N(a)|_{g_u+K}=u$. On the other hand, $i\in K_2$ was chosen such that $1< i_j< n_j^2$ for $1\leq j\leq d$. Since $N-2>k$ and $n_1,\ldots,n_d\geq 2$, the inclusion $g_u+F\subseteq K_m$ is derived. Altogether, \nameref{(D2)} follows for $\ws_S$ by the previous considerations. Hence, $\ws_S$ is a dictionary.

\vspace{.1cm}

(ii) Due to Theorem~\ref{Chap2-Theo-Shift+DictSpace}, there is a homeomorphism $\Phi:\DZd\to\SZd\big(\as^{\ZM^d}\big)$. The subset $\Xi_S\subseteq\as^{\ZM^d}$ is by definition the inverse with respect to $\Phi$ of the dictionary $\ws_S$. Consequently, $\Xi_S$ is a closed, $\ZM^d$-invariant subset of $\as^{\ZM^d}$ by (i) and Theorem~\ref{Chap2-Theo-Shift+DictSpace}.
\end{proof}

\medskip

Note that there are also non-primitive substitutions such that $\ws_S$ defines a dictionary, c.f. Section~\ref{Chap7-Sect-SierpinskiCarpetSubstitution}. In this case, $\Xi_S$ defines also a subshift. On the other hand, the following example shows that the set $\ws_S$ is not necessarily a dictionary if $S$ is not primitive.

\begin{example}
\label{Chap6-Ex-wsSNotDict}
Consider the alphabet $\as:=\{a,b\}$ and the substitution $S:\as\to\as^+$ defined by $a\mapsto ba\,,\; b\mapsto b$. The set $\ws_S$ does not contain a word that begins with the letter $a$. Thus, there does not exist an extension of the word $S(a)=ba\in\ws_S$ to the right. Hence, $\ws_S$ does not satisfy \nameref{(D2)}.
\end{example}

\begin{figure}[htb]
\centering
\includegraphics[scale=0.613]{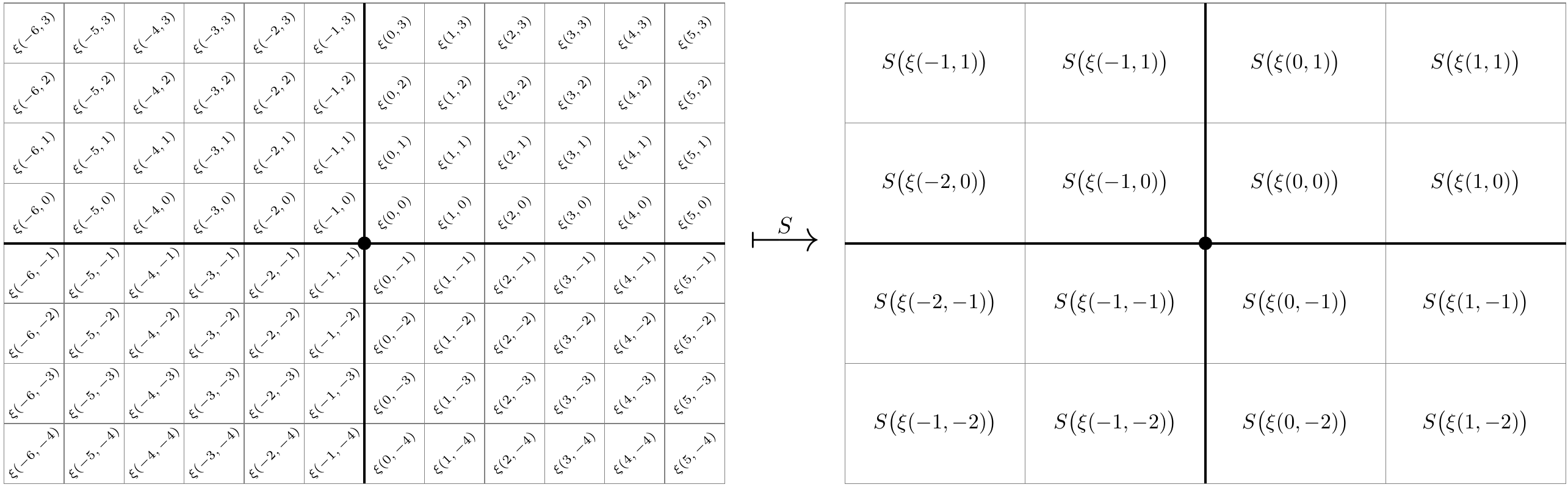}
\caption{Action of the block substitution on $\as^{\ZM^d}$ where the dot represents the origin.}
\label{Chap6-Fig-ActSubst}
\end{figure}

A block substitution naturally extends to a continuous map $S:\as^{\ZM^d}\to\as^{\ZM^d}$, c.f Figure~\ref{Chap6-Fig-ActSubst}.

\begin{lemma}
\label{Chap6-Lem-SubstitutionContinuous}
Let $S:\as\to\PaZd$ be a block substitution over the alphabet $\as$. Then the substitution $S$ extends to a continuous map $S:\as^{\ZM^d}\to\as^{\ZM^d}$ by acting on each letter separately. Furthermore, the substitution $S$ extends to a map $S:\PaZd\to\PaZd$. If $d=1$, the assertion holds for every substitution $S:\as\to\as^+$.
\end{lemma}

\begin{proof}
It suffices to show that $S$ extends to a continuous map $S:\as^{\ZM^d}\to\as^{\ZM^d}$. That $S$ also extends to a map $S:\PaZd\to\PaZd$ is similarly treated. The case $S:\as\to\as^+$ for $d=1$ is left to the reader.

\vspace{.1cm}

For every letter $a\in\as$, $S(a)$ is represented by an element of $\as^K$ with $K:=\prod_{j=1}^d\{1,\ldots,n_j\}$ where $n_1,\ldots,n_d\geq 2$ are the parameters of the block substitution $S$. Let $\xi\in\as^{\ZM^d}$. Then $\eta:=S(\xi)$ is defined pointwise by 
$$
\eta
	\begin{pmatrix}
		i_1\cdot k_1\\
		\vdots\\
		i_d\cdot k_d 
	\end{pmatrix} \; 
	:= \; \Big(S\big(\xi(i)\big)\Big)(k)\,,
	\qquad
	k:=\begin{pmatrix}
		k_1\\
		\vdots\\
		k_d
	\end{pmatrix}\in K
	\,,\;
	i:=\begin{pmatrix}
		i_1\\
		\vdots\\
		i_d
	\end{pmatrix} \in\ZM^d
	\,.
$$
In detail, let $i\cdot k$ be the pointwise multiplication of the entries of the vectors $k\in K\subseteq\ZM^d$ and $i\in\ZM^d$. Then the letter at $i\cdot k$ of $\eta$ is defined by the letter in the $k$-th entry of the pattern $S\big(\xi(i)\big)\in\as^K$. Since $S(a)\,,\; a\in\as\,,$ are defined on the same support $K\subseteq\ZM^d$, the map $S:\as^{\ZM^d}\to\as^{\ZM^d}$ is well-defined. Thus, $S(\eta)$ defines an element in $\as^{\ZM^d}$. Then the continuity follows by the facts that $S$ is defined on each letter separately, the support of $S(a)$ is a compact set $K\subseteq\ZM^d$ and the topology on $\as^{\ZM^d}$ is defined by the clopen sets $\os(F,[v]):=\{\xi\in\as^{\ZM^d}\;|\; \xi|_F\in[v]\}$ for a compact set $F\subseteq\ZM^d$ and $[v]\in\as^{[F]}$.
\end{proof}

\medskip

The following assertion is an immediate consequence of the definition of $\ws_S$.

\begin{lemma}
\label{Chap6-Lem-AssDictInvSubstit}
Let $\as$ be an alphabet and $S:\as\to\PaZd$ be a block substitution. Then $S(v)$ is an element of $\ws_S$ for each $v\in\ws_S$. If $d=1$, then the same assertion holds for each substitution $S:\as\to\as^+$.
\end{lemma}

\begin{proof}
Let $v\in\ws_S$. Thus, there exist a $k_v\in\NM$ and a letter $a_v\in\as$ such that $v$ occurs in $S^{k_v}(a_v)$. Then $S(v)$ is a subpattern of $S^{k_v+1}(a_v)$. Hence, $S(v)$ is an element of $\ws_S$.
\end{proof}

\begin{definition}[$k$-periodic element with respect to $S$]
\label{Chap6-Def-KPeriodPoint}
Let $S:\as\to\PaZd$ be a block substitution over the alphabet $\as$. An element $\xi\in\as^{\ZM^d}$ is called {\em  $k$-periodic with respect to the substitution $S$} if $S^k(\xi)=\xi$ for $k\in\NM$ and $\ws(\xi)\subseteq\ws_S$. 
\end{definition}

The existence of a $k$-periodic point for a $k\in\NM$ is proven in the following Proposition~\ref{Chap6-Prop-ExKPeriodPoint} for primitive block substitutions. This assertion is well-known, see e.g. \cite[Proposition~V.1]{Queffelec87}, \cite[Page~9]{Fogg02} for $d=1$ and \cite{LuPl87}, \cite[Lemma~2.1]{Sol98} for $d\neq 1$.

\begin{proposition}[\cite{Queffelec87,Fogg02}]
\label{Chap6-Cor-ExKPeriodPoint-d=1}
Let $S:\as\to\as^+$ be a primitive substitution on the alphabet $\as$ with at least two letters. Then there exist a $\xi\in\as^\ZM$ and a $k\in\NM$ such that $\xi$ is $k$-periodic with respect to the substitution $S$.
\end{proposition}

\begin{figure}[htb]
\centering
\includegraphics[scale=0.9]{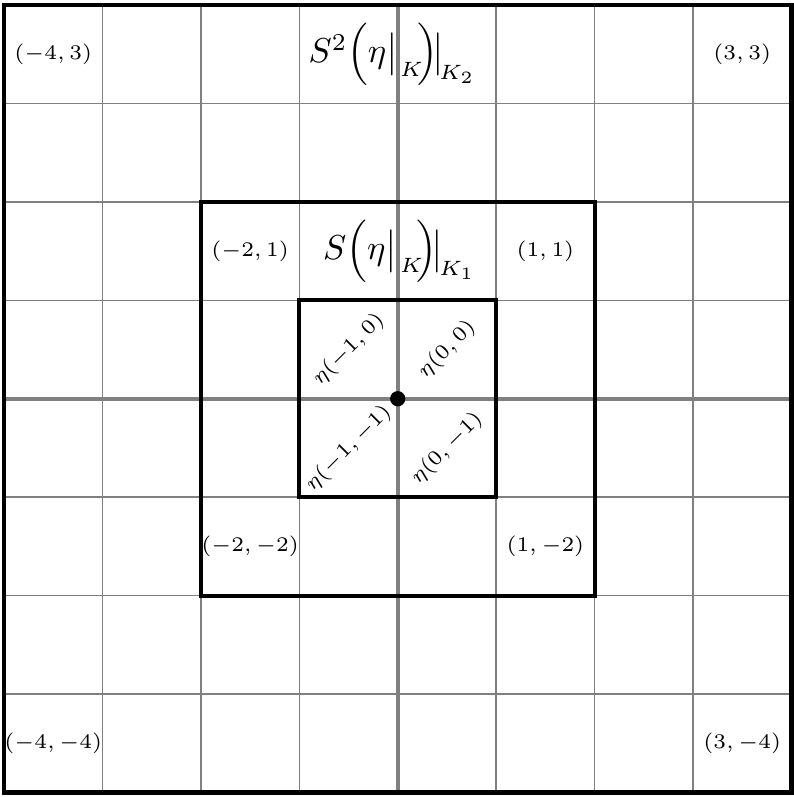}
\caption{The growth of $S^{k\cdot n}(\eta|_K)$ in $n\in\NM$ where $k=1$.}
\label{Chap6-Fig-ExKPeriodPoint}
\end{figure}

\begin{proposition}[\cite{LuPl87}, Existence of $k$-periodic elements with respect to $S$]
\label{Chap6-Prop-ExKPeriodPoint}
Let $S:\as\to\PaZd$ be a block substitution over the alphabet $\as$. Then there exist a $\xi\in\Xi_S$ and a $k\in\NM$ such that $S^k(\xi)=\xi$. 
\end{proposition}

\begin{proof}
Let $S:\as\to\PaZd$ be a block substitution with corresponding parameters $n_1,\ldots,n_d\geq 2$. Consider the compact set $K:=\prod_{j=1}^d\{-1,0\}$. Since $\as^K$ is finite, there exist a $k\in\NM$ and a pattern $v:K\to\as$ with $[v]\in\as^{[K]}\cap\ws_S$ such that, for every $\eta\in\as^{\ZM^d}$ with $[\eta|_K]=[v]$, the equality $S^k(\eta)|_K=\eta|_K$ holds. Let $\eta\in\Xi_S$ be such that $[\eta|_K]=[v]$ which exists by Lemma~\ref{Chap2-Lem-ExInfWord}. An induction over $n\in\NM_0$ leads to the equation $S^{k\cdot n}(\eta)|_{K_n}=S^{k\cdot (n+1)}(\eta)|_{K_n}$ and the fact that $S^{k\cdot n}(\eta)|_{K_n}$ is an element of $\ws_S$ where $K_n:=\prod_{j=1}^d\big\{-2^{k\cdot n},\ldots,2^{k\cdot n}-1\big\}$, c.f. Figure~\ref{Chap6-Fig-ExKPeriodPoint}.

\vspace{.1cm}

More precisely, the base case $n=0$ follows as $K_0=K$. Suppose now that $S^{k\cdot n}(\eta)|_{K_n}=S^{k\cdot (n+1)}(\eta)|_{K_n}\in\ws_S$ holds for $n\in\NM_0$ (induction hypothesis). A letter is supported on a block where all the parameters are equal to one. By applying the substitution to the letter $a\in\as$, the support growths in any direction $j\in\{1,\ldots,d\}$ at least by one. Thus, the support of $S^k(a)$ growths in each direction at least by $2^k$ as $S$ act letterwise. Consequently, $S^{k\cdot(n+1)}(\eta)|_{K_{n+1}}$ is a subpattern of $S\big(S^{k\cdot n}(\eta)|_{K_{n}}\big)$. Similarly, $S^{k\cdot(n+2)}(\eta)|_{K_{n+1}}$ is a subpattern $S\big(S^{k\cdot(n+1)}(\eta)|_{K_{n}}\big)$. Hence, $S^{k\cdot(n+1)}(\eta)|_{K_{n+1}} = S^{k\cdot(n+2)}(\eta)|_{K_{n+1}}$ is derived by the induction hypothesis. Since $S(v)$ is an element of $\ws_S$ if $v\in\ws_S$ by Lemma~\ref{Chap6-Lem-AssDictInvSubstit}, it follows that $S^{k\cdot (n+1)}(\eta)|_{K_{n+1}}=S^{k\cdot (n+2)}(\eta)|_{K_{n+1}}\in\ws_S$ by using that $S^{k\cdot(n+1)}(\eta)|_{K_{n+1}}$ is a subpattern of $S\big(S^{k\cdot n}(\eta)|_{K_{n}}\big)$ and \nameref{(D1)}. This finishes the induction.

\vspace{.1cm}

The sequence $(K_n)_{n\in\NM}$ of compact sets defines an exhausting sequence of $\ZM^d$. Thus, the sequence $\eta_n:=S^{k\cdot n}(\eta)\,,\; n\in\NM\,,$ defines a convergent sequence in $\as^{\ZM^d}$ by Lemma~\ref{Chap2-Lem-XiDictClosed}. Denote the limit $\lim_{n\to\infty}\eta_n$ by $\xi$. The map $S:\as^{\ZM^d}\to\as^{\ZM^d}$ is continuous, c.f. Lemma~\ref{Chap6-Lem-SubstitutionContinuous}. Consequently, the equalities $S^k(\xi)=\lim_{n\to\infty} S^{k\cdot (n+1)}(\xi)=\xi$ are derived. 

\vspace{.1cm}

Due to Lemma~\ref{Chap6-Lem-AssDictInvSubstit}, $S^n(\eta)$ is an element of $\Xi_S$ for $n\in\NM$ as $\eta\in\Xi_S$. Thus, $\xi:=\lim_{n\to\infty}\eta_n$ is contained in $\Xi_S$ since $\Xi_S$ is closed.
\end{proof}

\medskip

Primitive substitutions are of particular interests as the associated subshifts are minimal. This result is well-known, see e.g. \cite[Theorem~V.2]{Queffelec87}, \cite[Proposition~5.5]{Queffelec10}, \cite[Proposition~1.2.4]{Fogg02} for $d=1$ and \cite[Lemma~2.2]{Sol98} for $d\neq 1$.

\begin{proposition}[\cite{Queffelec87,Sol98}]
\label{Chap6-Prop-PrimSubstSubshZdMin}
Let $S:\as\to\PaZd$ be a primitive block substitution over the alphabet $\as$ with at least two letters. Then the associated subshift $\Xi_S$ is a minimal subshift and, hence, topologically transitive. In particular, $\Xi_S=\overline{\Orb(\xi)}$ holds for every $k$-periodic point $\xi\in\Xi_S$ for $k\in\NM$. If $d=1$, the assertion holds for all primitive substitutions $S:\as\to\as^+$.
\end{proposition}

\begin{proof}
In view of Remark~\ref{Chap6-Rem-BlockNotNec-d=1}, the assertion for primitive substitutions $S:\as\to\as^+$ follows the same lines as for primitive block substitution, c.f. \cite{Queffelec87,Queffelec10,Fogg02}. Thus, only the proof for a primitive block substitution is presented here. The proof is organized as follows. 
\begin{itemize}
\item[(i)] The equality $\ws(\xi)=\ws_S$ holds for every $k$-periodic element $\xi\in\as^\ZM$ with respect to $S$. This leads to $\Xi_S=\overline{\Orb(\xi)}$.
\item[(ii)] The subshift $\overline{\Orb(\xi)}=\Xi_S$ is minimal.
\end{itemize}

(i): Due to Proposition~\ref{Chap6-Prop-ExKPeriodPoint}, there exists at least one $k$-periodic element $\xi\in\Xi_S$ with respect to the substitution $S$ for a $k\in\NM$. Let $\xi\in\Xi_S$ be such a $k$-periodic element with respect to the substitution $S$ for $k\in\NM$. The inclusion $\ws(\xi)\subseteq\ws_S$ holds by Definition~\ref{Chap6-Def-KPeriodPoint}. The converse inclusion is proved as follows.

\vspace{.1cm}

Let $[v]\in\ws_S$ be a pattern. Then there exist a letter $a_v\in\as$ and an $l_v\in\NM$ such that $v$ is a subpattern of $S^{l_v}(a_v)$ by definition of $\ws_S$. Since $S$ is primitive, there is an $l_0\in\NM$ satisfying that $S^{l_0}(a)$ contains the letter $b$ for every pair of letters $a,b\in\as$. Hence, the pattern $v$ is a subpattern of $S^{k\cdot(l_0+l_v)}\big(\xi(0)\big)$. Since $\xi$ is $k$-periodic with respect to the substitution $S$, the equation $\xi=S^{k\cdot(l_0+l_v)}(\xi)$ follows. Consequently, $S^{k\cdot(l_0+l_v)}\big(\xi(0)\big)$ is a subpattern of $\xi$ leading to $[v]\in\ws(\xi)$ since $\ws(\xi)$ satisfies \nameref{(D1)}, c.f. Lemma~\ref{Chap2-Lem-xiDict}. Altogether, the equation $\ws(\xi)=\ws_S$ is derived implying $\Xi_S=\overline{\Orb(\xi)}$ by Theorem~\ref{Chap2-Theo-Shift+DictSpace}. Thus, $\Xi_S$ is topologically transitive.

\vspace{.1cm}

(ii): Let $\xi\in\as^{\ZM^d}$ be a $k$-periodic element with respect to the substitution $S$ for $k\in\NM$. According to (i), the equation $\ws(\xi)=\ws_S$ is concluded. Due to Proposition~\ref{Chap2-Prop-MinDictConst}, it suffices to show that $\ws(\xi)=\ws(\eta)$ for every $\eta\in\Xi_S$. The inclusion $\ws(\eta)\subseteq\ws(\xi)=\ws_S$ is derived by $\ws_S=\bigcup_{\eta\in\Xi_S}\ws(\eta)$. 

\vspace{.1cm}

For the converse inclusion, let $v\in\ws_S$ and $l_0\,,\, l_v\in\NM$ be such that $v$ is a subpattern of $S^{k\cdot(l_0+l_v)}\big(a\big)$ for each letter $a\in\as$. Due to the block structure of $S$, there is a $K_v\subseteq\ZM^d$ such that $S^{k\cdot(l_0+l_v)}\big(a\big)\in\as^{K_v}$ for all $a\in\as$. Hence, the support is independent of the letter $a\in\as$. Thus, $\xi|_{K_v+j}$ contains the pattern $v$ for every $j\in\ZM^d$ by using $S^{k\cdot (l_0+l_v)}(\xi)=\xi$. Consider a limit point $\eta:=\lim_{n\to\infty}\alpha_{j_n}(\xi)$ for a sequence $(j_n)_{n\in\NM}$ in $\ZM^d$, namely $\eta$ is an arbitrary element of $\Xi_S$ by (i). According to the definition of the topology on $\as^{\ZM^d}$, there is an $n_0\in\NM$ such that $\eta|_{K_v} = \alpha_{j_n}(\xi)|_{K_v}$ holds for all $n\geq n_0$. Since $\alpha_{j_n}(\xi)|_{K_v}$ is equal to $\xi|_{K_v-j_n}$ for $n\in\NM$, the previous considerations imply that $v$ occurs in $\eta|_{K_v}$. Thus, $v\in\ws(\eta)$ follows.
\end{proof}

\medskip

In the one-dimensional case, the minimality of the subshift $\Xi_S$ immediately implies that $\Xi_S$ is periodically approximable. For $d\neq 1$, a similar assertion does not exist in such a large generality, c.f. Section~\ref{Chap8-Sect-PerApprox}.

\begin{corollary}
\label{Chap6-Cor-PrimSubsPerAppr}
Let $\as$ be an alphabet with at least two letters and $S:\as\to\as^+$ be a primitive substitution. Then the associated subshift $\Xi_S$ is periodically approximable. More precisely, for each $k_0\in\NM$, there exists a strongly periodic $\eta\in\as^\ZM$ satisfying $\ws(\eta)\cap\as^{k_0}=\ws_S\cap\as^{k_0}$.
\end{corollary}

\begin{proof}
This follows by Proposition~\ref{Chap5-Prop-MinAllPath} and Proposition~\ref{Chap6-Prop-PrimSubstSubshZdMin}.
\end{proof}

\medskip

Note that there exist substitutions that are not primitive while they have a fixed point $\xi$ defining a minimal subshifts, see the following example.

\begin{example}
\label{Chap6-Ex-Chacon sequence}
Consider the alphabet $\as:=\{a,b\}$ and the substitution $S:\as\to\as^+$ defined by
$a\mapsto aaba\,,\; b\mapsto b$. Since $S^k(b)=b$ for every $k\in\NM$, this substitution is not primitive. There exists at least one fixed point as the length $|S^n(a)|$ tends to infinity. For instance, a fixed point is defined via $\eta:=\lim_{n\to\infty}S^n(\xi)$ where $\xi:=a^\infty$. This limit exists as $S(a)$ begins and ends with the letter $a$. It is proven in \cite[page~134]{Fogg02} that $\Xi_S:=\overline{\Orb(\eta)}$ is a minimal subshift. According to Proposition~\ref{Chap5-Prop-MinAllPath}, the subshift $\Xi_S$ is periodically approximable. This subshift $\Xi$ is related to the Chacon substitution, c.f. \cite{Cha69,Fer95}.
\end{example}

Recall that a subshift is called aperiodic if it contains a non-periodic element. It is shown in \cite[Theorem~1.1]{Sol98} that a tiling is aperiodic if and only if the substitution rule is injective. Specifically, the map $S:\as\to\PaZd$ is invertible on the image. The table substitution satisfies this condition, c.f. \cite[Example~2]{Sol98} and \cite{Fra05}. 

\begin{proposition}[\cite{Sol98,Fra05}]
\label{Chap6-Prop-BasPropTableSubst}
Let $\Xi_S$ be the subshift defined via the table substitution. Then $\Xi_S$ is minimal and aperiodic and so $\Xi_S$ is completely aperiodic.
\end{proposition}

\begin{proof}
Proposition~\ref{Chap6-Prop-PrimSubstSubshZdMin} leads to the minimality of $\Xi_S$ as the table substitution is a primitive block substitution with $l_0=2$ satisfying that $y$ occurs in $S^{l_0}(y)$ for every pair of letters $x,y\in\as$. The map $S:\as\to\textit{Pat}_{\ZM^2}(\as)$ is invertible on the image, i.e., it is injective. Thus, the subshift $\Xi_S$ is aperiodic by \cite[Theorem~1.1]{Sol98}. Altogether, $\Xi_S$ is completely aperiodic by Proposition~\ref{Chap2-Prop-AperSubs}~(c) since $\ZM^d$ is an abelian group.
\end{proof}

\section{Strongly periodic approximations for primitive block substitutions}
\label{Chap6-Sect-PerApprPrimBlockSubst}

This section deals with the strongly periodic approximation of subshifts arising by primi\-tive (block) substitutions. We define the sequence of strongly periodic subshifts $(\Xi_n)_{n\in\NM}$ via the orbit of $S^n(\eta)$ for $n\in\NM$ where $\eta\in\as^{\ZM^d}$ is a suitable strongly periodic element, c.f. Theorem~\ref{Chap6-Theo-PerApprZdPrimSubst}. In order to do so, we follow the strategy \nameref{(6.I)}--\nameref{(6.IV)} given in the introduction of this chapter.

\vspace{.1cm}

First, it is necessary to check whether $S^n(\eta)$ is still strongly periodic for $n\in\NM$ if $\eta\in\as^{\ZM^d}$ is strongly periodic corresponding to step \nameref{(6.II)} of the strategy.

\begin{lemma}
\label{Chap6-Lem-SnETAStronglyPeriodic}
Let $\eta\in\as^{\ZM^d}$ be a strongly periodic element and $S:\as\to\PaZd$ be a block substitution. Then $S^n(\eta)$ is strongly periodic for each $\eta\in\as^{\ZM^d}$. If $d=1$, then the same assertion holds for every substitution $S:\as\to\as^+$.
\end{lemma}

\begin{proof}
Due to Lemma~\ref{Chap6-Lem-ReprStrongPer}, there exist a block $K\subseteq\ZM^d$ and a pattern $v\in\as^K$ such that $\eta=v^\infty$. Let $n\in\NM$. Since $S$ is a block substitution, $S^n(v)$ is a pattern defined on a new block $K_n\subseteq\ZM^d$. With this at hand, the equality $S^n(\eta)=\big(S^n(v)\big)^\infty$ is derived as $S:\as^{\ZM^d}\to\as^{\ZM^d}$ acts letterwise, c.f. Lemma~\ref{Chap6-Lem-SubstitutionContinuous}. Thus, $S^n(\eta)$ is strongly periodic by Lemma~\ref{Chap6-Lem-ReprStrongPer}.

\vspace{.1cm}

The case $d=1$ is similarly treated by keeping Remark~\ref{Chap6-Rem-BlockNotNec-d=1}~(i) in mind.
\end{proof}

\medskip

According to Corollary~\ref{Chap2-Cor-DictExhSeq}, the convergence of dictionaries in $\DZd$ is characterized by the convergence along an exhausting sequence $(K_m)_{m\in\NM}$ of the group $\ZM^d$ independent of the choice of the sequence. More precisely, a sequence of dictionaries $(\ws_n)_{n\in\NM}$ converges to $\ws\in\DZd$ if and only if, for each $m\in\NM$, there exists an $n_0\in\NM$ such that 
$$
\ws_n\cap\as^{[K_m]} \; 
	= \; \ws\cap\as^{[K_m]}\,, 
	\qquad n\geq n_0\,.
$$
For $m\in\NM$, consider the block $K_m:=\prod_{j=1}^d\{-2^{m}+1,\ldots,2^{m}\}\subseteq\ZM^d$. In analogy with Section~\ref{Chap5-Sect-SymbDynSystZM}, the sequence of compact sets $(K_m)_{m\in\NM}$ defines an exhausting sequence for $\ZM^d$.  With this at hand, step \nameref{(6.III)} of the strategy is treated.

\begin{lemma}
\label{Chap6-Lem-GrowthAcceptPatternSubst}
Let $S:\as\to\PaZd$ be a block substitution over the alphabet $\as$ with at least two letters. Suppose that there exists a strongly periodic $\eta\in\as^{\ZM^d}$ satisfying $\ws(\eta)\cap\as^{(K]}\subseteq\ws_S\cap\as^{[K]}$ where $K:=\prod_{j=1}^d \{0,1\}\subseteq\ZM^d$. Then, for each $m\in\NM$, the inclusion $\ws\big(S^n(\eta)\big)\cap\as^{[K_m]}\subseteq\ws_S\cap\as^{[K_m]}$ holds for all $n\geq m$. If $d=1$ and $S:\as\to\as^+$ is a primitive substitution, then $\ws\big(S^{l_0+n}(\eta)\big)\cap\as^m\subseteq\ws_S\cap\as^m$ holds for all $n\geq m$ and $m\in\NM$ where $l_0\in\NM$ is the parameter of the primitivity of $S$.
\end{lemma}

\begin{proof}
Let $[u]\in\ws(\eta)\cap\as^{[K]}$ with representative $u\in\as^K$. Since $S:\as\to\PaZd$ is a block substitution (Definition~\ref{Chap6-Def-BlockPrimSubst}), $S(a)$ is extended in each direction $1\leq j\leq d$ at least by one. Thus, $S(u)$ is represented by an element of $\as^{F_1}$ with $K_1\subseteq F_1\subseteq\ZM^d$. With the same argument, an induction over $m\in\NM$ implies that $S^m(u)$ is represented by an element of $\as^{F_m}$ with $K_m \subseteq F_m\subseteq\ZM^d$ c.f. proof of Proposition~\ref{Chap6-Prop-ExKPeriodPoint}. Thus, for every $[w]\in\ws\big(S^m(\eta)\big)\cap\as^{[K_m]}$, there exists a $[u]\in\ws(\eta)\cap\as^{[K]}$ such that $[w]$ is a subpattern of $S^m(u)$.

\vspace{.1cm} 

Let $m\in\NM$. The previous considerations imply that each $[w]\in\ws\big(S^m(\eta)\big)\cap\as^{[K_m]}$ is a subpattern of $S^m(u)$ for a suitable $u:K\to\as$ with $[u]\in\ws(\eta)\cap\as^{[K]}\subseteq\ws_S\cap\as^{[K]}$. Recall that $\ws_S$ is defined by all patterns that are subpatterns of $S^l(a)$ for an $a\in\as$ and $l\in\NM$. Hence, there is a letter $a_u\in\as$ and an $l_u\in\NM$ such that $[u]$ is a subpattern of $S^{l_u}(a_u)$. Thus, $[w]$ is a subpattern of $S^{l_u+m}(a_u)$ by the previous considerations implying $[w]\in\ws_S$. Altogether, the inclusion $\ws\big(S^m(\eta)\big)\cap\as^{[K_m]}\subseteq\ws_S\cap\as^{[K_m]}$ is derived for each $m\in\NM$. Then Lemma~\ref{Chap6-Lem-AssDictInvSubstit} leads to the inclusion
$\ws\big(S^n(\eta)\big)\cap\as^{[K_m]}\subseteq\ws_S\cap\as^{[K_m]}$ for each $m\in\NM$ and $n\geq m$.

\vspace{.1cm}

The case of a primitive substitution $S:\as\to\as^+$ is similarly treated. Specifically, $l_0\in\NM$ satisfies that $b$ occurs in $S^{l_0}(a)$ for every pair $a,b\in\as$ by primitivity of $S$. Then $|S^{l_0+n}(a)|\geq n$ holds for each $n\in\NM$. With this at hand, the desired result is derived by following the same lines.
\end{proof}

\medskip

Next, a sufficient condition is provided so that $\Xi_S$ is periodically approxi\-mable for a primitive block substitution $S$.

\begin{theorem}[Periodic approximations, \cite{BeBeNi16}]
\label{Chap6-Theo-PerApprZdPrimSubst}
Let $S:\as\to\PaZd$ be a primitive block substitution over the alphabet $\as$ with at least two letters such that the associated subshift $\Xi_S$ is aperiodic. Suppose that there exists a strongly periodic $\eta\in\as^{\ZM^d}$ satisfying $\ws(\eta)\cap\as^{[K]}\subseteq\ws_S\cap\as^{[K]}$ where $K\subseteq\ZM^d$ fulfills that $\prod_{j=1}^d \{0,1\}\subseteq K$. Then the subshift $\Xi_S$ is periodically approximable and the sequence of strongly periodic subshifts $\Xi_n:=\Orb\big(S^n(\eta)\big)\,,\; n\in\NM\,,$ converges to $\Xi_S$ in the Hausdorff-topology on $\SZd\big(\as^{\ZM^d}\big)$.
\end{theorem}

\begin{proof}
According to Proposition~\ref{Chap6-Prop-AssDictSubshSubstZd}, the associated dictionary $\ws_S$ satisfies \nameref{(D1)}. Thus, $\ws(\eta)\cap\as^{[\tilde{K}]}\subseteq\ws_S\cap\as^{[\tilde{K}]}$ is deduced for $\tilde{K}:=\prod_{j=1}^d \{0,1\}\subseteq K$. Lemma~\ref{Chap6-Lem-SnETAStronglyPeriodic} implies that $\eta_n:=S^n(\eta)\in\as^{\ZM^d}$ is strongly periodic. The corresponding strongly periodic subshift $\Xi_n:=\Orb(\eta_n)$ satisfies $\ws(\eta_n)=\ws(\Xi_n)$ by Corollary~\ref{Chap2-Cor-DictOrbitSubshift}.

\vspace{.1cm}

Recall that the blocks $K_m:=\prod_{j=1}^d\{-2^{m}+1,\ldots,2^{m}\}\subseteq\ZM^d\,,\; m\in\NM\,,$ define an exhausting sequence of $\ZM^d$. Hence, it suffices to prove that, for each $m\in\NM$, there exists an $n(m)\in\NM$ such that 
$$
\ws(\eta_n)\cap\as^{[K_m]} \;
	= \;  \ws_S\cap\as^{[K_m]}\,,
	\qquad
	n\geq n(m)\,,
$$
by the following reason: In this case, the sequence of dictionaries $\big(\ws(\eta_n)\big)_{n\in\NM}$ converges to $\ws_S\in\DZd$ in the local pattern topology, c.f. Corollary~\ref{Chap2-Cor-DictExhSeq}. Consequently, the sequence of strongly periodic subshifts $\big(\Xi_n\big)_{n\in\NM}$ converges to the subshift $\Xi_S$ in $\SZd\big(\as^{\ZM^d}\big)$ by Theorem~\ref{Chap2-Theo-Shift+DictSpace}. Thus, $\Xi_S$ is periodically approximable.

\vspace{.1cm}

Let $m\in\NM$. Then $\ws(\eta_n)\cap\as^{[K_m]} \subseteq \ws_S\cap\as^{[K_m]}$ is derived for $n\geq m$ by Lemma~\ref{Chap6-Lem-GrowthAcceptPatternSubst}. For the converse inclusion, let $[w]\in\ws_S\cap\as^{[K_m]}$. According to definition of $\ws_S$, there is an $l_w\in\NM$ and a letter $a_w\in\as$ such that $[w]$ is a subpattern of $S^{l_w}(a_w)$. The substitution $S$ is primitive and so there is an $l_0\in\NM$ satisfying that $b$ occurs in $S^{l_0}(a)$ for every pair of letters $a,b\in\as$. Consequently, $[w]$ is a subpattern of $S^{l_0+l_w}\big(\eta(0)\big)$. Hence, $[w]$ is an element of $\ws\big(\eta_{l_0+l_w}\big)$. If $n\geq l_0+l_w$, then the letter $a_w$ occurs in $S^{n-l_w}(\eta)$ implying $[w]$ occurs in $\eta_n$ for $n\geq l_0+l_w$. Altogether, the desired equation $\ws(\eta_n)\cap\as^{[K_m]} = \ws_S\cap\as^{[K_m]}$ is derived for $n\geq n(m)$ where
$$
n(m) \; 
	:= \; \max\Big\{
		m\,,\; 
		\max\Big\{
			l_0+l_w
			\;\big|\; 
			[w]\in\ws_S\cap\as^{[K_m]}
		\Big\}
	\Big\}\,.
$$
Note that this maximum exists since the set $\as^{[K_m]}$ is finite.
\end{proof}

\begin{remark}
\label{Chap6-Rem-PerApprZdConvWord}
(i) Theorem~\ref{Chap6-Theo-PerApprZdPrimSubst} does not assert the convergence of $\big(S^n(\eta)\big)_{n\in\NM}$ to a $\xi\in\Xi_S$. It is possible that the sequence $\big(S^n(\eta)\big)_{n\in\NM}$ itself is not convergent, i.e., it can have more than one limit point. Theorem~\ref{Chap6-Theo-PerApprZdPrimSubst} only asserts the convergence of the corresponding dictionaries and subshifts. Due to compactness of $\as^{\ZM^d}$, a convergent subsequence $\big(S^{n_k}(\eta)\big)_{k\in\NM}$ can always be extracted. Then the limit $\lim_{k\to\infty}S^{n_k}(\eta)$ is necessarily an element of $\Xi_S$ by the convergence of the dictionaries $\lim_{k\to\infty}\ws\big(S^{n_k}(\eta)\big)=\ws_S$ in the local pattern topology. 

\vspace{.1cm}

(ii) Theorem~\ref{Chap6-Theo-PerApprZdPrimSubst} provides a general strategy to approximate substitutional subshifts $\Xi\in\SZd\big(\as^{\ZM^d}\big)$ by different subshifts $\Xi_n\in\SZd\big(\as^{\ZM^d}\big)\,,\; n\in\NM\,,$ which might be not strongly periodic. Specifically, let $\eta\in\as^{\ZM^d}$ (which is not necessarily strongly periodic) be so that $\ws(\eta)\cap\as^{[K]}\subseteq\ws_S\cap\as^{[K]}$ for a $K\subseteq\ZM^d$ satisfying $\prod_{j=1}^d \{0,1\}\subseteq K$. Then $\lim_{n\to\infty}\Xi_n=\Xi$ is concluded where $\Xi_n:=\overline{\Orb(S^n(\eta)}\,,\; n\in\NM$.
\end{remark}

In view of Remark~\ref{Chap6-Rem-BlockNotNec-d=1}, the condition that $S$ is a block substitution can be dropped if $d=1$.

\begin{corollary}[Periodic approximations for $d=1$]
\label{Chap6-Cor-PerApprZdPrimSubst-d=1}
Let $\as$ be an alphabet with at least two letters and $S:\as\to\as^+$ be a primitive substitution such that the associated subshift $\Xi_S$ is aperiodic. Consider an $\eta\in\as^\ZM$ and a $k_0\geq 2$ such that $\ws(\eta)\cap\as^{k_0} \subseteq \ws(\Xi_S)\cap\as^{k_0}$. Then the sequence of subshifts $(\Xi_n)_{n\in\NM}$ defined by $\Xi_n:=\overline{\Orb\big( S^n(\eta) \big)}$ converges to $\Xi_S$ in the Hausdorff-topology on $\SZ\big(\as^\ZM\big)$.
\end{corollary}

\begin{proof}
Let $S:\as\to\as^\ZM$ be a primitive substitution. Then $\Xi_S$ is periodically approximable by Corollary~\ref{Chap6-Cor-PerApprZdPrimSubst-d=1}. Hence, the existence of a strongly periodic $\eta\in\as^\ZM$ with $\ws(\eta)\cap\as^{k_0} \subseteq \ws(\Xi_S)\cap\as^{k_0}$ is guaranteed for every $k_0\geq 2$. Following the lines of the proof of Theorem~\ref{Chap6-Theo-PerApprZdPrimSubst}, the desired statement follows since Lemma~\ref{Chap6-Lem-SnETAStronglyPeriodic} and Lemma~\ref{Chap6-Lem-GrowthAcceptPatternSubst} are valid even if $S:\as\to\as^+$ is not a block substitution.
\end{proof}

\begin{corollary}
\label{Chap6-Cor-SubstitApprClosPath}
Let $\as$ be an alphabet with at least two letters and $S:\as\to\as^+$ be a primitive substitution with associated subshift $\Xi\in\SZ\big(\as^\ZM\big)$. Consider the de Bruijn graph $\gs_{k_0}$ of order $k_0\in\NM$ associated with $\ws(\Xi)$. Then the strongly periodic element $\eta_\wp\in\as^\ZM$ associated with a closed path $\wp$ of $\gs_{k_0}$ defines a sequence of strongly periodic subshifts $\Xi_n:=\Orb\big(S^n(\eta_\wp)\big)\,,\; n\in\NM\,,$ that converges to $\Xi\in\SZ\big(\as^\ZM\big)$.
\end{corollary}

\begin{proof}
Let $\wp$ be a closed path of $\gs_{k_0}$ for a $k_0\in\NM$. According to Lemma~\ref{Chap5-Lem-VerSetPerWor}, $\eta_\wp$ satisfies the inclusion $\ws(\eta_\wp)\cap\as^{k_0+1}\subseteq \es_{k_0} = \ws(\Xi)\cap\as^{k_0+1}$. Then Corollary~\ref{Chap6-Cor-PerApprZdPrimSubst-d=1} finishes the proof since $k_0+1\geq 2$.
\end{proof}

\medskip

The chapter is finished with a sufficient condition for a dictionary $\ws\in\mathfrak{D}_{\ZM^2}(\as)$ so that $\Xi(\ws)$ satisfies \nameref{(6.I)} of the strategy. More precisely, the existence of a strongly periodic $\eta\in\as^{\ZM^2}$ over a block $K:=\{0,1\}\times\{0,1\}\subseteq\ZM^2$, i.e., $\eta=v^\infty$ with $v\in\as^K$, is characterized for a dictionary $\ws$ such that $\ws(\eta)\cap\as^{[K]}\subseteq\ws\cap\as^{[K]}$.

\begin{proposition}
\label{Chap6-Prop-SuffPerApprZ2SymmPatt}
Let $\as$ be an alphabet and $\ws\in\mathfrak{D}_{\ZM^2}(\as)$ be a dictionary. Then the following assertions are equivalent. 
\begin{description}
\item[(i)] There exists a strongly periodic element $\eta:=v^\infty$ such that $v\in\as^{K}$ and $\ws(\eta)\cap\as^{[K]}\subseteq\ws\cap\as^{[K]}$ for $K:=\{0,1\}\times\{0,1\}\subseteq\ZM^2$.
\item[(ii)] There are letters $a,b,c,d\in\as$ (that may coincide) satisfying
$$
U \; := \;
\left\{
	\begin{pmatrix}
		a & b\\
		c & d
	\end{pmatrix}\,,\;
	\begin{pmatrix}
		c & d\\
		a & b
	\end{pmatrix}\,,\;
	\begin{pmatrix}
		b & a\\
		d & c
	\end{pmatrix}\,,\;
	\begin{pmatrix}
		d & c\\
		b & a
	\end{pmatrix}
\right\}\;
	\subseteq\ws\,.
$$
\end{description}
\end{proposition}

\begin{proof}
Let $v:=\begin{pmatrix}
		a & b\\
		c & d
\end{pmatrix}\in\as^{K}$ for $K:=\{0,1\}\times\{0,1\}\subseteq\ZM^2$. Then $\eta$ is represented by
\begin{gather*}
\eta \; := v^\infty \; = \;
\begin{array}{ccc}
& \vdots & \\
\ldots &
\begin{array}{cc|cc||cc|cc}
	a & b & a & b &		a & b & a & b\\
	c & d & c & \textcolor{gray}{d} &		\textcolor{gray}{c} & \textcolor{gray}{d} & \textcolor{gray}{c} & d\\
	\hline
	a & b & a & \textcolor{gray}{b} &		\textcolor{gray}{a} & \textcolor{gray}{b} & \textcolor{gray}{a} & b\\
	c & d & c & \textcolor{gray}{d} &		\textcolor{gray}{c} & \textcolor{gray}{d} & \textcolor{gray}{c} & d\\
	\hline				\hline
	a & b & a & \textcolor{gray}{b} &		\textcolor{gray}{a} & \textcolor{gray}{b} & \textcolor{gray}{a} & b\\
	c & d & c & d &		c & d & c & d\\
	\hline
	a & b & a & b &		a & b & a & b\\
	c & d & c & d &		c & d & c & d\\
\end{array} 
& \ldots \\
& \vdots & 
\end{array}
\in\as^{\ZM^2}\,.
\end{gather*}
Due periodicity of $\eta$ and since $v$ is supported on $K$, all patterns with support $[K]$ of $\eta$ occur in the pattern $\eta|_{F}$ with $F:=\{-1,0,1,2\}\times\{-1,0,1,2\}$. For convenience of the reader, the pattern $\eta|_{F}$ is printed in gray in the above representation. With this at hand, the equality $\ws(\eta)\cap\as^{[K]}=U$ is derived leading to the equivalence of (i) and (ii).
\end{proof}

\begin{remark}
\label{Chap6-Rem-SuffPerApprZ2SymmPatt}
(i) The table subshift and the Sierpinski carpet subshift satisfy this condition, c.f. Section~\ref{Chap7-Sect-TableSubst} and Section~\ref{Chap7-Sect-SierpinskiCarpetSubstitution}.

\vspace{.1cm}

(ii) Proposition~\ref{Chap6-Prop-SuffPerApprZ2SymmPatt} leads to a sufficient conditions of substitutional subshift over the group $\ZM^2$ to be periodically approximable in combination with Theorem~\ref{Chap6-Theo-PerApprZdPrimSubst}. 

\vspace{.1cm}

(iii) The assumption in Proposition~\ref{Chap6-Prop-SuffPerApprZ2SymmPatt} is essentially a constraint on the local symmetry of patterns in $\ws$. More precisely, the pattern
\begin{center}
\begin{tikzpicture}[>=stealth]
\draw[hellgrau, very thick] (-0.3,-0.3)--(0.3,-0.3);
\draw[hellgrau, very thick] (0,0)--(0,-0.6);
\draw[black] (-0.3,-0.05) node {$a$};
\draw[black] (0.3,0) node {$b$};
\draw[black] (-0.3,-0.65) node {$c$};
\draw[black] (0.3,-0.6) node {$d$};

\draw[hellgrau,very thick,->] (-0.3,0.2) to [out=45,in=170]++(0.8,0.25) to[out=-10,in=90]++(0.4,-0.5) to[out=-90,in=20] (0.5,-0.6);
\end{tikzpicture}
\end{center}
can be reflected at the gray axes so that the reflected pattern still belongs to $\ws$ as well as the rotation of the pattern around the center by $180^\circ$. 

\vspace{.1cm}

(iv) The assertion of Proposition~\ref{Chap6-Prop-SuffPerApprZ2SymmPatt} extends to blocks $K:=\{0,\ldots,n\}\times\{0,\ldots,n\}$. More precisely, there exists a strongly periodic element $\eta:=v^\infty$ such that $v\in\as^{K}$ and $\ws(\eta)\cap\as^{[K]}\subseteq\ws\cap\as^{[K]}$ if and only if the pattern $v$ satisfies certain symmetry properties in the dictionary $\ws$, c.f. (iii).

\vspace{.1cm}

(v) The assertion of Proposition~\ref{Chap6-Prop-SuffPerApprZ2SymmPatt} only provides a sufficient condition so that a subshift over the group $\ZM^2$ satisfies \nameref{(6.I)} while it is not a necessary condition. Specifically, there might exist subshifts $\Xi\in\as^{\ZM^2}$ so that there is no pattern $v\in\as^{[K]}$ fulfilling $\ws(v^\infty)\cap\as^{[K]}\subseteq\ws(\Xi)\cap\as^{[K]}$ where $K:=\{0,1\}\times\{0,1\}$.  Meanwhile, there could exists a pattern $u\in\as^{[F]}$ with block $[F]\neq[K]$ such that $\ws(u^\infty)\cap\as^{[K]}\subseteq\ws(\Xi)\cap\as^{[K]}$. In this case, different constraints on the local symmetries of the patterns are deduced on the dictionary $\ws$.

\vspace{.1cm}

(vi) Proposition~\ref{Chap6-Prop-SuffPerApprZ2SymmPatt} provides a class of subshifts of finite type 
$$
\Xi(K,U) \;
	:= \; \big\{
		 \xi\in\as^{\ZM^2} \;|\;
		 \ws(\xi)\cap\as^{[K]}\subseteq U
	\big\}
	\in\SG\big(\as^{\ZM^2}\big)
$$
so that $\Xi(K,U)$ contains periodic elements. More precisely, suppose there exists a primitive block substitution $S:\as\to\mathit{Pat}_{\ZM^2}(\as)$ such that $\ws_S$ satisfies the constraint of Proposition~\ref{Chap6-Prop-SuffPerApprZ2SymmPatt}. If $U=\ws_S\cap\as^{[F]}$ for a compact $F\in\ks(\ZM^2)$, then the corresponding subshift of finite type $\Xi(F,U)$ contains a periodic element by Theorem~\ref{Chap6-Theo-PerApprZdPrimSubst}, c.f. Section~\ref{Chap5-Sect-SymbDynSystZM}, Remark~\ref{Chap2-Rem-SubsFinTypPerSpectr}~(ii) and Section~\ref{Chap8-Sect-PerApprox}.
\end{remark}

\cleardoublepage


\chapter{Examples}
\label{Chap7-Examples}
\stepcounter{section}
\setcounter{section}{0}

Several examples are presented in this chapter to show how the results of this thesis apply to specific models. The main focus is on the existence of periodic approximations for pattern equivariant Schr\"odinger operators so that the spectra converge. In order to give a feeling for the structure of the de Bruijn graphs, some of the associated de Bruijn graphs are computed and sketched for the one-dimensional examples.

\medskip

The convergence of the spectra of periodic Schr\"odinger operators to aperiodic Schr\"odinger operators is proven in specific examples of substitutional systems if $d=1$, c.f. \cite[Proposition~5]{Sut87}. In \cite{BoGh93}, this result was extended to a class of primitive substitutions. Later \cite{LiTaWeWu02} provides a generalization of the work \cite{BoGh93} to all non-periodic primitive substitutions while only a semi-continuity of the spectra is proven and not the convergence spectra. In \cite[Proposition~3.9]{DaLiQu15} the convergence of the spectra is used to prove that the point spectrum of specific aperiodic Schr\"odinger operator is empty. Furthermore, \cite{LiQuYa16} used this convergence to estimate the behavior of the norms of the transfer matrices. These convergence results rely on transfer matrices and trace maps techniques. In addition to the convergence of the spectra, a monotone convergence is proven there. 

\medskip

In Chapter~\ref{Chap6-HigherDimPerAppr}, the convergence (and not only the semi-continuity) of the spectra is shown for all primitive substitutions for $d=1$ without using the methods of transfer matrices and trace maps. Furthermore, the developed theory is applicable to more general operators like pattern equivariant Schr\"odinger operators. Additionally, this approach applies to subshifts of $(\as^{\ZM^d},\ZM^d,\alpha)$ defined by substitutions if the subshift exhibits local symmetries of the patterns, c.f. Theorem~\ref{Chap6-Theo-PerApprZdPrimSubst} and Proposition~\ref{Chap6-Prop-SuffPerApprZ2SymmPatt}. The table substitution and the Sierpinski carpet substitution defined on the group $\ZM^2$ are analyzed in this chapter. We show that they satisfy these local symmetries. Altogether, some of the known results are confirmed and extended while the developed theory is independent of the dimension of the underlying system.

\medskip

Most of the examples considered here are defined via a primitive (block) substitution. Recall that a primitive block substitution $S:\as\to\PaZd$ or a primitive substitution $S:\as\to\as^+$ defines a dictionary $\ws_S$ and a subshift $\Xi_S$, c.f. Definition~\ref{Chap6-Def-AssDictSubshSubstZd} and Proposition~\ref{Chap6-Prop-AssDictSubshSubstZd}. The subshift $\Xi_S$ is minimal and, hence, topologically transitive, c.f. Proposition~\ref{Chap6-Prop-PrimSubstSubshZdMin}. Additionally, the subshift $\Xi_S$ is aperiodic if the substitution $S:\as\to\PaZd$ is injective, c.f.  \cite[Theorem~1.1]{Sol98}.

\medskip

Let $p:\as^\ZM\to\CM$ and $q:\as^\ZM\to\RM$ be pattern equivariant functions. Throughout this chapter, we consider, for $\xi\in\as^\ZM$, the bounded self-adjoint operators $J_\xi:\ell^2(\ZM)\to\ell^2(\ZM)$ defined by
\begin{equation}
\label{Chap7-Eq-JacobiOper}
\!\!\!\!
(J_\xi\psi)(m) 
	:= \left( 
			p\big(\alpha_{-m}(\xi)\big) \!\cdot \psi(m-1) + \overline{p\big(\alpha_{-(m+1)}(\xi)\big)} \!\cdot \psi(m+1)
		\right)
		+ q\big(\alpha_{-m}(\xi)\big) \!\cdot \psi(m)
\end{equation}
for $\psi\in\ell^2(\ZM)$ and $m\in \ZM$. We prove convergence results of the spectra of $J_\Xi$ for subshifts $\Xi\in\SZ(\as)$ by using Corollary~\ref{Chap4-Cor-CharSubshiftConvSpectrZM}, Theorem~\ref{Chap6-Theo-PerApprZdPrimSubst} and Corollary~\ref{Chap6-Cor-SubstitApprClosPath}. Clearly, the corresponding results hold also if the Jacobi operators $J_\xi\,,\; \xi\in\as^\ZM\,,$ is replaced by every other pattern equivariant Schr\"odinger operator, c.f. Theorem~\ref{Chap4-Theo-CharSubshiftConvSpectr}. Since the operators $J_\xi\,,\; \xi\in\as^\ZM\,,$ are self-adjoint and bounded, the spectrum $\sigma(J_\xi)$ is a compact subset of $\RM$ for $\xi\in\as^\ZM$. 

\medskip

For the sake of presentation, two different colors are used for the edges of the de Bruijn graphs. It is important to note that these colors do not have a meaning for the properties of the graphs or the specific edges. In the figures, the finite words corresponding to the vertices are printed in boldface and gray. Furthermore, the finite words associated with the edges are printed in the color of the corresponding edge.

\section{Full shift (Bernoulli shift)}
\label{Chap7-Sect-FullShift}

Let $\as$ be an alphabet. Then the {\em full shift}, also called {\em Bernoulli shift} in the literature, is the subshift containing each two-sided infinite word over the alphabet $\as$. Let $(\gs_k)_{k\in\NM}$ be the de Bruijn graphs associated with the full shift. More precisely, it is the dynamical system $(\as^\ZM,\ZM,\alpha)$. This example is interesting since every other subshift is contained in the full shift. Hence, every dictionary $\ws\in\DZ$ is contained in $\ws\big(\as^\ZM\big)=\PaZ$. Thus, the de Bruijn graph $\gs_k^\ws$ of order $k\in\NM$ associated with a dictionary $\ws\in\DZ$ is a subgraph of the de Bruijn graph $\gs_k$ of order $k$.

\medskip

\begin{figure}[htb]
\centering
\includegraphics[scale=0.9]{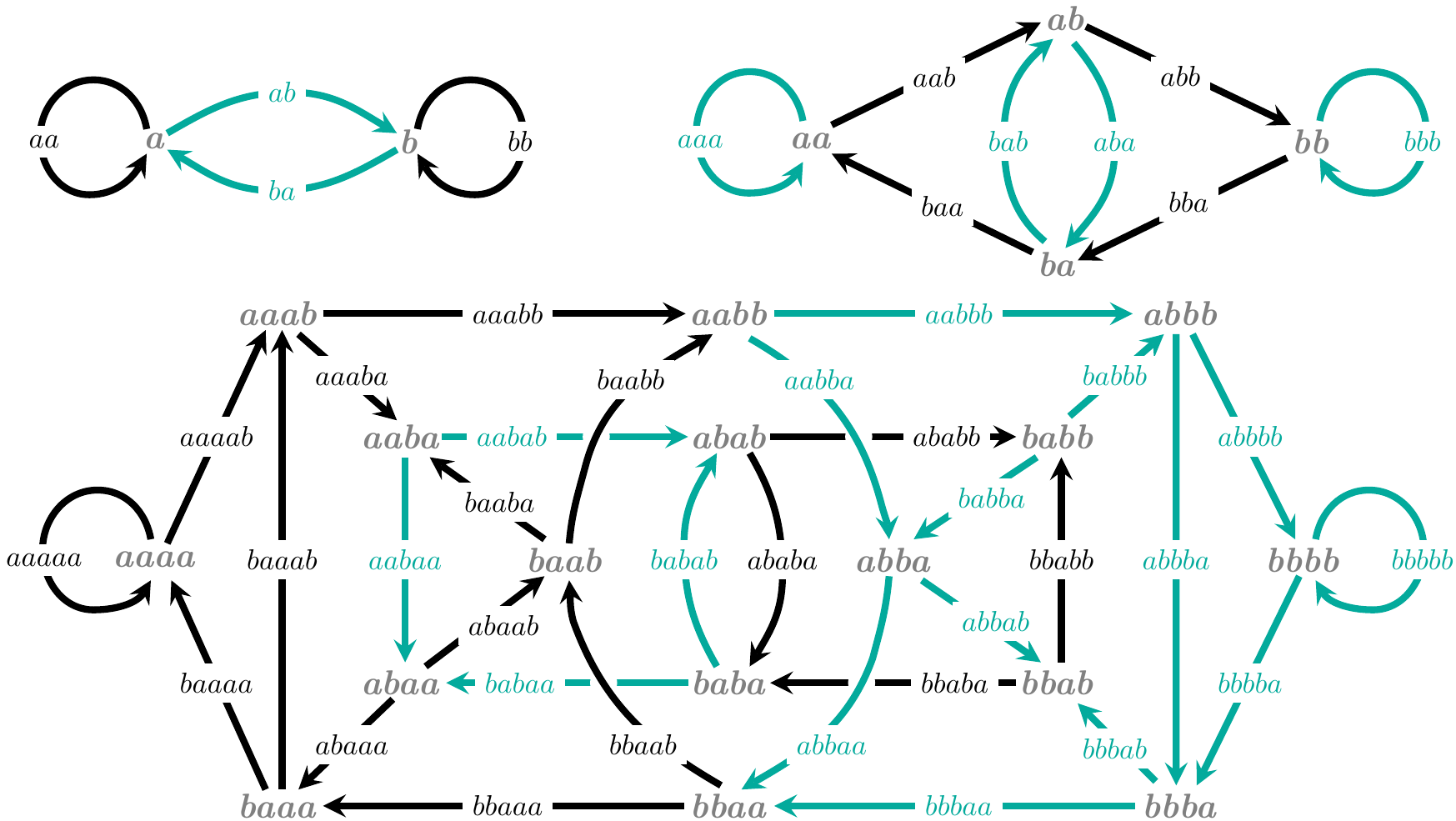}
\caption{The de Bruijn graphs of order $1$, $2$ and $4$ associated with the full shift.}
\label{Chap7-Fig-FullShift}
\end{figure}

The de Bruijn graph $\gs_4$ shows that the de Bruijn graphs are not planar in general, c.f. Figure~\ref{Chap7-Fig-FullShift}. Since each finite word over the alphabet $\as$ occurs in the full shift, the subword complexity function is given by $\comp(k):=\big(\sharp\as\big)^k$ for $k\in\NM$ and $\ws:=\ws(\as^\ZM)$. The full shift is aperiodic but not minimal. It is not difficult to check that the de Bruijn graphs associated with the full shift are all strongly connected. Hence, the full shift is periodically approximable by Theorem~\ref{Chap5-Theo-ExPerAppr}. Consider a strongly periodic sequence $(\Xi_n)_{n\in\NM}$ of subshifts tending to $\as^\ZM\in\SZ\big(\as^\ZM\big)$. Thanks to Corollary~\ref{Chap4-Cor-CharSubshiftConvSpectrZM}, the equality $\lim_{n\to\infty}\sigma\big(J_{\Xi_n}\big)=\sigma\big(J_{\as^\ZM}\big)$ is derived for every family of Jacobi operators $(J_\xi)_{\xi\in\as^\ZM}$ defined in Equation~\ref{Chap7-Eq-JacobiOper} where the limit is taken with respect to the Hausdorff metric on $\ks(\RM)$.

\section{Sturmian subshifts}
\label{Chap7-Sect-SturmSeq}

A minimal subshift $\Xi\in\SZ\big(\as^\ZM\big)$ is called {\em Sturmian} if $\sharp(\ws(\Xi)\cap\as^k)=k+1$ holds for all $k\in\NM$, c.f. \cite[Chapter 6]{Fogg02}. This is the least complex situation for aperiodic subshifts since a sequence with associated dictionary $\ws$ is strongly periodic if and only if there is a $k\in\NM$ satisfying $\comp(k)<k+\sharp\as-1$, c.f. \cite{MoHe40}. A Sturmian sequence has always an alphabet of exactly two letters which immediately follows from its definition. Thus, let $\as:=\{a,b\}$ be the alphabet. 

\medskip

Let $\Xi\in\SZ\big(\as^\ZM\big)$ be a subshift defined via a Sturmian subshift with dictionary $\ws:=\ws(\Xi)$. Then, the subword complexity function satisfies $\comp(k)=k+1$ for $k\in\NM$. Consequently, among the $k+1$ words of length $k$, only one can have two extensions on its right and all the other words of length $k$ being uniquely extendable to the right. Similarly, at most one word of length $k$ has two extensions to the left. Hence, the de Bruijn graph $\gs_k$ of order $k\in\NM$ admits exactly $k+1$ vertices and only one of them has two outgoing edges while only one admits two incoming edges. These two vertices might actually coincide. Consequently, the de Bruijn graph $\gs_k$ of order $k$ has either two branching vertices or only one branching vertex for every $k\in\NM$, c.f. Proposition~\ref{Chap5-Prop-BranPoiSubwComplex}.

\medskip

One of the most prominent example for a Sturmian sequence is the {\em Fibonacci subshift}, c.f. Example~\ref{Chap2-Ex-NonPeriodic-Fibonacci}, Example~\ref{Chap5-Ex-PerApprCutNotConv} and Example~\ref{Chap5-Ex-Fibonacci} This sequence was already known in the year 100 by {\sc Nikomachos von Gerasa}, c.f. \cite[Footnote~86]{Lan98}. The mathematician {\sc Leonardo de Pisa} (called Fibonacci) theoretically studied the evolution of rabbits by this sequence in 1202 in his famous \textit{Liber abbaci}, c.f. \cite[Chapter~XII]{Boncompagni57}. The Fibonacci sequence is defined via the primitive substitution $S:\as\to\as^+\,,\; a\mapsto ab\,,\; b\mapsto a$. Note that $S$ is not a block substitution. There exist de Bruijn graphs associated with the Fibonacci subshift having one or two branching vertices, c.f. Figure~\ref{Chap5-Fig-Fibonacci} and Proposition~\ref{Chap5-Prop-BranPoiSubwComplex}. A $2$-periodic point with respect to the substitution, i.e., $S^2(\xi)=\xi$, is given by the limit
$$
\lim_{n\to\infty} S^{2n}\big((b|a)^\infty\big) \;
	= \; \ldots \, b\, a\, b\, a\, a\, b\, a\, a\, b\, a\, b\, a\, a\, b\, a\, a\, b\, |\,  a\, b\, a\, a\, b\, a\, b\, a\, a\, b\, a\, a\, b\, a\, b\, a\, a\, \ldots
	\;\in\as^\ZM\,.
$$

\begin{proposition}[Periodic approximation of the Fibonacci operator]
\label{Chap7-Prop-ConvFibonacci}
Let $\Xi_S$ be the subshift generated by the Fibonacci substitution $S:\as\to\as^+$ over the alphabet $\as:=\{a,b\}$ and $J_\xi:\ell^2(\ZM)\to\ell^2(\ZM)\,,\; \xi\in\as^\ZM\,,$ be a Jacobi operator defined in Equation~\ref{Chap7-Eq-JacobiOper}. Consider the strongly periodic configurations $\eta_n^a:=\big( S^n(a) \big)^\infty=S^n\big(a^\infty\big)$ and $\eta_n^b:=\big( S^n(b) \big)^\infty=S^n\big(b^\infty\big)$ for $n\in\NM$. Then the equalities 
$$
\lim_{n\to\infty} \sigma\big( J_{\eta_n^a} \big) \;
	= \; \lim_{n\to\infty} \sigma\big( J_{\eta_n^b} \big) \;
	= \; \sigma\big( J_\xi \big) \;
	= \; \sigma_{ess}\big(J_\xi\big) \;
	= \; \sigma\big(J_{\Xi_S}\big)\,,
	\qquad
	\xi\in\Xi_S\,,
$$
hold where the limit is taken with respect to the Hausdorff metric on $\ks(\RM)$.
\end{proposition}

\begin{proof}
Since the substitution $S:\as\to\as^+$ is primitive, the subshift $\Xi_S$ is minimal, c.f. Proposition~\ref{Chap6-Prop-PrimSubstSubshZdMin}. Thus, the equations $\sigma\big(H_\xi\big)=\sigma_{ess}\big(H_\xi\big)=\sigma\big(H_{\Xi_S}\big)\,,\;\xi\in\Xi_S\,,$ follow by Theorem~\ref{Chap2-Theo-ConstSpectrMinimal}. According to Lemma~\ref{Chap6-Lem-SnETAStronglyPeriodic}, the elements $\eta_n^a, \eta_n^b\in\as^\ZM$ are strongly periodic and their period is bounded from below by $n+1$ for $n\in\NM$, c.f. Corollary~\ref{Chap5-Cor-PerGrowth}.

\vspace{.1cm}

Due to Corollary~\ref{Chap4-Cor-CharSubshiftConvSpectrZM}, it suffices to check that the sequences of subshifts $\Xi_n^a:=\Orb(\eta_n^a)\,,$ $n\in\NM\,,$ and $\Xi_n^b:=\Orb(\eta_n^b)\,,\; n\in\NM\,,$ converge to $\Xi_S$ in $\SZ\big(\as^\ZM\big)$, i.e., $\lim_{n\to\infty}\Xi_n^a=\Xi_S=\lim_{n\to\infty}\Xi_n^b$. Let $(\gs_k)_{k\in\NM}$ be the de Bruijn graphs associated with the dictionary $\ws_S$. The path $\wp_a:=(aa)$ consisting of only one edge in the de Bruijn graph $\gs_1$ is closed, c.f. Figure~\ref{Chap5-Fig-Fibonacci}. The associated strongly periodic word to $\wp_a$ is given by $\eta_0^a:=a^\infty$. Then Corollary~\ref{Chap6-Cor-SubstitApprClosPath} implies the desired convergence $\lim_{n\to\infty}\Xi_n^a=\Xi_S$. Since $S^n(b)=S^{n-1}(a)$ holds for $n\in\NM$, the equations $\Xi_n^b=\Xi_{n-1}^a\,,\; n\in\NM\,,$ are derived implying $\lim_{n\to\infty}\Xi_n^b = \Xi_S$ by the previous considerations.
\end{proof}

\medskip
 
\begin{figure}[htb]
\centering
\includegraphics[scale=0.97]{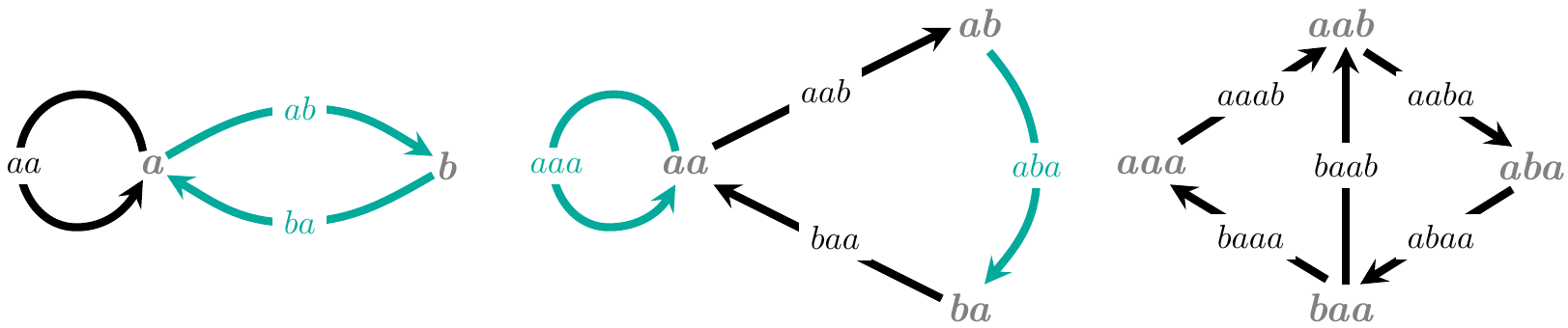}
\caption{The de Bruijn graphs of order $1$, $2$ and $3$ associated with the Silver mean subshift.}
\label{Chap7-Fig-SilverMean}
\end{figure}

Another example of a Sturmian sequence is provided by the {\em Silver mean subshift} obtained from the primitive substitution $S:\as\to\as^+\,,\; a\mapsto aab\,,\; b\mapsto a$. The Silver-mean substitution is primitive but $S$ is not a block substitution. The induced subshift $\Xi_S$ is minimal and aperiodic by Proposition~\ref{Chap6-Prop-PrimSubstSubshZdMin} and \cite[Theorem~1.1]{Sol98}, see \cite[Example~4.5]{BaakeGrimm13} as well. A $2$-periodic point with respect to the substitution $S$, i.e., $S^2(\xi)=\xi$, is given by the limit
$$
\lim_{n\to\infty} S^{2n}\big((a|b)^\infty\big) \;
	= \; \ldots \, a\, a\, b\, a\, a\, b\, a\, a\, a\, b\, a\, a\, b\, a\, a\, a\, b\, |\,  a\, a\, b\, a\, a\, b\, a\, a\, a\, b\, a\, a\, b\, a\, a\, a\, b\, \ldots
	\;\in\as^\ZM\,.
$$
Like in the Fibonacci case, there exist de Bruijn\- graphs having one branching vertex and two branching vertices, c.f. Figure~\ref{Chap7-Fig-SilverMean}. 

\begin{proposition}[Periodic approximation of the Silver mean operator]
\label{Chap7-Prop-ConvSilverMean}
Let $\Xi_S$ be the subshift generated by the Silver mean substitution $S:\as\to\as^+$ over the alphabet $\as:=\{a,b\}$ and $J_\xi:\ell^2(\ZM)\to\ell^2(\ZM)\,,\; \xi\in\as^\ZM\,,$ be a Jacobi operator defined in Equation~\ref{Chap7-Eq-JacobiOper}. Consider the strongly periodic configurations $\eta_n^a:=\big( S^n(a) \big)^\infty=S^n\big(a^\infty\big)$ and $\eta_n^b:=\big( S^n(b) \big)^\infty=S^n\big(b^\infty\big)$ for $n\in\NM$. Then the equalities 
$$
\lim_{n\to\infty} \sigma\big( J_{\eta_n^a} \big) \;
	= \; \lim_{n\to\infty} \sigma\big( J_{\eta_n^b} \big) \;
	= \; \sigma\big( J_\xi \big) \;
	= \; \sigma_{ess}\big(J_\xi\big) \;
	= \; \sigma\big(J_{\Xi_S}\big)\,,
	\qquad
	\xi\in\Xi_S\,,
$$
hold where the limit is taken with respect to the Hausdorff metric on $\ks(\RM)$.
\end{proposition}

\begin{proof}
Let $\Xi_S$ be the Silver mean subshift with corresponding dictionary $\ws_S$. The assertion is similarly proved like Proposition~\ref{Chap7-Prop-ConvFibonacci} with the following additional remarks. The path $\wp_a:=(aa)$ in the de Bruijn graph $\gs_1$ associated with $\ws_S$ is closed with associated periodic word $\eta_0^a=a^\infty$, c.f. Figure~\ref{Chap7-Fig-SilverMean}. Furthermore, the equation $S^n(b)=S^{n-1}(a)$ holds for $n\in\NM$. Hence, $\lim_{n\to\infty}\Orb(\eta_n^b)=\lim_{n\to\infty}\Orb(\eta_n^a)=\Xi_S$ is derived by Corollary~\ref{Chap6-Cor-SubstitApprClosPath}. Then Corollary~\ref{Chap4-Cor-CharSubshiftConvSpectrZM} and Theorem~\ref{Chap2-Theo-ConstSpectrMinimal} lead to the desired result.
\end{proof}

\begin{figure}[htb]
\centering
\includegraphics[scale=0.71]{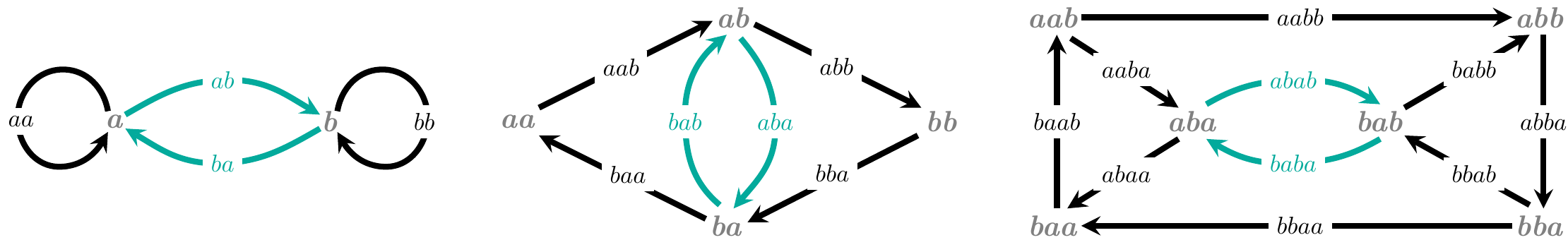}
\caption{The de Bruijn graphs of order $1$, $2$ and $3$ associated with the Prouhet-Thue-Morse subshift.}
\label{Chap7-Fig-ThueMorsGraph2}
\end{figure}

\section{The Prouhet-Thue-Morse subshift}
\label{Chap7-Sect-ThueMorse}

The first appearance of the Prouhet-Thue-Morse sequence goes back to {\sc Prouhet} \cite{Prou51} in 1851 without mentioning the sequence explicitly. Later {\sc Thue} \cite{Thu06,Thu12} introduced the sequence in order to study the combinatorial properties of words. {\sc Morse} \cite{Mors21} used the sequence to produce examples of geodesics on Riemannian manifolds of negative curvature independently which were neither periodic nor almost periodic. In the literature, the sequence is usually called {\em Thue-Morse sequence}.

\medskip

The Prouhet-Thue-Morse subshift is defined by the primitive block substitution $S:\as\to\as^+\,,\; a\mapsto ab\,,\, b\mapsto ba$ over the alphabet $\as=\{a,b\}$. A $2$-periodic point of the substitution, i.e., $S^2(\xi)=\xi$, is given by the limit
$$
\lim_{n\to\infty}S^{2n}\big(b^\infty\big) \;
	= \; \ldots \, b\, a\, a\, b\, a\, b\, b\, a\, a\, b\, b\, a\, b\, a\, a\, b\,|\, b\, a\, a\, b\, a\, b\, b\, a\, a\, b\, b\, a\, b\, a\, a\, b\, \ldots 
	\;\in\as^\ZM\,.
$$
The Prouhet-Thue-Morse subshift $\Xi_S$ is minimal and aperiodic by Proposition~\ref{Chap6-Prop-PrimSubstSubshZdMin} and \cite[Theorem~1.1]{Sol98}, c.f. \cite[Proposition~5.1.2]{Fogg02}. Note that the Prouhet-Thue-Morse subshift is not a pisot type substitution since the incidence matrix has eigenvalues $1$ and $0$, c.f. \cite[page~11]{Fogg02}. The reader is referred to \cite[Chapter~5]{Fogg02} and \cite[Section~4.6]{BaakeGrimm13} for a detailed discussion on the properties of the Prouhet-Thue-Morse subshift.

\begin{proposition}[Periodic approximation of the Prouhet-Thue-Morse operator]
\label{Chap7-Prop-ConvThueMorse}
Let $\Xi_S$ be the subshift generated by the Prouhet-Thue-Morse substitution $S:\as\to\as^+$ over the alphabet $\as:=\{a,b\}$ and $J_\xi:\ell^2(\ZM)\to\ell^2(\ZM)\,,\; \xi\in\as^\ZM\,,$ be a Jacobi operator defined in Equation~\ref{Chap7-Eq-JacobiOper}. Consider the strongly periodic configurations $\eta_n^a:=\big( S^n(a) \big)^\infty=S^n\big(a^\infty\big)$ and $\eta_n^b:=\big( S^n(b) \big)^\infty=S^n\big(b^\infty\big)$ for $n\in\NM$. Then the equalities 
$$
\lim_{n\to\infty} \sigma\big( J_{\eta_n^a} \big) \;
	= \; \lim_{n\to\infty} \sigma\big( J_{\eta_n^b} \big) \;
	= \; \sigma\big( J_\xi \big) \;
	= \; \sigma_{ess}\big(J_\xi\big) \;
	= \; \sigma\big(J_{\Xi_S}\big)\,,
	\qquad
	\xi\in\Xi_S\,,
$$
hold where the limit is taken with respect to the Hausdorff metric on $\ks(\RM)$.
\end{proposition}

\begin{proof}
Let $\Xi_S$ be the Prouhet-Thue-Morse subshift with corresponding dictionary $\ws_S$. The assertion is similarly proved like Proposition~\ref{Chap7-Prop-ConvFibonacci} with the following additional remarks. The paths $\wp_a:=(aa)$ and $\wp_b:=(bb)$ in the de Bruijn graph $\gs_1$ associated with $\ws_S$ are closed paths. The associated periodic words are denoted by $\eta_0^a=a^\infty$ and $\eta_0^b=b^\infty$, c.f. Figure~\ref{Chap7-Fig-ThueMorsGraph2}. Hence, $\lim_{n\to\infty}\Orb(\eta_n^a)=\lim_{n\to\infty}\Orb(\eta_n^b)=\Xi_S$ is derived by Corollary~\ref{Chap6-Cor-SubstitApprClosPath}. Then Corollary~\ref{Chap4-Cor-CharSubshiftConvSpectrZM} and Theorem~\ref{Chap2-Theo-ConstSpectrMinimal} lead to the desired result.
\end{proof}

\section{The Period-doubling subshift}
\label{Chap7-Sect-PerDoubl}

The Period-doubling sequence has been studied by several authors during the past decades \cite{BeBoGh91,BoGh93,HoKnSi95,Dam98PerDoub,Dam98,Dam00}. It is a factor of the Prouhet-Thue-Morse subshift, c.f. \cite[Section~4.5]{BaakeGrimm13} and refer\-ences therein for more details. 

\medskip

The Period-doubling subshift is defined via the primitive block substitution $S:\as\to\as^+\,,\; a\mapsto ab \,,\; b\mapsto aa\,,$ over the alphabet $\as:=\{a,b\}$. A $2$-periodic point of the substitution $S$, i.e., $S^2(\xi)=\xi$, is given by the limit
$$
\lim_{n\to\infty} S^{2n}\big((b|a)^\infty\big) \;
	:= \; \ldots\, a\, b\, a\, a\, a\, b\, a\, b\, a\, b\, a\, a\, a\, b\, a\, b\, |\, a\, b\, a\, a\, a\, b\, a\, b\, a\, b\, a\, a\, a\, b\, a\, a\, \ldots 
	\;\in\as^\ZM\,.
$$
The associated subshift $\Xi_S\in\SZ\big(\as^\ZM\big)$ is minimal and aperiodic by Proposition~\ref{Chap6-Prop-PrimSubstSubshZdMin} and \cite[Theorem~1.1]{Sol98}.

\begin{figure}[htb]
\centering
\includegraphics[scale=0.836]{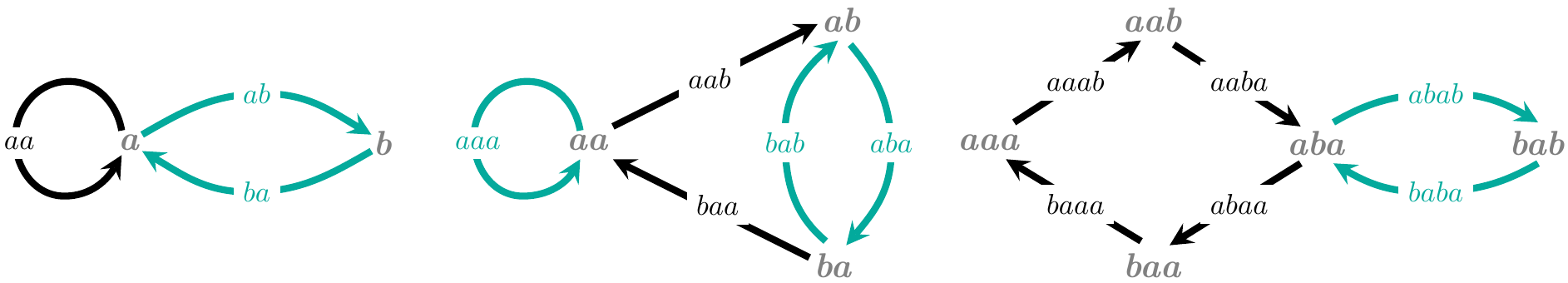}
\caption{The de Bruijn graphs of order $1$, $2$ and $3$ associated with the Period-doubling subshift.}
\label{Chap7-Fig-PeriodDoubling}
\end{figure}

\begin{proposition}[Periodic approximation of the Period-Doubling operator]
\label{Chap7-Prop-ConvPeriod-Doubling}
Let $\Xi_S$ be the subshift generated by the Period-doubling substitution $S:\as\to\as^+$ over the alphabet $\as:=\{a,b\}$ and $J_\xi:\ell^2(\ZM)\to\ell^2(\ZM)\,,\; \xi\in\as^\ZM\,,$ be a Jacobi operator defined in Equation~\ref{Chap7-Eq-JacobiOper}. Consider the strongly periodic configurations $\eta_n^a:=\big( S^n(a) \big)^\infty=S^n\big(a^\infty\big)$ and $\eta_n^b:=\big( S^n(b) \big)^\infty=S^n\big(b^\infty\big)$ for $n\in\NM$. Then the equalities 
$$
\lim_{n\to\infty} \sigma\big( J_{\eta_n^a} \big) \;
	= \; \lim_{n\to\infty} \sigma\big( J_{\eta_n^b} \big) \;
	= \; \sigma\big( J_\xi \big) \;
	= \; \sigma_{ess}\big(J_\xi\big) \;
	= \; \sigma\big(J_{\Xi_S}\big)\,,
	\qquad
	\xi\in\Xi_S\,,
$$
hold where the limit is taken with respect to the Hausdorff metric on $\ks(\RM)$.
\end{proposition}

\begin{proof}
Let $\Xi_S$ be the Period-doubling subshift with corresponding dictionary $\ws_S$. The assertion is similarly proved like Proposition~\ref{Chap7-Prop-ConvFibonacci} with the following additional remarks. The path $\wp_a:=(aa)$ in the de Bruijn graph $\gs_1$ associated with $\ws_S$ is closed with associated periodic word $\eta_0^a=a^\infty$, c.f. Figure~\ref{Chap7-Fig-PeriodDoubling}. Furthermore, the equation $S^n(b)=S^{n-1}(aa)$ holds for $n\in\NM$. The strongly periodic word $(aa)^\infty$ corresponds to the closed path $\wp_{aa}:=(aaa)$ consisting of one edge in the de Bruijn graph $\gs_2$ of order $2$, c.f. Figure~\ref{Chap7-Fig-PeriodDoubling}. Hence, the equalities $\lim_{n\to\infty}\Orb(\eta_n^a) = \lim_{n\to\infty}\Orb(\eta_n^b) = \Xi_S$ is derived by Corollary~\ref{Chap6-Cor-SubstitApprClosPath}. Then Corollary~\ref{Chap4-Cor-CharSubshiftConvSpectrZM} and Theorem~\ref{Chap2-Theo-ConstSpectrMinimal} lead to the desired result.
\end{proof}

\section{The Rudin-Shapiro subshift}
\label{Chap7-Sect-GolRudSha}

The Rudin-Shapiro subshift was first introduced by {\sc Shapiro} \cite{Sha51} and later by {\sc Rudin} \cite{Rud59} for estimations in harmonic analysis. {\sc Golay} \cite{Gol49} presented a predecessor in 1949. The reader is referred to \cite[Section~2.2]{Fogg02}, \cite[Section~4.7]{BaakeGrimm13} and references therein for more details. There are two representations of the sequence, one with two letters and one with four.

\medskip

\begin{figure}[htb]
\centering
\includegraphics[scale=0.82]{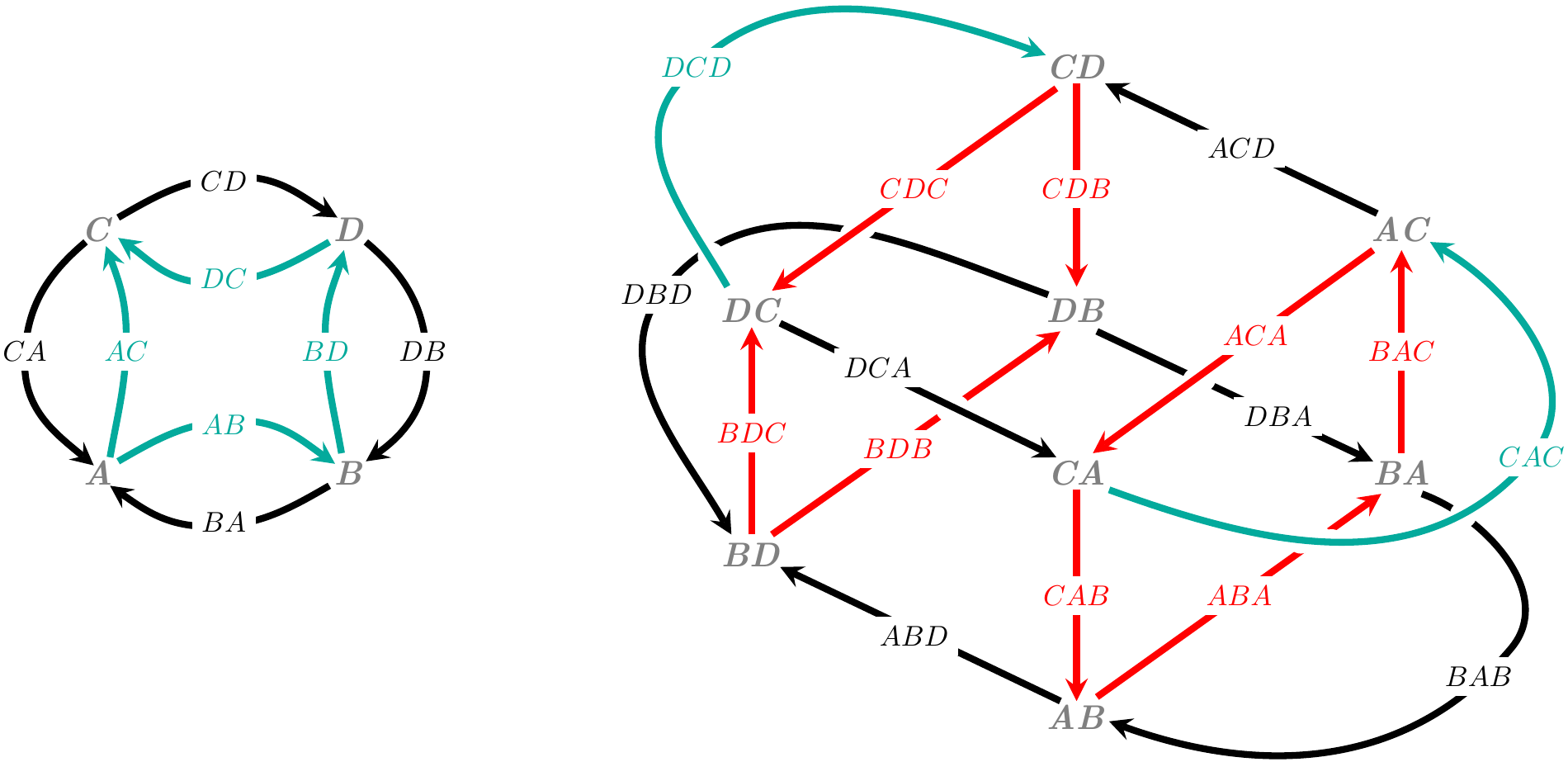}
\caption{The de Bruijn graphs of order $1$ and $2$ associated with the Rudin-Shapiro subshift.}
\label{Chap7-Fig-GolRudSha}
\end{figure}

The map $S':\as^+\to\as^+$ defined by
$$
aa\mapsto aaab
	\,,\qquad 
	ab\mapsto aaba
		\,,\qquad
		ba\mapsto bbab
			\,,\qquad
			bb\mapsto bbba\,,
$$
over the alphabet $\as:=\{a,b\}$ generates a subshift. Now, consider the primitive block substitution $S:\bs\to\bs^+$ over the alphabet $\bs=\{A,B,C,D\}$ defined by
$$
A\mapsto AB\,,
   \hspace{.7cm}
    B\mapsto AC\,,
     \hspace{.7cm}
      C\mapsto DB\,,
       \hspace{.7cm}
        D\mapsto DC\,.
$$
Then there exists a unique subshift $\Xi_S\in\SZ\big(\bs^\ZM\big)$ defined via the fixed points of the substitution $S$, c.f. Definition~\ref{Chap6-Def-AssDictSubshSubstZd} and Proposition~\ref{Chap6-Prop-PrimSubstSubshZdMin}. The substitution $S$ is induced by the map $S':\as^+\to\as^+$ by using $A:=aa\,,\; B:=ab\,,\; C:=ba\,,\; D:=bb$. A fixed point of the substitution $S$ is given by the limit $\lim_{n\to\infty}S^n\big((B|A)^\infty\big)$ represented by
$$
\ldots A\, B\, A\, C\, A\, B\, D\, B\, D\, C\, D\, B\, A\, B\, D\, B\, |\, A\, B\, A\, C\, A\, B\, D\, B\, A\, B\, A\, C\, D\, C\, A\, C \ldots 
	\in\bs^\ZM\,.
$$
The Rudin-Shapiro subshift $\Xi_S$ is minimal and aperiodic by Proposition~\ref{Chap6-Prop-PrimSubstSubshZdMin} and \cite[Theorem~1.1]{Sol98}. Note that the de Bruijn graph of order $2$ associated with the Rudin-Shapiro sequence is not planar, c.f. Figure~\ref{Chap7-Fig-GolRudSha}.

\begin{figure}[htb]
\centering
\includegraphics[scale=0.7]{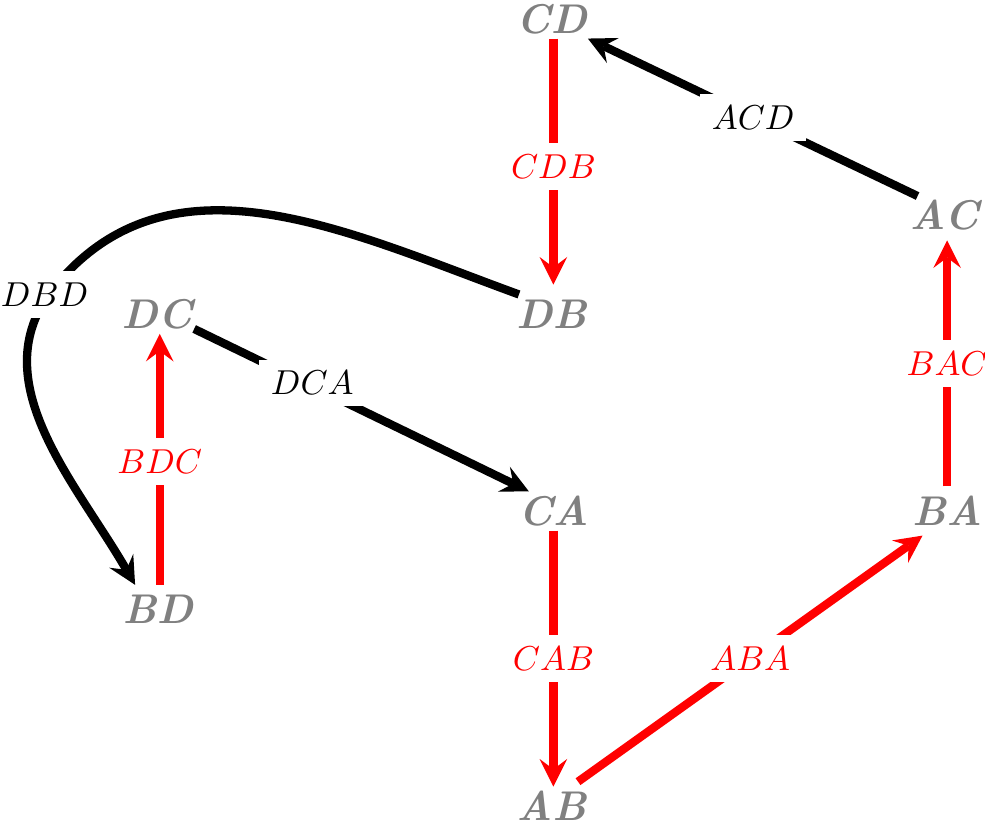}
\caption{The closed path $\wp= \big( DCA, CAB, ABA, BAC, ACD, CDB, DBD, BDC \big)$ in the de Bruijn graph of order $2$ associated with the Rudin-Shapiro subshift.}
\label{Chap7-Fig-ClosPathRudSha}
\end{figure}

\begin{remark}
\label{Chap7-Rem-RudShapClosPathNotElDict}
Let $\ws$ be a dictionary with associated de Bruijn graphs $(\gs_k)_{k\in\NM}$. According to Definition~\ref{Chap5-Def-AssPerWor}, a closed path $\wp:=(e_1,\ldots,e_l)$ for $l\in\NM$ in a de Bruijn graph $\gs_k$ of order $k\in\NM$ gives rise to a strongly periodic words $\eta:=v^\infty\in\as^\ZM$ satisfying $\ws(\eta)\cap\as^k\subseteq\vs_k$ where $v\in\as^l$. In general, the word $v$ is not an element of $\ws$ even if $\wp$ is a path that visit no vertex and no edge twice. For instance, consider the path, given in Figure~\ref{Chap7-Fig-ClosPathRudSha},
$$
\wp \;
	= \; \big( DCA, CAB, ABA, BAC, ACD, CDB, DBD, BDC \big)
$$
in the de Bruijn graph $\gs_2$ of order $2$ associated with the dictionary of the Rudin-Shapiro-sequence $\ws_S$. This closed path defines a strongly periodic word $\eta:=v^\infty$ where $v$ is the finite word $DCABACDB\in\as^8$, c.f. Definition~\ref{Chap5-Def-AssPerWor}. Then the subword $DCABA$ of $v$ can only arise by applying the substitution $S$ to the word $DAB$. Furthermore, $DA$ does never occur as a subword of the Rudin-Shapiro sequence. Thus, $v\not\in\ws_S$ follows.
\end{remark}

\begin{proposition}[Periodic approximation of the Rudin-Shapiro operator]
\label{Chap7-Prop-ConvRudin-Shapiro}
Let $\Xi_S$ be the subshift generated by the Rudin-Shapiro substitution $S:\bs\to\bs^+$ over the alphabet $\bs:=\{A,B,C,D\}$ and $J_\xi:\ell^2(\ZM)\to\ell^2(\ZM)\,,\; \xi\in\as^\ZM\,,$ be a Jacobi operator defined in Equation~\ref{Chap7-Eq-JacobiOper}. Consider the strongly periodic configurations $\eta_n^x:=\big( S^n(x) \big)^\infty=S^n\big(x^\infty\big)$ for the letter $x\in\bs$ and $n\in\NM$. Then the equalities 
$$
\lim_{n\to\infty} \sigma\big( J_{\eta_n^x} \big) \;
	= \; \sigma\big( J_\xi \big) \;
	= \; \sigma_{ess}\big(J_\xi\big) \;
	= \; \sigma\big(J_{\Xi_S}\big)\,,
	\qquad
	x\in\bs\,,\;
	\xi\in\Xi_S\,,
$$
hold where the limit is taken with respect to the Hausdorff metric on $\ks(\RM)$.
\end{proposition}

\begin{proof}
Let $\Xi_S$ be the Rudin-Shapiro subshift with dictionary $\ws_S$. The assertion is similarly proved like Proposition~\ref{Chap7-Prop-ConvFibonacci} with the following additional considerations. For $n\in\NM$, the identities 
\begin{align*}
S^n(A) \; 
	&= \; S^{n-1}(AB)\,,\qquad
S^n(B) \; 
	= \; S^{n-1}(AC)\,,\\
S^n(C) \; 
	&= \; S^{n-1}(DB)\,,\qquad
S^n(D) \; 
	= \; S^{n-1}(DC)\,,
\end{align*}
are derived. Consider the strongly periodic words $\eta_A:= (AB)^\infty\,,\; \eta_B:= (AC)^\infty\,,\; \eta_C:= (DB)^\infty$ and $\eta_D:= (DC)^\infty$. Then $\eta_x$ is the associated periodic word of the closed path $\wp_x$ in the de Bruijn graph $\gs_1$ of order $1$ for $x\in\bs$ where $\wp_A := (AB, BA)\,,\; \wp_B := (AC, CA)\,,\; \wp_C := (DB, BD)$ and $\wp_D := (DC, CD)$, c.f. Figure~\ref{Chap7-Fig-GolRudSha}. Consequently, the equation $\lim_{n\to\infty}\Orb(\eta_n^x)=\Xi_S$ is derived for each letter $x\in\bs$ by Corollary~\ref{Chap6-Cor-SubstitApprClosPath}. Then Corollary~\ref{Chap4-Cor-CharSubshiftConvSpectrZM} and Theorem~\ref{Chap2-Theo-ConstSpectrMinimal} lead to the desired result.
\end{proof}

\section{The one-defect}
\label{Chap7-Sect-OneDefect}

\begin{figure}[htb]
\centering
\includegraphics[scale=0.86]{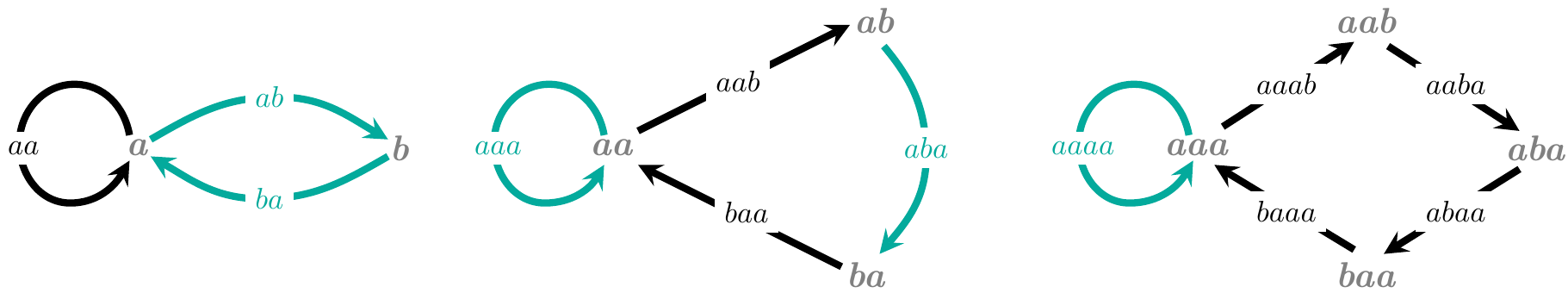}
\caption{The de Bruijn graphs of order $1$, $2$ and $3$ associated with the one-defect sequence.}
\label{Chap7-Fig-OneDefect}
\end{figure}

The {\em one-defect} subshift is defined over the alphabet $\as:=\{a,b\}$. The two-sided infinite word $\xi\in\as^\ZM$ is given by
$$
\xi \;
	:= \; \ldots a\, a\, a\, a\, a\, a\, a\, a\, a\, a\, a\, a\, a\, | \, b\, a\, a\, a\, a\, a\, a\, a\, a\, a\, a\, a\, a\, \ldots \as^\ZM\,, 
$$
i.e., $\xi(0):=b$ and $\xi(j):= a$ for all $j\in\ZM\setminus\{0\}$. The induced subshift $\Xi:=\overline{\Orb(\xi)}$ is equal to $\Orb(\xi)\cup\{ \eta\}$ where $\eta:=a^\infty$ is a periodic word satisfying $\alpha_1(\eta)=\eta$. The dictionary $\ws(\Xi)$ of $\Xi$ is equal to $\ws(\xi)$ by Corollary~\ref{Chap2-Cor-DictOrbitSubshift}. Since $\xi$ is not periodic, the subshift $\Xi$ is aperiodic. Clearly, $\eta$ satisfies $\Orb(\eta):=\{\eta\}$. Thus, the subshift $\Xi$ is not minimal. The associated de Bruijn graph $\gs_k$ of order $k\in\NM$ consists of two circles that attach at one vertex $v:=a^k\in\vs_k$, c.f. Figure~\ref{Chap7-Fig-OneDefect}. 

\begin{proposition}[Periodic approximation of the one-defect operator]
\label{Chap7-Prop-ConvOne-Defect}
Let $\Xi$ be the one-defect subshift and $J_\xi:\ell^2(\ZM)\to\ell^2(\ZM)\,,\; \xi\in\as^\ZM\,,$ be a Jacobi operator defined in Equation~\ref{Chap7-Eq-JacobiOper}. Consider the strongly periodic configurations $\eta_n:=(ba^n)^\infty$ for $n\in\NM$. Then $\Xi$ is periodically approximable by the subshifts $\Xi_n:=\Orb(\eta_n)\,,\; n\in\NM\,,$ and the equality
$$
\lim_{n\to\infty} \sigma\big( J_{\eta_n} \big) \;
	= \; \sigma\big(J_{\Xi}\big)
$$
holds where the limit is taken with respect to the Hausdorff metric on $\ks(\RM)$.
\end{proposition}

\begin{proof}
Let $\Xi_S$ be the one-defect subshift with dictionary $\ws_S$. For $n\in\NM$, the word $a^n$ is occurs infinitely often in $\xi$ to the left and to the right. Consequently, $\Xi$ is periodically approximable by Corollary~\ref{Chap5-Cor-SuffCondPerAppr}. By definition of $\eta_n$, the equality $\ws(\eta_n)\cap\as^{n}=\ws(\Xi)\cap\as^{n}$ is derived for each $n\in\NM$. Consequently, the dictionaries $\big(\ws(\Xi_n)\big)_{n\in\NM}$ converge in the local pattern topology to $\ws_S$ since $\ws(\Xi_n)=\ws(\eta_n)$ by Corollary~\ref{Chap2-Cor-DictOrbitSubshift}. Hence, $\lim_{n\to\infty}\Xi_n=\Xi$ is concluded by Theorem~\ref{Chap2-Theo-Shift+DictSpace}. Then Corollary~\ref{Chap4-Cor-CharSubshiftConvSpectrZM} and Theorem~\ref{Chap2-Theo-ConstSpectrMinimal} lead to the desired result.
\end{proof}

\section{The table substitution on \texorpdfstring{$\ZM^2$}{ZM2}}
\label{Chap7-Sect-TableSubst}

The table tiling in $\RM^2$ was first introduced by {\sc Solomyak} \cite{Sol98}. The tiling arises from an inflation map and a finite number of prototiles. According to \cite{Fra08}, this tiling of $\RM^2$ can be represented by a block substitutions on $\ZM^2$ over a finite alphabet, c.f. Example~\ref{Chap6-Ex-TableSubst}. The table tiling was studied by several authors, c.f. \cite{Sol98,Rob99,Fra05,Fra08,BaGaeGr13,BaakeGrimm13}.

\medskip

\begin{figure}[htb]
\centering
\includegraphics[scale=2.5]{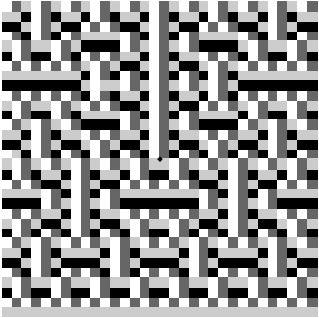}
\caption{The structure of the $2$-periodic element $\xi\in\Xi_S$ with respect to the table substitution $S$. A white box represents the letter $a$, the light gray box represents the letter $b$, the dark gray box represents the letter $c$ and the black box represents the letter $d$. The little black dot in the middle fixes the origin.}
\label{Chap7-Fig-TableTiling}
\end{figure}

Recall the substitution rule for the table tiling, c.f. Example~\ref{Chap6-Ex-TableSubst}. Let $\as:=\{a,b,c,d\}$ be the alphabet with four letters. The {\em table substitution} $S:\as\to\textit{Pat}_{\ZM^2}(\as)$ is defined by
$$
a\; \overset{S}{\longmapsto} \;
	\begin{matrix}
		b & a\\
		d & a
	\end{matrix}\;,
	\qquad
b\; \overset{S}{\longmapsto} \;
	\begin{matrix}
		a & c\\
		b & b
	\end{matrix}\;,
	\qquad
c\; \overset{S}{\longmapsto} \;
	\begin{matrix}
		c & b\\
		c & d
	\end{matrix}\;,
	\qquad
d\; \overset{S}{\longmapsto} \;
	\begin{matrix}
		d & d\\
		a & c
	\end{matrix}\; .
$$
This is a primitive block substitution with parameters $n_1=n_2=2$ for the block and $l_0=2$ for the primitivity. The subshift $\Xi_S$ associated with $S$ is defined via its dictionary $\ws_S$ or by a $k$-periodic element $\xi\in\as^{\ZM^2}$ with respect to the substitution $S$, c.f. Proposition~\ref{Chap6-Prop-AssDictSubshSubstZd} and Proposition~\ref{Chap6-Prop-PrimSubstSubshZdMin}. In detail, if $S$ is applied to the letter $b$, this leads to
\begin{equation}
\label{Chap7-Eq-TableSubstS4(b)} \tag{$\spadesuit$}
b\;\overset{S}{\longmapsto}\;
	\begin{array}{c||c}
		a & c \\
		\hline\hline
		b & b
	\end{array}
	\;\overset{S}{\longmapsto}\;
		\begin{array}{cc||cc}
			b & a & 			c & b\\
			d & a & 			c & d\\
			\hline			\hline
			a & c & 			a & c\\
			b & b & 			b & b		
		\end{array}
			\;\overset{S}{\longmapsto}\;
				\begin{array}{cc|cc||cc|cc}
					a & c & b & a &		c & b & a & c\\
					b & b & d & a &		c & d & b & b\\
					\hline
					d & d & b & a &		c & b & d & d\\
					a & c & d & a &		c & d & a & c\\
					\hline				\hline
					b & a & c & b &		b & a & c & b\\
					d & a & c & d &		d & a & c & d\\
					\hline
					a & c & a & c &		a & c & a & c\\
					b & b & b & b &		b & b & b & b\\
				\end{array}
				\;,
\end{equation}
where the double lines fix the origin. This immediately leads to the existence of two $2$-periodic elements with respect to $S$ as follows. Consider the patterns
$$
u \;
	:= \; \begin{array}{c||c}
		a & c\\
		\hline \hline
		b & b
	\end{array}\;,\qquad
v \;
	:= \; \begin{array}{c||c}
		a & c\\
		\hline \hline
		c & a
	\end{array}\;.
$$
Then, let $\eta\in\Xi_S$ be such that $\eta|_K=u$ for $K=\{-1,0\}\times\{-1,0\}\subseteq\ZM^2$, i.e.,
$$
\eta(-1,0) \; 
	= \; a\,,\qquad
\eta(0,0) \;
	= \; c\,,\qquad
\eta(0,-1) \;
	= \; b\,,\qquad
\eta(-1,-1) \;
	= \; b\,.
$$
The existence of such an $\eta$ is guaranteed by Lemma~\ref{Chap2-Lem-ExInfWord}. Then $S^2(\eta)$ fulfills $S^2(\eta)|_K=u$. Thus, the sequence $\big(S^{2n}(\eta)\big)_{n\in\NM}$ is convergent in $\as^{\ZM^2}$ by using Lemma~\ref{Chap2-Lem-XiDictClosed} since $S^n(\eta)|_K=u\in\ws_S$ and $S$ is a block substitution. Denote the limit point by $\xi:=\lim_{n\to\infty} S^{2n}(\eta)$. Due to Lemma~\ref{Chap6-Lem-SubstitutionContinuous}, $\xi\in\Xi_S$ defines a $2$-periodic element with respect to the substitution $S$. The pattern around the origin of the $2$-periodic element $\xi$ with respect to $S$ is given in Figure~\ref{Chap7-Fig-TableTiling}. Similarly, another fixed point is defined by replacing the pattern $u$ by the pattern $v$.

\medskip

The table subshift is minimal and aperiodic, c.f. Proposition~\ref{Chap6-Prop-BasPropTableSubst}. It turns out that the table subshift is periodically approximable by using Proposition~\ref{Chap6-Prop-SuffPerApprZ2SymmPatt}. More precisely, the existence of a pattern $v\in\ws_S$ is proven such that $v$ satisfies the symmetry described in Remark~\ref{Chap6-Rem-SuffPerApprZ2SymmPatt}~(iii).

\begin{proposition}[Table subshift periodically approximable]
\label{Chap7-Prop-TableTilingPerAppr}
The table subshift $\Xi_S\in\mathcal{I}_{\ZM^2}\big(\as^{\ZM^2}\big)$ arising by the table substitution $S:\as\to\textit{Pat}_{\ZM^2}(\as)$ is periodically approximable. In particular, let $\eta^1:=v_1^\infty\in\as^{\ZM^2}$ and $\eta^2:=v_2^\infty\in\as^{\ZM^2}$ be defined by
$$
v_1 \; 
	:= \;\begin{matrix}
			b & d\\
			d & b
		\end{matrix}\,,
\qquad
v_2 \; 
	:= \;\begin{matrix}
			a & c\\
			c & a
		\end{matrix}\,.
$$
For $j\in\{1,2\}$, the sequence of strongly periodic elements $\big( S^n(\eta^j)\big)_{n\in\NM}$ defines a sequence of strongly periodic subshifts $\Xi_n:=\Orb\big(S^n(\eta_j)\big)\,,\; n\in\NM\,,$ that converges to $\Xi_S$.
\end{proposition}

\begin{proof}
Let $j\in\{1,2\}$. The element $\eta^j:=v_j^\infty\in\as^{\ZM^2}$ is strongly periodic by Lemma~\ref{Chap6-Lem-ReprStrongPer}. Lemma~\ref{Chap6-Lem-SnETAStronglyPeriodic} implies that $S^n(\eta^j)\,,\; n\in\NM\,,$ are also strongly periodic elements. Due to Theorem~\ref{Chap6-Theo-PerApprZdPrimSubst}, it suffices to show that $\eta^j:=v_j^\infty$ satisfies $\ws(\eta_j)\cap\as^{[K]}\subseteq\ws_S\cap\as^{[K]}$ for $K:=\{0,1\}\times\{0,1\}$. Consider the sets
$$
U_1 \; 
	:= \; \left\{
		\begin{matrix}
			b & d\\
			d & b
		\end{matrix}\;,
		\qquad
		\begin{matrix}
			d & b\\
			b & d
		\end{matrix}
	\right\}\;,
\qquad
U_2 \; 
	:= \; \left\{
		\begin{matrix}
			a & c\\
			c & a
		\end{matrix}\;,
		\qquad
		\begin{matrix}
			c & a\\
			a & c
		\end{matrix}
	\right\}\,.
$$
The set $U_j$ is exactly the set defined in Proposition~\ref{Chap6-Prop-SuffPerApprZ2SymmPatt}~(ii) associated with the pattern $v_j$. 

\vspace{.1cm}

Both patterns in $U_1$ and the pattern $v_2$ occur in $S^4(b)$ by (\ref{Chap7-Eq-TableSubstS4(b)}). In addition, the pattern $u:=
\begin{matrix}
	d & d\\
	a & c
\end{matrix}$ is contained in $S^4(b)$ and $
\begin{matrix}
	c & a\\
	a & c
\end{matrix}$ appears in $S(u)$. Altogether, the inclusions $U_1, U_2\subseteq\ws_S$ are derived. Consequently, Proposition~\ref{Chap6-Prop-SuffPerApprZ2SymmPatt} and Theorem~\ref{Chap6-Theo-PerApprZdPrimSubst} lead to the desired result.
\end{proof}

\begin{proposition}[Periodic approximation of the table operator]
\label{Chap7-Prop-SchrOpTableTilingPerAppr}
Let $S:\as\to\textit{Pat}_{\ZM^2}(\as)$ be the table substitution with associated subshift $\Xi_S$. Consider a pattern equivariant Schr\"odinger operator $H_\xi:\ell^2\big(\ZM^2\big)\to\ell^2\big(\ZM^2\big)$, $\xi\in\as^{\ZM^2}$, defined by
$$
\begin{array}{rl}
(H_\xi\psi)(g) \; 
	:= \; &\left( 
			\sum\limits_{h\in K} 
				p_h\big(\alpha_{-g}(\xi)\big) \cdot \psi(g-h) + 
				\overline{p_h\big(\alpha_{-(g+h)}(\xi)\big)} \cdot \psi(g+h)
		\right)\\[0.2cm]
	&+\; p_e\big(\alpha_{-g}(\xi)\big) \cdot \psi(g)
\end{array}
$$
for $g\in\ZM^2$ and $\psi\in\ell^2\big(\ZM^2\big)$ with corresponding finite set $K\subseteq\ZM^2$ and pattern equivariant functions $p_h:\as^{\ZM^2}\to\CM\,,\; h\in K\,,$ and $p_e:\as^{\ZM^2}\to\RM$ for the neutral element $e:=(0,0)$ of $\ZM^2$. Define the strongly periodic elements $\eta^1_n:=S^n\big(v_1^\infty\big)\,,\; n\in\NM\,,$ and $\eta^2_n:=S^n\big(v_2^\infty\big)\,,\; n\in\NM\,,$ for the patterns
$$
v_1 \; 
	:= \;\begin{matrix}
			b & d\\
			d & b
		\end{matrix}\,,
\qquad
v_2 \; 
	:= \;\begin{matrix}
			a & c\\
			c & a
		\end{matrix}\in\ws_S\,.
$$
Then the equalities 
$$
\lim_{n\to\infty} \sigma\big( H_{\eta^j_n} \big) \;
	= \; \sigma\big( H_\xi \big) \;
	= \; \sigma_{ess}\big(H_\xi\big) \;
	= \; \sigma\big(H_{\Xi_S}\big)\,,
	\qquad
	\xi\in\Xi_S\,,\; j\in\{1,2\}\,,
$$
hold where the limit is taken with respect to the Hausdorff metric on $\ks(\RM)$.
\end{proposition}

\begin{proof}
The pattern equivariant Schr\"odinger operator is self-adjoint and so the corresponding spectrum is contained in $\RM$, c.f. Lemma~\ref{Chap2-Lem-PatEqSchrOp}. The subshifts $\Xi_n:=\Orb(\eta_n)\,,\; n\in\NM\,,$ and $\Xi_S$ are minimal. Hence, the equations $\sigma\big(H_{\eta_n}\big)=\sigma\big(H_{\Xi_n}\big)$ and $\sigma\big(H_{\xi}\big)=\sigma_{ess}\big(H_{\xi}\big)=\sigma\big(H_{\Xi_S}\big)$ are derived for $\xi\in\Xi_S$ by Theorem~\ref{Chap2-Theo-ConstSpectrMinimal}. According to Proposition~\ref{Chap7-Prop-TableTilingPerAppr}, the sequence $\big(\Xi_n\big)_{n\in\NM}$ of subshifts converges to $\Xi_S$ in the Hausdorff-topology of $\mathcal{I}_{\ZM^2}\big(\as^{\ZM^2}\big)$. Thus, Theorem~\ref{Chap4-Theo-CharSubshiftConvSpectr} leads to the desired result.
\end{proof}

\section{The Sierpinski carpet substitution on \texorpdfstring{$\ZM^2$}{ZM2}}
\label{Chap7-Sect-SierpinskiCarpetSubstitution}

The Sierpinski carpet goes back to {\sc Sierpi{\'{n}}ski} \cite{Sie16} and it is intensively studied in fractal geometry. In \cite[Example~3]{Fra08}, a substitutional subshift over the group $\ZM^2$ and an alphabet $\as$ with two letters is defined in view of the Sierpinski carpet construction. This substitution is analyzed in this section for the existence of strongly periodic approxi\-mations. Since the Sierpinski carpet substitution is not primitive, Theorem~\ref{Chap6-Theo-PerApprZdPrimSubst} cannot be applied directly. However, a similar strategy like in Section~\ref{Chap6-Sect-PerApprPrimBlockSubst} can be followed here. Thus, this section provides also a guideline to study other substitutional systems that do not arise by primitive substitutions.

\begin{definition}[Sierpinski carpet substitution]
\label{Chap7-Def-SierpinskiCarpetSubstitution}
Consider the alphabet $\as:=\{a,b\}$ with two letters. The block substitution $S:\as\to\textit{Pat}_{\ZM^2}(\as)$ defined by
$$
a\; \overset{S}{\longmapsto} \;
	\begin{matrix}
		a & a & a\\
		a & a & a\\
		a & a & a
	\end{matrix}\;,
	\qquad
b\; \overset{S}{\longmapsto} \;
	\begin{matrix}
		b & b & b\\
		b & a & b\\
		b & b & b
	\end{matrix}\;,
$$
is called {\em Sierpinski carpet substitution}. 
\end{definition}

Since only the letter $a$ occurs in $S^k(a)$ for each $k\in\NM$, the substitution is not primitive. Nevertheless, the subshift $\Xi_S$ and the dictionary $\ws_S$ associated with $S$ is well-defined. More precisely, they satisfy the requirements of Definition~\ref{Chap2-Def-SubdynSyst} and Definition~\ref{Chap2-Def-Dictionary}, respectively.

\begin{proposition}
\label{Chap7-Prop-SierpinskiCarpetDictSubs}
Let $S:\as\to\textit{Pat}_{\ZM^2}(\as)$ be the Sierpinski carpet substitution. Then the associated dictionary $\ws_S$ with the Sierpinski carpet substitution defined in Definition~\ref{Chap6-Def-AssDictSubshSubstZd} is a dictionary, i.e., $\ws_S$ satisfies \nameref{(D1)} and \nameref{(D2)}. Furthermore, the associated subshift $\Xi_S:=\big\{\xi\in\as^{\ZM^2}\;|\; \ws(\xi)\subseteq\ws_S \big\}\subseteq\as^{\ZM^2}$ with the Sierpinski carpet substitution is a subshift, i.e., $\Xi_S$ is closed and $\ZM^2$-invariant.
\end{proposition}

\begin{proof}
Recall that $\ws_S$ is defined by all $v\in\textit{Pat}_{\ZM^2}(\as)$ such that there is a $k_v\in\NM$ and a letter $a_v\in\as$ satisfying $v$ is a subpattern of $S^k(a_v)$. Thus, it is clear that $\ws_S$ satisfies \nameref{(D1)}. Condition \nameref{(D2)} is verified as follows. Consider the pattern $S^2(b)\in\as^{[K]}$ with $K:=\{0,1,\ldots, 8\}\times\{0,1,\ldots,8\}$. Then the letter $a$ and $b$ occur in the interior of the pattern, i.e., there exist an $i^a,i^b\in\ZM^2$ such that $0<i_j^a,i_j^b<8$ for $j\in\{1,2\}$. With this at hand and the fact that the Sierpinski carpet substitution is a block substitution, the proof that $\ws_S$ fulfills \nameref{(D2)} follows exactly the same lines of Lemma~\ref{Chap6-Lem-AssDictInvSubstit}. Hence, $\ws_S$ is a dictionary. Then Theorem~\ref{Chap2-Theo-Shift+DictSpace} implies that $\Xi_S$ is a subshift, i.e., it is closed and $\ZM^2$-invariant.
\end{proof}

\medskip

The subshift $\Xi_S$ induced by the Sierpinski carpet substitution is called {\em Sierpinski carpet subshift}. Similar to the primitive case, the Sierpinski carpet subshift is topologically transitive, i.e., it has a $1$-periodic element with respect to $S$ that generates the whole subshift $\Xi_S$.

\begin{proposition}
\label{Chap7-Prop-SierpinskiCarpetSubshift}
The Sierpinski carpet subshift $\Xi_S$ is not minimal and contains a periodic element $\eta\in\Xi_S$. Furthermore, there exists a non-periodic $\xi\in\Xi_S$ satisfying $S(\xi)=\xi$ and $\overline{\Orb(\xi)}=\Xi_S$, i.e., $\Xi_S$ is aperiodic and topologically transitive.
\end{proposition}

\begin{figure}[htb]
\centering
\includegraphics[scale=0.5]{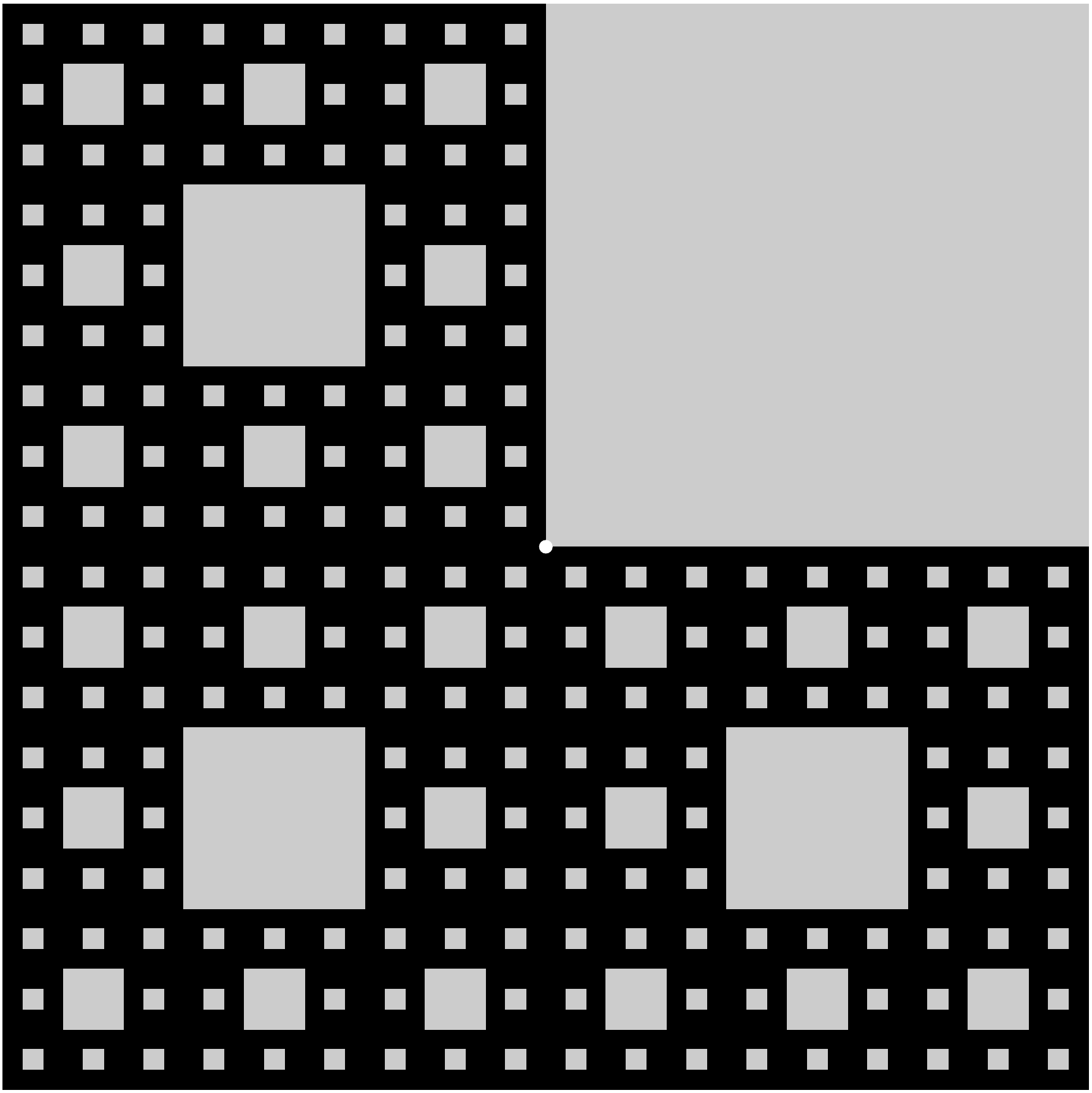}
\caption{The structure of the $1$-periodic element $\xi\in\Xi_S$ with respect to the Sierpinski carpet substitution $S$ defined in the proof of Proposition~\ref{Chap7-Prop-SierpinskiCarpetSubshift}. A light gray box represents the letter $a$ and the black box represents the letter $b$. The white dot in the middle fixes the origin.}
\label{Chap7-Fig-SierpinskiCarpet}
\end{figure}

\begin{proof}
The strongly periodic element $\eta\in\as^{\ZM^d}$ defined by $\eta(j)=a\,,\; j\in\ZM^2\,,$ is contained in $\Xi_S$ since arbitrary large patterns are contained in $\ws_S$ that only consists of the letter $a$. Then $\Orb(\eta)=\{\eta\}\subseteq\Xi_S$ follows. On the other hand, the dictionary $\ws_S$ contains patterns where the letter $b$ occurs. Then Lemma~\ref{Chap2-Lem-ExInfWord} implies that there exists a $\xi\in\Xi_S$ satisfying that it contains the letter $b$. Thus, $\Xi_S$ is not minimal by the previous considerations.

\vspace{.1cm}

It is left to show that there exists a non-periodic element $\xi$ in $\Xi_S$ such that $\Xi_S=\overline{\Orb(\xi)}$. Consider an $\omega\in\Xi_S\subseteq\as^{\ZM^2}$ satisfying
$$
\omega|_K \;
	= \; \begin{array}{c||c}
		b & a\\
		\hline \hline
		b & b
	\end{array}\,,\;
	\qquad
	K:=\{-1,0\}\times\{-1,0\}\subseteq\ZM^2\,,
$$ 
where the double lines fix the origin. The existence of $\omega$ is guaranteed by Lemma~\ref{Chap2-Lem-ExInfWord} since $\ws_S$ defines a dictionary. The pattern $[\omega|_K]$ is an element of $\ws_S$ as it occurs in $S(b)$. The application of the substitution $S$ to $\omega$ leads to 
$$
\begin{array}{c||c}
	b & a \\
	\hline\hline
	b & b
\end{array}
	\;\overset{S}{\longmapsto}\;
		\begin{array}{ccc||ccc}
			b & b & \textcolor{gray}{b} & 			a & a & a\\
			b & a & 	\textcolor{gray}{b} &			a & a & a\\
			b & b & \textcolor{gray}{b} & 			a & a & a\\
			\hline				\hline
			\textcolor{gray}{b} & \textcolor{gray}{b} & \textcolor{gray}{b} &
			 			\textcolor{gray}{b} & \textcolor{gray}{b} & \textcolor{gray}{b}\\
			b & a & \textcolor{gray}{b} & 			b & a & b\\
			b & b & \textcolor{gray}{b} & 			b & b & b\\
		\end{array}\;.
$$
Then $\big(S^n(\omega)\big)_{n\in\NM}$ defines a convergent sequence due to definition of the Sierpinski carpet substitution by using Lemma~\ref{Chap2-Lem-XiDictClosed}. Its limit point $\xi:=\lim_{n\to\infty} S^n(\omega)$ satisfies $\ws(\xi)\subseteq\ws_S$ and
$$
\xi(i_1,i_2) \;
	= \; \begin{cases}
		a\,,\quad \text{if } i_1,i_2\geq 0\,,\\
		b\,,\quad \text{if } i_1= -1 \text{ or } i_2=-1\,,
	\end{cases}
	\qquad
	(i_1,i_2)\in\ZM^2\,,
$$
c.f. Figure~\ref{Chap7-Fig-SierpinskiCarpet}. Specifically, the axes are labeled with the letter $b$ and the first quadrant of the plane is labeled with the letter $a$. Thus, each translation $\alpha_j(\xi)$ by a $j\in\ZM^2$ is not equal to $\xi$ besides $j=(0,0)$. Hence, $\xi$ is not periodic. Due to Lemma~\ref{Chap6-Lem-SubstitutionContinuous} and the definition of $\xi$, $S(\xi)=\xi$ is deduced, i.e., $\xi$ is a $1$-periodic element with respect to the Sierpinski carpet substitution. Since the limit point $\xi$ satisfies $S(\xi)=\xi$ and $\xi$ contains both letters $a$ and $b$, the equality $\ws(\xi)=\ws_S$ is derived. Consequently, the equality $\overline{\Orb(\xi)}=\Xi_S$ holds by Theorem~\ref{Chap2-Theo-Shift+DictSpace}.
\end{proof}

\medskip

According to Proposition~\ref{Chap7-Prop-SierpinskiCarpetSubshift}, the Sierpinski carpet subshift is not minimal but defined via a $1$-periodic element with respect to the substitution $S$. It turns out that the same strategy \nameref{(6.I)}-\nameref{(6.IV)} of Chapter~\ref{Chap6-HigherDimPerAppr} can be followed to show that the Sierpinski carpet subshift is periodically approximable. More precisely, the symmetry described in Proposition~\ref{Chap6-Prop-SuffPerApprZ2SymmPatt} and Remark~\ref{Chap6-Rem-SuffPerApprZ2SymmPatt}~(iii) is used to define a strongly periodic element satisfying~\nameref{(6.I)}.

\begin{proposition}[Sierpinski carpet subshift periodically approximable]
\label{Chap7-Prop-SierpinskiCarpetPerAppr}
The subshift $\Xi_S$ over $\as^{\ZM^2}$ arising by the Sierpinski carpet substitution is periodically approximable. In particular, let $\eta:=v^\infty$ be a strongly periodic element defined by the pattern 
$$
v \; 
	:= \;\begin{matrix}
			b & a\\
			b & b
		\end{matrix}\in\ws_S\,.
$$
Then the sequence of strongly periodic subshifts $\Xi_n:=\Orb\big(S^n(\eta)\big)\,,\; n\in\NM\,,$ converges to $\Xi_S$ in the Hausdorff-topology of $\textit{Sub}_{\ZM^2}\big(\as^{\ZM^2}\big)$.
\end{proposition}

\begin{proof}
The proof is along the lines of Proposition~\ref{Chap7-Prop-TableTilingPerAppr}. The pattern $v$ satisfies the condition of Proposition~\ref{Chap6-Prop-SuffPerApprZ2SymmPatt}~(ii), i.e., the set
$$
U \;
	:=\; \left\{
		\begin{matrix}
			b & b\\
			b & a
		\end{matrix}\;,
		\qquad
		\begin{matrix}
			b & b\\
			a & b
		\end{matrix}\;,
		\qquad
		\begin{matrix}
			a & b\\
			b & b
		\end{matrix}\;,
		\qquad
		\begin{matrix}
			b & a\\
			b & b
		\end{matrix}
	\right\}
$$
is contained in $\ws_S$ since the elements of $U$ are subpatterns of $S(b)$, c.f. Definition~\ref{Chap7-Def-SierpinskiCarpetSubstitution}. Hence, $\eta:=v^\infty$ satisfies $\ws(\eta)\cap\as^{[K]}\subseteq\ws_S\cap\as^{[K]}$ for $K:=\{0,1\}\times\{0,1\}$ by Proposition~\ref{Chap6-Prop-SuffPerApprZ2SymmPatt}. With this at hand, the desired result is concluded by following the lines of Theorem~\ref{Chap6-Theo-PerApprZdPrimSubst} as follows.

\vspace{.1cm}

For $n\in\NM$, $S^n(\eta)$ is strongly periodic by Lemma~\ref{Chap6-Lem-ReprStrongPer} and Lemma~\ref{Chap6-Lem-SnETAStronglyPeriodic}. Let $(K_m)_{m\in\NM}$ be the exhausting sequence of $\ZM^2$ defined by $K_m:=\prod_{j=1}^d\{-2^m+1,\ldots,2^m\}\,,\; m\in\NM$. Lemma~\ref{Chap6-Lem-GrowthAcceptPatternSubst} leads to
$$
\ws\big(S^n(\eta)\big)\cap\as^{[K_m]} \;
	\subseteq \; \ws_S\cap\as^{[K_m]}\,,
	\qquad	
	m\in\NM\,,\; n\geq m\,.
$$
Let $m\in\NM$ and $[w]\in\ws_S\cap\as^{[K_m]}$ with representative $w:K_M\to\as$. According to the definition of $\ws_S$, there exist an $l_w\in\NM$ and an $x_w\in\as$ such that $w$ appears in $S^{l_w}(x_w)$. Then $w$ is a subpattern of $S^{l_w}(v)$ since both letters $a$ and $b$ occur in the pattern $v$. Thus, $[w]\in\ws\big(S^{l_w}(\eta)\big)$ follows as $S^{l_w}(\eta)=\big(S^{l_w}(v)\big)^\infty$, c.f. Lemma~\ref{Chap6-Lem-SnETAStronglyPeriodic}. Due to the fact that $S(b)$ contains both letters $a$ and $b$, the pattern $S^l(v)$ contains also both letters $a$ and $b$ for all $l\in\NM$. Consequently, the same argument as before implies that $[w]$ is an element of $\ws\big(S^m(\eta)\big)$ for all $m\geq l_w$. The maximum $n(m):=\max\big\{ m,\max\big\{ l_w\;|\; [w]\in\ws_S\cap\as^{[K_m]} \big\} \big\}$ exists as $\as^{[K_m]}$ is finite. Then the previous considerations yield
$$
\ws\big(S^n(\eta)\big)\cap\as^{[K_m]} \;
	= \; \ws_S \cap\as^{[K_m]}\,,
	\qquad m\in\NM\,,\;  n\geq n(m)\,.
$$
Since $(K_m)_{m\in\NM}$ is an exhausting sequence of $\ZM^2$, the sequence $\big(\ws\big(S^n(\eta)\big)\big)_{n\in\NM}$ of diction\-aries converges to $\ws_S$ in the local pattern topology, c.f. Corollary~\ref{Chap2-Cor-DictExhSeq}. Hence, Theo\-rem~\ref{Chap2-Theo-Shift+DictSpace} leads to the convergence of the strongly periodic subshifts $\big(\Xi_n\big)_{n\in\NM}$ to $\Xi_S$. Consequently, $\Xi_S$ is periodically approximable.
\end{proof}

\begin{remark}
\label{Chap7-Rem-SierpinskiCarpetPerAppr}
Note that the second part of the proof of Proposition~\ref{Chap7-Prop-SierpinskiCarpetPerAppr} follows the lines of Theorem~\ref{Chap6-Theo-PerApprZdPrimSubst}. This is possible as all letters occur in the pattern $v$. However, Theorem~\ref{Chap6-Theo-PerApprZdPrimSubst} cannot be applied directly since the Sierpinski carpet substitution is not primitive.
\end{remark}

With this at hand, periodic approximations for pattern equivariant Schr\"odinger operators associated with the Sierpinski carpet subshift exist.

\begin{proposition}[Periodic approximation of the Sierpinski carpet operator]
\label{Chap7-Prop-SchrOpSierpinskiCarpetPerAppr}
Let $S:\as\to\textit{Pat}_{\ZM^2}(\as)$ be the Sierpinski carpet substitution with associated subshift $\Xi_S$. Consider a pattern equivariant Schr\"odinger operator $H_\xi:\ell^2\big(\ZM^2\big)\to\ell^2\big(\ZM^2\big)$, $\xi\in\as^{\ZM^2}$, defined by
\begin{align*}
(H_\xi\psi)(g) \; 
	:= \; &\left( 
			\sum\limits_{h\in K} 
				p_h\big(\alpha_{-g}(\xi)\big) \cdot \psi(g-h) + 
				\overline{p_h\big(\alpha_{-(g+h)}(\xi)\big)} \cdot \psi(g+h)
		\right)\\[0.2cm]
	&+ p_e\big(\alpha_{-g}(\xi)\big) \cdot \psi(g)
\end{align*}
for $g\in\ZM^2$ and $\psi\in\ell^2\big(\ZM^2\big)$ with corresponding finite set $K\subseteq\ZM^2$ and pattern equivariant functions $p_h:\as^{\ZM^2}\to\CM\,,\; h\in K\,,$ and $p_e:\as^{\ZM^2}\to\RM$ for the neutral element $e:=(0,0)$ of $\ZM^2$. Define the strongly periodic element $\eta_n:=S^n\big(v^\infty\big)$ for the pattern 
$$
v \; 
	:= \;\begin{matrix}
			b & a\\
			b & b
		\end{matrix}\in\ws_S\,.
$$
Then the equality
$$
\lim_{n\to\infty} \sigma\big( H_{\eta_n} \big) \;
	= \; \sigma\big(H_{\Xi_S}\big)
$$
holds where the limit is taken with respect to the Hausdorff metric on $\ks(\RM)$.
\end{proposition}

\begin{proof}
The pattern equivariant Schr\"odinger operator is self-adjoint and so the spectrum is contained in $\RM$, c.f. Lemma~\ref{Chap2-Lem-PatEqSchrOp}. According to Proposition~\ref{Chap7-Prop-SierpinskiCarpetPerAppr}, the subshifts $\Xi_n:=\Orb\big(S^n(\eta)\big)\,,\; n\in\NM\,,$ converge to $\Xi_S$ in the Hausdorff-topology of $\mathcal{I}_{\ZM^2}\big(\as^{\ZM^2}\big)$. In addition, the equation $\sigma\big( H_{\eta_n}\big)=\sigma\big( H_{\Xi_n}\big)$ holds for $n\in\NM$ since $\Xi_n$ is minimal, c.f.  Theorem~\ref{Chap2-Theo-MinCharConstSpectr}. Altogether, Theorem~\ref{Chap4-Theo-CharSubshiftConvSpectr} leads to the desired result.
\end{proof}

\cleardoublepage


\chapter{Discussion}
\label{Chap8-Disc}

In this chapter, open questions and future projects based on this work are discussed.

\section{Groupoids}
\label{Chap8-Sect-Groupoids}

Throughout this work it was assumed that all the considered groupoids are \'etale. The \'etale property is equivalent to the property that the $r$-fibers of the groupoid are discrete if the groupoid has a Haar system, c.f. \cite[Proposition~I.2.8]{Renault80} or Lemma~\ref{Chap2-Lem-CharUnitSpaceOpen} and Lemma~\ref{Chap2-Lem-FibreDiscrete}. According to Theorem~\ref{Chap2-Theo-UnitGroupoidCalgebra}, the reduced and the full $C^\ast$-algebra have a unit if the groupoid is \'etale and has compact unit space. As can be seen from Remark~\ref{Chap3-Rem-NecessIdentityContSpect}, it is necessary that the groupoid $C^\ast$-algebras have a unit for the application to the continuity of the spectra. Consequently, it is natural to ask if the assertion of Theorem~\ref{Chap2-Theo-UnitGroupoidCalgebra} is optimal, i.e., is the \'etale property and the compactness of the unit space necessary to insure that the associated reduced/full groupoid $C^\ast$-algebra has a unit? Or stated differently:

\begin{question}
\label{Chap8-Que-UnitSpace}
Is the unit space of $\Gamma$ open and compact if the reduced $C^\ast$-algebra $\CG^\ast_{red}(\Gamma)$ is unital?
\end{question}

From the philosophical point of view, the answer is yes. This was hinted in \cite{Exe14} without a proof and it was suggest to me in personal communication by {\sc J. Renault}. Furthermore, {\sc M. Landstad} \cite{Lan16} proved the assertion under the additional assumption that $\Gamma:=X\rtimes_\alpha G$ is a transformation group groupoid by using that the reduced group $C^\ast$-algebra $\CG^\ast_{red}(G)$ is unital if and only if the group $G$ is discrete. If the answer of Question~\ref{Chap8-Que-UnitSpace} is yes, the equivalence of the following assertions (i)-(iv) is immediately derived for a topological groupoid $\Gamma$ with left-continuous Haar system $\mu$.
\begin{itemize}
\item[(i)] The full $C^\ast$-algebra $\CG^\ast_{full}(\Gamma)$ is unital.
\item[(ii)] The reduced $C^\ast$-algebra $\CG^\ast_{red}(\Gamma)$ is unital.
\item[(iii)] The unit space $\Gamma^{(0)}$ is open and compact.
\item[(iv)] The groupoid $\Gamma$ is \'etale and has compact unit space $\Gamma^{(0)}$.
\end{itemize}
The implications (iv)$\Rightarrow$(i)$\Rightarrow$(ii) are proven in Theorem~\ref{Chap2-Theo-UnitGroupoidCalgebra} which are well-known in the literature. The equivalence of (iii) and (iv) is stated as an exercise in \cite{Renault09} and is proven in \cite{Tho10}, c.f. Lemma~\ref{Chap2-Lem-CharUnitSpaceOpen} as well. Meanwhile, the implication (ii)$\Rightarrow$(iii) is stated in Question~\ref{Chap8-Que-UnitSpace}.

\medskip

In order to show the assertion of Question~\ref{Chap8-Que-UnitSpace}, the connection between the unit $\mathpzc{1}\in\CG^\ast_{red}(\Gamma)$ and the unit space $\Gamma^{(0)}$ is the support of the unit $\mathpzc{1}$. More precisely, if the unit $\mathpzc{1}$ is represented as a continuous function $\mathpzc{1}:\Gamma^{(1)}\to\CM$, then it follows that the support $\supp(\mathpzc{1})$ is equal to $\Gamma^{(0)}$ and $\mathpzc{1}=\chi_{\Gamma^{(0)}}$. This implies that $\Gamma^{(0)}$ is compact and open. Altogether, the main task is to verify that the unit $\mathpzc{1}$ is represented by a continuous function on the groupoid $\Gamma$.

\medskip

Note that this work mainly focuses on groupoids arising by dynamical systems. Clearly, Theorem~\ref{Chap4-Theo-FieldGroupContCAlgebras} applies to more general groupoids than transformation group groupoids. The next Section~\ref{Chap8-Sect-DeloneTiling} deals with topological groupoids that do not arise from a dynamical system while they are interesting from the point of view of mathematical physics.

\section{Delone sets and tilings}
\label{Chap8-Sect-DeloneTiling}

The coloring of a group by an alphabet, i.e., a symbolic dynamical system $(\as^G,G,\alpha)$, models a specific case of solids in mathematical physics. More precisely, suppose that the group $G$ is a discrete subgroup of a second countable, locally compact, Hausdorff, (abelian) group $\mathfrak{X}$. The typical example for such a space is $\mathfrak{X}=\RM^d$. Then $G$ forms a lattice in $\mathfrak{X}$ and it represents the positions of the atoms in the space $\mathfrak{X}$ of the solid. In general, a solid is just described by a closed subset $\mathcal{D}$ of $\mathfrak{X}$ such that there is a minimal distance between two elements of $\mathcal{D}$ and no large gaps, c.f. also discussion in Chapter~\ref{Chap1-Intro}. 

\medskip

More precisely, a set $\mathcal{D}\subseteq\mathfrak{X}$ is called a {\em Delone set} if there exist an open $U\subseteq\mathfrak{X}$ and a compact $K\subseteq\mathfrak{X}$ such that $\sharp (y+U)\cap\mathcal{D}\leq 1$ for every $y\in\mathfrak{X}$ and $\bigcup_{y\in\mathcal{D}}y+K=\mathfrak{X}$. Delone sets are intensively studied in the literature where $\mathfrak{X}$ is an abelian locally compact, Hausdorff group, c.f. \cite{Moo97,Lag99,BeHeZa00,Sch00,LaPl03,BaLe04,BaLeMo07}. The set of Delone sets $Del_{U,K}$ with fixed parameters $U,K\subseteq\mathfrak{X}$ is naturally equipped with a topology that can be described in different ways. Either the topology is defined by the vague topology on the set of Dirac combs $\delta_{\mathcal{D}}$ where $\delta_{\mathcal{D}}:=\sum_{x\in\mathcal{D}}\delta_x$ gives constant weight to each point or it is described by a base for the topology which is called the {\em local rubber topology}, see e.g. \cite{BeHeZa00,BaLe04}. The local rubber topology is exactly the Fell-topology on $\cs(\mathfrak{X})$ which is also called Chabauty-Fell-topology. By Theorem~\ref{Chap3-Theo-VietFellHausMetricEquiv}, this topology agrees with the topology induced by the Hausdorff metric restricted to open bounded sets, see also \cite{BeHeZa00,LeSt03-Delone,BaLe04} for the case $\mathfrak{X}=\RM^d$. Whenever the Delone sets are of finite local complexity, another characterization of the convergence for Delone sets is provided in \cite[Proposition~2.8]{Bec12}. The space $Del_{U,K}$ equipped with the Fell-topology turns out to be compact, c.f. \cite[Theorem~1.5]{BeHeZa00}. Recent developments \cite{Yok05,BjHaPo16} study also Delone sets in non-commutative spaces $\mathfrak{X}$. 

\medskip

A Delone set naturally induces a groupoid defined by the translation by $\mathfrak{X}$. More precisely, the closure $\Xi_{\mathcal{D}}:=\overline{\{\mathcal{D}-x\;|\; x\in\mathcal{D} \}}$ with respect to the Fell-topology on $\cs(\mathfrak{X})$ is the unit space. This set $\Xi_{\mathcal{D}}$ is called the {\em transversal} of the Delone set $\mathcal{D}$. Then the set of arrows is defined by the pairs $(\tilde{\mathcal{D}},x)\in\Xi_{\mathcal{D}}\times\mathfrak{X}$ such that the translated set $\tilde{\mathcal{D}}-x$ contains the origin, i.e., $x\in\tilde{\mathcal{D}}$. Furthermore, the composition is defined by $(\tilde{\mathcal{D}},x)\circ(\tilde{\mathcal{D}}-x,y):=(\tilde{\mathcal{D}},x+y)$ and the inversion is given by $(\tilde{\mathcal{D}},x)^{-1}=(\tilde{\mathcal{D}}-x,-x)$. Since only those Delone sets are considered that contain the origin, the groupoid is not a transformation group groupoid even if the action is defined via the group action of $\mathfrak{X}$. Specifically, condition (ii) of Definition~\ref{Chap2-Def-DynSyst} does not hold since there are $z_1,z_2\in\mathfrak{X}$ such that $z_1+z_2\in\mathcal{D}$ while $z_2\not\in\mathcal{D}$. The Delone groupoid is an \'etale groupoid with compact unit space. An associated Schr\"odinger operator defined on the Delone graph turns out to be an element of the corresponding reduced groupoid $C^\ast$-algebra. Note that these operators are defined on different Hilbert spaces.

\medskip

Consider now the set $Del_{U,K}^0\subseteq\cs(\mathfrak{X})$ of all Delone sets $\mathcal{D}$ with the same parameters $U,K\subseteq\mathfrak{X}$ and $0\in\mathcal{D}$ where $0$ denotes the neutral element of the abelian group $\mathfrak{X}$. Let $\mathcal{I}\big(Del_{U,K}^0\big)$ be the set of all closed, invariant subsets of $Del_{U,K}^0$. Here the invariance is coming from the groupoid action and the translation by $\mathfrak{X}$, respectively. Then it is natural to ask for the behavior of the spectra of a family of Schr\"odinger operator indexed by $Y\in\mathcal{I}\big(Del_{U,K}^0\big)$ if $Y$ varies. It turns out that the tools of Chapter~\ref{Chap3-SpectAppr} and Chapter~\ref{Chap4-ToolContBehavSpectr} are applicable. More precisely, Theorem~\ref{Chap3-Theo-ContFieldCALgContSpectr} and Theorem~\ref{Chap4-Theo-FieldGroupContCAlgebras} can be used to characterize the continuous behavior of the spectra by the Fell-continuity of $Y\in\mathcal{I}\big(Del_{U,K}^0\big)$, c.f. \cite{BeBeNi16}. Like in the case of dynamical systems, the concept of universal groupoids is the key ingredient.

\medskip

Note that Delone sets are usual identified with the so called Voronoi tiling in the case $\mathfrak{X}=\RM^d$. Conversely, a tiling is identified with a Delone set by puncture each tile with its barycenter. According to Theorem~\ref{Chap2-Theo-Shift+DictSpace} and Remark~\ref{Chap5-Rem-SpaceBruijnGraphs}, the space of the symbolic dynamical system $(\as^G,G,\alpha)$ is homeomorphic to the space of dictionaries equipped with the local pattern topology and the space of de Bruijn graphs $\bs(\as)$ equipped with the product topology. Here the last homeomorphism only exists in the case $G=\ZM$. These connections are strongly used in Chapter~\ref{Chap5-OneDimCase} and Chapter~\ref{Chap6-HigherDimPerAppr}. Thus, it is desirable to extend this notation of a dictionary space to Delone sets in $\RM^d$ in a future project. In this case, the alphabet is replaced by the prototiles, i.e., the equivalence class of the tiles with respect to the translation by elements of $\RM^d$. Note that finite local complexity corresponds in this case to a finite alphabet.

\section{Fell-topology}
\label{Chap8-Sect-FellTop}

In view of Theorem~\ref{Chap4-Theo-CharDynSystConvSpectr}, the Fell-topology on $\SG(X)\subseteq\cs(X)$ for a dynamical system $(X,G,\alpha)$ is the correct topology in terms of continuous behavior of the spectra. According to Proposition~\ref{Chap2-Prop-SpaDynSyst}, the space $\SG(X)$ is a second-countable, compact, Hausdorff space. In the specific case that $X$ is totally disconnected, $\SG(X)$ is totally disconnected as well, c.f. Corollary~\ref{Chap2-Cor-SGTotaDisco}. Apart from this, the topological structure of $\SG(X)$ is not understood very deeply. Some special cases show that the topology depends on the interplay between the topology on $X$ and the dynamics as described next.

\medskip

The space $\cs(X)$ is compact by \cite[Theorem~1]{Fel62} implying that every $Y\in\cs(X)$ can be approximated by a sequence $Y_n\in\cs(X)\,,\; n\in\NM\,,$ of finite subsets of $X$. Consequently, $\cs(X)$ does not have isolated points if $X$ does so. On the other hand, the topological properties of $\SG(X)$ strongly depends on the dynamics induced by $G$. For instance, $\SZ\big(\as^\ZM\big)$ has isolated points like all the periodic subshifts as well as some aperiodic subshifts, c.f. Example~\ref{Chap5-Ex-DeBruijnNotStrongConnect} and Proposition~\ref{Chap5-Prop-PerPointIsolated}. On the other hand, all elements of $\SG(X)$ are non-isolated if there exists an invariant metric with respect to the group action, c.f. Theorem~\ref{Chap2-Theo-FinOrbDenseInvMetr}. The existence of a sequence of (a finite union of) strongly periodic dynamical subsystems converging to a dynamical subsystem $Y\in\SG(X)$ boils down to the question if $Y$ can be approximated by finite, $G$-invariant subsets of $X$, i.e., $Y\subseteq X$ is finite and $Y=\{\alpha_g(y)\;|\; y\in Y\,,\; g\in G\}$. Consequently, the first question is to understand whether $Y\in\SG(X)$ is isolated or not. 

\medskip

Altogether, the combination of the Fell-topology with the dynamics creates interesting topological properties of the space $\SG(X)$ that need to be studied in more detail.

\section{Periodic approximation}
\label{Chap8-Sect-PerApprox}

Periodic approximations are interesting by several reasons. First, the spectra of periodic operators can be analyzed with the Floquet-Bloch theory, c.f. Section~\ref{Chap4-Sect-ApprPerAppr}. Furthermore, some of the known spectral results (in dimension one) explicitly use periodic approxi\-mations to estimate, for instance, the Hausdorff dimension of the spectrum, c.f. the discussion in Chapter~\ref{Chap1-Intro} and Chapter~\ref{Chap7-Examples}. Also the choice of periodic boundary conditions is helpful for numerical calculations where it is necessary to know that the corresponding spectra converge, c.f. \cite{BeSi91}. Additionally, the rate of convergence is faster than the finite box approximation by experiences, see e.g. \cite[Section~2.1.3, Section~4.1]{RiethThesis95}. Philosophically, this is due to the fact that the boundary effects are washed out by the periodic boundary conditions. 

\medskip

According to Theorem~\ref{Chap4-Theo-PeriodicApproximations}, each generalized Schr\"odinger operators associated with a dynamical subsystem $Y$ can be approximated by periodic ones if the dynamical subsystem $Y$ is periodically approximable. Thus, it is interesting to find sufficient conditions for dynamical subsystems being periodically approximable. In addition, a subshift is periodically approximable only if the associated subshifts of finite type contain strongly periodic elements, c.f. Page~\ref{Page-GeneralExStrongPerSubshFinitType} and Page~\pageref{Page-ExStrongPerSubshFinitType} for more details. The existence of strongly periodic elements in a subshift of finite type is a current question in dynamical systems \cite{Pia08,Fio09,Hoc09,SiCo12,Coh14,CaPe15} which is also related to the existence of Wang tiles \cite{Wan61,Ber66,CuKa95}.

\medskip

As shown in Section~\ref{Chap5-Sect-PerApprSubs}, a large class of subshifts of the symbolic dynamical system $(\as^\ZM,\ZM,\alpha)$ is contained in the set of periodically approximable subshifts $\mathpzc{PA}_{\ZM}\big(\as^\ZM\big)$, c.f. Proposition~\ref{Chap5-Prop-MinAllPath}, Proposition~\ref{Chap5-Prop-SuffCondPerAppr} and Corollary~\ref{Chap5-Cor-SuffCondPerAppr}. Such a general result remains open for higher dimensional systems.

\medskip

In case of substitutional subshifts $\Xi\in\SZd\big(\as^{\ZM^d}\big)$, a sufficient condition is provided by Theo\-rem~\ref{Chap6-Theo-PerApprZdPrimSubst}. More precisely, the existence of a strongly periodic element of $\as^{\ZM^d}$ satisfying that all patterns up to a certain size appear in $\Xi$ implies that $\Xi$ is periodically approxi\-mable. The iterative application of the substitution to this strongly periodic el\-ement defines a strongly periodic approximation of $\Xi$. In Theorem~\ref{Chap6-Theo-PerApprZdPrimSubst}, it is assumed that $\Xi$ is defined via a primitive block substitution. The condition of being primitive can be relaxed, c.f. Section~\ref{Chap7-Sect-SierpinskiCarpetSubstitution}. Also the requirement to have a block substitution is not necessary. In detail, it suffices that the substitution extends to a map $S:\as^{\ZM^d}\to\as^{\ZM^d}$ and that $S^n(a)$ growths eventually with $n\in\NM$ in each direction $j\in\{1,\ldots,d\}$, c.f. Remark~\ref{Chap6-Rem-BlockNotNec-d=1}.

\medskip

Altogether, the main task is to show the existence of a strongly periodic element such that all patterns up to a fixed size appear in $\Xi$ whenever $\Xi$ is defined by a substitution. The existence of such an element corresponds to the existence of a torus (``a generalized closed path'') in the associated Anderson-Putnam complex, c.f. Section~\ref{Chap8-Sect-APComplex}. Proposition~\ref{Chap6-Prop-SuffPerApprZ2SymmPatt} provides a sufficient condition for the case $d=2$. Based on this result, the general philosophy is that such a strongly periodic element has constraints on the local symmetry of the patterns associated with the subshift $\Xi$, c.f. Remark~\ref{Chap6-Rem-SuffPerApprZ2SymmPatt}. This meets also the physical intuition that quasicrystals have local symmetries.

\medskip

The result of Theorem~\ref{Chap6-Theo-PerApprZdPrimSubst} seems to extend to substitutions on tilings of $\RM^d$ by following the same philosophy. Such a substitution consists of an inflation map and a subdivision of the inflated tiles, see e.g. \cite{LuPl87,AnPu98,Sol98,Rob99} for background on substitutional tilings in $\RM^d$.

\medskip

Substitutional subshifts and tilings represent a large class of quasicrystals that are studied in the literature. The class of model sets defined via a cut-and-project scheme is also of particular interest. Both classes have some overlap. More precisely, there are model sets that are also defined via a substitution. For instance, the Fibonacci sequence is such an example. Model sets are defined by cutting a strip in a high dimensional space (containing $\mathfrak{X}$) and then project the lattice points contained in the strip to $\mathfrak{X}$. This construction promises to be helpful for proving the existence of strongly periodic approximations since it has an analog of an angle of the corresponding strip. Meanwhile, it is well-known that rational angles induce periodic Delone sets. Thus, it is interesting to study the induced topology of $\mathcal{I}\big(Del_{U,K}^0\big)$ on the set of angles. Summing up, I suspect that Delone sets arising by a cut-and-project scheme provide a class of periodically approximable subshifts which is also a part of future projects. As a first step the relation to Diophantine approximations analyzed in \cite{HaKoWa15,Hay16,HaKo16,HaKoSaWa16} can be studied for the existence of periodic approximations in cut-and-project schemes.

\begin{conjecture}
\label{Chap8-Conj-CutAndProj}
The transversal of a Delone set $\mathcal{D}\in Del_{U,K}^0$ is periodically approximable if  $\mathcal{D}$ is defined via a cut-and-project scheme.
\end{conjecture}

\section{The role of the Anderson-Putnam complex}
\label{Chap8-Sect-APComplex}

Recall that tilings and Delone sets of $\RM^d$ can be treated on the same ground. The class of Delone sets with finite local complexity corresponds to tilings consisting of a finite number of prototiles. Consider a tiling $T$ of $\RM^d$ with finitely many prototiles. Then $\Gamma_0$ is the associated Anderson-Putnam complex (AP-complex) with the tiling $T$, c.f. \cite{AnPu98,BeGa03,BeBeGa06}. In the one-dimensional case, this AP-complex turns out to be a de Bruijn graph associated with a subshift $\Xi\in\SZ\big(\as^\ZM\big)$, c.f. Remark~\ref{Chap5-Rem-GAP}. The AP-complex describes the rule for the prototiles to glue them together. If the AP-complex has a torus, the torus defines a strongly periodic tiling by tiles of $T$. By passing to supertiles, a sequence of AP-complexes $(\Gamma_k)_{k\in\NM}$ is associated with $T$. This corresponds to the sequence of de Bruijn graphs for the one-dimensional situation. In case of a tiling arising by a substitution, the sequence can be defined by iteratively applying the substitution to the periodic tiling corresponding to a torus in the AP-complex, c.f. Remark~\ref{Chap6-Rem-StrategySubstitution} and Chapter~\ref{Chap6-HigherDimPerAppr}. 

\medskip

It is desirable to analyze the existence of tori in the Anderson-Putnam complex and its relation to the existence of periodic approximations for these systems. Specifically, this could lead to connections between the cohomology of the tiling and the property of being periodically approximable.

\section{H\"older-continuous behavior of the spectra}
\label{Chap8-Sect-HolderCont}

Let $(\ts,d)$ be a complete metric space. According to Theorem~\ref{Chap3-Theo-CharHolContSpect}, the spectra of a self-adjoint field of operators $(A_t)_{t\in\ts}$ vary H\"older-continuous in $t\in\ts$ if and only if the norms $\ts\ni t\mapsto\|p(A_t)\|\in\RM\,,\; p\in\Pt(M)\,,$ behave uniform H\"older-continuous for each $M>0$.

\medskip

If only the continuity of the norms is demanded, {\sc N. Landsman} and {\sc B. Ramazan} \cite{LaRa99} provide a general tool to verify the continuous behavior of the spectra in large generality by combining Theorem~\ref{Chap3-Theo-P2ContEquivContSpect} and Theorem~\ref{Chap4-Theo-TransGroupContFieldCAlg}. Such a tool is not available for the H\"older-continuity. In a current work with {\sc J. Bellissard}, we develop such a tool in the same spirit by using general techniques of $C^\ast$-algebras.

\medskip

More precisely, the continuous behavior of the spectra of a generalized Schr\"odinger operator strongly relies on the fact that the corresponding element $\fz$ of the groupoid $C^\ast$-algebra is a continuous function on the groupoid. For indeed, the proof of Theorem~\ref{Chap4-Theo-LanRamContFieldGroupoid} is based on the continuity of $\fz$ and that the norms $\|\fz\|_{red}$ and $\|\fz\|_{full}$ can be represented by an infimum and a supremum of functions of $\fz$ respectively. 

\begin{conjecture}
\label{Chap8-Conj-HolCont}
Let $(\Gamma,\ts,p)$ be a continuous field of groupoids. Consider a metric $d$ on $\ts$ such that $(\ts,d)$ is a complete metric space. Then the norm $\ts\ni t\mapsto\|\fz_t\|\in[0,\infty)$ is H\"older-continuous if $\fz\in\Cc_c\big(\Gamma^{(1)}\big)$ is H\"older-continuous in a suitable sense.
\end{conjecture}

In view of Remark~\ref{Chap3-Rem-CharHolContSpect}, it would not be surprising that there exists in general an $\fz\in\Cc_c\big(\Gamma^{(1)}\big)$ that is $\alpha$-H\"older-continuous while the norm $\ts\ni t\mapsto\|\fz_t\|\in[0,\infty)$ is only $\alpha/2$-H\"older-continuous. Specifically, the rate of convergence might decrease in general. Combining Conjecture~\ref{Chap8-Conj-HolCont} with Theorem~\ref{Chap3-Theo-CharHolContSpect}, this leads to a tool to prove the H\"older-continuous behavior of the spectra associated with $\fz$.

\cleardoublepage


\appendix

\chapter{Topology a reminder}
\label{App1-Topology}

This section provides a short reminder for standard concepts and results in topology and related subjects. The reader is referred to \cite{Kuratowski66,Kuratowski68,Querenburg2001} for a more detailed discussion.

\section{Topological spaces and metric spaces}
\label{App1-Sect-TopologyMetric}

The framework for topology are introduced and the connections with metric spaces are discussed.

\begin{definition}[Topology]
\label{App1-Def-Topology}
Let $X$ be a set and $\tau$ be a family of subsets of $X$. Then $\tau$ is called {\em topology} if 
\begin{description}
\item[(T1)\label{(T1)}] $X,\emptyset\in \tau$;
\item[(T2)\label{(T2)}] the union $\bigcup_{\iota\in I} U_\iota$ is contained in $\tau$ for every $U_\iota\in\tau\,,\; \iota\in I\,,$ and an index set $I$;
\item[(T3)\label{(T3)}] the finite intersection $\bigcap_{j=1}^N U_j$ is contained in $\tau$ for each $N\in\NM$ and $U_j\in\tau\,,\; 1\leq j\leq N$.
\end{description}
\end{definition}

In general, it is convenient to work with a smaller subfamily $\bs$ of a topology $\tau$ that encodes all properties of $\tau$. Specifically, let $\bs$ be a family of subsets of $X$ satisfying
\begin{description}
\item[(B1)\label{(B1)}] the union $\bigcup_{U\in\bs}U$ is equal to $X$;
\item[(B2)\label{(B2)}] for each $U,V\in\bs$, the intersection $U\cap V$ is equal to a union of elements of $\bs$.
\end{description}
Then there exists a unique topology $\tau(\bs)$ on $X$ with $\bs\subseteq\tau(\bs)$ defined by the intersection of all topologies containing $\bs$. In this case, every element of $\tau(\bs)$ is represented by a union of elements of $\bs$.

\begin{definition}[Base for a topology]
\label{App1-Def-BaseTopology}
Let $\bs$ be a family of subsets satisfying \nameref{(B1)} and \nameref{(B2)} and $\tau(\bs)$ be the associated unique topology. Then $\bs$ is called a {\em base for the topology $\tau(\bs)$}. The family $\bs_x:=\{U\in\bs\;|\; x\in U\}$ is called a {\em neighborhood base of $x\in X$}.
\end{definition}

A set $X$ equipped with a topology $\tau$ is called a {\em topological space $(X,\tau)$}. Whenever there is no confusion, the topological space $(X,\tau)$ is denoted by $X$. An element of $\tau$ is called an {\em open set}. A set $F\subseteq X$ is called {\em closed} whenever it is the complement of an open set, i.e., $F=X\setminus O$ for an $O\in\tau$. A subset $V\subseteq X$ is called a {\em neighborhood of a set $F\subseteq X$} if there is an open set $O\in\tau$ such that $F\subseteq O\subseteq V$. Let $Y$ be a subset of a topological space $X$ with topology $\tau$. Then $\tau_Y:=\{Y\cap U\;|\; U\in T\}$ defines a topology on $Y$. This topology is called {\em subspace topology on $Y$ of $X$}. Let $X$ be a set and $\tau_1,\tau_2$ be two topologies on $X$. Then $\tau_1$ is called {\em finer} than $\tau_2$ if, for each open set $U\subseteq X$ in the topology $\tau_2$, the set $U$ is open in the topology $\tau_1$, i.e., $\tau_2\subseteq\tau_1$. 
 
\begin{lemma}
\label{App1-Lem-TopFiner}
Let $X$ be a set and $\tau_1$ and $\tau_2$ be two topologies on $X$ where $\bs_2$ is a base for the topology $\tau_2$. If, for all $x\in X$ and each $V_2(x)\in\bs_2$ of $x$, there exists a $V_1(x)\in\tau_1$ such that $x\in V_1(x)\subseteq V_2(x)$, then the topology $\tau_1$ is finer than the topology $\tau_2$. 
\end{lemma}

\begin{proof}
Let $U\in\tau_2$ be open and $x\in U$. Then there is a $V_2(x)\in\bs_2$ such that $V_2(x)\subseteq U$. Hence, there exists a $V_1(x)\in\tau_1$ satisfying $x\in V_1(x)\subseteq V_2(x)\subseteq U$ by assumption. Consequently, $U$ is equal to the union $\bigcup_{x\in U} V_1(x)$ which is an open set in $\tau_1$. Thus, $U\in\tau_1$ follows.
\end{proof}

\begin{definition}[Convergent sequence]
\label{App1-Def-ConvergenceSequence}
Let $(X,\tau)$ be a topological space with base $\bs$. Then a sequence $(x_n)_{n\in\NM}$ in $X$ is called {\em convergent} to $x\in X$ if, for all $U\in\bs_x$, there exists an $n_0\in\NM$ such that $x_n\in U$ for $n\geq n_0$ where $\bs_x:=\{U\in\bs\;|\; x\in U\}$ is the neighborhood base of $x$ induced by the base $\bs$.
\end{definition}

Let $(X,\tau)$ be a topological space. Then $(X,\tau)$ is called {\em second-countable} if there exists a base $\bs$ for the topology $\ts$ such that $\bs$ is countable. Furthermore, $(X,\tau)$ is called {\em first countable} if every element $x\in X$ has a countable neighborhood base. The topological space $(X,\tau)$ is called {\em Hausdorff} if every pair of elements $x,y\in X$ with $x\neq y$ can be separated by open sets, i.e., there exist a neighborhood $U_x$ of $x$ and a neighborhood $U_y$ of $y$ such that $U_x\cap U_y=\emptyset$. Clearly, each singleton $\{x\}\subseteq X$ is closed if the topology on $X$ is Hausdorff since $X\setminus\{x\}$ is equal to the union $\bigcup_{y\neq x} U_y$ where $U_y$ is an open neighborhood of $y$ such that $x\not\in U_y$.

\begin{definition}[Normal space]
\label{App1-Def-NormalSpace}
A topological space $(X,\tau)$ is called {\em normal} if closed subsets can be separated by open sets, i.e., for each pair of disjoint, closed sets $F_1,F_2\subseteq X$, there are neighborhoods $U_1$ of $F_1$ and $U_2$ of $F_2$ such that $U_1\cap U_2=\emptyset$. \end{definition}

A subset $K\subseteq X$ of a topological space $(X,\tau)$ is called {\em compact} if every open cover of $K$ has a finite subcover of $K$. The topological space $(X,\tau)$ is called {\em locally compact} if each element $x\in X$ has a base of compact neighborhood, i.e., there exists a compact neighborhood $K_x\subseteq X$ such that $K_x\subseteq U_x$ for every open neighborhood $U_x$ of $x$.

\begin{definition}[Complete metric space]
\label{App1-Def-Metric}
Let $X$ be a set. A map $d:X\times X\to[0,\infty)$ is called a metric if, for all $x,y,z\in X$, the relations
\begin{description}
\item[(i)] $x=y\; \Leftrightarrow\; d(x,y)=0$,
\item[(ii)] $d(x,y)=d(y,x)$,
\item[(iii)] $d(x,y)\leq d(x,z) + d(z,y)$,
\end{description}
hold. A sequence $(x_n)_{n\in\NM}$ in $(X,d)$ is called a {\em Cauchy sequence} if, for each $\varepsilon>0$, there exists an $n(\varepsilon)>0$ such that $d(x_n,x_m)<\varepsilon$ for all $n,m\geq n(\varepsilon)$. Then $(X,d)$ is called a {\em complete metric space} if every Cauchy sequence converges in $X$ with respect to $d$.
\end{definition}

Let $(X,d)$ be a complete metric space. The family $\bs_d:=\{B_r(x)\;|\; x\in X\,,\; r>0\}$ satisfies \nameref{(B1)} and \nameref{(B2)} and so it defines a unique topology $\tau_d$ on $X$. This topology is called the {\em induced topology on $X$ by the metric $d$}.

\medskip

A topological space $(X,\tau)$ is called metrizable if there exists a metric $d:X\times X\to[0,\infty)$ such that $(X,d)$ is a complete metric space and the induced topology $\tau_d$ by the metric $d$ is equal to $\tau$.

\begin{proposition}
\label{App1-Prop-LocCompHausSecCountImplNormal}
Let $(X,\tau)$ be a second countable, locally compact, Hausdorff space. Then $(X,\tau)$ is metrizable and so $X$ is a normal space.
\end{proposition}

\begin{proof}
According to \cite[Satz 10.15]{Querenburg2001}, a locally compact, Hausdorff space $X$ is second countable if and only if $X$ is metrizable and countable at infinity. Then \cite[Satz 10.9]{Querenburg2001} and \cite[Satz 10.2]{Querenburg2001} imply that every metrizable space $X$ is normal.
\end{proof}

\begin{remark}
\label{App1-Rem-metrizable}
Note that a metric on a set is a more detailed concept than a topology. Specifically, a metrizable space $(X,\tau)$ can have different metrics that define the same topology. These metrics are not required to be equivalent at all. There might exists also a metric so that $(X,d)$ does not define a complete metric space while the metric $d$ induces the topology $\tau$.
\end{remark}

\section{Continuous functions}
\label{App1-Sect-ContFunct}

Certain properties of functions between topological spaces are studied in the following. Furthermore, it is proven that suitable maps between topological spaces transfer topologi\-cal properties. 

\medskip

Let $(X,\tau_X)$ and $(Y,\tau_Y)$ be topological spaces. A map $f:X\to Y$ is called {\em continuous} whenever the preimage $f^{-1}(V)$ is open in $(X,\tau_X)$ for each $V\in\tau_Y$. Furthermore, $f:X\to Y$ is called {\em open} if the image $f(U)$ is an open set in $(Y,\tau_Y)$ for each $U\in\tau_X$. Similarly, a map $f:X\to Y$ is called {\em closed} whenever $f(F)\subseteq Y$ is closed for every closed subset $F\subseteq X$. 

\begin{proposition}
\label{App1-Prop-ContOpenSurjMapPresSecCountable}
Let $(X,\tau_X)$ and $(Y,\tau_Y)$ be topological spaces. Let $f:X\to Y$ be a surjective, open and continuous map.
\begin{itemize}
\item[(i)] Then $Y$ is second-countable if $X$ is second-countable.
\item[(ii)] Then $Y$ is locally compact if $X$ is locally compact.
\end{itemize}
\end{proposition}

\begin{proof}
(i): Let $f:X\to Y$ be an open, continuous and surjective map where $X$ is second-countable. Consider a countable base $\bs_X=\{V_n\in\tau_X\;|\;n\in\NM\}$ of the topology $\tau_X$. It suffices to prove that $\bs_Y:=\{f(V_n)\;|\; n\in\NM\}$ is a base for the topology $\tau_Y$ on $Y$. 

\vspace{.1cm}

Then $\bs_Y\subseteq\tau_Y$ follows since $f$ is an open map. Assertion~\nameref{(B1)} is satisfied since $f$ is surjective and $\bs_X$ is a base for the topology $\tau_X$. Let $U,V\in\bs_Y$. Thus, there are $n_0,n_1\in\NM$ such that $f(V_{n_0})=U$ and $f(V_{n_1})=V$. Consider their intersection $U\cap V$ which is open by Definition~\ref{App1-Def-BaseTopology}~(iii). The map $f$ is surjective and so it has a right inverse. More precisely, the equality $f\big(f^{-1}(U\cap V)\big) = U\cap V$ is deduced. For each $y\in O:=f^{-1}(U\cap V)$, there is a $V(y)\in\bs_X$ satisfying $y\in V(y)\subseteq O$. Hence, the identity $O = \bigcup_{y\in O} V(y)$ follows. Consequently,
$$
U\cap V \; 
	= \; f\big(\bigcup_{y\in O} V(y)\big)
	= \; \bigcup_{y\in O} f\big(V(y)\big)
$$
where the second inequality is an elementary statement. This leads to the desired property \nameref{(B2)}. Thus, $\bs_Y$ defines a countable base for the topology $\tau(\bs_Y)$.

\vspace{.1cm}

It is left to prove that $\tau(\bs_Y)=\tau_Y$ while $\tau(\bs_Y)\subseteq\tau_Y$ holds since $f$ is an open map. For the converse inclusion, let $U\in\tau_Y$. Then $f^{-1}(U)$ is open by continuity of $f$. Since $\bs_X$ is a base for $\tau_X$, $f^{-1}(U)$ is equal to the union of elements of $\bs$. Hence, there is an index set $I\subseteq\NM$ and there are $U_k\in\bs_X$ for $k\in I$ such that $f^{-1}(U)=\bigcup_{k\in I} U_k$ holds. Then $f\big(f^{-1}(U)\big) = \bigcup_{k\in I} f(V_k)\,,\; f(V_k)\in\bs_Y\,,$ is derived. The equation $f\big(f^{-1}(U)\big)=U$ holds by surjectivity of $f$. Hence, $U$ is equal to the union of elements of $\bs_Y$ implying $U\in\tau(\bs_Y)$.

\vspace{.1cm}

(ii): Let $y\in Y$ and $U_y$ be an open neighborhood of $y$. Then there is an $x\in X$ such that $f(x)=y$ since $f$ is surjective. Furthermore, $f^{-1}(U_y)$ is an open neighborhood of $x$ as $f$ is continuous. Due to the fact that $X$ is locally compact, there exists a compact neighborhood $K_x\subseteq f^{-1}(U_y)$ of $x$. Then $f(K_x)$ is a compact neighborhood of $f(x)$ since $f$ is continuous and open. The image $f(K_x)$ is contained in $U_y$ by construction. Altogether, $Y$ is locally compact.
\end{proof}

\medskip

A map $f:X\to Y$ between two topological spaces is called {\em proper} if the preimage of a compact set is compact.

\begin{proposition}
\label{App1-Prop-ContClosMapIsProper}
A continuous, closed map $f:X\to Y$ is proper if the preimage $f^{-1}(\{y\})$ is compact for every $y\in Y$.
\end{proposition}

\begin{proof}
Let $f: X \to Y$ be a continuous, closed map such that $f^{-1}(\{y\})\subseteq X$ is compact for every $y \in Y$. It has to be shown that $f^{-1}(K)$ is compact for each compact $K\subseteq Y$. Let $K$ be a compact subset of $Y$ and $\us:=\{ U_{\iota} \;|\; \iota\in I \}$ be an open cover of $f^{-1}(K)$. 

\vspace{.1cm}

The family $\us$ is also an open cover of $f^{-1}(\{y\})$ for each $y\in K$. By assumption, there exist, for $y\in K$, a finite $I_y\subseteq I$ satisfying $f^{-1}(\{y\}) \subseteq \bigcup_{\iota\in I_y} U_{\iota}$. Then the set $X \setminus \bigcup_{\iota\in I_y} U_{\iota}$ is closed implying that its image through $f$ is closed in $Y$ since $f$ is a closed map. Consequently,
$$
V_y \;
	:= \; Y \setminus f\left(X \setminus \bigcup_{\iota\in I_y} U_{\iota}\right) 
$$
is open in $Y$ and contains $y$. Thus, the inclusion $K \subseteq \bigcup_{y \in K} V_y$ follows. Since $K$ is compact, there exist $y_1,\dots, y_n\in K$ satisfying $K \subseteq \bigcup_{j =1}^n V_{y_j}$. Due to construction, the union $I_K:=\bigcup_{j=1}^n I_{y_j}$ is finite as a finite union of finite sets. Thus, the inclusions
$$
f^{-1}(K) \;
	\subseteq \; f^{-1}\left(\bigcup_{j=1}^n V_{y_j}\right) \;
	\subseteq \; \bigcup_{j=1}^n f^{-1}\left(V_{y_j}\right) \;
	\overset{(\blacklozenge)}{\subseteq} \bigcup_{\iota\in I_K} U_{\iota}
$$ 
are concluded. Hence, there is a finite subcover of $\us$ for $f^{-1}(K)$. It is left to verify the inclusion $(\blacklozenge)$ which is proved as follows. 

\vspace{.1cm}

It suffices to show $f^{-1}(V_{y_j})\subseteq\bigcup_{\iota\in I_{y_j}} U_{\iota}$. For a set $F\subseteq X$, the inclusion $F\subseteq f^{-1}(f(F))$ holds implying $X\setminus f^{-1}(f(F))\subseteq X\setminus F$. Consequently, the inclusions
$$
f^{-1}\big(V_{y_j}\big) \; 
	= \; f^{-1}\left( Y\setminus f\left(X\setminus \bigcup_{\iota\in I_{y_j}} U_{\iota}\right) \right) \;
	\subseteq\; X\setminus f^{-1}\left(f\left(X\setminus \bigcup_{\iota\in I_{y_j}} U_{\iota}\right)\right) \;
	\subseteq \; \bigcup_{\iota\in I_{y_j}} U_{\iota}
$$
are derived proving $(\blacklozenge)$.
\end{proof}

\section{Urysohn's Lemma and Tietze's Extension Theorem}
\label{App1-Sect-UrysohnTietze}

In this section, the Lemma of Urysohn and Tietze's Extension Theorem are presented.

\medskip

Urysohn's Lemma states that the topology of a normal space can be described by the supports of continuous functions. This fact is used in Section~\ref{Chap4-Sect-ContSpecGenSchrOp}, c.f. Figure~\ref{Chap4-Fig-Concept} on Page~\pageref{Chap4-Fig-Concept}. The {\em support} $\supp(\phi)$ of a function $\phi:X\to\CM$ is defined by the closure $\overline{\{ x\in X\;|\; \phi(x)\neq 0 \}}$. 

\begin{proposition}[\cite{Ury25}, Urysohn's Lemma]
\label{App1-Prop-LemmaUrysohn}
Let $X$ be a second countable, locally compact, Hausdorff space. Consider closed subsets $A,B\subseteq X$. Then there exists a continuous map $\phi:X\to[0,1]$ such that $\phi|_A\equiv 1$ and $\phi|_B\equiv 0$. The function $\phi$ can be chosen with compact support if, additionally, $A$ is compact.
\end{proposition}

\begin{proof}
Thanks to Proposition~\ref{App1-Prop-LocCompHausSecCountImplNormal}, the space $X$ is normal. Then the existence of a continuous function $\phi:X\to[0,1]$ with $\phi|_A\equiv 1$ and $\phi|_B\equiv 0$ is proven in \cite[Satz, Anhang~III.25.]{Ury25}.

\vspace{.1cm}

Suppose now that the set $A$ is compact. Then, for each $x\in A$, there exists an open neighborhood $U_x$ of $x$ such that $\overline{U_x}$ is compact and $\overline{U_x}\cap B=\emptyset$. Hence, there are $x_1,\ldots,x_n\in A$ satisfying
$$
A\subseteq \underbrace{\bigcup_{j=1}^n U_{x_j}}_{=:U} \subseteq \underbrace{\bigcup_{j=1}^n \overline{U_{x_j}}}_{=:K} 
$$
by the compactness of $A$. The set $K$ is compact since it is the finite union of compact sets. Due to \cite[Satz, Anhang~III.25.]{Ury25}, there is a continuous function $\phi:X\to[0,1]$ such that $\phi|_A\equiv 1$ and $\phi|_{X\setminus U}\equiv 0$. Thus, the support of $\phi$ is contained in $U\subseteq K$. Hence, $\phi$ has compact support. Since $B\subseteq X\setminus U$ holds, the desired result follows.
\end{proof}

\medskip

Tietze's Extension Theorem asserts the existence of a continuous extensions of a continu\-ous function on a closed subset of a topological space. This is used in Chapter~\ref{Chap4-ToolContBehavSpectr}.

\begin{theorem}[\cite{Tie15}, Tietze's Extension Theorem]
\label{App1-Theo-TietzesThm}
Let $X$ be a normal space. Consider a continuous function $\fz:K\to\CM$ defined on a closed subset $K\subseteq X$. Then there exists a continuous $\gz:X\to\CM$ such that the restriction $\gz|_K$ is equal to $\fz$. If, additionally, $K$ is compact, the function $\gz$ can be chosen with compact support.
\end{theorem}

\begin{proof}
The existence of a continuous $\gz:X\to\CM$ satisfying $\gz|_K=\fz$ is proven in \cite[Satz~I.8.5.4]{Schubert75} by using that $\CM$ is homeomorphic to $\RM^2$. 

\medskip

Suppose now that $K\subseteq X$ is compact. Due to Proposition~\ref{App1-Prop-LemmaUrysohn}, there exists a $\phi:X\to[0,1]$ with compact support such that $\phi|_K\equiv 1$ since $X$ is normal. Then the function $\gz\cdot\phi:X\to\CM$ defined by pointwise multiplication is continuous with compact support $\supp(\gz\cdot\phi)=\supp(\phi)\cap\supp(\gz)$ and $\gz|_K=\fz$.
\end{proof}

\cleardoublepage
\markboth{Ehrenw\"ortliche Erkl\"arung}{Lebenslauf}
  \textbf{\large Pers\"ohnliche Daten}

\medskip

\begin{tabular}{ll}
	Name: & 
		Siegfried Beckus\\
	Geburtsdatum: & 
		11. April 1988\\
	Geburtsort: & 
		Erfurt\\
	Staatsangeh\"origkeit: & 
		deutsch\\
	Email: & 
		siegfried.beckus@uni-jena.de
\end{tabular}

\vspace{1cm}

\textbf{\large Beruflicher Werdegang}

\medskip

\begin{tabular}{ll}
	seit 03/2012: & 
		wissenschaftlicher Mitarbeiter am Institut f\"ur Mathematik,\\
	&
		Friedrich-Schiller-Universit\"at Jena\\
	02/2014 -- 04/2014 &
		Research collaborator (Student internship J1)\\
	&
		Advisor: Jean Bellissard,\\
	&
		Georgia Institute of Technology, Atlanta, USA
\end{tabular}

\vspace{1cm}

\textbf{\large Bildungsweg}

\medskip

\begin{tabular}{ll}
	seit 03/2012: & 
		Promotionsstudium,\\
	&
		Friedrich-Schiller-Universit\"at Jena,\\
	&
		Betreuer: Prof. Dr. Daniel H. Lenz\\
	08/2009 -- 12/2009&
		Austauschsemester,\\
	&
		Faculty of Science, University of Joensuu, Finland\\
	10/2006 -- 02/2012: &
		Studium der Mathematik mit Nebenfach Physik,\\
	&
		Friedrich-Schiller-Universit\"at Jena,\\
	&
		Diplom 02/2012,\\
	&
		Betreuer: Prof. Dr. Daniel H. Lenz\\
	1998 -- 2006: &
		Heinrich-Hertz Gymnasium,\\
	& 
		Abitur 2006	
\end{tabular}

\vspace{4cm}

\setlength{\parindent}{0mm} Jena, am 6. Oktober 2016\hfill Siegfried Beckus
  \relax 

Hiermit erkl\"are ich, dass
\vspace{0.5cm}

\begin{itemize}
\item[-] 
mir die Promotionsordnung der Fakult\"at f\"ur Mathematik und Informatik der Fried\-rich-Schiller-Universit\"at Jena bekannt ist,
\item[-] 
ich die Dissertation selbst angefertigt habe, keine Textabschnitte oder Ergebnisse eines Dritten oder eigene Pr\"ufungsarbeiten ohne Kennzeichnung \"ubernommen und alle von mir benutzten
Hilfsmittel, pers\"onliche Mitteilungen und Quellen in meiner Arbeit
angegeben habe,
\item[-] 
ich die Hilfe eines Promotionsberaters nicht in Anspruch genommen habe und dass
Dritte weder unmittelbar noch mittelbar geldwerte Leistungen von mir
f\"ur Arbeiten erhalten haben, die im Zusammenhang mit dem Inhalt der
vorgelegten Dissertation stehen,
\item[-]
ich die Dissertation noch nicht als Pr\"ufungsarbeit f\"ur eine staatliche oder andere wissenschaftliche Pr\"ufung eingereicht habe,
\item[-]
ich die gleiche, eine in wesentlichen Teilen \"ahnliche oder eine andere
Abhandlung nicht bei einer anderen Hochschule als Dissertation eingereicht
habe.
\end{itemize}
\vspace{3cm}

\setlength{\parindent}{0mm} Jena, am 6. Oktober 2016\hfill Siegfried Beckus, Verfasser
\cleardoublepage

\printglossary[title=List of symbols,toctitle=List of symbols,type=symbolslist, style=super3colheader]

\bibliographystyle{amsalphasort}
\bibliography{references}

\end{document}